\documentclass[11pt,a4paper]{amsart}
\IfFileExists{lmodern.sty}{\usepackage[T1]{fontenc}\usepackage{lmodern}}{}
\usepackage{microtype,mathtools,amssymb,mathrsfs,booktabs,array,enumitem,longtable,needspace}
\microtypesetup{expansion=false}
\usepackage[textwidth=158mm,textheight=243mm,centering]{geometry}
\usepackage[hidelinks]{hyperref}
\hypersetup{pdftitle={Vector balancing in convex order},
 pdfauthor={Eren Ercan},
 pdfsubject={Convex-order realization, Schatten discrepancy, trace-class quotients, martingale transport and Blackwell privacy},
 pdfkeywords={Komlos, Beck--Fiala, Schatten, Reis-Rothvoss, volume ratio, Talagrand, convex order, q-Bass, martingale transport, Blackwell, privacy}}
\allowdisplaybreaks[2]
\numberwithin{equation}{section}
\newtheorem{theorem}{Theorem}[section]
\newtheorem{lemma}[theorem]{Lemma}
\newtheorem{proposition}[theorem]{Proposition}
\newtheorem{corollary}[theorem]{Corollary}
\theoremstyle{definition}
\newtheorem{definition}[theorem]{Definition}
\theoremstyle{remark}
\newtheorem{remark}[theorem]{Remark}
\newtheorem{example}[theorem]{Example}
\newcommand{\R}{\mathbb R}\newcommand{\Z}{\mathbb Z}
\newcommand{\E}{\mathbb E}\newcommand{\Prb}{\mathbb P}
\newcommand{\Pp}{\mathbb P}
\newcommand{\Cov}{\operatorname{Cov}}\newcommand{\Var}{\operatorname{Var}}
\newcommand{\Ent}{\operatorname{Ent}}\newcommand{\KL}{D}
\newcommand{\tr}{\operatorname{tr}}\newcommand{\diag}{\operatorname{diag}}
\newcommand{\Lip}{\operatorname{Lip}}\newcommand{\supp}{\operatorname{supp}}
\newcommand{\conv}{\operatorname{conv}}\newcommand{\law}{\mathcal L}
\newcommand{\esssup}{\operatorname*{ess\,sup}}
\newcommand{\cx}{\preceq_{\mathrm{cx}}}
\newcommand{\norm}[1]{\lVert#1\rVert}
\newcommand{\ip}[2]{\langle#1,#2\rangle}\newcommand{\dd}{\,d}
\newcommand{\calS}{\mathcal S}

\newcommand{\G}{\mathcal G}
\newcommand{\trans}{^{\mathsf T}}\newcommand{\ri}{\mathfrak i}
\newcommand{\ind}{\mathbf1}\newcommand{\id}{\mathrm{id}}

\DeclareMathOperator{\dist}{dist}
\DeclareMathOperator{\TV}{TV}
\newcommand{\arxiv}[1]{\href{https://arxiv.org/abs/#1}{arXiv:#1}}

\newcommand{\Cstar}{6.8383231851485}
\newcommand{\gammaprof}{\frac{167}{200}}
\newcommand{\Sh}{\mathsf H}

\newcommand{\one}{\mathbf1}
\newcommand{\Diag}{\operatorname{Diag}}
\newcommand{\1}{\mathbf1}

\newcommand{\cP}{\mathcal P}

\newcommand{\ran}{\operatorname{ran}}
\newcommand{\lin}{\operatorname{span}}
\newcommand{\rank}{\operatorname{rank}}\newcommand{\cprank}{\operatorname{cp\text{-}rank}}
\newcommand{\argmax}{\operatorname*{arg\,max}}
\newcommand{\Law}{\mathcal L}
\newcommand{\sym}{\operatorname{sym}}

\newcommand{\Bern}{\operatorname{Bern}}\newcommand{\Unif}{\operatorname{Unif}}

\newcommand{\sgn}{\operatorname{sgn}}
\newcommand{\cF}{\mathcal F}

\newcommand{\N}{\mathbb N}\newcommand{\C}{\mathbb C}

\newcommand{\eps}{\varepsilon}
\DeclareMathOperator{\MMD}{MMD}

\newcommand{\PP}{\mathbb P}
\newcommand{\KLD}[2]{D\!\left(#1\,\middle\|\,#2\right)}
\newcommand{\Pois}{\operatorname{Pois}}\newcommand{\Bin}{\operatorname{Bin}}\newcommand{\Ber}{\operatorname{Bern}}
\newcommand{\St}{\operatorname{St}}\newcommand{\Per}{\operatorname{Per}}\newcommand{\diver}{\operatorname{div}}
\newcommand{\M}{\mathcal M}\newcommand{\MM}{\mathcal M_{\mathrm M}}
\newcommand{\SM}{\mathcal S}\newcommand{\SMm}{\mathcal S_{\mathrm M}}

\newcommand{\cd}{\preceq_{\mathrm{cd}}}\newcommand{\st}{\preceq_{\mathrm{st}}}

\DeclareMathOperator{\aff}{aff}\DeclareMathOperator{\diam}{diam}
\DeclareMathOperator{\intt}{int}
\newcommand{\Lp}{\mathcal L}\newcommand{\Up}{\mathcal U}
\newcommand{\Ga}{\gamma}\newcommand{\Pol}{\mathscr P}\newcommand{\op}{\mathrm{op}}

\DeclareMathOperator{\MCov}{MCov}\DeclareMathOperator{\bary}{bary}

\DeclareMathOperator{\prox}{prox}

\DeclareMathOperator{\dom}{dom}

\DeclareMathOperator{\Sym}{Sym}
\DeclareMathOperator{\vr}{vr}
\DeclareMathOperator{\ovr}{ovr}
\newcommand{\spnorm}[2]{\lVert #1\rVert_{S_{#2}}}
\newcommand{\kfnorm}[2]{\lVert #1\rVert_{(#2)}}

\title{Vector balancing in convex order}
\author{Eren Ercan}
\email{eren321@gmail.com}
\date{28 September 2026}
\begin{document}
\begin{abstract}
We construct Koml\'os signing laws with discrepancy below $6.84$, independent Gaussian reference blocks and exponentially many balanced signings. One law preserves hard constraints and exact conditional means while its reference controls every joint convex cost.

We resolve both Reis--Rothvoss Schatten conjectures. For $n$ symmetric $n\times n$ matrices, one prescribed-mean signing law bounds all centered Schatten discrepancies with sharp powers of $n$ and a universal constant. One Gaussian event of measure $e^{-O(n)}$ controls every unitarily invariant norm at the square-function radius. We prove Reis's volume-ratio conjecture: every $r$-dimensional quotient of the $a\times b$ trace-class unit ball has volume ratio $O(\sqrt{1+\min(a,b)/r})$.

We resolve Nutz--Wang--Zhang's directional and supermartingale conjectures and the higher-dimensional Bass-flow conjecture. For square-integrable laws, full convex reference support is the universal finite-potential criterion in Tschiderer's $q$-Bass problem. Every driver, even atomic, admits one extended potential and shift for all maximal-covariance martingale optimizers of each irreducible pair with full terminal affine span.

For Talagrand's problem, five copies of a compact balanced Gaussian half-mass set contain compact convex output of mass at least $3/4$, also in Wiener space. Cheeger's bound is asymptotically exact throughout subexponential Dirichlet spectral ranges of mixed weighted products.

For standard Borel sources, a Blackwell-greatest perfectly private full-data release exists exactly when, after discarding conditional point masses, the rest coincide or share a two-point support. Under conditional atomlessness, existence is equivalent to independence. Covariance-matched analytic uniformly log-concave sources can have a greatest useful-data-only private release with infinite information, while no greatest full-data release exists.
\end{abstract}
\maketitle
\clearpage

\begingroup
\makeatletter\let\addcontentsline\@gobblethree\makeatother
\section{Introduction}
\endgroup
\pdfbookmark[1]{1. Introduction}{R27-introduction}
A collection of vectors can have many balanced signings with very different statistical behavior. Selecting a probability law on those signings makes it possible to require exact means and to control objectives involving both the signs and their signed sums. The support condition must hold for every output, while the distributional guarantee must hold for the same law under all subsequent convex objectives.

We construct such laws. The initial reference contains independent physical and coefficient variables. Coordinate rearrangements select a discrete output whose expectation under every convex function is bounded by the corresponding expectation of the complete reference. Its bounded physical support gives hard discrepancy; its Gaussian comparison and information cost govern the distribution of the feasible outputs. The same construction can preserve the marginals of a sampling scheme and the exact constraints of an integer rounding problem.

For integrable random vectors, $U\cx V$ means
$\E F(U)\le\E F(V)$ for every finite convex function $F$, allowing infinite expectations on the right. Equivalently, there is a coupling with $\E[V\mid U]=U$ \cite{Strassen}. We use natural logarithms, write $G_n$ for a standard Gaussian vector, and call a sign law symmetric when it is invariant under multiplication by $-1$.
\begin{theorem}[Hard vector balancing in joint convex order]
\label{VB-thm:main}
For every $A\in\R^{m\times n}$ with Euclidean column norms at most one,
where $m,n\ge1$, there is a symmetric random vector
$\sigma\in\{-1,1\}^n$ such that
\begin{equation}\label{VB-eq:main-hard}
 \|A\sigma\|_\infty<\Cstar\qquad\text{at every outcome},
\end{equation}
and
\begin{equation}\label{VB-eq:main-joint}
 (A\sigma,\sigma)\cx(G,H),\qquad
 G\sim N(0,6.5I_m),\quad H\sim N(0,3.23I_n),\quad G\perp H.
\end{equation}
The same law satisfies
\begin{equation}\label{VB-eq:main-entropy}
 \Sh(\sigma)>n\log(1.22946),\qquad
 \Sh(\sigma)\ge n\log2-0.973141\|A\|_F^2.
\end{equation}
\end{theorem}
Guo--Fang--Lu~\cite{GFL} proved the Koml\'os conjecture and its square-root Beck--Fiala consequence. Theorem~\ref{VB-thm:main} strengthens the explicit radius and constructs a whole law with the joint Gaussian comparison and both entropy guarantees. These distributional conclusions persist under every subsequent convex objective.

\noindent\textit{Reading routes.}
The construction and sharp limits are in
\hyperref[R11-sec:construction]{Sections~\ref*{R11-sec:construction}--\ref*{R11-sec:joint}}.
\hyperref[R6-part:applications]{Part~\ref*{R6-part:applications}} gives exact-marginal sampling and polynomial-bit rounding; Sections~\ref{R27-sec:Schatten} and~\ref{R27-sec:volume} prove the Schatten and volume-ratio conjectures.
Martingale stability starts at \hyperref[R20-thm:intro-stability]{Theorem~\ref*{R20-thm:intro-stability}}, with canonical transport and $q$-Bass theory in
\hyperref[R11-sec:canonical]{Sections~\ref*{R11-sec:canonical}--\ref*{R12-sec:bass}}.
For Talagrand convexity and Cheeger spectra, see \hyperref[R19-thm:intro-Talagrand]{Theorems~\ref*{R19-thm:intro-Talagrand}} and \hyperref[R22-thm:intro-fivefold]{\ref*{R22-thm:intro-fivefold}}, and \hyperref[R11-sec:spectra]{Section~\ref*{R11-sec:spectra}}.
\hyperref[R6-part:information]{Parts~\ref*{R6-part:information}--\ref*{R11-part:experiments}} treat entropy and experiments, including the privacy discontinuities in \hyperref[R19-eq:intro-private-budget]{(\ref*{R19-eq:intro-private-budget})} and \hyperref[R22-thm:intro-privacy]{Theorem~\ref*{R22-thm:intro-privacy}}, and the full-data classification in Section~\ref{R27-sec:full-data}. The \hyperref[R11-tab:construction]{named-results guide} gives each conclusion its own proof route.

The entropy bound gives more than $1.22946^n$ distinct hard-balanced signings under the law of Theorem~\ref{VB-thm:main}. Theorem~\ref{MAIN-signing} sharpens the radius columnwise to
\begin{equation}\label{R19-eq:intro-radius}
 R_A=\max_{a_j\ne0}\left\{C\|a_j\|_2-\frac{167}{200}
                  \frac{\sum_i|a_{ij}|^3}{\|a_j\|_2^2}\right\},
 \qquad C=\Cstar,
\end{equation}
with $\|A\sigma\|_\infty<R_A$ when $A\ne0$. In a set system where each point belongs to at most $t\ge1$ sets, this gives the Beck--Fiala bound
\[
 \operatorname{disc}(A)<C\sqrt t-\frac{167}{200}.
\]
The same construction retains the joint reference and counting bound after the corresponding rescaling. Theorem~\ref{ENERGY-thm:main} also retains the physical Cram\'er cost in the weighted signing count.

\Needspace{10\baselineskip}
\noindent\textit{Explicit existential comparison, 27 September 2026.}
All bounds below concern arbitrary matrices with Euclidean column norms at most one.
\begin{center}
\begin{tabular}{@{}ll@{}}
\toprule
Construction & Universal discrepancy bound \\\midrule
Guo--Fang--Lu~\cite{GFL} & $3\sqrt{2\pi}=7.519884\ldots$ \\
Akbas--Sra~\cite[Appendix A]{R19-AkbasSraExposition} & $C_0=7.51500306616\ldots$ \\
Karingula--Lovett, version 2~\cite{KL} & $36$ \\
Theorem~\ref{MAIN-signing} & $C=\Cstar$, improved columnwise by \eqref{R19-eq:intro-radius} \\
\bottomrule
\end{tabular}
\end{center}
Kunisky's examples require a universal constant at least $1+\sqrt2$~\cite{R19-Kunisky}. The computational comparison below separates algorithmic discrepancy constants from the polynomial-bit rounding results.

After selecting the law, choose any finite family of vectors $(u_j,v_j)$ and thresholds $c_j$. Then
\begin{equation}\label{VB-eq:later-loss}
 \E\max_j(u_j^{\mathsf T}A\sigma+v_j^{\mathsf T}\sigma-c_j)_+
 \le\E\max_j(u_j^{\mathsf T}G+v_j^{\mathsf T}H-c_j)_+.
\end{equation}
Every sample retains the hard constraint. The query is chosen after constructing the distribution and before drawing the sample. The sampling section gives separate information bounds for sample-dependent queries.

\subsection{The mechanism and its necessary cost}
The proof retains the identity
\begin{equation}\label{eq:opening-regression}
 \E[(X,T)\mid\sigma]=(A\sigma,\sigma)
\end{equation}
with the prescribed independent law of $(X,T)$, at the critical overlap height for existence. Keeping a convex cost's epigraph fixed gives a feasible signing whose cost is at most the reference expectation. If a finite list admitted no common law, separation would produce a nonnegative combination violating this same inequality. Compactness then gives one law for every convex cost (Section~\ref{R11-sec:construction}).

At the critical height, the vertical coordinate absorbs the cost without changing the horizontal height assumptions. A uniform interval of length $2L$ above the graph adds mean height $L$ and has vertical variation $1/L$. The section estimate recovers the entire added height, and the rearranged horizontal means converge inside the same fixed epigraph. Thus the limiting selection retains the convex inequality at the existence threshold.

For independent Gaussian reference blocks with covariances $V,D>0$, the graph relation between $\sigma$ and $A\sigma$ lets us minimize over all representations of a query:
\begin{equation}\label{eq:opening-minimum}
 \min_u\{u^{\mathsf T}Vu+(t-A^{\mathsf T}u)^{\mathsf T}D(t-A^{\mathsf T}u)\}
 =t^{\mathsf T}(D^{-1}+A^{\mathsf T}V^{-1}A)^{-1}t.
\end{equation}
The vector comparison therefore gives
\begin{equation}\label{VB-eq:main-precision}
 \sigma\cx N\!\left(0,[3.23^{-1}I_n+6.5^{-1}A^{\mathsf T}A]^{-1}\right).
\end{equation}
In sampling, the same minimization charges the independent coefficient reference for the part of an integrand outside its chosen features.

Let $U_n$ denote the uniform law on $\{-1,1\}^n$. The least possible scalar Gaussian scale is
$\kappa=\sqrt{\pi/2}$. At that endpoint, a martingale coupling gives
\[
 1=\E\bigl|\E[\kappa G_i\mid\sigma]\bigr|=\E|\kappa G_i|.
\]
Equality in conditional Jensen forces $\sigma_i=\sgn G_i$ for every coordinate. Hence the whole endpoint law is the independent uniform sign law, which assigns positive probability to the extreme sum of the row $(1,\ldots,1)$. The following theorem measures the cost of approaching this endpoint while imposing a uniform hard discrepancy bound.
\Needspace{15\baselineskip}
\begin{theorem}[Optimal discrepancy--Gaussian tradeoff]
\label{VB-thm:tradeoff}
For $\varepsilon>0$, let $\mathcal R_G(\varepsilon)$ be the infimum
of the radii $R$ such that, for every $m,n\ge1$ and every
$A\in\R^{m\times n}$ with Euclidean column norms at most one,
there is a law on $\{-1,1\}^n$ satisfying
\[
 \|A\sigma\|_\infty\le R\quad\text{almost surely},\qquad
 \sigma\cx(1+\varepsilon)\kappa G_n.
\]
Then
\begin{equation}\label{VB-eq:tradeoff-rate}
 \mathcal R_G(\varepsilon)=\Theta(\varepsilon^{-1/2})
 \qquad(\varepsilon\downarrow0).
\end{equation}
The upper construction retains independent bounded references and an
energy-sensitive entropy bound. The lower bound applies to every
signing law and already follows from one-row matrices.
\end{theorem}
The upper construction has height reserve of order $\varepsilon$, and the shear consumes height of order $R^{-2}$. These quantities give the exponent and direct the scalar optimization in Section~\ref{R11-sec:joint}.
\Needspace{15\baselineskip}
\begin{theorem}[Sharp entropy deficit at the Gaussian endpoint]
\label{VB-thm:entropy-endpoint}
Let
\[
 \Delta_G(\varepsilon)=\sup_{n\ge1}
 \sup_{\substack{\mu\ \text{a law on }\{-1,1\}^n\\
                 \mu\cx(1+\varepsilon)\kappa G_n}}
       \frac{D(\mu\|U_n)}n.
\]
Then $\Delta_G(0)=0$, and
\begin{equation}\label{VB-eq:endpoint-entropy}
 \Delta_G(\varepsilon)\le
 \sqrt{\frac\pi3}\sqrt{\frac{\varepsilon}{1+\varepsilon}},
 \qquad
 \Delta_G(\varepsilon)\sim\sqrt{\frac\pi3}\sqrt\varepsilon
 \quad(\varepsilon\downarrow0).
\end{equation}
For each positive $\varepsilon$, uniform binary code laws attain the
supremum asymptotically as their lengths grow.
\end{theorem}
The code laws attain the Gaussian-dominated optimum; hard-balanced laws inherit its entropy bound with their additional discrepancy constraint. The source-specific theorem gives the full frontier.

\paragraph{\textbf{Gaussian hitting and convex-order laws.}}
A sharp transfer connects the construction to classical vector balancing. A compact set meeting every bounded open convex set of Gaussian measure at least $1/2$ supports a law dominated by the standard Gaussian in convex order, at optimal scale one (Theorem~\ref{GC-thm:hitting}, Proposition~\ref{R10-prop:Gaussian-scale}). For symmetric supports this strengthens the subgaussian connection of~\cite[Theorem~3.4]{DGLN}. Applied to Banaszczyk's theorem~\cite{Banaszczyk}, it gives $A\sigma\cx N(0,25I_m)$ for unit columns. Theorem~\ref{VB-thm:main} reduces the comparison variance to $6.5$ and adds the hard bound, independent coefficient reference and entropy guarantee. The complementary extraction principle turns positive hitting of a bounded open set by every bounded dominated law into a convex subset of Gaussian measure at least $1/2$ (Proposition~\ref{GC-prop:extraction}), leading to the Talagrand results.

\subsection{Sharp Talagrand convexity, with no dilation}
Talagrand's question asks whether finitely many sums of a large Gaussian set must contain a large convex set. He describes the proposed conclusion as ``such an extraordinary fact if true''~\cite[Section~2]{Tal26}. Hua--Song--Tudose establish the Gaussian representation and the earlier existence result~\cite{HST}. Exact convex extraction gives the following sharp unweighted form.
\begin{theorem}[Talagrand convexity with the optimal three summands]
\label{R19-thm:intro-Talagrand}
If Borel sets $A_1,A_2,A_3\subset\R^d$ satisfy
$\gamma_d(A_1)+\gamma_d(A_2)+\gamma_d(A_3)>2$, then
\[
 K\subset A_1+A_2+A_3,\qquad \gamma_d(K)>\tfrac12
\]
for a convex body $K$. Symmetric inputs admit symmetric $K$. In particular, $\gamma_d(A)>2/3$ suffices for $K\subset A+A+A$. Three is the least universal number of unweighted summands for this conclusion.
\end{theorem}
Mazhar proves the equal-input, undilated half-mass conclusion in~\cite[Theorem~1.2]{R20-Mazhar}. Theorem~\ref{GC-thm:three} gives the mixed-input statement above, its strict output inequality and the compact-input endpoint. Talagrand's two-summand obstruction proves optimality of the summand number. The difficult step is to preserve membership in the original sum while extracting convexity. Compact inner approximation retains strict mass slack, and a convex test with median at most its mean turns positive hitting by every dominated compact law into the required convex subset.

The higher-mass target has its own explicit answer. If $A=-A$ and $\gamma_d(A)>2/3$, then
\begin{equation}\label{R19-eq:intro-fourfold}
 K=-K\subset A+A+A+A,\qquad
 \gamma_d(K)>\beta_4:=2\Phi\!\left(2\Phi^{-1}(3/4)\right)-1
 >0.822656449.
\end{equation}
Here $\Phi$ is the standard normal distribution function. This proves Talagrand's Problem~2.1 with four summands: balanced inputs of mass at least $3/4$ produce convex output of mass greater than $3/4$ (Corollary~\ref{R15-cor:fourfold}). Mazhar's constructive theorem gives six summands for the balanced higher-mass problem~\cite{R20-Mazhar}. Here four give the stated existence conclusion at input threshold $2/3$. The optimal three-summand statement concerns the half-mass target.

The extraction criterion is itself exact for a general fixed reference; its dimension-uniform log-concave threshold is $1/e$ (Theorem~\ref{EXT-thm:fixed}). The weighted representation gives a convex body of mass $>1/2$ inside $2\operatorname{conv}_3(A)$, where $\operatorname{conv}_3(A)$ is the union of convex hulls of at most three points of $A$ (Corollary~\ref{R19-cor:conv-three}). The unweighted theorem also gives the $k=3$, $\varepsilon=1/4$ instance of Talagrand's Conjecture~1.1 as stated by Johnston, with the stronger undilated sum and output mass $>1/2$~\cite{R19-Johnston}. Johnston's separate dilation-free Question~1.2 retains its own scope. For every fixed input mass $\alpha>1/2$, Corollary~\ref{R19-cor:explicit-counts} gives explicit dimension-independent summand counts. The compact-input endpoint in Theorem~\ref{GC-cor:compact-endpoint} treats mass $2/3$ separately.

On a separable Gaussian Banach space, the output is a compact convex set of mass at least $1/2$. In particular, any Borel family of paths of probability $>2/3$ has a threefold pathwise sum containing such a family (Theorem~\ref{GC-thm:banach}, Corollary~\ref{GC-cor:wiener}). Finite-rank Gaussian approximation keeps each selected sum inside a common compact set, so exact membership survives the passage to path space. A weighted Hilbert-space realization also answers Green's product-space Problem~54 with three summands and output mass at least one half (Corollary~\ref{R20-cor:green}).

\subsection{Five summands below half Gaussian mass}
Balancedness permits a smaller input threshold. A set $A$ is balanced when $tA\subset A$ for every real $|t|\le1$; the radial condition matters as well as symmetry. Write $A^{+k}$ for the actual $k$-fold Minkowski sum, without any exterior dilation.
\begin{theorem}[Five summands below half Gaussian mass]
\label{R22-thm:intro-fivefold}
In every separable Banach space with a centered Radon Gaussian measure, a balanced Borel set $A$ of mass at least $5/12$ admits a compact symmetric convex set $C$ such that
\begin{equation}\label{R22-eq:intro-fivefold}
 C\subset A+A+A+A+A,\qquad \gamma(C)\ge\Phi(1/16)>0.524917.
\end{equation}
For compact balanced input of mass at least $1/2$, five summands give output mass at least $3/4$. These are subsets of the actual Minkowski sums, including pathwise sums for Wiener measure.
\end{theorem}
The decisive step is a sharp Gaussian chord bound: compact star-shaped input of mass $\alpha$ satisfies
\[
 \rho_\alpha B_H\subset A-A,\qquad
 \rho_\alpha=2\Phi^{-1}((1+\alpha)/2).
\]
Here $B_H$ is the closed Cameron--Martin unit ball, equal to $B_2^d$ for standard Gaussian measure on $\R^d$. For balanced $A$, this ball lies in $A+A$. The weighted theorem needs only radius
$r_\alpha=[\Phi^{-1}(1-\alpha/2)-\Phi^{-1}(\alpha)]_+$ to raise the other two input masses sufficiently. After taking a compact limit, it places a half-mass convex core $K$ inside $A+\tfrac12A+\tfrac12A+r_\alpha B_H$. The unused chord radius $d_\alpha=\rho_\alpha-r_\alpha$ gives the stronger containment
\[
 K+d_\alpha B_H\subset A^{+5},\qquad
 \gamma(K+d_\alpha B_H)\ge\Phi(d_\alpha).
\]
The Cameron--Martin ball is compact in the original Banach space. At $\alpha=5/12$, the strict inequality $d_\alpha>1/16$ survives compact approximation of Borel inputs. At compact input mass $1/2$, $d_\alpha=\Phi^{-1}(3/4)$ gives the three-quarter output (Theorem~\ref{R23-thm:five-enlargement} and Corollary~\ref{R23-cor:five-mass}). Thus the same construction retains an entire Gaussian enlargement, which determines the output-mass bound.

At compact balanced input mass $1/2$, the chord argument and Mazhar's threefold theorem already imply a fivefold half-mass conclusion. The weighted calculation gives both the lower input threshold in \eqref{R22-eq:intro-fivefold} and the fivefold three-quarter-mass endpoint. Six copies of compact balanced half-mass input give the larger output mass $\beta_4>0.822656449$, and $6j$ copies give $2\Phi(2j\Phi^{-1}(3/4))-1$ (Corollary~\ref{R22-cor:sixfold}). Theorem~\ref{R19-thm:intro-Talagrand} retains its distinct arbitrary-input and optimal-three-summand scope.

\subsection{The exact first-order certificate barrier}
The smaller radius depends on retaining the full overlap profile. For a probability density $f$ supported in $[-R,R]^m$, extend $f$ by zero and let $|D_vf|$ be the total variation of its distributional derivative in direction $v$. Among all such densities of bounded variation, including dependent ones, define
\[
 \tau_m(R)=\inf_f\sup_{\|v\|_2=1}|D_vf|.
\]
Theorem~\ref{M-first-order} proves the exact limit
\begin{equation}\label{R19-eq:intro-firstorder}
 \lim_{m\to\infty}\tau_m(R)=\frac{\sqrt{2\pi}}R.
\end{equation}
For the uniform auxiliary on $(-3,3)$, the sufficient first-order overlap certificate consequently has sharp dimension-uniform radius $3\sqrt{2\pi}$. Changing the initial density alone cannot lower that certificate's universal radius. The exact overlap retains how translation losses accumulate over a finite interval; replacing those losses by their derivative bound discards the information used to reach \eqref{R19-eq:intro-radius}.

Optimizing the full diffuse overlap is a separate extremal problem. Theorem~\ref{R7-thm:method-constant} permits dependent references chosen before the direction. Both its diffuse threshold and the common-auxiliary threshold lie between $6.8383231851$ and $6.8383231852$. The distinct optimization classes and their common numerical enclosure are specified in that theorem.

\subsection{What the retained law makes possible}\label{R25-sec:intro-applications}
The joint precision bound already gives a simultaneous norm consequence. For square matrices with $\ell_p$ column norms at most one, $p\ge2$, one signing satisfies the Reis--Rothvoss conjectured dimension order for every $q\ge p$ at once. For rectangular matrices the same conclusion holds at $p=2$ (Corollary~\ref{R19-cor:RR}). The precision formula selects one outcome with small quadratic energy and interpolates with its hard $\ell_\infty$ bound. The dimension orders in these named cases are already consequences of constant-discrepancy results and interpolation; the conclusion here gives the exact spectral energy and one signing for every target norm. Section~\ref{R19-sec:Lebesgue} proves both statements.

A second consequence starts with independent inputs $X_i$ and a function whose change in coordinate $i$ is at most $c_i$. Its Doob increments lie in predictable centered intervals. Completing their extremal two-point laws by Gaussian increments gives the sharp full-convex-order comparison
\begin{equation}\label{R19-eq:intro-bounded-differences}
 f(X)-\E f(X)\cx N\!\left(0,\frac\pi8\sum_i c_i^2\right).
\end{equation}
A single fair bit forces $\pi/8$; Theorem~\ref{SP-thm:mcdiarmid} covers every finite convex test. This same scalar coefficient governs unconstrained integer rounding. For adapted $0\le X_i\le c$ with $m_i=\E[X_i\mid\mathcal F_{i-1}]$ and deterministic predictable mean budget $\sum_i m_i\le L$, the companion sharp upper law is
\[
 \sum_i(X_i-m_i)\cx c\{\operatorname{Pois}(L/c)-L/c\}
\]
(Theorem~\ref{SP-thm:poisson}). Completion couples each reference to the whole original history. The Gaussian theorem uses predictable interval widths, and the Poisson theorem uses predictable means; neither replaces these budgets by quadratic variation.

The passage from signing to sampling requires exact marginals. A dyadic pairing of local choices preserves each prescribed local distribution. The same joint reference then controls every later integrand, including its residual after projection onto a feature space.
\Needspace{14\baselineskip}
\begin{theorem}[One sampling law for all integrands]\label{VB-thm:sampling}
Let $P$ be a probability law on a standard Borel space, and let $\Phi$ map measurably into a real separable Hilbert space with essential diameter at most $d>0$. For every $N$, there is an exchangeable law with $Z_i\sim P$ exactly such that, writing $\widehat P_Ng=N^{-1}\sum_i g(Z_i)$,
\begin{equation}\label{VB-eq:sampling}
 \Var(\widehat P_Ng)\le\inf_v\left\{
 \frac{2.038d^2}{N^2}\|v\|^2+
 \frac{5.569}{N}\Var_P(g-\langle v,\Phi\rangle)\right\}
 \quad(g\in L^2(P)).
\end{equation}
For finite-dimensional $\Phi$ the same law has $\|\widehat P_N\Phi-P\Phi\|_\infty<6.84d/N$ almost surely. For a canonical reproducing-kernel feature map it has $\sup_{\|f\|\le1}\E|\widehat P_Nf-Pf|^2\le2.038d^2/N^2$.
\end{theorem}
The same law retains Gaussian-process, mixed exponential and information bounds. Polynomial-exact local formulas give one equal-weight Sobolev rule across all lower smoothness, with exact uniform marginals, hard features and optimal scalar rates. Fixed-integrand and expected squared worst-case errors are distinguished in Section~\ref{R11-sec:sampling}. History-preserving composition also covers centered L\'evy references with only a first moment (Section~\ref{R11-sec:composition}).

Simultaneous exchange families retain hard features under strongly base-orderable or totally unimodular constraints. Unconstrained integer rounding has the optimal Gaussian covariance factor $\pi/8$ and an exact expected polynomial-bit rational sampler. For arbitrary integral base polyhedra, nearest-cell geometry gives Gaussian rounding in the minimal face and an expected polynomial-bit value-oracle sampler with positive variance slack. Table~\ref{R9-tab:runtime} separates these polynomial-bit algorithms from limiting existence and specifies their input models.

The matrix construction also permits an arbitrary starting vector $y\in[-1,1]^m$. If symmetric $d\times d$ matrices satisfy $\|\sum_i A_i^2\|\le s^2$, Theorem~\ref{R19-thm:matrix-partial} gives one law with $\E X=y$, more than $m-d$ coordinates in $\{-1,1\}$ at every outcome, hard matrix error $O(s)$ and a Gaussian comparison for $X-y$. The hard conclusion of Sachdeva--Thudi--Zhao's Conjecture~21 also follows from Akbas--Sra's full-coloring theorem by hereditary dyadic rounding. The result here preserves the prescribed mean and a Gaussian law comparison throughout that rounding, using the same analytic input. Hereditary signing controls the dyadic increments; their squared Gaussian scales sum to $4/3$.

\noindent\textit{Computational comparison.}
For unit Euclidean columns, deterministic bounds are $8272$ in Guo--Fang--Lu~\cite{GFLalg}, below $135$ in Akbas--Sra's arithmetic model~\cite[Theorem B.1]{R19-AkbasSraExposition}, and below $99$ with polynomial bit complexity for rational inputs in Li~\cite{R19-Li}. Akbas--Sra give rational-arithmetic and faster square-root/linear-algebra operation bounds. Table~\ref{R9-tab:runtime} specifies the polynomial-bit models of the integer and base-polyhedron samplers described above.

For spectral approximation, hard accuracy and unbiasedness hold under one sparse law. If $A_i\succeq0$ and $\sum_iA_i=I_d$, every weighted sum lies between $(1-\delta)I_d$ and $(1+\delta)I_d$, and each original weight has mean one. Each sample uses $O(d/(\varepsilon\delta^2))$ nonzero weights, with Poisson information cost $O(d\sqrt\varepsilon)$ independent of thinning depth (Theorem~\ref{FC-thm:poisson-sparse}). The continuous version retains the exact expected measure.

\subsection{The two Schatten conjectures and trace-class volume ratios}
A matrix law can control an entire scale of norms, so the exponent may be chosen after the law. The next two results resolve Conjectures~2 and~3 of Reis--Rothvoss~\cite{R19-RR}. Here $\|B\|_{S_p}$ is the $\ell_p$ norm of the singular values of $B$, and $S_\infty$ is the operator norm.
\begin{theorem}[One prescribed-mean law for every Schatten norm]\label{R27-thm:intro-Schatten}
For symmetric $n\times n$ matrices $A_1,\ldots,A_n$ and every $y\in[-1,1]^n$, there is a law on signs with $\E\sigma=y$ such that every outcome satisfies
\begin{equation}\label{R27-eq:intro-Schatten}
 \left\|\sum_i(\sigma_i-y_i)A_i\right\|_{S_p}
 \le C n^{\max\{1/2,1/p\}}\max_i\|A_i\|_{S_p}
 \qquad(1\le p\le\infty).
\end{equation}
The constant is universal, independent of $p$, and the powers of $n$ are optimal. The law may have at most $n+1$ support points.
\end{theorem}
For $p\ge2$ this gives the conjectured $C\sqrt n$ bound. Theorem~\ref{RRS-thm:general-law} gives the dimension-dependent rectangular and complex extensions. These are the simultaneous norm guarantees of this law; the joint Gaussian and entropy guarantees in Theorem~\ref{R3-thm:matrix} belong to the earlier reference-law construction.

\begin{theorem}[One Gaussian event for every unitarily invariant norm]\label{R27-thm:intro-square-function}
For the same matrices, put $V=(\sum_iA_i^2)^{1/2}$. There is a single measurable $K\subset\R^n$ with $\gamma_n(K)\ge e^{-Cn}$ such that
\begin{equation}\label{R27-eq:intro-square-function}
 \left\|\sum_i x_iA_i\right\|_{\mathcal N}\le\|V\|_{\mathcal N}
 \quad(x\in K)
 \quad\text{for every unitarily invariant norm }\mathcal N.
\end{equation}
\end{theorem}
The Gaussian event is proved first. A positive multiplier gives $A_i=HU_i+U_iH$ and $\sum_iU_i^2\preceq I$. Complementarity makes $T=H+\sum_iU_iHU_i\succeq H$, while row contraction gives $T^2\preceq V^2$. Hence $H\preceq V$: one operator-norm event controls every partial singular-value sum. Akbas--Sra's small-ball estimate~\cite{AkbasSraMatrix} then gives the stated Gaussian measure.

To obtain full signs, the event must remain large enough as the active set shrinks. The variance-sensitive form of Akbas--Sra's interpolation estimate gives this hereditary bound. Partial coloring preserves an arbitrary separating functional, so $y$ lies in the convex hull of the simultaneously good signings. Section~\ref{R27-sec:Schatten} proves this passage to the exact-mean law and gives the common-polar comparison with Cadilhac--Ricard~\cite{R27-CadilhacRicard}.

\paragraph{\textbf{Reis's trace-class volume-ratio conjecture.}}
The matrix density of Section~\ref{R11-sec:sparse} also controls volume. For a full-dimensional symmetric body $L\subset\R^r$ with largest inscribed ellipsoid $E_L$, set $\vr(L)=(|L|/|E_L|)^{1/r}$. Theorem~\ref{SVR-thm:main} proves, for every surjective $M:\R^{a\times b}\to\R^r$,
\begin{equation}\label{R27-eq:intro-volume}
 \vr(MB_{S_1}^{a\times b})\le C\sqrt{1+\min\{a,b\}/r}.
\end{equation}
Thus every quotient with dimension at least the smaller matrix dimension has bounded volume ratio; $a=b=r=d$ proves Reis's Conjecture~1~\cite{R27-ReisZonotopes}. John contacts bound both operator variances by $r$ and the undilated Gram trace by $r\min(a,b)$. The Fisher determinant gives a lower bound for the polar's volume through entropy, and Blaschke--Santal\'o turns it into \eqref{R27-eq:intro-volume}. Section~\ref{R27-sec:volume} gives the determinant refinement and its spectrahedral consequence using the Akbas--Sra density.

\subsection{From exact allocation to deterministic dynamics}
The deterministic realization problem starts after a law has been selected. To preserve an entire continuous output law and finitely many conditional moments, a simplex of source means divides the target measure into barycentric shares. Mean-independent factors realize these shares exactly. This allocation proves the continuous realization theorem. The related
finite-moment problem has a unique cost-minimizing polynomial partition
under positive feasibility; Section~\ref{R11-sec:polynomial} proves
its optimal selecting degree and exact-moment regularization.

\paragraph{\textbf{Exact Gaussian dynamics on prescribed observables.}}
For every Gaussian contraction and every finite degree $k$, one Gaussian-preserving deterministic map reproduces its polynomial action at every iterate. All multitime polynomial expectations of total degree at most $k$ agree with the Gaussian autoregression, while the last observation can be an arbitrary bounded test. The same map has relative-entropy contraction coefficient $1$ and infinite $L^p\to L^q$ norm for $p<q$ (Theorems~\ref{WC-thm:gaussian}--\ref{WC-thm:multitime}, Corollary~\ref{WC-cor:no-hyper}). Exact conditional moments make the construction possible. Keeping only four orders and letting the time step vanish gives stationary diffusion limits, including Brownian motion on every closed connected smooth manifold (Theorem~\ref{WC-thm:compact-diffusion}).

\subsection{Nutz--Wang--Zhang transport: approximation and the component label}
Directional allocation preserves order and both marginals, giving $W_\infty$ density and resolving Conjecture~5.4. Supermartingale allocation preserves the inequality and its equality region, giving $W_1$ density and $W_p$ density under finite $p$-moments, resolving Conjecture~5.3. For a countable initial law in arbitrary finite dimension, the martingale approximation is $W_\infty$-close and retains any prescribed finite list of bounded costs, resolving that clause of Conjecture~5.1. These are separate theorems in Section~\ref{R11-sec:monge}.

For continuous initial laws, full relative span or a terminal-recoverable common affine label permits the conditional realizations to be joined. Corollary~\ref{R19-cor:separated-labels} proves a geometric sufficient condition: terminal mass avoids all overlaps of the relevant component closures. A measurable singleton-section argument recovers the label from the terminal point; conditional atomlessness and the full-relative-span seed then give backward-Monge density.

The local realizations have a sharp global limitation.
\begin{theorem}[Atomless components with a forced backward bit]\label{R15-thm:intro-boundary}
There are compactly supported atomless marginals in $\R^4$, with full affine
span and a unique martingale coupling, such that every irreducible component
has an atomless conditional terminal law and permits a constant backward
map. Globally, however,
\[
 \law(X\mid Y)=\tfrac12\delta_{h(Y)}+\tfrac12\delta_{k(Y)},
 \qquad h(Y)\ne k(Y)\quad\text{almost surely}.
\]
Thus no backward-Monge coupling exists: one fair bit is necessary and sufficient.
\end{theorem}
Theorem~\ref{BM-thm:main} gives every Bernoulli parameter. A shared terminal boundary leaves two possible component labels at the same observed point, forcing the split. The NWZ consequence uses $\nu_I=\law(Y\mid I(X)=I)$; terminal recoverability $I(X)=J(Y)$ is stronger (Corollary~\ref{BM-cor:NWZ}). In this same example, every strict contraction of the initial law or strict inflation of the terminal law restores backward determinism (Corollary~\ref{BM-cor:quantitative}).

\subsection{Stability when the marginals and reference vary}
In dimensions at least two, small perturbations of both marginals can destroy approximation of a prescribed martingale coupling and continuity of the optimal value. Br\"uckerhoff--Juillet prove this instability even for perturbations with full support~\cite{R20-BJ}. Real-line stability is established by Backhoff-Veraguas--Pammer and Wiesel~\cite{R20-BP,R20-Wiesel}. In higher dimensions the condition that restores approximation here is an exact positive product component. It leaves enough mass to repair the conditional means while retaining the new terminal marginal.

Write $\mathcal M(\mu,\nu)$ for the martingale couplings, and $\mathcal P_p(\R^d)$ for probability laws with finite $p$th moment. The following theorem permits the entire initial coupling of the old and new problems to be prescribed.
\begin{theorem}[All-marginal adapted stability]\label{R20-thm:intro-stability}
Let $p\ge1$, let $\mu$ be compactly supported, and let $\nu\in\mathcal P_p(\R^d)$ have full affine span. Suppose
\[
 s\in\mathcal M(\mu,\nu),\qquad s\ge\eta\,\mu\otimes\nu
 \quad\text{for some }\eta>0.
\]
If $\mu_n\cx\nu_n$ and $(\mu_n,\nu_n)\to(\mu,\nu)$ in $W_p$, then for every $\pi\in\mathcal M(\mu,\nu)$ and every coupling $\xi_n$ of $\mu,\mu_n$ with $\int|x-x'|^p\,d\xi_n\to0$, there are $\pi_n\in\mathcal M(\mu_n,\nu_n)$ such that
\[
 \int W_p^p(\pi_x,(\pi_n)_{x'})\,\xi_n(dx,dx')\longrightarrow0.
\]
Here $\pi_x$ and $(\pi_n)_{x'}$ are the conditional terminal laws. Both marginals and every conditional mean are exact at each $n$.
\end{theorem}
This is Theorem~\ref{AS-thm:main}. The product condition is equivalent to convex order after some strict outward dilation of $\mu$ about its mean. It holds for irreducible pairs whose initial support is compactly inside the terminal convex hull (Lemma~\ref{AS-lem:product}). The approximating pairs may be unbounded or reducible. The proof first groups the initial states into finitely many cells. A bounded family of centered terminal functions repairs every row mean, and a final correction repairs the terminal marginal without disturbing those means. The product component keeps these signed corrections nonnegative.

\Needspace{6\baselineskip}
Optimizers can therefore be stable even when every input to their variational problem changes. For a reference law $q$, write
\[
 \mathcal P_q(\mu,\nu)
 =\sup_{\pi\in\mathcal M(\mu,\nu)}
       \int\MCov(\pi_x,q)\,\mu(dx),
\]
where $\MCov(r,q)$ is the largest $\E\langle Y,W\rangle$ over couplings with laws $r,q$.
\begin{corollary}[Both marginals and the driving law]\label{R20-cor:intro-three-laws}
Under the preceding hypotheses with $p=2$, suppose $q_n\to q$ in $W_2$ and $q$ is absolutely continuous. Then the values $\mathcal P_{q_n}(\mu_n,\nu_n)$ converge to $\mathcal P_q(\mu,\nu)$. Every sequence of asymptotically optimal couplings converges to the unique limiting optimizer in adapted $W_2$, and satisfies the displayed conditional-kernel convergence along every prescribed $\xi_n$. The approximating references $q_n$ may be atomic.
\end{corollary}
Exact feasibility gives the lower bound for the optimal value. The reverse passage must also recover the conditional laws: initial points merging in the limit can conceal different terminal kernels. Averaging those kernels preserves every martingale constraint. Maximal covariance is strictly concave in the terminal law when the limiting driver is absolutely continuous, so an unresolved mixture would improve the objective on averaging. Approximate optimality therefore forces one conditional kernel at each limiting initial point. Lemma~\ref{AS-lem:lift} turns this collapse into the adapted convergence in Corollary~\ref{AS-cor:three-laws}.

This conclusion requires no terminal log-concavity. Any compact-source convex-order pair with a full-dimensional second-moment terminal law enters the product-component regime after an arbitrarily small strict terminal dilation. The original signing law is therefore retained while both its terminal reference and its driving law are approximated. Example~\ref{AS-ex:upper} exhibits the hidden-kernel phenomenon and its strictly positive Gaussian objective loss.

The same correction preserves an entire prescribed history and gives finite-horizon, non-Markov adapted approximation. It also gives ordinary Wasserstein Hausdorff convergence of feasible sets and continuity for continuous costs of the stated moment growth. For weak costs, the corresponding convexity or concavity in the conditional law is retained (Corollaries~\ref{AS-cor:history} and~\ref{AS-cor:optimization}). Example~\ref{AS-ex:upper} shows that adapted upper Hausdorff convergence can still fail, even in this interior regime.

\subsection{Sharp regularity of stretched Brownian motion}
With the endpoints fixed, terminal curvature determines the regularity of the variationally selected martingale. The objective chooses one coupling; the contraction theorem then controls its conditional maps and its entire future path. This conclusion concerns the optimizer itself, including reducible endpoint pairs.

Let $\nu(dy)=Z^{-1}e^{-V(y)}dy$, with $V-|\cdot|^2/(2L^2)$ convex, and let $\mu\cx\nu$. Write $\mathcal P(\mu,\nu)$ for the maximum of $\E\langle Y,G\rangle$ over martingale pairs $(X,Y)$ with laws $(\mu,\nu)$ and $G\sim\gamma_d$ independent of $X$. The conditional Brenier maps $T_X$ of its unique optimizer are $L$-Lipschitz without irreducibility, and satisfy the sharp deficit bound
\begin{equation}\label{R11-eq:intro-deficit}
 \mathcal P(\mu,\nu)-\E\langle Y,G\rangle
 \ge\frac1{2L}\E|Y-T_X(G)|^2,
 \qquad G\sim\gamma_d,\ G\perp X.
\end{equation}

For a competing martingale $M_t=X+\int_0^t\sigma_s\,dB_s$ with $\law(M_1)=\nu$ and the canonical martingale $M^*$ on the same $(X,B)$, put
$\Delta=\mathcal P(\mu,\nu)-\E\int_0^1\tr\sigma_t\,dt$.
Then
\begin{equation}\label{R21-eq:intro-same-driver}
 \E\int_0^1\|\sigma_t-\sigma_t^*\|_{\mathrm{HS}}^2\,dt
 =\E|M_1-M_1^*|^2\le 2L\Delta.
\end{equation}
The factor $2L$ is sharp for this fixed-driver inequality.
Conditional laws at every stopping time satisfy logarithmic Sobolev
and Poincar\'e inequalities; the complete future-path law satisfies
Gaussian isoperimetry and transportation with sharp coefficient
$2L^2(1-\tau)$ in the uniform metric. Section~\ref{R11-sec:canonical}
proves these statements for the optimizer itself.

The sharp canonical deficit also quantifies a change of driver.
For a uniformly log-concave terminal law as above, every optimizer
$\pi^q$ for a centered reference $q$ satisfies
\begin{equation}\label{R13-eq:intro-driver-rate}
 \left(\int W_2^2(\pi_x^q,\pi_x^\gamma)\,\mu(dx)\right)^{1/2}
 \le 2L W_2(q,\gamma_d).
\end{equation}
The comparison uses the same initial state, includes atomic approximations and keeps both endpoint laws exact; optimization error is also allowed. Thus a finite driving reference approximates the canonical coupling of the original signing law without changing any signing probability (Corollary~\ref{R13-cor:signing-reference}). An upper log-density curvature bound gives non-Gaussian conditional contraction and sharp deficit for arbitrary initial laws (Theorem~\ref{R14-thm:reference-curvature}); a common global potential requires irreducibility.

\subsection{Gaussian convex order and bounded-time embedding}
The canonical contraction has a direct probabilistic consequence. Whitening the terminal covariance and taking the heat extension of each conditional map produces a martingale with a deterministic instantaneous covariance cap. Completing its Brownian driver gives, in every dimension,
\[
 \mu\cx N(0,Q)
 \quad\Longleftrightarrow\quad
 \begin{gathered}
 M_0\sim\mu,\quad M_1\sim N(0,Q),\quad M\text{ continuous martingale},\\
 0\preceq d\langle M\rangle_t/dt\preceq Q.
 \end{gathered}
\]
The Gaussian has covariance $Q\succ0$. Every stopping-time residual retains the corresponding conditional Gaussian comparison. In dimension one, time change gives the bounded-time characterization
\[
 \mu\cx N(0,1)
 \quad\Longleftrightarrow\quad
 X+B_\tau\sim N(0,1),\quad X\sim\mu,\quad\tau\le1,
\]
where $B$ is Brownian in a filtration containing $X$ at time zero and permitting auxiliary randomness; hence $B$ is independent of $X$. Corollaries~\ref{GC-cor:gauss-dynamic} and~\ref{GC-cor:skorokhod} give the constructions. The conditional-mean identity $\E[X+B_\tau\mid X]=X$ gives necessity.

The starting law is $\mu$ and the terminal law is Gaussian. Classical bounded-time Skorokhod embedding starts at a point and prescribes the target; Lipschitz Gaussian quantiles give sufficient bounded-time conditions in that direction~\cite{R20-AS,R20-AHS}. For example, a symmetric two-point law $\{-a,a\}$ with $0<a\le\sqrt{2/\pi}$ lies below the Gaussian in convex order, so the displayed characterization evolves it to a Gaussian by time one. Brownian motion started at zero cannot reach that two-point terminal law by a deterministic time bound, since its first exit from $(-a,a)$ has unbounded support.

Existence of a suitable weighted Gaussian representation also follows by combining Hua--Song--Tudose's finite strict-slack conditional log-concavity with Caffarelli contraction, their half-sum representation and finite approximation~\cite{HST}. The regularity theorem above identifies the canonical optimal realization and proves its sharp deficit and path inequalities. Conditional log-concavity is used in its finite strict-slack role; the bounded-covariance proof follows canonical contraction, heat interpolation and Brownian completion.

\subsection{Arbitrary references, all optimizers and exact \texorpdfstring{$q$}{q}-Bass existence}\label{R25-sec:intro-general-reference}
A prescribed driver raises a compatibility question: the optimal conditional allocations must fit together while preserving the whole terminal law. For an irreducible convex-order pair with finite second moments and full terminal affine span, every second-moment reference $q$, including an atomic one, has one convex potential and one measurable shift representing \emph{every} optimizer:
\[
 Y\in\partial v(a(X)+W),\qquad W\sim q,\quad W\perp X,
 \qquad \E[Y\mid X]=X.
\]
Theorem~\ref{QB-thm:arbitrary} proves conditional dual attainment in this same generality. Conditional Jensen gaps give a common normalized limit even when its separate global integrals diverge; recovering competitors with their means exact completes the limiting variational problem.

If $q$ charges no set of Hausdorff dimension at most $d-1$, the subgradient relation becomes
\begin{equation}\label{R13-eq:intro-qbass}
 Y=\nabla v(a(X)+W),\qquad W\sim q,\quad W\perp X,
 \qquad \E[Y\mid X]=X.
\end{equation}
Theorem~\ref{QB-thm:main} therefore gives a common extended $q$-Bass potential. The marginals may be singular or unbounded, and the potential is proper lower-semicontinuous convex. Theorem~\ref{R24-thm:q-Bass-universal} gives the exact universal finite-potential criterion, already detected by one Gaussian terminal test:
\[
 \begin{gathered}
 \overline{\conv}\supp q=\R^d\quad\Longleftrightarrow\quad
 (\delta_0,\gamma_d)\text{ admits a globally finite subgradient representation}\\
 \Longleftrightarrow\quad\text{one globally finite potential represents all optimizers}\\
 \text{for every irreducible second-moment pair with full terminal affine span}.
 \end{gathered}
\]
This holds for every $q\in\mathcal P_2(\R^d)$, allowing atoms, arbitrary latent laws and randomized subgradient allocations. Under the preceding nondegeneracy hypothesis it settles the universal finite-potential form of Tschiderer's $q$-Bass problem. Compact terminal support gives finiteness for an individual pair, while extended potentials retain existence for arbitrary reference support.

For the Gaussian test, couple two latent allocations over the same driver. Their supporting-plane gaps are nonnegative and have expected sum zero because both terminal means are zero. The gaps therefore vanish, so the entire latent mixture uses one translated subgradient coupling. Gaussian terminal support makes its conjugate finite everywhere; a half-space containing the reference support would then force the original potential to be infinite outside that half-space. This explains both the full-support criterion and the uniform-reference obstruction. For $q$ equivalent to Lebesgue measure, existence of a globally finite common-gradient representation is equivalent to irreducibility, and the latent translation is unique up to a common shift.

\subsection{Positive noise, atomic targets and Bass calibration}
Calibration requires a different regularity property: curvature of the optimized cost as the latent variable moves. The target allocation itself changes with that variable, so its mass constraints must be retained in the second variation. Let $q,\nu$ have finite second moments, let $q$ have a continuous everywhere-positive density, and define $F_{\nu,q}(Z)=\MCov(\law(Z+\xi),\nu)$ for independent $\xi\sim q$. Theorem~\ref{BF-thm:maincurvature} gives the exact criterion
\[
 \nu\text{ has full affine span}
 \quad\Longleftrightarrow\quad
 \begin{gathered}
 F_{\nu,q}\text{ is strongly convex modulo translations}\\
 \text{on every bounded $L^\infty$ latent set}.
 \end{gathered}
\]
The modulus is independent of the latent probability space. Atomic targets are included: their distributional transport Hessians retain curvature concentrated on cell boundaries, which additive noise preserves in the optimized second variation.

For Gaussian noise, confinement gives a bounded latent trajectory, so this local curvature yields exponential Bass-flow convergence in every dimension. The pair is irreducible, the initial support is compactly inside the terminal convex hull and initialization is bounded; the target may be singular and unbounded (Theorem~\ref{BF-thm:flow}). Corollary~\ref{BF-cor:published} resolves Backhoff--Pammer--Schachermayer's conjecture~\cite{BF-BPS}. Under terminal curvature, endpoint couplings and conditional-mean residuals converge exponentially with both marginals exact at finite time. Every reference inflation $c>1$ makes the original signing law eligible for this calibration and for a common $q$-Bass potential with any admissible driver, preserving all signing probabilities and hard constraints (Corollary~\ref{QB-cor:signing}).

\subsection{Asymptotic spectral equality and reference geometry}
The ground-state score identity determines an entire asymptotic spectral range. Briani--Buttazzo--Prinari prove that the sharp dimension-dependent first-eigenvalue Cheeger ratio tends to its lower-bound constant as dimension grows~\cite{R20-BBP}. The following conclusion holds along every product sequence in a fixed family and for every subexponential eigenvalue index. For arbitrary products of $N$ factors drawn from a fixed finite family of weighted bounded connected Lipschitz domains, put
$E_N=\sum_i\lambda_{1,2}(\Omega_i,w_i)$. Then
\begin{equation}\label{R11-eq:intro-spectrum}
 h_D(D_N)\sim2\sqrt{E_N},\qquad
 \lambda_{k_N,p}(D_N)^{1/p}\sim\frac2p\sqrt{E_N}
 \quad\text{if }\log k_N=o(N),\quad1<p<\infty.
\end{equation}
Here the perimeter is the full Dirichlet perimeter and the eigenvalues form the genus sequence. Cheeger's lower inequality becomes asymptotically exact throughout every subexponential spectral range. For the cube $(-1,1)^N$, this gives $h_D\sim\pi\sqrt N$ and $\lambda_{k_N,p}^{1/p}\sim\pi\sqrt N/p$. The factor domains may be nonconvex and their proportions need not converge. The summed score field gives the lower bound. For the upper bound, $\lceil\log_2 k_N\rceil$ modified factors encode $k_N$ disjoint trials, at a cost negligible relative to the order-$N$ ground-state energy.

Optimizing a reference density leads to a spectral problem with an exact interpretation. The Smirnov--Vershynin Fisher-information capacity satisfies
\[
 I(K)^2=4\Lambda(K)
 =4\sup_{M\succ0,\,\tr M=1}\lambda_1(-\operatorname{div}(M\nabla);K)
\]
for every nonempty bounded open $K$, including nonconvex and disconnected domains. The Sobolev relaxation attains the optimum. Products take the maximum capacity. For convex domains, inverse capacity is concave under Minkowski interpolation and balls minimize capacity at fixed volume (Corollary~\ref{R20-cor:capacity}). This identification retains the operator norm in the existing capacity~\cite{R20-SV} and its exact support normalization.

The full score law governs further retained geometry: compact log-concave densities with equal Fisher matrices and all three-coordinate score marginals can have diagonal widths differing by order $\sqrt m$ (Theorem~\ref{S-thm:r7-compact-pair}). Even all polynomial score moments can agree while sharp support radii differ (Theorem~\ref{R9-thm:score-moment-indeterminacy}).

\subsection{Exact entropy data and reconstruction of experiments}
\paragraph{\textbf{Exact entropy data identify the normalized source.}}
The exact entropy optimum for a symmetric scalar source is the convex conjugate of $\E\log\cosh(tW)$, attained with the entire product source and every conditional mean preserved. Its exact values on positive scales accumulating at one positive point determine the normalized source, including atoms, under only a first moment (Theorem~\ref{SP-thm:entropy-inverse}). Analyticity, convex duality and Mellin--Fourier inversion recover the law even when moments leave retained geometry undetermined.

Fixed common and coordinate noises admit an exact all-dimensional feasibility and entropy classification. Gaussian noises have an explicit phase boundary and linear-arithmetic coupling. If the common noise may vary, cancellation gives every empirical-mean limit. Arbitrarily small smooth coordinate-noise changes can force consensus while each fixed dimension still permits independence.

For a retained statistical experiment, conditional quantile allocation determines every feasible independent-source extension. One Laplacian comparison controls its full nonlinear information potential at every amplitude. Complete positivity gives the exact finite latent-realization criterion. 

\paragraph{\textbf{One simulator for every bounded decision problem.}}
Refining the target creates new calibration tests without changing the observed information: every finite refinement level is strict, their completion is Blackwell order, and a quantitative comparison gives one simulator for all bounded decision problems. 

\paragraph{\textbf{A directed cost profile reconstructs the entire experiment.}}
Varying the source weights of one fixed directed cycle cost then recovers the full posterior law, including its boundary atoms, with a sharp inverse modulus on the closed simplex (Theorems~\ref{CYC-thm:boundary} and~\ref{CYC-thm:sharp}). For $q$ states, write $J_\nu(a)$ for the mean cost of transporting source weights $a$ to the posterior under that fixed cycle, and let $\Delta_q$ be the probability simplex. Then
\[
 d_{\rm cyc}(\nu,\nu')=\sup_{a\in\Delta_q^\circ}|J_\nu(a)-J_{\nu'}(a)|,
 \qquad W_1(\nu,\nu')\le C_q d_{\rm cyc}(\nu,\nu')^{1/q}.
\]
The exponent $1/q$ is optimal even for smooth interior densities with the same prior; the upper bound also covers every boundary atom. Moving mass from the last source state to the first recovers the joint posterior distribution function from a directional derivative of the cost. Recovering its measure requires $q$ derivatives, which explains the sharp inverse exponent.

For odd alphabets, symmetrizing the cost loses identifiability even for smooth interior examples.

\subsection{Exact privacy capacity and the Gaussian budget transition}
In the full-data model, an encoder observes the private variable $Z$ and useful variable $Y$. Conditional allocation preserves exact privacy while approaching full finite-alphabet utility. Conditional atomlessness of $Y$ given $Z$ is equivalent to the utility supremum $H(p)$ for every prescribed finite output law $p$, under exact independence of the release and $Z$ (Corollary~\ref{SP-cor:privacy-iff}). For correlated Gaussian scalars with $0<|\rho|<1$, the finite-alphabet value is $\log M$, unattained for each $M\ge2$, disproving Rassouli--G\"und\"uz Conjecture~1.

For every prescribed $b>0$, the construction also permits one deterministic uniform release $U=F(Z,Y)$ satisfying
\[
 U\perp Z,\qquad I(U;Y)=\infty,\qquad I(Z;U\mid Y)\le b,
 \qquad \Prb\{U=g(Y)\}\ge e^{-b}
\]
for a measurable predictor $g$ (Theorem~\ref{SP-thm:infinite}, Corollary~\ref{R15-cor:private-prediction}). The whole-output predictor $g(Y)$ may itself carry private information. For the correlated Gaussian pair above, define
\[
 \mathcal U(b)=\sup\{I(U;Y):U\perp Z,\ I(Z;U\mid Y)\le b\},
\]
where full-data mechanisms with standard Borel outputs are allowed. The exact budget profile is
\begin{equation}\label{R19-eq:intro-private-budget}
 \mathcal U(0)=0,\qquad \mathcal U(b)=\infty\quad(b>0)
\end{equation}
(Corollary~\ref{SP-cor:gaussian-budget}); each positive-budget value is realized by the single release above.

\subsection{Smooth data laws with discontinuous exact privacy}
Zero conditional-information budget makes the release a channel from $Y$ alone. Its utility can change from zero to infinity under smooth covariance-preserving changes of the data law. Write $\mathcal U_P$ for the preceding utility when $(Z,Y)$ has joint law $P$.

\begin{theorem}[Privacy instability at every correlated Gaussian pair]
\label{R22-thm:intro-privacy}
Fix $0<|\rho|<1$ and $P_0=N(0,\Sigma_\rho)$ with
$\Sigma_\rho=\left(\begin{smallmatrix}1&\rho\\\rho&1\end{smallmatrix}\right)$.
For every $0<a\le1/2$ there is a centered law $P_a$ with exactly this covariance and private marginal $N(0,1)$, having a positive real-analytic density $e^{-V_a}$ and
\begin{equation}\label{R22-eq:intro-privacy-curvature}
 (1-a^2)\Sigma_\rho^{-1}\preceq\nabla^2V_a
 \preceq(1+a^2)\Sigma_\rho^{-1}\quad\hbox{everywhere}.
\end{equation}
These laws converge to $P_0$ in $W_2$ and in both relative-entropy directions, while
\[
 \mathcal U_{P_a}(0)=\infty\quad(a>0),\qquad \mathcal U_{P_0}(0)=0.
\]
For each $a>0$, a deterministic uniform release from $Y$ alone attains infinite utility with exact privacy and whole-output prediction probability one.
\end{theorem}
The example is explicit. Take independent $Z,G\sim N(0,1)$ and $T_a$ uniform on $[-a,a]$, put $c=\sqrt{1-\rho^2}$, and set
\begin{equation}\label{R22-eq:intro-private-mask}
 Y_a=\rho Z+c\bigl(\sqrt{1-a^2/3}\,G+T_a\bigr),\qquad
 U_a=\left\{\frac{Y_a}{2ac}\right\},
\end{equation}
where braces denote fractional part. Conditional on $Z,G$, reduction modulo one makes $U_a$ uniform; it is therefore independent of $Z$. The variance correction fixes the full covariance exactly. The posterior-variance formula for the Gaussian mixture proves \eqref{R22-eq:intro-privacy-curvature} and both relative-entropy limits (Theorem~\ref{R22-thm:privacy-instability}).

The noise characteristic function has exact zeros at nonzero multiples of $\pi/a$, which escape to infinity as the Gaussian is approached. The classical uniform-dither cancellation~\cite{R22-Schuchman} explains the release. Wiener's translation theorem~\cite{R22-Wiener} gives a precise converse: for additive positive continuous noise and a private law with full support, zero-budget utility vanishes exactly when the noise characteristic function has no real zeros. Covariance-matched analytic examples attain every value $\log m$, $m\ge2$, and an arbitrarily small exponential tilt changes the infinite value to zero. These results identify both failures of semicontinuity within the same smooth uniformly log-concave class (Theorem~\ref{R22-thm:privacy-Fourier} and Corollaries~\ref{R22-cor:privacy-logm}--\ref{R22-cor:privacy-tilt}).

The uniform release in \eqref{R22-eq:intro-private-mask} also determines the complete decision order of the zero-budget private channels. Every such kernel factors through $U_a$:
\begin{equation}\label{R23-eq:intro-private-factor}
 \law(V\mid Y=y)=L_{\{y/(2ac)\}}
 \quad\hbox{for a Markov kernel }L.
\end{equation}
Thus one deterministic release is Blackwell-greatest among useful-data-only mechanisms: postprocessing it reproduces every private channel from $Y$, for every bounded decision problem about $Y$. The converse follows by deconvolving the Gaussian and differentiating the resulting constant moving average; privacy forces periodicity of the entire kernel (Theorem~\ref{R23-thm:private-greatest}).

The smooth finite-capacity useful-data-only examples have a different order. Each admits two deterministic releases attaining $\log m$ that are incomparable maximal private experiments, with no common private upper bound (Theorem~\ref{R23-thm:private-maximal}). 

\paragraph{\textbf{Greatest useful-data-only experiments.}}
For finite useful data, a greatest useful-data-only private experiment exists exactly when the feasible posterior polytope is a simplex (Proposition~\ref{R23-prop:privacy-simplex}). This criterion follows from the Rassouli--G\"und\"uz posterior representation and Blackwell comparison, connecting exact privacy capacity to the comparison theory developed in the final sections.

\subsection{The complete full-data Blackwell classification}
Allowing the encoder to observe both $Z$ and $Y$ changes which experiments can be greatest. We compare releases by postprocessing about $Y$, retaining exact privacy $U\perp Z$. Set $\mu_z=\law(Y\mid Z=z)$.
\begin{theorem}[Greatest full-data private releases]\label{R27-thm:intro-full-data}
For arbitrary standard Borel $(Z,Y)$, a greatest perfectly private full-data release exists exactly when, after discarding the conditional point masses, either the remaining $\mu_z$ almost surely coincide or they are almost surely supported on one common two-point set. Conditional point masses may have arbitrary locations. Under conditional atomlessness, existence is equivalent to $Y\perp Z$.
\end{theorem}
The conflict already occurs between linear and quadratic Gaussian prediction. For $Y=rZ+sE$, $0<|r|<1$, $s^2=1-r^2$ and independent standard Gaussians $Z,E$, both $E$ and $\sgn(Z)E$ are private. Optimal prediction of $Y$ forces a release to be Blackwell-equivalent to $E$, while $\sgn(Z)E$ predicts $Y^2$ strictly better. Hence they have no common private upper bound (Theorem~\ref{fdp:thm:gaussian}). The classification also gives a greatest useful-data-only release and no greatest full-data release for the smooth dependent sources of Theorem~\ref{R22-thm:intro-privacy}.

For the general classification, conditional quantiles optimize each scalar task. Simultaneous optimality for a countable determining family of ternary observables forces pairwise conditional couplings onto the diagonal or one common two-point support. Private refinement exposes these couplings without finite-state assumptions (Theorem~\ref{fdp:thm:quadratic}); Theorem~\ref{fdp:thm:classification} constructs the greatest mechanisms in both cases. Vector quadratic objectives give Wasserstein barycenters of $(\mu_z)$, with the all-optimizer correspondence and optimal-transport/fair-regression credit in Theorem~\ref{fdp:thm:vector}.

\subsection{Relation to previous work}
Banaszczyk's Gaussian-measure theorem gives $O(\sqrt{\log n})$ Koml\'os discrepancy~\cite{Banaszczyk}; constructive vector balancing relates this geometry to subgaussian signing laws and polynomial-time rounding~\cite{DGLN,BDGL}. Bansal--Jiang's affine spectral-independence method gives a polynomial-time $\widetilde O((\log n)^{1/4})$ bound for unit Euclidean columns and $O(\sqrt k)$ Beck--Fiala discrepancy when $k=\Omega((\log n)^2)$~\cite{BJ-ASI}. Their stochastic-calculus exposition states the Koml\'os bound as $O((\log n)^{1/4}(\log\log n)^{7/4})$~\cite[Theorem~1.1]{BJ-Exposition}; the logarithmic refinement in \cite{ErcanQuarter} builds on that framework. Affine spectral-independence controls how the evolving row constraints interfere with a covariance-selected random walk. The construction here encodes all columns in the initial full shear; coordinate rearrangement then retains independent reference blocks through the choice of signs. Guo--Fang--Lu develop directional-variation and polynomial-time constant-discrepancy constructions~\cite{GFL,GFLalg}; Karingula--Lovett prove an elementary constant bound~\cite{KL}. The joint law in Theorem~\ref{VB-thm:main}, its energy-sensitive entropy and the unrestricted near-Gaussian tradeoff state the simultaneous guarantees pursued here.

Akbas--Sra's reciprocal/Fisher and small-ball estimates~\cite{AkbasSraBSB,AkbasSraMatrix} are the common analytic inputs to the reference laws, Schatten theorems and volume-ratio theorem. The Parseval-frame laws retain exact means, a Gaussian reference and information bounds; the classical rank-one partition theorem is due to Marcus--Spielman--Srivastava~\cite{MSS}. Spectral sparsification~\cite{FC-BSS,FC-SHS} and randomized quadrature~\cite{KriegNovak,KunschRudolf,Kunsch} provide the approximation benchmarks. The corresponding sections explain the additional requirement here: a whole law supported on accurate approximations, with each original weight or sample marginal preserved exactly and a common information or convex-order bound.

Strassen~\cite{Strassen} represents fixed convex order by a martingale. Purification~\cite{DWW,Edwards1987,KhanRath2009} removes randomization under finite integral constraints. Sections~\ref{R11-sec:allocation}--\ref{R11-sec:monge} retain the entire output law and conditional means while making the terminal-to-initial map deterministic, proving the stated Nutz--Wang--Zhang approximation results~\cite{NWZ2024}. Brown's weak-operator approximation theorem~\cite{Brown1966} gives the background for Section~\ref{R11-sec:dynamics}, which preserves exact finite-observable identities at all iterates and in diffusion limits.

Hua--Song--Tudose's representation and Talagrand existence theorem~\cite{HST} underlie exact Gaussian extraction. Mazhar~\cite{R20-Mazhar} proves the equal-input threefold conclusion and constructive higher-mass theorem. The mixed inputs, compact endpoints, fourfold bound and Gaussian Banach extension are stated separately here. The finite-point formulation follows Johnston~\cite{R19-Johnston}; Borell's general inequality~\cite{R20-Borell} improves the integer summand counts. Song~\cite{Song} develops the representation and amplification connections.

For optimal martingale couplings, the existence and Bass representation of stretched Brownian motion are established in~\cite{CM-BBHK,CM-Bass}. Section~\ref{R11-sec:canonical} combines this structure with contraction~\cite{Caffarelli,Kolesnikov} to turn terminal curvature into a sharp same-driver deficit and conditional path geometry. Tschiderer's $q$-Bass duality and primal uniqueness~\cite{QB-Tsch} underlie the common-potential theorem; Acciaio--Marini~\cite{QB-AM} prove one-dimensional existence in the first-moment regime. The theorem here is all-dimensional in the second-moment regime, with a separate atomic-reference dual and an explicit potential-domain criterion. Backhoff--Pammer--Schachermayer's gradient flow~\cite{BF-BPS} is then calibrated using the additive-noise strong-convexity theorem.

Shannon's rate--distortion problem~\cite{R9Shannon} provides the operational setting for the entropy frontier. Exact product-source retention and conditional calibration identify its full transform here; exact entropy data then recover the reference law. Blackwell's comparison theorem~\cite{Blackwell} motivates a different inverse question for posterior distributions: which independent-source tests recover the decision order, and which costs recover the entire experiment? The refinement hierarchy and directed inverse answer these questions in Sections~\ref{R11-sec:Blackwell}--\ref{R11-sec:inverse}.

Rassouli--G\"und\"uz already establish unbounded Gaussian utility in the full-data privacy model~\cite[Corollary~6.1]{RG2021}. The finite-alphabet optimum, conditional-atomlessness criterion and single infinite-output construction with an arbitrarily small conditional-information budget are distinct conclusions here (Section~\ref{R11-sec:privacy}).

\subsection{Guide to the named results}
Tables~\ref{R11-tab:construction}--\ref{R11-tab:information} give every named conclusion, its scope and a direct proof route. Transport starts from prescribed marginals; the spectral and information theories take domains, sources and posterior laws as inputs. Their hypotheses determine which constructions can be combined.

\begingroup
\setlength{\tabcolsep}{4pt}
\renewcommand{\arraystretch}{1.12}
\begin{longtable}{@{}>{\raggedright\arraybackslash}p{.20\textwidth}>{\raggedright\arraybackslash}p{.60\textwidth}>{\raggedright\arraybackslash}p{.16\textwidth}@{}}
\caption{Construction, exact constraints and sampling.}\label{R11-tab:construction}\\
\toprule Problem or question & Mathematical conclusion & Statement and proof \\ \midrule
\endfirsthead
\multicolumn{3}{@{}l}{\textit{Construction, exact constraints and sampling (continued)}}\\[4pt]
\toprule Problem or question & Mathematical conclusion & Statement and proof \\ \midrule
\endhead
\midrule\multicolumn{3}{r@{}}{\textit{Continued on the next page}}\\\endfoot
\bottomrule\endlastfoot
\textbf{Koml\'os and Beck--Fiala} & Hard discrepancy below $C=\Cstar$, with a column-sensitive improvement, independent Gaussian reference blocks and more than $1.22946^n$ hard-balanced signings. The same law loses at most $0.973141\|A\|_F^2$ entropy. Incidence column degree $t$ gives $C\sqrt t-167/200$. & Theorems \ref{VB-thm:main}, \ref{MAIN-signing}\par pp.~\pageref{VB-thm:main}, \pageref{MAIN-signing}\\[7pt]
\textbf{Gaussian discrepancy endpoint} & \emph{Sharp order for all constructions.} The universal hard radius under $\sigma\cx(1+\varepsilon)\sqrt{\pi/2}\,G$ is $\Theta(\varepsilon^{-1/2})$. A one-row example gives the lower bound. & Theorem~\ref{VB-thm:tradeoff}\par p.~\pageref{VB-thm:tradeoff}\\[7pt]
\textbf{Gaussian endpoint entropy} & \emph{Sharp asymptotic constant.} The maximal entropy deficit per coordinate is $\sim\sqrt{\pi/3}\sqrt\varepsilon$ at Gaussian scale $(1+\varepsilon)\sqrt{\pi/2}$, attained asymptotically by code laws. & Theorem~\ref{VB-thm:entropy-endpoint}\par p.~\pageref{VB-thm:entropy-endpoint}\\[7pt]
\textbf{First-order directional-variation barrier} & Over all zero-extended $BV$ probability densities supported in $[-R,R]^m$, $\lim_m\inf_f\sup_{\|v\|=1}|D_vf|=\sqrt{2\pi}/R$. The uniform $(-3,3)$ auxiliary's first-order overlap certificate has sharp universal radius $3\sqrt{2\pi}$, including dependent densities. & Theorem~\ref{M-first-order}\par p.~\pageref{M-first-order}\\[7pt]
Universal-reference constant & The diffuse overlap problem permits dependent competitors selected before the direction. Its critical radius and the common-auxiliary threshold are separately enclosed between $6.8383231851$ and $6.8383231852$. & Theorem~\ref{R7-thm:method-constant}\par p.~\pageref{R7-thm:method-constant}\\[7pt]
Reis--Rothvoss Conjecture 1 & The named cases are square matrices for all $p\ge2$ and rectangular matrices for $p=2$. The precision-energy formula gives one signing for every $q\ge p$, with power $n^{1/2-1/p+1/q}$ in the square case and a rank-sensitive constant. A separate extension covers $1\le p\le2$. The hard dimension orders also follow from constant discrepancy and interpolation. & Theorem~\ref{R19-thm:Lebesgue}; Corollary~\ref{R19-cor:RR}\par p.~\pageref{R19-cor:RR}\\[7pt]
\textbf{Reis--Rothvoss Conjecture 2: Schatten signing} & For $n$ symmetric $n\times n$ matrices and $y\in[-1,1]^n$, one law with $\E\sigma=y$ satisfies $\|\sum_i(\sigma_i-y_i)A_i\|_{S_p}\le C n^{\max(1/2,1/p)}\max_i\|A_i\|_{S_p}$ at every outcome, simultaneously for $1\le p\le\infty$. The constant is universal and both powers are sharp. Support at most $n+1$; general dimensions incur the stated logarithmic factor. & Theorems \ref{RRS-thm:square-law}, \ref{RRS-thm:general-law}\par p.~\pageref{RRS-thm:square-law}\\[7pt]
\textbf{Reis--Rothvoss Conjecture 3: Gaussian square function} & One event of Gaussian measure $e^{-O(n)}$ controls every unitarily invariant norm at radius $\|(\sum A_i^2)^{1/2}\|_{\mathcal N}$ for $n$ symmetric $n\times n$ matrices. The stronger common weak singular-value majorization holds; general size $d$ has measure $\ge8^{-n}e^{-c_{\rm AS}d}$. & Theorem~\ref{RRS-thm:gaussian}\par p.~\pageref{RRS-thm:gaussian}\\[7pt]
\textbf{Reis's trace-class volume-ratio conjecture} & Every surjective real $M:\R^{a\times b}\to\R^r$ satisfies $\vr(MB_{S_1})\le C\sqrt{1+\min(a,b)/r}$. The case $a=b=r=d$ proves the conjecture. The proof retains a determinant bound and the dual outer-volume bound for spectrahedra. & Theorems \ref{SVR-thm:main}, \ref{SVR-thm:determinant}; Corollary~\ref{SVR-cor:Reis}\par p.~\pageref{SVR-thm:main}\\[7pt]
\textbf{Sharp bounded differences} & For independent inputs and coordinate widths $c_i$, $f-\E f\cx N(0,\frac\pi8\sum_i c_i^2)$, with optimal $\pi/8$ for all convex tests. \par\emph{Sharp adaptive Poisson law.} If $0\le X_i\le c$ and $\sum_i\E[X_i\mid\mathcal F_{i-1}]\le L$ pathwise, the centered sum is dominated by $c\{\operatorname{Pois}(L/c)-L/c\}$, the least universal upper law. & Theorems \ref{SP-thm:mcdiarmid}, \ref{SP-thm:poisson}\par p.~\pageref{SP-thm:mcdiarmid}\\[7pt]
\textbf{Gaussian integer rounding} & \emph{Sharp constructive bound.} Unbiased integer rounding with optimal covariance factor $\pi/8$ and an exact expected polynomial-bit rational sampler. Its input model is given in the theorem. & Theorem~\ref{thm:exact-gaussian-integer}\par p.~\pageref{thm:exact-gaussian-integer}\\[7pt]
\textbf{Integral base polyhedra} & \emph{General constrained sampler.} Arbitrary integral submodular $f$: unbiased nearest-cell rounding in $B(f)$, every tight face retained, Gaussian covariance $\frac\pi3P_H$, and an expected polynomial-bit value-oracle sampler with positive variance slack. & Theorem~\ref{BASE-thm:main}\par p.~\pageref{BASE-thm:main}\\[7pt]
\textbf{Exact-marginal sampling} & \emph{One law for every later integrand.} Exact individual marginals, a joint Hilbert reference and an $N^{-2}$ feature cost plus the residual variance divided by $N$. & Theorem~\ref{VB-thm:sampling}\par p.~\pageref{VB-thm:sampling}\\[7pt]
Sobolev integration & \emph{Optimal rates with added constraints.} One equal-weight law for all $s\le s_0$: exact uniform marginals and polynomial moments, hard measurable features, and optimal confidence and mean-square rates. & Theorem~\ref{thm:quadrature}\par p.~\pageref{thm:quadrature}\\[7pt]
\textbf{Matrix Spencer; Weaver} & \emph{Simultaneous guarantees.} Hard operator control, near-endpoint Gaussian comparison and entropy loss $O(d\sqrt\varepsilon)$ hold together, with optional scalar constraints. Uses Akbas--Sra; the rank-one existence theorem is due to MSS. & Theorem~\ref{R3-thm:matrix}; Corollary~\ref{R3-cor:Weaver}\par p.~\pageref{R3-thm:matrix}\\[7pt]
Sachdeva--Thudi--Zhao Conjecture 21 & Arbitrary-start matrix partial coloring retains exact mean $y$, more than $m-d$ saturated coordinates, hard error $K_\varepsilon s$, and Gaussian covariance $\frac{2\pi}{3}(1+\varepsilon)^2P_I$ when $\|\sum_iA_i^2\|\le s^2$. The Akbas--Sra full-coloring input already gives the hard conclusion; hereditary dyadic rounding here retains the mean and law comparison. & Theorem~\ref{R19-thm:matrix-partial}\par p.~\pageref{R19-thm:matrix-partial}\\[7pt]
Spectral graph signing & \emph{Same-sign constraints.} Signed adjacency norm $O_\varepsilon(\sqrt\Delta)$ and all signed degrees $O_\varepsilon(1)$, with the joint reference and entropy guarantee. & Corollary~\ref{R7-cor:graph}\par p.~\pageref{R7-cor:graph}\\[7pt]
\textbf{Unbiased spectral sparsifiers} & \emph{Entire sparse law.} Every outcome is a spectral approximation, each original weight is unbiased, the support has linear dimension order, and the Poisson information cost is independent of thinning depth. A continuous version retains the exact expected measure. & Theorem~\ref{FC-thm:poisson-sparse}\par p.~\pageref{FC-thm:poisson-sparse}\\[7pt]
\end{longtable}

\begin{longtable}{@{}>{\raggedright\arraybackslash}p{.20\textwidth}>{\raggedright\arraybackslash}p{.60\textwidth}>{\raggedright\arraybackslash}p{.16\textwidth}@{}}
\caption{Geometry, transport and spectral consequences.}\label{R11-tab:geometry}\\
\toprule Problem or question & Mathematical conclusion & Statement and proof \\ \midrule
\endfirsthead
\multicolumn{3}{@{}l}{\textit{Geometry, transport and spectral consequences (continued)}}\\[4pt]
\toprule Problem or question & Mathematical conclusion & Statement and proof \\ \midrule
\endhead
\midrule\multicolumn{3}{r@{}}{\textit{Continued on the next page}}\\\endfoot
\bottomrule\endlastfoot
\textbf{Talagrand Problem 2.1} & For symmetric $A$ with $\gamma_d(A)>2/3$, four summands contain a symmetric convex body of mass $>0.822656449$, answering the balanced $3/4$-input, $3/4$-output problem. Compact input admits the non-strict endpoint. Mazhar's constructive result uses six summands; four give existence here. & Corollary~\ref{R15-cor:fourfold}\par p.~\pageref{R15-cor:fourfold}\\[7pt]
\textbf{Talagrand convexity} & Three Borel sets with total Gaussian mass $>2$ have an undilated sum containing a convex body of mass $>1/2$. Three is optimal for the universal unweighted half-mass conclusion. The equal-input conclusion is also proved by Mazhar; the mixed inputs, strict output and compact endpoint are explicit here. Uses HST and exact extraction. & Theorems \ref{R19-thm:intro-Talagrand}, \ref{GC-thm:three}\par p.~\pageref{GC-thm:three}\\[7pt]
\textbf{Talagrand convexity below half mass} & Balanced Borel input of mass $\ge5/12$: five undilated copies contain compact symmetric convex output of mass $\ge\Phi(1/16)>0.524917$ in every separable Gaussian Banach space, including Wiener space. An explicit Cameron--Martin enlargement remains inside the fivefold sum. \par\emph{Half-mass input.} Five copies of compact balanced input give output mass $\ge3/4$; six give $\ge\beta_4>0.822656449$, and $6j$ give the full Gaussian $S$-profile. Explicit counts cover every positive compact input mass. & Theorems \ref{R22-thm:fivefold}, \ref{R23-thm:five-enlargement}; Corollaries \ref{R23-cor:five-mass}, \ref{R22-cor:sixfold}, \ref{R22-cor:positive-mass}\par p.~\pageref{R23-thm:five-enlargement}\\[7pt]
\textbf{Wiener and Gaussian Banach spaces} & Three input masses summing to $>2$ give a compact convex subset of the actual sum, of measure at least $1/2$; weighted sums are also covered. This includes pathwise Wiener sums. \par\emph{Green Problem 54.} In Gaussian product space, three copies of Borel input of mass $>2/3$ give product-compact convex output of mass $\ge1/2$. & Theorem~\ref{GC-thm:banach}; Corollaries \ref{GC-cor:wiener}, \ref{R20-cor:green}\par p.~\pageref{GC-thm:banach}\\[7pt]
Johnston's finite-point formulation & For $\gamma_d(A)>2/3$, $2\operatorname{conv}_3(A)$ contains a convex body of mass $>1/2$. The unweighted theorem gives the stronger $k=3$, $\varepsilon=1/4$ instance of Talagrand's Conjecture~1.1 as stated by Johnston; the dilation-free Question~1.2 is distinct. Explicit counts cover every fixed input threshold $\alpha>1/2$. & Corollaries \ref{R19-cor:conv-three}, \ref{R19-cor:explicit-counts}\par p.~\pageref{R19-cor:conv-three}\\[7pt]
\textbf{Gaussian hitting and convex-order laws} & A compact support hitting every bounded open convex set of Gaussian mass $\ge1/2$ carries a law $\mu\cx\gamma_d$. Scale one is sharp, even for finite symmetric supports. Banaszczyk's theorem then gives comparison variance $25$ for unit columns; Theorem~\ref{VB-thm:main} gives $6.5$ with a hard bound, independent coefficient reference and entropy. & Theorem~\ref{GC-thm:hitting}; Proposition~\ref{R10-prop:Gaussian-scale}\par pp.~\pageref{GC-thm:hitting}, \pageref{R10-prop:Gaussian-scale}\\[7pt]
Fixed-reference extraction & \emph{Exact threshold.} Positive hitting by all compactly supported dominated laws forces a convex subset of mass $\kappa(\nu)$. For Gaussians this is $1/2$; the dimension-uniform log-concave value is sharply $1/e$. & Theorem~\ref{EXT-thm:fixed}\par p.~\pageref{EXT-thm:fixed}\\[7pt]
Gaussian convex order in continuous time & Gaussian domination is equivalent to a continuous martingale from the dominated law to that Gaussian, with its covariance as instantaneous cap. Stopping-time residuals retain conditional Gaussian comparison. \par\emph{Bounded-time Skorokhod embedding.} In dimension one, $\mu\cx N(0,1)$ iff a Brownian motion started with law $\mu$ reaches $N(0,1)$ by a stopping time $\tau\le1$ in the stated enlarged filtration. & Corollaries \ref{GC-cor:gauss-dynamic}, \ref{GC-cor:skorokhod}; Theorem~\ref{GC-thm:residual}\par p.~\pageref{GC-cor:gauss-dynamic}\\[7pt]
\textbf{All-marginal interior stability} & A compact initial law, full-span terminal law and an exact positive product component give adapted $W_p$ approximation as both marginals vary, along every prescribed initial coupling. The approximating pairs may be unbounded or reducible. Finite-horizon, continuous-cost and weak-cost consequences retain their stated hypotheses. & Theorem~\ref{R20-thm:intro-stability}; Theorem~\ref{AS-thm:main}; Corollaries \ref{AS-cor:history}, \ref{AS-cor:optimization}\par p.~\pageref{AS-thm:main}\\[7pt]
\textbf{Stability of both endpoints and driver} & Under the preceding interior condition in $\mathcal P_2$, both marginals and the reference may vary in $W_2$. If the limiting driver is absolutely continuous, every asymptotically optimal coupling converges in adapted $W_2$ to the unique limiting optimizer, along every prescribed initial coupling. Atomic driver approximants are allowed; terminal curvature is unnecessary. & Corollaries \ref{R20-cor:intro-three-laws}, \ref{AS-cor:three-laws}\par p.~\pageref{AS-cor:three-laws}\\[7pt]
\textbf{Stretched Brownian motion} & \emph{Sharp stability.} Uniform terminal log-concavity gives $L$-Lipschitz conditional maps and sharp same-input squared-error deficit constant $1/(2L)$, without irreducibility. Centered, possibly atomic $q$ gives conditional-kernel error $\le2L W_2(q,\gamma_d)$. Upper-curvature bounds give non-Gaussian conditional contraction, sharp deficit and reference variation for arbitrary initial laws. & Theorems \ref{CM-thm:contraction}, \ref{CM-thm:deficit}, \ref{R14-thm:reference-curvature}; Corollaries \ref{R13-thm:reference-rate}, \ref{R15-cor:variation}\par p.~\pageref{CM-thm:contraction}\\[7pt]
\textbf{Conditional future-path geometry} & \emph{Sharp whole-path guarantee.} At every stopping time, Gaussian isoperimetry and $T_2$ with sharp coefficient $2L^2(1-\tau)$ hold in uniform path distance. Conditional one-time laws also satisfy log-Sobolev and Poincar\'e inequalities. & Theorems \ref{CM-thm:conditional}, \ref{CM-thm:path}\par p.~\pageref{CM-thm:path}\\[7pt]
\textbf{General-reference $q$-Bass existence} & \emph{Exact universal existence criterion.} For every second-moment reference, one globally finite subgradient potential represents all optimizers of every irreducible second-moment pair with full terminal affine span iff the reference has full convex support. Nondegeneracy gives $q$-Bass gradients. The single pair $(\delta_0,\gamma_d)$ detects necessity. Extended potentials cover arbitrary support; compact targets give global finiteness for individual pairs. \par\emph{Fixed signing law.} Every strict terminal inflation permits any admissible driver while retaining every signing probability. & Theorems \ref{QB-thm:main}, \ref{R24-thm:q-Bass-universal}; Proposition~\ref{R13-prop:finite-domain}; Corollary~\ref{QB-cor:signing}\par pp.~\pageref{QB-thm:main}, \pageref{R24-thm:q-Bass-universal}\\[7pt]
\textbf{General-reference dual attainment} & \emph{Atomic references included.} For an irreducible pair and any reference with finite second moments, one shifted subdifferential represents every optimizer. For $q$ equivalent to Lebesgue measure, a globally finite common-gradient representation exists exactly for irreducible pairs. Specifying a reference martingale gives an optimal path realization. & Theorems \ref{QB-thm:arbitrary}, \ref{R13-thm:reference-process}\par p.~\pageref{QB-thm:arbitrary}\\[7pt]
\textbf{Noise-regularized optimal transport} & \emph{Exact geometric criterion.} Continuous positive-density noise with a finite second moment makes maximal covariance strongly convex on bounded $L^\infty$ latent sets modulo translations iff the second-moment target has full affine span. Atomic targets are included; no martingale constraint is imposed. & Theorem~\ref{BF-thm:maincurvature}\par p.~\pageref{BF-thm:maincurvature}\\[7pt]
\textbf{Bass gradient flow} & \emph{Backhoff--Pammer--Schachermayer conjecture resolved.} Exponential convergence in every dimension from bounded initialization, for irreducible pairs with compact interior initial support. Singular and unbounded second-moment targets are allowed. & Theorem~\ref{BF-thm:flow}; Corollary~\ref{BF-cor:published}\par p.~\pageref{BF-thm:flow}\\[7pt]
\textbf{NWZ Conjectures 5.4 and 5.3} & \emph{Scalar conjectures resolved.} Backward-Monge directional couplings are dense in $W_\infty$ with exact marginals and order. Supermartingale couplings are dense in $W_1$, and in $W_p$ under finite $p$-moments. & Theorems \ref{FP-thm:directional-density}, \ref{FP-thm:supermartingale-density}\par p.~\pageref{FP-thm:directional-density}\\[7pt]
\textbf{NWZ Conjecture 5.1: countable initial law} & \emph{Specified higher-dimensional clause resolved.} In every finite dimension, $W_\infty$ backward-Monge martingale approximation retains both marginals and any prescribed finite list of bounded measurable costs. & Theorem~\ref{FP-thm:discrete-density}\par p.~\pageref{FP-thm:discrete-density}\\[7pt]
Full-span and common-fiber density & Full-relative-span conditional supports give $W_\infty$ backward-Monge approximation. A measurable full-span seed gives weak and $W_p$ density under the stated moments. For common fibers, terminal atomlessness and recoverability of the label are retained; avoidance of shared closures gives a measurable sufficient criterion. & Theorems \ref{SP-thm:multidimensional}, \ref{WC-thm:fiber-martingale}; Corollary~\ref{R19-cor:separated-labels}\par p.~\pageref{R19-cor:separated-labels}\\[7pt]
\textbf{NWZ 5.1: component disintegration} & \emph{Boundary obstruction.} With $\nu_I=\law(Y\mid I(X)=I)$ atomless, a unique compact martingale in $\R^4$ still forces a backward Bernoulli choice. Each component admits a constant backward map. Every strict contraction of the initial law or inflation of the terminal law restores backward determinism in this example. The component label need not be terminal-recoverable. & Theorem~\ref{BM-thm:main}; Corollaries \ref{BM-cor:NWZ}, \ref{BM-cor:quantitative}\par p.~\pageref{BM-thm:main}\\[7pt]
Continuous deterministic realization & \emph{Exact full output law.} Atomless source and full conditional affine span: deterministic realization preserves the entire output law and finite conditional moments, with arbitrary $W_\infty$ pair approximation. Strict contractions of martingale targets are covered in every dimension. & Theorems \ref{SP-thm:main}, \ref{SP-thm:contraction}\par p.~\pageref{SP-thm:main}\\[7pt]
Canonical polynomial partitions & \emph{Uniqueness and optimal degree.} Positive feasibility gives a unique cost-minimizing partition and gauge-fixed dual; selecting degree $r+1$ is optimal in dimension at least three. Regularization retains all moments and has an exact interface expansion. & Theorem~\ref{R9-thm:canonical-moments}\par p.~\pageref{R9-thm:canonical-moments}\\[7pt]
Deterministic diffusion approximation & \emph{Exact finite observations and path limits.} Exact invariant maps matching four conditional moment orders converge in path space to stationary OU and compact elliptic diffusions. Finite Gaussian polynomial action also holds at every horizon. \par\emph{All-horizon Gaussian records.} Every prescribed finite degree agrees exactly, including multitime tests, while the entropy contraction coefficient is $1$ and every improving $L^p\to L^q$ norm is infinite. & Theorems \ref{WC-thm:gaussian}, \ref{WC-thm:multitime}, \ref{WC-thm:compact-diffusion}; Corollary~\ref{WC-cor:no-hyper}\par p.~\pageref{WC-thm:compact-diffusion}\\[7pt]
\textbf{Cheeger inequality and Dirichlet spectra} & \emph{Equality throughout a subexponential spectral range.} For mixed weighted products, $h_D\sim2\sqrt{E_N}$ and $\lambda_{k_N,p}^{1/p}\sim(2/p)\sqrt{E_N}$ whenever $\log k_N=o(N)$. Nonconvex factor domains and arbitrary proportions are allowed. & Theorem~\ref{R9-thm:mixed-weighted}\par p.~\pageref{R9-thm:mixed-weighted}\\[7pt]
Smirnov--Vershynin information capacity & For every nonempty bounded open domain, $I(K)^2=4\Lambda(K)$ with an exact matrix-valued Dirichlet dual and an attained Sobolev relaxation of the smooth-density infimum. Cube and ball values and the product maximum rule are explicit. Convex domains also have Brunn--Minkowski concavity of $1/I$. Disjoint unions of $s$ components in dimension $d$ have the sharp factor $1/\sqrt{\min(s,d)}$ bound. The existing thinning theorem uses destination $2K$. & Corollary~\ref{R20-cor:capacity}; Theorems \ref{GEOM-thm:product-ball}, \ref{GEOM-thm:bm}; Proposition~\ref{R22-prop:components}\par p.~\pageref{R20-cor:capacity}\\[7pt]
Geometry beyond Fisher information & \emph{Complete-score dependence.} Equal Fisher matrices and all three-coordinate score marginals can coexist with retained widths differing by order $\sqrt m$. Even all scalar score moments can agree while sharp support radii differ. & Theorems \ref{S-thm:r7-compact-pair}, \ref{R9-thm:score-moment-indeterminacy}\par p.~\pageref{S-thm:r7-compact-pair}\\[7pt]
Compact-orbit comparison & \emph{Exact criterion.} An expected orbit support function characterizes comparison for every compact orthogonal action. Matrix spectra and singular values give the existing real and complex criteria; specialized maximizers retain their endpoint rigidity. & Theorem~\ref{R15-thm:orbit}; Section~\ref{R11-sec:orbits}\par p.~\pageref{R11-sec:orbits}\\[7pt]
\end{longtable}

\begin{longtable}{@{}>{\raggedright\arraybackslash}p{.20\textwidth}>{\raggedright\arraybackslash}p{.60\textwidth}>{\raggedright\arraybackslash}p{.16\textwidth}@{}}
\caption{Entropy, dependence and statistical experiments.}\label{R11-tab:information}\\
\toprule Problem or question & Mathematical conclusion & Statement and proof \\ \midrule
\endfirsthead
\multicolumn{3}{@{}l}{\textit{Entropy, dependence and statistical experiments (continued)}}\\[4pt]
\toprule Problem or question & Mathematical conclusion & Statement and proof \\ \midrule
\endhead
\midrule\multicolumn{3}{r@{}}{\textit{Continued on the next page}}\\\endfoot
\bottomrule\endlastfoot
\textbf{Complete privacy capacity} & In the full-data model, conditional atomlessness characterizes full finite-alphabet utility. For every $b>0$, one uniform private release has infinite information, $I(Z;U\mid Y)\le b$, and whole-output prediction probability $\ge e^{-b}$. \par\emph{Sharp Gaussian budget threshold.} For scalar Gaussians with $0<|\rho|<1$, utility is zero at $b=0$ and infinite for $b>0$. & Corollaries \ref{SP-cor:privacy-iff}, \ref{R15-cor:private-prediction}, \ref{SP-cor:gaussian-budget}; Theorem~\ref{SP-thm:infinite}\par p.~\pageref{SP-cor:privacy-iff}\\[7pt]
\textbf{Rassouli--G\"und\"uz Conjecture 1} & \emph{Exact capacity and nonattainment.} In the full-data observation model, the perfectly private $M$-alphabet utility supremum is $\log M$, unattained for each $M\ge2$ and $0<|\rho|<1$. & Theorem~\ref{FP-thm:gaussian-nonattainment}\par p.~\pageref{FP-thm:gaussian-nonattainment}\\[7pt]
\textbf{Smooth-data privacy instability} & At every nondegenerate correlated Gaussian scalar pair, analytic uniformly log-concave laws with exactly the same covariance converge in $W_2$, both relative-entropy directions and uniform potential Hessian, while zero-budget utility is infinite and the Gaussian value is zero. One deterministic uniform output from $Y$ alone attains infinity. \par\emph{Exact values and both discontinuity directions.} Fourier nonvanishing characterizes zero utility for the stated additive-noise model. Smooth covariance-matched laws attain every $\log m$; small tilts destroy infinite utility. & Theorems \ref{R22-thm:privacy-instability}, \ref{R22-thm:privacy-Fourier}; Corollaries \ref{R22-cor:privacy-logm}, \ref{R22-cor:privacy-tilt}\par p.~\pageref{R22-thm:privacy-instability}\\[7pt]
\textbf{Blackwell-greatest useful-data-only releases} & In the useful-data-only model, every zero-budget private kernel in the smooth infinite-capacity examples factors through one deterministic modulo release. The smooth $\log m$-capacity examples admit incomparable maximal releases and no greatest private experiment. \par\emph{Finite useful data.} A greatest private experiment exists exactly when the feasible posterior polytope is a simplex; all maximal experiments are its vertex-supported posterior laws. & Theorems \ref{R23-thm:private-greatest}, \ref{R23-thm:private-maximal}; Proposition~\ref{R23-prop:privacy-simplex}\par p.~\pageref{R23-thm:private-greatest}\\[7pt]
\textbf{Greatest full-data private releases} & For standard Borel sources, discard deterministic conditional useful laws. A greatest full-data release exists iff the remaining laws coincide almost surely or share a common two-point support. Conditional atomlessness reduces the criterion to independence. Ternary quadratic tasks detect it. For $0<|\rho|<1$, two explicit Gaussian private releases have no common private upper bound. & Theorems \ref{fdp:thm:classification}, \ref{fdp:thm:quadratic}, \ref{fdp:thm:gaussian}\par p.~\pageref{fdp:thm:classification}\\[7pt]
\textbf{Blackwell comparison} & \emph{Strict hierarchy and completion.} Every successive finite refinement level is strict; all levels are exactly complete. For $q$ states, quantitative recovery has exponent $2/(q+1)$ and yields one approximate simulator for all bounded decision problems. & Theorems \ref{BW-thm:completion}, \ref{R8-thm:blackwell-rate}\par p.~\pageref{BW-thm:completion}\\[7pt]
\textbf{Posterior reconstruction} & \emph{Complete directed inverse.} Varying the source weights of one fixed directed cycle cost recovers the whole closed-simplex posterior law, including boundary atoms. On $q$ states the inverse has sharp exponent $1/q$ on the closed simplex, with sharper face and endpoint bounds. & Theorems \ref{CYC-thm:boundary}, \ref{CYC-thm:sharp}\par p.~\pageref{CYC-thm:boundary}\\[7pt]
Rate--distortion and exact calibration & \emph{Exact attained frontier.} The product-reference entropy optimum is the conjugate of $\E\log\cosh(tW)$, attained while preserving the full source and all prescribed conditional means. & Theorem~\ref{R9-thm:frontier}\par p.~\pageref{R9-thm:frontier}\\[7pt]
\textbf{Entropy determines the source} & \emph{Complete inverse.} Exact entropy measurements on scales accumulating at one positive point determine the whole normalized symmetric law. Only a first moment is required; atoms are allowed. & Theorem~\ref{SP-thm:entropy-inverse}\par p.~\pageref{SP-thm:entropy-inverse}\\[7pt]
Nonlinear calibration sparsification & \emph{All amplitudes, one sparse reference.} A Laplacian comparison controls the entire calibration potential for every symmetric integrable source. $O(N/\varepsilon^2)$ latent states retain the exact target marginal and uniform nonlinear bounds. & Theorem~\ref{NLI-thm:nonlinear-sparse}\par p.~\pageref{NLI-thm:nonlinear-sparse}\\[7pt]
Exact latent realization & \emph{Characterization and exact complexity.} A reversible transition $P$ is a posterior two-step transition iff $D_pP$ is completely positive; its cp-rank is the minimum finite latent cardinality. Exact complexity can be quadratic when approximation is linear. & Theorem~\ref{NLI-thm:cp}\par p.~\pageref{NLI-thm:cp}\\[7pt]
Fixed common and coordinate noises & \emph{Complete feasibility boundary.} A directing-law factorization characterizes all-dimensional finite-alphabet feasibility and the attained entropy rate. For $aZ_0+bZ_i$, feasibility, maximal entropy and minimal dependence are explicit. & Theorems \ref{r3:thm:factorization}, \ref{r3:thm:gaussianboundary}\par p.~\pageref{r3:thm:factorization}\\[7pt]
Reference-induced dependence & \emph{Exact limit sets.} With varying common noise, every possible empirical-mean limit is characterized. Small smooth coordinate-density perturbations realize arbitrary closed permitted sets and can force consensus while fixed dimensions still permit independence. & Theorem~\ref{R8-thm:empirical-limits}\par p.~\pageref{R8-thm:empirical-limits}\\[7pt]
Local sampling dynamics & \emph{Exponential obstruction for stationary local processes.} Full-support, high-entropy laws with many independent marginals admit fast global sampling but force exponentially slow stationary local processes, including hidden-state lifts and arbitrary memory. The Gaussian and actual cosine references have separate constructions. & Section~\ref{R11-sec:local}\par p.~\pageref{R11-sec:local}\\[7pt]
\end{longtable}

\endgroup



\clearpage
\setcounter{tocdepth}{1}\tableofcontents

\clearpage
\part{Hard rounding and retained references}
\label{R6-part:construction}
The support must satisfy the hard constraint while conditional means retain the entire reference. We construct one output law with both properties, then determine the universal cost of making its Gaussian comparison sharp.

\section{Geometric rounding and calibrated references}
\label{R11-sec:construction}
Coordinate rearrangement finds a feasible alphabet point for each convex test. Keeping the test's epigraph fixed preserves its cost at the critical height; separation then gives one law for all tests. The shear $(X,T)\mapsto(X-AT,T)$ places every column in the initial density. Its heights survive rearrangement, so columns incur no repeated loss. The scalar calibration below and Appendix~\ref{app:overlap} verify those initial heights.

\subsection{Rounding a conditional mean}\label{sec:coordinate}\label{F-sec:joint}
Let $B\subset\R^d\times\prod_{j=1}^n I_j$ be open and convex, with
each $I_j$ a bounded open interval. Let $(Y,T)$ have an integrable
probability density $F$, supported on $B$, and finite mean $(p,\mu)$.
Extend densities by zero. For $z_{-j}$ denoting all coordinates except
$t_j$, define
\begin{equation}\label{eq:height} h_j(F)=\frac14\int\!\int_{I_j}\!\int_{I_j}
 \min\{F(z_{-j},s),F(z_{-j},t)\}\,ds\,dt\,dz_{-j}.
\end{equation}
A finite alphabet $\mathcal A_j\subset I_j$ has gap $g_j$ if every
open subinterval of $I_j$ of length greater than $g_j$ meets
$\mathcal A_j$. In particular, $\{-1,1\}\subset(-3,3)$ has gap two.

\begin{theorem}[Coordinate rounding]\label{thm:coordinate}
If $h_j(F)\ge g_j/2$ for every $j$, there is a law on
$N\in\prod_j\mathcal A_j$ and a coupling with $(Y,T)$ such that
\begin{equation}\label{eq:coordinate-conclusion}
 (p,N)\in B\quad\text{at every outcome},\qquad
 \E[Y\mid N]=p,\qquad \E[T\mid N]=N.
\end{equation}
Thus $(p,N)\cx(Y,T)$ and $\E N=\mu$. Equality in the height
hypotheses is allowed.
\end{theorem}

The factor $1/4$ in \eqref{eq:height} converts a squared superlevel
length into a first absolute moment. The next calculation explains
both the threshold and its preservation under successive operations.
Write $\mathfrak S_j$ for symmetric decreasing rearrangement in $t_j$
and for Steiner symmetrization of the supporting set. Layer cake gives
\begin{equation}\label{eq:rearrangement}
 \int\min\{a g^*,b h^*\}\ge\int\min\{ag,bh\},\qquad
 \|g^*-h^*\|_1\le\|g-h\|_1.
\end{equation}
Indeed, centered intervals maximize the intersection at each level.
The same computation gives
\begin{equation}\label{eq:height-rearrangement}
 h_j(\mathfrak S_iF)\ge h_j(F)\ (i\ne j),\qquad
 \E_{\mathfrak S_jF}|T_j|=h_j(F).
\end{equation}
Both sides of the second identity integrate one quarter of the squared
length of a superlevel set. Later rearrangements preserve the entire
marginal of an earlier coordinate. The supporting set stays open and
convex and gains the corresponding reflection symmetries.

We also need a vertical section estimate. Define
$\mathsf V_z(F)=\sup_{a\ne0}\|F-\tau_{ae_z}F\|_1/|a|$.
\begin{lemma}\label{lem:vertical}
If a finite-first-moment density on an open convex set has mean
$(p,m)$ in coordinates $(y,z)$ and
$0<\mathsf V_z(F)\le v<\infty$, the set contains
$\{(p,m+t):|t|<1/v\}$.
\end{lemma}
\begin{proof}
A nonnegative integrable line density $g$ has
$\|g\|_\infty\le\mathsf V(g)/2$. A probability density bounded by
$M$ on $(0,\infty)$ has mean at least $1/(2M)$, by integrating the
bound $\Prb(X\le x)\le Mx$. Thus a probability density $g$ on $(0,\infty)$ satisfies
\[
 \left(\int_0^\infty xg(x)\,dx\right)\mathsf V(g)\ge1.
\]
For a supporting inequality $a^{\mathsf T}y+bz<c$ with $b\ne0$,
the positive distance has variation at most $v/|b|$, since translation
distances contract under marginalization. Its mean is therefore at
least $|b|/v$. If $b=0$, its mean distance is strictly positive.
Every asserted point satisfies all supporting inequalities strictly;
separation proves the lemma. Uniform density on $(0,2L)$ attains the
constant, with mean $L$ and variation $1/L$.
\end{proof}

\paragraph{Retaining a cost at the critical height.}
The threshold has no strict reserve to spend on a perturbation of the horizontal density. We therefore keep the cost's epigraph fixed and vary only the probability placed inside it. A uniform vertical interval of length $2L$ adds mean height $L$ and has vertical variation $1/L$. The section estimate recovers that entire added height. Meanwhile the horizontal marginal approaches the unlifted density in $L^1$, so its rearranged mean converges inside the same fixed epigraph. Continuity there preserves the cost inequality at equality in the height assumptions.

\begin{proof}[Proof of Theorem~\ref{thm:coordinate}]
Rearrange in all $t_j$ and fold them by absolute value. The resulting
mean is $(p,M)$, where $M_j\ge h_j(F)\ge g_j/2$ by
\eqref{eq:height-rearrangement}. It belongs to
$B^*=\mathfrak S_n\cdots\mathfrak S_1B$: an integrable law supported
on an open convex set has its mean in that set, as separation shows.
Convexity and the separate reflection symmetries imply
$(p,g_1/2,\ldots,g_n/2)\in B^*$.

Reverse the symmetrizations. The last centered fibre contains $g_n/2$
and is open, so its length exceeds $g_n$. Its original fibre has the
same length and therefore contains $N_n\in\mathcal A_n$. Repeat at
the resulting base point in the preceding coordinate. This produces
an element of the finite set
\[
 \mathcal S=\{N\in\textstyle\prod_j\mathcal A_j:(p,N)\in B\}.
\]

To prove the convex-cost inequality for $\Phi$ with finite expectation, keep the fixed
open convex epigraph
\[
 D=\{(y,t,z):(y,t)\in B,\ z>\Phi(y,t)\}
\]
and give it density
\begin{equation}\label{eq:epigraph-density}
 F_L(y,t,z)=F(y,t)g_L(z-\Phi(y,t)),\qquad
 g_L=(2L)^{-1}\mathbf1_{(0,2L)}.
\end{equation}
Its vertical mean is $\E\Phi+L$ and its vertical variation is $1/L$.
Horizontal rearrangement commutes with vertical translation and
contracts its $L^1$ distance, so this variation bound survives all
rearrangements and folding. Moreover,
\begin{equation}\label{eq:epigraph-limit}
 \|F_L-Fg_L\|_1
 =\E\min\{2,|\Phi(Y,T)|/L\}\longrightarrow0.
\end{equation}
The rearrangement contractions imply that the folded auxiliary mean
$M^{(L)}$ tends to $M$. Here the auxiliary intervals are bounded.
Lemma~\ref{lem:vertical} puts
\[
 (p,M^{(L)},\E\Phi+\varepsilon)\in D^*
 \qquad(0<\varepsilon<2L),
 \quad D^*=\mathfrak S_n\cdots\mathfrak S_1D.
\]

The projection of $D^*$ is $B^*$. To see this, the horizontal sections
of $D$ increase to $B$ as $z\to\infty$, and fibre lengths show that
Steiner symmetrization commutes with increasing unions. Thus $D^*$
is the epigraph of a finite convex function on $B^*$. Its lower
finiteness follows from the minimum of $\Phi(y,\cdot)$ on the compact
closure of the auxiliary box. This function is continuous at
$(p,M)\in B^*$. Passing to the limit in \eqref{eq:epigraph-limit}
therefore gives
$(p,M,\E\Phi+\varepsilon)\in D^*$ for every $\varepsilon>0$.
Lower $M$ to $(g_j/2)_j$ and reverse the symmetrizations inside this
same epigraph, with $z$ fixed. The finite set $\mathcal S$ consequently
satisfies
\begin{equation}\label{eq:one-convex-cost}
 \min_{N\in\mathcal S}\Phi(p,N)\le\E\Phi(Y,T).
\end{equation}

For a finite list of convex tests, separation of attainable expectation
deficits from the negative orthant would give a nonnegative linear
combination contradicting \eqref{eq:one-convex-cost}. Some law on
$\mathcal S$ therefore satisfies the whole list. Compactness of its
probability simplex gives one law for all tests. The affine tests fix
its mean, and Strassen's theorem gives \eqref{eq:coordinate-conclusion}.
\end{proof}

The threshold is sharp: for the sign alphabet, uniform
$T$ on $(-2a,2a)$ has height $a$ and $\E|T|=a<1$ when $a<1$;
the convex test $|t|$ excludes every sign law. More generally, rescale
this example into an alphabet gap.

\subsubsection{The affine constraint}
Let $X$ have a density $f$ on an open convex $K\subset\R^m$, mean
$p$, and finite first moment. Independently let $T_j\sim q_j$ be
supported on $I_j$, with means $\mu_j$. Define
\begin{equation}\label{F-eq:H-weighted} H_f(v;q)=\frac14\iint\!\int
 \min\{q(s)f(y+sv),q(t)f(y+tv)\}\,dy\,ds\,dt.
\end{equation}
\begin{corollary}[Affine rounding]\label{thm:affine}
If $H_f(v_j;q_j)\ge g_j/2$ for every column of $V$, there is a law
on $N\in\prod_j\mathcal A_j$ with
\begin{equation}\label{eq:affine-identity} \E[X\mid N]=p+V(N-\mu),\qquad \E[T\mid N]=N,
 \qquad p+V(N-\mu)\in K.
\end{equation}
If the references and alphabets are centrally symmetric, the law may
be chosen symmetric.
\end{corollary}
\begin{proof}
The shear $(Y,T)=(X-VT,T)$ has density
$F(y,t)=f(y+Vt)\prod_jq_j(t_j)$ on
$B=\{(y,t):t\in\prod_jI_j,\ y+Vt\in K\}$ and mean
$(p-V\mu,\mu)$. Changing $y$ by $\sum_{i\ne j}t_iv_i$ in
\eqref{eq:height} gives exactly $h_j(F)=H_f(v_j;q_j)$.
Apply Theorem~\ref{thm:coordinate} and undo the shear. Averaging a law
with its reflection proves the symmetric assertion.
\end{proof}

All column conditions are evaluated on the original product reference.
The single shear retains the other coefficient variables while each
coordinate is rounded.

\subsubsection{A quantitative inner copy of the fractional section}
\label{F-sec:interior-means}
For the sign alphabet $\mathcal A_j=\{-1,1\}\subset(-3,3)$,
strict height gives a quantitative inclusion for the whole section. Let
$P_p=\{t:(p,t)\in B\}$ and
$\mathcal S_p=P_p\cap\{-1,1\}^n$ in the coordinate theorem.
\begin{theorem}[A homothetic section in the signing polytope]
\label{F-thm:homothetic-hull}\label{F-thm:strict-mean}
Put $h=\min_j h_j(F)>1$. For $h^{-1}\le c<1$,
\begin{equation}\label{F-eq:homothetic-hull}
 c\mu+(1-c)P_p\subset\operatorname{conv}\mathcal S_p.
\end{equation}
In particular, $\mu$ is an interior point of this convex hull. For the
full shear, the assertion becomes
\[
 c\mu+(1-c)\{t\in(-3,3)^n:p+V(t-\mu)\in K\}
 \subset\operatorname{conv}\mathcal T_\mu,
\]
where $\mathcal T_\mu=\{\sigma:p+V(\sigma-\mu)\in K\}$.
No boundedness of $K$ is needed. The factor $1-h^{-1}$ in
\eqref{F-eq:homothetic-hull} is optimal, for each $1<h\le3/2$.
\end{theorem}
\begin{proof}
Fix $t_0\in P_p$ and contract the entire source about $(p,t_0)$ by $c$.
Convexity keeps its support in $B$. Its mean is
$(p,c\mu+(1-c)t_0)$ and its coordinate heights are exactly $c h_j(F)$:
the two sliced integrations leave one additional factor $c$ after the
density Jacobian cancels. The endpoint theorem gives a law on
$\mathcal S_p$ with the displayed mean. Taking $t_0$ throughout $P_p$
proves the inclusion. Since $\mu\in P_p$, it also proves the interior.

For sharpness take $d=0$, $n=1$, $B=(-3,3)$, and a uniform density on
$(3-4h,3)$. Its height is $h$ and its mean is $\mu=3-2h$.
The right endpoint of $c\mu+(1-c)B$ is $3-2hc$, which exceeds one
as soon as $c<h^{-1}$. Thus no larger homothetic factor works uniformly.
\end{proof}
For example, if $\mu+rB_2^\circ\subset P_p$, the good-signing polytope
contains $\mu+(1-h^{-1})rB_2^\circ$.

\begin{corollary}[Linear objectives and an entropy-maximizing law]
\label{F-cor:mean-objective}
Every nonzero coefficient objective has a good signing whose value
exceeds its value at $\mu$. A law with mean $\mu$ has a representation
with at most $n+1$ outcomes, with rational probabilities for rational
$\mu$. Among all laws on $\mathcal S_p$ with this mean, the unique
one of largest Shannon entropy has probabilities proportional to
$e^{\theta^{\mathsf T}\sigma}$, for a unique $\theta\in\R^n$.
\end{corollary}
\begin{proof}
The interior gives the objective assertion and full affine span.
Carath\'eodory's theorem gives the support bound; an affinely independent
support gives rational weights by a rational linear system.
The function $\log\sum_{\sigma\in\mathcal S_p}
 e^{\theta^{\mathsf T}(\sigma-\mu)}$ is strictly convex. An inner ball
at $\mu$ makes it coercive. Its critical equation is the prescribed
mean equation, and relative entropy proves the maximizing property.
\end{proof}

The law selected in Corollary~\ref{F-cor:mean-objective} has the stated
support and mean. The reference-compatible law is constructed by
Theorem~\ref{thm:coordinate}; the entropy maximization here is over all
feasible laws with that mean.

\subsubsection{The terminal section and recovery of signs}
The same lifted set will be used to retain regions of starting points.
For the sign alphabet, write
\begin{equation}\label{R6-eq:full-shear}
 B=\{(y,t):t\in(-3,3)^n,\ y+Vt\in K\},\qquad
 F(y,t)=f(y+Vt)\prod_jq_j(t_j).
\end{equation}
For $1\le\eta<3$, define the buffered Banaszczyk transform
\[
 \calS_v^{(\eta)}K=
 ((K-\eta v)\cap(K+\eta v))+(-(3-\eta),3-\eta)v,
 \qquad \calS_v=\calS_v^{(1)}.
\]
If $B^*=\mathfrak S_n\cdots\mathfrak S_1B$, the fibre identity is
\begin{equation}\label{R6-eq:full-shear-section}
 \{y:(y,\eta_1,\ldots,\eta_n)\in B^*\}
 =\calS_{v_n}^{(\eta_n)}\cdots\calS_{v_1}^{(\eta_1)}K.
\end{equation}
An open fibre in $(-3,3)$ reaches $\eta$ after centering precisely when
it contains $[a-\eta,a+\eta]$ for some $|a|<3-\eta$.
The two endpoints require $w\pm\eta v\in K$, and the position of the
interval permits the horizontal displacement $av$. This proves the
identity for one coordinate; fixing the other coordinates and repeating
proves \eqref{R6-eq:full-shear-section}.

For backward recovery, write $y=w+tv$ in $\calS_v^{(\eta)}K$.
Choose $\epsilon\in\{-1,1\}$ with $|t+\epsilon|\le1$. Convexity gives
\[
 y+\epsilon v+[-(\eta-1),\eta-1]v\subset K.
\]
Thus each terminal point determines a signing whose image lies in $K$. Applying the
full shear to a prefix gives the same assertion in every column order.
Section~\ref{S-sec:score} will place two explicitly determined convex
regions in this terminal section. The conditional-mean law of
Corollary~\ref{thm:affine} and those retained regions come from the same
sequence of rearrangements.


\subsection{A reference on the cube}\label{sec:reference}
For $R>0$, let $X_R\in\R^m$ have independent coordinates with density
\begin{equation}\label{MAIN-eq:cosine} f_{R,1}(x)=R^{-1}\cos^2\!\left(\frac{\pi x}{2R}\right)
 \mathbf1_{(-R,R)}(x),\qquad f_{R,m}=f_{R,1}^{\otimes m}.
\end{equation}
Its support imposes the desired hard coordinate bound. The question
in Corollary~\ref{thm:affine} is now whether each initial height is
at least one.

\subsubsection{The first-order estimate}
For the uniform auxiliary $q_0=\tfrac16\mathbf1_{(-3,3)}$, put
$\omega_f(v)=\|f-f(\cdot+v)\|_1$.
Since $\int\min(f,f(\cdot+v))=1-\omega_f(v)/2$,
\begin{equation}\label{F-eq:uniform-height}
 H_f(v;q_0)=\frac32-\frac1{24}\int_0^6(6-t)\omega_f(tv)\,dt
 \ge\frac32-\frac32\|\partial_v f\|_{\TV}.
\end{equation}
Here $\|\partial_v f\|_{\TV}=\sup_{t\ne0}\omega_f(tv)/|t|$,
which includes boundary jumps of a zero extension.

For $X\sim f_{R,m}$, the variables
$U_i=\tan(\pi X_i/(2R))$ are independent, with density
$2/[\pi(1+u^2)^2]$. They have the representation
$U_i=G_i/\sqrt{W_i}$, where the $G_i$ are standard Gaussians and
the $W_i$ are independent chi-square variables with three degrees of
freedom, independent of the $G_i$. The gamma integral proves the
representation and $\E W_i^{-1}=1$. Thus
\begin{equation}\label{F-eq:cosine-score-L1}
 \|\partial_vf_{R,m}\|_1
 =\frac\pi R\sqrt{\frac2\pi}
       \E\sqrt{\sum_i v_i^2/W_i}
 \le\frac{\sqrt{2\pi}}R\|v\|_2.
\end{equation}
The normalized all-ones direction and the central limit theorem show
that the dimension-uniform coefficient is sharp. Equations
\eqref{F-eq:uniform-height}--\eqref{F-eq:cosine-score-L1} and the affine
rounding theorem already give
\begin{equation}\label{eq:elementary-joint}
 \|A\sigma\|_\infty<3\sqrt{2\pi},\qquad
 (A\sigma,\sigma)\cx(X_{3\sqrt{2\pi}},U),
 \quad U_j\sim\Unif(-3,3),\quad U\perp X.
\end{equation}

\subsubsection{The weighted overlap}
The exact overlap in \eqref{F-eq:H-weighted} retains the size of
each coordinate of a column. Set
\begin{equation}\label{VB-eq:constants}
 C=\Cstar,\qquad \gamma=\gammaprof,\qquad
 S_3(v)=\frac{\sum_i|v_i|^3}{\|v\|_2^3}\quad(v\ne0),
\end{equation}
and, for $A\ne0$, put
\begin{equation}\label{MAIN-eq:radius}\csname ltx@label\endcsname{MAIN-eq:radius-intro}
 R_A=\max_{a_j\ne0}\|a_j\|_2\bigl(C-\gamma S_3(a_j)\bigr).
\end{equation}
There is a fixed even density $q=q_\dagger$ on $(-3,3)$ for which
\begin{equation}\label{eq:scalar-input}
 H_{f_{C-\gamma S_3(u),m}}(u;q)>1\qquad(\|u\|_2=1).
\end{equation}
Its degree-ten logarithm is given in \eqref{R-eq:htrial};
Appendix~\ref{app:overlap} derives \eqref{eq:scalar-input} from the
cosine likelihood transform and gives the finite scalar certificates.
The density is independent of $A,m,n$.

\begin{theorem}[The bounded joint reference]\label{MAIN-signing}
For $A\ne0$ and every $R>0$ such that
$H_{f_{R,m}}(a_j;q)\ge1$ for all $j$, a symmetric sign law satisfies
\begin{equation}\label{MAIN-stationarity}
 \E[X_R\mid\sigma]=A\sigma,\qquad
 \E[T\mid\sigma]=\sigma,\qquad \|A\sigma\|_\infty<R,
 \quad \law(X_R,T)=f_{R,m}\otimes q^{\otimes n}.
\end{equation}
In particular $R=R_A$ is admissible. For unit columns,
$R_A\le C-\gamma/\sqrt m$; if each column has at most $k$ nonzero
coordinates, $\sqrt m$ may be replaced by $\sqrt k$.
\end{theorem}
\begin{proof}
Use the sign alphabet in Corollary~\ref{thm:affine}. The scalar input
and scaling give height greater than one at $R_A$; translation
overlap for the log-concave cosine density increases with $R$.
The zero-column height is greater than one as shown in
Appendix~\ref{app:overlap}. Finally H\"older's inequality gives
$S_3(u)\ge m^{-1/2}$, or $k^{-1/2}$ on a support of size $k$.
\end{proof}

The numerical radius enters only through the displayed height
hypothesis. Every conclusion below uses the chosen reference at its
admissible radius. Unequal row budgets follow by replacing $A$ by
$D^{-1}A$ for a positive diagonal $D$ and then undoing that map.
For a $0$--$1$ incidence matrix of column degree at most $t$,
\eqref{MAIN-eq:radius} gives discrepancy less than $C\sqrt t-\gamma$,
a square-root bound in the Beck--Fiala setting \cite{BeckFiala}.

We record the reference quantities needed later. Let $T_1\sim q$ and
put
\begin{equation}\label{eq:reference-constants} v_0=\frac13-\frac2{\pi^2},\quad
 \alpha=\frac\pi2\left(\frac12-\frac2{\pi^2}\right)^2,\quad
 w=\E T_1^2<2.784448,\quad
 \tau=\frac\pi2(\E|T_1|)^2<3.23.
\end{equation}
Then
\begin{equation}\label{eq:reference-gaussian}
 (X_R)_i\cx N(0,\alpha R^2),\qquad T_1\cx N(0,\tau),
 \qquad \Var((X_R)_i)=v_0R^2.
\end{equation}
The Gaussian variances are optimal for the two scalar laws, as shown
in Appendix~\ref{app:reference-inequalities}. The cosine law also has
its exact variance as a subgaussian parameter:
\begin{equation}\label{eq:cosine-laplace}
 \E e^{t(X_R)_i}
 =\frac{\pi^2\sinh(Rt)}{Rt(\pi^2+R^2t^2)}
 =\prod_{k=2}^\infty\left(1+\frac{R^2t^2}{\pi^2k^2}\right)
 \le e^{v_0R^2t^2/2}.
\end{equation}
The inequality is $\log(1+x)\le x$ and
$\sum_{k\ge2}(\pi k)^{-2}=1/6-1/\pi^2$.
Thus Gaussian convex-order covariance and the smaller variance proxy
for exponential moments remain explicitly distinct.

\subsection{Numerical calibration and column-dependent radius}
A different optimality problem fixes a universal reference before an unknown vector direction. Its competitors may be fully dependent. The diffuse value is characterized exactly, and the two thresholds satisfy
\begin{equation}\label{VB-eq:reference-thresholds}
 6.8383231851<R_{\rm diff}\le R_{\rm com}<6.8383231852.
\end{equation}
Theorem~\ref{R7-thm:method-constant} states this reference-before-direction optimization and its quantifiers. The first-order directional-variation certificate has sharp radius $3\sqrt{2\pi}$ over all admissible BV densities (Theorem~\ref{M-first-order}). Equation~\eqref{MAIN-eq:radius} gives the column-dependent radius; the full likelihood profile sharpens this cubic bound.

The explicit existential comparison follows Theorem~\ref{VB-thm:main}; Section~\ref{R25-sec:intro-applications} gives the separate computational comparison. The exact rational integer sampler and positive-slack base-polyhedron sampler have the polynomial-bit input models in Table~\ref{R9-tab:runtime}.


\section{The joint law: Gaussian projection, entropy and optimal discrepancy}
\label{R11-sec:joint}
The graph relation gives equivalent reference representations of a query. Minimizing their variance strengthens the joint comparison and later controls sampling residuals. Replica energy bounds the entropy cost. At the least Gaussian scale, conditional Jensen forces independent fair signs; the one-row obstruction and bounded auxiliary determine the cost of approaching that endpoint.

\subsection{Linear relations and independent references}\label{sec:projection}
The graph relation in \eqref{MAIN-stationarity} determines how the two
independent reference blocks combine. The same linear calculation explains the sampling estimates. Projecting the Gaussian onto that relation adds their precision matrices; this retains information that separate bounds on the two marginals would discard.

\begin{proposition}[Gaussian projection]\label{prop:projection}\label{r3:thm:short}\label{GAUSS-prop:subspace}
Let $V,D$ be positive definite and let $A$ be a real matrix. Put
$Q=(D^{-1}+A^{\mathsf T}V^{-1}A)^{-1}$. For any integrable $Z$,
\begin{equation}\label{eq:graph-equivalence}
 (AZ,Z)\cx N(0,\diag(V,D))
 \quad\Longleftrightarrow\quad Z\cx N(0,Q).
\end{equation}
The same covariance follows from the variational identity
\begin{equation}\label{eq:quadratic-infimum}
 t^{\mathsf T}Qt=
 \min_u\{u^{\mathsf T}Vu+(t-A^{\mathsf T}u)^{\mathsf T}
                      D(t-A^{\mathsf T}u)\}.
\end{equation}
More generally, if $LZ=0$ and $\Sigma\succeq0$, then
\[
 Z\cx N(0,\Sigma)\ \Longleftrightarrow\
 Z\cx N(0,\Sigma-\Sigma L^{\mathsf T}
                  (L\Sigma L^{\mathsf T})^\dagger L\Sigma).
\]
\end{proposition}
\begin{proof}
Set $P=I-\Sigma L^{\mathsf T}(L\Sigma L^{\mathsf T})^\dagger L$.
Then $PZ=Z$ and $P\Sigma P^{\mathsf T}$ is the displayed covariance.
Applying $P$ proves the forward comparison; adding independent centered
Gaussian noise proves its converse. Take $L(x,z)=x-Az$ to obtain
\eqref{eq:graph-equivalence}. Completing the square gives
\eqref{eq:quadratic-infimum}; its minimizer is
$(V+ADA^{\mathsf T})^{-1}ADt$.
\end{proof}

For scalar covariance blocks, every representation
$t=A^{\mathsf T}u+(t-A^{\mathsf T}u)$ splits the query between
independent sources. The orthogonal projection selects the smallest
reference covariance. In particular, the law in
Theorem~\ref{MAIN-signing} satisfies simultaneously
\begin{align}
 \sigma&\cx N(0,[\tau^{-1}I+(\alpha R^2)^{-1}A^{\mathsf T}A]^{-1}),
                   \label{eq:sign-gaussian}\\
 \Cov(\sigma)&\preceq
 [w^{-1}I+(v_0R^2)^{-1}A^{\mathsf T}A]^{-1}.
                   \label{eq:sign-covariance}
\end{align}
Tensorize \eqref{eq:reference-gaussian} for the first assertion. For
the second apply conditional variance directly to
\eqref{eq:opening-regression} and use
\eqref{eq:quadratic-infimum}. The second assertion concerns the actual
covariance and uses the exact reference variances.

\subsubsection{Information and the number of feasible signings}\label{INFO-sec}
Let $\psi_X(u)=\log\E e^{\langle u,X_R\rangle}$ and
$I_X=\psi_X^*$, where $f^*(x)=\sup_u(\langle u,x\rangle-f(u))$.
Set
\begin{equation}\label{INFO-dq}
 d_q=\sup_{a\in\R}\{a-\E\log\cosh(aT_1)\}.
\end{equation}
For a fair sign $\epsilon$ with $\E[T_1\mid\epsilon]=\epsilon$,
write $b(t)=\E[\epsilon\mid T_1=t]$. Its mutual information is
$\E\mathfrak i(b(T_1))$, where
\[
 \mathfrak i(s)=\tfrac{1+s}{2}\log(1+s)+
               \tfrac{1-s}{2}\log(1-s),\qquad |s|\le1.
\]
Since $\mathfrak i^*=\log\cosh$ and $\E[T_1b(T_1)]=1$,
Fenchel's inequality gives $I(\epsilon;T_1)\ge d_q$.
Equality is attained by $b(t)=\tanh(a_*t)$, where
$\E[T_1\tanh(a_*T_1)]=1$. Existence and uniqueness follow from
$\E|T_1|>1$ and strict monotonicity.

Independence of the references adds these scalar costs. Indeed,
relative-entropy superadditivity against their product law gives
\[
 \Sh(\sigma)\ge I(\sigma;X_R,T)
 \ge I(\sigma;X_R)+\sum_j I(\sigma_j;T_j).
\]
The conditional physical mean is $A\sigma$, so the Gibbs variational
inequality bounds the first term below by $\E I_X(A\sigma)$.
Consequently
\begin{equation}\label{INFO-entropy}
 \Sh(\sigma)\ge nd_q+\E I_X(A\sigma),\qquad
 \sum_{\|As\|_\infty<R}e^{-I_X(As)}\ge e^{nd_q}.
\end{equation}
The second inequality is the finite Gibbs inequality applied to the
signing law. The scalar estimate in Appendix~\ref{app:reference-inequalities}
is $d_q>0.2065768178565815>\log(1.22946)$.

For $R\le C_0=6.8383231852$, the scalar enclosures give
$\alpha R^2<6.5$ and $\tau<3.23$. Tensorization and addition of
independent Gaussian noise give \eqref{VB-eq:main-joint}, while
Proposition~\ref{prop:projection} gives \eqref{VB-eq:main-precision}.
The squared-column-mass bound will be obtained in
Section~\ref{ENERGY-sec} while retaining the same reference; that
law proves all the assertions of Theorem~\ref{VB-thm:main}.
If $A=0$, independent fair signs suffice, because conditioning
$\sqrt{\pi/2}\,G_n$ on its coordinate signs gives those signs.

\subsubsection{Other convex supports}
The hard bound originates in the support of $X$, so the affine theorem
also applies to a convex body before any Gaussian comparison is made.
\begin{corollary}\label{thm:body}
Let $K\subset\R^m$ be open and convex, with
$\theta=\gamma_m(K)>0$, and set $p=\E[G\mid G\in K]$ for a standard
Gaussian $G$. If $\|v_j\|_2\le\theta\sqrt{\pi/2}/3$, there is a
mean-zero sign law with
\[
 p+V\sigma\in K,\qquad
 (V\sigma,\sigma)\cx N(0,\diag(I_m,\tfrac{9\pi}{8}I_n)).
\]
The two Gaussian blocks are independent. Its coefficient covariance
in convex order is
$[8/(9\pi)I_n+V^{\mathsf T}V]^{-1}$.
\end{corollary}
\begin{proof}
Every line section of $\varphi_m\mathbf1_K$ is a Gaussian density
restricted to an interval. Its variation is twice its maximum, so
for $f=\varphi_m\mathbf1_K/\theta$,
$\|\partial_vf\|_{\TV}\le\sqrt{2/\pi}\|v\|_2/\theta$.
Equation~\eqref{F-eq:uniform-height} and Corollary~\ref{thm:affine}
give $(V\sigma,\sigma)\cx(X-p,U)$, with
$X\sim(G\mid G\in K)$ and independent $U_j\sim\Unif(-3,3)$.
Harg\'e's theorem gives $X-p\cx G$ \cite{Harge}: the centered law has
log-concave density ratio $c e^{-\langle p,x\rangle}\mathbf1_{K-p}$
relative to Gaussian measure. The scalar comparison
$U_j\cx N(0,9\pi/8)$ follows from the two-crossing argument of
Appendix~\ref{app:reference-inequalities}. Tensorize and apply
Proposition~\ref{prop:projection}.
\end{proof}


\subsection{Energy-sensitive signing laws and counts}
\label{ENERGY-sec}
The complete source can be retained while replacing the ambient number
of columns in the entropy budget by their squared Euclidean mass.
The key step is to couple many signed copies through one signing,
then select one copy. This preserves the physical and coefficient
references as well as the hard constraint. The replica entropy argument of \cite{AkbasSraBSB} bounds the information
spent in selecting a copy. The conditional source identity below also
retains the full product reference.

Put $C_0=6.8383231852$, and let $q=q_\dagger$ be the even density with the
exact polynomial coefficients in \eqref{R-eq:htrial}. Write
\begin{equation}\label{ENERGY-eq:constants}
 d_q=\sup_{t\in\R}\{t-\E\log\cosh(tT_1)\},\quad
 \delta_q=\log2-d_q,\quad
 L_0=\frac{6455525558018172689}{31250000000000000000}.
\end{equation}
The scalar coefficient estimate in Section~\ref{INFO-sec} gives $d_q\ge L_0$. This section uses
the uncontracted universal cosine--auxiliary reference. The profile-dependent radius in Theorem~\ref{MAIN-signing}
applies to the law constructed there. Replicas here use the same
uncontracted reference for every signed copy.

For $r>0$, let $X_r$ be the product cosine reference on $(-r,r)^m$ and
let $I_r$ be its Cram\'er transform. Thus
\[
 \log\E e^{\langle u,X_r\rangle}
 =\sum_{i=1}^m\log\frac{\pi^2\sinh(ru_i)}
                   {ru_i(\pi^2+r^2u_i^2)},
\]
with the continuous value used when $u_i=0$.

Let $U_n$ be the uniform law on signs. For positive squared column
weights $w$, $F_B(w)$ denotes the fractional configuration optimum
with capacity $B$, defined in \eqref{PACK-eq:LP}; put
$\eta_B(w)=\sum_i1/\lfloor B/w_i\rfloor$.
Zero weights are omitted from both expressions.

\begin{theorem}[One calibrated law with an energy-dependent entropy deficit]
\label{ENERGY-thm:main}
Let $A\in\R^{m\times n}$ have columns $a_i$, and let
$r\ge C_0\max_i\|a_i\|_2$, $r>0$. Put
\[
 B=(r/C_0)^2,\qquad w_i=\|a_i\|_2^2,\qquad
 \tau=\sum_iw_i=\|A\|_F^2,
\]
and omit zero weights in $F_B$ and $\eta_B$. There is a symmetric law
$\mu$ on signs with
\begin{equation}\label{ENERGY-eq:reference}
 \|A\sigma\|_\infty<r,\qquad
 \E[X_r\mid\sigma]=A\sigma,\quad \E[T\mid\sigma]=\sigma,
 \qquad \law(X_r,T)=\law(X_r)\otimes q^{\otimes n},
\end{equation}
and
\begin{align}
 D(\mu\Vert U_n)+\E_\mu I_r(A\sigma)
 &\le\delta_q F_B(w)\le\delta_q\eta_B(w)\notag\\
 &\le\min\left\{\delta_q n,
            2\delta_q C_0^2\frac\tau{r^2}\right\}.
 \label{ENERGY-eq:budget}
\end{align}
In particular, for normalized columns at $r=C_0$,
\begin{equation}\label{ENERGY-eq:numeric}
 D(\mu\Vert U_n)+\E_\mu I_{C_0}(A\sigma)<0.973141\,\|A\|_F^2
 \quad\hbox{when }A\ne0.
\end{equation}
The weighted number of good signings satisfies
\begin{equation}\label{ENERGY-eq:weighted-count}
 \sum_{\|As\|_\infty<r}e^{-I_r(As)}
       \ge 2^n\exp[-\delta_q F_B(w)].
\end{equation}
\end{theorem}

\subsubsection{The exact replica entropy identity}\label{REP-sec}
Fix permutations $P_1,\ldots,P_k$ of $N$ coordinates. Independently choose
$kN$ fair signs $\xi_{\ell i}$, put
$G_\ell=\operatorname{diag}(\xi_{\ell1},\ldots,\xi_{\ell N})P_\ell$,
and write $\Xi$ for the array of masks. Conditional on $\Xi$, choose any
law for $Y\in\{-1,1\}^N$. Set $X_\ell=G_\ell Y$.
For random variables $X_1,\ldots,X_k$, their total correlation is
\[
 \operatorname{TC}(X_1,\ldots,X_k)
 =D\!\left(\law(X_1,\ldots,X_k)\,\middle\Vert\,
                  \bigotimes_{\ell=1}^k\law(X_\ell)\right).
\]

\begin{theorem}[Exact replica entropy identity]\label{REP-thm:identity}
With the preceding notation,
\begin{align}
 &\sum_{\ell=1}^kD(\law(X_\ell)\Vert U_N)
   +\operatorname{TC}(X_1,\ldots,X_k)
   +I(Y;X_1,\ldots,X_k)\notag\\
 &\hspace{25mm}=H(Y)-H(Y\mid\Xi)=I(Y;\Xi).
 \label{REP-eq:identity}
\end{align}
Consequently, if a nonnegative function $\mathcal L$ on signs satisfies
\begin{equation}\label{REP-eq:conditional-info}
 H(Y\mid\Xi=\xi)\ge Nd+
       \sum_{\ell=1}^k\E[\mathcal L(X_\ell)\mid\Xi=\xi]
 \quad\hbox{for every mask array }\xi,
\end{equation}
then a uniformly selected replica has law $\mu$ satisfying
\begin{equation}\label{REP-eq:budget}
 D(\mu\Vert U_N)+\E_\mu\mathcal L\le\frac Nk(\log2-d).
\end{equation}
The same inequality holds after deleting fixed zero-column coordinates,
provided $\mathcal L$ depends only on the remaining coordinates.
\end{theorem}
\begin{proof}
For each fixed $Y=y$, the map from the masks to
$(G_1y,\ldots,G_ky)$ is a bijection of $\{-1,1\}^{kN}$.
Thus the full map $(\Xi,Y)\mapsto(X_1,\ldots,X_k,Y)$ is bijective and
\[
 H(X_1,\ldots,X_k,Y)=kN\log2+H(Y\mid\Xi).
\]
Expanding relative entropy against $U_N^{\otimes k}\otimes\law(Y)$
gives
\[
 D(\law(X_1,\ldots,X_k)\Vert U_N^{\otimes k})
   +I(Y;X_1,\ldots,X_k)
 =H(Y)-H(Y\mid\Xi).
\]
Relative entropy to a product measure decomposes into the sum of
marginal relative entropies and total correlation, proving
\eqref{REP-eq:identity}. Average \eqref{REP-eq:conditional-info}, use
$H(Y)\le N\log2$, and drop the two nonnegative information terms.
Convexity of $D(\cdot\Vert U_N)$ proves the assertion for the uniform
mixture. Coordinate projection contracts relative entropy to the
corresponding uniform product measure.
\end{proof}

The identity records the two terms discarded by the estimate. Dependence
among different replicas and information retained about the shared signing
both reduce the sum of their entropy deficits. No independence of the
output replicas is assumed.

\paragraph{The source marginal after selecting a replica}
Let $A$ have $n\le N$ columns, padded by zero columns to form $\widehat A$.
Let $q$ be an even probability law and $\nu$ a centered law on $\R^m$,
both with finite first moments. Suppose that, for every mask array,
there is a joint coupling of $Y,Z_1,\ldots,Z_k,T$ such that
\begin{equation}\label{REP-eq:source-input}
 \law(Z_1,\ldots,Z_k,T\mid\Xi)
       =\nu^{\otimes k}\otimes q^{\otimes N},\qquad
 \E[Z_\ell\mid Y,\Xi]=\widehat A G_\ell Y,
 \quad \E[T\mid Y,\Xi]=Y.
\end{equation}
This is an equality of source laws for each mask array, even though the
coupling to $Y$ can depend on the masks.

\begin{proposition}[Retention of the complete selected reference]
\label{REP-prop:reference}
Choose $J$ uniformly in $\{1,\ldots,k\}$, independently of all preceding
variables, and let $\sigma$ be the first $n$ coordinates of $G_JY$.
Put $Z=Z_J$ and let $T'$ be the first $n$ coordinates of $G_JT$. Then
\begin{equation}\label{REP-eq:source-output}
 \law(Z,T')=\nu\otimes q^{\otimes n},\qquad
 \E[Z\mid\sigma]=A\sigma,\qquad \E[T'\mid\sigma]=\sigma.
\end{equation}
If every conditional law of $Y$ is symmetric, the law of $\sigma$ is
symmetric. If every replica has hard support in a prescribed set,
that support is retained as well.
\end{proposition}
\begin{proof}
For fixed $(\Xi,J)$, signed permutation invariance gives
$G_JT\sim q^{\otimes N}$. This variable is independent of $Z_J\sim\nu$,
so the selected pair has the same product distribution for every
$(\Xi,J)$. Mixing leaves that distribution unchanged. Conditional on
$(Y,\Xi,J)$, its mean is $(\widehat A G_JY,G_JY)$ after the designated
coordinate projection. Both coordinates of this mean depend on the
retained signing through $(A\sigma,\sigma)$. The tower property proves
\eqref{REP-eq:source-output}. The support and symmetry statements follow
for each conditional law before mixing.
\end{proof}

The same argument applies when additional source coordinates depend on
the mask array, as long as the displayed independent blocks have the
fixed conditional marginal in \eqref{REP-eq:source-input}. In that case
the conclusion retains precisely those blocks. This distinction will
be useful for the matrix barrier, whose coefficient density varies
with the replica array.

\begin{corollary}[Information amplification with one retained reference]
\label{REP-cor:amplification}
Assume \eqref{REP-eq:source-input}, symmetry, finite cumulants, and
$\E_q|T_1|>1$. Define
\[
 d_q=\sup_{a\in\R}\{a-\E_q\log\cosh(aT_1)\},\qquad
 \delta_q=\log2-d_q,\qquad I_\nu=(\log\E_\nu e^{\langle\cdot,Z\rangle})^*.
\]
The selected law has the joint reference \eqref{REP-eq:source-output}
and satisfies
\begin{equation}\label{REP-eq:amplified-cost}
 D(\law(\sigma)\Vert U_n)+\E I_\nu(A\sigma)
       \le\delta_q\frac Nk.
\end{equation}
\end{corollary}
\begin{proof}
The information inequality for independent source blocks gives, for
fixed masks,
\[
 H(Y\mid\Xi=\xi)\ge Nd_q+
        \sum_{\ell=1}^k\E[I_\nu(\widehat A G_\ell Y)\mid\Xi=\xi].
\]
Indeed, each auxiliary coordinate contributes at least $d_q$ by its
stationary fair-sign moment, and each physical block contributes its
Cram\'er rate at its conditional mean. This is the argument of
Section~\ref{INFO-sec}, applied conditionally to the larger product.
Theorem~\ref{REP-thm:identity} and Proposition~\ref{REP-prop:reference}
now give both conclusions for the same selected law.
\end{proof}

\subsubsection{Exact packing of the squared column norms}\label{PACK-sec}
The factor $N/k$ in \eqref{REP-eq:amplified-cost} is a packing problem.
Allowing zero columns lets us enforce the column budget exactly.

Fix positive weights $w_1,\ldots,w_n$ and a capacity $B\ge\max_iw_i$.
Treat the weights as item types, permitting repeated occurrences of a type
in one bin. A configuration is a nonzero vector $c\in\mathbb Z_{\ge0}^n$ with
$\sum_i c_iw_i\le B$. Denote the finite collection of configurations by
$\mathcal C_B$. Define the fractional bin-packing optimum
\begin{equation}\label{PACK-eq:LP}
 F_B(w)=\min\left\{\sum_{c\in\mathcal C_B}\lambda_c:
   \lambda_c\ge0,\quad
   \sum_{c\in\mathcal C_B}\lambda_cc_i=1\quad(1\le i\le n)\right\}.
\end{equation}
For an empty list of positive weights we set $F_B=0$. The equality
form in \eqref{PACK-eq:LP} has the usual fractional bin-packing value:
excess occurrences in a covering solution can be removed from bins,
because every subconfiguration remains feasible.

\begin{theorem}[Exact replica realization of fractional bin packing]
\label{PACK-thm:exact}
The number $F_B(w)$ is attained as $N/k$ by a finite arrangement of $k$
copies of each of the $n$ items in $N$ bins, with these properties:
each bin has weight at most $B$, and the $k$ copies can be labelled by
replicas so that each replica uses each original item exactly once and
uses any bin at most once. Equivalently, each replica is a permutation
of the original column array padded by $N-n$ zero columns.
Every such arrangement satisfies $N/k\ge F_B(w)$.
\end{theorem}
\begin{proof}
The constraint matrix of \eqref{PACK-eq:LP} and its right-hand side are
integral. An optimal basic feasible solution therefore has rational
coordinates, even when the original weights are real. Choose $k$ to
clear all their denominators, and make $k\lambda_c$ bins of type $c$.
Each item type occurs exactly $k$ times, and $N=kF_B(w)$. Multiplying
$k$ further ensures that each bin contains at most $k$ occurrences and
that $N\ge n$.

Form a bipartite multigraph with item types on the left, bins on the
right, and one edge for each occurrence. Left degrees are $k$ and right
degrees are at most $k$. Its edges split into $k$ matchings. For
completeness, add vertices and edges to make a $k$-regular bipartite
multigraph with equal part sizes. Hall's condition follows because the
$k|S|$ edges from any set $S$ enter vertices of total degree at most
$k|N(S)|$. A perfect matching exists. Remove it and repeat. Restricting
the resulting edge coloring to the original graph partitions its edges into $k$
matchings, each meeting every original left vertex once. The colors
are the required replica labels. Finally, an arbitrary replica
arrangement yields a feasible fractional solution by assigning to each
configuration its number of bins divided by $k$.
\end{proof}

\begin{proposition}[Explicit bounds and polynomial-size packing]
\label{PACK-prop:bounds}
Put $\tau=\sum_iw_i$, $r_i=\lfloor B/w_i\rfloor$, and
$\eta_B(w)=\sum_i r_i^{-1}$. Then
\begin{equation}\label{PACK-eq:bounds}
 \frac\tau B\le F_B(w)\le\eta_B(w)
       \le\min\left\{n,\frac{2\tau}{B}\right\}.
\end{equation}
For every $\varepsilon>0$ there is an explicit arrangement using
$k=\lceil n/\varepsilon\rceil$ replicas and
\begin{equation}\label{PACK-eq:finite-size}
 N=\sum_i\left\lceil\frac{k}{r_i}\right\rceil
       \le\min\{k\eta_B(w)+n,kn\}
\end{equation}
shared coordinates. Its ratio is at most $\eta_B(w)+\varepsilon$.
The arrangement is constructed in time polynomial in $n$, $k$ and
the bit length of rational input weights and capacity.
\end{proposition}
\begin{proof}
Sum the capacity inequalities against a feasible fractional solution
to obtain $\tau\le B\sum_c\lambda_c$. The configurations $r_ie_i$,
used with weights $1/r_i$, give $F_B\le\eta_B$. The inequalities
$r_i\ge1$ and $\lfloor x\rfloor\ge x/2$ for $x\ge1$ give the remaining
bounds.
For the explicit arrangement, treat each item type separately. Place
its occurrences labelled $1,\ldots,k$ in consecutive bins, with at
most $r_i$ occurrences per bin. Each replica then uses exactly one
bin of this type, and bins belonging to different types are distinct.
There are $\lceil k/r_i\rceil$ bins of type $i$, which gives
\eqref{PACK-eq:finite-size}; the last bound uses
$\lceil k/r_i\rceil\le k$ for each $i$. One can cap $r_i$ at $k$ when
constructing the bins, so tiny positive weights do not cause a large
intermediate enumeration.
\end{proof}

When every weight equals $B/r$ for an integer $r$, equality holds in
$F_B=\tau/B=n/r$. If every weight is slightly larger than $B/2$, every
bin holds at most one item, so $F_B=n$ and $F_B/(\tau/B)$ approaches two.
Thus the factor two in the energy bound is sharp for this packing
problem. Mixed configurations can improve the explicit homogeneous
bound: for weights $(9/25,16/25)$ and $B=1$, one bin containing both
items gives $F_B=1$, whereas $\eta_B=3/2$.

Zero-weight original items can be assigned distinct unused positions in
each replica. Replacing $N$ by $\max\{N,n\}$ preserves the bound
$N\le k\eta_B+n$, where $n$ now counts all original columns.

\begin{proof}[Proof of Theorem~\ref{ENERGY-thm:main}]
First assume every column is nonzero. Theorem~\ref{PACK-thm:exact}
provides $k$ replicas and $N$ bins with $N/k=F_B(w)$. Let $P_\ell$
be the padded permutation corresponding to replica $\ell$, and generate
the independent sign masks. Stack the $k$ matrices
$\widehat A G_\ell$ vertically to form $\mathcal A$.
Every column of $\mathcal A$ has squared Euclidean norm at most $B$:
its squared norm is exactly the load of the corresponding bin, because
sign masks do not change the load. The universal calibration applied
to $\mathcal A/\sqrt B$ gives a symmetric law for $Y$ with independent
physical cosine blocks of radius $C_0$ and auxiliary law $q^{\otimes N}$.
Multiply the physical blocks by $\sqrt B=r/C_0$. The resulting source
has $k$ independent cosine blocks of radius $r$, and
\[
 \E[Z_\ell\mid Y,\Xi]=\widehat A G_\ell Y,
 \qquad\E[T\mid Y,\Xi]=Y.
\]
Each replica has hard discrepancy strictly below $r$.
Corollary~\ref{REP-cor:amplification} gives the first inequality in
\eqref{ENERGY-eq:budget}, together with the exact reference
\eqref{ENERGY-eq:reference}. Proposition~\ref{PACK-prop:bounds} gives
the remaining inequalities.

If some columns vanish, apply the preceding construction to the
nonzero columns and append independent fair signs. Each appended sign
can be coupled to an independent $T_1\sim q$ using
\[
 \Pr(\epsilon=s\mid T_1=t)
 =\frac12\left(1+\frac{s\operatorname{sgn}(t)}{\E|T_1|}\right).
\]
The assumption $\E|T_1|>1$ makes this a probability kernel and gives
$\E[T_1\mid\epsilon]=\epsilon$. This preserves the product source and
adds zero relative entropy. When $A=0$, use independent fair signs
throughout. Finally, the finite Gibbs variational inequality gives
\[
 \log\sum_{\|As\|_\infty<r}e^{-I_r(As)}
 \ge H(\mu)-\E_\mu I_r(A\sigma),
\]
which proves \eqref{ENERGY-eq:weighted-count}.
The strict numerical bound follows from
$2(\log2-L_0)<0.973141$, certified by a rational remainder bound for
$\log2$.
\end{proof}

\subsubsection{Consequences of the same law}
Let $v_r=r^2(1/3-2/\pi^2)$,
$b_r=r\sqrt{\pi/2}(1/2-2/\pi^2)$,
$w_q=\E_qT^2$ and $b_q=\sqrt{\pi/2}\E_q|T|$.
Proposition~\ref{prop:projection} applied to
\eqref{ENERGY-eq:reference} gives
\begin{equation}\label{ENERGY-eq:precision}
 \Cov(\sigma)\preceq(w_q^{-1}I+v_r^{-1}A^{\mathsf T}A)^{-1},\qquad
 \sigma\cx N(0,[b_q^{-2}I+b_r^{-2}A^{\mathsf T}A]^{-1}).
\end{equation}
For $r=C_0$ the scalar bounds already proved in
Section~\ref{sec:projection} give $b_r^2<6.5$ and $b_q^2<3.23$.
Together with $d_q>\log(1.22946)$ and
\eqref{ENERGY-eq:numeric}, this proves all the assertions of
Theorem~\ref{VB-thm:main} for one law.

\begin{corollary}[Beck--Fiala with a complete reference]
\label{ENERGY-cor:BF}
Let $M$ be the incidence matrix of a set system in which each element
belongs to at most $t\ge1$ sets. There is a symmetric sign law with
\[
 \|M\sigma\|_\infty<C_0\sqrt t,\qquad
 (M\sigma,\sigma)\cx
 (N(0,6.5tI),N(0,3.23I)),
\]
where the Gaussian blocks are independent. The same law obeys
\[
 D(\law(\sigma)\|U_n)
 <0.973141\,\frac{\#\{\text{incidences}\}}t
\]
when $M\ne0$, and has entropy greater than $n\log(1.22946)$.
\end{corollary}
\begin{proof}
Apply Theorem~\ref{ENERGY-thm:main} to $A=M/\sqrt t$ and $r=C_0$,
then scale the physical block by $\sqrt t$.
\end{proof}
The square-root discrepancy order in the Beck--Fiala problem follows
from the Koml\'os conclusion; the same construction retains the joint
reference and the incidence-sensitive information bound.

\begin{corollary}[Nearly independent random coordinate blocks]
\label{ENERGY-cor:blocks}
For the law of Theorem~\ref{ENERGY-thm:main}, choose a uniform
$k$-element subset $S\subseteq[n]$, independently of $\sigma$. Then
\[
 \E_S D(\law(\sigma_S)\|U_k)
 \le\frac{k}{n}\delta_qF_B(w),\qquad
 \E_S\TV(\law(\sigma_S),U_k)
 \le\sqrt{\frac{k\delta_qF_B(w)}{2n}}.
\]
\end{corollary}
\begin{proof}
In the entropy chain rule for the coordinates in their natural order,
conditioning on earlier selected coordinates loses at most as much
entropy as conditioning on all earlier coordinates. Averaging gives
$\E_S H(\sigma_S)\ge(k/n)H(\sigma)$. Subtract from $k\log2$,
apply \eqref{ENERGY-eq:budget}, and use Pinsker and Cauchy--Schwarz.
\end{proof}

Uniform rejection from the cube has acceptance probability at least
$\exp[-\delta_qF_B(w)]$. For rational $A,r$, testing the hard condition
costs polynomial time per trial. Its expected number of trials is at
most $\exp(2\delta_q C_0^2\|A\|_F^2/r^2)$, hence polynomial when
$\|A\|_F^2/r^2$ is logarithmic in the input length. This algorithm
samples the uniform law on the good signings. The reference-coupled
law uses the replica construction, whose explicit packing has
polynomial overhead relative to a sampler for the underlying joint law.
The exact configuration optimum has a finite realization of possibly
larger size.


\subsection{The optimal discrepancy--Gaussian tradeoff}\label{R3-sec:height}\label{VB-subsec:endpoint}
The Gaussian scale of a fair sign is at least $\kappa=\sqrt{\pi/2}$,
because $\E|\sigma_i|=1$. We determine the order of the hard discrepancy
needed to approach this endpoint. The one-row obstruction below determines the necessary order without
any hypothesis on how the law is constructed. The remaining subsections
attain that order while retaining the independent reference.

\subsubsection{The one-row obstruction}
\begin{theorem}[A necessary inverse-square-root radius]
\label{R3-thm:radius-lower}\label{VB-eq:one-row-width}
Suppose $R\ge0$ and $\varepsilon>0$ have the following property:
for every matrix with Euclidean column norms at most one, some signing
law satisfies $\|A\sigma\|_\infty\le R$ and
$\sigma\cx(1+\varepsilon)\kappa G_n$. Then
\begin{equation}\label{R3-eq:radius-lower}
 R^2\ge\frac{\pi(1+\varepsilon)}{6144\varepsilon}-\frac38.
\end{equation}
Thus the exponent $\varepsilon^{-1/2}$ in
Theorem~\ref{R3-thm:hard-Gaussian} is optimal.
\end{theorem}
\begin{proof}
Choose an even $n$ with $8R^2\le n\le8R^2+2$, taking $n=2$ if
$R=0$, and let $A$ be the single row of $n$ ones. Put
$\mathcal V=\{y\in\{-1,1\}^n:|\sum y_i|\le R\}$.
If it is empty the required property already fails. Its support
function $h$ equals $n$ on $\mathcal V$, so any such law gives
$n\le(1+\varepsilon)\kappa\E h(G)$.

Let $S=\sum_i\operatorname{sgn}G_i$ and $a_0=\E|G_1|=\sqrt{2/\pi}$.
To bring the unconstrained maximizer into $\mathcal V$, at least
$k\ge(|S|-R)_+/2$ signs in the majority class must change.
Conditionally on the signs, the magnitudes in that class are
independent half-normal variables. Their distribution function is at
most $a_0t$, so the expected sum of their $k$ smallest values is at
least $k(k+1)/(2a_0(M+1))$, where $M\le n$ is the class size.
Consequently
\[
 n a_0-\E h(G)\ge
       \frac{\E(|S|-R)_+^2}{4a_0(n+1)}.
\]
Since $\E S^2=n$ and $\E S^4\le3n^2$, the Paley--Zygmund
inequality gives $\Prb(|S|\ge\sqrt{n/2})\ge1/12$.
On this event $(|S|-R)^2\ge n/8$. Therefore
\[
 \frac{\varepsilon}{1+\varepsilon}
 \ge\frac{n a_0-\E h(G)}{n a_0}
 \ge\frac1{384a_0^2(n+1)}
 =\frac{\pi}{768(n+1)}.
\]
Use $n+1\le8R^2+3$ and rearrange. The displayed lower constant is
explicit; the assertion of optimality concerns the exponent.
\end{proof}

\subsubsection{Excess height pays for the physical shear}
The lower bound used only the support of a sign law and a Gaussian
support-function test. To attain its order, we choose a coefficient
source whose first absolute moment is $1+\varepsilon$. This gives
an excess height $\varepsilon$ before the shear. The physical Fisher
cost scales as $R^{-2}$, so a quadratic height inequality can retain
the threshold one exactly when $R$ is of order
$\varepsilon^{-1/2}$. We now compute that inequality and its optimal
coefficient.

For a density $f$ with $\partial_w\sqrt f\in L^2$, put
$\mathcal J_w(f)=4\int|\partial_w\sqrt f|^2$. For an even density
$q$ on $(-a,a)$, write
\begin{equation}\label{R3-eq:heightdef}
 H_f(w;q)=\frac14\int_{-a}^a\int_{-a}^a\int
 \min\{q(s)f(x-sw),q(t)f(x-tw)\}\,dx\,ds\,dt.
\end{equation}
The threshold in Theorem~\ref{thm:affine} is one.
A decreasing auxiliary has $H_f(0;q)=\E_q|T|$.
The next calculation spends its excess height quadratically in the
physical score. That is the source of the inverse-square-root law.

\subsubsection{The sharp quadratic coefficient and its extremizer}
The auxiliary density determines the complete second-order height loss.
Its coefficient has an exact defect identity, which explains the
Gaussian shape and the extra cost of a bounded support.

Let $q(t)=Z^{-1}e^{-V(t)}\mathbf1_{\{|t|<a\}}$, where
$0<a\le\infty$, $V$ is even and twice continuously differentiable on
$(-a,a)$, $V'(t)>0$ for $t>0$, and $V''$ is nonnegative and nondecreasing
on $(0,a)$. Assume $q$ is a probability density with finite first
absolute moment $m=\E_q|T|$. Put
\begin{equation}\label{R4-eq:height-functional}
 K(q)=\int_0^a\frac{t^2q(t)}{V'(t)}\,dt
     =\int_0^a\frac{t^2q(t)^2}{-q'(t)}\,dt<\infty.
\end{equation}
We use the height and directional Fisher information defined above.

\begin{theorem}[A general sharp quadratic height inequality]
\label{R4-thm:general-height}
For every density with finite directional Fisher information,
\begin{equation}\label{R4-eq:general-height}
 H_f(w;q)\ge m-K(q)\mathcal J_w(f).
\end{equation}
Its coefficient is exact:
\begin{equation}\label{R4-eq:height-variation}
 \lim_{\lambda\to0}
 \frac{H_f(0;q)-H_f(\lambda w;q)}{\lambda^2}
 =K(q)\mathcal J_w(f).
\end{equation}
Moreover $T\cx\sqrt{\pi/2}\,mG_1$.
\end{theorem}
\begin{proof}
Use the difference $u=t-s=2h>0$ and midpoint $d=(s+t)/2$.
For fixed physical values $F=f(x)$ and $G=f(x-uw)$, put
$M=F+G$, $z=(F-G)/M$. Let $g_u(z)$ be the loss in the midpoint
integral divided by $M$; thus $g_u(0)=g_u'(0)=0$.
The crossing of $(1+z)q(d-h)$ and $(1-z)q(d+h)$ satisfies
\[
 z=-\tanh A_h(d),\qquad
 A_h(d)=\frac{V(d+h)-V(d-h)}2.
\]
At an interior crossing, differentiation of the two integrals gives
\begin{equation}\label{R4-eq:general-crossing}
 g_u''(z)=
 \frac{\sqrt{q(d-h)q(d+h)}}{A_h'(d)(1-z^2)^{3/2}}.
\end{equation}
For $d\ge0$, convexity and evenness give
$V(d+h)+V(d-h)\ge2V(h)$.
Also $A_h'(0)=V'(h)$ and
$A_h''(d)=(V''(d+h)-V''(d-h))/2\ge0$,
because $V''$ is even and nondecreasing in absolute value.
The numerator in \eqref{R4-eq:general-crossing} is therefore at most
$q(h)$ and its first denominator is at least $V'(h)$. Symmetry
handles $d<0$. If the crossing leaves the allowed midpoint interval,
$g_u''=0$; the first derivative remains continuous at that transition.
Integrating twice yields
\[
 0\le g_u(z)\le\frac{q(h)}{V'(h)}
                       (1-\sqrt{1-z^2}).
\]
Since $M(1-\sqrt{1-z^2})=(\sqrt F-\sqrt G)^2$, the directional
Sobolev translation bound gives
\[
 \int M g_u(z)\,dx
 \le\frac{q(u/2)}{V'(u/2)}\frac{u^2}{4}\mathcal J_w(f).
\]
The height integral contributes a further factor $1/2$ for $u>0$.
Changing variables $u=2h$ gives precisely $K(q)$.

At $z=0$, equality holds in the second derivative in
\eqref{R4-eq:general-crossing}. The directional difference quotient
of $\sqrt f$ converges in $L^2$ to its weak derivative. Its squares
are uniformly integrable; on $\{f=0\}$ that weak derivative vanishes
almost everywhere. The ratio of $g_u(z)$ to
$1-\sqrt{1-z^2}$ therefore converges to $q(u/2)/V'(u/2)$ in the
weighted difference-quotient integral. The preceding bound, integrated
against $u^2/8$, is integrable by \eqref{R4-eq:height-functional}.
Dominated convergence proves \eqref{R4-eq:height-variation}.

Finally $W(s)=V(\sqrt s)$ is convex: the inequality
$tV''(t)\ge V'(t)$ follows by integration of the nondecreasing
function $V''$. The logarithm of the ratio of $q$ to the centered
Gaussian density of scale $\sqrt{\pi/2}\,m$ is concave as a function
of $t^2$. The two densities have the same mass and first absolute
moment. Their difference thus has the sign pattern $-,+,-$ on the
positive half-line, allowing empty intervals and equality. A single
crossing would contradict equality of the first moments. Integrating
the difference twice, with zero mass and first moment, gives
$\E(|T|-s)_+\le\E(\sqrt{\pi/2}\,m|G_1|-s)_+$ for all $s\ge0$.
Symmetry and the stop-loss characterization prove convex order.
\end{proof}

The class includes truncated Gaussians and the parabolic auxiliary.
For a truncated Gaussian of variance parameter $v$,
$K(q)=vm/2$. For $q(t)=3(a^2-t^2)/(4a^3)$ on $(-a,a)$,
$K(q)=a^3/16$. The truncated Gaussian will give the explicit signing family below.

\subsubsection{A sharp lower bound and an exact defect identity}
The next result uses only monotonicity and Gaussian convex order; its
assumptions are broader than those of Theorem~\ref{R4-thm:general-height}.
Suppose $q$ is even, decreasing and continuously differentiable on
$(0,a)$, continuous at zero, with $-q'(t)>0$ almost everywhere and
finite $K(q)$. Write $q_a=q(a-)$ when $a<\infty$, and $q_a=0$ when
$a=\infty$. Assume the fundamental theorem of calculus gives
$\int_0^a(-q')=q(0)-q_a$.

\begin{theorem}[Gaussian optimality for the quadratic coefficient]
\label{R4-thm:height-optimality}
If $T\sim q$ and $T\cx\sqrt{\pi/2}\,mG_1$, where $m=\E|T|>0$,
then
\begin{equation}\label{R4-eq:height-barrier}
 K(q)\ge\frac{\pi m^3}{4}.
\end{equation}
Equality holds exactly for the centered Gaussian of variance
$\pi m^2/2$, on the full real line. More precisely, with
$d=q(0)-q_a$ and $c=m/(2d)$,
\begin{equation}\label{R4-eq:height-defect}
 K(q)=\frac{m^2}{4d}
       +\int_0^a\frac{(tq(t)+c q'(t))^2}{-q'(t)}\,dt,
 \qquad d\le q(0)\le\frac1{\pi m}.
\end{equation}
When $a<\infty$, the explicit stronger bound is
\begin{equation}\label{R4-eq:compact-height-barrier}
 K(q)\ge\frac{\pi m^3}{4}
       \left(1+\exp\left\{-\frac{2a^2}{\pi m^2}\right\}\right).
\end{equation}
\end{theorem}
\begin{proof}
For a symmetric continuous density at zero, its absolute stop-loss has
expansion
\[
 \E(|T|-s)_+=m-s+q(0)s^2+o(s^2).
\]
The Gaussian in the statement has the same $m$ and density
$1/(\pi m)$ at zero. Convex order therefore gives
$q(0)\le1/(\pi m)$.

Put $r(t)=-q'(t)$ and $d=\int_0^a r(t)\,dt$.
Because $\int_0^a t q(t)\,dt=m/2$, expansion of the square gives
\[
 \int_0^a\frac{(tq-cr)^2}{r}
 =K(q)-cm+c^2d=K(q)-\frac{m^2}{4d}.
\]
This proves \eqref{R4-eq:height-defect} and
\eqref{R4-eq:height-barrier}. Equality requires
$d=q(0)=1/(\pi m)$, $q_a=0$ and $-q'=tq/c$ almost everywhere.
The resulting Gaussian has variance $c=\pi m^2/2$. It vanishes at
the endpoint only when $a=\infty$, and it indeed attains equality.

For the finite-support strengthening, denote the square integral in
\eqref{R4-eq:height-defect} by $D$. Integration of the first-order
linear equation with integrating factor $e^{t^2/(2c)}$ gives
\[
 \left|q_a-q(0)e^{-a^2/(2c)}\right|
 \le\frac{\sqrt{dD}}c.
\]
To obtain this bound, multiply $tq-cr$ by
$e^{-(a^2-t^2)/(2c)}/c$ and use Cauchy--Schwarz with measure
$r(t)\,dt$; its exponential factor is at most one.
Let $x=\pi m d\le1$, $u=q_a/d$, and $A=a^2/(\pi m^2)$.
The constraint $q(0)=d(1+u)\le1/(\pi m)$ gives $u\le x^{-1}-1$.
Consequently
\[
 K(q)\ge\frac{\pi m^3}{4x}
 \left[1+\left(1-\frac{1-e^{-Ax}}x\right)_+^2\right]
 \ge\frac{\pi m^3}{4}
 \left[\frac1x+\frac{(x-1+e^{-A})_+^2}{x^3}\right].
\]
For fixed $e\in(0,1)$, the function
$x^{-1}+(x-1+e)_+^2x^{-3}$ decreases on $(0,1]$.
Where its second term is positive, the negative of its derivative,
multiplied by $x^4$, is
$x^2-2x(x-1+e)+3(x-1+e)^2>0$.
Its minimum is thus $1+e^2$, at $x=1$. Taking $e=e^{-A}$ proves
\eqref{R4-eq:compact-height-barrier}.
\end{proof}

\subsubsection{A bounded family near the Gaussian endpoint}
For $v,a>0$ put
\begin{equation}\label{R3-eq:truncated}
 q_{v,a}(t)=Z_{v,a}^{-1}e^{-t^2/(2v)}\mathbf1_{|t|<a},\quad
 Z_{v,a}=\int_{-a}^a e^{-t^2/(2v)}dt,\quad m_{v,a}=\E|T|.
\end{equation}
Its potential has constant second derivative, and
$K(q_{v,a})=vm_{v,a}/2$. Thus
\begin{equation}\label{R3-eq:Gaussian-height}
 H_f(w;q_{v,a})\ge m_{v,a}
       -\frac{vm_{v,a}}2\mathcal J_w(f).
\end{equation}
The coefficient is sharp by Theorem~\ref{R4-thm:general-height}.
Set $\kappa=\sqrt{\pi/2}$. For $0<\varepsilon<1/2$, let $v_\varepsilon$
be the unique positive solution of
\begin{equation}\label{R3-eq:v-root}
 m_{v_\varepsilon,3}=1+\varepsilon,
 \qquad
 m_{v,3}=\frac{2v(1-e^{-9/(2v)})}{\int_{-3}^3e^{-t^2/(2v)}dt},
\end{equation}
and put
\begin{equation}\label{R3-eq:new-radius}
 q_\varepsilon=q_{v_\varepsilon,3},\qquad
 R_\varepsilon=\pi\sqrt{\frac{v_\varepsilon(1+\varepsilon)}{2\varepsilon}}.
\end{equation}
The mean in \eqref{R3-eq:v-root} increases strictly from $0$ to $3/2$:
its derivative is a positive multiple of
$\Cov_{q_{v,3}}(|T|,T^2)$. Thus the parameter is unambiguous.

\begin{theorem}[Hard discrepancy at a near-endpoint Gaussian scale]
\label{R3-thm:hard-Gaussian}
Let $0<\varepsilon<1/2$, let $A\in\R^{m\times n}$ have columns $a_j$,
and let $r\ge R_\varepsilon\max_j\|a_j\|_2$, $r>0$.
There is one symmetric signing law $\mu$ satisfying
\begin{gather}
 \|A\sigma\|_\infty<r,\qquad
 \E[X_r\mid\sigma]=A\sigma,\qquad\E[T\mid\sigma]=\sigma,
 \label{R3-eq:hard-means}\\
 (X_r,T)\sim\law(X_r)\otimes q_\varepsilon^{\otimes n}.
 \label{R3-eq:hard-source}
\end{gather}
Write $\nu_\varepsilon=\E_{q_\varepsilon}T^2$,
$v_r=r^2(1/3-2/\pi^2)$ and
$b_r=\kappa r(1/2-2/\pi^2)$. The same law has
\begin{align}
 \Cov(\sigma)&\preceq
  \left(\nu_\varepsilon^{-1}I_n+v_r^{-1}A^{\mathsf T}A\right)^{-1},
  \label{R3-eq:newcov}\\
 \sigma&\cx N\!\left(0,
  \left(\frac{I_n}{\kappa^2(1+\varepsilon)^2}
                   +b_r^{-2}A^{\mathsf T}A\right)^{-1}\right)
 \cx(1+\varepsilon)\kappa G_n.
 \label{R3-eq:newcx}
\end{align}
With $B=(r/R_\varepsilon)^2$ and $w_j=\|a_j\|_2^2$,
\begin{equation}\label{R3-eq:newenergy}
 D(\mu\Vert U_n)+\E I_r(A\sigma)
 \le\delta_{q_\varepsilon}F_B(w),\qquad
 \delta_{q_\varepsilon}\le
       \pi\sqrt{\frac{q_\varepsilon(0)\varepsilon}{3}}.
\end{equation}
As $\varepsilon\downarrow0$,
\begin{equation}\label{R3-eq:radiusasymptotic}
 R_\varepsilon\sim\frac{\pi\sqrt{v_{\mathrm{tr}}/2}}{\sqrt\varepsilon},\qquad m_{v_{\mathrm{tr}},3}=1.
\end{equation}
\end{theorem}
\begin{proof}
Inequality~\eqref{R3-eq:Gaussian-height} gives height at least
$1+\varepsilon-v_\varepsilon(1+\varepsilon)\pi^2/(2R_\varepsilon^2)=1$.
Apply Theorem~\ref{thm:affine} to the full shear,
then use Corollary~\ref{REP-cor:amplification} and
Theorem~\ref{PACK-thm:exact} exactly as in the proof of
Theorem~\ref{ENERGY-thm:main}.

For completeness, an even density proportional to $e^{-W(t^2)}$, with
$W$ convex on its interval of support, obeys
$T\cx\kappa\E|T|G_1$. To check this, compare the densities of $|T|$
and $\kappa\E|T||G_1|$. Their log ratio is concave as a function of
$t^2$, so the density difference has the sign pattern $-,+,-$ or a
limiting version of it. Equal mass and equal first moment then imply
the convex order, by the secant-line test. Evenizing a convex function
on $\R$ completes the scalar comparison. This applies to both the
truncated Gaussian and the cosine density. Tensorizing, and applying
the precision left inverse to the target pair $(A\sigma,\sigma)$,
proves \eqref{R3-eq:newcx}. Conditional variance gives
\eqref{R3-eq:newcov}.

The density of $|T|$ is bounded by $2q_\varepsilon(0)$. Hence
\[
 \delta_{q_\varepsilon}
 \le\inf_{t>0}\left\{t\varepsilon+
       \frac{\pi^2q_\varepsilon(0)}{12t}\right\}
 =\pi\sqrt{q_\varepsilon(0)\varepsilon/3},
\]
using $\int_0^\infty\log(1+e^{-u})du=\pi^2/12$.
Continuity and strict monotonicity of \eqref{R3-eq:v-root} give the
asymptotic. 
\end{proof}

For cosine physical densities, the quadratic criterion has radius
$R=\pi\sqrt{K(q)/\varepsilon}$ when $\E|T|=1+\varepsilon$.
The exact defect identity therefore separates the cost of the density
shape from its endpoint trace. In the Gaussian-dominated class with
support in $(-3,3)$, it implies
\begin{equation}\label{R7-eq:quadratic-barrier}
 \sqrt\varepsilon R\ge
 \frac{\pi^{3/2}(1+\varepsilon)^{3/2}}2
 \sqrt{1+e^{-18/(\pi(1+\varepsilon)^2)}}.
\end{equation}
This is a lower bound for the universal quadratic certificate; the
exact likelihood criterion of Section~\ref{sec:reference} retains
higher-order information. The following perturbation strictly improves
the truncated-Gaussian coefficient at fixed mean and Gaussian scale.
\subsubsection{A boundary perturbation improves the truncated Gaussian}
The coefficient problem also identifies a perturbation which lowers the explicit
radius. Fix $a>0$ and $0<m<a/2$. For sufficiently small $\beta\ge0$,
define
\begin{equation}\label{R4-eq:boundary-family}
 q_\beta(t)=Z_\beta^{-1}e^{-\alpha_\beta t^2}
       (1-t^2/a^2)^\beta\mathbf1_{\{|t|<a\}},
 \qquad \E_{q_\beta}|T|=m.
\end{equation}
The parameter $\alpha_\beta>0$ is uniquely determined near $\beta=0$.

\begin{proposition}[Strict improvement by an endpoint perturbation]
\label{R4-prop:boundary-improvement}
Every sufficiently small positive $\beta$ gives a density in
Theorem~\ref{R4-thm:general-height}, with the same Gaussian comparison
scale $\sqrt{\pi/2}\,m$, and
\begin{equation}\label{R4-eq:boundary-improvement}
 K(q_\beta)=K(q_0)
 -\frac{q_0(a-)}{4\alpha_0^2}\,\beta\log(1/\beta)
 +O(\beta).
\end{equation}
Consequently the truncated Gaussian is strictly suboptimal in this
compactly supported class at fixed first absolute moment and fixed
Gaussian comparison scale.
\end{proposition}
\begin{proof}
Differentiating the first moment with respect to $\alpha$ gives
$-\operatorname{Cov}(|T|,T^2)<0$. Differentiation with respect to
$\beta$ is valid because $\log(1-t^2/a^2)$ is integrable at the
endpoints. The implicit function theorem gives
$\alpha_\beta=\alpha_0+O(\beta)$ and
$\|q_\beta-q_0\|_1=O(\beta)$. The potential
$\alpha_\beta t^2-\beta\log(1-t^2/a^2)$ has nonnegative
second derivative increasing on the positive half-interval, so the
height and comparison theorem applies.

Using
$-q_\beta'(t)/q_\beta(t)=2t(\alpha_\beta+\beta/(a^2-t^2))$
in the definition of $K$, replace $\alpha_\beta,q_\beta$ by
$\alpha_0,q_0$ with an $O(\beta)$ error. The reciprocal denominator
and its $\alpha$ derivative are uniformly bounded. With
$h=a^2-t^2$, the remaining loss is
\[
 K(q_0)-K(q_\beta)
 =\frac{\beta}{4\alpha_0}\int_0^{a^2}
       \frac{q_0(\sqrt{a^2-h})}{\alpha_0h+\beta}\,dh+O(\beta).
\]
Since $q_0(\sqrt{a^2-h})=q_0(a-)+O(h)$, its constant term gives
$\beta q_0(a-)(4\alpha_0^2)^{-1}\log(1/\beta)+O(\beta)$,
and the remainder is $O(\beta)$. This proves the expansion.
\end{proof}


\subsection{One signing for simultaneous Lebesgue-norm bounds}
\label{R19-sec:Lebesgue}
The hard support bound and the covariance bound have different uses here. The covariance bound selects a single signing with small total energy. Since that signing remains in the hard support, interpolation gives all larger norms for the same outcome.

The named dimension orders below follow from constant Koml\'os bounds and interpolation in the indicated regimes. The additional conclusion is the exact precision energy, its rank-sensitive formula, and simultaneous control of every target norm by one selected signing.

\begin{theorem}[Energy-sensitive simultaneous norm bounds]
\label{R19-thm:Lebesgue}
Set $c=6.84$, $a^2=6.5$ and $b^2=3.23$. If the columns of $A\in\R^{m\times n}$ have Euclidean norm at most one, put
\[
 M=(b^{-2}I_n+a^{-2}A^{\mathsf T}A)^{-1},\qquad
 \mathcal E(A)=\operatorname{tr}(AMA^{\mathsf T}).
\]
There is one $s\in\{-1,1\}^n$ such that $\|As\|_\infty<c$ and, simultaneously for $2\le q<\infty$,
\begin{equation}\label{R19-eq:energy-norms}
 \|As\|_q\le c^{1-2/q}\mathcal E(A)^{1/q}.
\end{equation}
If $r=\operatorname{rank}A>0$ and $F=\|A\|_F^2$, then
\begin{equation}\label{R19-eq:energy-rank}
 \mathcal E(A)\le\frac{a^2b^2Fr}{a^2r+b^2F}
 \le\frac{a^2b^2}{a^2+b^2}\,n<2.158n.
\end{equation}
If instead the columns have $\ell_p$ norm at most one for some $1\le p\le2$, there is one signing with $\|As\|_\infty<c$ and
\begin{equation}\label{R19-eq:small-p}
 \|As\|_q\le c^{1-p/q}b^{p/q}n^{1/q}
 \qquad(p\le q<\infty).
\end{equation}
\end{theorem}
\begin{proof}
Use the law in \eqref{VB-eq:main-precision}. The convex function $x\mapsto\|Ax\|_2^2$ gives $\E\|A\sigma\|_2^2\le\mathcal E(A)$. Choose one outcome $s$ with energy at most this mean. Its hard bound and
$\sum_j|(As)_j|^q\le\|As\|_\infty^{q-2}\|As\|_2^2$
prove \eqref{R19-eq:energy-norms} for every $q\ge2$ without choosing a new outcome.

If $s_1,\ldots,s_r$ are the positive singular values of $A$, then
\[
 \mathcal E(A)=\sum_{j=1}^r\frac{a^2b^2s_j^2}{a^2+b^2s_j^2}.
\]
Concavity in $s_j^2$ proves the first bound in \eqref{R19-eq:energy-rank}. The fraction is increasing in both $F$ and $r$, and $F,r\le n$, proving the second. The zero matrix has energy zero and satisfies all conclusions directly.

For $p\le2$, the columns also have Euclidean norm at most one. If $u_j$ is row $j$ of $A$, the covariance comparison and Lyapunov's inequality give
\[
 \E|u_j\sigma|^p\le b^p\|u_j\|_2^p
 \le b^p\sum_i|a_{ji}|^p.
\]
Summing over $j$ gives $\E\|A\sigma\|_p^p\le b^pn$. Choose one outcome at most this mean and interpolate with the hard bound, obtaining \eqref{R19-eq:small-p}.
\end{proof}

\begin{corollary}[Reis--Rothvoss: square matrices and the rectangular $p=2$ case]
\label{R19-cor:RR}
Let $2\le p\le\infty$ and $A\in\R^{n\times n}$ have columns of $\ell_p$ norm at most one. Set $\eta=6.5\cdot3.23/(6.5+3.23)$. One signing satisfies, simultaneously for all $2\le q\le\infty$,
\begin{equation}\label{R19-eq:RR-square}
 \|As\|_q\le 6.84^{1-2/q}\eta^{1/q}
           n^{1/2-1/p+1/q},
\end{equation}
with the usual limiting interpretation at $q=\infty$. In particular this proves Conjecture~1 of Reis--Rothvoss~\cite[Section~6]{R19-RR} for every square matrix. For $p=2$, the same conjecture holds for every rectangular matrix with $n\le m$.
\end{corollary}
\begin{proof}
In the square case, every Euclidean column norm is at most $n^{1/2-1/p}$. Apply Theorem~\ref{R19-thm:Lebesgue} to the matrix divided by this factor, and then use $\mathcal E\le\eta n$. This proves \eqref{R19-eq:RR-square}. For $q\ge p$, the conjectured factor is $\sqrt{\min(p,\log(2m/n))}$. When $m=n$, it is $\sqrt{\log2}$; the coefficient in \eqref{R19-eq:RR-square} is at most $6.84$, so a universal conjecture constant $6.84/\sqrt{\log2}<8.3$ suffices. In the rectangular $p=2$ case, apply \eqref{R19-eq:energy-norms} directly and use $\min(2,\log(2m/n))\ge\log2$.
\end{proof}
The square conclusion covers all $p\ge2$. For rectangular $p>2$, the row dimension enters the Euclidean normalization; the corollary's unrestricted rectangular assertion is precisely the $p=2$ case. The matrix-valued partial-coloring consequence is proved separately in Theorem~\ref{R19-thm:matrix-partial}.


\Needspace{12\baselineskip}
\part{Exact marginals, constraints and sampling}
\label{R6-part:applications}
Can successive choices preserve complete marginals as well as hard constraints? Conditional reference completion permits adaptation; exchange geometry and dyadic pairing then enforce the discrete constraints and individual sampling laws.

\section{Composition and polynomial-bit integer rounding}
\label{R11-sec:composition}
Each realized increment is coupled to a reference with that conditional mean. Completing the unused convolution parameter makes the terminal reference law fixed even when used parameters depend on the complete history. Centered integrable L\'evy references require no second or exponential moments. The Gaussian specialization gives optimal $\pi/8$ integer rounding and an exact expected polynomial-bit rational sampler; base polyhedra require the exchange geometry of the next section.

\subsection{Completion of convolution reference families}\label{SP-sec:convolution}

Adaptive rounding chooses its next move from the current state. To retain
a prescribed reference, couple each move to a conditionally available
source and add the unused part of an additive parameter budget.
The original history is unchanged. The completed source has a fixed law
because its transform factors into the same total parameter on every
history. Gaussian covariance, Poisson intensity and stable spectral
measure are instances of this calculation.

\subsubsection{Families with moment-generating functions}
Let $(\mathcal T,+,0)$ be a measurable additive parameter space with a measurable family of centered integrable laws $\nu_\theta$ on $\R^d$, such that
\[
 \nu_0=\delta_0,\qquad \nu_{\theta+\vartheta}=\nu_\theta*\nu_\vartheta.
\]
Only parameters and residual parameters that occur below need be included in $\mathcal T$. Assume the family is defined as a measurable probability kernel. Let
$M_\theta(t)=\int e^{\langle t,x\rangle}\nu_\theta(\dd x)$.

\begin{theorem}[Completion with a deterministic convolution parameter]\label{SP-thm:convolution-mgf}
Let $(\mathcal F_i)_{i=0}^r$ be an original history in standard Borel spaces, with integrable martingale differences $\Delta_i$. Suppose $\theta_i$ is $\mathcal F_{i-1}$-measurable and
\[
 \Law(\Delta_i\mid\mathcal F_{i-1})\cx\nu_{\theta_i}.
\]
Suppose a measurable residual parameter $\theta_0^*$ satisfies
\[
 \theta_0^*+\sum_{i=1}^r\theta_i=\Theta
\]
for a deterministic $\Theta$, and $M_\Theta$ is finite on a neighbourhood of zero. Then an extension preserving the entire original history carries $Z_\Theta$ such that
\begin{equation}\label{SP-eq:convolution-output}
 \Law(Z_\Theta\mid\mathcal F_0)=\nu_\Theta,
 \qquad \E[Z_\Theta\mid\mathcal F_r]=\sum_{i=1}^r\Delta_i.
\end{equation}
In particular $Z_\Theta$ is independent of $\mathcal F_0$, and the sum is convex-order dominated by $\nu_\Theta$, also conditionally on $\mathcal F_0$.
\end{theorem}
\begin{proof}
Choose measurable conditional martingale couplings using \cite{LeskelaVihola}. Given the full original history, draw each reference $Z_i$ from its conditional kernel at $(\mathcal F_{i-1},\Delta_i)$ using fresh mutually independent randomization. Then
\[
 \E[Z_i\mid\mathcal F_r]=\Delta_i,
 \quad
 \Law(Z_i\mid\mathcal F_{i-1},Z_1,\ldots,Z_{i-1})=\nu_{\theta_i}.
\]
The second identity follows because the added past randomization carries no information about the next original increment beyond the original past. It does not require that $\mathcal F_i$ be generated by the martingale sum alone.

Put $S=\sum_iZ_i$ and conditionally draw $Z_0^*\sim\nu_{\theta_0^*}$, independently of the added references given the whole history. Its conditional mean is zero. Thus $Z_\Theta=S+Z_0^*$ has the conditional mean in \eqref{SP-eq:convolution-output}.

Fix $t$ in a neighbourhood on which $M_\Theta(t)<\infty$. Centering gives $M_\theta(t)\ge1$. Factorization with the residual therefore makes every encountered moment-generating function finite. In the augmented filtration,
\[
 L_k(t)=\exp\left(\left\langle t,\sum_{i\le k}Z_i\right\rangle
                         -\sum_{i\le k}\log M_{\theta_i}(t)\right)
\]
is a positive martingale of conditional mean one given $\mathcal F_0$. Finite-step conditional integration proves its integrability inductively. Consequently
\begin{align*}
 \E[e^{\langle t,Z_\Theta\rangle}\mid\mathcal F_0]
 &=\E[e^{\langle t,S\rangle}M_{\theta_0^*}(t)\mid\mathcal F_0]\\
 &=M_\Theta(t)\E[L_r(t)\mid\mathcal F_0]=M_\Theta(t).
\end{align*}
Uniqueness of moment-generating functions proves the fixed conditional law. One may first use a countable dense set of $t$ and then continuity, so that all identities hold outside a single null set. Every original joint law was preserved by construction.
\end{proof}

\begin{example}[Integer parameters and matrix references]\label{SP-ex:convolution}
For a fixed centered law $\rho$ with an exponential moment, take $\nu_k=\rho^{*k}$, $k\in\mathbb N_0$. Predictably allocated integer budgets $k_i$ with $\sum k_i\le N$ complete to exactly $\rho^{*N}$. The family can be multitype, with a separate integer budget for each source distribution.

For a matrix example, in the Euclidean space of symmetric $d\times d$ matrices let
\[
 \nu_k=\Law\left(\sum_{j=1}^k(g_jg_j^{\mathsf T}-I_d)\right),
 \qquad g_j\sim N(0,I_d)\text{ independently}.
\]
These are centered Wishart references with integer degrees of freedom. Their moment-generating functions are finite near zero, so predictable allocations with total degrees at most $N$ complete to the exact centered $N$-degree reference in full matrix convex order. No assertion that matrix Wishart laws are infinitely divisible is needed for this integer-parameter example.
\end{example}

\subsubsection{Centered infinitely divisible laws with only a first moment}
A parameter is now $\theta=(Q,\nu)$, where $Q\succeq0$ and $\nu$ is a L\'evy measure satisfying
\begin{equation}\label{SP-eq:levy-integrability}
 \int_{\|x\|\le1}\|x\|^2\nu(\dd x)
 +\int_{\|x\|>1}\|x\|\nu(\dd x)<\infty.
\end{equation}
Let $\mathcal L_\theta$ be the centered infinitely divisible law with characteristic exponent
\begin{equation}\label{SP-eq:levy-exponent}
 \psi_\theta(t)=-\frac12t^{\mathsf T}Qt
   +\int(e^{i\langle t,x\rangle}-1-i\langle t,x\rangle)\nu(\dd x).
\end{equation}
The convention compensates all jumps. Condition \eqref{SP-eq:levy-integrability} gives an integrable centered law. Addition of parameters corresponds to convolution. Parameter domination means Loewner domination of $Q$ and measure domination of $\nu$.

\begin{theorem}[L\'evy reference completion]\label{SP-thm:levy}
In the history setting above, suppose
\[
 \Law(\Delta_i\mid\mathcal F_{i-1})\cx\mathcal L_{\theta_i},
 \qquad \sum_i\theta_i\preceq\Theta
\]
pathwise for a deterministic parameter satisfying \eqref{SP-eq:levy-integrability}. There is an extension of the original history carrying $Z_\Theta\sim\mathcal L_\Theta$, independent of $\mathcal F_0$, with
\[
 \E[Z_\Theta\mid\mathcal F_r]=\sum_i\Delta_i.
\]
No exponential or second moment of the reference is required.
\end{theorem}
\begin{proof}
Use the same conditional couplings and the residual parameter
$\Theta-\sum_i\theta_i$. For fixed $t$, the process
\[
 L_k(t)=\exp\left(i\left\langle t,\sum_{i\le k}Z_i\right\rangle
                         -\sum_{i\le k}\psi_{\theta_i}(t)\right)
\]
is a complex martingale. Its modulus is bounded by
$\exp[-\Re\psi_\Theta(t)]$, because $-\Re\psi$ is a nonnegative additive functional on the parameter cone. Conditional integration of the residual gives
\[
 \E[e^{i\langle t,Z_\Theta\rangle}\mid\mathcal F_0]
 =e^{\psi_\Theta(t)}\E[L_r(t)\mid\mathcal F_0]
 =e^{\psi_\Theta(t)}.
\]
Characteristic-function uniqueness identifies the conditional law.

The variables used in the conditional-mean calculation are integrable. Here is a uniform bound which also justifies integration over random parameters. Split any $\theta\preceq\Theta$ at jump size one. Its Gaussian plus compensated small-jump part has first absolute moment bounded by the square root of the trace of its covariance, which is at most the corresponding bound for $\Theta$. The large-jump compound Poisson part, after centering, has first absolute moment at most $2\int_{\|x\|>1}\|x\|\nu_\Theta(\dd x)$. This bounds $\E\|Z_\theta\|$ uniformly over the parameter interval. The history-preserving mean identities now follow exactly as before.
\end{proof}

\begin{corollary}[Stable reference budgets]\label{SP-cor:stable}
Fix $1<\alpha<2$. Let $S_{\alpha,q}$ be the symmetric isotropic stable law with characteristic function $e^{-q\|t\|^\alpha}$. If
\[
 \Law(\Delta_i\mid\mathcal F_{i-1})\cx S_{\alpha,q_i},
 \qquad q_i\ge0\text{ predictable},\qquad\sum_iq_i\le q,
\]
then the sum is dominated by $S_{\alpha,q}$ through a coupling retaining the complete history. The same statement holds for a fixed symmetric stable spectral measure, with predictable scalar multiples; more general predictable spectral measures can be used if their sum is dominated as a measure by a deterministic one.
\end{corollary}
\begin{proof}
Stable laws in this range are integrable and their L\'evy measures satisfy \eqref{SP-eq:levy-integrability}. Apply Theorem~\ref{SP-thm:levy}.
\end{proof}

\subsubsection{Exact Poisson domination under a predictable mean budget}
\begin{lemma}[Bernoulli to Poisson]\label{SP-lem:bernoulli-poisson}
For $0\le p\le1$,
\[
 \operatorname{Bern}(p)-p\cx\operatorname{Pois}(p)-p.
\]
Among centered unit-jump Poisson references $\operatorname{Pois}(\lambda)-\lambda$, the parameter $\lambda=p$ is minimal for $0<p<1$.
\end{lemma}
\begin{proof}
Let $N\sim\operatorname{Pois}(p)$. Set $B=1$ whenever $N>0$, and randomize on $\{N=0\}$ so that $\Pp(B=1)=p$. This is possible since $1-e^{-p}\le p$. Then $\E[N\mid B=0]=0$ and $\E[N\mid B=1]=p/p=1$, proving the martingale coupling. If $\lambda<p$, the convex function $(-x-\lambda)_+$ has positive expectation at $\operatorname{Bern}(p)-p$ and zero expectation at the proposed Poisson reference.
\end{proof}

\begin{remark}[The deterministic endpoint of the scalar lemma]
For $p=1$ the centered Bernoulli is zero and the minimal Poisson parameter is zero. The domination in Lemma~\ref{SP-lem:bernoulli-poisson} remains valid at both endpoints. This endpoint is separated explicitly because its zero-probability lower atom cannot be used as a support obstruction.
\end{remark}

\begin{theorem}[Adaptive bounded variables and a Poisson reference]\label{SP-thm:poisson}
Let $0\le X_i\le c$ be adapted, with $c>0$, and put
$m_i=\E[X_i\mid\mathcal F_{i-1}]$. If $\sum_i m_i\le L$ pathwise, where $L$ is deterministic, then
\begin{equation}\label{SP-eq:poisson-bound}
 \sum_i(X_i-m_i)\cx c\{\operatorname{Pois}(L/c)-L/c\}.
\end{equation}
The reference can be coupled with the whole original history so that its conditional mean is the displayed sum. In the class determined only by $c,L$ and an unrestricted number of increments, this reference is sharp as an upper law: it is a limit of admissible sums.
\end{theorem}
\begin{proof}
Conditional chord interpolation bounds $X_i/c-m_i/c$ in convex order by $\operatorname{Bern}(m_i/c)-m_i/c$. Lemma~\ref{SP-lem:bernoulli-poisson} then gives the centered Poisson reference of rate $m_i/c$ and jump size $c$. Its predictable rates sum to at most $L/c$, so either completion theorem applies. For sharpness take independent $X_i=cB_i$, $1\le i\le n$, with $B_i\sim\operatorname{Bern}(L/(cn))$ and $n\ge L/c$. Their centered sum converges in law and first moment to the reference. Any fixed integrable convex-order upper law for the whole class must therefore dominate this limit.
\end{proof}

The conditional comparison in \eqref{SP-eq:poisson-bound} controls every convex function, including stop-loss losses. Its deterministic budget bounds the sum of conditional means; predictable quadratic variation is a different constraint. The completion couples the reference to the entire original history and retains the conditional-mean identity with that history. General semimartingale comparisons are developed in~\cite{KopferRuschendorf}.


\subsection{Adaptive Gaussian completion and an exact integer sampler}
\label{sec:exact-gaussian}

\subsubsection{Composition under a predictable covariance bound}
Let $(\mathcal F_i)_{i=0}^r$ be the filtration generated by a history in
standard Borel spaces, and let $\Delta_i\in\R^d$ be martingale
differences. Conditional convex order means convex order between the
regular conditional laws, almost surely. The following statement uses a
bound on the covariances of the Gaussian envelopes.

\begin{theorem}[Gaussian comparison for adaptive increments]
\label{thm:gaussian-composition}
Suppose $Q_i\succeq0$ is $\mathcal F_{i-1}$-measurable and
\begin{equation}\label{eq:conditional-gaussian-increments}
 \law(\Delta_i\mid\mathcal F_{i-1})\cx N(0,Q_i)
 \quad\hbox{almost surely},\qquad
 \sum_{i=1}^r Q_i\preceq Q
\end{equation}
for a deterministic $Q\succeq0$. There is an extension of the entire
original history carrying $G_Q\sim N(0,Q)$ such that
\begin{equation}\label{eq:gaussian-composition-coupling}
 G_Q\perp\mathcal F_0,\qquad
 \E[G_Q\mid\mathcal F_r]=\sum_{i=1}^r\Delta_i.
\end{equation}
In particular $\sum_i\Delta_i\cx N(0,Q)$, conditionally on
$\mathcal F_0$ as well as unconditionally. Every joint law of the
original history is preserved by the extension.
\end{theorem}
The Gaussian envelope used at a step depends on the earlier history, so its accumulated covariance is random. The pathwise bound leaves a positive semidefinite residual covariance. Drawing that residual after the history is complete makes the total convolution parameter exactly $Q$ on every path. The conditional characteristic-function calculation below gives the fixed Gaussian law while its conditional mean retains the entire original sum.
\begin{proof}
Apply Theorem~\ref{SP-thm:levy} to the parameter pair $(Q_i,0)$ and
deterministic total $(Q,0)$. The residual source is
$N(0,Q-\sum_iQ_i)$, drawn conditionally on the complete history.
The compensated characteristic function in that theorem becomes
\[
 \exp\!\left(i\left\langle t,\sum_{j\le k}Z_j\right\rangle
       +\frac12\sum_{j\le k}t^{\mathsf T}Q_jt\right).
\]
Its modulus is bounded by $\exp(t^{\mathsf T}Qt/2)$. Integrating the
residual cancels the random accumulated covariance and leaves
$\exp(-t^{\mathsf T}Qt/2)$ on every initial history. The same construction
retains all original joint laws and gives
\eqref{eq:gaussian-composition-coupling}.
\end{proof}

The matrices $Q_i$ may depend on earlier increments and need not commute.
The conclusion is a comparison of the entire vector law. A deterministic
upper bound on the accumulated Gaussian envelope is sufficient; no
independence of the original increments is required.

\subsubsection{Exact integer rounding}
\begin{theorem}[Exact Gaussian integer rounding]\label{thm:exact-gaussian-integer}\label{MAIN-exact}
Let $A=[a_1\ \cdots\ a_n]\in\R^{m\times n}$ and $b\in\R^n$.
Set $J=\{j:b_j\notin\Z\}$, and order the columns of $A_J$ arbitrarily.
Let $g_1,\ldots,g_r$, $r=|J|$, be their unnormalized Gram--Schmidt
residuals, including zero residuals. There is a law on at most $2^r$
integer vectors $N$ such that
\begin{equation}\label{eq:exact-gaussian-main}
 \E N=b,\qquad N_j=b_j\ (j\notin J),\qquad
 A(N-b)\cx N\!\left(0,\frac\pi8\sum_{i=1}^r g_i g_i^{\mathsf T}\right).
\end{equation}
In particular, if $\|a_j\|_2\le1$, then
\begin{equation}\label{eq:exact-gaussian-projection}
 A(N-b)\cx N\!\left(0,\frac\pi8 P_{\operatorname{range}(A_J)}\right)
 \cx N(0,\tfrac\pi8 I_m).
\end{equation}
The constant $\pi/8$ is optimal. For rational $A,b$, the law has an exact
sampler with at most $r$ Bernoulli choices and expected polynomial bit
complexity in the input length.
\end{theorem}

We first record the exact one-dimensional comparison, including its
uniform constant. Write $\varphi$ and $\Phi$ for the standard normal density
and distribution function.

\begin{lemma}[Bernoulli--Gaussian comparison]\label{lem:bernoulli-gaussian-exact}
For $0<p<1$, put
\[
 t(p)=\frac{p(1-p)}{\varphi(\Phi^{-1}(p))}.
\]
Then $\operatorname{Bern}(p)-p\cx t(p)G$, where $G$ is standard normal.
This scale is the smallest possible, and
\begin{equation}\label{eq:bernoulli-scale-sharp}
 t(p)\le\sqrt{\pi/8},
\end{equation}
with equality precisely at $p=1/2$. Set $t(0)=t(1)=0$.
\end{lemma}
\begin{proof}
Let $c=\Phi^{-1}(1-p)$ and $B=\ind_{\{G>c\}}$. Gaussian integration gives
$\E[G\ind_{\{G>c\}}]=\varphi(c)$; consequently
$\E[t(p)G\mid B]=B-p$.
Conversely, in any martingale coupling between $B-p$ and $sG$,
\[
 p(1-p)=s\E[G\ind_{\{B=1\}}]\le s\varphi(c).
\]
The last inequality follows by placing a set of Gaussian mass $p$ on the
upper tail; more generally the same bound holds for randomized sets by
integrating their conditional probabilities. This proves optimality.

For completeness, $t(p)\le\sqrt{\pi/8}$ is equivalent to $\varphi(x)\ge4\varphi(0)\Phi(x)\Phi(-x)$. Both sides are even. On $x\ge0$ let their difference be $F(x)$. Then
$F(0)=\lim_{x\to\infty}F(x)=0$ and
\[
 F'(x)=\varphi(x)h(x),\quad
 h(x)=4\varphi(0)(2\Phi(x)-1)-x,\quad
 h'(x)=\frac4\pi e^{-x^2/2}-1.
\]
The derivative $h'$ changes sign once, while $h(0)=0$ and $h(x)\to-\infty$.
Thus $F$ first increases and then decreases to zero, and is positive on
$(0,\infty)$.
\end{proof}

\begin{corollary}[Sharp comparison for bounded martingale increments]
\label{cor:bounded-martingale-gaussian}
Let $a_i\in\R^d$ and $p_i\in[0,1]$ be predictable, and let
$\xi_i$ satisfy
\[
 \E[\xi_i\mid\mathcal F_{i-1}]=0,\qquad
 -p_i\le\xi_i\le1-p_i.
\]
If $\sum_i a_ia_i^{\mathsf T}\preceq Q$ pathwise for a deterministic
$Q\succeq0$, then
\begin{equation}\label{eq:bounded-martingale-gaussian}
 \sum_i a_i\xi_i\cx N(0,\tfrac\pi8 Q).
\end{equation}
The conditional-history coupling \eqref{eq:gaussian-composition-coupling}
also holds. More precisely, $\pi Q/8$ can be replaced by any deterministic
upper bound for $\sum_i t(p_i)^2a_ia_i^{\mathsf T}$. The universal factor
$\pi/8$ cannot be reduced, even for one scalar increment.

The same convex-order bound holds for any stopped sum whose accumulated
matrix budget is bounded by $Q$, including an unbounded stopping time
when the stopped sum is defined by its $L^2$ limit.
\end{corollary}
\begin{proof}
A centered variable in $[-p,1-p]$ is convex-order dominated by
$\operatorname{Bern}(p)-p$: the secant joining the two endpoint values
majorizes every convex function on that interval, and its expectation
has the Bernoulli endpoint weights. Lemma~\ref{lem:bernoulli-gaussian-exact}
therefore gives the conditional envelope
$N(0,t(p_i)^2a_ia_i^{\mathsf T})$. Apply
Theorem~\ref{thm:gaussian-composition} and
$t(p_i)^2\le\pi/8$. For one fair Bernoulli increment, the convex test
$|x|$ gives the matching lower bound $\pi/8$.

For a stopping time $T$, replace $a_i$ by
$\ind_{\{T\ge i\}}a_i$, which is predictable, and first stop at a
finite deterministic index. The martingale has
\[
 \E\left\|\sum_{i\le n\wedge T}a_i\xi_i\right\|_2^2
 =\sum_{i\le n}\E[\ind_{\{T\ge i\}}\|a_i\|_2^2
                      \E(\xi_i^2\mid\mathcal F_{i-1})]
 \le\tfrac14\tr Q.
\]
It therefore converges in $L^2$. Convex order is closed under first-moment
convergence against a fixed integrable upper law, giving the limit.
The sharper $t(p_i)^2$ budget instead gives the bound $\tr Q$ for the
corresponding Gaussian envelope budget and the same argument.
\end{proof}

The hypothesis concerns interval widths, or the sharper Gaussian
variances $t(p_i)^2$, rather than only the conditional second moments.
For $p\downarrow0$, the optimal Gaussian variance of
$\operatorname{Bern}(p)-p$ is asymptotic to
$1/(2\log(1/p))$, while its variance is asymptotic to $p$.
Thus replacing the envelope budget by a universal multiple of the
predictable quadratic variation would fail already at one step.

\begin{proof}[Proof of Theorem~\ref{thm:exact-gaussian-integer}]
Work on $J$ and leave all other coordinates fixed. Denote its ordered
columns again by $a_1,\ldots,a_r$. Gram--Schmidt elimination gives an
upper-triangular matrix $D=[d_1\ \cdots\ d_r]$ with diagonal entries one
such that
\begin{equation}\label{eq:triangular-orthogonalization}
 A_JD=[g_1\ \cdots\ g_r],\qquad
 \ip{g_i}{g_j}=0\quad(i\ne j).
\end{equation}
One explicit recursion, valid also for rank-deficient matrices, is
\[
 d_i=e_i-\sum_{k<i:g_k\ne0}
       \frac{\ip{a_i}{g_k}}{\|g_k\|_2^2}d_k,
 \qquad g_i=A_Jd_i.
\]
Thus $d_i$ is supported on $\{1,\ldots,i\}$ and $(d_i)_i=1$.

Initialize $y=b_J$. For $i=r,r-1,\ldots,1$, set
$p_i=y_i-\lfloor y_i\rfloor$, draw a Bernoulli variable $B_i$ of parameter
$p_i$, and replace
\begin{equation}\label{eq:nearest-plane-update}
 y\longleftarrow y+(B_i-p_i)d_i.
\end{equation}
The update makes coordinate $i$ integral and leaves all later coordinates
unchanged. Each update has conditional mean zero. The final vector is
therefore integral and has mean $b_J$. With $\epsilon_i=B_i-p_i$, it obeys
\begin{equation}\label{eq:rounding-error-factorization}
 N_J-b_J=D\epsilon,\qquad
 A(N-b)=\sum_i\epsilon_i g_i.
\end{equation}
Theorem~\ref{thm:gaussian-composition}, applied to the coordinate increments $\epsilon_i e_i$ in reverse order with deterministic envelope $\pi e_ie_i^{\mathsf T}/8$, proves the
stronger coefficient-space comparison
\begin{equation}\label{eq:coefficient-gaussian-exact}
 N_J-b_J\cx N(0,\tfrac\pi8 DD^{\mathsf T}),
\end{equation}
and hence \eqref{eq:exact-gaussian-main}. Since
$\|g_i\|_2\le\|a_i\|_2\le1$ and the nonzero $g_i$ are orthogonal,
$\sum_i g_i g_i^{\mathsf T}\preceq P_{\operatorname{range}(A_J)}$.
Adding independent Gaussian noise gives
\eqref{eq:exact-gaussian-projection}.

Every outcome satisfies the explicit finite displacement estimate
\begin{equation}\label{eq:triangular-window}
 |N_{J_k}-b_{J_k}|<\sum_{i=k}^r|D_{ki}|.
\end{equation}
For rational input, Gram--Schmidt can be performed without square roots.
Its coefficients are ratios of minors of rational Gram matrices; standard
determinant bounds give polynomial bit lengths after clearing input
denominators. The entries of $D^{-1}$ have polynomial bit lengths by the
same bounds for a unit triangular rational matrix. At step $i$, the
current coordinate can alternatively be computed from the already chosen
integers by
\[
 y_i=b_{J_i}-\sum_{j>i}(D^{-1})_{ij}(N_{J_j}-b_{J_j}).
\]
The estimate \eqref{eq:triangular-window} bounds the bit lengths of those
integers polynomially in the input length. Thus the probabilities used
in \eqref{eq:nearest-plane-update} also have polynomial bit lengths.
A Bernoulli probability $a/q$ is sampled exactly by drawing a uniform
integer in $\{0,\ldots,q-1\}$ using binary rejection sampling and comparing
it with $a$; the expected number of trials is less than two. There are
$r$ such choices and polynomially many rational arithmetic operations.

Finally, for $A=I_m$ and $b=(1/2,\ldots,1/2)$, every integer output has
$\E|N_i-1/2|\ge1/2$. Domination by $N(0,s^2I_m)$ therefore requires
$1/2\le s\sqrt{2/\pi}$, or $s^2\ge\pi/8$.
\end{proof}

\begin{corollary}[Ellipsoidal Gaussian references]\label{cor:ellipsoid-exact}
Let $Q\succeq0$ and suppose $a_j\in\operatorname{range}(Q)$ and
$a_j^{\mathsf T}Q^\dagger a_j\le1$ for every initially nonintegral
coordinate. Then unbiased integer rounding is possible with
\[
 A(N-b)\cx N(0,\tfrac\pi8 Q),
\]
while fixing every initially integral coordinate.
\end{corollary}
\begin{proof}
Apply Theorem~\ref{thm:exact-gaussian-integer} to the columns
$Q^{\dagger/2}a_j$ on $\operatorname{range}(Q)$ and map back by $Q^{1/2}$.
\end{proof}

\begin{proposition}[Rigidity at the universal endpoint]\label{prop:gaussian-endpoint-rigidity}
Suppose $N\in\Z^m$ has mean $\tfrac12\mathbf1$ and
$N-(1/2,\ldots,1/2)\cx N(0,\pi I_m/8)$.
Then the coordinates of $N$ are independent fair Bernoulli variables.
In any martingale coupling with this Gaussian, $N_i$ is almost surely the
indicator that its $i$th Gaussian coordinate is positive.
\end{proposition}
\begin{proof}
In a martingale coupling write $\Delta=N-(1/2,\ldots,1/2)$ and
$Z\sim N(0,\pi I_m/8)$, with $\E[Z\mid\Delta]=\Delta$. The inequalities $\tfrac12\le\E|\Delta_i|\le\E|Z_i|=\tfrac12$ are equalities. Thus $|\Delta_i|=1/2$ almost surely. Equality in conditional
Jensen for the absolute value implies that $Z_i$ has the sign of
$\Delta_i$ almost surely. Since Gaussian coordinates have no atom at zero,
$N_i=\ind_{\{Z_i>0\}}$. Independence follows from the Gaussian product law.
\end{proof}

\begin{remark}
Theorem~\ref{thm:exact-gaussian-integer} gives the optimal variance,
an explicit matrix-dependent displacement window, and an expected
polynomial-time integer sampler. The bound \eqref{eq:triangular-window}
may grow when columns are nearly dependent. The hard-supported signing
law of Theorem~\ref{MAIN-signing} and the nearest-cell laws require the
separate constructions given in their respective sections.
\end{remark}

\subsubsection{Gaussian comparison for an entire rounding path}
Theorem~\ref{thm:gaussian-composition} preserves the original adapted history. Lifting each increment into a path space therefore controls convex losses of the entire trajectory.

\begin{proposition}[Trajectory comparison]\label{R7-prop:trajectory}
Let $\Delta_i\in\R^d$, $1\le i\le r$, be adapted to a filtration $(\mathcal F_i)$ and suppose, conditionally on $\mathcal F_{i-1}$,
\[
 \Delta_i\cx N(0,Q_i)
\]
for deterministic positive semidefinite matrices $Q_i$. Put $M_k=\sum_{i\le k}\Delta_i$. Let $G_1,\ldots,G_r$ be independent with $G_i\sim N(0,Q_i)$. Then
\begin{equation}\label{R7-eq:path-order}
 (M_1,\ldots,M_r)\cx
 (G_1,G_1+G_2,\ldots,G_1+\cdots+G_r).
\end{equation}
The coupling can retain the full original history, as in Theorem~\ref{thm:gaussian-composition}.
\end{proposition}
\begin{proof}
Let $L_i:\R^d\to(\R^d)^r$ place its argument in blocks $i,i+1,\ldots,r$ and zero in earlier blocks. The lifted increments satisfy
\[
 L_i\Delta_i\cx N(0,L_iQ_iL_i^{\mathsf T})
 \quad\text{conditionally on }\mathcal F_{i-1}.
\]
Their sum is $(M_1,\ldots,M_r)$. Apply the composition theorem with deterministic budget
\[
 K=\sum_{i=1}^rL_iQ_iL_i^{\mathsf T}.
\]
The Gaussian on the right side of \eqref{R7-eq:path-order} has precisely this covariance.
\end{proof}

Thus any integrable convex loss of the entire trajectory is controlled by the Gaussian path. Examples include
\[
 \max_{1\le k\le r}\|M_k\|,
 \qquad \max_{k,t\in T}\langle t,M_k\rangle,
 \qquad \exp\left(\lambda\max_{1\le k\le r}\|M_k\|\right)
 \quad(\lambda\ge0).
\]
The statement avoids reducing every path loss to a union of separate terminal tests. Predictable random $Q_i$ also work if the \emph{lifted} sum $\sum_iL_iQ_iL_i^{\mathsf T}$ is bounded by a deterministic path-space covariance. A bound on $\sum_iQ_i$ alone controls the terminal sum, not the entire path-space covariance asserted here.

\subsubsection{Computational scope}\label{NEW-sec:construction}
The construction retains exact probability laws, while its algorithmic
implementations depend on how those laws and the scalar primitives are
presented. Table~\ref{R9-tab:runtime} records the guarantees used here.
The code examples in Section~\ref{SP-sec:codes} concern local heat-bath
dynamics and do not change these implementation statements.

\begin{table}[!htbp]
\caption{Constructions and computational guarantees.}
\label{R9-tab:runtime}
\centering\small
\begin{tabular}{@{}>{\raggedright\arraybackslash}p{.29\textwidth}>{\raggedright\arraybackslash}p{.65\textwidth}@{}}
\toprule
Construction & Available realization and input scope \\
\midrule
Joint hard signing law & Finite support selected by convex separation;
scalar overlap certificates are explicit. A general polynomial-time
kernel is not asserted.\\[4pt]
Exact rational integer rounding & Expected polynomial bit complexity
for rational matrix and mean, with the explicit displacement window of
Theorem~\ref{thm:exact-gaussian-integer}.\\[4pt]
Integral base-polyhedron rounding & Expected polynomial bit complexity with an integral submodular value oracle and polynomial submodular minimization. The nearest cell and tight faces are exact; the Gaussian covariance has any prescribed positive slack. See Theorem~\ref{BASE-thm:main}.\\[4pt]
Prescribed marginals and Hilbert sampling & Finite dyadic constructions
use conditional signing calls; the general law follows by fixed-marginal
compactness. Local polynomial formulas are finite.\\[4pt]
Exchange and capacity constraints & The shared signing construction
uses the exchange decompositions stated in Section~\ref{R7-sec:exchange}.
Their availability is part of the input.\\[4pt]
Poisson spectral sparsification & A hard-supported finite count law
with exact expectations; the argument proves existence and entropy
bounds, without a polynomial-time sampling claim.\\[4pt]
Canonical moment partitions & A coercive finite-dimensional convex
dual with unique multipliers and exact-moment soft approximations.
No uniform conditioning or numerical complexity bound is imposed.\\[4pt]
Backward-Monge maps & Measurable partitions and quantile maps, with
exact marginals and the stated Wasserstein approximation. General
measurable input is not a finite algorithmic encoding.\\[4pt]
Exact private encoders & Deterministic conditional-quantile correction;
the Gaussian formulas give explicit error rates. Evaluating arbitrary
conditional distributions requires its own input representation.\\
\bottomrule
\end{tabular}
\end{table}


\subsection{Bounded differences in full convex order}\label{SP-sec:bounded-differences}
The bounded-increment comparison gives a Gaussian reference for every convex loss of a function of independent variables. Its optimal full-convex-order variance is determined already by one fair bit.

\begin{theorem}[Bounded differences in full convex order]\label{SP-thm:mcdiarmid}
Let $X_1,\ldots,X_n$ be independent and let changing coordinate $i$ change a real-valued integrable function $f$ by at most $c_i$. Then
\[
 f(X)-\E f(X)\cx N\left(0,\frac\pi8\sum_i c_i^2\right).
\]
The coefficient $\pi/8$ is optimal for domination against all convex tests.
\end{theorem}
\begin{proof}
The Doob martingale obtained by revealing coordinates has centered increments in predictable intervals of length at most $c_i$. For $0<p<1$, put
\[
 t(p)=\frac{p(1-p)}{\varphi(\Phi^{-1}(p))}.
\]
If $G$ is standard normal and $B=\ind_{\{G>\Phi^{-1}(1-p)\}}$, Gaussian tail integration gives
$\E[t(p)G\mid B]=B-p$. Thus every centered variable in $[-p,1-p]$ is convex-order dominated by $t(p)G$, by chord interpolation followed by this coupling. The scalar inequality $t(p)\le\sqrt{\pi/8}$ is proved in Lemma~\ref{lem:bernoulli-gaussian-exact}; its proof uses differentiation of
$\varphi(x)-4\varphi(0)\Phi(x)\Phi(-x)$. Gaussian completion gives the displayed sum bound.

For a single fair Bernoulli and $c_1=1$, the convex test $|x|$ requires
$1/2\le s\sqrt{2/\pi}$ for any dominating $N(0,s^2)$, giving $s^2\ge\pi/8$.
\end{proof}

Pinelis's normal-domination theorem \cite{Pinelis} uses a smaller Gaussian variance for a restricted test class; Theorem~\ref{SP-thm:mcdiarmid} determines the optimal variance for all convex tests. The usual Hoeffding exponential-moment bound has its separate, smaller parameter. In vector form, Corollary~\ref{cor:bounded-martingale-gaussian} and the completion theorem control every norm and support function under a deterministic accumulated Gaussian-envelope matrix. The interval-width and reference budgets specify the hypotheses of those comparisons.


\section{Exchange constraints and Gaussian rounding on base polyhedra}
\label{R11-sec:constraints}
Swap rounding on matroid bases preserves marginals and controls concentration and negative dependence~\cite{BASE-CVZ}. A Gaussian comparison of the entire error requires feasible paired moves whose joint signing also retains the hard features.

For an arbitrary integral base polyhedron, the nearest unit cell is a translated matroid base polytope. Symmetric exchanges of complementary bases charge covariance to the minimal face. They preserve every tight face and yield the stated expected polynomial-bit value-oracle sampler with positive variance slack. Actual covariance, negative correlation and convex-order domination retain their separate conclusions.

\subsection{Exact constraints from simultaneous exchanges}
\label{R7-sec:exchange}

\subsubsection{Complementary feasibility and a common Gram bound}
A difference between two feasible points may admit many local moves. The construction requires complementary feasibility for every subset of those moves. This permits one signing to control all current pairs simultaneously. A matrix bound on their Gram sums then makes the reference at each merge independent of the current list.

Put $R_m=C-\gamma/\sqrt m$. The constants $\alpha,\tau$ are
those of \eqref{eq:reference-constants}. The next theorem makes both
the hard support and the Gaussian coefficient budget independent of
the evolving list of feasible points.

\begin{theorem}[Simultaneous-exchange rounding]\label{R7-thm:exchange}
Let $V\subset\R^n$ be finite, let $x\in\conv V$, and put $H=\operatorname{span}(V-V)$. Let $A:\R^n\to\R^m$. Suppose there are $d\ge0$ and a positive semidefinite operator $B$ supported on $H$ such that every pair $u,v\in V$ admits a decomposition
\[
 u-v=\sum_a c_a
\]
with the following properties:
\begin{equation}\label{R7-eq:all-subsets}
 u-\sum_{a\in J}c_a\in V,
 \qquad v+\sum_{a\in J}c_a\in V
 \quad\text{for every subset }J;
\end{equation}
\begin{equation}\label{R7-eq:exchange-budgets}
 \sum_a c_ac_a^{\mathsf T}\preceq B,
 \qquad \|Ac_a\|_2\le d\quad\text{for every }a.
\end{equation}
There is one random $Z\in V$ with $\E Z=x$ for which, writing $E=Z-x$,
\begin{equation}\label{R7-eq:exchange-joint}
 (AE,E)\cx(U_{d,m},G_B),\qquad
 U_{d,m}=d\sum_{k\ge1}2^{-k}X_{R_m}^{(k)},\qquad
 G_B\sim N(0,\tau B/2),\qquad U_{d,m}\perp G_B.
\end{equation}
The $X_{R_m}^{(k)}$ are independent copies. Every outcome satisfies
\[
 \|A(Z-x)\|_\infty<dR_m\quad(d>0),
\]
and $A(Z-x)=0$ when $d=0$.

If $x$ is an average of $2^K$ members of $V$, the series can be truncated at $K$. The hard radius is then $dR_m(1-2^{-K})$, the coefficient covariance is $\tau B(1-2^{-K})/2$, and the physical Gaussian-envelope variance is
\[
 \frac{\alpha d^2R_m^2}{3}(1-4^{-K}).
\]
The assertion can be restricted to any exposed face of $\conv V$ containing $x$.
\end{theorem}

\begin{proof}
First take $d>0$ and a list of $q=2^K$ points averaging to $x$. Pair the list. For each pair, property \eqref{R7-eq:all-subsets} says that
\[
 w=\frac{u+v}{2}+\frac12\sum_a\sigma_a c_a
\]
is feasible for every choice of signs. Stack the move columns for every pair in a single matrix $D$. Replacing each pair by its selected point changes the average by $D\sigma/q$. By \eqref{R7-eq:exchange-budgets}, $DD^{\mathsf T}\preceq qB/2$. Apply Theorem~\ref{MAIN-signing} to $AD/d$ across all pairs at once. Its columns have Euclidean norm at most one. Enlarge the physical cosine radius to $R_m$ and replace the auxiliary product law by its product Gaussian envelope. Conditional on the current list, the joint average increment is dominated in convex order by independent blocks
\[
 \left(\frac d q X_{R_m},\ N\left(0,\frac{\tau B}{2q}\right)\right).
\]
For the coefficient block, independent Gaussian padding uses $DD^{\mathsf T}\preceq qB/2$. The resulting source law depends only on the merge level, not on the current list.

Conditional Jensen composition now uses independent reference increments at the successive levels $q=2^k$. The physical sums are $d\sum_{k=1}^K2^{-k}X_{R_m}^{(k)}$. The coefficient covariances sum to $\tau B(1-2^{-K})/2$. The physical Gaussian variances sum as the geometric series $\alpha d^2R_m^2\sum_{k=1}^K4^{-k}$. Conditional means of the increments vanish, so the terminal mean remains $x$.

For arbitrary convex-combination weights, approximate them by dyadic weights. The starting means converge to $x$, the outputs remain in the same finite set $V$, and the reference sums converge in first moment. Passing to a subsequence of the joint martingale couplings gives \eqref{R7-eq:exchange-joint}. Every coordinate of $U_{d,m}$ lies strictly between $-dR_m$ and $dR_m$ almost surely. Its conditional mean at each positive-probability output retains the strict inequality. The strict bound follows from this conditional-mean statement at each output atom. If $d=0$, apply the argument with positive $d$ tending to zero; the coefficient reference is unchanged and the feature error vanishes.

Finally let a linear functional expose a face containing $u,v$. Each complementary pair in \eqref{R7-eq:all-subsets} sums to $u+v$. Both points are feasible and their functional values sum to twice the maximum; hence both are in the same face. This permits restriction to the face before running the proof.
\end{proof}

\subsubsection{Matroid bases and totally unimodular systems}
A matroid is \emph{strongly base-orderable} if any two bases admit a
bijection fixing their intersection such that exchanging any subset of
matched elements leaves both complementary sets as bases.
For a totally unimodular matrix $M$, a \emph{circuit} is a nonzero
kernel vector of inclusion-minimal support, normalized so that its
nonzero entries are $\pm1$. These familiar exchange decompositions
give the two applications of the common theorem below.

\begin{corollary}[Exact feasible outputs]\label{thm:matroid-hard}
\label{MAIN-matroid}\label{TU-thm:rounding}
Let $x$ belong to either a strongly base-orderable matroid base polytope
or $P=\{x\in[0,1]^n:Mx=b\}$ with totally unimodular $M$ and integral $b$.
Let $F$ be its minimal face and $H=\operatorname{span}(F-F)$.
There is a vertex $Z$ of $F$ with $\E Z=x$ and the same-law comparison
\eqref{R7-eq:exchange-joint}, with the following parameters:
\[
\begin{array}{c|c|c}
 &B&d\\ \hline
\text{strongly base-orderable bases}&2P_H&
 \max\|a_e-a_f\|_2\\[1mm]
\text{totally unimodular system}&\ell P_H&
 \max_c\|Ac\|_2.
\end{array}
\]
In the first row the maximum runs over possible single exchanges in
$F$. In the second row $c$ ranges over face-compatible circuits and
$\ell=\max_c|\operatorname{supp}c|$. Empty maxima are zero.
Every output has hard feature error less than $dR_m$ when $d>0$.
The corresponding coefficient Gaussian is $N(0,\tau P_H)$ or
$N(0,\tau\ell P_H/2)$, independent of the bounded physical reference.
\end{corollary}
\begin{proof}
For matroid bases, the matching gives disjoint moves $e_i-e_j$.
All subsets are feasible by strong base orderability, and
\[
 \sum_{(i,j)}\langle z,e_i-e_j\rangle^2
 \le2\sum_{(i,j)}(z_i^2+z_j^2)\le2\|z\|^2.
\]
Each move lies in $H$, so its Gram sum is at most $2P_H$.
An exposing functional for $F$ puts both complementary outputs in
$F$: their sum is the original pair and both attain the same maximal
functional value.

For the second case, take vertices $u,v$ and $w=u-v\in\{-1,0,1\}^n$.
Every nonzero kernel vector has a conformal circuit: intersect its
orthant and coordinate support with the kernel, then choose an
extreme ray. Its support is minimal; total unimodularity gives a
primitive representative with entries $\pm1$. Subtract this circuit
from $w$. Every coordinate it uses becomes zero and no sign changes.
Iteration gives a disjoint conformal decomposition $w=\sum_a c_a$.
Every subset exchanged from $u$ to $v$ stays in $[0,1]^n$ and satisfies
the equalities. Its complementary point does too. The same face
argument applies, and disjointness gives
\[
 \sum_a\langle z,c_a\rangle^2
 \le\sum_a|\operatorname{supp}c_a|
                  \sum_{i\in\operatorname{supp}c_a}z_i^2
 \le\ell\|P_Hz\|^2.
\]
Theorem~\ref{R7-thm:exchange} proves both conclusions.
\end{proof}

For a partition matroid this preserves all prescribed integer quotas
and all individual means. For assignments the circuits are alternating
cycles in a bipartite graph. For network flows they are signed cycles
of the oriented incidence matrix. Their feature sums give $d$ directly;
when each column has norm at most one, $d\le\ell$ suffices for the
unimodular case. The exchange-specific value can be much smaller.

\begin{corollary}[Integral capacities and inequalities]
\label{R7-cor:capacities}
Let $M$ be totally unimodular and let $l,u,b$ be integral. For any
$x\in[l,u]$ satisfying $Mx\le b$, there is an integral $Z$ with
\[
 \E Z=x,\qquad MZ\le b,\qquad
 Z_j\in\{\lfloor x_j\rfloor,\lceil x_j\rceil\}.
\]
Every initially integral coordinate and every tight face is preserved.
The preceding hard and joint-reference bounds hold for the circuits
of the augmented nearest-integer cell. Their constants depend on its
circuit geometry and feature sums, independently of the magnitudes
of the capacities.
\end{corollary}
\begin{proof}
Introduce slack $s=b-Mx\ge0$ and the totally unimodular equality
matrix $[M\ I]$. Intersect with the cell whose coordinates are between
the floors and ceilings of $(x,s)$. Translate by the lower corner;
fix integral coordinates. The resulting $0$--$1$ polytope is integral.
Extend $A$ by zero on the slack coordinates and apply the preceding
corollary, then project to the original coordinates. Nonnegative slack
preserves feasibility. If an inequality is tight at $x$, its expected
slack is zero and all output slacks are nonnegative, so it stays tight.
The minimal-face assertion also follows directly from the exchange
theorem.
\end{proof}

On an equality face, features differing by $LM$ have identical errors:
\begin{equation}\label{R7-eq:feature-quotient}
 (A+LM)(Z-x)=A(Z-x),\qquad (A+LM)c=Ac.
\end{equation}
Thus row and column potentials disappear from assignment-cycle sums,
and vertex potentials disappear from flow-cycle sums. This identity
can be used before bounding features by individual coefficients.
Exact circuit bounds themselves are unchanged.


\subsection{Gaussian rounding on integral base polyhedra}
\label{BASE-sec:main}
Let $E=\{1,\ldots,n\}$. For $z\in\R^E$ and $S\subseteq E$, write $z(S)=\sum_{i\in S}z_i$. A function $f:2^E\to\Z$ is submodular if
\[
 f(S)+f(T)\ge f(S\cup T)+f(S\cap T),
 \qquad S,T\subseteq E.
\]
Assume $f(\varnothing)=0$, and define its base polyhedron
\begin{equation}\label{BASE-eq:base}
 B(f)=\{z\in\R^E:z(S)\le f(S)\ (S\subseteq E),\ z(E)=f(E)\}.
\end{equation}
The function need not be monotone. Given $x\in B(f)$, put
\[
 a^0=\lfloor x\rfloor,\qquad p=x-a^0,\qquad J=\{i:0<p_i<1\},\qquad m=|J|.
\]
All floors and ceilings are coordinatewise. Let
\begin{equation}\label{BASE-eq:face}
 P=B(f)\cap\prod_i[\lfloor x_i\rfloor,\lceil x_i\rceil],
 \qquad F=\text{the minimal face of $P$ containing $x$},\qquad
 H=\lin(F-F).
\end{equation}
Write $P_H$ for orthogonal projection onto $H$. In particular, $H$ is supported on $J$ and lies in the zero-sum hyperplane. The sampler need not compute this projection. Integral input is returned unchanged.

The Gaussian reference $P_HG_n$ may be singular; the convex-order definition remains the one in the introduction.

\begin{theorem}[Gaussian rounding on an integral base polyhedron]\label{BASE-thm:main}
For every $x\in B(f)$ there is an integer random vector $N$ such that
\begin{equation}\label{BASE-eq:main-hard}
 \E N=x,\qquad N\in B(f),\qquad
 N_i\in\{\lfloor x_i\rfloor,\lceil x_i\rceil\}
 \quad\text{at every outcome},
\end{equation}
and
\begin{equation}\label{BASE-eq:main-cx}
 N-x\cx\sqrt{\pi/3}\,P_HG_n.
\end{equation}
Every inequality of $P$ tight at $x$ is therefore preserved as an equality.

If $2^Kx$ is integral, a finite $K$-stage sampler gives the stronger comparison
\begin{equation}\label{BASE-eq:dyadic}
 N-x\cx\sqrt{\frac\pi3(1-4^{-K})}\,P_HG_n.
\end{equation}
For rational $x$ and any rational $\varepsilon>0$, a finite exact sampler gives~\eqref{BASE-eq:main-hard} and
\begin{equation}\label{BASE-eq:epsilon}
 N-x\cx\sqrt{\frac{(1+\varepsilon)\pi}{3}}\,P_HG_n.
\end{equation}
With a value oracle for the integral submodular function, its expected bit complexity is polynomial in $n$, the input and oracle-value bit lengths, and $\log^+(1/\varepsilon)$. The oracle algorithm uses polynomial submodular minimisation. A dyadic input needs at most $mK/2$ fair binary choices; a general rational input needs at most $m+mK/2$ binary choices.
\end{theorem}
The endpoint~\eqref{BASE-eq:main-cx} for a general real input is a finite-alphabet compactness conclusion. The finite rational sampler has the explicit finite-depth bound in Theorem~\ref{BASE-thm:finite}; its variance slack is stated separately from that endpoint. Compact cut and coverage functions admit exact minimum-cut implementations. For an arbitrary function given by its complete value table, the complexity is polynomial in that table size, which is $2^n$.

\begin{corollary}[Half-integral points]\label{BASE-cor:half}
If every nonintegral coordinate of $x$ has fractional part $1/2$, one finite stage gives
\begin{equation}\label{BASE-eq:half}
 N-x\cx\frac{\sqrt\pi}{2}P_HG_n.
\end{equation}
The scalar $\sqrt\pi/2$ is optimal uniformly over this class. The finite sampler uses at most $m/2$ fair binary choices.
\end{corollary}
This includes every feasible half-integral point of a matroid base polytope.
Its exchanges have the form $e_i-e_j$. General transportation fibres can
have longer circuits; their separate circuit-dependent guarantees are
stated in Section~\ref{R7-sec:exchange}.

Every additional real matrix $A$ inherits from the same law
\begin{equation}\label{BASE-eq:allmaps}
 A(N-x)\cx\sqrt c\,AP_HG_n,
 \qquad \E e^{t^{\mathsf T}A(N-x)}\le
 e^{c\norm{P_HA^{\mathsf T}t}_2^2/2},
\end{equation}
where $c$ is the applicable covariance coefficient. These comparisons do not impose an extra almost-sure bound on an arbitrary map $A$.

The construction works inside the nearest integer cell, which is a translated matroid base polytope. Two bases on a tight face can be exchanged symmetrically while remaining on that face: the supporting functional forces its value on the exchange direction to vanish. During a merge, each exchanged coordinate pair is then fixed and never reused. The resulting pathwise bound on the sum of exchange matrices lies in $H$, so Gaussian completion preserves the exact mean and gives the projected reference $P_HG_n$. This is why the covariance depends on the minimal face while the sampler never needs to compute its projection.

\subsubsection{A unit box in a base polyhedron is a translated matroid base}
\label{BASE-sec:box}
The reduction retains the whole family of submodular inequalities. It is needed again at every dyadic scale, since the point being rounded may lie on new tight faces.

\begin{lemma}[Residual rank formula]\label{BASE-lem:rank}
Let $h:2^E\to\Z$ be submodular, with $h(\varnothing)=0$ and $h(S)\ge0$ for every $S$. Define
\begin{equation}\label{BASE-eq:rho}
 \rho(U)=\min_{T\subseteq E}\{h(T)+|U\setminus T|\}.
\end{equation}
Then $\rho$ is a matroid rank function and
\begin{equation}\label{BASE-eq:independent-polytope}
 \{z\in[0,1]^E:z(S)\le h(S)\ \forall S\}
 =\{z\ge0:z(S)\le\rho(S)\ \forall S\}.
\end{equation}
If this set contains $z$ with $z(E)=h(E)$, then $\rho(E)=h(E)$ and its section of that total is the base polytope of $\rho$.
\end{lemma}
\begin{proof}
The function $\rho$ is integer-valued, monotone, normalised, and at most $|U|$. Monotonicity follows termwise in~\eqref{BASE-eq:rho}; normalisation uses $h\ge0$ and $h(\varnothing)=0$; the bound uses $T=\varnothing$. Adding one element changes each term by at most one, hence $\rho(U+i)-\rho(U)\le1$.

Take minimisers $T_1,T_2$ for $U,V$. Submodularity of $h$ and the elementary cardinality inequality
\[
 |U\setminus T_1|+|V\setminus T_2|
 \ge |(U\cup V)\setminus(T_1\cup T_2)|
       +|(U\cap V)\setminus(T_1\cap T_2)|
\]
prove submodularity of $\rho$. The rank axioms follow. Explicitly, declare $I$ independent when
$\rho(I)=|I|$. The unit increment bound gives deletion, and diminishing
returns gives augmentation: if no element of $J\setminus I$ augmented
$|I|<|J|$, then $\rho(I\cup J)=\rho(I)<\rho(J)$, a contradiction.

For $z$ in the left side of~\eqref{BASE-eq:independent-polytope}, $z(U)\le z(T)+|U\setminus T|\le h(T)+|U\setminus T|$. Minimise in $T$. Conversely, $\rho(S)\le h(S)$ and $\rho(\{i\})\le1$ give all left-side constraints. Finally, a feasible point of total $h(E)$ implies $\rho(E)\ge h(E)$, while $T=E$ gives the reverse bound. The standard greedy description of a matroid base polytope identifies its vertices with the base indicators \cite{BASE-CVZ}.
\end{proof}

\begin{corollary}[Scaled adjacent boxes]\label{BASE-cor:box}
Suppose $y\in B(f)$, $t$ is a positive integer, and $a=\lfloor ty\rfloor$. Put $h_t(S)=tf(S)-a(S)$. Then $h_t$ satisfies Lemma~\ref{BASE-lem:rank}, and
\begin{equation}\label{BASE-eq:scaled-box}
 \{z\in[0,1]^E:a+z\in tB(f)\}=B(\rho_t),
 \qquad \rho_t(U)=\min_T\{tf(T)-a(T)+|U\setminus T|\}.
\end{equation}
Coordinates with zero residual may be fixed at zero by restriction of this matroid. In that restriction the total rank remains $tf(E)-a(E)$.
\end{corollary}
\begin{proof}
Feasibility of $y$ and $a\le ty$ imply $tf(S)-a(S)\ge0$ for all $S$. The residual $ty-a$ is feasible in the box and has total $tf(E)-a(E)$, so apply the lemma. After fixing zero coordinates, this same point witnesses that the rank on the remaining coordinates is still the prescribed total.
\end{proof}
This proves integrality of the adjacent-box intersection using an explicit rank formula. It also explains why the construction extends to finite integral base polyhedra whose coordinates need not be binary.

\begin{lemma}[Complementary bases at a midpoint]\label{BASE-lem:partition}
If $\1_S/2$ belongs to a matroid base polytope on $S$, then $|S|$ is even and $S$ partitions into two bases. Such a partition can be found with a polynomial number of independence queries.
\end{lemma}
\begin{proof}
Its rank is $|S|/2$. For every $T\subseteq S$, feasibility gives $|T|\le2\rho(T)$. Edmonds' partition theorem says precisely that this condition permits a partition into two independent sets \cite{BASE-Edmonds}. Their total size is twice the rank, so both are bases. The augmenting-path matroid-partition algorithm is polynomial in independence queries; see \cite{BASE-Terao} for the classical bound and subsequent improvements.
\end{proof}
The augmenting-path algorithm maintains two independent sets. Each arc records an allowed one-element exchange, so complementary bases are obtained without enumerating all bases.

\subsubsection{The Gaussian cost of symmetric basis exchange}
\label{BASE-sec:merge}
Let $\phi,\Phi$ be the standard Gaussian density and distribution function. For $0<p<1$, put
\begin{equation}\label{BASE-eq:tau}
 \tau(p)=\frac{p(1-p)}{\phi(\Phi^{-1}(p))},\qquad \tau(0)=\tau(1)=0.
\end{equation}
Lemma~\ref{lem:bernoulli-gaussian-exact} gives
\begin{equation}\label{BASE-eq:Bernoulli}
 B-p\cx\tau(p)G_1,\qquad B\sim\operatorname{Bernoulli}(p),
 \qquad \tau(p)^2\le\pi/8.
\end{equation}
The tail coupling and its sharpness are proved in
Lemma~\ref{lem:bernoulli-gaussian-exact}. In particular, a fair sign has
Gaussian reference variance $\pi/2$.

We state the probabilistic assembly used throughout.
\begin{lemma}[Adaptive covariance completion]\label{BASE-lem:completion}
Let $(\xi_t)$ be the increments of a finite martingale tree. Suppose,
conditionally on its past, $\xi_t\cx N(0,Q_t)$, where $Q_t$ is determined
before that step. If $\sum_tQ_t\preceq C$ on every path for a deterministic
$C\succeq0$, then $\sum_t\xi_t\cx N(0,C)$.
\end{lemma}
\begin{proof}
Apply Theorem~\ref{thm:gaussian-composition}. Its coupling preserves the
entire discrete path, including all feasibility decisions.
\end{proof}
A bound on the expected covariance alone would not suffice for this argument. All bounds below hold for every path.

Let $B_0,B_1$ be bases of one matroid and let $0<p<1$. Start at
$z=p\1_{B_0}+(1-p)\1_{B_1}$. Choose $i\in B_0\setminus B_1$. Symmetric basis exchange gives $j\in B_1\setminus B_0$ for which both $B_0-i+j$ and $B_1-j+i$ are bases. With probability $p$, replace $B_1$ by $B_1-j+i$; otherwise replace $B_0$ by $B_0-i+j$. Iterate until the bases coincide. This is the classical \texttt{MergeBases} procedure \cite{BASE-CVZ}.

\begin{theorem}[A complete two-base Gaussian reference]\label{BASE-thm:merge}
The merge returns a base $B$ with
\[
 \E\1_B=p\1_{B_0}+(1-p)\1_{B_1}.
\]
If $H$ is the direction space of any face containing the two initial bases, then
\begin{equation}\label{BASE-eq:merge-bound}
 \1_B-p\1_{B_0}-(1-p)\1_{B_1}
 \cx\sqrt2\,\tau(p)P_HG_n.
\end{equation}
It takes at most $|B_0\setminus B_1|$ binary choices and polynomially many rank queries. For fixed $p$, $\sqrt2\tau(p)$ is the smallest scalar valid uniformly over all matroids and base pairs when the reference is projected onto their minimal face.
\end{theorem}
\begin{proof}
With $d=e_i-e_j$, the weighted current vector changes by $(1-p)d$ with probability $p$ and by $-pd$ with probability $1-p$. Its conditional mean is unchanged and its exact Gaussian covariance is $\tau(p)^2dd^{\mathsf T}$ by~\eqref{BASE-eq:Bernoulli}.

After either update, the two bases agree at both $i$ and $j$. Those coordinates are never used in later exchanges of this merge. Thus all the exchange vectors along a path have disjoint supports and squared length two. They belong to $H$: if both bases maximise a supporting linear functional, the two valid cross-exchanged bases force its value on $e_i-e_j$ to be zero. Consequently
\begin{equation}\label{BASE-eq:pair-Gram}
 \sum_{\text{exchanges}}dd^{\mathsf T}\preceq2P_H
\end{equation}
on every path, even when the pair choices depend on earlier coin outcomes. Lemma~\ref{BASE-lem:completion} proves~\eqref{BASE-eq:merge-bound}. Agreement increases at every step, giving the stated count. Testing candidates for $j$ gives polynomial rank-query complexity.

For necessity, use a rank-one matroid on two elements and the mean $(p,1-p)$. The only feasible law is the Bernoulli choice between its two bases. The unit contrast $(1,-1)/\sqrt2$ has error $\sqrt2(B-p)$, whose exact Gaussian scale is $\sqrt2\tau(p)$. This proves uniform optimality.
\end{proof}
The new bound is a bound for the entire classical merge. Treating its exchanges separately without~\eqref{BASE-eq:pair-Gram} would introduce the rank. The covariance can be completed to any larger deterministic matrix; no Gaussian variable needs to be sampled by the terminal algorithm.

\subsubsection{Coarsening the grid while preserving every base inequality}
\label{BASE-sec:dyadic}
Suppose $2^Ky$ is integral and $y\in P$. At a level $k\ge1$, put
\[
 t=2^{k-1},\qquad a=\lfloor ty\rfloor,\qquad r=ty-a.
\]
All coordinates of $r$ are zero or $1/2$. Let $S$ be the half-coordinate set. Corollary~\ref{BASE-cor:box} and Lemma~\ref{BASE-lem:partition} give complementary bases $B_0,B_1$ on $S$ in the residual matroid, with
$r=(\1_{B_0}+\1_{B_1})/2$. Merge these bases with $p=1/2$ and set
\begin{equation}\label{BASE-eq:coarsen}
 y'=(a+\1_B)/t.
\end{equation}
The new point belongs to $B(f)$, has the same conditional mean as $y$, and lies on the $2^{-(k-1)}$ grid. A coordinate with even numerator stays fixed; an odd numerator changes by $\pm1$ at denominator $2^k$. Hence all stages remain inside the original adjacent-integer box.

\begin{proposition}[Pathwise covariance at one level]\label{BASE-prop:level}
Every elementary exchange of level $k$ changes the current physical vector by
$\pm2^{-k}(e_i-e_j)$. Its Gaussian covariance is
\[
 Q_{k,ij}=\frac\pi2\,4^{-k}(e_i-e_j)(e_i-e_j)^{\mathsf T}.
\]
The sum for that level is at most $\pi4^{-k}P_H$ on every path.
\end{proposition}
\begin{proof}
The factor $t^{-1}$ in~\eqref{BASE-eq:coarsen} scales the midpoint increments $\pm d/2$ to $\pm2^{-k}d$. Theorem~\ref{BASE-thm:merge} gives the level sum by disjointness. All intermediate weighted vectors are in $P$ and all their positive-probability transitions are unbiased. A supporting equality at the initial mean is consequently preserved at every transition: nonnegative slack with zero conditional expectation is zero on both branches. Every exchange direction therefore belongs to the initial minimal-face space $H$.
\end{proof}

Summing over levels gives the deterministic budget
\begin{equation}\label{BASE-eq:sum}
 \sum_{k,ij}Q_{k,ij}\preceq
 \pi\sum_{k=1}^K4^{-k}P_H
 =\frac\pi3(1-4^{-K})P_H.
\end{equation}
Lemma~\ref{BASE-lem:completion} proves~\eqref{BASE-eq:dyadic}. At most $m/2$ disjoint pairs occur at each level. The case $K=1$ proves Corollary~\ref{BASE-cor:half}; its uniform lower bound is the rank-one half-pair in Theorem~\ref{BASE-thm:merge}.

\paragraph{Rational inputs and an explicit initial cost}
For arbitrary rational $x$, write
\begin{equation}\label{BASE-eq:initial}
 x=2^{-K}(a+r),\qquad a=\lfloor2^Kx\rfloor,\qquad 0\le r_i<1.
\end{equation}
The residual belongs to the base polytope in Corollary~\ref{BASE-cor:box}. Round $r$ to an integral base indicator $M$ by ordinary two-coordinate pipage steps. At a nonvertex point, choose $e_i-e_j$ feasible in both directions in its minimal face, and move to the two endpoints with their mean-preserving probabilities. Each endpoint lies in a proper face, so at most $m-1$ nontrivial moves are needed; we use $m$ as a uniform bound. The existence of such an exchange direction follows either from the matroid base-edge description or the tight-set argument below.

For completeness, tight rank sets form a lattice under union and intersection. A maximal chain in this lattice partitions the ground set into blocks. Every tight set is a union of whole blocks: otherwise its intersection with a chain set, united with the preceding chain set, would insert a strict intermediate tight set. Differences of consecutive chain indicators are the block indicators, so these span all tight-set incidence vectors. A direction $e_i-e_j$ is feasible in both directions exactly when $i,j$ belong to the same block, after coordinate-bound tight sets are included. If every block is a singleton, all coordinates are fixed and the point is a vertex. Otherwise one such pair exists. At a maximal endpoint at least one new independent tight equality appears, reducing the face dimension. This also justifies the move count.

Let the two step lengths be $\alpha,\beta>0$. Since the residual remains in the unit cube and the direction has unit entries, $\alpha+\beta\le1$. The exact scalar Gaussian variance is
$(\alpha+\beta)^2\tau(\alpha/(\alpha+\beta))^2$. After scaling by $2^{-K}$, it is at most $\pi4^{-K}/8$ on $dd^{\mathsf T}$ and hence at most $\pi4^{-K}P_H/4$. The initial sum is at most $\pi m4^{-K}P_H/4$.

\begin{theorem}[Exact finite-depth reference]\label{BASE-thm:finite}
Round~\eqref{BASE-eq:initial} to $Y=2^{-K}(a+M)$ as above and then coarsen all $K$ levels. The resulting law satisfies~\eqref{BASE-eq:main-hard} and
\begin{equation}\label{BASE-eq:cK}
 N-x\cx\sqrt{c_K}\,P_HG_n,\qquad
 c_K=\pi\left[\frac{1-4^{-K}}3+\frac{m4^{-K}}4\right].
\end{equation}
The second term is omitted when initial gridding is unnecessary. It takes at most $m+mK/2$ binary choices. All probabilities and state coordinates are rational.
\end{theorem}
\begin{proof}
The initial stage and the dyadic stages are one finite martingale tree. Their pathwise covariance bounds add, so Lemma~\ref{BASE-lem:completion} gives~\eqref{BASE-eq:cK}. Feasibility, adjacent-integer support, and the exact mean have been preserved at every step.
\end{proof}
Choose $4^K\ge3m/(4\varepsilon)$. Then
\[
 c_K=\frac\pi3\left[1+\left(\frac{3m}4-1\right)4^{-K}\right]
 \le\frac{(1+\varepsilon)\pi}3.
\]
This proves the finite rational bound of Theorem~\ref{BASE-thm:main}. For real input, the same finite trees with real probabilities exist for every $K$. Their terminal laws lie on a fixed finite subset of the adjacent box. A convergent subsequence has mean $x$ and the same hard support, and its Gaussian references converge in first moment to $N(0,(\pi/3)P_H)$. Convex order passes to the limit by first testing convex Lipschitz functions and then increasing maxima of affine functions. This proves~\eqref{BASE-eq:main-cx}.

\subsubsection{Oracle complexity and compact minimum-cut implementations}
\label{BASE-sec:oracles}
Every residual rank query has the form
\begin{equation}\label{BASE-eq:rank-query}
 \rho_t(U)=|U|+\min_T\{tf(T)-a(T)-|U\cap T|\}.
\end{equation}
The expression minimised is submodular. Thus one submodular minimisation computes a rank value; comparison with $|U|$ answers independence. Polynomial value-oracle minimisation of integral submodular functions is classical \cite{BASE-IFF}.

For an initial pipage step at residual $r$, the maximal feasible positive length in direction $e_i-e_j$ is
\begin{equation}\label{BASE-eq:step-query}
 \beta_{ij}=\min\left\{1-r_i,r_j,
 \min_{T\ni i,\,T\not\ni j}\bigl(tf(T)-a(T)-r(T)\bigr)\right\}.
\end{equation}
The negative length is $\beta_{ji}$. Fixing membership and nonmembership preserves the submodular minimisation problem on the remaining elements. Searching all pairs finds a two-sided feasible direction. The formula includes the coordinate bounds and proves exact maximality. It uses the original function $f$, avoiding a nested minimisation over $\rho_t$.

With a standard shortest augmenting-path partition algorithm, each level uses $O(n^3)$ independence tests and $O(n^2)$ tests for its merges. The initial stage uses $O(n^3)$ constrained minimisations. Hence $O((K+1)n^3)$ submodular minimisation calls suffice. The rank identities and the matroid partition theorem determine these oracle calls directly.

Let $Q$ be a common denominator of the input point. During the initial stage all residual coordinates and endpoint distances remain multiples of $1/Q$, since each minimised slack is an integer minus a sum of such coordinates. The initial probabilities therefore have polynomial bit length. Later states use at most $K$ additional dyadic bits. A rational coin with denominator $d$ is sampled by rejection from a uniform bit string of length $\lceil\log_2d\rceil$; acceptance probability is at least $1/2$. This proves the expected bit bound in Theorem~\ref{BASE-thm:main}. The construction does not compute $\pi$, a Gaussian quantile, or a covariance square root.

\paragraph{Directed-cut base polyhedra}
Let a directed graph have integral capacities $c_{ij}\ge0$ and integral modular coefficients $b_i$. Define
\begin{equation}\label{BASE-eq:cut}
 f(S)=b(S)+\sum_{i\in S,\,j\notin S}c_{ij}.
\end{equation}
A query $\min_S[tf(S)+v(S)]$ is an exact minimum cut. An original arc has capacity $tc_{ij}$. A positive unary coefficient $v_i+tb_i$ becomes an arc from $i$ to the sink; a negative one becomes an arc from the source to $i$, with the corresponding negative constant added to the cut value. Forced membership or nonmembership is imposed by an arc with capacity greater than the sum of all finite capacities. All capacities are rational and have polynomial bit length.

If $x-b$ is the net outgoing balance of a feasible fractional flow, then $x\in B(f)$. Conversely, the cut inequalities and the total equality are the capacitated-flow feasibility conditions, by the max-flow/min-cut theorem. Therefore the sampler rounds node balances to adjacent integers, preserving their exact expectations and every cut capacity. A second integral maximum-flow calculation recovers a feasible integral flow for those rounded balances. The guarantee concerns the means of the node balances; it does not prescribe the mean of each recovered edge flow.

\paragraph{Coverage allocations}
Let resource $a$ have an integer capacity $c_a\ge0$ and a nonempty set $E_a$ of eligible bins. Put
\begin{equation}\label{BASE-eq:coverage}
 f(S)=\sum_a c_a\1_{\{E_a\cap S\ne\varnothing\}}.
\end{equation}
This is an integral monotone submodular function. Its base polytope is the set of bin totals obtained by distributing every resource among its eligible bins. The equality follows from the bipartite flow criterion: all subset upper inequalities bound how much resource can enter those bins, and the total uses every resource.

For minimisation of $tf(S)+v(S)$, introduce one node for each resource. Add a large-capacity arc from eligible bin $i$ to resource $a$, and an arc of capacity $tc_a$ from $a$ to the sink. Selecting any eligible bin on the source side forces the resource node onto that side and pays its capacity. The same unary construction handles $v$. This proves a compact exact minimum-cut oracle for~\eqref{BASE-eq:rank-query} and~\eqref{BASE-eq:step-query}.

The sampled adjacent-integer bin totals can be realised by an integral bipartite flow. Every resource is used in full, every eligibility restriction holds, and each bin total has the original prescribed mean. Again, the edge-level allocation means are not prescribed by this reconstruction.

Both compact classes can be implemented with a rational Edmonds--Karp minimum-cut routine. Their running times are polynomial in the compact graph or incidence description and rational bit length, rather than in an explicit table of all subsets. For a general submodular function represented by its complete table, each minimisation scans the table, whose size is therefore part of the input.

\subsubsection{An endpoint construction using only explicit bases and rank queries}
\label{BASE-sec:tree}
A second access model avoids submodular minimisation. Suppose a matroid is given by a rank oracle and the input includes an explicit decomposition
\begin{equation}\label{BASE-eq:decomposition}
 x=\sum_{i=1}^M\lambda_i\1_{B_i},\qquad
 \lambda_i>0,\quad\sum_i\lambda_i=1,
\end{equation}
with rational weights. Such decompositions are standard inputs to swap rounding \cite{BASE-CVZ}. The input includes all $M$ terms of the decomposition, and the running time depends on $M$.

For any fixed binary merge tree with leaf weights $\lambda_i$, let $w_v$ be an internal node's total weight and $p_v$ the fraction in its first child. Merging its independently constructed child bases with the classical operation gives a finite martingale. Theorem~\ref{BASE-thm:merge}, scaled by $w_v$, proves
\begin{equation}\label{BASE-eq:tree-cost}
 N-x\cx\left(2\sum_{v\text{ internal}}w_v^2\tau(p_v)^2\right)^{1/2}P_HG_n.
\end{equation}
The node weights and ratios are deterministic. The nodes can be processed bottom-up, so their global weighted increments have mean zero; covariance completion applies to that one history. A sequential tree can have a large sum in~\eqref{BASE-eq:tree-cost}. A compressed balanced tree has a uniform bound.

\begin{theorem}[Finite balanced-tree endpoint]\label{BASE-thm:tree}
Given~\eqref{BASE-eq:decomposition}, choose $K\ge1$ and put $\delta=2^{-K}$. There is a finite rank-oracle sampler with
\begin{equation}\label{BASE-eq:tree-bound}
 \E N=x,\quad N\text{ a base indicator},\qquad
 N-x\cx\sqrt{\widetilde c_K}\,P_HG_n,
\end{equation}
where
\begin{equation}\label{BASE-eq:tree-coefficient}
 \widetilde c_K=\pi\left[\frac{1-2^{-K}}2+\frac{M-1}{4}4^{-K}\right].
\end{equation}
Here $H$ is the direction space of the minimal base-polytope face containing $x$. The algorithm has $O(M(K+1))$ two-base merges. Choosing $2^K\ge M$ gives the finite endpoint $\widetilde c_K\le\pi/2$, with expected polynomial bit complexity in the explicit decomposition and rank-oracle model. No approximation of the weights is needed.
\end{theorem}
\begin{proof}
Partition $[0,1]$ into consecutive intervals of lengths $\lambda_i$, labelled by bases $B_i$. Subdivide dyadically to depth $K$. A dyadic cell contained in one label interval is a deterministic leaf and is not subdivided further. Within a mixed depth-$K$ cell, merge the labels intersecting that cell sequentially with their exact relative intersection lengths. Above those leaves, merge the two children of each mixed cell with equal weights.

At depth $d<K$, at most $2^d$ merges of global weight $2^{-d}$ occur. Their uniform Gaussian cost is at most
\[
 \frac\pi4\sum_{d=0}^{K-1}2^d4^{-d}
 =\frac\pi2(1-2^{-K}).
\]
This uses $\tau(1/2)^2=\pi/8$ and the factor two in~\eqref{BASE-eq:tree-cost}. In a mixed terminal cell of mass $\delta$, every sequential merge has global weight at most $\delta$. Each costs at most $\pi\delta^2/4$. The number of such merges, summed over terminal cells, is at most $M-1$, since each is charged to a distinct boundary between original label intervals inside that cell. This gives the second term of~\eqref{BASE-eq:tree-coefficient}.

At every depth only cells containing an original interval boundary can be mixed, so at most $M-1$ of them need work. Thus the actual number of merges is $O(M(K+1))$, even if the complete balanced tree would be very large. Exact interval intersections and relative weights are rational with bit length polynomial in the original weights and $K$. Each merge takes polynomial rank queries and exact rational coins. If $2^K\ge M$, then $(M-1)2^{-K}<1$, and~\eqref{BASE-eq:tree-coefficient} is less than $\pi/2$. The mean and face are preserved by every merge.
\end{proof}
The two main access models have different guarantees. Submodular minimisation gives variance arbitrarily close to $\pi/3$ in finite rational time. An explicit decomposition and a rank oracle alone give a finite variance at most $\pi/2$ with the simple compressed implementation. The one-merge theorem attains its sharper $2\tau(p)^2$ whenever the input has just two bases and that direct construction is chosen.

\subsubsection{Marginal variances, negative correlation, and simultaneous losses}
\label{BASE-sec:covariance}
All elementary moves in both constructions are centred transfers between two coordinates. This retains the classical negative-correlation calculation for pipage and swap rounding \cite{BASE-CVZ}, alongside the new full Gaussian reference.

\begin{proposition}[An exact covariance Laplacian]\label{BASE-prop:covariance}
For the adjacent-integer sampler put $Y=N-a^0$ and $p=x-a^0$. There are nonnegative weights $w_{ij}$ such that
\begin{equation}\label{BASE-eq:Laplacian}
 \Cov(Y)=\sum_{i<j}w_{ij}(e_i-e_j)(e_i-e_j)^{\mathsf T},
 \qquad \sum_{j\ne i}w_{ij}=p_i(1-p_i).
\end{equation}
Consequently
\begin{equation}\label{BASE-eq:actualcov}
 \Cov(N-x)\preceq2P_H\diag(p_i(1-p_i))P_H.
\end{equation}
For every $S\subseteq E$,
\begin{equation}\label{BASE-eq:negative}
 \E\prod_{i\in S}Y_i\le\prod_{i\in S}p_i,
 \qquad
 \E\prod_{i\in S}(1-Y_i)\le\prod_{i\in S}(1-p_i).
\end{equation}
The same assertions hold for the explicit-base tree sampler with $a^0=0$ and $p=x$. They pass to endpoint limits of the finite laws.
\end{proposition}
\begin{proof}
Every elementary martingale increment is $D_t(e_i-e_j)$ with conditional scalar variance $v_t\ge0$. Martingale increments are orthogonal in second moment. Sum their expected conditional covariance matrices, putting
$w_{ij}=\E\sum_{t:\{i_t,j_t\}=\{i,j\}}v_t$. The terminal coordinate is Bernoulli of mean $p_i$, so its diagonal variance is exactly $p_i(1-p_i)$, proving~\eqref{BASE-eq:Laplacian}. For any $u$, the inequality $(u_i-u_j)^2\le2u_i^2+2u_j^2$ gives the diagonal bound; projection by $P_H$ is valid because the error lies in $H$.

At an elementary transfer, a product over $S$ is affine in $D_t$ if at most one affected coordinate lies in $S$. If both lie in $S$, its conditional expectation decreases by
$v_t$ times the product of the remaining coordinates. These lie in $[0,1]$ throughout. The same calculation applies to the complementary coordinates. Thus both products are supermartingales, proving~\eqref{BASE-eq:negative}. Compact finite terminal support permits all expectations to pass to a limit.
\end{proof}
The factor two in~\eqref{BASE-eq:actualcov} is attained by the half-pair example. This is an estimate for the actual covariance of the same law. Its marginal variances are not substituted for the Gaussian convex-order covariance: rare Bernoulli variables already exclude a universal such substitution, as follows already from the sharp scalar scale in Lemma~\ref{lem:bernoulli-gaussian-exact}.

For any norm $\norm{\cdot}_{\mathcal K}$ and any further map $A$,~\eqref{BASE-eq:allmaps} yields
\[
 \E\norm{A(N-x)}_{\mathcal K}
 \le\sqrt c\,\E\norm{AP_HG_n}_{\mathcal K}.
\]
The law is independent of the chosen norm. A deterministic outcome selected for that particular norm may depend on the norm. Feasible spanning trees, integer allocation loads, and network balances can therefore be compared for multiple convex objectives through one sampled law.


\section{Sampling laws chosen before the integrand}
\label{R11-sec:sampling}
Optimal Sobolev orders and feature-based integration are established in~\cite{KriegNovak,KunschRudolf,Kunsch,Bach,DM}. Here $N$ equal-weight points must also retain exact marginals and hard features under one law for later integrands. Dyadic pairing preserves each local distribution, while Gaussian projection separates feature cost from residual variance. Polynomial-exact formulas cancel low frequencies; a Hilbert sum covers all smoothness below the design order. Fixed-integrand error and expected worst-case kernel error keep their different orders.

\subsection{Prescribed means with independent references}\label{MEAN-sec}
The signing law applies to every subset of the columns. A dyadic
rounding scheme can therefore use it conditionally at each level while
preserving the prescribed coordinate means. The complete level
references are chosen from deterministic data. Summing independent
copies of these references retains a joint convex-order comparison for
the final error.
\subsubsection{Independent references at every prescribed mean}\label{r3:sec:dyadic}
For every column $a_j\ne0$, set
\[
 R_j=C\|a_j\|_2-\gamma\frac{\sum_i|a_{ij}|^3}{\|a_j\|_2^2},
\]
and set $R_j=0$ for a zero column. Write $R_I=\max_{j\in I}R_j$,
with $R_\varnothing=0$. Theorem~\ref{MAIN-signing} constructs, for every
$I$, a symmetric law $S^I$ with
\begin{equation}\label{r3:assumption:signing}
 (A_IS^I,S^I)\cx(X_{R_I},T_I),\qquad X_{R_I}\perp T_I.
\end{equation}
Here $X_0=0$, $T_j\sim q_\dagger$, $w=\E T_j^2$, and
\[
 v_0=\frac13-\frac2{\pi^2},\qquad
 \alpha=\frac\pi2\left(\frac12-\frac2{\pi^2}\right)^2,\qquad
 \tau=\frac\pi2(\E|T_j|)^2.
\]
The scalar comparisons are \eqref{eq:reference-gaussian}.

Fix prescribed means $p=(p_1,\ldots,p_n)\in[0,1]^n$. For $k\ge1$, define
\begin{equation}\label{r3:eq:dyadicdata}
 r_{j,k}=2\operatorname{dist}(2^{k-1}p_j,\mathbb Z),\qquad
 J_k=\{j:r_{j,k}>0\},\qquad R_k=R_{J_k}.
\end{equation}
Thus $r_{j,k}\in[0,1]$. A coordinate of exact dyadic denominator $2^K$ belongs to $J_1,\ldots,J_K$ and to no later $J_k$. A nondyadic coordinate in $(0,1)$ belongs to every $J_k$. All deterministic coordinates can be omitted from the formulas below.

Take independent copies of all reference variables and set
\begin{equation}\label{r3:eq:dyadicreference}
 \widehat X=\sum_{k\ge1}2^{-k}X_{R_k}^{(k)},\qquad
 \widehat T_j=\sum_{k:j\in J_k}2^{-k}T_j^{(k)}.
\end{equation}
These series converge absolutely almost surely and in $L^1$. The physical block and the coefficient block are independent, and the coordinates within each block are independent.

\begin{theorem}[Prescribed means and independent references]\label{r3:thm:dyadic}\label{MAIN-means}\label{NEW-thm:means}
For the hereditary reference \eqref{r3:assumption:signing}, there is a random vector $Z\in\{0,1\}^n$ with $\mathbb EZ=p$ such that, for $E=Z-p$,
\begin{equation}\label{r3:eq:jointbiased}
 (AE,E)\cx(\widehat X,\widehat T).
\end{equation}
Equivalently, the reference can be coupled to $Z$ so that
\[
 \mathbb E[\widehat X\mid Z]=A(Z-p),\qquad
 \mathbb E[\widehat T\mid Z]=Z-p.
\]
If $B_A(p)=\sum_{k\ge1}2^{-k}R_k>0$, every outcome of positive probability satisfies
\begin{equation}\label{r3:eq:biasedhard}
 \|A(Z-p)\|_\infty<B_A(p)\le\max_jR_j.
\end{equation}
If $B_A(p)=0$, then $A(Z-p)=0$ almost surely.

Suppose that the Gaussian-envelope hypotheses above hold. Put
\[
 v_G=\alpha\sum_{k\ge1}4^{-k}R_k^2,
 \qquad
 D_G=\tau\operatorname{diag}\left(
               \sum_{k:j\in J_k}4^{-k}\right).
\]
On the nondeterministic coordinates, when $v_G>0$,
\begin{equation}\label{r3:eq:biasedgaussian}
 Z-p\cx N(0,Q_G),\qquad
 Q_G=\left(D_G^{-1}+v_G^{-1}A^{\mathsf T}A\right)^{-1}.
\end{equation}
Zero physical blocks are handled by omitting the corresponding constraint.
\end{theorem}

The uniform consequence of \eqref{r3:eq:biasedgaussian} is
\begin{equation}\label{r3:eq:uniformbiasedgaussian}
 Z-p\cx N\!\left(0,
 \frac13\left(\tau^{-1}I+(\alpha R^2)^{-1}A^{\mathsf T}A\right)^{-1}\right),
 \qquad R=\max_jR_j.
\end{equation}
For biased signs $S=2Z-1$ with prescribed mean $b=2p-1$, multiply the error by $2$. The hard radius becomes $2B_A(p)$ and the Gaussian covariance becomes $4Q_G$. When all means have dyadic denominator dividing $2^K$, the hard radius in binary coordinates is at most $(1-2^{-K})R$, and the factor $1/3$ in \eqref{r3:eq:uniformbiasedgaussian} can be replaced by $(1-4^{-K})/3$. In particular, $p_j=1/2$ recovers the original signing scale exactly.

\begin{lemma}\label{r3:lem:dyadicprob}
In unbiased dyadic rounding of a scalar $p\in[0,1]$, the probability of being active at level $k$ is $2\operatorname{dist}(2^{k-1}p,\mathbb Z)$. Moreover,
\begin{equation}\label{r3:eq:dyadicenergy}
 \sum_{k\ge1}4^{-k}\,2\operatorname{dist}(2^{k-1}p,\mathbb Z)=p(1-p).
\end{equation}
These statements hold for each coordinate of an arbitrarily dependent unbiased dyadic rounding.
\end{lemma}
\begin{proof}
At mesh $2^{-k}$, the coordinate lies at the two adjacent grid points bracketing $p$. Its mean is $p$, so the two probabilities are determined. Exactly one of those grid points has odd numerator unless the coordinate is already on the coarser mesh. Summing its probability gives the stated triangular function.

Start at a random point on a mesh $2^{-K}$ obtained by unbiased rounding between the two grid points bracketing $p$. A level-$k$ active increment has square $4^{-k}$ and conditional mean zero. Orthogonality of martingale increments therefore gives
\[
 \sum_{k=1}^K4^{-k}r_k=p(1-p)-\operatorname{Var}(p^{(K)}).
\]
The final variance on the right is at most $4^{-K}/4$. Letting $K$ tend to infinity proves \eqref{r3:eq:dyadicenergy}. The argument uses only conditional unbiasedness, so dependence between different coordinates does not affect it.
\end{proof}

\begin{proof}[Proof of Theorem~\ref{r3:thm:dyadic}]
First take a dyadic initial point. At level $k$, let $I_k$ be the coordinates with odd numerator and apply the signing law from \eqref{r3:assumption:signing} to $A_{I_k}$. Add $2^{-k}S^{I_k}$ on the active coordinates. Each coordinate remains in $[0,1]$, moves to the coarser grid, and preserves its conditional mean.

Every possible active set is contained in $J_k$. Conditional on the preceding choices, the joint increment of the physical and coefficient errors is dominated by $2^{-k}(X_{R_k},T_{J_k})$, with zeros on the omitted coefficient coordinates. Indeed, contraction of a centred reference gives the smaller physical scale $R_{I_k}$, and independent centred coordinates can fill $J_k\setminus I_k$. These operations preserve the product law of the larger reference. Applying conditional Jensen's inequality successively shows that the sum of the adaptive increments is dominated by the sum of independent level references. This proves the finite-series version of \eqref{r3:eq:jointbiased}.

For general $p$, first round each coordinate unbiasedly onto a mesh $2^{-K}$, and then perform the preceding construction. The initial error tends uniformly to zero. The level probabilities and the sets $J_k$ are still those of \eqref{r3:eq:dyadicdata}. The output laws have a convergent subsequence because they are supported on a finite set. The reference partial sums converge in $L^1$; hence their couplings are tight and their first moments are uniformly integrable. Passing the martingale identity to a subsequential limit gives \eqref{r3:eq:jointbiased} and the displayed conditional expectations.

Each coordinate of $\widehat X$ is strictly between $-B_A(p)$ and $B_A(p)$ almost surely when this number is positive. Its conditional expectation on every positive-probability output is also strictly between the same endpoints. This establishes the strict hard bound without taking a limit of strict inequalities.

Finally, independent Gaussian comparisons tensorize and add. They dominate $(\widehat X,\widehat T)$ by the independent Gaussian blocks of covariance $v_G I$ and $D_G$. Proposition~\ref{r3:thm:short} gives \eqref{r3:eq:biasedgaussian}. The uniform statements follow by summing the geometric series and enlarging the Gaussian covariance.
\end{proof}

\subsubsection{A marginal-sensitive covariance for the same law}
The probabilities in Lemma~\ref{r3:lem:dyadicprob} retain information which is lost by filling every possible active coordinate with Gaussian noise. To express it, put $D_p=\operatorname{diag}(p_j(1-p_j))$ and define
\begin{equation}\label{r3:eq:variancebudget}
 v_A(p)=v_0\sum_{k\ge1}4^{-k}\int_0^\infty
       \min\left\{1,\sum_{j:R_j^2>t}r_{j,k}\right\}\,dt.
\end{equation}
This scalar is determined by the prescribed means and the column radii. If $R_{(1)}^2\ge\cdots\ge R_{(n)}^2$ and $R_{(n+1)}=0$, its level-$k$ integral is
\[
 \sum_{\ell=1}^n(R_{(\ell)}^2-R_{(\ell+1)}^2)
           \min\left\{1,\sum_{j=1}^{\ell}r_{(j),k}\right\}.
\]
A truncation after $K$ levels has a nonnegative remainder at most
$v_0 R^2/(3\cdot4^K)$. Thus this quantity is computable to a specified absolute accuracy by sorting the radii once and evaluating finite sums.

\begin{theorem}\label{r3:thm:variance}
The law in Theorem~\ref{r3:thm:dyadic} can be chosen so that
\[
 \operatorname{Cov}(AE,E)\preceq
       \operatorname{diag}\bigl(v_A(p)I,wD_p\bigr).
\]
Consequently, on the nondeterministic coordinates and when $v_A(p)>0$,
\begin{equation}\label{r3:eq:marginalcov}
 \operatorname{Cov}(Z-p)\preceq
 Q_V:=\left((wD_p)^{-1}+v_A(p)^{-1}A^{\mathsf T}A\right)^{-1}.
\end{equation}
In particular,
\[
 \operatorname{Cov}(Z-p)\preceq wD_p,
 \qquad
 v_A(p)\le\frac{v_0 R^2}{3}.
\]
All these estimates hold for the same law as the hard bound and the joint-reference comparison.
\end{theorem}
\begin{proof}
Conditional on the past, the covariance of a level-$k$ joint increment is at most
\[
 4^{-k}\operatorname{diag}
       \bigl(v_0 R_{I_k}^2I,w\operatorname{diag}(\mathbf1_{j\in I_k})\bigr).
\]
Different martingale increments are orthogonal in $L^2$. The sum of the expected coefficient blocks is $wD_p$ by \eqref{r3:eq:dyadicenergy}. For the physical block, the layer-cake identity and a union bound give
\[
 \mathbb E R_{I_k}^2
 =\int_0^\infty\mathbb P\bigl(I_k\cap\{j:R_j^2>t\}\ne\varnothing\bigr)\,dt
 \le\int_0^\infty\min\left\{1,\sum_{j:R_j^2>t}r_{j,k}\right\}\,dt.
\]
This proves the joint covariance estimate at finite mesh. Its limit holds because the output is bounded and the initial mesh error tends uniformly to zero. Apply the second-moment calculation following Proposition~\ref{r3:thm:short} to obtain \eqref{r3:eq:marginalcov}.
\end{proof}

The integrand in \eqref{r3:eq:variancebudget} is the largest possible union probability compatible with the displayed marginal activation probabilities. For a fixed level, its bound on the expected maximum is attained by placing events of those probabilities successively around a circle, in decreasing order of the radii. This establishes sharpness of that intermediate bound among arbitrary couplings of the masks; it does not assert that every such coupling is realized by the rounding construction.

\begin{corollary}[A three-point integer window]\label{LAT-cor:integer}
For a nonzero matrix $A$ with unit Euclidean column bounds and any
$b\in\R^n$, there is an integer vector $N$ of mean $b$, fixing every
integral coordinate, such that
\begin{equation}\label{LAT-eq:integer}
 \|A(N-b)\|_\infty<R_A/2<C/2,\qquad |N_j-b_j|<3/2.
\end{equation}
Writing $J=\{j:b_j\notin\Z\}$, the same law satisfies
\begin{equation}\label{LAT-eq:integer-cx}
 (A(N-b),(N-b)_J)\cx\tfrac12(X_{R_A},T_J),\qquad X_{R_A}\perp T_J.
\end{equation}
Each active coordinate takes at most three integer values; at a
half-integer these are the two nearest integers.
\end{corollary}
\begin{proof}
Delete integral coordinates and use Corollary~\ref{thm:affine} with
physical reference $X_{R_A}/2$, coefficient reference $b+T_J/2$,
intervals $(b_j-3/2,b_j+3/2)$ and their integer alphabets of gap one.
Halving both references halves the overlap, so its height exceeds
$1/2$ by \eqref{eq:scalar-input}. The affine conditional means give
\eqref{LAT-eq:integer}--\eqref{LAT-eq:integer-cx}; the interval length
bounds the number of values. The zero active matrix has zero physical
error and can use independent scalar roundings. A lattice $g\Z$
similarly gives the physical radius $gR_A/2$.
\end{proof}


\subsection{Couplings with prescribed marginals}\label{sec:coupling}
\label{BAL-sec}
Let $P_j$ be probability measures on standard Borel spaces $\Omega_j$,
and let $\Phi_j:\Omega_j\to\mathcal H$ be strongly measurable maps
of essential diameter at most $d>0$. Here $\mathcal H$ is real and
separable. Put
\begin{equation}\label{eq:sampling-constants} A_0=\frac{v_0C^2}{3}<2.038,\qquad B_0=\frac{\alpha C^2}{3},\qquad
 b=2w<5.569,\qquad \beta=\frac1{2\kappa},\quad
 \kappa=\frac{\log2}{I_q(1)},
\end{equation}
where $I_q=(\log\E e^{(\cdot)T_1})^*$.
Appendix~\ref{BAL-app:scalar} gives $\beta=0.13421456\ldots$.
An isonormal process is an $L^2$-linear centered Gaussian family
$(W(u))_{u\in\mathcal H}$ satisfying
$\E[W(u)W(v)]=\langle u,v\rangle$.

\begin{theorem}[Coupling the marginals]\label{thm:coupling}
There is one coupling $\mu$ of the $P_j$, with
$S=\sum_j(\Phi_j(Z_j)-P_j\Phi_j)$, which has an extension carrying
an isonormal process such that
\begin{equation}\label{eq:coupling-Gaussian}
 \E[\sqrt{B_0}\,d\,W(u)\mid Z]=\langle u,S\rangle
 \qquad(u\in\mathcal H).
\end{equation}
For all $g_j\in L^2(P_j)$, and for all $u\in\mathcal H$ and
nonnegative integrable $f_j$, the same law satisfies
\begin{align}
 \Var_\mu\sum_jg_j(Z_j)
 &\le\inf_{u\in\mathcal H}\left\{A_0d^2\|u\|^2+
 b\sum_j\Var_{P_j}(g_j-\langle u,\Phi_j\rangle)\right\},
                  \label{eq:coupling-variance}\\
 \E_\mu e^{\langle u,S\rangle}\prod_j f_j(Z_j)^\beta
 &\le e^{A_0d^2\|u\|^2/2}\prod_j(P_jf_j)^\beta.
                  \label{eq:coupling-mixed}
\end{align}
In finite dimension it retains the bounded reference
\begin{equation}\label{eq:coupling-bounded}
 U_d=d\sum_{k\ge1}2^{-k}X_k,\qquad
 \E[U_d\mid Z]=S,\qquad \|S\|_\infty<dC,
 \quad X_k\stackrel{\rm iid}{\sim}f_{C,m}.
\end{equation}
In Hilbert space, any prescribed finite set of orthonormal coordinates
may retain $|\langle e_i,S\rangle|\le dC$.
For nonnegative weights $\pi_j$, $d$ is replaced by
$\max_j\pi_j\operatorname{diam}(\Phi_j)$ after scaling
$\Phi_j,g_j$ by $\pi_j$.
\end{theorem}

The two geometric series
\begin{equation}\label{eq:two-series}
 \sum_{k\ge1}2^{-k}=1,\qquad \sum_{k\ge1}4^{-k}=\frac13
\end{equation}
explain the bounded and quadratic scales. The coefficient reference
has an additional role: its conditional means control arbitrary
observables of the sample, including the residual in
\eqref{eq:coupling-variance}.

\subsubsection{Pairing and the conditional means}
First take finitely supported marginals. Their probabilities form
vectors $p_i$ in groups $G$, each of total one; the argument also
allows any integral total $r_G$. Associate a feature $a_i\in\R^m$
to each option, with $\|a_i-a_j\|_2\le d$ within a group, and let
$A=(a_i)_i$.
For dyadic $p$ of denominator $2^K$, visit levels $k=K,\ldots,1$.
The odd numerators in each group occur in even number. Pair them and
let $E_k$ have columns $e_i-e_j$ for the disjoint pairs. The update
\begin{equation}\label{eq:pair-update}
 x\longmapsto x+2^{-k}E_k\epsilon_k,
 \qquad B_k=d^{-1}AE_k
\end{equation}
uses the signing law of Theorem~\ref{MAIN-signing} for $B_k$.
Each possible update stays in $[0,1]$, preserves the group sums and
moves to the coarser grid. Its conditional mean is zero.
The terminal $Y$ has exact totals and $\E Y=p$.

At each history take a fresh joint source $(X_k,Q_k)$, with law
$f_{C,m}\otimes q^{\otimes s_k}$, and use its signing kernel. Future
steps use only the past signs. Hence the local source residuals have
zero conditional mean given the whole sign history, and therefore
\begin{equation}\label{eq:source-sums}
 \E[U_{d,K}\mid Y]=A(Y-p),\qquad
 \E[V_K\mid Y]=Y-p,
 \quad
 U_{d,K}=d\sum_{k=1}^K2^{-k}X_k,\quad
 V_K=\sum_{k=1}^K2^{-k}E_kQ_k.
\end{equation}
The physical sources $X_k$ are independent. Their support gives
$\|A(Y-p)\|_\infty<dC(1-2^{-K})$.
For real $p$, approximate by dyadic probabilities with the same group
totals and pass to a coupling limit. Since every coordinate of $U_d$
lies strictly inside $(-dC,dC)$, so does its conditional mean. This
retains the strict support bound in \eqref{eq:coupling-bounded}.

\subsubsection{The residual variance is a graph quadratic form}
The source covariance is determined by the matrix
\[
 L=\E\sum_{k=1}^K4^{-k}E_kE_k^{\mathsf T}.
\]
It is a graph Laplacian supported within the groups. The scalar
coordinate martingale starts at $p_i$ and ends in $\{0,1\}$, so
orthogonality of increments gives $L_{ii}=p_i(1-p_i)$.
Independence at each step and martingale orthogonality between steps
also give
\begin{equation}\label{eq:source-covariances}
 \Cov(U_{d,K})=\frac{d^2v_0C^2(1-4^{-K})}{3}I,\qquad
 \Cov(V_K)=wL,\qquad \Cov(U_{d,K},V_K)=0.
\end{equation}
Every edge lies within a group, so for arbitrary group constants $c_G$,
\[
 r^{\mathsf T}Lr\le2\sum_G\sum_{i\in G}p_i(1-p_i)(r_i-c_G)^2.
\]
For a group of total one, use $p_i(1-p_i)\le p_i$ and minimize over
$c_G$ to obtain $L_G\preceq2(\diag p_G-p_Gp_G^{\mathsf T})$.
Now \eqref{eq:source-sums} gives, for every $t,u$,
\begin{equation}\label{eq:paired-query}
 t^{\mathsf T}(Y-p)
 =\E[u^{\mathsf T}U_{d,K}+(t-A^{\mathsf T}u)^{\mathsf T}V_K\mid Y].
\end{equation}
Apply conditional variance and \eqref{eq:source-covariances}. If $t_i$
is the value of $g_G$ on option $i$, the result is exactly
\eqref{eq:coupling-variance}. This is the calculation
\eqref{eq:opening-regression} with the coefficient quadratic form
replaced by a graph Laplacian.

\subsubsection{The exponential inequality}
Write $\psi=\log\E e^{(\cdot)T_1}$ and
$\psi_X(u)=\log\E e^{u^{\mathsf T}X_C}$.
The scalar rate inequality $I_q(s)\ge\kappa^{-1}\mathfrak i(s)$ on
$[-1,1]$, proved in Appendix~\ref{BAL-app:scalar}, implies the local
bound
\begin{equation}\label{eq:local-product}
 \log\E e^{(B^{\mathsf T}u+r)^{\mathsf T}\epsilon}
 \le\psi_X(u)+\kappa^{-1}\sum_i\log\cosh(\kappa r_i).
\end{equation}
Indeed, conditional Jensen first gives $\psi_X(u)+\sum_i\psi(s_i)$;
clipping the coefficient test on $\epsilon_i\in\{-1,1\}$ adds
$\sum_i|r_i-s_i|$. Minimize over $s_i$ and take convex conjugates.

The scalar potential
\[
 b_{x,2\kappa}(r)=(2\kappa)^{-1}\log(1-x+xe^{2\kappa r})-xr
\]
vanishes when $x\in\{0,1\}$. Fix $t,u$ and write $r=t-A^{\mathsf T}u$.
For a paired step of size $a=2^{-k}$, set
$z_i=a(e^{2\kappa r_i}-1)/(1-x_i+x_ie^{2\kappa r_i})$.
Since both updated coordinates are feasible, $|z_i|<1$.
For an oriented pair $(i,j)$, the residual increment plus the potential
change equals
\[
 \frac{\log(1-z_i^2)+\log(1-z_j^2)}{4\kappa}
 +\frac{\operatorname{artanh}z_i-\operatorname{artanh}z_j}{2\kappa}
 \epsilon.
\]
Convexity of $\log\cosh$ gives
\[
 \kappa^{-1}\log\cosh\!\left(\frac{x-y}{2}\right)
 \le\frac{\log\cosh x+\log\cosh y}{2\kappa},
\]
which cancels the constant term when
$x=\operatorname{artanh}z_i$, $y=\operatorname{artanh}z_j$.
Thus \eqref{eq:local-product} bounds the conditional exponential of
$t^{\mathsf T}\Delta+\sum_i[b_{x_i+\Delta_i,2\kappa}(r_i)-b_{x_i,2\kappa}(r_i)]$
by $\exp\psi_X(dau)$. Telescoping to the terminal vertex gives
\begin{equation}\label{eq:quota-exponential}
 \log\E e^{t^{\mathsf T}(Y-p)}
 \le\log\E e^{u^{\mathsf T}U_d}
        +\sum_i b_{p_i,2\kappa}(t_i-a_i^{\mathsf T}u).
\end{equation}
The dyadic limit preserves this inequality.

For a categorical group and positive values $f_G(i)$, divide them by
$P_Gf_G$ and denote the result by $\widetilde f_G(i)$.
Use $t_i=a_i^{\mathsf T}u+\beta\log\widetilde f_G(i)$ in
\eqref{eq:quota-exponential}. The resulting group factor is
$\prod_{i\in G}(1-p_i+p_i\widetilde f_G(i))^\beta$.
Its bases have sum $|G|$, so their product is at most one.
This proves the mixed inequality with
$\log\E e^{u^{\mathsf T}U_d}$ in place of $A_0d^2\|u\|^2/2$.
Equation~\eqref{eq:cosine-laplace} and \eqref{eq:two-series} give
\eqref{eq:coupling-mixed}; zero values follow by a limit.

\subsubsection{A single limit for all observables}
For bounded finite-dimensional features, partition their ranges into
sets of diameter tending to zero and use their conditional means on
the cells. Couple the cell labels, then lift independently using the
original conditional marginals. Each marginal remains exactly $P_j$.
Conditional Jensen preserves the product inequality because
$0<\beta\le1$. Given the cell labels, the residuals
$g_j-\E[g_j\mid\text{cell}_j]$ are independent and centered, so their
variance is absorbed by $b\ge1$. The feature error tends uniformly to
zero. We retain all limiting observables in the following passage.

For separable $\mathcal H$, center the features, so their norms are
at most $d$, and project onto the first $r$ basis directions, always
including the prescribed hard coordinates. Use an increasing sequence
of finite partitions that resolves these coordinates and a countable
generating algebra of each $\Omega_j$.
To pass all tests to one law, embed each space by a Borel injection into
$[0,1]$, its full Hilbert feature, and the sequence of bounded rational
simple tests from this algebra. This is a Borel injection into a Polish
space, with Borel image. The pushed-forward marginals are fixed;
their couplings are tight. A subsequence converges and can be pulled
back to the original spaces.

Expectations involving a fixed finite number of bounded coordinate
tests pass to the limit by clipping continuous coordinate projections
at their bounds. Conditional finite-feature errors tend to zero in
$L^2$ of each fixed marginal. These facts pass
\eqref{eq:coupling-variance} and \eqref{eq:coupling-mixed} first to
simple tests and finite-support vectors. Density in marginal $L^2$
extends variance; bounded $L^1$ approximation, the inequality
$|x^\beta-y^\beta|\le|x-y|^\beta$, and truncation extend the product
inequality. Norm approximation extends the Hilbert vectors.
The selected hard coordinates pass as closed support conditions.
In finite dimension, keeping the bounded source itself retains the
strict bound by its conditional-mean identity.

It remains to retain the Gaussian process in this same limit.
For every finite projection, \eqref{eq:reference-gaussian} and
\eqref{eq:two-series} give
$U_d\cx N(0,B_0d^2I)$. Couple that Gaussian to the whole sample using
a martingale kernel conditional on $S$, and append an independent
Gaussian tail. The sample marginal and the infinite Gaussian-sequence
marginal are fixed. Include this sequence in the tightness argument.
Testing against bounded functions of the augmented sample and using
uniform integrability of each Gaussian coordinate passes all finite
conditional means to the limit. A monotone-class argument gives them
for the full sample sigma-field. If $(G_i)$ is the limiting independent
Gaussian sequence, put $W(u)=\sum_i\langle u,e_i\rangle G_i$ in
$L^2$. Conditional expectation is an $L^2$ contraction, so this proves
\eqref{eq:coupling-Gaussian} for every $u$. This completes the proof of
Theorem~\ref{thm:coupling}.

\subsubsection{Sampling and the covariance operator}\label{sec:hilbert}
For identical marginals and maps, average the law over sample
permutations. The sum is unchanged, and the bounds are invariant under
this averaging. Taking $g_j=g$ in \eqref{eq:coupling-variance} and
dividing by $N^2$ proves Theorem~\ref{VB-thm:sampling}.
If $C_\Phi=\Cov_P(\Phi)$, the substitution
$g=\langle t,\Phi\rangle$ and the minimization
\begin{equation}\label{eq:Hilbert-minimum}
 \inf_u\{\lambda\|u\|^2+\langle t-u,C_\Phi(t-u)\rangle\}
 =\lambda\langle t,C_\Phi(C_\Phi+\lambda I)^{-1}t\rangle
\end{equation}
give
\begin{equation}\label{eq:kernel-covariance}
 \Cov(\widehat P_N\Phi-P\Phi)
 \preceq\frac{A_0d^2}{N^2}C_\Phi(C_\Phi+\lambda_N I)^{-1},
 \qquad\lambda_N=\frac{A_0d^2}{bN}.
\end{equation}
Thus the same least-squares calculation as
\eqref{eq:quadratic-infimum} acts on the feature covariance.
For a reproducing-kernel Hilbert space with canonical feature map,
\begin{align}
 \sup_{\|f\|_{\mathcal H}\le1}\E|\widehat P_Nf-Pf|^2
 &\le\frac{A_0d^2}{N^2},\qquad
 \log\E e^{t(\widehat P_Nf-Pf)}
 \le\frac{A_0d^2t^2\|f\|_{\mathcal H}^2}{2N^2},
                 \label{eq:kernel-errors}\\
 \E\MMD_k(\widehat P_N,P)^2
 &\le\frac{A_0d^2}{N^2}
 \tr[C_\Phi(C_\Phi+\lambda_NI)^{-1}].\label{eq:MMD}
\end{align}
The last line takes the supremum inside the expectation; the
trace retains the effective dimension.

\subsubsection{Conditioning and exact quotas}
The variational formula for relative entropy applied to
\eqref{eq:coupling-mixed} gives, for any alternative law $\nu$,
\begin{equation}\label{eq:coupling-entropy}
 D(\nu\Vert\mu)\ge
 \frac{\|\E_\nu S\|^2}{2A_0d^2}
 +\beta\sum_jD(\nu_j\Vert P_j).
\end{equation}
Optimize independently over $u$ and $\log f_j$ to obtain the formula.
For an observation $W$ of an identically distributed sample, and a
measurably selected unit-norm kernel function $f_W$, it follows that
\begin{equation}\label{eq:adaptive-kernel}
 \big|\E[\widehat P_Nf_W-Pf_W]\big|
 \le\frac{\sqrt{2A_0}\,d}{N}
 \sqrt{I(W;Z)-\beta\sum_j I(W;Z_j)}.
\end{equation}
Apply \eqref{eq:coupling-entropy} to the conditional law given $W$,
then use Cauchy--Schwarz. This is the information-selection estimate
of \cite{RussoZou} with the marginal-information term retained.

The finite construction above also proves unbiased binary rounding
with prescribed integral group totals: \eqref{eq:source-sums} and
\eqref{eq:quota-exponential} hold for every such input, with all frozen
coordinates fixed. Dyadic depth $K$ uses exactly $K$ conditional signing
calls on disjoint pairs. General probabilities and the Hilbert-space
laws use the limiting argument just given; these existence statements
do not assert a polynomial-time sampler for that limit.


\subsection{Local exactness and Sobolev quadrature}\label{sec:quadrature}
Let $\Omega=[0,1]^D$ with uniform measure $P$. Use the restriction
norm from the Bessel-potential space
\[
 \|F\|_{H^s(\R^D)}^2=\int(1+|\xi|^2)^s|\widehat F(\xi)|^2\,d\xi,
 \qquad
 \|f\|_{H^s(\Omega)}=\inf_{F|_\Omega=f}\|F\|_{H^s(\R^D)},
\]
with unitary Fourier transform. Continuous representatives are used
when $s>D/2$. Write $Q_Nf=N^{-1}\sum_{i=1}^Nf(Z_i)$ and
$E_Nf=Q_Nf-Pf$.

\begin{theorem}[Quadrature with hard feature constraints]\label{thm:quadrature}
Fix $s_0>D/2$ and a bounded measurable
$\Psi:\Omega\to\R^m$ of essential Euclidean diameter at most one.
For every integer $N\ge N_0(D,s_0)$, there is one exchangeable law
of exactly $N$ points, each uniform on $\Omega$, such that
\begin{equation}\label{eq:quad-hard-feature}
 \|Q_N\Psi-P\Psi\|_\infty\le L_{D,s_0}/N
 \quad\text{almost surely}.
\end{equation}
For each $D/2<s\le s_0$ and every finite family $f_1,\ldots,f_k\in H^s$,
\begin{equation}\label{eq:quad-Gaussian}
 ((E_Nf_i)_{i=1}^k,E_N\Psi)\cx(G_f,H_\Psi),
\end{equation}
where the two Gaussian blocks are independent and
\[
 \Cov(G_f)=L_{D,s,s_0}N^{-1-2s/D}
                 (\langle f_i,f_j\rangle_{H^s})_{i,j},\qquad
 \Cov(H_\Psi)=L_{D,s_0}N^{-2}I_m.
\]
The same law satisfies
\begin{align}
 \sup_{\|f\|_{H^s}\le1}\E|E_Nf|^2
 &\le L_{D,s,s_0}N^{-1-2s/D} &&(0\le s\le s_0),
                         \label{eq:quad-mse}\\
 \sup_{\|f\|_{H^s}\le1}|E_Nf|
 &\le L_{D,s,s_0}N^{-s/D} &&(D/2<s\le s_0)
                         \label{eq:quad-hard}
\end{align}
on every outcome for the second line. It integrates every polynomial
of degree at most $r=\lceil s_0\rceil-1$ exactly.
All constants are independent of $m,\Psi$. The law is chosen before
$s$, the integrands and the convex test.
\end{theorem}

A measurable partition $B_1,\ldots,B_m$ can be encoded by
$\Psi=2^{-1/2}(\mathbf1_{B_1},\ldots,\mathbf1_{B_m})$.
Then \eqref{eq:quad-hard-feature} bounds the integer count error in
every part by a constant independent of the number of parts, under the
same law as \eqref{eq:quad-Gaussian}.

\subsubsection{Randomizing a polynomial formula}
A uniformly rotated tetrahedron projected onto a fixed axis gives four
uniform points on $[-1,1]$ with average zero and average square $1/3$
at every rotation. For higher moments, use a degree-$r$ polynomial
design on $(S^2)^D$, whose dimension is $2D$. The theorem of
Etayo, Marzo and Ortega-Cerd\`a \cite{EMO} constructs such a design with
every cardinality $q\ge K(D,r)$, with $K(D,r)\le C_Dr^{2D}$ for
$r\ge1$. For $r=0$ take $K=1$.

Write its points as $v_i=(v_{i1},\ldots,v_{iD})$ and choose independent
Haar rotations $R_1,\ldots,R_D\in SO(3)$. Set
\[
 X_i=(\langle e_3,R_1v_{i1}\rangle,\ldots,
                         \langle e_3,R_Dv_{iD}\rangle).
\]
A spherical zone between heights $a,b$ has area $2\pi(b-a)$; hence
every $X_i$ is exactly uniform on $[-1,1]^D$. For each fixed rotation
array, composition with these coordinate projections preserves
polynomial degree. The design identity therefore integrates every
polynomial of degree at most $r$ exactly. Affine maps transport the
construction to rectangles.

Write $N=\sum_jq_j$ with $K\le q_j<2K$ and partition the cube into
rectangles $B_j$ of volumes $\pi_j=q_j/N$ and aspect ratio at most
five. Such a partition is obtained by dividing the volume list into
two almost equal cardinality lists and repeatedly splitting the
longest side in their volume ratio. The comparable requested volumes
keep each splitting fraction in $[1/5,4/5]$.
For $h=N^{-1/D}$ we have $|B_j|\asymp h^D$ and
$\operatorname{diam}(B_j)\asymp h$.
Let $Q_j^{\omega_j}$ be the local polynomial formula and $\rho_j$
its product-Haar law. Each point has uniform marginal $P_j$ on $B_j$,
so
\begin{equation}\label{eq:local-formula}
 \E_{\rho_j}Q_j^{\omega_j}g=P_jg,\qquad
 \Var_{\rho_j}(Q_j^{\omega_j}g)\le\Var_{P_j}(g).
\end{equation}
The variance bound is Jensen applied to the average of
$(g-P_jg)^2$. We will couple the arrays $\omega_j$, keeping every
$\rho_j$ fixed.

\subsubsection{One Hilbert sum contains all the smoothness classes}
The coupling theorem accepts orthogonal sums of feature spaces.
We therefore include ordinary Sobolev kernels at the orders
\[
 s_\ell=D/2+(s_0-D/2)2^{-\ell},\qquad \ell\ge0,
\]
with summable weights. Put
\begin{equation}\label{eq:sobolev-sum-kernel}
 k_{s,h}(\xi)=h^{D-2s}(1+|\xi|^2)^{-s},\qquad
 \lambda_h(\xi)=\sum_{\ell\ge0}4^{-\ell}k_{s_\ell,h}(\xi).
\end{equation}
For fixed $h>0$, this is positive and integrable. It defines a
translation-invariant kernel and a Hilbert space $\mathcal H_h$ with
norm $\int|\widehat F|^2/\lambda_h$, restricted to $\Omega$.
The direct-sum interpretation is the exact formula
\begin{equation}\label{eq:sobolev-decomposition}
 \|F\|_{\mathcal H_h(\R^D)}^2
 =\inf_{F=\sum_{\ell\ge0}F_\ell}
       \sum_{\ell\ge0}4^\ell h^{2s_\ell-D}
                                  \|F_\ell\|_{H^{s_\ell}}^2.
\end{equation}
Pointwise Cauchy--Schwarz in Fourier coordinates proves the inequality;
equality is obtained with
$\widehat F_\ell=4^{-\ell}k_{s_\ell,h}\widehat F/\lambda_h$.
Thus the norm again minimizes a sum of quadratic costs over all
representations of the same observable.

\begin{lemma}\label{lem:sobolev-sum}
The averaged feature map $\omega_j\mapsto Q_j^{\omega_j}\Phi_h$ has
diameter bounded by $L_{D,s_0}$ in $\mathcal H_h$, uniformly in $h,j$.
For every $D/2<s\le s_0$,
\begin{equation}\label{eq:sobolev-sum-embedding}
 \|f\|_{\mathcal H_h}^2\le L_{D,s,s_0}h^{2s-D}\|f\|_{H^s}^2.
\end{equation}
For every $0\le s\le s_0$ there is $v\in\mathcal H_h$ with
\begin{equation}\label{eq:sobolev-approximation}
 \|v\|_{\mathcal H_h}^2\le Lh^{2s-D}\|f\|_{H^s}^2,
 \qquad \|f-v\|_{L^2(P)}^2\le Lh^{2s}\|f\|_{H^s}^2.
\end{equation}
\end{lemma}
\begin{proof}
Two local formulas agree on all moments through degree $r$. Taylor's
formula for $e^{i\xi\cdot x}$ about the center of their rectangle gives
\[
 |(Q_j^\omega-Q_j^{\omega'})e^{i\xi\cdot x}|
 \le L\min\{1,(h|\xi|)^{r+1}\}.
\]
For $D/2<s\le s_0$, splitting the radial integral at $1$ and $h^{-1}$
shows
\begin{equation}\label{eq:sobolev-local-integral}
 \int k_{s,h}(\xi)\min\{1,(h|\xi|)^{2r+2}\}\,d\xi
 \le L_{D,s_0}\left(1+\frac1{2s-D}\right).
\end{equation}
The upper tail is $|S^{D-1}|/(2s-D)$. The lower integral is bounded
because $D+2r+2-2s_0>0$. Sum \eqref{eq:sobolev-local-integral}
with weights $4^{-\ell}$; since $2s_\ell-D=(2s_0-D)2^{-\ell}$,
the series converges. This proves the diameter bound. The same tail
calculation without the cancellation factor proves integrability of
$\lambda_h$ for each fixed $h$.

If $s=\theta s_\ell+(1-\theta)s_{\ell+1}$, then
$k_{s,h}=k_{s_\ell,h}^{\theta}k_{s_{\ell+1},h}^{1-\theta}$.
The weighted arithmetic--geometric mean inequality gives
\begin{equation}\label{eq:kernel-interpolation}
 \lambda_h\ge4^{-\ell}k_{s_\ell,h}
              +4^{-(\ell+1)}k_{s_{\ell+1},h}
 \ge4^{-(\ell+1-\theta)}k_{s,h}.
\end{equation}
Integration and the restriction norm prove
\eqref{eq:sobolev-sum-embedding}. Every $s>D/2$ lies between two of
the displayed orders. The constants may diverge as $s\downarrow D/2$.

For \eqref{eq:sobolev-approximation}, take an $H^s$ extension of $f$
with norm at most twice its restriction norm, and keep its Fourier
frequencies $|\xi|\le h^{-1}$. On this ball, the term $\ell=0$ gives
\[
 \lambda_h^{-1}\le h^{2s_0-D}(1+|\xi|^2)^{s_0}
 \le Lh^{2s-D}(1+|\xi|^2)^s.
\]
The norm estimate follows, and Parseval bounds the discarded tail by
$Lh^{2s}\|f\|_{H^s}^2$.
\end{proof}

\subsubsection{Coupling the local formulas}
Apply Theorem~\ref{thm:coupling} to the marginals $\rho_j$ and maps
\[
 \omega_j\longmapsto\pi_j
 \bigl(Q_j^{\omega_j}\Phi_h,\ Q_j^{\omega_j}\Psi\bigr)
 \quad\text{in }\mathcal H_h\oplus\R^m.
\]
Their diameter is at most $L/N$: the first component is controlled
by Lemma~\ref{lem:sobolev-sum}, and the second lies in the convex
hull of a set of diameter one. The finite coordinates give
\eqref{eq:quad-hard-feature}. The Gaussian processes on the two
orthogonal summands are independent. The bounded inclusion
$H^s\to\mathcal H_h$ in \eqref{eq:sobolev-sum-embedding} compares
all finite Gram matrices, so adding independent Gaussian noise gives
\eqref{eq:quad-Gaussian}. The scale is
\begin{equation}\label{eq:rate-factorization}
 N^{-2}h^{2s-D}=N^{-1-2s/D}.
\end{equation}

The residual-variance bound retains the exact moments. Let
$\mathcal P_r(\mathcal B)$ denote functions polynomial of degree at
most $r$ separately on the $B_j$. For any $g\in L^2(P)$,
\begin{equation}\label{eq:quad-residual}
 \Var(Q_Ng)\le L\inf_{v\in\mathcal H_h,\ p\in\mathcal P_r(\mathcal B)}
 \left\{\frac{\|v\|_{\mathcal H_h}^2}{N^2}
       +\frac{\|g-v-p\|_{L^2(P)}^2}{N}\right\}.
\end{equation}
Indeed, local exactness removes $p$ from each centered error, and
\eqref{eq:local-formula} bounds its residual variance by
$\sum_j b\pi_j^2\Var_{P_j}(g-v-p)\le LN^{-1}\|g-v-p\|_2^2$.
Use \eqref{eq:sobolev-approximation} in \eqref{eq:quad-residual} to
obtain \eqref{eq:quad-mse}, including $s\le D/2$.
Taking $v=0$ and optimizing the constant polynomial gives the general
bound $\Var(Q_Ng)\le L\Var_P(g)/N$.

There are $q_j$ output points from block $j$ and $\pi_j=q_j/N$,
so the global formula has exactly $N$ equal positive weights. A uniform
random permutation makes its points exchangeable and each marginal
is $\sum_j(q_j/N)P_j=P$. All formula values and bounds are unchanged.
Exact marginals also make the rule unbiased for all integrable $g$.

For the hard Sobolev bound, local polynomial approximation gives
\begin{equation}\label{eq:local-Sobolev}
 \inf_{p\in\mathcal P_{\lceil s\rceil-1}}
 \|f-p\|_{L^\infty(B_j)}
 \le L_{D,s}h^{s-D/2}|f|_{H^s(B_j)}\qquad(s>D/2).
\end{equation}
On a unit rectangle, subtract the $L^2$ projection onto the indicated
polynomials. The higher-order Poincar\'e inequality controls the
remaining Sobolev norm by its order-$s$ seminorm, and Sobolev embedding
gives its supremum norm. Rescaling proves \eqref{eq:local-Sobolev};
the aspect ratios bound its constants. For noninteger $s=k+\theta$
use the Slobodeckij seminorm of derivatives of order $k$, whose
nullspace is the polynomials of degree at most $k$. For integer $s=k$
the nullspace has degree at most $k-1$. This is the
Bramble--Hilbert approximation argument \cite{BrambleHilbert}.
The squared local seminorms sum to at most $L\|f\|_{H^s(\Omega)}^2$,
since their derivative integrals or nonnegative fractional double
integrals are bounded by those of a whole-space extension.
Local exactness and Cauchy--Schwarz now give, at every outcome,
\[
 |E_Nf|\le Lh^{s-D/2}\sum_j|B_j||f|_{H^s(B_j)}
 \le Lh^s\|f\|_{H^s(\Omega)}.
\]
This proves \eqref{eq:quad-hard} and completes the proof of
Theorem~\ref{thm:quadrature}.

\subsubsection{Confidence and sharpness}
For a randomized rule $Q$, let $e(Q,H^s,\delta)$ be the least error
threshold exceeded with probability at most $\delta$, uniformly on the
unit ball. Let $e^{\rm ran}(N,H^s,\delta)$ be its infimum over all
algorithms using at most $N$ function values, allowing adaptivity and
nonlinearity.
\begin{corollary}\label{cor:confidence}
For $D/2<s\le s_0$ and $0<\delta<1/4$, the law just constructed satisfies
\begin{equation}\label{eq:optimal-confidence}
 e(Q_N,H^s,\delta)\asymp e^{\rm ran}(N,H^s,\delta)
 \asymp N^{-s/D}\min\left\{1,
                    \sqrt{\frac{\log(1/\delta)}N}\right\}.
\end{equation}
The constants may depend on $D,s,s_0$. The exponent in
\eqref{eq:quad-mse} is optimal for every $s\ge0$, and the hard
worst-case exponent in \eqref{eq:quad-hard} is optimal for $s>D/2$.
\end{corollary}
\begin{proof}
The Gaussian comparison gives Chernoff tails at variance scale
$N^{-1-2s/D}$. Combining them with \eqref{eq:quad-hard} proves the
upper bound without changing the law with $\delta$.

For completeness, the standard disjoint-bump lower bound gives all the
real orders here. Partition the cube into $M=q^D$ small cubes, where
$q=\lceil(4N)^{1/D}\rceil$ and $h_*=q^{-1}$. Put one scaled copy
$\varphi_j$ of a fixed positive smooth bump strictly inside each cube.
For every $s\ge0$,
\[
 \Big\|\sum_jc_j\varphi_j\Big\|_{H^s}
 \le Lh_*^{D/2-s}\|(c_j)\|_2.
\]
Disjoint derivative supports prove this at integer orders; interpolation
of the whole-space $L^2$ and $H^k$ estimates proves real orders.
Thus $f_\epsilon=a h_*^s\sum_{j=1}^M\epsilon_j\varphi_j$ has norm
at most one for fixed small $a>0$. Each sign contributes
$b_*=c h_*^{s+D}\asymp N^{-1-s/D}$ to its integral.
Under independent fair signs, a transcript of at most $N$ function
values leaves at least $M-N\ge3M/4$ signs independent and unobserved.
The conditional error is a shifted sum of these signs times $b_*$.

The binomial mass estimate, summed over $\lfloor\sqrt{k}\rfloor$
consecutive values, gives a tail at least
$c\exp[-C t^2/k]$ for $\sqrt{k}\lesssim t\le ck$.
The central mass bound gives a fixed positive tail below a small
multiple of $\sqrt{k}$, uniformly in the shift; the all-equal event
has probability $2^{-k}$ and gives the saturated lower bound of order $k$.
These are the binomial estimates used in
\cite[Lemma~A.1]{KunschRudolf}. They prove the lower order in
\eqref{eq:optimal-confidence}, including $\log(1/\delta)\gtrsim N$.
The remaining conditional variance is at least $c b_*^2M$, proving
the mean-square lower bound. For a fixed set of $N$ points, take all
unhit bumps with positive sign. Their integral is at least
$cN^{-s/D}$, while the rule sees only zero. This proves the hard
worst-case lower bound.
\end{proof}

If $W$ is a standard Borel observation of the sample and $f_W$ is a
measurably chosen unit-norm $H^s$ function, relative-entropy duality and
the Gaussian exponential bound give
\[
 |\E E_Nf_W|\le L N^{-s/D}
       \min\{1,\sqrt{I(W;Z)/N}\}.
\]
Under the product law of $W$ and $Z$, the logarithmic moment generating
function is at most $Lt^2N^{-1-2s/D}/2$. The variational inequality
and optimization in $t$ give the second term, and
\eqref{eq:quad-hard} gives the first. Selection from a fixed family of
$M$ functions has $I(W;Z)\le\log M$.

The theorem also applies to periodic Sobolev functions by multiplying
a periodic extension by a fixed smooth cutoff. A smooth bi-Lipschitz
change of variables, with bounded derivatives and inverse derivatives
through order $\lceil s_0\rceil+1$, transfers it to the corresponding
pushforward measure. The composition operator is bounded in the
stated Sobolev range, so the Gaussian Gram bounds and the errors
transfer together. Exact polynomials become the pullbacks of those
polynomials. These changes preserve the number and weights of points.


\section{Entire laws of matrix partitions and sparse approximations}
\label{R11-sec:sparse}
Akbas--Sra's reciprocal/Fisher estimates give the matrix reference~\cite{AkbasSraBSB,AkbasSraMatrix}; retaining its independent coefficient block adds exact selection probabilities and joint convex-order and entropy guarantees. The rank-one partition existence input is the MSS resolution of Weaver's Kadison--Singer formulation~\cite{MSS}.

Sparse spectral approximation is developed in~\cite{FC-BSS,FC-SHS,FC-BKLM}. Here every original summand must have expected weight one. Conditioning independent thinning on success would change these means, so each conditional partition is selected with its means already correct. Information charges are summable over the full thinning history, and repeated copies give a Poisson count reference. The continuous version produces hard norm discretizations with exact expected measure in the setting of~\cite{FC-CM,FC-CD}. Finite algorithms and limiting existence retain their separate input models.

\subsection{Matrix partitions with near-endpoint reference laws}
\label{R3-sec:matrix}
For an equal-norm Parseval frame $u_1,\ldots,u_n\in\R^d$,
Corollary~\ref{R3-cor:Weaver} constructs a fair half-selection $N$
with
\[
 \left\|\sum_iN_i u_iu_i^{\mathsf T}-\tfrac12I_d\right\|
   =O\!\left(\sqrt{\frac d{n\varepsilon}}\right),\qquad
 n\log2-H(N)=O(d\sqrt\varepsilon),
\]
under the near-endpoint Gaussian comparison for $N-\tfrac12\mathbf1$.
Here $\sum_i u_iu_i^{\mathsf T}=I_d$ and $\|u_i\|^2=d/n$.
The entropy budget depends on the matrix dimension even when $n/d$
is large. The proof obtains this dependence by combining a matrix
Fisher bound with the squared-energy replica argument.

Matrix norms in this section are operator norms. We use the following
analytic input: there are universal constants $\eta,L>0$ such that, if
$\|\sum_i B_i^2\|\le\eta/t^2$, the density
\[
 p_t(x)\propto\exp\left[-\|x\|_2^2/(2t^2)
 -32\Phi\left(\sum_i x_iB_i\right)
 +28\tr\left(\sum_i x_iB_i\right)^2\right],
 \]
on $\|\sum_i x_iB_i\|<1$, with $\Phi(H)=-\log\det(I-H^2)$,
has Fisher information matrix bounded by
\begin{equation}\label{R3-eq:external-matrix-Fisher}
 J(p_t)\preceq t^{-2}I+L(\tr B_iB_j)_{i,j}.
\end{equation}
This is the reciprocal/Fisher estimate of Akbas and Sra
\cite[Lemmas 2.4--2.5 and 3.3]{AkbasSraBSB}, derived there from
\cite[Proposition 3.8]{AkbasSraMatrix}. The constants below retain
this dependence explicitly. Its extension by zero is $C^2$ with finite Fisher information; Section~\ref{R27-sec:volume} uses this regularity to turn the same matrix bound into a support-volume estimate.

Use the truncated Gaussian parameters of
\eqref{R3-eq:v-root}--\eqref{R3-eq:new-radius}, and set
\begin{equation}\label{R3-eq:new-matrix-params}
 h_\varepsilon=\frac{\varepsilon}{v_\varepsilon(1+\varepsilon)},
 \quad L_* =\max\{L,\pi^2\},\quad
 v_* =\frac{\eta h_\varepsilon}{8},\quad
 q_* =\frac{h_\varepsilon}{8L_*},\quad
 t_* =\sqrt{2/h_\varepsilon}.
\end{equation}
These symbols $v_*,q_*$ are matrix variance and energy thresholds;
$q_\varepsilon$ continues to denote the scalar probability density.
A directional Fisher bound at most $h_\varepsilon$ gives height at
least $1+\varepsilon/2$ by inequality~\eqref{R3-eq:Gaussian-height}.

Fix a basis and write each real symmetric matrix as $A_i=D_i+R_i$,
where $D_i$ is diagonal and $R_i$ has zero diagonal. Let $\Gamma$
have columns $\diag D_i$, and put
\[
 \nu=\left\|\sum_i R_i^2\right\|,
 \quad w_i=\tr A_i^2,\quad q_{\max}=\max_i w_i.
\]
Let $Z_\circ$ be the product cosine law of radius one in dimension $d$,
let $I_\circ$ be its Cram\'er transform, and put
$b_\circ=\sqrt{\pi/2}(1/2-2/\pi^2)$.

\begin{theorem}[Matrix reference laws with inverse-square-root loss]
\label{R3-thm:matrix}
Suppose $\nu\le v_*r^2/4$ and $q_{\max}\le q_*r^2/4$, with $r>0$.
Then one symmetric signing law has
\begin{gather}
 \left\|\sum_i\sigma_iA_i\right\|<r/2,
 \qquad \E[Z_\circ\mid\sigma]=-4\Gamma\sigma/r,
 \qquad \E[T\mid\sigma]=\sigma,
 \label{R3-eq:matrix-source}\\
 (Z_\circ,T)\sim\law(Z_\circ)\otimes q_\varepsilon^{\otimes n}.
 \notag
\end{gather}
For $B=q_*r^2/4$, the same law satisfies
\begin{align}
 D(\mu\Vert U_n)+\E I_\circ(4\Gamma\sigma/r)
    &\le\delta_{q_\varepsilon} F_B(w),
       \label{R3-eq:matrix-energy}\\
 \sigma&\cx N\!\left(0,
   \left[\frac{I_n}{\kappa^2(1+\varepsilon)^2}
             +\frac{16}{r^2b_\circ^2}\Gamma^{\mathsf T}\Gamma
   \right]^{-1}\right).
       \label{R3-eq:matrix-Gaussian}
\end{align}
\end{theorem}
\begin{proof}
At base scale assume $\nu\le v_*$ and $q_{\max}\le q_*$.
Apply \eqref{R3-eq:external-matrix-Fisher} to $2R_i$ and $t_*$:
$4\nu\le\eta h_\varepsilon/2=\eta/t_*^2$.
For independent $X\sim p_{t_*}$ and $Z_\circ$, use the sheared
physical source $(X,2\Gamma X+Z_\circ)$ in directions $(e_i,0)$.
The directional Fisher information is at most
\[
 h_\varepsilon/2+4L\|R_i\|_F^2+4\pi^2\|D_i\|_F^2
       \le h_\varepsilon/2+4L_*w_i\le h_\varepsilon.
\]
Independence makes the cross score mean zero. The quadratic height
bound and Theorem~\ref{thm:affine} give conditional means
$(\sigma,0)$ for the sheared source and $\sigma$ for the auxiliary.
The open support gives
$\|\sum_i\sigma_iR_i\|<1/2$ and $\|\Gamma\sigma\|_\infty<1/2$.

For the stated scale, use the exact replica packing of $w_i$ with
capacity $B$. Scale each signed direct-sum matrix by $2/r$.
Each column has squared Frobenius norm at most $q_*$, and the
variance of its off-diagonal square sum is at most $v_*$.
The preceding base construction therefore applies to every mask
array. Its barrier-density block can depend on the masks; its cosine
and auxiliary product marginal is fixed and independent of them.
Selecting a replica and applying Corollary~\ref{REP-cor:amplification}
gives \eqref{R3-eq:matrix-source} and \eqref{R3-eq:matrix-energy}.
The scalar Gaussian comparisons and the precision left inverse give
\eqref{R3-eq:matrix-Gaussian}.
\end{proof}

\begin{corollary}[Positive semidefinite partitions]
\label{R3-cor:Weaver}
Suppose $A_i\succeq0$, $\sum_iA_i=I_d$ and
$\tr A_i\le\theta\le1$. Define
\begin{equation}\label{R3-eq:Kepsilon}
 K_\varepsilon^2=\max\{16/v_*,4/q_*\}.
\end{equation}
There is a law on $N\in\{0,1\}^n$ such that
\begin{gather}
 \E N_i=1/2,\qquad
 \left\|\sum_i N_iA_i-\tfrac12I_d\right\|
              <\frac{K_\varepsilon}{4}\sqrt\theta,
       \label{R3-eq:Weaver-hard}\\
 N-\tfrac12\mathbf1\cx
             (1+\varepsilon)\sqrt{\pi/8}\,G_n,
 \qquad D(\law(N)\Vert\operatorname{Bernoulli}(1/2)^{\otimes n})
          \le\frac{8\delta_{q_\varepsilon}}
                       {q_*K_\varepsilon^2}\,d.
       \label{R3-eq:Weaver-cx}
\end{gather}
It retains the independent product source in
\eqref{R3-eq:matrix-source} with $\sigma=2N-\mathbf1$ and
$r=K_\varepsilon\sqrt\theta$.
As $\varepsilon\downarrow0$,
\[
 K_\varepsilon=O(\varepsilon^{-1/2}),\qquad
 \frac{8\delta_{q_\varepsilon}}{q_*K_\varepsilon^2}
       =O(\sqrt\varepsilon).
\]
The constants depend only on the external $\eta,L$.
\end{corollary}
\begin{proof}
Since $A_i^2\preceq\theta A_i$, one has
$\sum_i A_i^2\preceq\theta I_d$, $\sum_iw_i\le\theta d$, and
$q_{\max}\le\theta^2\le\theta$.
Also $\sum_iD_i^2\preceq\theta I_d$ and
$R_i^2\preceq2(A_i^2+D_i^2)$, so $\nu\le4\theta$.
The definition of $K_\varepsilon$ ensures the hypotheses of
Theorem~\ref{R3-thm:matrix}. Put $N=(\sigma+\mathbf1)/2$ and use
$F_B(w)\le2\sum_iw_i/B$.
Finally $v_\varepsilon\to v_{\mathrm{tr}}$,
$h_\varepsilon\sim\varepsilon/v_{\mathrm{tr}}$, and
$\delta_{q_\varepsilon}=O(\sqrt\varepsilon)$.
\end{proof}
For rank-one summands, this is the Kadison--Singer/Weaver partition
setting \cite{MSS}. The simultaneous Gaussian and information
conclusions come from the retained reference. The dependence on the distance from the sharp Gaussian endpoint is
$\varepsilon^{-1/2}$, while the entropy deficit depends on the
ambient matrix dimension $d$ rather than the number $n$ of summands.

\subsubsection{Arbitrary starting points and matrix partial coloring}
A partial coloring must start from a prescribed vector. Symmetric signing of every subfamily lets us round that vector dyadically while preserving its mean. The following lemma records both the hard and distributional costs of this passage.

\begin{lemma}[Dyadic rounding from hereditary Gaussian signing]
\label{R19-lem:dyadic-start}
Let $B_1,\ldots,B_m$ belong to a normed real vector space. Suppose that for every $J\subset\{1,\ldots,m\}$ there is a symmetric signing law $\epsilon^J$ with
\[
 \left\|\sum_{i\in J}\epsilon_i^J B_i\right\|\le h,
 \qquad \epsilon^J\cx\tau G_J,
\]
where $G_J$ is standard Gaussian on the coordinate subspace indexed by $J$. For every $y\in[-1,1]^m$ there is a law on $X\in\{-1,1\}^m$ satisfying
\[
 \E X=y,\qquad
 \left\|\sum_i(X_i-y_i)B_i\right\|\le2h,
 \qquad X-y\cx N(0,\tfrac43\tau^2I_m).
\]
\end{lemma}
\begin{proof}
First assume $z=(y+\mathbf1)/2$ has dyadic denominator $2^N$. At stage $k=N,N-1,\ldots,1$, let $J$ be the coordinates whose numerator at denominator $2^k$ is odd. Add $2^{-k}\epsilon^J$ on those coordinates. Each coordinate remains in $[0,1]$, the new vector has denominator $2^{k-1}$, and the increment has conditional mean zero. Its feature norm is at most $2^{-k}h$. The final $Z\in\{0,1\}^m$ thus has mean $z$ and feature displacement at most $h\sum_{k=1}^N2^{-k}<h$.

For $X=2Z-\mathbf1$, the stage increment is conditionally dominated by $2^{1-k}\tau G_J$, hence by $2^{1-k}\tau G_m$ after independent Gaussian completion on the missing coordinates. Iterated conditional convex comparison, or the convolution rule of Section~\ref{R11-sec:composition}, gives covariance at most
$\tau^2\sum_{k=1}^N4^{1-k}I_m\preceq\tfrac43\tau^2I_m$.
The feature bound is $2h$. For general $y$, use dyadic approximations and pass to a subsequence of laws on the finite cube. The exact means, closed hard bound and each finite convex comparison pass to the limit.
\end{proof}

The hard arbitrary-start conclusion in Sachdeva--Thudi--Zhao's Conjecture~21 already follows from Akbas--Sra's full-coloring estimate and hereditary dyadic rounding. The following theorem performs that passage with the prescribed mean and the Gaussian conditional comparison retained at every step. The analytic matrix estimate remains the Akbas--Sra input.

\begin{theorem}[Sachdeva--Thudi--Zhao matrix partial coloring with exact mean]
\label{R19-thm:matrix-partial}
Fix $\varepsilon>0$ and $K_\varepsilon$ from \eqref{R3-eq:Kepsilon}. Let symmetric $d\times d$ matrices $A_1,\ldots,A_m$ satisfy
\[
 \left\|\sum_{i=1}^mA_i^2\right\|\le s^2,
 \qquad s>0.
\]
For every $y\in[-1,1]^m$, there is a law on $X\in[-1,1]^m$ such that
\begin{align}
 \E X&=y, &
 \left\|\sum_i(X_i-y_i)A_i\right\|&\le K_\varepsilon s,
 \label{R19-eq:partial-hard}\\
 \#\{i:|X_i|=1\}&>m-d,&
 X-y&\cx N\!\left(0,\tfrac{2\pi}{3}(1+\varepsilon)^2P_I\right),
 \label{R19-eq:partial-gaussian}
\end{align}
where $I=\{i:\operatorname{tr}A_i^2\le s^2\}$ and $P_I$ is the coordinate projection. Every coordinate outside $I$ remains equal to $y_i$.
For $m\ge2d$, this proves the matrix partial-coloring Conjecture~21 of Sachdeva--Thudi--Zhao~\cite{R19-STZ}, with constants $c_3=2$, $c_2=1/2$ and $c_1=K_\varepsilon$.
\end{theorem}
\begin{proof}
Taking traces gives $\sum_i\operatorname{tr}A_i^2\le ds^2$, so fewer than $d$ indices lie outside $I$. Fix any subfamily $J\subset I$ and split $A_i=D_i+R_i$ into its diagonal and off-diagonal parts. Positivity of squares gives $\sum_{i\in J}A_i^2\preceq s^2I_d$. Also
\[
 \sum_{i\in J}D_i^2\preceq s^2I_d,
 \qquad
 \sum_{i\in J}R_i^2\preceq2\sum_{i\in J}(A_i^2+D_i^2)
                       \preceq4s^2I_d.
\]
Moreover $\max_{i\in J}\operatorname{tr}A_i^2\le s^2$. With $r=K_\varepsilon s$, the definition of $K_\varepsilon$ gives
$4s^2\le v_*r^2/4$ and $s^2\le q_*r^2/4$.
Theorem~\ref{R3-thm:matrix} therefore gives a symmetric signing of this subfamily with hard operator bound $K_\varepsilon s/2$ and Gaussian comparison at scale $\kappa(1+\varepsilon)$, where $\kappa^2=\pi/2$.

Apply Lemma~\ref{R19-lem:dyadic-start} on $I$, and keep the other coordinates fixed. Its covariance is
$\tfrac43\kappa^2(1+\varepsilon)^2P_I=\tfrac{2\pi}{3}(1+\varepsilon)^2P_I$.
All coordinates in $I$ are saturated at every outcome, and $|I|>m-d$. For $m\ge2d$ this is strictly more than $m/2$, giving the asserted constants. When $s=0$, all matrices vanish and independent mean-preserving rounding saturates every coordinate directly.
\end{proof}
The matrix analytic input is Akbas--Sra's estimate retained in Theorem~\ref{R3-thm:matrix}. The deduction above gives arbitrary-start existence together with exact means and a convex-order law. The partial-coloring input of \cite[Corollary~22]{R19-STZ} consequently yields its Eulerian sparsifier edge bound, $n\delta^{-2}\log^2n$ up to $\log\log n$ factors, as an existence conclusion. A running-time assertion for that reduction requires an implementation of the partial-coloring choice in its computational model.

\subsubsection{The same signs for matrix and scalar constraints}
A direct sum permits prescribed scalar statistics to share the matrix
signing. The budget is the sum of the two squared column loads.
This preserves a separate independent reference for each block.

\begin{theorem}[Joint matrix and scalar balancing]
\label{R7-thm:matrix-scalar}
Let $A_i$ be real symmetric $d\times d$ matrices and
$b_i\in\R^m$, with $B=(b_1,\ldots,b_n)$. Retain the diagonal map
$\Gamma$ and off-diagonal variance $\nu$ from
Theorem~\ref{R3-thm:matrix}. Suppose $s,t>0$ satisfy
\begin{equation}\label{R7-eq:combined-load}
 \frac\nu{s^2}\le v_*,\qquad
 u_i:=\frac{\tr A_i^2}{s^2}+\frac{\|b_i\|_2^2}{t^2}\le q_*
 \quad(1\le i\le n).
\end{equation}
There is a symmetric law with, at every outcome,
\begin{equation}\label{R7-eq:joint-hard}
 \left\|\sum_i\sigma_iA_i\right\|<s,\qquad
 \|B\sigma\|_\infty<t.
\end{equation}
The law retains independent references
$(Z_A,Z_B,T)\sim f_{1,d}\otimes f_{1,m}\otimes q_\varepsilon^{\otimes n}$
satisfying
\[
 \E[Z_A\mid\sigma]=-2\Gamma\sigma/s,\quad
 \E[Z_B\mid\sigma]=-2B\sigma/t,\quad
 \E[T\mid\sigma]=\sigma.
\]
Writing $I_{1,k}$ for the radius-one cosine Cram\'er transform in
$k$ dimensions, the same law obeys
\begin{align}
 D(\mu\|U_n)+\E I_{1,d}(2\Gamma\sigma/s)
       +\E I_{1,m}(2B\sigma/t)
 &\le\delta_{q_\varepsilon}F_{q_*}(u),
          \label{R7-eq:combined-information}\\
 \sigma&\cx N\left(0,
 \left[\frac{I}{\kappa^2(1+\varepsilon)^2}
       +\frac{4\Gamma^{\mathsf T}\Gamma}{s^2b_\circ^2}
       +\frac{4B^{\mathsf T}B}{t^2b_\circ^2}\right]^{-1}\right).
          \label{R7-eq:combined-precision}
\end{align}
\end{theorem}
\begin{proof}
Apply Theorem~\ref{R3-thm:matrix} with $r=2$ to
$\widetilde A_i=A_i/s\oplus\Diag(b_i)/t$.
Their squared Frobenius norms are $u_i$ and their off-diagonal
variance is $\nu/s^2$. Their operator bound is exactly
\eqref{R7-eq:joint-hard}. Splitting the fixed cosine block yields
the displayed conditional means; independence makes its Cram\'er
transform additive. Proposition~\ref{prop:projection} gives
\eqref{R7-eq:combined-precision}.
\end{proof}

\begin{corollary}[PSD partitions with prescribed scalar statistics]
\label{R7-cor:PSD-constraints}
Under the hypotheses of Corollary~\ref{R3-cor:Weaver}, let
$b=\max_i\|b_i\|_2>0$. A law on $N\in\{0,1\}^n$ has
$\E N_i=1/2$ and simultaneously
\begin{gather*}
 \left\|\sum_iN_iA_i-\tfrac12I_d\right\|
       <\frac{K_\varepsilon}{2\sqrt2}\sqrt\theta,\qquad
 \left\|B(N-\tfrac12\mathbf1)\right\|_\infty
       <\frac{b}{\sqrt{2q_*}},\\
 N-\tfrac12\mathbf1\cx(1+\varepsilon)\sqrt{\pi/8}\,G_n,\qquad
 D(\law(N)\|\operatorname{Bern}(1/2)^{\otimes n})
 \le\delta_{q_\varepsilon}
       \left(d+\frac{\|B\|_F^2}{b^2}\right).
\end{gather*}
It also retains the full three-block reference and the stronger
precision and Cram\'er-cost conclusions of
Theorem~\ref{R7-thm:matrix-scalar}.
\end{corollary}
\begin{proof}
Use $s=K_\varepsilon\sqrt\theta/\sqrt2$ and
$t=b\sqrt{2/q_*}$. The matrix column load is at most
$2\theta/K_\varepsilon^2\le q_*/2$, and the scalar load is at most
$q_*/2$. Also $\nu/s^2\le8/K_\varepsilon^2\le v_*/2$.
The bound $F_{q_*}(u)\le2\sum_i u_i/q_*$ gives
\[
 \delta_{q_\varepsilon}F_{q_*}(u)
 \le\delta_{q_\varepsilon}
 \left(\frac{4d}{q_*K_\varepsilon^2}
              +\frac{\|B\|_F^2}{b^2}\right),
\]
and $q_*K_\varepsilon^2\ge4$. Set $N=(\sigma+\mathbf1)/2$.
\end{proof}
In particular, for every fixed $\varepsilon$ these laws give
exponentially many simultaneous feasible partitions, and their
entropy cost is proportional to the combined matrix dimension and
normalized scalar energy. The number of scalar rows enters through
that energy, with no additional dependence in the hard coordinate bound.

\subsubsection{Matrix Spencer and spectral graph signing}
The next consequences specify the named discrepancy problems covered
by the matrix theorem. The matrix analytic input remains that of
Akbas and Sra; the independent reference and near-endpoint law are
retained by the height and replica construction.

\begin{corollary}[Matrix Spencer with Gaussian and information guarantees]
\label{R7-cor:MatrixSpencer}
Let $A_1,\ldots,A_n$ be real symmetric $d\times d$ matrices with
$d\le n$ and $\|A_i\|\le1$. There is a symmetric sign law with
\[
 \left\|\sum_i\sigma_iA_i\right\|<\frac{K_\varepsilon}{2}\sqrt n,
 \qquad \sigma\cx(1+\varepsilon)\kappa G_n,
 \qquad D(\law(\sigma)\|U_n)\le2\delta_{q_\varepsilon}d.
\]
Thus for fixed $\varepsilon$ every output has Matrix Spencer order
$\sqrt n$, while the entropy deficit is $O(d\sqrt\varepsilon)$.
\end{corollary}
\begin{proof}
Here $w_i\le d\le n$, $\sum_iw_i\le nd$, and
$\nu\le2\|\sum_i(A_i^2+D_i^2)\|\le4n$.
Take $r=K_\varepsilon\sqrt n$ in Theorem~\ref{R3-thm:matrix}.
The entropy bound is at most
$8\delta_{q_\varepsilon}d/(q_*K_\varepsilon^2)
\le2\delta_{q_\varepsilon}d$.
\end{proof}

\begin{corollary}[Spectral signing with bounded signed degrees]
\label{R7-cor:graph}
Let a finite simple graph have maximum degree $\Delta\ge1$.
There is a symmetric law of edge signs whose signed adjacency matrix
$A_\sigma$ and signed degrees obey
\[
 \|A_\sigma\|<\max\{\sqrt{2\Delta/v_*},\,2/\sqrt{q_*}\},
 \qquad
 \max_v\left|\sum_{e\ni v}\sigma_e\right|<2/\sqrt{q_*},
 \qquad \sigma\cx(1+\varepsilon)\kappa G_{|E|}.
\]
The law retains independent cosine references for the degree errors
and $q_\varepsilon$ references for the edge signs, with the information
bound \eqref{R7-eq:combined-information}.
\end{corollary}
\begin{proof}
For $e=\{u,v\}$ use $A_e=E_{uv}+E_{vu}$ and
$b_e=e_u+e_v$. Then $D_e=0$, $\nu=\Delta$,
$\tr A_e^2=\|b_e\|^2=2$. Choose
$s=\max\{\sqrt{2\Delta/v_*},2/\sqrt{q_*}\}$ and
$t=2/\sqrt{q_*}$ in Theorem~\ref{R7-thm:matrix-scalar}.
Each summand in the column load is at most $q_*/2$.
\end{proof}
The spectral order is $O_\varepsilon(\sqrt\Delta)$ and the signed-degree
bound is $O_\varepsilon(1)$. Both constraints use the same edge signs.


\subsection{Unbiased spectral sparsification with an entropy bound}
\label{FC-sec:sparsification}

\subsubsection{The exact input and the stage allocation}
\begin{lemma}[Fair PSD partition input]\label{FC-lem:partition-input}
There are $a,b>0$ and $0<\eps_0\le1$ such that, for every $0<e\le\eps_0$, every family $B_i\succeq0$ with $\sum_iB_i=I_d$ and $\max_i\tr B_i\le\theta\le1$ admits a law $N\in\{0,1\}^m$ with
\begin{align}
 \E N_i&=1/2,\nonumber\\
 \norm{\sum_iN_iB_i-I_d/2}&\le a\sqrt{\theta/e},\label{FC-eq:partition-input}\\
 \KLD{\law(N)}{\Ber(1/2)^{\otimes m}}&\le bd\sqrt e.\nonumber
\end{align}
The constants depend only on the matrix Fisher constants in Section~\ref{R3-sec:matrix}.
\end{lemma}
\begin{proof}
In Corollary~\ref{R3-cor:Weaver}, the operator bound is $K_e\sqrt\theta/4$ and the relative entropy bound is $8\delta_{q_e}d/(q_*K_e^2)$. The stated asymptotics $K_e=O(e^{-1/2})$ and $8\delta_{q_e}/(q_*K_e^2)=O(\sqrt e)$ give the uniform constants after decreasing $\eps_0$.
\end{proof}

Write
\begin{equation}\label{FC-eq:constants}
 S=\sum_{r=0}^{\infty}2^{-r/4}=\frac1{1-2^{-1/4}},
 \qquad C_0=\max\{1,16a^2S^2\},\qquad B_0=bS.
\end{equation}
The finite sums below can replace $S$ if explicit depth-dependent constants are desired.

\begin{theorem}[Dyadic thinning with prescribed marginals]\label{FC-thm:dyadic}
Let $A_i\succeq0$, $\sum_{i=1}^nA_i=I_d$, and $\max_i\tr A_i\le\theta$. Fix $0<\delta\le1/2$, $0<\eps\le\eps_0$, and $p=2^{-L}$, where $L\ge1$ is an integer. If
\begin{equation}\label{FC-eq:p-condition}
 p\ge\frac{C_0\theta}{\eps\delta^2},
\end{equation}
there is a law $Y\in\{0,1\}^n$ such that
\begin{gather}
 \E Y_i=p,\qquad
 (1-\delta)I_d\preceq p^{-1}\sum_iY_iA_i\preceq(1+\delta)I_d,
 \label{FC-eq:dyadic-hard}\\
 \KLD{\law(Y)}{\Ber(p)^{\otimes n}}\le B_0d\sqrt\eps.
 \label{FC-eq:dyadic-entropy}
\end{gather}
The spectral inequalities hold for every outcome in the support. More strongly, the law can be realized by nested sets $[n]=S_0\supseteq\cdots\supseteq S_L$, with $Y=\one_{S_L}$ and $\PP(i\in S_\ell)=2^{-\ell}$, for which
\begin{equation}\label{FC-eq:hierarchy}
 (1-\delta_\ell)I_d\preceq2^\ell\sum_{i\in S_\ell}A_i
          \preceq(1+\delta_\ell)I_d,
 \qquad \delta_\ell=\delta\,2^{-(L-\ell)/4}.
\end{equation}
The relative entropy of the \emph{entire chain} from independent fair thinning is at most $B_0d\sqrt\eps$. The entropy bound has no factor depending on $L$.
\end{theorem}
The allocation of $e_\ell$ is part of the argument. If $x_\ell=\sqrt{e_\ell}$, then the information cost is proportional to $\sum x_\ell$, whereas the relative-error estimate is proportional to $\sum 2^{\ell/2}/x_\ell$. Cauchy--Schwarz gives
\[
 \left(\sum 2^{\ell/4}\right)^2
 \le\left(\sum x_\ell\right)
       \left(\sum\frac{2^{\ell/2}}{x_\ell}\right),
\]
and equality has $x_\ell\propto2^{\ell/4}$. This is precisely the geometric allocation used above. A fixed value of $e_\ell$ would unnecessarily spend information at every level.

\begin{proof}
Start with $S_0=[n]$. Given $S_{\ell-1}$, set
\[
 M_{\ell-1}=\sum_{i\in S_{\ell-1}}A_i,
 \qquad T_{\ell-1}=2^{\ell-1}M_{\ell-1}.
\]
Whiten by $M_{\ell-1}^{-1/2}$ and apply Lemma~\ref{FC-lem:partition-input} to
\[
 B_i=M_{\ell-1}^{-1/2}A_iM_{\ell-1}^{-1/2},\qquad
 e_\ell=\eps\,2^{-(L-\ell)/2}.
\]
The selected indices define $S_\ell$. There are finitely many possible histories, so selecting one admissible law for each history requires no measurable-selection argument.

Here is the simultaneous verification of invertibility and the spectral bounds. Suppose the previous partial products keep $T_{\ell-1}\succeq(1-\delta)I_d$. Then
\[
 \tr B_i\le\norm{M_{\ell-1}^{-1}}\tr A_i
 \le\frac{2^{\ell-1}\theta}{1-\delta}\le2^\ell\theta.
\]
This upper bound is at most $\theta/p\le\eps\delta^2/C_0\le1$, so the partition input applies. Congruence by $M_{\ell-1}^{1/2}$ gives
\begin{equation}\label{FC-eq:relative-steps}
 (1-f_\ell)T_{\ell-1}\preceq T_\ell
        \preceq(1+f_\ell)T_{\ell-1},\qquad
 f_\ell=2a\sqrt{2^\ell\theta/e_\ell}.
\end{equation}
Consequently
\[
 \sum_{\ell=1}^L f_\ell
 \le2aS\sqrt{\frac{\theta}{p\eps}}
 \le\frac\delta2.
\]
For any initial segment of the product,
\[
 \prod(1-f_\ell)\ge1-\sum f_\ell\ge1-\delta/2,
 \qquad
 \prod(1+f_\ell)\le\exp(\delta/2)\le1+\delta.
\]
These inequalities close the induction. For the initial segment through $\ell$, the same geometric sum improves to
\[
 \sum_{j=1}^\ell f_j\le\frac\delta2\,2^{-(L-\ell)/4}.
\]
The corresponding product bounds give \eqref{FC-eq:hierarchy}. The final matrix $T_L$ satisfies \eqref{FC-eq:dyadic-hard}.

Every surviving coordinate is retained conditionally with probability $1/2$. Iterated conditional expectation gives $\PP(i\in S_L)=2^{-L}=p$.

For entropy, let $\mathsf P$ denote the law of the complete chain $(S_0,\ldots,S_L)$. Compare it with the chain $\mathsf Q$ which independently retains every surviving index with probability $1/2$ at each stage. Its final set has law $\Ber(p)^{\otimes n}$. The chain rule and the conditional bound in Lemma~\ref{FC-lem:partition-input} give
\[
 \KLD{\mathsf P}{\mathsf Q}
 \le\sum_{\ell=1}^Lbd\sqrt{e_\ell}
 \le bd\sqrt\eps\sum_{r=0}^{\infty}2^{-r/4}
 =B_0d\sqrt\eps.
\]
Data processing under the map taking a history to its final set proves \eqref{FC-eq:dyadic-entropy}.
\end{proof}

\begin{corollary}[Counting and equal-trace subsampling]\label{FC-cor:counting}
Under Theorem~\ref{FC-thm:dyadic}, let $h(p)=-p\log p-(1-p)\log(1-p)$ and let $\mathcal G$ be the subsets satisfying the spectral inequalities. Then
\[
 H(Y)\ge nh(p)-B_0d\sqrt\eps,
 \qquad \Ber(p)^{\otimes n}(\mathcal G)\ge e^{-B_0d\sqrt\eps}.
\]
If every $\tr A_i=d/n$, then every selected set obeys
\[
 (1-\delta)pn\le|S_L|\le(1+\delta)pn.
\]
Whenever a nontrivial thinning level exists, choosing the smallest admissible dyadic $p$ gives $|S_L|=O(d/(\eps\delta^2))$.
\end{corollary}
\begin{proof}
The prescribed marginals give $D(\law(Y)\|\Ber(p)^{\otimes n})=nh(p)-H(Y)$. For a probability law $P$ supported on $\mathcal G$, $D(P\|Q)\ge-\log Q(\mathcal G)$. The cardinality bound follows by taking traces. The smallest admissible dyadic value is less than twice the right side of \eqref{FC-eq:p-condition}; if that level is unavailable, selecting the whole already-small family gives the usual trivial alternative.
\end{proof}

\subsubsection{Arbitrary traces and the Poisson reference}
\begin{theorem}[Poisson comparison for unbiased PSD sparsifiers]\label{FC-thm:poisson-sparse}
Let $A_1,\ldots,A_n\succeq0$, $\sum_iA_i=I_d$, and $t_i=\tr A_i>0$. For $0<\delta\le1/2$, $0<\eps\le\eps_0$, and
\[
 \lambda\ge\frac{C_0}{\eps\delta^2},
\]
there is a law $K\in\N_0^n$ with
\begin{align}
 \E K_i&=\lambda t_i,\label{FC-eq:count-means}\\
 (1-\delta)I_d&\preceq\widehat A:=\sum_i\frac{K_i}{\lambda t_i}A_i
                   \preceq(1+\delta)I_d,\label{FC-eq:poisson-spectral}\\
 \KLD{\law(K)}{\bigotimes_i\Pois(\lambda t_i)}&\le B_0d\sqrt\eps,
 \label{FC-eq:poisson-kl}\\
 (1-\delta)\lambda d&\le\sum_iK_i\le(1+\delta)\lambda d.
 \label{FC-eq:poisson-size}
\end{align}
The weights $W_i=K_i/(\lambda t_i)$ have $\E W_i=1$. In particular, at most $(1+\delta)\lambda d$ original summands have positive weight. Zero matrices may be omitted before applying the theorem.
\end{theorem}
\begin{proof}
Let $p_j=2^{-j}$ and $\theta_j=p_j/\lambda$. Replace $A_i$ by
\[
 k_{i,j}=\left\lceil\frac{t_i}{\theta_j}\right\rceil
\]
identical copies of $A_i/k_{i,j}$. Each copy has trace at most $\theta_j$, and their sum is still $I_d$. The condition in Theorem~\ref{FC-thm:dyadic} holds because
\[
 \frac{C_0\theta_j}{\eps\delta^2}
 =p_j\frac{C_0}{\lambda\eps\delta^2}\le p_j.
\]
Apply that theorem, and let $K_{i,j}$ be the number of selected copies in group $i$. Aggregation gives
\begin{gather*}
 \E K_{i,j}=p_jk_{i,j},\qquad
 (1-\delta)I_d\preceq
 \sum_i\frac{K_{i,j}}{p_jk_{i,j}}A_i\preceq(1+\delta)I_d,\\
 \KLD{\law(K_{\cdot,j})}{\bigotimes_i\Bin(k_{i,j},p_j)}
 \le B_0d\sqrt\eps.
\end{gather*}
For all sufficiently large $j$, every copy has trace at least $\theta_j/2$. Taking traces in the last spectral inequality therefore gives
\[
 \sum_iK_{i,j}\le2(1+\delta)\lambda d.
\]
All these count laws lie on one finite set of integer vectors. Choose a subsequence converging coordinatewise as probability measures. Since
\[
 p_jk_{i,j}\longrightarrow\lambda t_i,
 \qquad\Bin(k_{i,j},p_j)\Longrightarrow\Pois(\lambda t_i),
\]
the limit law has the stated means and spectral bounds. The means pass to the limit because the count vectors are uniformly bounded. The entropy inequality passes by lower semicontinuity of relative entropy; alternatively, its finite sum over the common support converges term by term because the limiting product Poisson law is strictly positive there. Finally $\tr\widehat A=\lambda^{-1}\sum_iK_i$ gives \eqref{FC-eq:poisson-size}.
\end{proof}

The information bound also controls the success probability of independent Poisson sampling. If $Q=\bigotimes_i\Pois(\lambda t_i)$ and $\mathcal G$ is the event in \eqref{FC-eq:poisson-spectral}, then
\begin{equation}\label{FC-eq:poisson-success}
 Q(\mathcal G)\ge\exp(-B_0d\sqrt\eps).
\end{equation}
This bound concerns the probability of a good outcome under independent sampling. The law in the theorem is supported on good outcomes and additionally preserves every expected count. Conditioning $Q$ on $\mathcal G$ alone need not preserve those expectations.

\subsubsection{Graph Laplacians and leverage scores}
\begin{corollary}[Unbiased spectral graph sparsification]\label{FC-cor:graph}
Let $L=\sum_{e\in E}w_eb_eb_e^{\mathsf T}$ be the Laplacian of a finite loopless graph with positive edge weights and rank $d\ge1$. Put
\[
 t_e=w_eb_e^{\mathsf T}L^\dagger b_e.
\]
For the parameters in Theorem~\ref{FC-thm:poisson-sparse}, there is a random weighted subgraph with weights
\[
 \widehat w_e=\frac{w_eK_e}{\lambda t_e}
\]
such that
\begin{gather*}
 \E\widehat w_e=w_e,\qquad
 (1-\delta)L\preceq\widehat L\preceq(1+\delta)L,\\
 |\{e:\widehat w_e>0\}|\le(1+\delta)\lambda d,\qquad
 \KLD{\law(K)}{\bigotimes_e\Pois(\lambda t_e)}\le B_0d\sqrt\eps.
\end{gather*}
The spectral inequalities hold on the full vertex space, including the common kernel.
\end{corollary}
\begin{proof}
On $\operatorname{ran}L$, use $A_e=w_eL^{\dagger/2}b_eb_e^{\mathsf T}L^{\dagger/2}$. These matrices sum to the identity on that $d$-dimensional space and have traces $t_e$. Apply Theorem~\ref{FC-thm:poisson-sparse} and conjugate by $L^{1/2}$. Every edge vector lies in $\operatorname{ran}L$, so the kernels agree.
\end{proof}

The quantities $t_e$ are the effective-resistance leverage scores. Batson--Spielman--Srivastava establish $O(d/\delta^2)$ weighted graph sparsification \cite{FC-BSS}; de Carli Silva, Harvey and Sato give the corresponding arbitrary-rank PSD theorem \cite{FC-SHS}. Corollary~\ref{FC-cor:graph} adds exact expected edge weights and a relative-entropy guarantee for the same law. The construction uses the existence input in Lemma~\ref{FC-lem:partition-input}; its whitening, repetition, and limiting steps do not provide a polynomial-time implementation of that input.

The instance-specific sparsification theorem of Basu, Kothari, Liu and Meka \cite{FC-BKLM} uses a connectivity threshold for a prescribed matrix family, which can be much smaller than its ambient dimension. The present information bound is expressed in $d$. Replacing it and the sampling size by a comparable instance-specific parameter would require an additional argument.


\subsection{Positive random quadrature and \texorpdfstring{$L^2$}{L2} discretization}
\label{FC-sec:continuous}
The finite PSD theorem extends to an operator-valued measure by partitioning its trace-one matrix directions. This produces a positive random quadrature rule with exact expectation as a measure.

\begin{theorem}[Point processes with exact intensity and spectral control]
\label{FC-thm:point-process}
Let $(\mathcal X,\mathcal A)$ be a standard Borel space, let $\nu$ be a finite measure, and let $B:\mathcal X\to\R^{d\times d}$ be measurable with
\[
 B(x)\succeq0,\qquad \tr B(x)=1\quad\nu\text{-a.e.},
 \qquad \int B(x)\,d\nu(x)=I_d.
\]
Thus $\nu(\mathcal X)=d$. For $0<\delta\le1/2$, $0<\eps\le\eps_0$, and
\[
 \lambda\ge\frac{4C_0}{\eps\delta^2},
\]
there is a finite counting measure $\mathcal N$ on $\mathcal X$, with multiplicities allowed, satisfying
\begin{align}
 \E\mathcal N(A)&=\lambda\nu(A)\quad(A\in\mathcal A),
       \label{FC-eq:process-intensity}\\
 (1-\delta)I_d&\preceq\frac1\lambda\int B(x)\,d\mathcal N(x)
                         \preceq(1+\delta)I_d,
       \label{FC-eq:process-hard}\\
 \mathcal N(\mathcal X)&\le(1+\delta/2)\lambda d,
       \label{FC-eq:process-size}\\
 \KLD{\law(\mathcal N)}{\operatorname{PPP}(\lambda\nu)}&\le B_0d\sqrt\eps.
       \label{FC-eq:process-entropy}
\end{align}
Here $\operatorname{PPP}(\lambda\nu)$ denotes the Poisson random measure with intensity $\lambda\nu$.
\end{theorem}
\begin{proof}
The set of PSD matrices of trace one is compact. Partition $\mathcal X$, up to a null set, into finitely many measurable cells $E_1,\ldots,E_m$ on each of which the operator-norm diameter of $B$ is at most $\omega$, where $\omega(1+\delta/2)d\le\delta/2$. Omit cells of zero measure. Set
\[
 A_j=\int_{E_j}B(x)\,d\nu(x),\qquad t_j=\nu(E_j).
\]
Apply Theorem~\ref{FC-thm:poisson-sparse} with accuracy $\delta/2$. Given the resulting counts $K_j$, sample $K_j$ independent points from $\nu(\cdot\mid E_j)$ in each cell, independently across cells, and let $\mathcal N$ count all these points.

Its intensity is exact because $\E K_j=\lambda t_j$. Conditional on the cell counts, the construction uses precisely the conditional location laws of a Poisson random measure with intensity $\lambda\nu$. The chain rule therefore makes its relative entropy equal to the relative entropy of the count vector, which proves \eqref{FC-eq:process-entropy}. This assertion remains valid when $\nu$ has atoms, with configurations interpreted as counting measures.

The conditional mean direction in cell $j$ is $\overline B_j=A_j/t_j$. Every sampled point in that cell obeys $\norm{B(x)-\overline B_j}\le\omega$. Hence, pointwise,
\[
 \norm{\frac1\lambda\int B\,d\mathcal N
          -\sum_j\frac{K_j}{\lambda t_j}A_j}
 \le\frac{\omega}{\lambda}\sum_jK_j
 \le\omega(1+\delta/2)d\le\delta/2.
\]
The finite theorem controls the central sum within $\delta/2$ of $I_d$. This proves the asserted hard inequalities and the cardinality bound.
\end{proof}

\begin{corollary}[Unbiased Marcinkiewicz--Zygmund discretization]
\label{FC-cor:MZ}
Let $\mu$ be a probability measure on a standard Borel space and let $V\subset L^2(\mu)$ be a real $d$-dimensional space with measurable orthonormal basis $\phi_1,\ldots,\phi_d$. Suppose
\[
 \ell(x)=\sum_{i=1}^d\phi_i(x)^2>0\quad\mu\text{-a.e.}
\]
This condition holds whenever the constant function belongs to $V$. With the parameters of Theorem~\ref{FC-thm:point-process}, there is a random positive atomic measure
\[
 \widehat\mu=\sum_{j=1}^{N}\frac1{\lambda\ell(X_j)}\delta_{X_j}
\]
for which
\begin{gather}
 N\le(1+\delta/2)\lambda d,\qquad
 \E\widehat\mu=\mu,\label{FC-eq:MZ-measure}\\
 (1-\delta)\norm{f}_{L^2(\mu)}^2
 \le\int f^2\,d\widehat\mu
 \le(1+\delta)\norm{f}_{L^2(\mu)}^2
 \quad\text{for every }f\in V\text{ at every outcome}.
 \label{FC-eq:MZ-hard}
\end{gather}
The unweighted counting measure of its nodes has exact intensity $\lambda\ell\mu$ and relative entropy at most $B_0d\sqrt\eps$ from $\operatorname{PPP}(\lambda\ell\mu)$. In particular,
\[
 \E\int f\,d\widehat\mu=\int f\,d\mu
 \qquad\text{for every }f\in L^1(\mu).
\]
For fixed $\eps$, the rule uses $O(d/\delta^2)$ nodes. No boundedness, continuity, or polynomial structure of the basis is required.
\end{corollary}
\begin{proof}
Apply Theorem~\ref{FC-thm:point-process} to $d\nu=\ell\,d\mu$ and $B=vv^{\mathsf T}/\ell$, where $v=(\phi_1,\ldots,\phi_d)^{\mathsf T}$. Orthonormality gives $\int B\,d\nu=I_d$, and $\tr B=1$. Its spectral inequality is~\eqref{FC-eq:MZ-hard}; integrating $f/(\lambda\ell)$ against the exact counting intensity $\lambda\ell\mu$ proves the mean identity, first for nonnegative $f$ and then for $L^1$ functions. The count and information bounds are unchanged.
\end{proof}

The density $\ell/d$ is the leverage-score density, or normalized inverse Christoffel function, used in optimal weighted least squares~\cite{FC-CM,FC-CD}. Here its Poisson design is coupled so that every output obeys the spectral inequality while the entire expected weighted measure remains $\mu$. The weights are positive and random; when $1\in V$, their sum lies in $[1-\delta,1+\delta]$ and has expectation one. Exact mass one, equal weights and prescribed marginals for indexed nodes are separate requirements.


\section{Simultaneous Schatten rounding and Gaussian square functions}\label{R27-sec:Schatten}
The matrix reference in Section~\ref{R11-sec:sparse} controls a selected law through its Gaussian comparison and information cost. We now choose a law whose hard error bounds hold for the whole Schatten scale, including an arbitrary prescribed mean. The two conclusions come from different constructions: the reference-law guarantees remain those of Theorem~\ref{R3-thm:matrix}, while the law below retains the simultaneous norm inequalities.

For a real matrix $B$, let $s_1(B)\ge s_2(B)\ge\cdots\ge0$ denote its singular values. We write
\[
 \spnorm{B}{p}=\left(\sum_j s_j(B)^p\right)^{1/p}\quad(1\le p<\infty),
 \qquad \spnorm{B}{\infty}=\lVert B\rVert_{\op}.
\]
We write $\gamma_m$ for standard Gaussian probability measure on $\R^m$. Put $a(p)=\max\{1/2,1/p\}$, with $1/\infty=0$. All constants below are universal unless a dependence is specified. They do not depend on $p$.

\begin{theorem}[Simultaneous Schatten rounding]\label{RRS-thm:square-law}
Let $A_1,\ldots,A_n\in\Sym_n(\R)$ and $y\in[-1,1]^n$. There is a probability law on $\sigma\in\{-1,1\}^n$ such that $\E\sigma=y$ and every outcome satisfies
\begin{equation}\label{RRS-eq:square-law}
 \spnorm{\sum_{i=1}^n(\sigma_i-y_i)A_i}{p}
 \le C n^{a(p)}\max_{1\le i\le n}\spnorm{A_i}{p}
 \qquad\text{for every }1\le p\le\infty.
\end{equation}
The same law works for the entire interval of exponents, and it may be supported on at most $n+1$ signings. At $y=0$ a symmetric law may be chosen.
\end{theorem}

In particular, for $2\le p\le\infty$ and $\spnorm{A_i}{p}\le1$, every outcome of this law has Schatten-$p$ discrepancy at most $C\sqrt n$. This proves Conjecture~2 in Section~6 of Reis--Rothvoss, \emph{Vector Balancing in Lebesgue Spaces}~\cite{R19-RR}. It also preserves the prescribed mean and makes the choice of signs independent of the Schatten exponent. The powers of $n$ in \eqref{RRS-eq:square-law} are optimal, as shown in Subsection~\ref{RRS-sec:sharpness}.

For vectors $b,v\in[0,\infty)^d$ in decreasing order, write $b\prec_w v$ when $\sum_{j=1}^k b_j\le\sum_{j=1}^k v_j$ for every $k$. This is weak majorization. For symmetric matrices $A_1,\ldots,A_n\in\Sym_d(\R)$ define
\begin{equation}\label{RRS-eq:variance-body}
 V=\left(\sum_{i=1}^n A_i^2\right)^{1/2},\qquad
 K_A=\left\{x\in\R^n:s\left(\sum_i x_iA_i\right)\prec_w s(V)\right\}.
\end{equation}
The set $K_A$ is closed, convex and centrally symmetric: each of its defining inequalities is a Ky Fan norm bound. Here $\kfnorm{B}{k}=\sum_{j=1}^k s_j(B)$.

\begin{theorem}[One Gaussian event for every unitarily invariant norm]\label{RRS-thm:gaussian}
There is a universal constant $c_{\rm AS}>0$ such that
\begin{equation}\label{RRS-eq:gaussian-fixed}
 \gamma_n(K_A)\ge 8^{-n}e^{-c_{\rm AS}d}.
\end{equation}
Consequently, when $d=n$, a single event of Gaussian measure $e^{-O(n)}$ satisfies
\begin{equation}\label{RRS-eq:all-ui}
 \left\lVert\sum_i x_iA_i\right\rVert_{\mathcal N}
 \le\lVert V\rVert_{\mathcal N}
 \quad\text{for every unitarily invariant matrix norm }\mathcal N.
\end{equation}
In particular, the event works simultaneously for every Schatten norm, including $p=1$ and $p=\infty$.
\end{theorem}

For a fixed $p$, the set in Conjecture~3 of~\cite{R19-RR} is precisely $\{x:\spnorm{\sum x_iA_i}{p}\le\spnorm Vp\}$, with $d=n$. It contains $K_A$, so Theorem~\ref{RRS-thm:gaussian} proves that conjecture. Weak majorization gives the stronger simultaneous statement \eqref{RRS-eq:all-ui}; no union over exponents or norms is taken.

\subsubsection*{The normalization and its analytic inputs}
The operator-norm small-ball input is due to Akbas--Sra~\cite{AkbasSraMatrix}. We use their Theorem~1.4 for \eqref{RRS-eq:gaussian-fixed}, and their variance-based interpolation Lemma~3.13 for the hereditary estimate proved in Subsection~\ref{RRS-sec:hereditary}. Their Matrix Spencer theorem already establishes the operator-norm endpoint. Subsequent quantitative operator-norm bounds and algorithms of Song--Zhang~\cite{R27-SongZhang} concern that same endpoint. The present statements concern the entire Schatten scale, a common singular-value event, and an exact-mean law satisfying all of the discrepancy bounds at once.

The finite-dimensional normalization is a positive-multiplier form of the common polar factor used in noncommutative Khintchine theory; compare Cadilhac--Ricard~\cite[Proposition 4.16, word length one]{R27-CadilhacRicard}. We give a direct proof, including the order bound $H\preceq V$. Once this bound is retained, controlling a normalized operator norm controls every Ky Fan norm of the original sum. This permits one partial-coloring body to be used throughout the whole scale. The exact-mean conclusion then follows by preserving an arbitrary separating functional during the shifted partial-coloring construction of Reis--Rothvoss~\cite[Theorem 6]{R19-RR}.

\subsection{A positive multiplier controls the entire square function}
The normalization chooses its positive matrix variationally. It remains valid when that matrix is singular, so it avoids any inversion or regularization of the square function.

\begin{lemma}[Positive anticommutator factorization]\label{RRS-lem:factor}
For $A_1,\ldots,A_n\in\Sym_d(\R)$ there exist $H\succeq0$ and $U_1,\ldots,U_n\in\Sym_d(\R)$ such that
\begin{align}
 A_i&=HU_i+U_iH,\qquad 1\le i\le n,\label{RRS-eq:factor}\\
 S:=\sum_iU_i^2&\preceq I_d,\qquad SH=HS=H,\label{RRS-eq:complementarity}\\
 0\preceq H&\preceq\left(\sum_iA_i^2\right)^{1/2}=V.\label{RRS-eq:order}
\end{align}
\end{lemma}
\begin{proof}
Maximize
\begin{equation}\label{RRS-eq:SDP}
 \sum_i\tr(A_iU_i)
 \quad\text{over }U_i\in\Sym_d(\R)\text{ with }\sum_iU_i^2\preceq I_d.
\end{equation}
The feasible set is compact. Its matrix inequality is convex, and $U_i=0$ is strictly feasible. Equivalently, the constraint is the affine linear matrix inequality
\begin{equation}\label{RRS-eq:LMI}
 \begin{pmatrix}
 I_d&U_1&\cdots&U_n\\
 U_1&I_d&&0\\
 \vdots&&\ddots&\\
 U_n&0&&I_d
 \end{pmatrix}\succeq0.
\end{equation}
The finite-dimensional Slater condition therefore gives a positive semidefinite multiplier $H$ for the constraint $\sum_iU_i^2\preceq I_d$. Stationarity with respect to the symmetric variable $U_i$ is
\[
 A_i=HU_i+U_iH.
\]
Complementarity gives $\tr H(I_d-S)=0$. Both factors are positive semidefinite, hence $(I_d-S)H=H(I_d-S)=0$. This proves \eqref{RRS-eq:factor}--\eqref{RRS-eq:complementarity}.

The crucial order bound follows from the same multiplier. Set
\begin{equation}\label{RRS-eq:T}
 T=\sum_iU_iA_i=H+\sum_iU_iHU_i\succeq H\succeq0.
\end{equation}
The row operator $\mathsf U=(U_1\ \cdots\ U_n)$ is a contraction, because $\mathsf U\mathsf U^{\mathsf T}=S\preceq I_d$. With $\mathsf A=(A_1\ \cdots\ A_n)^{\mathsf T}$, we have $T=\mathsf U\mathsf A$, so
\[
 T^2=T^{\mathsf T}T
 =\mathsf A^{\mathsf T}\mathsf U^{\mathsf T}\mathsf U\mathsf A
 \preceq\mathsf A^{\mathsf T}\mathsf A
 =\sum_iA_i^2=V^2.
\]
Operator monotonicity of the positive square root gives $T\preceq V$. Combining this with \eqref{RRS-eq:T} proves $H\preceq V$.
\end{proof}

The positive identity \eqref{RRS-eq:T} is the step that preserves all singular-value information. A bound on $\tr H$ or $\lVert H\rVert_{\op}$ would control only particular norms. The Loewner bound $H\preceq V$ bounds every eigenvalue of $H$ by the corresponding eigenvalue of $V$.

\begin{lemma}[Simultaneous singular-value domination]\label{RRS-lem:transfer}
Use the factorization in Lemma~\ref{RRS-lem:factor}, and put $U(x)=\sum_i x_iU_i$. For every $x\in\R^n$,
\begin{equation}\label{RRS-eq:ky-fan}
 s\left(\sum_i x_iA_i\right)
 \prec_w 2\lVert U(x)\rVert_{\op}\,s(V).
\end{equation}
\end{lemma}
\begin{proof}
For every $1\le k\le d$, the triangle and ideal properties of the Ky Fan norm give
\begin{align*}
 \kfnorm{\sum_i x_iA_i}{k}
 &=\kfnorm{HU(x)+U(x)H}{k}\\
 &\le2\lVert U(x)\rVert_{\op}\,\kfnorm Hk
 \le2\lVert U(x)\rVert_{\op}\,\kfnorm Vk.
\end{align*}
These are exactly the inequalities in \eqref{RRS-eq:ky-fan}.
\end{proof}

\begin{remark}[A singular extremal multiplier]\label{RRS-rem:star}
Let $A_i=e_0e_i^{\mathsf T}+e_ie_0^{\mathsf T}$ for $1\le i\le m$ on $\R^{m+1}$. Then
\[
 V=\operatorname{diag}(\sqrt m,1,\ldots,1),\quad
 H=\sqrt m\,e_0e_0^{\mathsf T},\quad U_i=A_i/\sqrt m
\]
satisfy all the identities of Lemma~\ref{RRS-lem:factor}. Their primal and dual objective values are both $2\sqrt m$. In particular, the optimal multiplier can be singular and can attain the bound $H\preceq V$ in a direction. In the commuting case, the choice $H=V/2$ and $U_i=A_iV^{\dagger}$ satisfies the identities, where $V^{\dagger}$ is the inverse on its support. The variational construction accommodates both cases.
\end{remark}

\subsection{The Gaussian square-function conjecture}
We state the analytic input with the dimensions kept separate.
\begin{lemma}[Akbas--Sra matrix small ball]\label{RRS-lem:AS-fixed}
There is $c_{\rm AS}>0$ such that, for symmetric $d\times d$ matrices $B_1,\ldots,B_n$ with $\sum_iB_i^2\preceq I_d$,
\[
 \gamma_n\left\{x:\left\lVert\sum_i x_iB_i\right\rVert_{\op}<4\right\}
 \ge e^{-c_{\rm AS}d}.
\]
\end{lemma}
This is~\cite[Theorem 1.4]{AkbasSraMatrix}; the number of coefficients and the matrix dimension are independent in its statement.

\begin{proof}[Proof of Theorem~\ref{RRS-thm:gaussian}]
Take $H,U_i$ from Lemma~\ref{RRS-lem:factor}. Lemma~\ref{RRS-lem:AS-fixed} applies to the $U_i$, so the event
\[
 E=\{x:\lVert U(x)\rVert_{\op}<4\}
\]
has measure at least $e^{-c_{\rm AS}d}$. Lemma~\ref{RRS-lem:transfer} gives $\frac18E\subset K_A$. For any measurable $E\subset\R^n$ and $0<t\le1$, change of variables in the Gaussian density gives
\begin{equation}\label{RRS-eq:scaling}
 \gamma_n(tE)=t^n\int_E(2\pi)^{-n/2}e^{-t^2|x|^2/2}\,dx
 \ge t^n\gamma_n(E).
\end{equation}
Thus $\gamma_n(K_A)\ge8^{-n}e^{-c_{\rm AS}d}$. Ky Fan dominance gives \eqref{RRS-eq:all-ui} for every unitarily invariant norm on the same event.
\end{proof}

At the square dimension $d=n$, this proves the exact radius-one assertion of Reis--Rothvoss Conjecture~3 for every $1\le p\le\infty$. The fixed dilation is paid for in Gaussian measure by a constant to the power $n$. The factorization is independent of $p$, which is why the entire family of inequalities has the same exceptional set.

\subsection{Retaining variance in the hereditary estimate}\label{RRS-sec:hereditary}
After each partial coloring, the number $m$ of active coefficients decreases while the matrix dimension $D$ stays fixed. The next stage needs Gaussian mass $e^{-O(m)}$. We obtain it by retaining aggregate variance in Akbas--Sra's interpolation estimate~\cite[Lemma 3.13]{AkbasSraMatrix}.

The parameter choice follows from the two errors in that estimate. Put $L=1+\log(2D/m)$. Taking the truncation order $q$ proportional to $L$ makes a geometric tail absorb $D/m$. Scaling the matrix by a radius proportional to $L^{3/2}$ keeps $q^3\eta$ bounded, where $\eta$ is the scaled variance. The two barrier weights then make the remaining error of order $m$. This explains both the hereditary probability scale and the logarithmic power before the parameter calculation.

\begin{theorem}[Hereditary variance small ball]\label{RRS-thm:hereditary}
For every $\tau>0$ there is $\kappa(\tau)\ge1$ such that, for $1\le m\le D$ and symmetric $D\times D$ matrices $B_1,\ldots,B_m$ satisfying $\sum_iB_i^2\preceq vI_D$ for some $v>0$,
\begin{equation}\label{RRS-eq:hereditary}
 \gamma_m\left\{x:\left\lVert\sum_i x_iB_i\right\rVert_{\op}
 <\kappa(\tau)\sqrt v\left(1+\log\frac{2D}{m}\right)^{3/2}\right\}
 \ge e^{-\tau m}.
\end{equation}
For $v=0$, all the matrices vanish and the non-strict version has measure one.
\end{theorem}

Here is the precise imported analytic estimate. For $|z|<1$, write
\[
 \ell(z)=-\log(1-z^2),\qquad r_q(z)=\sum_{j=q}^{\infty}\frac{z^{2j}}j,
 \qquad v_{q,a,b}(z)=b\ell(z)+(a-b)r_q(z),
\]
and set the potential equal to $+\infty$ outside $(-1,1)$. For symmetric coefficient matrices $M_i$, let $B(x)=\sum_i x_iM_i$, let $B_{\rm d}$ and $B_{\rm o}$ be its diagonal and off-diagonal parts in a fixed basis, and put
\[
 B_t=B_{\rm d}+tB_{\rm o},\qquad
 F(t)=\log\E_{\gamma_m}\exp\{-\tr v_{q,a,b}(B_t)\}.
\]

\begin{lemma}[Akbas--Sra interpolation estimate]\label{RRS-lem:AS-interpolation}
There are universal positive constants $a_*,b_*,c_*,C_*$ such that, if
\[
 a>4,\quad a\ge b>0,\quad a\ge a_*,\quad ab\ge b_*,\quad
 q\in\mathbb N,\quad q\ge2,\quad q^2\eta\le c_*,\quad \sum_iM_i^2\preceq\eta I_D,
\]
then $F\in C^2([0,1])$ and
\begin{equation}\label{RRS-eq:interpolation}
 F''(t)\ge-C_*bD\eta-C_*aqD(C_*q^3\eta)^q,
 \qquad 0\le t\le1.
\end{equation}
Moreover, $F'(0)=0$.
\end{lemma}
The regularity and inequality are~\cite[Lemmas 3.12--3.13]{AkbasSraMatrix}. The last identity follows directly: the derivative of the trace potential at a diagonal matrix is diagonal, so its inner product with $B_{\rm o}$ is zero. Lemma~\ref{RRS-lem:AS-interpolation} is used as an established analytic input; the parameter argument is given in full below.

\begin{proof}[Proof of Theorem~\ref{RRS-thm:hereditary}]
We may assume $v=1$ by homogeneity. Put
\[
 r=D/m,\qquad L=1+\log(2r),\qquad
 R=\kappa L^{3/2},\qquad M_i=B_i/R,\qquad
 \eta=\frac1{\kappa^2L^3}.
\]
Then $\sum_iM_i^2\preceq\eta I_D$. Fix $\delta=\min\{\tau,1\}/4$. Choose a constant $b_0\ge\max\{b_*,1\}$, then choose $a_0$ large enough, depending only on $\delta$, that
\begin{equation}\label{RRS-eq:a0}
 a_0\ge\max\{a_*,5,\sqrt{b_0}\},\qquad
 \frac{(C_*+2)b_0}{a_0}\le\delta.
\end{equation}
Choose $c\ge2$ sufficiently large that
\begin{equation}\label{RRS-eq:cchoice}
 e^2 4^{-c}\le1,\qquad 2a_0e^2 4^{-c}\le\delta.
\end{equation}
Finally choose $\kappa$ sufficiently large as specified below, and set
\begin{equation}\label{RRS-eq:parameters}
 q=\lceil cL\rceil,\qquad a=a_0\max\{1,r\eta\},\qquad b=b_0/a.
\end{equation}
The first four conditions of Lemma~\ref{RRS-lem:AS-interpolation} hold. Since $cL\le q\le2cL$,
\begin{equation}\label{RRS-eq:qeta}
 q^2\eta\le\frac{4c^2}{\kappa^2},\qquad
 C_*q^3\eta\le\frac{8C_*c^3}{\kappa^2}=:\theta.
\end{equation}
Increasing $\kappa$ makes $q^2\eta\le c_*$, $\eta\le1$, and $\theta<1$.

First consider the two terms in \eqref{RRS-eq:interpolation}. We have
\begin{equation}\label{RRS-eq:first-error}
 C_*bD\eta\le\frac{C_*b_0}{a_0}m\le\delta m.
\end{equation}
Also $ar\le a_0r^2\le a_0e^{2L}$, whence
\[
 \frac{C_*aqD(C_*q^3\eta)^q}{m}
 \le2C_*a_0cL\,e^{2L}\theta^{cL}.
\]
Choose $\kappa$ so large that $x=e^2\theta^c\le e^{-1}$ and $2C_*a_0cx\le\delta$. Since $Lx^L\le x$ for $L\ge1$, the last display is at most $\delta$. Consequently,
\begin{equation}\label{RRS-eq:Fsecond}
 F''(t)\ge-2\delta m.
\end{equation}

We next bound $F(0)$, where the matrix is diagonal. Write $B_0=\operatorname{diag}(\zeta_1,\ldots,\zeta_D)$ and $\Sigma=\operatorname{Cov}(\zeta)$. Then
\[
 \operatorname{Var}(\zeta_j)\le\eta,\qquad \tr\Sigma\le D\eta.
\]
For $|z|\le1/2$, the elementary estimates $\ell(z)\le2z^2$ and $r_q(z)\le 2\cdot4^{-q}/q$ imply
\[
 v_{q,a,b}(z)\le2bz^2+2a4^{-q}/q.
\]
On the event $E=\{\max_j|\zeta_j|\le1/2\}$, this gives
\[
 e^{F(0)}\ge e^{-2aD4^{-q}/q}\E[e^{-2b|\zeta|^2}\mathbf1_E].
\]
The Gaussian quadratic tilt has normalizing factor
\[
 \det(I+4b\Sigma)^{-1/2}\ge e^{-2b\tr\Sigma}\ge e^{-2bD\eta},
\]
and its covariance is $\Sigma^{1/2}(I+4b\Sigma)^{-1}\Sigma^{1/2}\preceq\Sigma$. The Gaussian strip correlation inequality therefore gives, when $2e^{-1/(8\eta)}\le1/2$,
\[
 \widetilde{\mathbb P}(E)\ge(1-2e^{-1/(8\eta)})^D
 \ge e^{-4D e^{-1/(8\eta)}}.
\]
The argument also applies to singular Gaussian covariance by a limit. Hence
\begin{equation}\label{RRS-eq:Fzero}
 F(0)\ge-2bD\eta-\frac{2aD}{q}4^{-q}-4D e^{-1/(8\eta)}.
\end{equation}
The first term is at most $\delta m$ by \eqref{RRS-eq:a0}. The second, divided by $m$, is at most
\[
 2a_0e^{2L}4^{-cL}\le2a_0e^2 4^{-c}\le\delta.
\]
Finally, for $\kappa\ge4$ and $L\ge1$,
\[
 4r e^{-1/(8\eta)}
 \le4e^{L-\kappa^2L^3/8}
 \le4e^{-\kappa^2/16}.
\]
Increasing $\kappa$ makes this at most $\delta$. Thus $F(0)\ge-3\delta m$. Equation~\eqref{RRS-eq:Fsecond} and $F'(0)=0$ give
\[
 F(1)\ge F(0)+\int_0^1(1-t)F''(t)\,dt\ge-4\delta m\ge-\tau m.
\]
The potential is nonnegative and infinite outside the operator-norm unit ball. Its partition function is therefore at most $\gamma_m\{\|B(x)\|_{\op}<1\}$. This proves \eqref{RRS-eq:hereditary}.
\end{proof}

\begin{corollary}[Dimension-sensitive simultaneous Gaussian event]\label{RRS-cor:general-gaussian}
For $A_i\in\Sym_d(\R)$, set $D=\max\{d,n\}$ and $L_{n,d}=1+\log(2D/n)$. There is a universal $c>0$ such that
\begin{equation}\label{RRS-eq:general-gaussian}
 \gamma_n(K_A)\ge\left(\frac{c}{L_{n,d}^{3/2}}\right)^n.
\end{equation}
\end{corollary}
\begin{proof}
Pad the matrices by zeros to dimension $D$, apply Lemma~\ref{RRS-lem:factor}, and use Theorem~\ref{RRS-thm:hereditary} with $m=n$, $v=1$ and $\tau=1$ for the normalized coefficients $U_i$. Lemma~\ref{RRS-lem:transfer} and the Gaussian scaling inequality \eqref{RRS-eq:scaling} give the result with $c=e^{-1}/(2\kappa(1))$.
\end{proof}

\subsection{Full rounding with exact mean}
At each active set we use one spectral body for every $p$. To prescribe the mean, fix an arbitrary linear functional and keep each partial-coloring increment in its kernel. The final rounding can only increase that functional. If the prescribed point lay outside the convex hull of the simultaneously good signs, a separating functional would contradict this construction. The proof below implements this argument with constants independent of the functional and of $p$.

\begin{lemma}[Square-function estimates across the Schatten scale]\label{RRS-lem:Sp-square}
For symmetric matrices $A_i$ and $V_S=(\sum_{i\in S}A_i^2)^{1/2}$,
\begin{equation}\label{RRS-eq:Sp-square}
 \spnorm{V_S}{p}\le
 \begin{cases}
 (\sum_{i\in S}\spnorm{A_i}{p}^{p})^{1/p},&1\le p\le2,\\
 (\sum_{i\in S}\spnorm{A_i}{p}^{2})^{1/2},&2\le p\le\infty.
 \end{cases}
\end{equation}
In particular, $\spnorm{V_S}{p}\le |S|^{a(p)}\max_i\spnorm{A_i}{p}$.
\end{lemma}
\begin{proof}
For $p\ge2$, use the triangle inequality in $S_{p/2}$, including the operator norm at $p=\infty$. For $1\le p<2$, put $t=p/2\in(0,1)$. Positive matrices satisfy $\tr(X+Y)^t\le\tr X^t+\tr Y^t$. To see this directly, use
\[
 z^t=\frac{\sin(\pi t)}{\pi}\int_0^\infty s^{t-1}\frac{z}{z+s}\,ds
\]
and the inverse-order inequalities $(X+Y+sI)^{-1}\preceq(X+sI)^{-1}$ and $(X+Y+sI)^{-1}\preceq(Y+sI)^{-1}$ after taking traces against $X$ and $Y$. Iterating gives
\[
 \tr\left(\sum_{i\in S}A_i^2\right)^{p/2}
 \le\sum_{i\in S}\tr|A_i|^p.
\]
The endpoint $p=2$ is equality.
\end{proof}

\begin{lemma}[Shifted partial coloring, Reis--Rothvoss]\label{RRS-lem:partial}
There is a universal $C_{\rm pc}$ such that the following holds. Let $K\subset\R^m$ be symmetric and convex with $\gamma_m(K)\ge e^{-m}$. Let $z\in[-1,1]^m$, and let $F\subset\R^m$ be a linear subspace with $\dim F\ge3m/4$. There is $x\in C_{\rm pc}K\cap F$ for which $z+x\in[-1,1]^m$ and at least $m/2$ coordinates of $z+x$ lie in $\{-1,1\}$.
\end{lemma}
This is~\cite[Theorem 6]{R19-RR} with $(\alpha,\beta,\gamma)=(1,3/4,1/4)$. An unbounded $K$ can be intersected with sufficiently large Euclidean balls and the measure exponent enlarged by an absolute constant, which only changes $C_{\rm pc}$.

\begin{theorem}[Simultaneous rounding in arbitrary matrix dimension]\label{RRS-thm:general-law}
Let $A_1,\ldots,A_n\in\Sym_d(\R)$, $y\in[-1,1]^n$, and
\[
 L_{n,d}=1+\log\frac{2\max\{d,n\}}n.
\]
There is a law on $\sigma\in\{-1,1\}^n$ with $\E\sigma=y$, supported on at most $n+1$ points, such that every outcome satisfies
\begin{equation}\label{RRS-eq:general-law}
 \spnorm{\sum_i(\sigma_i-y_i)A_i}{p}
 \le C L_{n,d}^{3/2}n^{a(p)}\max_i\spnorm{A_i}{p}
 \quad\text{for every }1\le p\le\infty.
\end{equation}
\end{theorem}
\begin{proof}
Pad by zeros so the matrix dimension is $D=\max\{d,n\}$. Fix an arbitrary $w\in\R^n$. We construct a signing satisfying \eqref{RRS-eq:general-law} and
\begin{equation}\label{RRS-eq:support-functional}
 \langle w,\sigma-y\rangle\ge0.
\end{equation}
Begin at $z=y$. At a stage with active set $S=\{i:|z_i|<1\}$ of size $m\ge4$, apply Lemma~\ref{RRS-lem:factor} to $(A_i)_{i\in S}$, obtaining $H_S,U_i^{S}$ with $H_S\preceq V_S$ and $\sum_{i\in S}(U_i^S)^2\preceq I_D$. Define
\[
 K_S=\left\{x\in\R^S:
 \left\lVert\sum_{i\in S}x_iU_i^S\right\rVert_{\op}
 \le\kappa(1)\left(1+\log\frac{2D}{m}\right)^{3/2}\right\}.
\]
Theorem~\ref{RRS-thm:hereditary} gives $\gamma_m(K_S)\ge e^{-m}$. The subspace
\[
 F_S=\{x\in\R^S:\langle w_S,x\rangle=0\}
\]
has dimension at least $m-1\ge3m/4$. Lemma~\ref{RRS-lem:partial} gives an increment $x\in C_{\rm pc}K_S\cap F_S$ that fixes at least half the active coordinates. Lemmas~\ref{RRS-lem:transfer} and~\ref{RRS-lem:Sp-square} give, on the same increment for all $p$,
\begin{equation}\label{RRS-eq:increment}
 \spnorm{\sum_{i\in S}x_iA_i}{p}
 \le 2C_{\rm pc}\kappa(1)
 \left(1+\log\frac{2D}{m}\right)^{3/2}
 m^{a(p)}M_p,\qquad M_p=\max_i\spnorm{A_i}{p}.
\end{equation}
The increment preserves $\langle w,z\rangle$ exactly.

Continue until at most three coordinates remain active. Round each remaining coordinate toward the sign of $w_i$, using either sign when $w_i=0$. This final operation increases $\langle w,z\rangle$ and has Schatten-$p$ error at most $6M_p$ for every $p$. Thus \eqref{RRS-eq:support-functional} holds.

To sum \eqref{RRS-eq:increment}, let $m_j$ be the nonzero active sizes and put $b_j=\lfloor\log_2(n/m_j)\rfloor$. Each partial-coloring stage at least halves $m_j$, so the $b_j$ are distinct increasing nonnegative integers. Writing $L_0=L_{n,d}$, we have
\[
 m_j^{a(p)}\le n^{a(p)}2^{-b_j/2},\qquad
 1+\log\frac{2D}{m_j}\le L_0+(b_j+1)\log2\le L_0(b_j+2).
\]
Since $a(p)\ge1/2$ for every $p$, it follows that
\begin{align*}
 \sum_jm_j^{a(p)}\left(1+\log\frac{2D}{m_j}\right)^{3/2}
 &\le n^{a(p)}L_0^{3/2}\sum_{b=0}^\infty2^{-b/2}(b+2)^{3/2}\\
 &\le n^{a(p)}L_0^{3/2}\sum_{b=0}^\infty2^{-b/2}(b+2)^2
 <95n^{a(p)}L_0^{3/2}.
\end{align*}
The last sum equals $(4-3r+r^2)/(1-r)^3$ for $r=2^{-1/2}$ and is less than $95$. This proves \eqref{RRS-eq:general-law} with a universal constant, uniformly in $w$ and $p$.

Let $G_y$ be the finite set of signings satisfying \eqref{RRS-eq:general-law} for all $p$. The construction shows that, for every $w$, there exists $\sigma\in G_y$ with $\langle w,\sigma\rangle\ge\langle w,y\rangle$. Finite-dimensional separation gives $y\in\operatorname{conv}(G_y)$. A convex representation of $y$ defines the required law; Carath\'eodory's theorem reduces its support to at most $n+1$ points. For $y=0$, the set $G_0$ is nonempty and symmetric, so the uniform law on $\{\sigma,-\sigma\}$ for any $\sigma\in G_0$ is a symmetric choice with support size two.
\end{proof}

Taking $d=n$ in Theorem~\ref{RRS-thm:general-law} proves Theorem~\ref{RRS-thm:square-law}, since $L_{n,n}=1+\log2$ is an absolute constant. Preserving a separating functional at every stage places the prescribed mean in the convex hull of the full signings. The polynomial-bit sampling results of Table~\ref{R9-tab:runtime} concern their separately specified rounding models.

\subsection{Optimal exponents and nonsymmetric matrices}\label{RRS-sec:sharpness}
\begin{proposition}[Optimal orders in the square case]\label{RRS-prop:sharp}
The powers $n^{1/p}$ for $1\le p\le2$ and $\sqrt n$ for $2\le p\le\infty$ in Theorem~\ref{RRS-thm:square-law} cannot be improved uniformly in $n$, even for diagonal matrices and $y=0$.
\end{proposition}
\begin{proof}
For $1\le p\le2$, use $A_i=e_ie_i^{\mathsf T}$. Then $\spnorm{A_i}{p}=1$, while every signed sum has Schatten-$p$ norm $n^{1/p}$.

For $2\le p\le\infty$, take a Sylvester Hadamard matrix $Q\in\{-1,1\}^{n\times n}$, for $n$ a power of two, and let $A_i$ be diagonal with diagonal equal to its $i$th column. For every signing $\sigma$, $\lVert Q\sigma\rVert_2=n$. Hence
\[
 \spnorm{\sum_i\sigma_iA_i}{p}
 =\lVert Q\sigma\rVert_p
 \ge n^{1/p-1/2}\lVert Q\sigma\rVert_2
 =\sqrt n\,n^{1/p}
 =\sqrt n\,\max_i\spnorm{A_i}{p}.
\]
The same calculation applies at $p=\infty$ with $n^{1/p}=1$.
\end{proof}

\begin{corollary}[Rectangular real and complex matrices]\label{RRS-cor:rectangular}
Let $A_1,\ldots,A_n$ be real or complex $a\times b$ matrices and $y\in[-1,1]^n$. The conclusion of Theorem~\ref{RRS-thm:general-law} holds with
\[
 L=1+\log\frac{2\max\{n,2(a+b)\}}n
\]
in place of $L_{n,d}$. In particular, for $n$ arbitrary $n\times n$ real or complex matrices, one exact-mean law satisfies \eqref{RRS-eq:square-law} for the whole Schatten scale with a universal constant.
\end{corollary}
\begin{proof}
Use the Hermitian dilations
\[
 \mathcal A_i=\begin{pmatrix}0&A_i\\A_i^*&0\end{pmatrix}.
\]
For finite $p$, their Schatten norms are $2^{1/p}\spnorm{A_i}{p}$, and their operator norms equal those of $A_i$. For complex matrices, realification doubles every singular-value multiplicity once more and produces real symmetric matrices of dimension $2(a+b)$. These multiplicity factors appear identically on both sides of \eqref{RRS-eq:general-law} and cancel. A common law on the coefficient signs therefore gives the claimed result for every $p$.
\end{proof}


\section{Trace-class quotients and Reis's volume-ratio conjecture}\label{R27-sec:volume}
For a centrally symmetric convex body $L\subset\R^r$ with nonempty interior, let $E_L$ be its John ellipsoid, the ellipsoid of largest volume contained in $L$. Its \emph{volume ratio} is
\[
 \vr(L)=\left(\frac{|L|}{|E_L|}\right)^{1/r}.
\]
Here $|\cdot|$ denotes $r$-dimensional Lebesgue measure. This quantity is invariant under invertible linear maps. Let $B_{S_1}^{a\times b}$ be the unit ball of the nuclear norm on real $a\times b$ matrices, with the Hilbert--Schmidt inner product. Put $k=\min\{a,b\}$.

\begin{theorem}[Volume ratios of trace-class quotients]\label{SVR-thm:main}
There are universal constants $c_0,c_1>0$ such that every surjective linear map $M:\R^{a\times b}\to\R^r$ satisfies
\begin{equation}\label{SVR-eq:main}
 \vr\bigl(MB_{S_1}^{a\times b}\bigr)
 \le \sqrt{c_0+c_1\frac{\min\{a,b\}}r}.
\end{equation}
Consequently, every $r$-dimensional quotient of $S_1^{a\times b}$ has uniformly bounded volume ratio whenever $r\ge\min\{a,b\}$. For $r\ge\theta\min\{a,b\}$, the bound depends only on $\theta>0$.
\end{theorem}

\begin{corollary}[Reis's volume-ratio conjecture]\label{SVR-cor:Reis}
There is a universal constant $C$ such that, for every positive integer $d$ and every surjective linear map $M:\R^{d\times d}\to\R^d$,
\begin{equation}\label{SVR-eq:Reis}
 \vr\bigl(MB_{S_1}^{d\times d}\bigr)\le C.
\end{equation}
\end{corollary}

Corollary~\ref{SVR-cor:Reis} answers Conjecture~1 in Section~7 of Reis's \emph{Optimal Vector Balancing for Zonotopes}~\cite{R27-ReisZonotopes}. Surjectivity states the full-dimensional hypothesis in that paper's definition of volume ratio. Theorem~\ref{SVR-thm:main} also specifies the dependence on the intrinsic image dimension when the matrix size and image dimension differ. The bound is uniform over all surjective linear maps $M$.

The Fisher bound \eqref{R3-eq:external-matrix-Fisher} also controls the volume needed to support the matrix reference. John contact points normalize both matrix variances and retain the smaller matrix dimension through their rank bounds. Differential entropy then converts the determinant of Fisher information into a volume estimate, which polarity transfers to the quotient. Thus the Akbas--Sra density used for signing gives a geometric theorem for every surjective map.

The determinant estimate retains the contribution of each coefficient direction. Write $K=L^\circ$ for the polar body, and let $E_K^{\rm out}$ be its minimum-volume containing ellipsoid. Define
\[
 \ovr(K)=\left(\frac{|E_K^{\rm out}|}{|K|}\right)^{1/r}.
\]
After putting $L=MB_{S_1}^{a\times b}$ in John position, so that $E_L=B_2^r$, set
\begin{equation}\label{SVR-eq:coefficients}
 A_i=M^*e_i,\qquad
 \mathcal A_i=\begin{pmatrix}0&A_i\\ A_i^{\mathsf T}&0\end{pmatrix},\qquad
 G_{ij}=\tr(\mathcal A_i\mathcal A_j)=2\tr(A_i^{\mathsf T}A_j).
\end{equation}
The adjoint here is computed after the change to John coordinates. Set
\[
 v=\left\|\sum_{i=1}^r\mathcal A_i^2\right\|_{\rm op}.
\]

\begin{theorem}[Determinant bound]\label{SVR-thm:determinant}
Let $\eta>0$ and $\Lambda>0$ be the universal constants in the Fisher estimate of Lemma~\ref{SVR-lem:Fisher}. In the notation above,
\begin{align}
 \vr(L)&\le\ovr(K)
 \le\det\left(\frac{v}{\eta r}I_r+\frac{\Lambda}{r}G\right)^{1/(2r)},
 \label{SVR-eq:determinant}\\
 v&\le r,\qquad \tr G\le 2rk.
 \label{SVR-eq:trace-bounds}
\end{align}
In particular, \eqref{SVR-eq:main} holds with $c_0=\eta^{-1}$ and $c_1=2\Lambda$.
\end{theorem}

\subsection{John contacts control both matrix variances}
The polar of $L=MB_{S_1}^{a\times b}$ has the exact description
\begin{equation}\label{SVR-eq:polar}
 K=\{x\in\R^r:\|M^*x\|_{\rm op}\le1\}.
\end{equation}
Indeed, the support function of $L$ is the operator norm of $M^*x$, by duality of the nuclear and operator norms. Surjectivity of $M$ makes $M^*$ injective, so $K$ is a bounded, full-dimensional convex body.

\begin{lemma}[Contact normalization]\label{SVR-lem:John}
In John coordinates for $L$, the matrices in \eqref{SVR-eq:coefficients} satisfy
\begin{equation}\label{SVR-eq:rowcolumn}
 \sum_i A_iA_i^{\mathsf T}\preceq rI_a,
 \qquad
 \sum_i A_i^{\mathsf T}A_i\preceq rI_b,
 \qquad
 \sum_i\|A_i\|_{\rm F}^2\le rk.
\end{equation}
Their symmetric dilations therefore satisfy \eqref{SVR-eq:trace-bounds}.
\end{lemma}
\begin{proof}
John's contact decomposition gives unit vectors $u_j\in L\cap K\cap S^{r-1}$ and positive weights $c_j$ with
\begin{equation}\label{SVR-eq:John}
 \sum_jc_ju_ju_j^{\mathsf T}=I_r,
 \qquad \sum_jc_j=r.
\end{equation}
We use the symmetric form recorded in~\cite[Theorem 6]{R27-ReisZonotopes}. Put $T_j=M^*u_j=\sum_i(u_j)_iA_i$. Since $u_j\in K$, equation~\eqref{SVR-eq:polar} gives $\|T_j\|_{\rm op}\le1$. Expanding the contact identity yields
\begin{align*}
 \sum_jc_jT_jT_j^{\mathsf T}
 &=\sum_{i,\ell}\left(\sum_jc_j(u_j)_i(u_j)_\ell\right)A_iA_\ell^{\mathsf T}
 =\sum_i A_iA_i^{\mathsf T},\\
 \sum_jc_jT_j^{\mathsf T}T_j&=\sum_i A_i^{\mathsf T}A_i.
\end{align*}
Both variance bounds follow from $T_jT_j^{\mathsf T}\preceq I_a$, $T_j^{\mathsf T}T_j\preceq I_b$, and $\sum_jc_j=r$. Moreover,
\[
 \sum_i\|A_i\|_{\rm F}^2
 =\sum_jc_j\|T_j\|_{\rm F}^2
 \le\sum_jc_j\,\rank(T_j)\le rk.
\]
Finally,
\[
 \sum_i\mathcal A_i^2
 =\begin{pmatrix}\sum_iA_iA_i^{\mathsf T}&0\\0&\sum_iA_i^{\mathsf T}A_i\end{pmatrix},
 \qquad \tr G=2\sum_i\|A_i\|_{\rm F}^2.
\]
This proves the claimed bounds.
\end{proof}

The normalization uses the complete contact identity. It controls the row and column variances simultaneously, although the individual coefficient matrices need not commute. The last line of the proof is also why the final estimate depends on $\min\{a,b\}$ rather than the larger dilation size $a+b$.

\subsection{Fisher information determines a support-volume bound}
For a probability density $p$ on $\R^r$, write
\[
 J(p)=\int_{\{p>0\}}\frac{\nabla p\,\nabla p^{\mathsf T}}p,
 \qquad h(p)=-\int p\log p.
\]
The first quantity is its Fisher information matrix; the second is its differential entropy.

A lower entropy bound limits how small the support can be. The matrix form matters: an affine change of variables turns the trace Fisher inequality into a determinant bound, retaining every coefficient direction. We use the density in \eqref{R3-eq:external-matrix-Fisher}, with $\Lambda$ denoting the constant called $L$ there; its required regularity is part of the same Akbas--Sra input.
\begin{lemma}[Regularity of the matrix reference; Akbas--Sra]\label{SVR-lem:Fisher}
For real symmetric $B_1,\ldots,B_r$ with
$\|\sum_iB_i^2\|_{\rm op}\le\eta/t^2$, let $p_t$ be the barrier density defined before \eqref{R3-eq:external-matrix-Fisher}, and put $H_{ij}=\tr(B_iB_j)$. Its extension by zero across $\|\sum_i x_iB_i\|_{\rm op}=1$ is $C^2$, belongs to $W^{1,1}(\R^r)$, and has finite Fisher information with
\begin{equation}\label{SVR-eq:Fisher}
 J(p_t)\preceq t^{-2}I_r+\Lambda H.
\end{equation}
\end{lemma}
This is \cite[Lemma 3.3]{AkbasSraBSB}; $\Lambda=64C_R$ is an admissible choice with the constants of its Lemma~2.4, derived from \cite[Proposition 3.8]{AkbasSraMatrix}.

\begin{lemma}[Entropy and the determinant of Fisher information]\label{SVR-lem:entropy}
Suppose $p$ is a bounded, compactly supported probability density in $W^{1,1}(\R^r)$ with finite positive-definite Fisher information matrix. Then
\begin{equation}\label{SVR-eq:Stam-matrix}
 h(p)\ge\frac r2\log(2\pi e)-\frac12\log\det J(p).
\end{equation}
If $p$ is supported on a measurable set $\Omega$ of finite volume, then
\begin{equation}\label{SVR-eq:volume-Fisher}
 |\Omega|\ge\frac{(2\pi e)^{r/2}}{\sqrt{\det J(p)}}.
\end{equation}
\end{lemma}
\begin{proof}
We record the entropy--Fisher calculation, including its matrix dependence. The Gaussian logarithmic Sobolev inequality~\cite{Gross1975} for $\gamma_\sigma=N(0,\sigma^2I_r)$ gives
\[
 D(p\|\gamma_\sigma)\le\frac{\sigma^2}{2}I(p\|\gamma_\sigma).
\]
The finite-Fisher Sobolev form follows by approximation from the smooth form. Bounded compact support makes $h(p)$ and the second moment finite. Integration by parts gives $\int x\cdot\nabla p=-r$. Consequently,
\begin{align*}
 D(p\|\gamma_\sigma)
 &=-h(p)+\frac r2\log(2\pi\sigma^2)
                 +\frac{\mathbb E\|X\|^2}{2\sigma^2},\\
 I(p\|\gamma_\sigma)
 &=\tr J(p)-\frac{2r}{\sigma^2}
                 +\frac{\mathbb E\|X\|^2}{\sigma^4}.
\end{align*}
The second-moment terms cancel. Choosing $\sigma^2=r/\tr J(p)$ yields
\begin{equation}\label{SVR-eq:Stam-trace}
 h(p)\ge\frac r2\log\left(\frac{2\pi e\,r}{\tr J(p)}\right).
\end{equation}
Let $Y=J(p)^{1/2}X$. Its Fisher information matrix is $I_r$, while
\[
 h(Y)=h(X)+\tfrac12\log\det J(p).
\]
Apply \eqref{SVR-eq:Stam-trace} to $Y$ to obtain \eqref{SVR-eq:Stam-matrix}. Finally, comparison with the uniform density on $\Omega$ gives $h(p)\le\log|\Omega|$, which proves \eqref{SVR-eq:volume-Fisher}.
\end{proof}

\begin{proposition}[Volume of an operator-norm section]\label{SVR-prop:sectionvolume}
Let $B_1,\ldots,B_r$ be linearly independent real symmetric matrices and put
\[
 \Omega=\left\{x:\left\|\sum_i x_iB_i\right\|_{\rm op}<1\right\},
 \qquad v=\left\|\sum_i B_i^2\right\|_{\rm op},
 \qquad H_{ij}=\tr(B_iB_j).
\]
Then
\begin{equation}\label{SVR-eq:sectionvolume}
 |\Omega|\ge
 \frac{(2\pi e)^{r/2}}
 {\sqrt{\det\bigl((v/\eta)I_r+\Lambda H\bigr)}}.
\end{equation}
\end{proposition}
\begin{proof}
Linear independence makes $\Omega$ bounded and gives $v>0$. Choose $t^2=\eta/v$ in Lemma~\ref{SVR-lem:Fisher}. The resulting density has compact support in $\overline\Omega$, is bounded, and satisfies
\[
 J(p_t)\preceq(v/\eta)I_r+\Lambda H.
\]
Its Fisher information is positive definite. Indeed, a zero quadratic form in a nonzero direction would give a zero weak derivative in that direction, so an integrable compactly supported density would be constant along almost every parallel line and hence vanish. Lemma~\ref{SVR-lem:entropy}, monotonicity of determinant on positive-definite matrices, and the zero volume of the boundary of a convex body give \eqref{SVR-eq:sectionvolume}.
\end{proof}

This is the step that converts the reference into a geometric statement. The directional bounds used for signing are entries of $J(p_t)$. Keeping the determinant instead bounds the volume required to support the whole density. The argument retains every coefficient direction before taking the trace.

\subsection{Proof of the quotient theorem and the dual statement}
\begin{proof}[Proof of Theorems~\ref{SVR-thm:main} and~\ref{SVR-thm:determinant}]
Put $L=MB_{S_1}^{a\times b}$ in John position and form $K$, $\mathcal A_i$, $v$, and $G$ as above. Symmetric dilation preserves the operator norm:
\[
 \left\|\sum_i x_i\mathcal A_i\right\|_{\rm op}=\|M^*x\|_{\rm op}.
\]
Thus Proposition~\ref{SVR-prop:sectionvolume} applies to $\operatorname{int}K$ and gives
\begin{equation}\label{SVR-eq:K-volume}
 |K|^{1/r}\ge
 \frac{\sqrt{2\pi e}}
      {\det\bigl((v/\eta)I_r+\Lambda G\bigr)^{1/(2r)}}.
\end{equation}
The polar of a John ellipsoid is the L\"owner ellipsoid of the polar body, so $E_K^{\rm out}=B_2^r$. Also
\begin{equation}\label{SVR-eq:ballvolume}
 |B_2^r|^{1/r}\le\sqrt{\frac{2\pi e}{r}}.
\end{equation}
For completeness, integrate $e^{-t\|x\|^2}\ge e^{-t}$ on $B_2^r$, and optimize $|B_2^r|\le e^t(\pi/t)^{r/2}$ at $t=r/2$. Combining \eqref{SVR-eq:K-volume} and \eqref{SVR-eq:ballvolume} proves
\[
 \ovr(K)\le\det\left(\frac{v}{\eta r}I_r+\frac\Lambda rG\right)^{1/(2r)}.
\]
The symmetric Blaschke--Santal\'o inequality~\cite{R27-BianchiKelly} gives $|L||K|\le|B_2^r|^2$, and therefore
\[
 \vr(L)=\left(\frac{|L|}{|B_2^r|}\right)^{1/r}
 \le\left(\frac{|B_2^r|}{|K|}\right)^{1/r}=\ovr(K).
\]
This proves the determinant estimate. Finally, apply the arithmetic--geometric mean inequality to its eigenvalues and use Lemma~\ref{SVR-lem:John}:
\begin{align*}
 \vr(L)^2
 &\le\frac{v}{\eta r}+\frac\Lambda{r^2}\tr G\\
 &\le\eta^{-1}+2\Lambda\frac{k}{r}.
\end{align*}
This is \eqref{SVR-eq:main}. Setting $a=b=r=d$ proves Corollary~\ref{SVR-cor:Reis}.
\end{proof}

\begin{corollary}[Symmetric matrices and spectrahedra]\label{SVR-cor:spectrahedra}
For a surjective map $M:\operatorname{Sym}_d\to\R^r$,
\[
 \vr\bigl(MB_{S_1}^{\operatorname{Sym}_d}\bigr)
 \le\sqrt{\eta^{-1}+\Lambda d/r}.
\]
Let $\mathcal K=\{x\in\R^r:\|\sum_i x_iB_i\|_{\rm op}\le1\}$ be bounded, where the $B_i$ are symmetric $d\times d$ matrices. Then
\begin{equation}\label{SVR-eq:ovr}
 \ovr(\mathcal K)\le\sqrt{\eta^{-1}+\Lambda d/r}.
\end{equation}
In particular, centrally symmetric spectrahedra of this form with matrix size proportional to their dimension have bounded outer volume ratio.
\end{corollary}
\begin{proof}
For symmetric matrices no dilation is required. John contacts give $\sum_iA_i^2\preceq rI_d$ and $\tr G\le rd$, so the same proof removes the factor two. The second claim applies that proof to the polar of $\mathcal K$, whose support function is $\|\sum_i x_iB_i\|_{\rm op}$. Equivalently, it is the image of the symmetric nuclear-norm ball under the adjoint of $x\mapsto\sum_i x_iB_i$.
\end{proof}

The spectrahedron is described by the two linear matrix inequalities $I_d\pm\sum_i x_iB_i\succeq0$. The conclusion compares the volume of the body with the volume of its minimum containing ellipsoid. For fixed $\theta>0$, the estimates hold for every image or section with $r\ge\theta d$.

\subsection{A second proof of the named conjecture with explicit constants}
The shared small-ball input, Lemma~\ref{RRS-lem:AS-fixed}, gives a second proof with explicit constants. For $C_i\in\operatorname{Sym}_D$ and $\sum_i C_i^2\preceq I_D$, its quantitative form is
\begin{equation}\label{SVR-eq:smallball}
 \mathbb P\left\{\left\|\sum_i g_iC_i\right\|_{\rm op}<4\right\}\ge\beta^D,
 \quad
 \beta=\frac34\left(\frac{1020\sqrt3-1671}{128}\right)^8>0,
\end{equation}
where the $g_i$ are independent standard Gaussians. The displayed exact choice follows from its rank-one pinching proof; $-\log\beta=2.614861\ldots$.

\begin{proposition}[Explicit dimension-comparable estimate]\label{SVR-prop:explicit}
In the rectangular setting,
\begin{equation}\label{SVR-eq:explicit}
 \vr(L)\le\ovr(L^\circ)\le4\sqrt e\,\beta^{-(a+b)/r}.
\end{equation}
Consequently, Corollary~\ref{SVR-cor:Reis} holds with $C<1232$. For symmetric $d\times d$ matrices with $r=d$, the same argument gives $C=4\sqrt e\,\beta^{-1}$.
\end{proposition}
\begin{proof}
Use the John-normalized dilations and apply \eqref{SVR-eq:smallball} to $C_i=\mathcal A_i/\sqrt r$, whose variance is at most the identity by Lemma~\ref{SVR-lem:John}. With $D=a+b$,
\[
 \gamma_r(4\sqrt r\,K)\ge\beta^D.
\]
The standard Gaussian density is at most $(2\pi)^{-r/2}$, so
\[
 |K|^{1/r}\ge\frac{\sqrt{2\pi}}{4\sqrt r}\,\beta^{D/r}.
\]
Equation~\eqref{SVR-eq:ballvolume} and Blaschke--Santal\'o prove \eqref{SVR-eq:explicit}. For $a=b=r=d$, the constant is $4\sqrt e\,\beta^{-2}=1231.535442\ldots<1232$. Symmetric matrices require no dilation and give the last assertion.
\end{proof}

The determinant argument proves the stronger dependence in Theorem~\ref{SVR-thm:main}: its cost is polynomial in $k/r$ and is independent of the longer matrix dimension. Proposition~\ref{SVR-prop:explicit} gives a second complete deduction of the dimension-comparable conjecture from a different stated analytic input.

\subsection{Intrinsic image dimension}
The full-dimensional convention matters. A rank-deficient map is covered by Theorem~\ref{SVR-thm:main} using its actual rank $r$. A dimension-free intrinsic-volume claim for arbitrary ranks would be false. To see this, choose $d=2^{r-1}$ and vectors $v_1,\ldots,v_d$ containing one representative of each antipodal pair of vertices of $[-1,1]^r$. The map
\[
 M(X)=\sum_{j=1}^d X_{jj}v_j
\]
sends the $d\times d$ nuclear-norm ball exactly onto $[-1,1]^r$. Indeed, its diagonal vectors fill the $\ell_1^d$ unit ball: $\sum_j|X_{jj}|\le\|X\|_*$ by nuclear/operator duality, and every $\ell_1$ vector is realized by a diagonal matrix. The image is therefore $\operatorname{conv}\{\pm v_j\}$. Its John ellipsoid is $B_2^r$, and its volume ratio grows like $\sqrt r$. The main theorem records the necessary matrix-size/image-dimension dependence.


\Needspace{12\baselineskip}
\part{Deterministic realization and transport}
\label{R11-part:transport}
We now realize a chosen law as a deterministic function of an atomless source. The allocation must retain the entire output measure and every prescribed conditional moment. Directional, supermartingale and component constraints require different local allocations.

\section{Allocating an entire output law}
\label{R11-sec:allocation}
The Dvoretzky--Wald--Wolfowitz theorem preserves finitely many integrals while eliminating randomization~\cite{DWW}; countable partitions require further allocation~\cite{Edwards1987,KhanRath2009}. Here a deterministic map must reproduce every measurable subset of a continuous output law, as well as the conditional moments. We keep that measure intact and partition its mass according to the finite-dimensional mean geometry.

In a small output cell, every required mean has barycentric coordinates in a simplex of available source means. These coordinates split the entire target measure into shares. A mean-independent factor realizes each share without altering its source mean. The interval example makes the allocation explicit before the general lemma. Cellwise construction preserves the source--cell masses and gives $W_\infty$ approximation of the pair law. For countable targets, fixed label masses control integrable Borel-cost tails, allowing an extreme minimizing kernel to be a partition. Transport and deterministic dynamics use these exact-law conclusions.

\subsection{Countable purification and deterministic reference laws}
\label{FP-sec:purification}
For probability measures $\mu,\nu$ on $\R^d$ with finite first moments,
$\M(\mu,\nu)$ denotes their martingale couplings:
$X\sim\mu$, $Y\sim\nu$ and $\E[Y\mid X]=X$.
A coupling is \emph{backward Monge} if $X=T(Y)$ for a Borel map $T$;
write $\MM(\mu,\nu)$ for this subclass.
The distance $W_p$ is the Wasserstein distance on the product Euclidean
space; $W_\infty$ is allowed to take the value $+\infty$.

The following formulation makes explicit the finite number of conditions attached to each output. The underlying purification theorem is classical, including countable extensions \cite{DWW,Edwards1987,KhanRath2009}.

\begin{theorem}[Countable purification with prescribed moments]\label{FP-thm:purification}\label{SP-lem:purify}\label{SP-sec:purification}
Let $(S,\nu)$ be an atomless standard probability space, and let $J$ be finite or countable. Suppose $q_j:S\to[0,1]$ are measurable, $\sum_jq_j=1$, and, for each $j$, finitely many functions $f_{j,r}\in L^1(\nu)$ are specified. Also prescribe measurable allowed sets $A_j$ with $q_j=0$ off $A_j$. There is a measurable partition $(E_j)_{j\in J}$, with $E_j\subseteq A_j$ modulo null sets, such that
\[
 \nu(E_j)=\int q_j\,d\nu,
 \qquad
 \int_{E_j}f_{j,r}\,d\nu=\int q_jf_{j,r}\,d\nu
 \quad\text{for all }j,r.
\]
\end{theorem}
\begin{proof}
Put $p_j=\int q_j\,d\nu$. In the product of weak-star compact unit balls of $L^\infty(\nu)$, consider the allocations $(a_j)$ satisfying nonnegativity, all finite inequalities $\sum_{j\in F}a_j\le1$, the masses $\int a_j=p_j$, all prescribed moments, and $a_j=0$ off $A_j$. These constraints are weak-star closed: inequalities can be tested against nonnegative $L^1$ functions, and every moment is paired with an $L^1$ function. The resulting set is nonempty, compact and convex. Since $\sum_jp_j=1$, monotone convergence shows that $\sum_ja_j=1$ almost everywhere for every allocation in this set. Thus no mass escapes to infinitely many labels.

Choose an extreme point by the Krein--Milman theorem. If it is not a partition, countability gives labels $j\ne k$, a number $\eta>0$, and a set $E$ of positive measure on which $a_j,a_k\ge\eta$. Collect the constant function and all moment functions for these two labels into a finite vector $F$ with $L$ coordinates. Split $E$ into $L+1$ sets of positive measure. Linear dependence of their $F$-integrals gives a nonzero bounded simple function $h$, supported on $E$, with $\int hF\,d\nu=0$. Scale it so that $\|h\|_\infty\le\eta$. Replacing $(a_j,a_k)$ by $(a_j+h,a_k-h)$ or $(a_j-h,a_k+h)$ preserves all constraints, including the sum and allowed sets. Both allocations are nonnegative; all other coordinates are unchanged. Their midpoint is the purported extreme point, a contradiction. Therefore each $a_j$ is an indicator and these indicators form the required partition.
\end{proof}

\begin{corollary}[Backward Monge realization of a discrete target]\label{FP-cor:Monge-existence}
Let $\mu$ be countably supported on $\R^d$, let $\nu$ be atomless with finite first moment, and suppose $\mu\cx\nu$. There is a Borel map $T$ with
\[
 T_\#\nu=\mu,
 \qquad \E[Y\mid T(Y)]=T(Y),\quad Y\sim\nu.
\]
One can also preserve finitely many additional reference moments separately at each target atom, whenever a coupling satisfying them exists.
\end{corollary}
\begin{proof}
Take a martingale coupling, write $\mu=\sum_jp_j\delta_{x_j}$, and disintegrate backwards to obtain $q_j(y)=P(X=x_j\mid Y=y)$. Apply Theorem~\ref{FP-thm:purification} with the coordinate functions of $y$. On the resulting set $E_j$, put $T(y)=x_j$. Then $\nu(E_j)=p_j$ and $\int_{E_j}y\,d\nu=p_jx_j$, which are exactly the two asserted identities.
\end{proof}

\begin{theorem}[$W_\infty$ density for countable initial marginals]\label{FP-thm:discrete-density}\label{R9-intro:transport}
Under the hypotheses of Corollary~\ref{FP-cor:Monge-existence}, every $\pi\in\M(\mu,\nu)$ and every $\eta>0$ admit a backward Monge martingale coupling $\pi_T$ with
\[
 W_\infty(\pi,\pi_T)\le\eta.
\]
The same conclusion holds for supermartingale couplings when $d=1$ and $\mu$ is countably supported. Measurable pointwise allowed sets for each atom can be retained in either assertion. For any fixed finite family of bounded measurable costs $c_r(x,y)$, the approximating map can also satisfy
\[
 \int c_r\,d\pi_T=\int c_r\,d\pi\quad\text{for every }r.
\]
\end{theorem}
This proves the countably supported initial-marginal clause of
Nutz--Wang--Zhang, Conjecture~5.1 \cite{NWZ2024}. The conclusion holds
in every finite dimension and retains the stated finite costs exactly.
\begin{proof}
Partition $\R^d$ into countably many half-open cubes $C_\ell$ of diameter at most $\eta$. On each nonzero restriction $\nu|_{C_\ell}$, purify the backward kernel $q_j$ while preserving the mass and coordinate moments separately for each $j$. Combine the partitions over $\ell$. Summing the coordinate moments preserves each target atom's conditional mean, hence also any supermartingale inequality satisfied by that mean. For each pair $(j,\ell)$, the old and new restrictions have the same mass, their first coordinate equals $x_j$, and their second coordinates lie in $C_\ell$. Couple these restrictions arbitrarily and sum the couplings. The product-space distance is at most $\eta$ almost surely. Allowed sets are included in each local purification. To retain the costs exactly, add the finitely many functions $c_r(x_j,y)$ to the moments for label $j$ on each source cell. Summing the preserved integrals over cells and labels proves the assertion.
\end{proof}

\begin{theorem}[A deterministic optimizer for measurable costs]
\label{FP-prop:extreme}\label{R9-thm:Borel-cost}
Let $\mu=\sum_jp_j\delta_{x_j}$ be countably supported on $\R^d$ and
let $\nu$ be atomless, with $\mu\cx\nu$ and finite first moments.
Every extreme point of $\M(\mu,\nu)$ is backward Monge.
Moreover, every Borel cost $c$ satisfying
\begin{equation}\label{R9-eq:cost-envelope}
 |c(x_j,y)|\le a_j+b(y),\qquad
 a_j\ge0,\quad \sum_jp_ja_j<\infty,\quad 0\le b\in L^1(\nu),
\end{equation}
has a backward Monge minimizer over $\M(\mu,\nu)$.
In particular this holds for every bounded measurable cost.
Measurable allowed sets and finitely many integrable moments for each
label may be imposed whenever the resulting feasible set is nonempty.
\end{theorem}
\begin{proof}
Disintegrate backwards and use the weak-star compact allocation set
in Theorem~\ref{FP-thm:purification}. The overlap perturbation for two
labels preserves their masses, means, allowed sets and all other
specified finite moment lists. Hence every extreme allocation is a
partition, also after those extra constraints are imposed.

It remains to show that the cost is continuous in this topology.
For each fixed label, $c(x_j,\cdot)\in L^1(\nu)$, so
$\int c(x_j,y)a_j^{\rm ker}(y)\,d\nu(y)$ is weak-star continuous in
the corresponding kernel coordinate. Finite sums are continuous.
For the tail, put $q_J=\sum_{j>J}a_j^{\rm ker}$; then
$0\le q_J\le1$ and $\int q_J\,d\nu=\sum_{j>J}p_j$.
For any $B>0$, uniformly over feasible kernels,
\[
 \sum_{j>J}\int |c(x_j,y)|a_j^{\rm ker}(y)\,d\nu
 \le \sum_{j>J}p_ja_j
       +B\sum_{j>J}p_j+\int_{\{b>B\}}b\,d\nu.
\]
First choose $B$ large and then $J$ large. The tail tends uniformly to
zero. The objective is therefore a continuous affine function on a
compact convex set. Its minimizing face has an extreme point, which
is an extreme feasible allocation and hence a deterministic decoder.
The same overlap argument without extra constraints proves the stated
extreme-point description of $\M(\mu,\nu)$.
\end{proof}
The countable target makes the stronger cost class possible: after the
first finitely many labels, the fixed remaining mass controls the whole
cost tail. No continuity or semicontinuity of $c$ is required.

\begin{corollary}[Deterministic realization of the joint signing law]\label{FP-cor:R8-deterministic}
Suppose Theorem~\ref{MAIN-signing} or its prescribed-mean extension yields a law of $\sigma\in\{-1,1\}^n$ and an atomless integrable reference $R$ with
\[
 \E[R\mid\sigma]=(A\sigma,\sigma).
\]
There is a measurable function $F$ of the reference alone for which
\[
 \law(F(R))=\law(\sigma),\qquad
 \E[R\mid F(R)]=(AF(R),F(R)).
\]
All hard constraints and all information bounds that depend only on the signing law are preserved. If $R=(G,H)$ has independent reference blocks, their full joint distribution remains unchanged. The same assertion applies to countable integer-valued targets.
\end{corollary}
The conclusion concerns the output law and conditional first moments. Purification can change the original joint reference--output distribution, its mutual information, and other conditional properties not among the imposed moments. The map is obtained by measurable allocation; this construction has no asserted polynomial-time implementation.

\subsubsection{A uniform factor preserving finitely many means}
The supermartingale argument uses the following consequence to retain a constant conditional mean over an interval of earlier coordinates.

\begin{lemma}[Mean-independent uniform factor]\label{FP-lem:mean-independent}\label{SP-lem:uniform}
Let $(S,\nu)$ be atomless and let $f:S\to\R^d$ be integrable. There is a measurable $V:S\to[0,1]$ with uniform law and
\[
 \E[f\mid V]=\int f\,d\nu.
\]
For any prescribed probability law $\alpha$ on a standard Borel space $E$, there is therefore a measurable $T:S\to E$ with $T_\#\nu=\alpha$ and the same constant conditional mean of $f$.
\end{lemma}
\begin{proof}
Apply Theorem~\ref{FP-thm:purification} to the two constant weights $1/2$, preserving the vector $(1,f)$. This divides $S$ into two sets with equal masses and equal halves of every $f$-integral. Repeat inside every resulting set. The first $k$ binary choices define $2^k$ sets, each of measure $2^{-k}$ and with $f$-integral $2^{-k}\int f\,d\nu$. Let $V$ be the binary expansion determined by this nested partition. It is uniform, and conditional expectations of $f$ given each finite prefix equal $\int f\,d\nu$. Martingale convergence in $L^1$ gives the identity conditional on $V$. Every standard Borel probability law is a measurable image of uniform measure on $[0,1]$: use a Borel embedding into that interval and its quantile map. Composing with this map proves the last assertion by the tower property.
\end{proof}
This lemma preserves a finite-dimensional mean. It does not assert independence of $V$ from the entire input $S$ or from $f$.


\subsection{Continuous output laws with exact conditional moments}
\label{R10-sec:continuous}
The spaces below are Polish, and measurable maps are taken in Borel versions on sets of full measure. On their product, $W_\infty$ uses the maximum of the coordinate metrics. It is an extended distance and requires no moment assumption.

\begin{theorem}[Conditional moments and backward deterministic approximation]
\label{SP-thm:main}\label{R10-intro:continuous}
Let $X,Y$ take values in Polish spaces, let $\nu=\law(Y)$ be atomless, and let $F(Y)\in L^1(\R^d)$. Suppose
\[
 \E[F(Y)\mid X]=g(X),
 \qquad
 \aff\supp\law(F(Y)\mid X=x)=\R^d
 \quad\text{for }\law(X)\text{-almost every }x.
\]
For every $\eta>0$, there is a Borel map $T$ such that
\begin{align*}
 \law(T(Y))&=\law(X),&
 \E[F(Y)\mid T(Y)]&=g(T(Y)),\\
 W_\infty\bigl(\law(X,Y),\law(T(Y),Y)\bigr)&\le\eta.
\end{align*}
The same assertion holds when the conditional affine spans range over a countable family of affine subspaces, each used with its relative dimension.
\end{theorem}

\subsubsection{A local construction for a continuous output law}
\label{SP-sec:local}
The following lemma converts finite-dimensional freedom in the source means into an exact realization of a whole output distribution. Its quantitative radius is convenient for countable localization; optimal constants are unnecessary for the applications.

For the scalar construction, let $F$ be uniform on $[-1,1]$, and
let $\alpha$ be any centered target law supported on
$[-1/12,1/12]$. The source pieces
$B_+=[-5/12,7/12]$ and $B_-=[-1,1]\setminus B_+$
have probabilities $1/2$ and means $1/12$ and $-1/12$.
Give the two target shares densities $1+12x$ and $1-12x$
relative to $\alpha$. Each is a probability law, their equal mixture
is $\alpha$, and their posterior weights give conditional mean
\[
 \frac{1+12x}{2}\frac1{12}
 +\frac{1-12x}{2}\left(-\frac1{12}\right)=x.
\]
A mean-independent factor on each source piece realizes its entire
target share. This produces the prescribed continuous law and its
conditional mean simultaneously. In higher dimension the two source
means become the vertices of a simplex.

\begin{lemma}[Local realization of conditional means]
\label{SP-lem:simplex}
Let $(S,\beta)$ be atomless, let $F:S\to\R^d$ be integrable with mean $m$, and let $\phi:S\to\R^d$ satisfy $\|\phi\|\le L$, where $L>0$. Assume the matrix
\[
 C=\int (F-m)(\phi-\E\phi)^{\mathsf T}\,d\beta
\]
has smallest singular value at least $a>0$. Let $\alpha$ be a probability law on a standard Borel space $E$, and let $g:E\to\R^d$ satisfy
\[
 \int g\,d\alpha=m,
 \qquad
 \|g(x)-m\|\le\frac{a}{4Ld}
 \quad\text{for }\alpha\text{-almost every }x.
\]
There is a measurable $T:S\to E$ such that
\[
 T_\#\beta=\alpha,
 \qquad
 \E[F\mid T]=g(T).
\]
\end{lemma}
\begin{proof}
Choose the vertices $u_0,\ldots,u_d$ of a regular simplex centered at zero, normalized by
\[
 |u_j|=1,\qquad \sum_{j=0}^d u_j=0,
 \qquad u_i\cdot u_j=-1/d\quad(i\ne j).
\]
The simplex contains the ball of radius $1/d$. Set $\delta=a/(4L)$, $v_j=m+\delta u_j$, and
\begin{equation}\label{SP-eq:simplex-alloc}
 q_j(s)=\frac1{d+1}
 \left[1+\delta(C^{-1}u_j)\cdot(\phi(s)-\E\phi)\right].
\end{equation}
Because $\|C^{-1}\|\le1/a$ and $\|\phi-\E\phi\|\le2L$, the term added to $1$ has absolute value at most $1/2$. Hence the $q_j$ are nonnegative and sum to one. Their moments are
\begin{equation}\label{SP-eq:simplex-moments}
 \int q_j\,d\beta=\frac1{d+1},
 \qquad
 \int Fq_j\,d\beta=\frac{m+\delta u_j}{d+1}=\frac{v_j}{d+1}.
\end{equation}
The second identity follows from $C(C^{-1}u_j)=u_j$; no symmetry of $C$ is assumed. Lemma~\ref{SP-lem:purify} gives a partition $B_0,\ldots,B_d$ with these same masses and $F$-integrals.

The range of $g$ lies in the simplex with vertices $v_j$. Let $w_j(x)$ be its barycentric coordinates there. They are measurable and satisfy
\[
 w_j\ge0,\qquad \sum_jw_j=1,\qquad
 \sum_j w_j(x)v_j=g(x).
\]
Integrating and using the affine independence of the $v_j$ shows that $\int w_j\,d\alpha=1/(d+1)$. Define probability measures
\[
 \alpha_j(dx)=(d+1)w_j(x)\alpha(dx).
\]
The normalized restriction of $\beta$ to $B_j$ is atomless and has $F$-mean $v_j$. Lemma~\ref{SP-lem:uniform} realizes $\alpha_j$ on this restriction by a map $T_j$ with $\E[F\mid T_j]=v_j$. Merge the maps $T_j$. The resulting output law is $\alpha$ because $(d+1)^{-1}\sum_j\alpha_j=\alpha$. For every bounded measurable scalar $h$ on $E$,
\[
 \E[Fh(T)]
 =\frac1{d+1}\sum_jv_j\int h\,d\alpha_j
 =\int h(x)g(x)\alpha(dx).
\]
This is the required conditional-mean identity.
\end{proof}

The construction uses only $d+1$ purified source pieces. A separate mean-independent factor within each piece realizes its share of the full output law. This is the step that allows a continuous target without imposing infinitely many moment constraints in a single application of purification.

\begin{lemma}[A bounded score detects full affine span]
\label{SP-lem:truncated}
Let $Z\in L^1(\R^d)$ have mean $m$ and full affine support. For $R\ge1$, set $\phi_R(Z)=Z\ind_{\{|Z|\le R\}}$ and
\[
 C_R=\E[(Z-m)\phi_R(Z)^{\mathsf T}].
\]
For some integer $R$ and some $a>0$, $\sym C_R\succeq2aI_d$.
\end{lemma}
\begin{proof}
Write
\[
 C_R=A_R+e_Rm^{\mathsf T},\quad
 A_R=\E[(Z-m)(Z-m)^{\mathsf T}\ind_{\{|Z|\le R\}}],\quad
 e_R=\E[(Z-m)\ind_{\{|Z|\le R\}}].
\]
The positive semidefinite matrices $A_R$ increase with $R$. Some $A_{R_0}$ is positive definite: otherwise their decreasing kernels have nontrivial intersection in finite dimension, giving a nonzero vector $u$ with $u\cdot(Z-m)=0$ almost surely. That would contradict full affine support. Moreover $e_R\to0$ by integrability and centering. Consequently $\sym C_R$ is positive definite for all sufficiently large $R$, with a fixed positive lower bound. Reducing that bound gives $2a$.
\end{proof}
Only a first moment is used. Truncation makes the auxiliary score bounded, and the possibly infinite second moment of $Z$ is never needed.

\subsubsection{Exact conditional moments with small displacement}
\label{SP-sec:global}
We prove Theorem~\ref{SP-thm:main}. The localization has two purposes: to make each prescribed conditional-mean profile small enough for Lemma~\ref{SP-lem:simplex}, and to retain the mass of every small source--target rectangle.

\begin{proof}[Proof of Theorem~\ref{SP-thm:main} for full affine span]
Write $\mu=\law(X)$ and $\pi=\law(X,Y)$, and fix regular conditional distributions. Null exceptional sets will be discarded. Put
\[
 \phi_R(y)=F(y)\ind_{\{|F(y)|\le R\}},\qquad
 C_R(x)=\E[(F(Y)-g(x))\phi_R(Y)^{\mathsf T}\mid X=x].
\]
By Lemma~\ref{SP-lem:truncated}, every admissible $x$ belongs to one of the countably many measurable sets specified by
\[
 R\in\N,\quad a\in\{1/n:n\in\N\},\quad
 \sym C_R(x)\succeq2aI_d.
\]
Choose the first such pair in a fixed enumeration and partition accordingly. Refine each part into countably many Borel sets $E_i$ such that, for its assigned $R_i,a_i$,
\begin{equation}\label{SP-eq:target-cells}
 \diam(E_i)\le\eta,
 \qquad
 \diam(g(E_i))\le h_i:=\frac{a_i}{4R_i d}.
\end{equation}
Countable small-diameter covers of the Polish space, combined with small cubes in the range of $g$, give such a partition. Retain cells with $p_i=\mu(E_i)>0$ and set
\[
 \alpha_i=p_i^{-1}\mu|_{E_i},\quad
 \beta_i=\law(Y\mid X\in E_i),\quad
 m_i=\int g\,d\alpha_i.
\]
Then $\int F\,d\beta_i=m_i$ and $|g(x)-m_i|\le h_i$ for $\alpha_i$-almost every $x$.

Let $\bar\phi_i(x)=\E[\phi_{R_i}(Y)\mid X=x]$. The cross-moment matrix of $F$ and $\phi_{R_i}$ under $\beta_i$ is
\begin{align}
 C_i
 &=\int (F-m_i)(\phi_{R_i}-\E_{\beta_i}\phi_{R_i})^{\mathsf T}\,d\beta_i\notag\\
 &=\int C_{R_i}(x)\alpha_i(dx)
   +\int(g(x)-m_i)\bar\phi_i(x)^{\mathsf T}\alpha_i(dx).
 \label{SP-eq:aggregate-matrix}
\end{align}
The final term has norm at most $R_i h_i\le a_i/4$. Therefore $\sym C_i\succeq a_iI_d$, which implies that the smallest singular value of $C_i$ is at least $a_i$.

Now partition the source space into countably many Borel sets $C_\ell$ of diameter at most $\eta$. Let $q_i(y)=\Pp(X\in E_i\mid Y=y)$. On each nonzero restriction $\nu|_{C_\ell}$, apply Lemma~\ref{SP-lem:purify} to these countable allocations. For label $i$, preserve mass and the finitely many scalar entries of
\begin{equation}\label{SP-eq:preserved-vector}
 F,\qquad \phi_{R_i},\qquad F\phi_{R_i}^{\mathsf T}.
\end{equation}
Every function is integrable, since $F\in L^1$ and $\phi_{R_i}$ is bounded. Combine the partitions over $\ell$ to obtain disjoint sets $B_i$. They satisfy
\begin{align}
 \nu(B_i\cap C_\ell)&=\pi(E_i\times C_\ell),\label{SP-eq:rectangle-masses}\\
 \nu(B_i)&=p_i,\quad
 \E_{\widetilde\beta_i}F=m_i,\quad
 C(F,\phi_{R_i};\widetilde\beta_i)=C_i,
 \qquad \widetilde\beta_i=p_i^{-1}\nu|_{B_i}.
 \label{SP-eq:retained-cross}
\end{align}
Here $C(\cdot,\cdot;\cdot)$ denotes the cross-moment matrix in Lemma~\ref{SP-lem:simplex}. The last equality follows by expanding that matrix using the moments in \eqref{SP-eq:preserved-vector}. Each $\widetilde\beta_i$ is atomless.

Apply Lemma~\ref{SP-lem:simplex} on $B_i$ with score $\phi_{R_i}$, bound $L=R_i$, target $\alpha_i$, and conditional-mean profile $g|_{E_i}$. Equations \eqref{SP-eq:target-cells} and \eqref{SP-eq:retained-cross} verify every hypothesis. We obtain $T_i:B_i\to E_i$ with
\[
 (T_i)_\#\widetilde\beta_i=\alpha_i,
 \qquad
 \E_{\widetilde\beta_i}[F\mid T_i]=g(T_i).
\]
The merged map $T$ has the prescribed full target law and the prescribed conditional means.

The old plan $\pi$ and new plan $\pi_T=\law(T(Y),Y)$ give equal mass to every rectangle $E_i\times C_\ell$ by \eqref{SP-eq:rectangle-masses}. Couple their restrictions within each rectangle, and sum these couplings. Both coordinate displacements are at most $\eta$. Hence $W_\infty(\pi,\pi_T)\le\eta$ for the maximum product metric.
\end{proof}

\begin{proof}[Extension to a countable family of affine spans]
Suppose the conditional affine spans belong to $H_1,H_2,\ldots$, where $H_j=b_j+L_j$ and $L_j$ is a linear subspace. Partition the target measurably according to its span. Membership in a fixed affine space is tested by conditional probabilities; full relative span is tested by the truncated matrices in Lemma~\ref{SP-lem:truncated}. On the part assigned to $H_j$ with positive dimension, express $F-b_j$ and $g-b_j$ in an orthonormal basis of $L_j$ and run the preceding argument in dimension $\dim L_j$. Include the allowed set $\{y:F(y)\in H_j\}$ in every purification label originating from this part. Thus the orthogonal coordinates remain fixed, while the relative-coordinate conditional means are preserved.

If $H_j$ is a singleton, $F$ equals that point on the allowed set and $g$ equals it on the corresponding target part. Purify the required masses inside each source cell and realize each target restriction by any atomless map. There is then no nonconstant moment to control. All labels remain countable and each has finitely many integrable moment constraints. The rectangle-mass argument is unchanged.
\end{proof}

\begin{corollary}[One nondegenerate reference coupling suffices]
\label{SP-cor:seed}
Fix $\mu,\nu,F,g$ as in Theorem~\ref{SP-thm:main}, without a support assumption on a given coupling $\pi$. Suppose the constraint set
\[
 \mathcal C=\{\gamma\in\Pi(\mu,\nu):\E_\gamma[F(Y)\mid X]=g(X)\}
\]
contains one coupling $\gamma_*$ whose conditional $F$-laws have full affine span. Then backward deterministic couplings are weakly dense in $\mathcal C$. If the fixed marginals have finite $p$-th moments for some $p\ge1$, they are dense in $W_p$.
\end{corollary}
\begin{proof}
For $\epsilon>0$, $\pi_\epsilon=(1-\epsilon)\pi+\epsilon\gamma_*$ belongs to $\mathcal C$, and its conditional supports contain those of $\gamma_*$. Theorem~\ref{SP-thm:main} gives a backward deterministic plan within $\epsilon$ of $\pi_\epsilon$ in $W_\infty$. Mixtures converge weakly to $\pi$. For finite $p$-th moments, couple the common $(1-\epsilon)$ part identically and couple the remaining parts arbitrarily. The resulting $p$-cost is $O(\epsilon)$ because the marginals are fixed. The triangle inequality gives $W_p$ convergence.
\end{proof}



\subsection{Finite-moment weak transport with arbitrary output dependence}
\label{WC-sec:weak}
Let $\alpha$ be a probability law on a Polish output space $E$, let $\nu$ be an atomless probability on a Polish source space $S$, and let $F\in L^p(\nu;\R^d)$, $p\ge1$. For a coupling $\pi\in\Pi(\alpha,\nu)$, write
\[
 g_\pi(x)=\E_\pi[F(Y)\mid X=x],\qquad m=\int F\,d\nu.
\]
Weak transport allows costs depending on conditional distributions; here the dependence is through a finite vector of their moments. The barycentric case is $F(y)=y$. See \cite{GJ2020,GRST} for the general weak-transport framework.

\begin{theorem}[Exact output law and prescribed moment attenuation]
\label{WC-thm:attenuate}
For every $\pi\in\Pi(\alpha,\nu)$ and every $0\le\lambda<1$, there is a measurable $T:S\to E$ with
\begin{equation}\label{WC-eq:attenuation}
 T_\#\nu=\alpha,
 \qquad
 \E[F(Y)\mid T(Y)=x]=m+\lambda(g_\pi(x)-m).
\end{equation}
The corresponding pair law may be chosen arbitrarily close in $W_\infty$ to $\lambda\pi+(1-\lambda)\alpha\otimes\nu$.
\end{theorem}
\begin{proof}
Work in the affine hull of the source feature law. If that hull is a singleton, the zero-dimensional case of Theorem~\ref{SP-thm:main}, applied to the mixed coupling, gives the output law and the requested displacement; the moment condition is automatic. Otherwise, the mixed coupling has full conditional affine support in the relative coordinates and has the profile displayed in \eqref{WC-eq:attenuation}. Apply Theorem~\ref{SP-thm:main}, retaining the original output space and its entire marginal.
\end{proof}

\begin{theorem}[Equality of finite-moment weak transport infima]
\label{WC-thm:weak-cost}
Suppose $c:E\times\R^d\to\R$ is measurable, continuous in its second variable for $\alpha$-almost every $x$, and
\[
 |c(x,z)|\le a(x)+C|z|^p,
 \qquad a\in L^1(\alpha).
\]
Then
\begin{equation}\label{WC-eq:weak-inf}
 \inf_{T_\#\nu=\alpha}\int c\bigl(x,\E[F(Y)\mid T(Y)=x]\bigr)\alpha(dx)
 =\inf_{\pi\in\Pi(\alpha,\nu)}\int c(x,g_\pi(x))\alpha(dx).
\end{equation}
The same assertion holds with additional measurable pointwise constraints $g_\pi(x)\in C(x)$ whenever every $C(x)$ is star-shaped about $m$ and the admissible set is nonempty. The infima are then taken over those constrained classes.
\end{theorem}
\begin{proof}
Fix an admissible $\pi$ and use Theorem~\ref{WC-thm:attenuate}. Its deterministic profiles $g_\lambda=m+\lambda(g_\pi-m)$ converge pointwise to $g_\pi$ and obey $|g_\lambda|^p\le 2^{p-1}(|g_\pi|^p+|m|^p)$. The right side is integrable by conditional Jensen. Dominated convergence proves convergence of the costs. No continuity in $x$ is used, because its law stays exact and the profiles are given explicitly at that same label. Star-shapedness retains the constraints at every $\lambda$. Taking infima proves one direction; inclusion of deterministic couplings gives the other.
\end{proof}
If $c$ is uniformly $L$-Lipschitz in its moment argument, this proof gives the explicit bound
\[
 |\operatorname{Cost}(T)-\operatorname{Cost}(\pi)|
 \le L(1-\lambda)\int|g_\pi-m|\,d\alpha.
\]
The same sequence retains any finite collection of such costs asymptotically. The theorem permits arbitrary measurable dependence on the output label and requires no convexity of the cost in the moments. Extended-valued hard constraints outside star-shaped feasible sets need a separate argument. In particular, the equality constraint $g(x)=x$ is not stable under this attenuation for a general nonconstant target; \eqref{WC-eq:weak-inf} does not prove the unrestricted exact backward-Monge martingale conjecture.

\begin{example}[A continuous cost with no deterministic minimizer]\label{R10-ex:weak-nonattainment}
Let $Y$ be uniform on $[0,1]$, let the output law be uniform on $[0,1]^2$, take $F(y)=y$, and put $c((a,b),z)=(z-a)^2$. The randomized output $(Y,V)$, with $V$ independent uniform, has cost zero. The deterministic infimum is also zero, but it is unattained.
\end{example}
\begin{proof}
Theorem~\ref{WC-thm:attenuate} gives an exact uniform-square output whose conditional mean is $1/2+\lambda(a-1/2)$, so its cost is $(1-\lambda)^2/12$. If a deterministic output $(A(Y),B(Y))$ attained zero, then $\E[Y\mid A,B]=A$. Since $Y$ and $A$ have the same second moment, $\E(Y-A)^2=0$ and $A=Y$ almost surely. The square output law would then require the function $B(Y)$ to be a nonconstant uniform variable independent of $Y$, a contradiction.
\end{proof}
This example explains why Theorem~\ref{WC-thm:weak-cost} gives equality of infima. The countable-target optimizer in Theorem~\ref{R9-thm:Borel-cost} uses compactness from the fixed masses of the individual labels.


\section{Order constraints and backward-Monge transport}
\label{R11-sec:monge}
A backward-Monge coupling has $X=T(Y)$, so every realization uses the terminal variable alone. Nutz--Wang--Zhang prove existence and weak density for scalar martingale couplings with atomless terminal marginal~\cite{NWZ2024}. They also ask for density under directional and supermartingale constraints, and for higher-dimensional extensions. Those constraints determine the local allocation that is needed here.

For directional transport, only the diagonal rectangles of a fine partition need an ordered allocation. For a supermartingale, the equality region $\E[Y\mid X]=X$ must retain exact martingale means, while the strict-slack region permits a controlled local change. This gives $W_\infty$ directional approximation and $W_1$ supermartingale approximation, with the stated $W_p$ extension. In several dimensions, the countable-initial-marginal theorem retains finite bounded Borel costs as well as both marginals. Strict contraction creates interior convex-order slack for arbitrary initial marginals; the common-fiber endpoint theorem uses its separate measurability hypothesis. Keeping these statements distinct identifies exactly which part of the higher-dimensional question each argument answers.

\subsection{Multidimensional martingales and strict contraction}\label{SP-sec:martingales}
The moment realization theorem now applies with $F(y)=y$ and $g(x)=x$. The following statements preserve the complete marginal laws; their respective support and approximation hypotheses are part of the conclusions about backward Monge transport.

\begin{theorem}[Backward Monge approximation in higher dimensions]
\label{SP-thm:multidimensional}
Let $\nu$ be atomless and let $\pi\in\M(\mu,\nu)$. If
\[
 \aff\supp\pi_x=\R^d\quad\text{for }\mu\text{-almost every }x,
\]
then for every $\eta>0$ there exists $\pi_T\in\MM(\mu,\nu)$ with $W_\infty(\pi,\pi_T)\le\eta$. More generally the same assertion holds if the conditional affine spans of $\pi$ belong to a countable family.

If some $\pi_*\in\M(\mu,\nu)$ has full-dimensional conditional supports, then $\MM(\mu,\nu)$ is weakly dense in $\M(\mu,\nu)$, and is $W_p$ dense whenever both marginals have finite $p$-th moments.
\end{theorem}
\begin{proof}
Take $F(y)=y$ and $g(x)=x$ in Theorem~\ref{SP-thm:main}, and use Corollary~\ref{SP-cor:seed} for the second assertion.
\end{proof}

\paragraph{The exact conjecture reached.}
Nutz--Wang--Zhang \cite[Conjecture~5.1]{NWZ2024} formulate higher-dimensional backward Monge density using atomlessness of the terminal laws after irreducible decomposition. Their discussion uses the maximal martingale coupling $\widehat\pi$ of the convex-paving theory of De~March--Touzi \cite{DMT2019}, with
\[
 I(x)=\operatorname{ri}\conv\supp\widehat\pi_x
\]
and every other conditional martingale support contained in the corresponding closed affine span. If $\dim I(x)=d$ almost everywhere, $\widehat\pi$ is the full-dimensional coupling needed in Theorem~\ref{SP-thm:multidimensional}. Thus the asserted density conclusion of that conjecture holds for all such pairs with atomless $\nu$.

The same argument works when the affine hulls of the maximal components belong to a countable family. Mix an arbitrary coupling with $\widehat\pi$; the affine hull of each mixed conditional support is then the hull of the maximal component. Apply the countable-affine version of the theorem. The countable-initial-marginal clause is proved directly by the countable allocation in Theorem~\ref{FP-thm:discrete-density}.

Existence on the real line also has the injective construction of Hobson--Norgilas \cite{HN2023}.

Theorem~\ref{WC-thm:fiber-martingale} below also allows uncountably varying lower-dimensional components when their common parameter is recoverable from both coordinates and the terminal laws are conditionally atomless. General convex pavings with shared terminal boundary mass retain an additional issue. In particular, the example in \cite[Remark~5.2]{NWZ2024}, with common first coordinate and a discrete conditional second-coordinate source, remains an obstruction to an unrestricted global-atomlessness statement.

\begin{corollary}[Equality of optimal transport values]
\label{SP-cor:cost}
Under the full-dimensional coupling hypothesis of Theorem~\ref{SP-thm:multidimensional}, every continuous cost $c$ with
\[
 |c(x,y)|\le a(x)+b(y),\qquad a\in L^1(\mu),\ b\in L^1(\nu),
\]
satisfies
\[
 \inf_{\pi\in\MM(\mu,\nu)}\int c\,d\pi
 =\inf_{\pi\in\M(\mu,\nu)}\int c\,d\pi.
\]
\end{corollary}
\begin{proof}
Use weak approximation and uniform integrability furnished by the fixed marginals and the displayed bound. To pass from bounded to unbounded $c$, truncate it continuously; the tails are uniformly integrable because $a(X)+b(Y)$ has fixed integrable marginal summands. The reverse inequality is inclusion of the admissible sets.
\end{proof}
This is equality of infima. The result does not assert that an optimizer for every cost is backward Monge.

\begin{theorem}[Strict homothetic contractions]
\label{SP-thm:contraction}\label{R10-eq:contraction}
Let $\mu\cx\nu$ be integrable laws on $\R^d$, with $\nu$ atomless and common mean $m$. For every $0\le\lambda<1$, put
\[
 S_\lambda(x)=m+\lambda(x-m),\qquad \mu_\lambda=(S_\lambda)_\#\mu.
\]
There exists a Borel $T_\lambda:\R^d\to\R^d$ with
\[
 (T_\lambda)_\#\nu=\mu_\lambda,
 \qquad \E[Y\mid T_\lambda(Y)]=T_\lambda(Y).
\]
If $\mu,\nu$ have finite $p$-th moments and $\pi\in\M(\mu,\nu)$ is prescribed, the maps may be selected as $\lambda\uparrow1$ so that
\[
 W_p\bigl((T_\lambda,\id)_\#\nu,\pi\bigr)\longrightarrow0.
\]
\end{theorem}
\begin{proof}
Given $\pi\in\M(\mu,\nu)$, Theorem~\ref{WC-thm:attenuate} with $F(y)=y$ constructs a map $A_\lambda$ with unchanged output law $\mu$ and conditional mean $S_\lambda(A_\lambda(Y))$. For $\lambda>0$, set $T_\lambda=S_\lambda\circ A_\lambda$. Injectivity of $S_\lambda$ gives the martingale identity; at $\lambda=0$ the constant map $m$ suffices.

For the approximation choose the pair law of $A_\lambda$ within $1-\lambda$ in $W_\infty$ of $\lambda\pi+(1-\lambda)\mu\otimes\nu$. This mixture tends to $\pi$ in $W_p$. The final contraction changes the $p$th power of the displacement by at most $(1-\lambda)^p\int |x-m|^p\mu(dx)$. The triangle inequality proves the claim.
\end{proof}
The fixed-source assertion holds for every strict contraction and in every dimension. The approximation changes the first marginal from $\mu$ to $\mu_\lambda$; exact preservation of $\mu$ at $\lambda=1$ is the additional issue in the full higher-dimensional conjecture.

\begin{corollary}[Arbitrarily small inflation of an atomless reference]
\label{SP-cor:inflation}
If $X\cx R$, $R$ is atomless and integrable, and $m=\E R$, then for every $c>1$ the unchanged target law $\law(X)$ has a backward Monge martingale realization from the reference law of $m+c(R-m)$.
\end{corollary}
\begin{proof}
Apply Theorem~\ref{SP-thm:contraction} with $\lambda=1/c$ to the source law of $R$. Affinely rescale both coordinates by $z\mapsto m+c(z-m)$. The first marginal becomes $\law(X)$ and the source is the asserted inflated reference.
\end{proof}



\subsection{Measurable families and uncountably many affine fibers}
\label{WC-sec:fibers}
The preceding continuous-output proof can be selected measurably over a parameter. It therefore also applies when the affine supports form an uncountable measurable family.

\begin{theorem}[Conditional realization with a parameter]
\label{WC-thm:parameter}
Let the parameter space be standard Borel with probability law $\xi$, and let $z\mapsto\pi_z$ be a probability kernel on a product of Polish spaces $E\times S$, with marginals $\alpha_z,\nu_z$, and measurable features $F(z,y)\in\R^d$. Suppose $F$ is integrable for the joint source law, $\nu_z$ is atomless for $\xi$-almost every $z$, and
\[
 \E_{\pi_z}[F(z,Y)\mid X]=g(z,X),
 \qquad
 \aff\supp\law_{\pi_z}(F(z,Y)\mid X=x)=\R^d
\]
for almost every $(z,x)$. For every measurable positive $\eta(z)$, there is a jointly measurable map $T(z,y)$ such that, almost surely in $z$,
\begin{align*}
 T(z,\cdot)_\#\nu_z&=\alpha_z,\\
 \E[F(z,Y)\mid z,T(z,Y)]&=g(z,T(z,Y)),\\
 W_\infty\bigl(\pi_z,(T(z,\cdot),\id)_\#\nu_z\bigr)&\le\eta(z).
\end{align*}
Maps and kernels can be taken Borel after modification on a set of zero joint source measure.
\end{theorem}
\begin{proof}
We give the selection argument, since fiberwise existence alone does not produce a measurable map. The space $\mathcal P(E\times S)$ is standard Borel. There is a jointly Borel version of regular conditional distributions as a function of a probability measure and its conditioning coordinate; it may be constructed by conditional expectations on a countable refining sequence of finite partitions and almost-sure limits. Fix such a disintegration $\kappa^y$ for each $\kappa\in\mathcal P(E\times S)$.

The set of backward deterministic couplings is Borel. Indeed choose a countable family of bounded measurable real functions separating points of $E$ (for example indicators of a countable separating algebra). The conditional source-to-output kernel is Dirac almost everywhere exactly when every conditional variance
\[
 \int\left[\int h(x)^2\kappa^y(dx)-\left(\int h(x)\kappa^y(dx)\right)^2\right]\kappa_S(dy)
\]
vanishes. Each displayed quantity is a Borel function of $\kappa$, and vanishing of all of them forces a probability kernel to be supported on a single point.

The prescribed marginal conditions are Borel. The conditional-moment condition is equivalent to the countable family of identities
\[
 \int h(x)(F(z,y)-g(z,x))\,\kappa(dx,dy)=0
\]
for indicators $h$ from a generating algebra, together with the corresponding integrability conditions. These are Borel conditions by truncation of the measurable integrands. Finally, $W_\infty(\kappa,\pi_z)\le\eta(z)$ is Borel: for each fixed radius the relation admits a coupling supported in the closed radius set, and is closed under weak limits by tightness. Thus the set of feasible $(z,\kappa)$ is Borel.

For almost every $z$ its section is nonempty by Theorem~\ref{SP-thm:main}. The Jankov--von Neumann selection theorem provides a universally measurable selected coupling $z\mapsto\kappa_z$ \cite{Kechris1995}. Relative to the fixed law $\xi$, replace this selector by a Borel version. Disintegrate $\kappa_z$ over $y$. The conditional measures are Dirac almost everywhere, and the inverse of the Borel embedding $x\mapsto\delta_x$ extracts their unique point as a jointly measurable $T(z,y)$. The marginal, moment and distance identities follow from the selection. Exceptional parameter and source sets may be assigned arbitrary values.
\end{proof}
The fixed finite feature dimension can vary measurably: partition the parameter space by that dimension. For a measurable field of affine subspaces, choose their nearest points to the origin and orthonormal bases of their direction spaces by measurable Gram--Schmidt applied to a fixed countable dense set. The theorem then applies in relative coordinates.

\begin{theorem}[Backward Monge density with a common parameter]
\label{WC-thm:fiber-martingale}
Let $\pi$ be a martingale coupling of integrable laws $\mu,\nu$ on $\R^d$. Suppose there are measurable maps $a,b$ into a standard Borel space such that $Z=a(X)=b(Y)$ under $\pi$. Disintegrate $\pi$ over $Z$ as $\pi_z\in\mathcal M(\mu_z,\nu_z)$. Assume that, for almost every $z$, the marginals lie in a measurable affine space $H_z$, the terminal law $\nu_z$ is atomless, and the conditional terminal supports of $\pi_z$ fill $H_z$ affinely for $\mu_z$-almost every starting point. Then $\pi$ has arbitrarily $W_\infty$-close backward Monge martingale couplings preserving both marginals and the common parameter.

More generally, fix the conditional marginals $(\mu_z,\nu_z)$ and suppose a measurable family of full-relative-span martingale seed couplings exists. Backward Monge martingale couplings preserving $Z$ are weakly dense among all martingale couplings with those conditional marginals and that common parameter, and are $W_p$ dense under finite $p$-th moments. Zero-dimensional fibers may be included without atomlessness; their coupling is the identity.
\end{theorem}
\begin{proof}
In positive-dimensional fibers, use relative orthonormal coordinates and apply Theorem~\ref{WC-thm:parameter} with $F(y)=y$ and $g(x)=x$ in those coordinates. The parameter is recovered from the terminal point as $z=b(y)$, so $T(y)=T(b(y),y)$ is a function of the unchanged terminal state alone. The law of its output conditional on $z$ is $\mu_z$, which is supported on $a^{-1}(z)$. Therefore the selected maps preserve the common parameter as well as the full marginals and martingale condition. Integrating the fiberwise displacement couplings gives the global $W_\infty$ bound. The required displacement couplings can themselves be selected by the same Borel coupling argument.

For the second assertion, mix the given conditional plan with an arbitrarily small positive fraction of the seed on each fiber. The mixed conditional support fills $H_z$. Apply the first assertion, then let the mixing fraction and approximation radius tend to zero. For $W_p$, couple the common part of the mixture identically; the remaining $p$-cost tends to zero because both global marginals have finite $p$-th moments. A zero-dimensional fiber has equal point masses as its two marginals and is treated by the identity map.
\end{proof}

\begin{corollary}[Terminal recovery away from shared component closures]
\label{R19-cor:separated-labels}
Let $\pi\in\mathcal M(\mu,\nu)$ have integrable marginals on $\R^d$, and let $Z=a(X)$ take values in a standard Borel space $S$. Suppose $C_z\subset\R^d$ are closed sets whose incidence relation $\{(z,y):y\in C_z\}$ is Borel, and $Y\in C_Z$ almost surely. Assume there is a Borel $\nu$-null set $N$ such that, for each $y\notin N$, at most one $z$ satisfies $y\in C_z$. Then a Borel map $b:\R^d\to S$ satisfies $b(Y)=a(X)$ almost surely.

Suppose in addition that the conditional marginals lie in measurable affine spaces $H_z$, that $\nu_z=\law(Y\mid Z=z)$ is atomless on positive-dimensional fibers, and that the conditional terminal supports of $\pi_z$ fill $H_z$ affinely. Then $\pi$ is a $W_\infty$ limit of backward-Monge martingale couplings with both marginals and $Z$ preserved. If a measurable full-relative-span seed exists instead, the same couplings are weakly dense among all plans with those conditional marginals, and $W_p$ dense under finite $p$-moments.
\end{corollary}
\begin{proof}
The relation
$\mathcal R=\{(y,z):y\notin N,\ y\in C_z\}$
is Borel and has at most one point in each section over $y$. The Lusin--Novikov theorem~\cite{Kechris1995} makes its projection Borel and its unique-point selector Borel there. Extend that selector by any fixed value of $S$ outside the projection. Since $Y\in C_{a(X)}$ and $\nu(N)=0$, the resulting map satisfies $b(Y)=a(X)$ almost surely. Theorem~\ref{WC-thm:fiber-martingale} gives the remaining assertions with this common parameter. Zero-dimensional fibers use their identity coupling.
\end{proof}
For a measurable convex paving with $Y\in\overline{I(X)}$, take $C_z$ to be the component closures. Avoidance of their overlaps recovers the component label; a measurable maximal seed whose conditional support fills the affine hull verifies the separate span assumption. Shared-boundary mass requires an allocation analysis, as the unique backward-bit example in Theorem~\ref{BM-thm:main} demonstrates. The condition above is sufficient and makes no claim that every overlap forces randomness.

\begin{corollary}[Common-coordinate families]
\label{WC-cor:common-coordinate}
Let $Z$ have any integrable law on $\R^r$, and consider marginals of $(Z,U)$ and $(Z,V)$ on $\R^{r+d}$. Suppose their conditional laws $(\mu_z,\nu_z)$ in the last $d$ coordinates admit a measurable full-dimensional martingale seed, and $\nu_z$ is atomless almost surely. Then backward Monge martingale couplings are weakly dense in the full martingale coupling class, and are $W_p$ dense when the marginals have finite $p$-th moments.
\end{corollary}
\begin{proof}
Every martingale coupling must preserve its first $r$ coordinates. Their marginals coincide; apply conditional Jensen to the strictly convex integrable function $z\mapsto\sqrt{1+|z|^2}$. Equality of its two expectations forces those coordinates to agree almost surely. Thus every martingale coupling belongs to the common-parameter class of Theorem~\ref{WC-thm:fiber-martingale}.
\end{proof}
Nutz--Wang--Zhang Conjecture~5.1 concerns backward-Monge density after irreducible decomposition~\cite{NWZ2024}. The full-dimensional and countable-affine-family cases are in Theorem~\ref{SP-thm:multidimensional}; a recoverable common label permits uncountably many lower-dimensional fibers here. Corollary~\ref{R19-cor:separated-labels} gives a geometric way to verify recoverability. General convex pavings can share terminal boundary mass~\cite{DMT2019}, which is precisely where that verification matters. Theorem~\ref{BM-thm:main} exhibits the unique martingale for which the shared label forces backward randomness.


The common-parameter theorem glues its local maps using a component label recoverable from the terminal state. At shared boundaries, the same terminal point can belong to different components, leaving that choice unresolved. The following example forces exactly a Bernoulli choice and loses the obstruction under every strict contraction or inflation. Interior stability then corrects both endpoint laws along every prescribed initial coupling. Its exact product reserve repairs conditional means and the terminal marginal simultaneously. For approximate optimizers, strict concavity then excludes a hidden mixture of terminal kernels at any limiting initial state. This gives convergence of conditional laws when both endpoints and the driver vary, with an absolutely continuous limiting driver and possibly atomic approximants.

\subsection{Atomless components with an unavoidable backward choice}
\label{BM-sec:boundary}
The component theorems preserve the allocation of terminal mass as well as its conditional means. At a shared boundary, that allocation can itself carry randomness. We now construct a unique martingale for which every component has an atomless terminal law and a constant backward map, but choosing between those local maps requires an independent Bernoulli variable. Thus the terminal-recoverability condition in Theorem~\ref{WC-thm:fiber-martingale} concerns a real part of the realization problem.

\begin{theorem}[A unique martingale with atomless components]\label{BM-thm:main}
For each $0<p<1$ there are centered, compactly supported, atomless laws $\mu_p,\nu$ on $\R^4$, both with full affine span, such that $\M(\mu_p,\nu)=\{\pi_p\}$. Every forward conditional law $(\pi_p)_x$ is atomless with affine span of dimension two. There are Borel maps $h,k:\mathbb S^1\to\R^4$ such that, writing $Y=(U,V)$,
\begin{equation}\label{BM-eq:reverse}
 \law(X\mid Y=(u,v))=p\delta_{h(v)}+(1-p)\delta_{k(u)},
 \qquad h(v)\ne k(u).
\end{equation}
The maximal irreducible component map $I$ satisfies
\begin{equation}\label{BM-eq:component}
 \nu_I:=\law_{\pi_p}(Y\mid I(X)=I)
 \quad\text{is atomless for }\law_{\mu_p}(I(X))\text{-almost every }I.
\end{equation}
Every component contains exactly one point of a Borel carrier of $\mu_p$, and almost every terminal point belongs to the relative boundaries of two used components. In particular $\MM(\mu_p,\nu)=\varnothing$. The same $\nu$ works for all $p$, and both its circle marginals are uniform.
\end{theorem}

The geometric mechanism is an equality in a second moment. If $Y=(Y_1,Y_2)$ lies on $\mathbb S^1\times\mathbb S^1$ and the initial point is $(b,v)$ with $|v|=1$, the martingale equation gives
\[
 \E[|Y_2-v|^2\mid X=(b,v)]=1-2|v|^2+|v|^2=0.
\]
That boundary coordinate freezes $Y_2$. A second family freezes $Y_1$. There are therefore only two possible predecessors of each terminal point. An extreme ordinary coupling will fix their probabilities.

\begin{lemma}[An extreme coupling on two circles]\label{BM-lem:extreme}
There exists an extreme point $\zeta\in\Pi(\sigma,\sigma)$, where $\sigma$ is uniform arclength on $\mathbb S^1$, whose conditional laws in both directions are atomless. If
\[
 b(v)=\E[U\mid V=v],\qquad c(u)=\E[V\mid U=u],\qquad (U,V)\sim\zeta,
\]
then $|b(v)|<1$, $|c(u)|<1$ almost everywhere, and every conditional support has affine span $\R^2$.
\end{lemma}
\begin{proof}
Kun's Theorem~1.4 gives a conditionally atomless coupling on standard Borel spaces which is the only coupling carried by a specified Borel relation~\cite{Kun2024}. It is extreme in the set of all couplings: each term of a convex decomposition must also be carried by that relation. Its two marginals are atomless, since a marginal atom would be a conditional atom on a positive-measure set in the reverse disintegration. Transport both marginal spaces by measure-preserving Borel isomorphisms to $(\mathbb S^1,\sigma)$; extremality and conditional atomlessness are retained. A line meets the circle in at most two points, so each atomless conditional law has full affine span. Equality in the triangle inequality for its mean would concentrate it at one point, proving the strict inequalities.
\end{proof}

\begin{proof}[Proof of Theorem~\ref{BM-thm:main}]
Use the coupling in Lemma~\ref{BM-lem:extreme} and set
\[
 h(v)=(b(v),v),\qquad k(u)=(u,c(u)),\qquad
 \nu=\zeta,\quad \mu_p=p\,h_\#\sigma+(1-p)\,k_\#\sigma.
\]
The two initial carriers are disjoint: they lie respectively in $B_2^\circ\times\mathbb S^1$ and $\mathbb S^1\times B_2^\circ$. Draw $(U,V)\sim\zeta$ and, independently, $B\sim\Bern(p)$; choose $X=h(V)$ if $B=1$ and $X=k(U)$ otherwise. The definitions of $b,c$ prove $\E[Y\mid X]=X$. The marginals are centered and atomless, and all variables are bounded.

Consider any other martingale with these marginals. The boundary calculation above forces $V=v$ when $X=h(v)$ and $U=u$ when $X=k(u)$. Its reverse kernel must therefore have the form
\[
 a(u,v)\delta_{h(v)}+(1-a(u,v))\delta_{k(u)},\qquad 0\le a\le1.
\]
Preservation of the initial marginal gives $\E[a\mid V]=p$ and $\E[a\mid U]=p$. Hence
\[
 \zeta_H=\frac ap\zeta,\qquad
 \zeta_K=\frac{1-a}{1-p}\zeta
 \quad\text{belong to }\Pi(\sigma,\sigma),\qquad
 \zeta=p\zeta_H+(1-p)\zeta_K.
\]
Extremality implies $\zeta_H=\zeta_K=\zeta$, so $a=p$ almost surely. This proves uniqueness and~\eqref{BM-eq:reverse}. The Borel relation in Kun's input has not been imposed on the martingale problem; the specified marginals and martingale equation alone forced this allocation.

Full affine span of $\nu$ follows by conditioning a putative affine relation $\alpha\cdot U+\beta\cdot V=d$ on $V$ and using conditional atomlessness on the circle: $\alpha=0$, and then $\beta=0$. If $\mu_p$ lay in a proper affine hyperplane, the conditional mean of its defining affine functional of $Y$ would be constant. More directly, equality of $\alpha\cdot b(v)+\beta\cdot v$ and $\alpha\cdot u+\beta\cdot c(u)$ to the same constant implies that $\alpha\cdot U+\beta\cdot V$ has both conditional expectations equal to that constant. Its variance in $L^2(\zeta)$ is then zero, by expanding against these two conditional expectations. The full-span assertion for $\nu$ again gives $\alpha=\beta=0$.

It remains to identify the components, including their boundary allocation. Put
\[
 C_v=\conv\supp\law(U\mid V=v),\qquad
 D_u=\conv\supp\law(V\mid U=u).
\]
They are two-dimensional compact convex sets, and their conditional barycenters lie in their interiors. Since the martingale is unique, the maximal conditional convex supports in the irreducible paving of De March--Touzi~\cite[Theorem~2.1]{DMT2019} are exactly
\begin{equation}\label{BM-eq:paving}
 I(h(v))=\operatorname{int}C_v\times\{v\},\qquad
 I(k(u))=\{u\}\times\operatorname{int}D_u.
\end{equation}
These sets are pairwise disjoint: different horizontal sets have different second coordinates, different vertical sets have different first coordinates, and the two families lie in $B_2^\circ\times\mathbb S^1$ and $\mathbb S^1\times B_2^\circ$, respectively. Each contains only its indicated point of that carrier. Thus conditioning on $I(X)$ is equivalent to conditioning on $X$, and~\eqref{BM-eq:component} is the atomless conditional law from Lemma~\ref{BM-lem:extreme}. Almost surely $u\in\supp\law(U\mid V=v)$ and $v\in\supp\law(V\mid U=u)$. Points of the circle are extreme in the disk, so $(u,v)$ lies on the relative boundaries of both components in~\eqref{BM-eq:paving}.
\end{proof}

\begin{corollary}[The component-disintegration formulation]\label{BM-cor:NWZ}
Atomlessness of the laws $\law_\pi(Y\mid I(X))$ after maximal irreducible-component disintegration does not imply existence, and hence does not imply density, of backward-Monge martingales. This gives a negative answer to the component-disintegrated formulation of Nutz--Wang--Zhang Conjecture~5.1~\cite[Section~5]{NWZ2024}, with the disintegration specified by~\eqref{BM-eq:component}.
\end{corollary}
\begin{proof}
The coupling of Theorem~\ref{BM-thm:main} is unique, so its component disintegration is canonical. All its component terminal laws are atomless, whereas its reverse kernel has two distinct atoms.
\end{proof}

The component label in this statement is $I(X)$. Requiring a terminal-measurable label $I(X)=J(Y)$ imposes an additional condition on shared boundary mass, as distinguished in~\cite[Section~2.2]{DMT2019}. It fails here: conditional on $Y=(u,v)$, the two distinct labels in~\eqref{BM-eq:paving} have probabilities $p$ and $1-p$. Thus the positive common-parameter theorem retains its stated terminal-recoverability hypothesis. Each local backward map is constant, but the global marginal constraint prevents those maps from being chosen deterministically.

\begin{corollary}[Exact reverse randomness and failure at the endpoint]\label{BM-cor:quantitative}
For the coupling above, its branch $B$ is independent of $Y$ and has entropy
$H(B\mid Y)=h_2(p)=-p\log p-(1-p)\log(1-p)$. Moreover
\begin{equation}\label{BM-eq:TV}
 \sup_T\Prb\{X=T(Y)\}=\max(p,1-p),\qquad
 \inf_T\TV\bigl(\pi_p,(T,\id)_\#\nu\bigr)=\min(p,1-p),
\end{equation}
where the extrema range over all Borel maps, without a marginal or martingale restriction. For every $0\le\lambda<1$, $\MM((\lambda\id)_\#\mu_p,\nu)$ is nonempty. Every choice of these couplings converges to $\pi_p$ in each finite-order Wasserstein distance as $\lambda\uparrow1$, but has total-variation distance one from $\pi_p$. Likewise every strict inflation $(c\id)_\#\nu$, $c>1$, admits a backward-Monge coupling from the unchanged $\mu_p$, whereas $c=1$ fails.
\end{corollary}
\begin{proof}
Formula~\eqref{BM-eq:reverse} gives the entropy and success probability. Disintegrating total variation over the common terminal marginal shows that its distance to a graph kernel is $1-\pi_p(\{T(y)\}\mid y)$; choosing the more probable branch attains~\eqref{BM-eq:TV}. The strict-contraction theorem gives the contracted realizations. Compactness makes every subsequential limit a martingale with marginals $\mu_p,\nu$, hence the unique $\pi_p$. Bounded support upgrades weak convergence to every finite $W_r$. The original initial carrier has $\max(|x_1|,|x_2|)=1$ and its contraction has that maximum $\lambda$, proving total variation one. Scaling a coupling from $c^{-1}\mu_p$ to $\nu$ by $c$ gives the terminal-inflation assertion.
\end{proof}

At $p=1/2$ the exact unresolved choice is one fair bit, even after observing the complete terminal state. The next argument works under an interior condition which produces a positive product component. That component permits a bounded correction of both marginals and every conditional mean; no such product component exists in this boundary example.


\subsection{Interior slack and changes of all marginals}\label{AS-sec:stability}
Br\"uckerhoff--Juillet show that unrestricted higher-dimensional martingale transport is unstable even under full-support perturbations~\cite{R20-BJ}. Real-line stability is established in~\cite{R20-BP,R20-Wiesel}. The hypothesis below gives an exact interior condition in every finite dimension; it allows both marginals to vary and preserves any prescribed initial coupling.

A perturbation of the terminal law changes both the available conditional means and the mass assigned to each initial state. An interior reserve permits these errors to be corrected together. The reserve used below is an exact measure inequality,
\begin{equation}\label{AS-eq:product}
 s\in\M(\mu,\nu),\qquad s\ge\varepsilon\mu\otimes\nu,
 \qquad\varepsilon>0.
\end{equation}
A bounded finite-rank correction will have zero row mass, the required row mean and exactly the required column marginal. Its positivity comes from~\eqref{AS-eq:product}.

For $1\le p<\infty$ and $\pi=\mu(dx)\pi_x(dy)$, define the two-time adapted Wasserstein distance by
\[
 AW_p(\pi,\pi')^p=\inf_{\xi\in\Pi(\mu,\mu')}
 \int\bigl(|x-x'|^p+W_p(\pi_x,\pi'_{x'})^p\bigr)\,d\xi.
\]
Write $m_p(r)=\int|y|^p\,r(dy)$. Here irreducibility means that for every Borel $A,B$ of positive $\mu$- and $\nu$-mass, some martingale coupling assigns positive mass to $A\times B$.

\begin{theorem}[All-marginal adapted approximation]\label{AS-thm:main}
Let $\mu\cx\nu$ in $\cP_p(\R^d)$, let $\nu$ have full affine span, and suppose $\mu$ is compactly supported and~\eqref{AS-eq:product} holds. In particular these assumptions hold if the pair is irreducible and
\begin{equation}\label{AS-eq:interior}
 \supp\mu\Subset\operatorname{int}\overline{\conv}\supp\nu.
\end{equation}
If $\mu_n\cx\nu_n$ and $W_p(\mu_n,\mu)+W_p(\nu_n,\nu)\to0$, then, for every $\pi\in\M(\mu,\nu)$ and every prescribed $\xi_n\in\Pi(\mu,\mu_n)$ with $\int|x-x'|^p\,d\xi_n\to0$, there are $\pi_n\in\M(\mu_n,\nu_n)$ such that
\begin{equation}\label{AS-eq:prescribed}
 \int W_p(\pi_x,(\pi_n)_{x'})^p\,\xi_n(dx,dx')\longrightarrow0.
\end{equation}
Thus $AW_p(\pi_n,\pi)\to0$, with both marginals and every conditional mean exact at each $n$. The approximating pairs may be reducible and have unbounded supports.
\end{theorem}

\begin{lemma}[Outward dilation and an exact product component]\label{AS-lem:product}
If $m$ is the common mean and $0<\varepsilon<1$, condition~\eqref{AS-eq:product} is equivalent to
\begin{equation}\label{AS-eq:dilation}
 \left(x\mapsto m+\frac{x-m}{1-\varepsilon}\right)_\#\mu\cx\nu.
\end{equation}
An irreducible pair satisfying~\eqref{AS-eq:interior} has this property for some $\varepsilon>0$.
\end{lemma}
\begin{proof}
Subtract $\varepsilon\mu\otimes\nu$ from $s$ and divide by $1-\varepsilon$. The remaining kernel has conditional mean $m+(x-m)/(1-\varepsilon)$; relabeling proves one direction. Conversely, pull a martingale for~\eqref{AS-eq:dilation} back to $x$ and mix its kernel with $\nu$ in proportions $1-\varepsilon,\varepsilon$.

We prove that the geometric hypotheses allow a small outward dilation. Write $U=\operatorname{int}\overline{\conv}\supp\nu$. On every compact $L\Subset U$, nonnegative convex functions satisfy
\begin{equation}\label{AS-eq:localconvex}
 \sup_L f\le C_L\int f\,d\nu.
\end{equation}
Indeed finitely many support points surround a neighborhood of $L$. Choose small positive-mass balls around them whose arbitrary representative points still surround $L$. The integral bound guarantees a point of bounded $f$ in each ball, and convexity bounds $f$ inside their convex hull. Applying this on a larger compact set also bounds the local Lipschitz constant.

There is a uniform $\delta>0$ such that
\begin{equation}\label{AS-eq:Jensen}
 \int f\,d\mu\le(1-\delta)\int f\,d\nu
 \quad(f\ge0\text{ convex},\ f(m)=0).
\end{equation}
Otherwise normalize $\int f_n\,d\nu=1$ and take a locally uniform convex limit using~\eqref{AS-eq:localconvex}. Compact initial support gives $\int f\,d\mu=1$. Extend $f$ closed-convexly to the terminal convex hull. The radial inequality
$f_n(y)\ge t^{-1}f_n(m+t(y-m))$, $t<1$, and Fatou give $\int f\,d\nu\le1$. Jensen gives equality. Choose bounded subgradients $a(x)$ on the compact initial support. Every martingale must then avoid the positive set of
$D(x,y)=f(y)-f(x)-a(x)\cdot(y-x)$.
If $f$ is not affine on $U$, this set has positive product mass and contains a positive-mass Borel rectangle: bound $y$, fix a positive rational gap, and partition the coefficients $(a(x),f(x)-a(x)\cdot x)$ into sufficiently small boxes. Irreducibility excludes such a rectangle. Thus $f$ is affine; nonnegativity and $f(m)=0$ force $f=0$, a contradiction.

Choose a compact neighborhood of $\conv\supp\mu$ inside $U$. By~\eqref{AS-eq:localconvex}, dilation by $1+\theta$ changes each normalized integral against $\mu$ by at most $C\theta\int f\,d\nu$. Taking $C\theta<\delta$ in~\eqref{AS-eq:Jensen}, and subtracting an affine support at $m$ from a general convex test, proves~\eqref{AS-eq:dilation} for $\varepsilon=\theta/(1+\theta)$ by Strassen's criterion.
\end{proof}

The interior step cannot be replaced by~\eqref{AS-eq:Jensen} alone. For $\nu=(\delta_{-1}+\delta_1)/2$ and $\mu=(1-\delta)\nu+\delta\delta_0$, that gap holds, but each initial endpoint forces the identical terminal endpoint and precludes~\eqref{AS-eq:product}.

\begin{proof}[Proof of Theorem~\ref{AS-thm:main}]
\emph{A positive finite-rank kernel.}
Choose a bounded measurable vector $h$ with $\int h\,d\nu=0$ and
\[
 C=\int y h(y)^{\mathsf T}\,\nu(dy)\quad\text{invertible}.
\]
One can take $h(y)=(y-m)\mathbf1_{\{|y-m|\le R\}}-\int(z-m)\mathbf1_{\{|z-m|\le R\}}\,\nu(dz)$. Then $C$ is a truncated second-moment matrix and is positive definite for large finite $R$. This uses only a first moment.

Mix $\pi$ with a small amount $\tau s$ of~\eqref{AS-eq:product}. Average the resulting terminal kernels over a finite partition into small continuity cells $C_i$, with weights $w_i>0$ and barycenters $b_i$. Their densities $r_i$ relative to $\nu$ satisfy
\[
 \tau\varepsilon\le r_i\le w_i^{-1},\quad
 \int r_i\,d\nu=1,\quad \int y r_i(y)\,\nu(dy)=b_i,
 \quad\sum_iw_i r_i=1.
\]
For $x\in C_i$ put
\begin{equation}\label{AS-eq:regular}
 k(x,y)=r_i(y)+h(y)^{\mathsf T}C^{-1}(x-b_i).
\end{equation}
The added term has zero mass, mean $x-b_i$ and zero column integral. A sufficiently small cell diameter preserves a uniform positive lower bound. Thus $k$ is a bounded, positive martingale density with both prescribed marginals.

These kernels approximate $\pi$ in $\int W_p^p\,d\mu$. Refining continuity partitions generate the initial sigma-field, so conditional expectations of a countable determining family of tests and of $|Y|^p$ converge. Their averaged terminal law is $\nu$, which ensures uniform integrability. The added density in~\eqref{AS-eq:regular} tends uniformly to zero as the mesh decreases, and a common-part coupling gives
\begin{equation}\label{AS-eq:commonpart}
 W_p(a\nu,b\nu)^p\le2^{p-1}\int|y|^p|a(y)-b(y)|\,\nu(dy).
\end{equation}
Finally let $\tau\downarrow0$. Fix one such regular kernel while perturbing the marginals.

\emph{Exceptional rows.}
Let $G$ be the bounded union of the positive-mass cells, chosen with $\mu(\partial G)=0$ and $\mu(G)=1$. Take any $\eta_n\in\M(\mu_n,\nu_n)$, put $a_n=\mu_n(G)$, and normalize its restriction to these rows. Its marginals $\widetilde\mu_n,\widetilde\nu_n$ converge in $W_p$ to $\mu,\nu$. Indeed $a_n\to1$; the removed terminal submeasure is dominated by $\nu_n$, has vanishing mass, and has vanishing $p$-moment by uniform integrability. Outside $G$ we will keep $\eta_n$ unchanged.

\emph{Move the functions and correct all constraints.}
Couple $Y\sim\nu$ and $Z\sim\widetilde\nu_n$ by $\beta_n$ with $\E|Y-Z|^p\to0$. Define
\begin{align*}
 r_i^n(z)&=\E[r_i(Y)\mid Z=z],& h_n(z)&=\E[h(Y)\mid Z=z],\\
 b_i^n&=\int z r_i^n(z)\,d\widetilde\nu_n,&
 C_n&=\int z h_n(z)^{\mathsf T}\,d\widetilde\nu_n.
\end{align*}
The density bounds survive, $\int h_n\,d\widetilde\nu_n=0$, $\sum_iw_i r_i^n=1$, and boundedness of $r_i,h$ gives $b_i^n\to b_i$, $C_n\to C$. Write $w_i^n=\widetilde\mu_n(C_i)$ and $m_n$ for the common normalized mean. Set
\begin{align}
 d_n(z)&=\sum_iw_i^n r_i^n(z)-1+
 h_n(z)^{\mathsf T}C_n^{-1}\left(m_n-\sum_iw_i^n b_i^n\right),\label{AS-eq:column}\\
 k_n(x',z)&=r_i^n(z)+h_n(z)^{\mathsf T}C_n^{-1}(x'-b_i^n)-d_n(z),
 \qquad x'\in C_i.\label{AS-eq:correction}
\end{align}
The column error $d_n$ has zero mass and first moment, and $\|d_n\|_\infty\to0$. Direct integration gives
\begin{equation}\label{AS-eq:identities}
 \int k_n(x',z)\,d\widetilde\nu_n=1,\qquad
 \int z k_n(x',z)\,d\widetilde\nu_n=x',\qquad
 \int k_n(x',z)\,d\widetilde\mu_n=1.
\end{equation}
The fixed small mesh and the original positive lower bound make $k_n$ positive for large $n$; its upper bounds are uniform. Use $k_n\widetilde\nu_n$ on $G$ and restore $\eta_n$ on the exceptional rows. This defines the required exact martingale marginals.

\emph{Respect the prescribed initial coupling.}
Under $\xi_n$, cell mismatches have probability tending to zero, because the cells are continuity sets. On a common cell, transporting $k(X,\cdot)\nu$ through $\beta_n$ costs at most a fixed multiple of $\int|y-z|^p\,d\beta_n$. Its transported density differs from~\eqref{AS-eq:correction} by
\[
 h_n(z)^{\mathsf T}
 \bigl[C_n^{-1}(X'-b_i^n)-C^{-1}(X-b_i)\bigr]-d_n(z).
\]
Its supremum is at most $C'|X-X'|+o(1)$ on these bounded cells. Equation~\eqref{AS-eq:commonpart} makes its integrated transport cost vanish. On mismatched good cells the uniform density bounds control the conditional moments. On exceptional rows the removed terminal moment tends to zero. Hence~\eqref{AS-eq:prescribed} holds for the fixed regular kernel. A diagonal choice over the regularization, using the $L^p$ triangle inequality, proves it for $\pi$.
\end{proof}

\subsubsection{Histories, feasible sets and optimization}
The prescribed coupling in~\eqref{AS-eq:prescribed} allows an already coupled past to remain unchanged. The correction uses the current mean as one observable of that past; it does not require the transition to be Markovian.

\begin{corollary}[History-indexed kernels and finite horizons]\label{AS-cor:history}
Let $H\sim\alpha$, $H_n\sim\alpha_n$ take values in a Polish space, and let $g,g_n$ be measurable $\R^d$-valued observables. Suppose a prescribed coupling satisfies
\[
 \E[d(H,H_n)^p+|g(H)-g_n(H_n)|^p]\to0.
\]
If $g_\#\alpha,\nu$ satisfy Theorem~\ref{AS-thm:main}, $\nu_n\to\nu$ in $W_p$, and $(g_n)_\#\alpha_n\cx\nu_n$, then every terminal kernel $\pi_h$ of mean $g(h)$ and marginal $\nu$ admits kernels $\pi^n_{h'}$ of mean $g_n(h')$ and marginal $\nu_n$ with $\E W_p^p(\pi_H,\pi^n_{H_n})\to0$.

Consequently, if each adjacent pair of $\boldsymbol\mu=(\mu_0,\ldots,\mu_T)$ satisfies Theorem~\ref{AS-thm:main} and $\mu_t^n\to\mu_t$ in $W_p$ with $\mu_t^n\cx\mu_{t+1}^n$, every martingale law $P$ with these one-time marginals has approximants $P_n$ with marginals $\boldsymbol\mu^n$ and $AW_p(P_n,P)\to0$. Here pathwise $AW_p$ minimizes the sum of coordinate $p$-costs over bicausal couplings. No Markov hypothesis is imposed.
\end{corollary}
\begin{proof}
Lift the input to $(H,g(H))$. A product-minorized martingale kernel lifts by depending only on $g(H)$. Average on finite continuity partitions of this lifted space which generate its sigma-field and have small $g$-diameter. Formulas~\eqref{AS-eq:regular}--\eqref{AS-eq:identities} apply with $x=g(h)$, $x'=g_n(h')$ and $b_i$ the cell average of $g$. A feasible exceptional-row kernel is obtained from any martingale for $(g_n)_\#\alpha_n,\nu_n$. The same estimates prove the first assertion. For the second, couple the initial states and apply the first assertion successively to the entire coupled histories, taking $g$ to be their last coordinate. Measurable conditional optimal couplings extend the past coupling. Each of their marginals depends only on its own history, so the extension is bicausal. Induction over the fixed finite horizon proves convergence.
\end{proof}

For optimizer convergence, the conditional law must be retained through the limit. Nearby initial points can carry different terminal kernels even when their joint laws converge. Record the kernel as a random measure:
\[
 \Lambda_n(dx,dr)=\mu_n(dx)\delta_{\pi_x^n}(dr).
\]
A limit may have the form $\mu(dx)\Lambda_x(dr)$, with several kernels $r$ at the same $x$. Its averaged kernel $\bar r_x=\int r\,\Lambda_x(dr)$ still has mean $x$ and retains the terminal marginal. A concave objective improves when this unresolved choice is averaged; strict concavity makes the improvement strict whenever $\Lambda_x$ is nonconstant. Convergence of the optimal values therefore forces the limiting random kernel to be a point mass. The lemma proves this implication and then upgrades it to convergence along any prescribed initial coupling.
\begin{lemma}[Limits of conditional kernels]\label{AS-lem:lift}
Suppose $\nu_n\to\nu$ in $W_p$ and the laws $\Lambda_n$ on $\R^d\times\cP_p(\R^d)$ have first marginals converging in $W_p$ to $\mu$ and averaged second components $\nu_n$. They are tight with uniformly integrable $p$-moments for the product Wasserstein metric. If
$\int|\bary(r)-x|^p\,\Lambda_n(dx,dr)\to0$, every limit $\mu(dx)\Lambda_x(dr)$ has $\bary(r)=x$ and averaged terminal law $\nu$. Its barycentric kernel $\bar r_x=\int r\,\Lambda_x(dr)$ is a martingale kernel. If a continuous concave objective with marginally controlled $p$-growth attains its limiting maximum only at $\pi^*$ and is strictly concave in $r$, every asymptotically maximizing sequence with convergent values has limit $\mu(dx)\delta_{\pi_x^*}(dr)$. For graph laws $\Lambda_n=\mu_n(dx)\delta_{\pi_x^n}(dr)$, the conditional kernels converge in integrated $W_p^p$ along any prescribed initial coupling of vanishing $p$-cost.
\end{lemma}
\begin{proof}
Uniform terminal tails and Markov's inequality, applied at a sequence of radii with summable tail bounds, confine the random measures to $W_p$-compact sets with arbitrarily high probability. A superlinear convex function integrable uniformly against $|Y|^p$ and conditional Jensen give uniform integrability of $m_p(r)$. Barycenters and bounded terminal tests pass to the limit. Write $\mathcal J_x(r)$ for the limiting objective at $x$. Concavity gives
\[
 \int\!\int\mathcal J_x(r)\,\Lambda_x(dr)\,\mu(dx)
 \le\int\mathcal J_x(\bar r_x)\,\mu(dx)
 \le\sup_{\pi\in\M(\mu,\nu)}\int\mathcal J_x(\pi_x)\,\mu(dx).
\]
The left side is the limiting optimal value by the growth and uniform-integrability assumptions. Both inequalities are therefore equalities. To see why the first equality forces a point mass, partition any nonconstant random kernel using a bounded continuous terminal test. Its conditional barycenters on the two parts are distinct, and strict two-point concavity makes the first inequality strict. Thus $\Lambda_x=\delta_{\bar r_x}$ almost everywhere; uniqueness of the optimizer gives $\bar r_x=\pi_x^*$. Restrict $x\mapsto\pi_x^*$ to compact Lusin sets. Graph convergence and convergence of the coupled initial states then imply conditional convergence in probability along any prescribed coupling. The uniform moment bound upgrades this to integrated $W_p^p$.
\end{proof}

\begin{corollary}[Feasible sets and continuous or weak costs]\label{AS-cor:optimization}
Under the finite-horizon hypotheses of Corollary~\ref{AS-cor:history}, the full sets of martingale laws converge in ordinary path-space $W_p$ Hausdorff distance. Minima and maxima of every continuous path cost bounded in absolute value by $C(1+\sum_t|x_t|^p)$ converge; subsequential limits of optimizers are optimal.
For the two-time hypotheses, if $C(x,r)$ is continuous, convex in $r$, and
\[
 |C(x,r)|\le A(1+|x|^p+m_p(r)),
\]
then $\inf_{\pi\in\M(\mu_n,\nu_n)}\int C(x,\pi_x)\,\mu_n(dx)$ converges to the corresponding infimum at $(\mu,\nu)$. The assertion for suprema holds for concave $C$.
\end{corollary}
\begin{proof}
The finitely many marginal sequences give compactness in ordinary $W_p$ and close the martingale equations. This proves the upper Hausdorff inclusion; adapted lower approximation and a finite covering of the compact limiting feasible set prove the lower inclusion. Uniform integrability passes continuous costs to the limit. For weak costs, Theorem~\ref{AS-thm:main} gives the upper value bound for the infimum. Lemma~\ref{AS-lem:lift} gives a limit of nearly minimizing kernels; averaging them and applying convexity gives the lower bound. Negate the cost for the concave case.
\end{proof}

\begin{corollary}[Simultaneous stability of endpoints and reference]\label{AS-cor:three-laws}
Under Theorem~\ref{AS-thm:main} with $p=2$, let also $q_n\to q$ in $W_2$, with $q$ absolutely continuous. The values
\[
 P_{q_n}(\mu_n,\nu_n)=\sup_{\pi\in\M(\mu_n,\nu_n)}
 \int\MCov(\pi_x,q_n)\,\mu_n(dx)
\]
converge to $P_q(\mu,\nu)$, where $\MCov(r,q)$ is the largest expected inner product over couplings of $r,q$. Every asymptotically optimal sequence converges to the unique limiting optimizer in $AW_2$, and satisfies~\eqref{AS-eq:prescribed} along every prescribed converging initial coupling. The approximating references may be atomic. These conclusions also hold at $(\mu,(m+c(\,\cdot-m))_\#\nu)$ for any compactly supported $\mu\cx\nu$, full-span $\nu\in\cP_2$ and $c>1$.
\end{corollary}
\begin{proof}
Maximal covariance is continuous in $W_2$, with
\[
 |\MCov(r,q_n)-\MCov(r,q)|\le\sqrt{m_2(r)}W_2(q_n,q).
\]
For strict concavity, let $r_0\ne r_1$ and choose their optimal transports $T_0,T_1$ from the absolutely continuous law $q$. For $0<t<1$, mixing the two transport plans is feasible from $q$ to $tr_0+(1-t)r_1$, so
\[
 \MCov(tr_0+(1-t)r_1,q)
 \ge t\MCov(r_0,q)+(1-t)\MCov(r_1,q).
\]
Equality would make this mixed plan optimal. Brenier uniqueness says that the optimal plan from $q$ is a map, whereas the mixed plan chooses between $T_0$ and $T_1$ at each source point. It can be a map only if $T_0=T_1$ $q$-almost everywhere, which would give $r_0=r_1$. Hence the inequality is strict. Theorem~\ref{AS-thm:main} gives the lower value bound; Lemma~\ref{AS-lem:lift} gives the upper bound and rules out every nontrivial limiting mixture. Its graph conclusion is precisely the asserted convergence of conditional kernels. For inflation, mix the dilated terminal kernel of any martingale from $\mu$ to $\nu$ with the independent inflated terminal law in proportions $1/c,1-1/c$. Its conditional mean is $x$ and it has the required product component.
\end{proof}

\begin{example}[Adapted upper Hausdorff convergence fails]\label{AS-ex:upper}
Let $\mu=\delta_0$, $A=(\delta_{-1}+\delta_1)/2$, $B=(\delta_{-2}+\delta_2)/2$, $\nu=(A+B)/2$, and $\mu_n=(\delta_{-1/n}+\delta_{1/n})/2$. For $n\ge2$, set
\[
 K_-^n=(1-n^{-1})A+n^{-1}B+\frac{\delta_{-2}-\delta_2}{4n},
 \qquad K_+^n=n^{-1}A+(1-n^{-1})B-\frac{\delta_{-2}-\delta_2}{4n}.
\]
These probabilities have means $-1/n,1/n$ and average $\nu$. They define martingales $\pi_n$ with the perturbed marginals, although the limiting pair is irreducible and strictly interior. Its unique martingale has kernel $\nu$, while
\[
 AW_p(\pi_n,\delta_0\otimes\nu)^p
 =n^{-p}+\tfrac12\bigl(W_p(K_-^n,\nu)^p+W_p(K_+^n,\nu)^p\bigr)
 \longrightarrow\tfrac12\bigl(W_p(A,\nu)^p+W_p(B,\nu)^p\bigr)>0.
\]
The two starting points coalesce, while the choice between $A$ and $B$ survives in the conditional law. For the standard Gaussian driver this choice has an explicit positive objective cost. Write $\varphi,\Phi$ for its density and distribution function, and let $z=\Phi^{-1}(3/4)$. Monotone transport gives
\[
 \MCov(A,\gamma_1)=2\varphi(0),\qquad
 \MCov(B,\gamma_1)=4\varphi(0),\qquad
 \MCov(\nu,\gamma_1)=2\varphi(0)+2\varphi(z).
\]
Indeed the transport to $\nu$ takes values $-2,-1,1,2$ on the four Gaussian quartiles. Continuity of maximal covariance therefore yields
\[
 P_{\gamma_1}(\delta_0,\nu)
 -\lim_{n\to\infty}\tfrac12\{\MCov(K_-^n,\gamma_1)+\MCov(K_+^n,\gamma_1)\}
 =2\varphi(z)-\varphi(0)>0.
\]
These feasible martingales stay a fixed amount below the limiting optimum. Approximate optimizers cannot retain their unresolved kernel choice. Every limiting martingale has adapted approximants, while ordinary Hausdorff convergence and adapted lower approximation remain distinct from adapted upper Hausdorff convergence.
\end{example}


\subsection{Directional couplings: a full \texorpdfstring{$W_\infty$}{W-infinity} approximation theorem}
\label{FP-sec:directional}
For equal-mass finite measures on $\R$, $\alpha\st\beta$ means $\alpha(({-\infty},t])\ge\beta(({-\infty},t])$ for every $t$. Equivalently, their increasing quantiles are ordered. This convention gives couplings $X\le Y$.

\begin{lemma}[Upper-quantile submeasure]\label{FP-lem:upper-slice}
Let $\nu$ be an atomless finite measure on an interval, and let $0\le m\le\nu(\R)$. Its uppermost submeasure $\nu^{\mathrm{up}}_m$ of mass $m$ stochastically dominates every submeasure of $\nu$ having mass $m$. If $\alpha\st\nu^{\mathrm{up}}_m$, there is a measurable map from $\nu^{\mathrm{up}}_m$ to $\alpha$ satisfying $T(y)\le y$ almost everywhere.
\end{lemma}
\begin{proof}
Choose the upper quantile cut, possible without splitting an atom. For each $t$, this submeasure has upper-tail mass $\min\{m,\nu((t,\infty))\}$, the largest possible value. This proves stochastic domination. The increasing quantile coupling is induced from the atomless source by its continuous distribution function. Ordered quantiles give $T(y)\le y$ almost surely.
\end{proof}

\begin{theorem}[Directional backward Monge density]\label{FP-thm:directional-density}
Let $\mu,\nu$ be probabilities on $\R$ with $\mu\st\nu$ and $\nu$ atomless. For every coupling $\pi$ supported on $\{(x,y):x\le y\}$ and every $\eta>0$, there is a Borel map $T$ such that
\[
 T_\#\nu=\mu,\qquad T(y)\le y\quad\nu\text{-a.e.},
 \qquad W_\infty\bigl(\pi,(T,\mathrm{id})_\#\nu\bigr)\le\eta.
\]
Thus Conjecture~5.4 of Nutz--Wang--Zhang holds, with $W_\infty$ density in place of weak density.
\end{theorem}
\begin{proof}
Choose a grid $I_k=[kh,(k+1)h)$ with $\sqrt2h\le\eta$. Let $\pi_{k\ell}$ be the restriction of $\pi$ to $I_k\times I_\ell$, with marginals $\mu_{k\ell},\nu_{k\ell}$ and common mass $m_{k\ell}$. Only $k\le\ell$ can have positive mass. For each $\ell$, the source restriction $\nu_\ell=\nu|_{I_\ell}$ decomposes as $\sum_{k\le\ell}\nu_{k\ell}$.

Assign to the diagonal label $k=\ell$ the uppermost submeasure $\widetilde\nu_{\ell\ell}$ of $\nu_\ell$ with mass $m_{\ell\ell}$. The old diagonal coupling gives $\mu_{\ell\ell}\st\nu_{\ell\ell}$, and Lemma~\ref{FP-lem:upper-slice} gives $\mu_{\ell\ell}\st\nu_{\ell\ell}\st\widetilde\nu_{\ell\ell}$. The remainder of $\nu_\ell$ is atomless. Partition it measurably into submeasures $\widetilde\nu_{k\ell}$, $k<\ell$, of the prescribed masses $m_{k\ell}$. Such a countable partition is obtained from a uniform quantile coordinate and intervals of those lengths.

For the diagonal pair, use the ordered quantile map from $\widetilde\nu_{\ell\ell}$ to $\mu_{\ell\ell}$. For $k<\ell$, use any quantile map from $\widetilde\nu_{k\ell}$ to $\mu_{k\ell}$. Its outputs satisfy $x\in I_k$, $y\in I_\ell$, hence $x<y$. The source pieces are mutually singular, so these countably many maps combine into a Borel map $T$. Its output marginal is $\sum_{k,\ell}\mu_{k\ell}=\mu$ and it satisfies the direction constraint.

The old and new couplings give exactly mass $m_{k\ell}$ to each rectangle $I_k\times I_\ell$. Couple their normalized restrictions independently inside that rectangle. Every such pair of points has distance at most $\sqrt2h$. Summing these couplings gives the asserted $W_\infty$ estimate.
\end{proof}

\begin{corollary}[No relaxation gap for directional transport]\label{FP-cor:directional-cost}
For every bounded continuous cost $c$, the infimum over directional couplings equals the infimum over directional backward Monge couplings. Every admissible plan can be approximated in value by maps with the same marginals and exact direction constraint.
\end{corollary}
\begin{proof}
Apply Theorem~\ref{FP-thm:directional-density}; $W_\infty$ convergence implies weak convergence, and the cost is bounded and continuous.
\end{proof}
The theorem approximates each admissible plan and its bounded continuous costs. The initial marginal is arbitrary. The upper-quantile choice keeps the only order-sensitive rectangles feasible.


\subsection{Supermartingale couplings: preserving the equality region}
\label{FP-sec:supermartingale}
For finite-first-moment laws on $\R$, write $\mu\cd\nu$ when
$\int\phi\,d\mu\le\int\phi\,d\nu$ for every convex decreasing
$\phi$. Equivalently there is a coupling with $\E[Y\mid X]\le X$.
Write $\SM(\mu,\nu)$ for these couplings and $\SMm(\mu,\nu)$ for
those induced by a backward map.

The mean-independent factor in Lemma~\ref{FP-lem:mean-independent} handles every region with a strictly positive supermartingale gap. The region on which the gap vanishes requires the one-dimensional martingale theory. We use the following two established inputs from \cite{NWZ2024}.

\begin{proposition}[Nutz--Wang--Zhang inputs]\label{FP-prop:NWZ-inputs}
The following statements hold on $\R$.
\begin{enumerate}[label=(\roman*),leftmargin=2em]
\item If $\alpha\cx\beta$ and $\beta$ is atomless, backward Monge martingale couplings are weakly dense among all martingale couplings with these marginals \cite[Theorem~2.3]{NWZ2024}.
\item If an atomless finite measure $\beta$ is written as a countable sum $\sum_{i\ge0}\beta_i$, there are mutually singular submeasures $\widetilde\beta_i$ summing to $\beta$, with the same individual masses, such that
\[
 \beta_0\cx\widetilde\beta_0,
 \qquad
 \int y\,d\widetilde\beta_i=\int y\,d\beta_i\quad(i\ge1).
\]
This is the mass-preserving construction of \cite[Lemma~4.10]{NWZ2024}, with index $0$ distinguished.
\end{enumerate}
\end{proposition}
The statements for finite subprobabilities follow by normalization. If the distinguished measure is zero, ordinary countable purification suffices for (ii). Individual mass preservation in (ii) also follows directly from the cited proof: it replaces the other pieces by atoms of their original masses and barycenters before taking shadows. In a restriction to an interval, every new piece stays in that interval.

\begin{lemma}[Uniform integrability with fixed marginals]\label{FP-lem:fixed-marginal-Wp}
If $\pi_n,\pi$ have fixed marginals with finite $p$th moments, $1\le p<\infty$, and $\pi_n$ converges weakly to $\pi$, then $W_p(\pi_n,\pi)\to0$.
\end{lemma}
\begin{proof}
The $p$th powers of the joint norms are uniformly integrable because the two marginals are fixed. For example, the sum $(|x|+|y|)^p$ outside a radius $R$ is bounded by a constant depending only on $p$ times $|x|^p\ind_{\{|x|>R/2\}}+|y|^p\ind_{\{|y|>R/2\}}$. Weak convergence together with this uniform integrability is the standard characterization of $W_p$ convergence.
\end{proof}

\begin{theorem}[Supermartingale backward Monge density]\label{FP-thm:supermartingale-density}
Let $\mu,\nu$ have finite first moments on $\R$, assume $\mu\cd\nu$, and suppose $\nu$ is atomless. Then backward Monge supermartingale couplings are $W_1$-dense in $\SM(\mu,\nu)$. If both marginals have finite $p$th moments, the density holds in $W_p$ for every $1\le p<\infty$. This proves Conjecture~5.3 of Nutz--Wang--Zhang.
\end{theorem}
\begin{proof}
Fix $\pi\in\SM(\mu,\nu)$ and a small $h>0$. Choose a measurable version
\[
 b(x)=\E_\pi[Y\mid X=x],\qquad d(x)=x-b(x)\ge0.
\]
Let $E_0=\{d=0\}$. Partition $\{d>0\}$ into countably many Borel sets $E_i$, $i\ge1$, as follows. For each integer $r$, intersect the set $2^{-r}\le d<2^{-r+1}$ with a half-open grid of length at most $\min\{h,2^{-r-1}\}$. Remove null pieces and enumerate the rest. Every $E_i$ lies in an interval $[a_i,a_i+h_i)$ with $h_i\le h$, and
\begin{equation}\label{FP-eq:gap-partition}
 b(x)\le a_i\le x\quad\text{for }\mu\text{-a.e. }x\in E_i.
\end{equation}
Indeed, on that interval $b(x)=x-d(x)\le a_i+h_i-2^{-r}\le a_i$.

Let $\pi_i=\pi|_{\{X\in E_i\}}$ have marginals $\mu_i,\nu_i$ and mass $m_i$. The restriction $\pi_0$ is a martingale coupling. For $i\ge1$ put
\[
 \beta_i=\frac1{m_i}\int y\,d\nu_i
 =\frac1{m_i}\int_{E_i}b(x)\,d\mu(x)\le a_i.
\]

Partition the later-coordinate line into intervals $C_\ell$ of length $h$. On each restriction $\nu|_{C_\ell}$ apply Proposition~\ref{FP-prop:NWZ-inputs}(ii) to the decomposition $\sum_{i\ge0}\nu_i|_{C_\ell}$, distinguishing $i=0$. It gives mutually singular pieces $\widetilde\nu_{i\ell}$ with
\begin{align}
 \widetilde\nu_{i\ell}(\R)&=\nu_i(C_\ell),\label{FP-eq:local-mass}\\
 \int y\,d\widetilde\nu_{i\ell}&=\int_{C_\ell}y\,d\nu_i\quad(i\ge1),\label{FP-eq:local-mean}\\
 \nu_0|_{C_\ell}&\cx\widetilde\nu_{0\ell}.\label{FP-eq:local-cx}
\end{align}
Put $\widetilde\nu_i=\sum_\ell\widetilde\nu_{i\ell}$. All these measures are atomless and mutually singular, and they sum to $\nu$.

For $i\ge1$, apply Lemma~\ref{FP-lem:mean-independent} to $Y\sim\widetilde\nu_i/m_i$, with prescribed output law $\mu_i/m_i$. It yields a map $T_i$ with
\[
 (T_i)_\#\widetilde\nu_i=\mu_i,
 \qquad\E[Y\mid T_i(Y)]=\beta_i\le a_i\le T_i(Y).
\]
The mean is $\beta_i$ because~\eqref{FP-eq:local-mean} preserves the full first moment. Thus each resulting $\widehat\pi_i=(T_i,\mathrm{id})_\#\widetilde\nu_i$ is a supermartingale subcoupling.

For the equality region, take martingale kernels transporting $\nu_0|_{C_\ell}$ to $\widetilde\nu_{0\ell}$, whose existence follows from~\eqref{FP-eq:local-cx}. They move the later coordinate by at most $h$. Composing them with $\pi_0$ gives a martingale coupling $\pi'_0$ from $\mu_0$ to $\widetilde\nu_0$ and a coupling of $\pi_0,\pi'_0$ moving only the later coordinate by at most $h$. By Proposition~\ref{FP-prop:NWZ-inputs}(i) and Lemma~\ref{FP-lem:fixed-marginal-Wp}, choose a backward Monge martingale coupling $\widehat\pi_0$ from $\mu_0$ to $\widetilde\nu_0$ at normalized $W_p$ distance at most $h$ from $\pi'_0$. Skip this step when $m_0=0$.

The maps just constructed have disjoint source pieces, so they combine into one map $T$. Its marginals are exactly $\mu,\nu$ and its conditional-mean inequality holds on each $E_i$, hence globally.

It remains to verify approximation to the original plan. For each $i\ge1$ and $\ell$, the restrictions of $\pi_i$ and $\widehat\pi_i$ to $\{Y\in C_\ell\}$ have equal mass by~\eqref{FP-eq:local-mass}; both first coordinates lie in $E_i$, of diameter at most $h$, and both second coordinates lie in $C_\ell$. Coupling the normalized restrictions arbitrarily gives distance at most $\sqrt2h$. The equality-region composition costs at most $h$, followed by at most $h$ in $W_p$. Taking the mixture of these couplings shows that the global $W_p$ distance is at most $2h$. Let $h\downarrow0$. For the asserted $W_1$ result only first moments were used.
\end{proof}

\begin{corollary}[Cost approximation under the supermartingale constraint]\label{FP-cor:super-cost}
For every bounded continuous cost $c$, the infimum over $\SM(\mu,\nu)$ equals the infimum over $\SMm(\mu,\nu)$. Both marginals and the supermartingale inequality remain exact along an approximating sequence.
\end{corollary}
\begin{proof}
Use Theorem~\ref{FP-thm:supermartingale-density} and bounded-continuous convergence.
\end{proof}

\begin{remark}[The equality and strict-slack regions]
The equality region uses the existing martingale theorem. On a region with positive gap, the extra space below $x$ makes a constant conditional mean admissible. The mean-independent factor realizes an arbitrary distribution of earlier coordinates inside that short region, while keeping this constant mean. The one distinguished convex-order piece in the splitting lemma lets the two constructions use mutually singular parts of the same unchanged later marginal. No higher-dimensional martingale-density assertion enters the argument.
\end{remark}


\section{Canonical partitions with prescribed polynomial moments}
\label{R11-sec:polynomial}
A selecting cost makes an exact-moment partition canonical. Positivity of the feasible kernel gives dual attainment, and the selecting degree prevents cells from sharing a positive-measure equality set. The theorem determines the optimal degree in dimension at least three and exact-moment entropy regularization. For quadratic interfaces, their geometry also determines the leading cost of residual randomness.

\subsection{Canonical partitions with exact polynomial moments}
\label{R9-sec:canonical}
Let $\Omega\subset\R^D$ be a bounded open convex box, and let $\rho$
be a probability measure with density positive almost everywhere on
$\Omega$. For $r\ge1$, let $v_r(z)$ list all monomials of total degree
at most $r$, including the constant. Prescribe vectors $b_i$,
$1\le i\le q$, whose constant coordinates are $p_i>0$, and assume
that a kernel $k^0$ satisfies
\begin{equation}\label{R9-eq:positive-moments}
 \sum_i k_i^0=1,\qquad k_i^0\ge\eta p_i\quad\rho\text{-a.e.},
 \qquad \int v_r k_i^0\,d\rho=b_i
\end{equation}
for some $\eta>0$. In particular, $\sum_i p_i=1$.

\begin{theorem}[Degree-optimal canonical moment partitions]
\label{R9-thm:canonical-moments}
Suppose $D\ge3$. Choose distinct real numbers $\lambda_i$ and put
\[
 P(z)=\sum_{j=1}^D z_j^{r+1},\qquad c_i(z)=\lambda_iP(z).
\]
Among all kernels satisfying the moments in
\eqref{R9-eq:positive-moments}, there is a unique minimizer of
$\mathcal C(k)=\int\sum_i c_i k_i\,d\rho$. It is deterministic:
$k_i=\one_{E_i}$, where, up to null sets,
\begin{equation}\label{R9-eq:polynomial-cells}
 E_i=\{z\in\Omega:(\theta_i-\theta_j)\cdot v_r(z)
                   >(\lambda_i-\lambda_j)P(z)\quad(j\ne i)\}.
\end{equation}
The multipliers are unique modulo addition of a common polynomial of
degree at most $r$; the gauge $\sum_i p_i\theta_i=0$ selects one.
Every cell is described by at most $q-1$ inequalities of degree $r+1$.
This selecting degree is the least possible for a universal construction
with uniqueness over all feasible kernels.

For every $\varepsilon>0$, minimizing
$\mathcal C(k)+\varepsilon I(Z;I)$ with the same exact moments has a
unique kernel. Its gauge-fixed dual converges to $\theta$, and its
kernel converges in $L^1(\rho)$ to the partition as
$\varepsilon\downarrow0$.
\end{theorem}

The source is bounded, so the selecting polynomial need not be even.
Its highest homogeneous part is irreducible for every degree $r+1\ge2$
in at least three variables. This gives degree $r+1$ without a parity
restriction. The general purification theorem and moment-constrained
transport duality have classical and contemporary forms
\cite{DWW,KhanRath2009,CarlierMalamutSylvestre}; the cell formula
and multiplier uniqueness here come from this particular cost.

\begin{proof}
Write $\mathcal H=\{\theta:\sum_i p_i\theta_i=0\}$ and consider
\[
 J_0(\theta)=\int\max_i\{\theta_i\cdot v_r-c_i\}\,d\rho
                 -\sum_i\theta_i\cdot b_i.
\]
For $w_i=\theta_i\cdot v_r$, the gauge gives $\sum_i p_iw_i=0$
pointwise. The positive feasible kernel implies
\begin{align*}
 J_0(\theta)
 &\ge\int\left(\max_iw_i-\sum_i k_i^0w_i\right)d\rho
                  -\max_i\|c_i\|_\infty\\
 &\ge\eta\int\max_i w_i\,d\rho-\max_i\|c_i\|_\infty.
\end{align*}
The last homogeneous functional is positive on every nonzero direction
in $\mathcal H$. Indeed, its integrand is nonnegative; if the integral
is zero, all $w_i$ vanish almost everywhere, and positivity of the
density makes each polynomial zero. Compactness of the unit sphere in
this finite-dimensional gauge gives coercivity of $J_0$ and a minimizer.

Pairwise score differences have nonzero leading term
$(\lambda_j-\lambda_i)P$, so their zero sets are $\rho$-null.
Differentiating the maximum therefore gives
$\int v_r\one_{E_i}\,d\rho=b_i$. The pointwise inequality
\[
 \sum_i(c_i-\theta_i\cdot v_r)k_i
 \ge\min_i(c_i-\theta_i\cdot v_r)
\]
shows that this partition minimizes cost. Equality forces the unique
winning label almost everywhere, proving uniqueness among all kernels.

It remains to identify the dual. The homogeneous polynomial $P$ is
irreducible over $\C$: its projective hypersurface is smooth because
its partial derivatives vanish simultaneously only at the zero vector;
two nonconstant homogeneous factors in at least three variables would
have a common projective zero and create a singularity. Any
factorization of a pairwise tie polynomial would factor its leading
homogeneous part, so that tie polynomial is also irreducible.
A nonempty smooth real patch of an irreducible real hypersurface is
Zariski dense in it. Consequently, a polynomial of degree at most $r$
that vanishes on such a patch is zero: otherwise the two polynomials
would have a common hypersurface component, forcing divisibility by
the degree-$r+1$ tie polynomial.

The graph joining labels along smooth active interfaces is connected.
To see this, join interior points of their positive-measure open cells
by a generic path in the convex box. Singular sets of an irreducible
tie hypersurface, and intersections of nonproportional tie
hypersurfaces, have codimension at least two and can be avoided.
Proportional triple ties would make the score whose $\lambda$ is
intermediate a convex combination of the other two scores everywhere;
that label could have no positive winning cell. Thus the path crosses
only regular two-label interfaces. Two optimal duals give the same
partition. On each active patch the difference of their corresponding
score differences is a degree-at-most-$r$ polynomial that vanishes
there. It is zero. Connectedness makes all parameter differences
common, and the gauge removes the common polynomial.

The entropy-regularized dual is
\[
 J_\varepsilon(\theta)=\varepsilon\int
 \log\sum_i p_i e^{(\theta_i\cdot v_r-c_i)/\varepsilon}\,d\rho
                 -\sum_i\theta_i\cdot b_i.
\]
Its difference from $J_0$ is uniformly bounded in absolute value by
$\varepsilon|\log p_{\min}|$. Hence minimizers exist and remain in a
common compact set as $\varepsilon\downarrow0$. The Hessian in any
nonzero gauge direction is the integral of a strictly positive
conditional variance, since a nonzero polynomial cannot vanish
almost everywhere. Thus the dual minimizer is unique. Its derivative
gives all moments exactly, and its kernel is the associated softmax.
Uniform convergence of the duals and uniqueness of the hard minimizer
give convergence of the multipliers. Away from the null tie sets the
softmax converges to the winning label, so dominated convergence gives
$L^1$ convergence.

Finally, if every $c_i$ has degree at most $r$, its integral against
$k_i$ is fixed by $b_i$. The positive feasible kernel in
\eqref{R9-eq:positive-moments} and a deterministic purification are
distinct minimizers when $q\ge2$. Thus degree at most $r$ cannot give
a universally unique minimizer; degree $r+1$ attains that lower bound.
\end{proof}

\subsubsection{Quadratic interfaces and the cost of residual randomness}
For first moments, the construction works already in dimension two.
Curvature also determines the second-order cost of keeping randomness
in a decoder while preserving the moments exactly.

\begin{theorem}[Quadratic dual rigidity]
\label{R9-thm:quadratic-rigidity}
Let $D\ge2$, let $\rho$ have a bounded density $f$ continuous on
$\overline\Omega$ and positive in $\Omega$, and assume
\eqref{R9-eq:positive-moments} with $v(z)=(1,z)$.
For distinct $\lambda_i$ and $Q\succ0$, take
$c_i(z)=\lambda_i z^{\mathsf T}Qz$.
The cost-minimizing partition and its gauge-fixed dual are unique.
Entropy regularization preserves masses and barycenters exactly and
converges to them. If $s_i=\theta_i\cdot v-c_i$ and
\[
 F_{ij}=\{z\in\Omega:s_i(z)=s_j(z)>s_l(z)\ (l\ne i,j)\},
\]
then the hard dual has Hessian
\begin{equation}\label{R9-eq:interface-Hessian}
 D^2J_0(\theta)[u,u]
 =\sum_{i<j}\int_{F_{ij}}
 \frac{((u_i-u_j)\cdot v(z))^2f(z)}{|\nabla(s_i-s_j)(z)|}
 \,d\mathcal H^{D-1}(z),
\end{equation}
which is positive definite on the gauge.
\end{theorem}
\begin{proof}
The coercivity, derivative and pointwise-duality arguments above use
only that the quadratic costs have higher degree than the affine
constraints. After the change of variables $z\mapsto Q^{1/2}z$, each
pairwise score difference is a definite quadratic plus an affine
function. Both cells have positive mass, so the difference takes both
signs; its zero set is a sphere of positive radius. A nonzero affine
function cannot vanish on an open spherical patch in dimension at
least two. The active-interface path argument therefore proves dual
uniqueness, including in dimension two. Softmax convergence follows
as before.

For a small change of the multipliers, the normal displacement of a
regular interface is $(u_i-u_j)\cdot v/|\nabla(s_i-s_j)|$.
Differentiating the cell moment integrals by coarea gives
\eqref{R9-eq:interface-Hessian}. Distinct tie spheres intersect in
codimension at least two; coincident active triple ties are excluded
by the intermediate-$\lambda$ argument. These exceptional sets have
zero interface measure. The integrals are finite, since the spheres
have positive radii and bounded surface area in the box. A zero
quadratic form makes every affine difference vanish on every active
patch; connectedness and the gauge then give $u=0$.
\end{proof}

\begin{theorem}[Exactly calibrated cost--entropy asymptotics]
\label{R9-thm:thermal-frontier}
Under the hypotheses of Theorem~\ref{R9-thm:quadratic-rigidity}, with
$q\ge2$, let $k^\varepsilon$ minimize
$\mathcal C(k)+\varepsilon I(Z;I)$ under the prescribed masses and
barycenters. Write $C_0$ for the deterministic optimum and put
\[
 \Lambda=\sum_{i<j}\int_{F_{ij}}
       \frac{f(z)}{|\nabla(s_i-s_j)(z)|}\,d\mathcal H^{D-1}(z)>0.
\]
Then, as $\varepsilon\downarrow0$,
\begin{align}
 \mathcal C(k^\varepsilon)-C_0
   &=\frac{\pi^2}{6}\Lambda\varepsilon^2+o(\varepsilon^2),
       \label{R9-eq:thermal-cost}\\
 h_\varepsilon:=H(I\mid Z)
   &=\frac{\pi^2}{3}\Lambda\varepsilon+o(\varepsilon).
       \label{R9-eq:thermal-entropy}
\end{align}
The least feasible cost excess under $H(I\mid Z)\ge h$ is
\begin{equation}\label{R9-eq:thermal-frontier}
 \frac{3}{2\pi^2\Lambda}h^2+o(h^2)\qquad(h\downarrow0).
\end{equation}
\end{theorem}
\begin{proof}
Absorb the priors into the constant multipliers:
$\alpha_i=\theta_i+\varepsilon\log p_i\,e_0$, with a common shift
to restore the gauge. Define the unweighted log-sum-exp functional
\[
 L_\varepsilon(\alpha)=\varepsilon\int
 \log\sum_i e^{(\alpha_i\cdot v-c_i)/\varepsilon}\,d\rho
             -\sum_i\alpha_i\cdot b_i.
\]
Then $J_\varepsilon=L_\varepsilon-\varepsilon H(p)$ and
$L_0=J_0$. We first prove, uniformly in a sufficiently small compact
neighborhood of the hard optimum,
\begin{equation}\label{R9-eq:uniform-layer}
 L_\varepsilon(\alpha)=L_0(\alpha)
             +\frac{\pi^2}{6}\Lambda(\alpha)\varepsilon^2
             +o(\varepsilon^2).
\end{equation}
Positive cell masses persist in that neighborhood. Every pairwise
sphere therefore has radius bounded below, and its gradient on nearby
levels is bounded away from zero. Its level-surface areas in the box
are uniformly bounded. Away from triple ties and the boundary, a
regular two-label interface with gap $g=s_i-s_j$ contributes
$\varepsilon\log(1+e^{-|g|/\varepsilon})$. Coarea, with
$g=\varepsilon t$, gives the coefficient
\[
 \int_{\R}\log(1+e^{-|t|})\,dt=\frac{\pi^2}{6}.
\]
Triple ties and intersections with box faces have zero hypersurface
measure. The inequality
\[
 0\le\log\sum_i e^{(s_i-\max_j s_j)/\varepsilon}
 \le\sum_{i<j}\log(1+e^{-|s_i-s_j|/\varepsilon})
\]
controls their neighborhoods by the same binary coarea integrals.
Compactness and continuity of the spherical level measures make the
omitted surface mass uniformly small. Exponential tails control large
$|g|/\varepsilon$. First take the limit off those neighborhoods and
then shrink them. This proves \eqref{R9-eq:uniform-layer} and
continuity of $\Lambda(\alpha)$.

The regularized minimizers converge to the unique hard dual. Comparing
$L_\varepsilon$ at those two minimizers in both directions shows that
its optimized second-order coefficient is unchanged. Finite-dimensional
duality therefore gives
\[
 V(\varepsilon):=\min_k\{\mathcal C(k)+\varepsilon I(Z;I)\}
 =C_0+\varepsilon H(p)
       -\frac{\pi^2}{6}\Lambda\varepsilon^2+o(\varepsilon^2).
\]
For $\varepsilon>0$, uniqueness and the envelope identity give
$V'(\varepsilon)=I(k^\varepsilon)$. Concavity bounds this derivative
between secants at $\varepsilon(1\pm a)$. Letting first
$\varepsilon\downarrow0$ and then $a\downarrow0$ yields
\[
 V'(\varepsilon)=H(p)-\frac{\pi^2}{3}\Lambda\varepsilon
                    +o(\varepsilon).
\]
Now $h_\varepsilon=H(p)-V'(\varepsilon)$ and
$\mathcal C(k^\varepsilon)=V(\varepsilon)-\varepsilon V'(\varepsilon)$,
which prove \eqref{R9-eq:thermal-cost}--\eqref{R9-eq:thermal-entropy}.
For any feasible kernel with $H(I\mid Z)\ge h$,
\[
 \mathcal C(k)-C_0\ge V(\varepsilon)-C_0-\varepsilon H(p)
                         +\varepsilon h.
\]
Optimize the quadratic lower bound in $\varepsilon$. The regularized
curve attains the matching upper bound; its conditional entropy is
continuous and increases through all sufficiently small positive
values. This proves \eqref{R9-eq:thermal-frontier}.
\end{proof}

The coefficient $\pi^2/6$ is the semi-discrete entropic-transport
boundary-layer constant of Altschuler, Niles-Weed and Stromme
\cite{AltschulerNilesWeedStromme}. Here the dual is recalibrated at
every temperature to retain the labelwise barycenters as well as the
masses. The positive interface Hessian makes that calibration compatible
with the same second-order coefficient. Higher-degree singular
interfaces are covered by the canonical partition theorem, while
Theorem~\ref{R9-thm:thermal-frontier} uses the stated regular quadratic
interfaces.


\section{Exact Gaussian observables and deterministic diffusion limits}
\label{R11-sec:dynamics}
Brown's theorem gives weak-operator approximation by measure-preserving maps~\cite{Brown1966}. Conditional allocation makes the action exact on a finite-dimensional invariant space and all its iterates. Matching four moment orders controls stationary diffusion limits~\cite{EK1986}. The deterministic transfer operator remains a coisometry, so exact finite-mode dynamics can coexist with failure of global hypercontractivity and strict entropy contraction.

\subsubsection{Conditional polynomial moments and Gaussian kernels}\label{SP-sec:history}
\begin{corollary}[Any fixed finite order of conditional moments]
\label{SP-cor:polynomials}
Let $X$ take values in a Polish space, let $Y\in\R^r$, and fix an integer $k\ge1$. Suppose $\E|Y|^k<\infty$ and the conditional law of $Y$ given $X=x$ is absolutely continuous with respect to Lebesgue measure for almost every $x$. For every $\eta>0$, there exists a Borel $T:\R^r\to\mathcal X$ with the entire law of $X$ and such that
\begin{align*}
 \E[Y^\alpha\mid T(Y)]&=g_\alpha(T(Y)),
 &g_\alpha(x)&=\E[Y^\alpha\mid X=x],
 &1\le|\alpha|&\le k,\\
 W_\infty\bigl(\law(T(Y),Y),\law(X,Y)\bigr)&\le\eta.
\end{align*}
Here $\alpha$ is a multi-index and $Y^\alpha=\prod_jY_j^{\alpha_j}$.
\end{corollary}
\begin{proof}
Let $F$ be the finite vector of nonconstant monomials of degrees at most $k$. An affine relation on the conditional support of $F(Y)$ would give a nonzero real polynomial vanishing on the conditional support of $Y$. A nonzero polynomial has a Lebesgue-null zero set. This follows by induction on the dimension: regard the polynomial as a polynomial in the final coordinate, exclude the common zero sets of its nonzero coefficient polynomials, and then use that a nonzero one-variable polynomial has finitely many roots. Absolute continuity rules out such a relation. The source marginal is atomless, so Theorem~\ref{SP-thm:main} applies.
\end{proof}

In particular, the corollary applies to every jointly Gaussian pair for which the conditional covariance of $Y$ given $X$ is positive definite. It preserves the Gaussian source and target marginals and any fixed finite collection of Gaussian conditional moments, although the resulting coupling is backward deterministic. The maps for successively larger orders need not coincide. Requiring all conditional polynomial moments simultaneously would, for moment-determinate Gaussian conditional laws, prescribe the original kernel and can preclude a backward deterministic representation.

\begin{corollary}[A whole Brownian history encoded in its later Gaussian reference]
\label{SP-cor:brownian}
Fix $0<t<T$, an integer $k\ge1$, and $\eta>0$. Let $W$ be a standard Brownian motion in $\R^r$. There is a Borel map
\[
 \mathcal T:\R^r\longrightarrow C([0,t],\R^r)
\]
such that, for $Y\sim N(0,TI_r)$, the path $B=\mathcal T(Y)$ has exactly the Brownian path law on $[0,t]$, and
\[
 \E[Y^\alpha\mid B]
 =\E[(B_t+\sqrt{T-t}\,G)^\alpha\mid B]
 \quad(1\le|\alpha|\le k),
\]
where $G\sim N(0,I_r)$ is independent on the right. Moreover
\[
 W_\infty\bigl(\law(B,Y),\law(W|_{[0,t]},W_T)\bigr)\le\eta
\]
for the maximum of the uniform path metric and Euclidean terminal metric. In particular $\E[Y\mid B]=B_t$.
\end{corollary}
\begin{proof}
Take $X=W|_{[0,t]}$ and terminal $Y=W_T$ in Corollary~\ref{SP-cor:polynomials}. The path space is Polish, and the terminal conditional law is $N(X_t,(T-t)I_r)$, which has a positive density because $t<T$. All required moments are finite.
\end{proof}
The full history up to $t$ is preserved in distribution. The strict time gap $T-t>0$ makes the conditional reference nondegenerate and is required by the allocation argument.



\subsection{Exact finite-dimensional restrictions of Markov operators}
\label{WC-sec:main}
Let $(S,\mu)$ be an atomless standard Borel probability space, equipped with a Polish topology. A Markov kernel $K(x,dy)$ with invariant law $\mu$ defines
\[
 Pf(x)=\int f(y)K(x,dy),\qquad \int Pf\,d\mu=\int f\,d\mu.
\]
For a measurable $\mu$-preserving map $T:S\to S$, its Perron--Frobenius operator $\Lp_T$ is specified by
\begin{equation}\label{WC-eq:transfer}
 \int \Lp_Tf(x)h(x)\mu(dx)=\int f(y)h(T(y))\mu(dy).
\end{equation}
Equivalently, $\Lp_Tf(T(Y))=\E[f(Y)\mid T(Y)]$ for $Y\sim\mu$. The Koopman operator is $\Up_Th=h\circ T$, and $\Lp_T=\Up_T^*$ on $L^2(\mu)$. Both operators are contractions on every $L^p$, $1\le p\le\infty$. We work with specified measurable kernels, so no operator-kernel representation issue arises.

Theorem~\ref{SP-thm:main} realizes this invariant self-coupling on the same probability space. The transfer operator describes conditioning an earlier source state on its later deterministic image, so an invariant observable space retains every iterate.

\begin{theorem}[Finite-dimensional interpolation by a deterministic map]
\label{WC-thm:interpolation}
Let $V\subset L^1(\mu)$ be a real finite-dimensional space and contain the constants. For every $0\le\lambda<1$, there is a $\mu$-preserving measurable map $T$ such that
\begin{equation}\label{WC-eq:interpolation}
 \Lp_T f=\lambda Pf+(1-\lambda)\int f\,d\mu\qquad(f\in V).
\end{equation}
If $1,f_1,\ldots,f_d$ is a basis of $V$ and
\begin{equation}\label{WC-eq:full-support}
 \aff\supp\law\bigl((f_1(Y),\ldots,f_d(Y))\mid X=x\bigr)=\R^d
\end{equation}
for $\mu$-almost every $x$ under $\pi(dx,dy)=\mu(dx)K(x,dy)$, the same statement holds with $\lambda=1$.

For any $\eta>0$, the map in the first assertion can moreover satisfy
\[
 W_\infty\bigl((T,\id)_\#\mu,\lambda\pi+(1-\lambda)\mu\otimes\mu\bigr)\le\eta.
\]
Under \eqref{WC-eq:full-support}, the comparison is directly with $\pi$. On the product, $W_\infty$ uses the maximum coordinate metric.
\end{theorem}
\begin{proof}
Choose the $f_i$ centered. Linear independence together with the constant function implies that their joint law under $\mu$ has full affine span. The kernel $K_\lambda(x,dy)=\lambda K(x,dy)+(1-\lambda)\mu(dy)$ has the same invariant law and contains a positive multiple of $\mu$. Hence its conditional feature laws have full affine span. Apply Theorem~\ref{SP-thm:main} with output law $\mu$, source law $\mu$, feature vector $F=(f_1,\ldots,f_d)$, and profile $g=P_\lambda F$. The resulting exact conditional-mean identities are \eqref{WC-eq:interpolation}, by \eqref{WC-eq:transfer}. The same argument without mixing applies under \eqref{WC-eq:full-support}. Constants are automatically preserved.
\end{proof}

\begin{corollary}[Exact spectral data and all iterates]
\label{WC-cor:iterate}
If $P(V)\subseteq V$, then the map above satisfies
\[
 \Lp_T^n|_V=P_\lambda^n|_V\qquad(n\ge0),
 \qquad P_\lambda=\lambda P+(1-\lambda)\E_\mu.
\]
In particular, for a centered eigenfunction $f\in V$, $Pf=af$ implies $\Lp_T^nf=(\lambda a)^nf$. Under \eqref{WC-eq:full-support}, the eigenvalue and all its iterates are retained without attenuation. The statement retains the full matrix of $P_\lambda|_V$, including its Jordan structure. In the undamped case this is the matrix of $P|_V$.
\end{corollary}
\begin{proof}
The space $V$ is invariant under $P_\lambda$. Repeatedly apply Theorem~\ref{WC-thm:interpolation} on that same space.
\end{proof}
The equalities are identities of functions. In particular, for $f\in V$ and every bounded measurable $h$,
\begin{equation}\label{WC-eq:all-tests}
 \E[f(Y)h(T^nY)]=\int P_\lambda^nf(x)h(x)\mu(dx).
\end{equation}
Thus the comparison retains all later-state tests against each selected earlier observable.

Brown proved weak-operator density of invertible measure-preserving transformations among Markov operators~\cite{Brown1966}. Theorem~\ref{WC-thm:interpolation} imposes exact identities on a specified finite-dimensional space of observables, including nonsimple functions, and these identities persist at every iterate. The construction uses measurable partitions of an atomless space; its conclusion specifies a measurable map, with no additional regularity or evaluation bound.

\subsubsection{Gaussian channels and exact moments at every time horizon}
\label{WC-sec:gaussian}
Write $\Ga_d=N(0,I_d)$, and let $\Pol_k(\R^d)$ denote the polynomials of total degree at most $k$. For a real matrix $A$ with $\|A\|_\op\le1$, define the Gaussian Markov operator
\begin{equation}\label{WC-eq:gaussian-channel}
 P_A f(x)=\E f\bigl(Ax+(I_d-AA^{\mathsf T})^{1/2}G\bigr),
 \qquad G\sim\Ga_d.
\end{equation}
It preserves $\Ga_d$ and each $\Pol_k$. This is the usual Gaussian operator associated with a contraction; its action on homogeneous Gaussian chaoses is the corresponding tensor-power action, with the transpose dictated by the convention in \eqref{WC-eq:gaussian-channel}.

\begin{theorem}[All Gaussian contractions at every finite moment order]
\label{WC-thm:gaussian}
For every contraction $A\in\R^{d\times d}$, every finite integer $k\ge1$, and every $\eta>0$, there is a $\Ga_d$-preserving measurable map $T$ such that
\begin{equation}\label{WC-eq:gaussian-exact}
 \Lp_Tf=P_Af\qquad(f\in\Pol_k(\R^d)).
\end{equation}
The pair $(T(Y),Y)$, $Y\sim\Ga_d$, can be chosen within $\eta$ in $W_\infty$ of the Gaussian coupling with conditional source law \eqref{WC-eq:gaussian-channel}. The result includes singular conditional covariance matrices.
\end{theorem}
\begin{proof}
First take a scalar $s\in[0,1)$. Under the Gaussian coupling with conditional law $N(sx,1-s^2)$, the vector $(Y,Y^2,\ldots,Y^k)$ has full affine support in $\R^k$: a nonzero polynomial cannot vanish on all of $\R$. Theorem~\ref{WC-thm:interpolation} with $\lambda=1$ produces a Gaussian-preserving scalar map $t_{s,k}$ matching all conditional moments up to degree $k$, and arbitrarily close to that coupling in $W_\infty$. For $s=1$, take the identity map.

Write a singular value decomposition $A=U\operatorname{diag}(s_1,\ldots,s_d)V^{\mathsf T}$. For $y\in\R^d$, put
\[
 T(y)=V\bigl(t_{s_1,k}((U^{\mathsf T}y)_1),\ldots,
 t_{s_d,k}((U^{\mathsf T}y)_d)\bigr).
\]
The input coordinates after rotation by $U^{\mathsf T}$ are independent standard Gaussians, and their scalar images are again independent standard Gaussians. Conditional source moments factorize over these coordinate pairs. Consequently every monomial of total degree at most $k$ has precisely the conditional moment of the independent Gaussian scalar channels. Rotating gives conditional mean $Ax$ and covariance $I-AA^{\mathsf T}$, and proves \eqref{WC-eq:gaussian-exact} for all polynomials. Choose each scalar pair displacement at most $\eta/\sqrt d$, couple the pairs independently, and rotate. Each full vector coordinate then moves by at most $\eta$.
\end{proof}

\begin{theorem}[Exact multitime polynomial moments]
\label{WC-thm:multitime}
Let $T$ be as in Theorem~\ref{WC-thm:gaussian}, and set $X_n=T^nG$ for $G\sim\Ga_d$. Let $(Z_n)_{n\ge0}$ be the stationary Gaussian autoregression
\begin{equation}\label{WC-eq:AR}
 Z_{n+1}=A^{\mathsf T}Z_n+(I-A^{\mathsf T}A)^{1/2}\xi_{n+1},
\end{equation}
with independent standard Gaussian innovations. For any $0\le n_1<\cdots<n_r$ and polynomials $q_1,\ldots,q_r$ whose degrees sum to at most $k$,
\begin{equation}\label{WC-eq:mixed-moments}
 \E\prod_{i=1}^r q_i(X_{n_i})=\E\prod_{i=1}^r q_i(Z_{n_i}).
\end{equation}
More generally, the last factor can be an arbitrary bounded measurable function if the degrees of all earlier polynomial factors sum to at most $k$.
\end{theorem}
\begin{proof}
The reverse transition of the stationary Gaussian chain \eqref{WC-eq:AR} is $P_A$. For the deterministic orbit, conditioning an earlier state on a later state is governed by powers of $\Lp_T$. Successively use
\[
 \E[q_1(X_0)q_2(X_n)q_3(X_m)]
 =\int q_3\,\Lp_T^{m-n}\bigl(q_2\Lp_T^nq_1\bigr)\,d\Ga_d.
\]
The analogous formula holds with $P_A$. At each step, the degree is at most the sum of the degrees already encountered. Thus every application lies in $\Pol_k$, where Corollary~\ref{WC-cor:iterate} gives exact equality. The final factor is integrated without another application of the operator. Stationarity handles $n_1>0$. Gaussian marginal moments and H\"older's inequality justify all integrals.
\end{proof}

\begin{corollary}[Hermite eigenfunctions and the Gaussian power spectrum]
\label{WC-cor:hermite}
For any $-1<\rho<1$ and finite $k$, there is a Gaussian-preserving map $T$ such that, for the probabilists' Hermite polynomials $H_j$,
\[
 \Lp_T H_j=\rho^jH_j\quad(0\le j\le k),\qquad
 \E[H_j(X_m)H_\ell(X_n)]=\ind_{\{j=\ell\}}j!\rho^{j(n-m)}
\]
whenever $n\ge m$, $j\le k$, and $\ell\ge0$. In particular, the scalar state has exactly the power spectral density
\begin{equation}\label{WC-eq:PSD}
 s_\rho(\theta)=\frac{1-\rho^2}{2\pi(1-2\rho\cos\theta+\rho^2)},\qquad -\pi\le\theta\le\pi.
\end{equation}
\end{corollary}
\begin{proof}
The generating function $\exp(tx-t^2/2)=\sum_{j\ge0}H_j(x)t^j/j!$ gives $P_\rho H_j=\rho^jH_j$ by one Gaussian integration. Apply \eqref{WC-eq:all-tests}, using truncation for the integrable polynomial test $H_\ell$, and Hermite orthogonality. Summing the absolutely convergent covariance series $\rho^{|n|}$ gives \eqref{WC-eq:PSD}.
\end{proof}
The invariant law, finite-order multitime moments and power spectrum are exact. The transition law has a deterministic direction. This is a finite-spectral realization in the inverse Frobenius--Perron setting; classical work on that problem prescribes invariant densities and studies their associated correlation structure \cite{IFPP1998,IFPP2020}. The map constructed here need not be piecewise smooth, finite-branch, ergodic or numerically tractable.

\begin{corollary}[Fixed-step approximation of the full Gaussian process]
\label{WC-cor:fdd}
Fix $A$, and choose such a map $T_k$ for each $k\to\infty$. For every fixed finite collection of times, the joint law of $(T_k^{n_i}G)_i$ converges weakly, and in every finite $W_p$, to the corresponding stationary Gaussian law in \eqref{WC-eq:AR}.
\end{corollary}
\begin{proof}
Gaussian marginals give tightness and uniform integrability of every fixed-degree monomial. Any subsequential limit has every joint moment of the stated Gaussian vector by Theorem~\ref{WC-thm:multitime}. Gaussian laws are moment determinate, including their possibly degenerate versions; for example every linear projection has Gaussian moments and an exponential-square moment in a neighbourhood of zero. The Cram\'er--Wold theorem identifies the limit. Uniformly bounded higher marginal moments give convergence in every finite $W_p$.
\end{proof}

\subsubsection{Exact finite spectra coexist with maximal information retention}
\label{WC-sec:limits}
The full transfer operator retains an exact right inverse: composing with $T$ and then conditioning on $T$ returns the original observable. Pulling a density back by $T$ therefore preserves its entire divergence when the transfer operator is applied. This mechanism gives coefficient $1$ even for maps that reproduce arbitrarily many Gaussian polynomial modes at every iterate.

\begin{proposition}[Coisometry and information contraction]
\label{WC-prop:coisometry}
For every measure-preserving map $T$ on an atomless probability space,
\[
 \Lp_T\Up_T=I.
\]
Consequently:
\begin{enumerate}[label=(\roman*),nosep]
\item $\|\Lp_T^n\|_{L^2_0\to L^2_0}=1$ for every $n\ge1$;
\item $\|\Lp_T\|_{L^p\to L^q}=\infty$ whenever $1\le p<q\le\infty$;
\item for every convex $\Phi$ with $\Phi(1)=0$ admitting a density of positive finite $\Phi$-divergence, the invariant-law divergence contraction coefficient of $\Lp_T$ equals one.
\end{enumerate}
In particular, the coefficient for relative entropy is one.
\end{proposition}
\begin{proof}
For every $h$, conditioning $h(T(Y))$ on $T(Y)$ gives $h(T(Y))$, which proves $\Lp_T\Up_T=I$. The map $\Up_T$ is an isometry preserving means. Choose a centered $h$ of $L^2$ norm one and apply $\Lp_T^n$ to $h\circ T^n$; contraction gives the reverse inequality and proves (i).

For a measurable set $E$ of mass $a$, put $f=\ind_E\circ T$. Then $\Lp_Tf=\ind_E$, and the norm ratio is $a^{1/q-1/p}$. Atomlessness permits sets with $a\downarrow0$, proving (ii).

For (iii), choose a density $h$ and set $f=h\circ T$. Then $f$ and $h$ have the same distribution under $\mu$, and $\Lp_Tf=h$. Therefore
\[
 \int\Phi(\Lp_Tf)\,d\mu=\int\Phi(f)\,d\mu.
\]
Conditional Jensen gives the upper bound one for every admissible density. Choosing a density with positive finite divergence proves equality. The case $\Phi(s)=s\log s$ is relative entropy.
\end{proof}

\begin{corollary}[Arbitrarily many exact Gaussian modes without hypercontractivity]
\label{WC-cor:no-hyper}
Fix $0<\rho<1$. For every finite $k$, Theorem~\ref{WC-thm:gaussian} gives a Gaussian-preserving map whose transfer operator equals $P_\rho$ on $\Pol_k$ and on all iterates of that space, while every improving $L^p\to L^q$ norm is infinite and its relative-entropy contraction coefficient is one.
\end{corollary}
The Ornstein--Uhlenbeck operator $P_\rho$ satisfies the classical Nelson--Gross estimate $\|P_\rho f\|_q\le\|f\|_p$ when $1<p\le q<\infty$ and $\rho\le\sqrt{(p-1)/(q-1)}$ \cite{Gross1975}. Thus exact agreement on any prescribed finite collection of Gaussian chaoses permits complete failure of global hypercontractivity. Proposition~\ref{WC-prop:coisometry} also precludes a strict uniform $L^2$ contraction on the mean-zero space. It makes no claim about spectral gaps on stronger, regularity-dependent Banach spaces.

\begin{proposition}[Singularity and the infinite-order endpoint]
\label{WC-prop:singular}
For $Y\sim\Ga_d$ and a Gaussian-preserving measurable map $T$, the pair $(T(Y),Y)$ has infinite mutual information and Hirschfeld--Gebelein--R\'enyi maximal correlation one. If $\|A\|_\op<1$, its law is mutually singular with the Gaussian pair of \eqref{WC-eq:gaussian-channel}. Hence the total variation distance between those pair laws is one, even when Theorem~\ref{WC-thm:gaussian} makes their $W_\infty$ distance arbitrarily small.

A single Gaussian-preserving map satisfies $\Lp_Tf=P_Af$ for \emph{all} polynomials if and only if $A$ is orthogonal.
\end{proposition}
\begin{proof}
The graph $\{(x,y):x=T(y)\}$ has full probability under the pair law and zero probability under $\Ga_d\otimes\Ga_d$, by atomlessness. This gives infinite mutual information. Any nonconstant centered square-integrable $h(T(Y))$ gives a correlation-one pair of functions of the two variables. A strictly contracting Gaussian channel has a density relative to the product law, and also gives the graph zero mass; this proves the singularity and total variation statements.

Polynomials are dense in $L^2(\Ga_d)$. Equality on all polynomials would give $\Lp_T=P_A$ on $L^2$. Since $\Lp_T\Lp_T^*=I$, the Gaussian operator would be a coisometry. Its adjoint is $P_{A^{\mathsf T}}$. On the vector of coordinate functions, $P_AP_{A^{\mathsf T}}\id=A^{\mathsf T}Ax$. Thus $A^{\mathsf T}A=I$, so $A$ is orthogonal. Conversely, for orthogonal $A$ take $T(y)=A^{\mathsf T}y$.
\end{proof}
The total-variation and information conclusions concern full distributions. Equation~\eqref{WC-eq:mixed-moments} concerns selected moment data. The preservation statements in the construction keep these two kinds of information distinct.

\begin{corollary}[A uniform distance from compact Gaussian operators]
\label{WC-cor:compact-distance}
For every measure-preserving $T$ on an atomless standard probability space and every compact operator $C$ on $L^2_0(\mu)$,
\[
 \|\Lp_T-C\|_{L^2_0\to L^2_0}\ge1.
\]
In particular, if $\|A\|_\op<1$, every realization in Theorem~\ref{WC-thm:gaussian} satisfies $\|\Lp_T-P_A\|_{L^2_0\to L^2_0}\ge1$, regardless of the number of polynomial orders matched.
\end{corollary}
\begin{proof}
Choose an orthonormal sequence $h_j$ in $L^2_0$ and put $f_j=h_j\circ T$. This is again orthonormal, hence weakly null, so $Cf_j\to0$ in norm by compactness. Meanwhile $\Lp_Tf_j=h_j$. Thus $\|(\Lp_T-C)f_j\|\to1$. For a strictly contracting Gaussian matrix $A$, singular value decomposition on the Gaussian chaoses shows that $P_A$ is Hilbert--Schmidt: its squared singular values sum to $\prod_{i=1}^d(1-s_i^2)^{-1}<\infty$, where $s_i$ are the singular values of $A$. The same holds on the mean-zero subspace.
\end{proof}

\begin{corollary}[Strong operator approximation]
\label{WC-cor:strong}
Every invariant Markov kernel on an atomless standard probability space is a strong operator limit, on every $L^p$ with $1\le p<\infty$, of Perron--Frobenius operators of measure-preserving maps. One sequence can be used for all these $p$.
\end{corollary}
\begin{proof}
Choose increasing finite-dimensional spaces spanned by constants and indicators from a countable generating algebra. Their union is dense in each finite $L^p$. Apply Theorem~\ref{WC-thm:interpolation} with $\lambda_n\uparrow1$. Exact agreement on every fixed finite span and the uniform contraction bound imply strong convergence by approximation.
\end{proof}
The last density statement also follows directly from simple-function partitions. Exact interpolation of nonsimple observables uses the continuous-feature construction; the resulting identities determine the spectral restrictions and path-law limits below.



\subsection{A diffusion limit with only four matched moment orders}
\label{WC-sec:diffusion}
The exact operator identities also give a deterministic approximation to a continuous-time diffusion. Here the matching order stays fixed while the time step decreases.

To identify a diffusion from four conditional moments, reverse each finite stationary orbit. Earlier reverse states are deterministic functions of the current state, so the reverse process has a Markov transition described by the transfer operator. The matched first and second moments identify its generator; the fourth moment controls both the Taylor remainder and path tightness. The unique reverse martingale problem then identifies the limit, and time reversal returns the required forward diffusion.

\begin{theorem}[Stationary Ornstein--Uhlenbeck limit]
\label{WC-thm:OU-limit}
Let $B\in\R^{d\times d}$ satisfy $B+B^{\mathsf T}\preceq0$, and put $D=-(B+B^{\mathsf T})$. For each $h>0$, choose a $\Ga_d$-preserving map $T_h$ satisfying
\[
 \Lp_{T_h}|_{\Pol_4}=P_{\exp(hB^{\mathsf T})}|_{\Pol_4}.
\]
Such maps exist by Theorem~\ref{WC-thm:gaussian}. Let $X^h$ linearly interpolate $T_h^nG$ at the times $nh$, where $G\sim\Ga_d$. On every fixed finite time interval,
\begin{equation}\label{WC-eq:OU-limit}
 X^h\ \Longrightarrow\ X\quad\text{in }C([0,S],\R^d),
\end{equation}
where $X$ is the stationary solution of
\[
 dX_t=BX_t\,dt+D^{1/2}dW_t,\qquad X_0\sim\Ga_d.
\]
The driving randomness of each approximating orbit is contained entirely in its initial Gaussian state.
\end{theorem}
\begin{proof}
The dissipativity assumption gives $\|e^{hB^{\mathsf T}}\|_\op\le1$. The Gaussian process in Theorem~\ref{WC-thm:multitime} is now the exact $h$-skeleton of the stationary Ornstein--Uhlenbeck process with drift $B$. In particular, its increment fourth moments agree with those of the deterministic orbit, for all pairs of grid times. Gaussian covariance formulas and $e^{tB}=I+tB+O(t^2)$ give, uniformly for $0\le i h,j h\le S+1$,
\begin{equation}\label{WC-eq:grid-tightness}
 \E|T_h^iG-T_h^jG|^4\le C_S|i h-j h|^2.
\end{equation}
For completeness, the increment covariance is $2I-e^{tB}-e^{tB^{\mathsf T}}$, with $t=|i-j|h$, and a centered Gaussian vector with covariance $C$ has fourth norm moment $(\tr C)^2+2\tr(C^2)$. This proves \eqref{WC-eq:grid-tightness}.

Linear interpolation retains a bound $\E|X^h_t-X^h_s|^4\le C'_S|t-s|^2$. For $|t-s|\ge h$, use the two adjacent grid endpoints and \eqref{WC-eq:grid-tightness}. For $|t-s|<h$, the increment is a sum of at most two fractions of adjacent increments, with total fraction $(t-s)/h$. Its fourth moment is at most a constant times $h^2((t-s)/h)^4\le C(t-s)^2$. The initial laws are fixed, so the Kolmogorov tightness criterion gives tightness in $C([0,S],\R^d)$.

It remains to identify the limit; fourth moments alone do not identify fixed-step joint laws. Reverse a finite stationary orbit. The reverse chain is Markov with transition $R_h=\Lp_{T_h}$. Indeed, its earlier states in reverse time are later deterministic iterates of its current state, so they provide no further information about the next reverse state. The Gaussian comparison kernel has conditional increment mean and covariance
\[
 (e^{hB^{\mathsf T}}-I)x,
 \qquad I-e^{hB^{\mathsf T}}e^{hB}.
\]
The first four conditional raw moments of $R_h$ agree with those of this kernel. For $f\in C_c^\infty(\R^d)$, Taylor's formula to order two therefore gives
\begin{equation}\label{WC-eq:generator}
 \left\|h^{-1}(R_hf-f)-\mathcal A^*f\right\|_{L^1(\Ga_d)}\longrightarrow0,
 \quad
 \mathcal A^*f=(B^{\mathsf T}x)\cdot\nabla f+\tfrac12D:\nabla^2f.
\end{equation}
To verify the remainder, the integrated third absolute increment moment is at most the three-quarter power of its integrated fourth moment, hence is $O(h^{3/2})$ by \eqref{WC-eq:grid-tightness}. The drift and covariance Taylor errors are $O(h^2)$ times Gaussian-integrable polynomials. Dividing by $h$ proves \eqref{WC-eq:generator}.

For the reversed discrete chain, the process
\[
 f(Z^h_{nh})-f(Z^h_0)-\sum_{j<n}(R_hf-f)(Z^h_{jh})
\]
is a martingale. Stationarity and \eqref{WC-eq:generator} show that replacing the sum by $\int_0^{nh}\mathcal A^*f(Z^h_t)dt$ changes its expectation against every bounded past test by $o(1)$, uniformly on bounded time intervals. Tightness, boundedness and continuity of $\mathcal A^*f$ permit passage to the limit. Every limit solves the martingale problem for the linear Ornstein--Uhlenbeck generator $\mathcal A^*$, with initial law $\Ga_d$. This martingale problem is unique: its solution is $e^{tB^{\mathsf T}}Z_0+\int_0^te^{(t-s)B^{\mathsf T}}D^{1/2}dW_s$, including when $D$ is singular. The standard martingale-problem convergence argument is the one developed in \cite{EK1986}; the estimates needed here are \eqref{WC-eq:grid-tightness}--\eqref{WC-eq:generator}.

The time reversal of this stationary Gaussian process has drift $B$: its covariance at positive lag is the transpose of the reverse covariance. Time reversal is continuous on continuous path space. Choosing the reversed horizon at a grid point changes the desired endpoint by at most $h$, which is negligible by the established tightness estimate. This proves \eqref{WC-eq:OU-limit}.
\end{proof}
In dimension one, $B=-1$ gives $dX_t=-X_tdt+\sqrt2\,dW_t$. Every approximating state is exactly standard normal at every grid time. The result is an existence theorem for measurable maps; it does not assert that the maps have the regularity required for a practical deterministic numerical integrator.

The mechanism extends to other polynomial Markov models when their stationary time-reversed kernels preserve the required polynomial spaces and suitable increment bounds identify a unique limiting martingale problem. Polynomial-process theory describes these finite moment spaces by matrix exponentials \cite{CKRT2012}; the exact deterministic restriction is furnished by Theorem~\ref{WC-thm:interpolation}.

\begin{theorem}[Compact elliptic diffusions from deterministic stationary maps]
\label{WC-thm:compact-diffusion}
Let $M\subset\R^d$ be a compact connected smooth submanifold without boundary, and let $(Z_t)$ be a stationary diffusion with smooth uniformly elliptic generator on $M$. Suppose its invariant law $\mu$ has a smooth strictly positive density. For every $h>0$ there exists a $\mu$-preserving measurable map $T_h:M\to M$ such that
\begin{equation}\label{WC-eq:compact-moments}
 \E[Y^\alpha\mid T_h(Y)=x]=\E[Z_0^\alpha\mid Z_h=x],
 \qquad |\alpha|\le4,\quad Y\sim\mu.
\end{equation}
For every family of maps satisfying these identities, the linear interpolation in $\R^d$ of $T_h^nY$ at times $nh$ converges in $C([0,S],\R^d)$ to $(Z_t)_{0\le t\le S}$ as $h\downarrow0$, for each finite $S$. Each approximating state at each grid time has exactly the invariant law $\mu$.
\end{theorem}
\begin{proof}
The time-reversed stationary diffusion has a smooth elliptic generator $\mathcal B$ and a strictly positive transition density on the connected manifold. Its time-$h$ kernel is $Q_h(x,dy)=\law(Z_0\mid Z_h=x)$. Positivity follows from the parabolic strong maximum principle. Let $V$ be the restrictions to $M$ of ambient polynomials of degree at most four, and choose a basis modulo its exact linear relations. For every $x$, equivalence of $Q_h(x,\cdot)$ with $\mu$ implies full affine span of the nonconstant basis features. The undamped part of Theorem~\ref{WC-thm:interpolation} therefore constructs $T_h$ satisfying \eqref{WC-eq:compact-moments}. No invariance of $V$ under $Q_h$ is required here.

Reverse a finite stationary orbit of $T_h$ and write its transition as $R_h=\Lp_{T_h}$. Smoothness and compactness give uniform diffusion increment bounds
\[
 \begin{aligned}
 \sup_x|\E[Z_0-x\mid Z_h=x]|&\le Ch,\\
 \sup_x\E[|Z_0-x|^2\mid Z_h=x]&\le Ch,\\
 \sup_x\E[|Z_0-x|^4\mid Z_h=x]&\le Ch^2.
 \end{aligned}
\]
for small $h$. These follow, for example, from a smooth ambient stochastic differential equation representing the compact diffusion, its bounded coefficients, and the second and fourth martingale moment estimates. Each expression is a polynomial in the source coordinate of degree at most four, with coefficients depending on $x$. Hence the same bounds hold for $R_h$.

The reversed chain decomposes into a vector martingale and a predictable drift whose increments have norm at most $Ch$. The conditional second moment of each martingale increment is at most $Ch$. For bounded stopping times and a further time interval of length $\theta$, the conditional expected squared martingale displacement is at most $C(\theta+h)$, while the drift displacement is at most $C(\theta+h)$. These bounds give Aldous tightness. Since every state lies in $M$, compact containment is automatic. Moreover, over a fixed horizon, the union bound and the fourth moment estimate give
\[
 \Pp\!\left(\max_{j\le S/h}|Y_{j+1}-Y_j|>\delta\right)
 \le C_S h\delta^{-4}\longrightarrow0.
\]
Thus all subsequential limits have continuous paths, and linear interpolation has the same limits in continuous path space.

For a smooth function on $M$, choose a smooth extension to a neighbourhood of $M$ and then a compactly supported extension on $\R^d$. Taylor expansion to order two, the exact first and second conditional moments, and the uniform fourth moment bound give
\[
 \left\|h^{-1}(R_hf-f)-\mathcal B f\right\|_\infty\longrightarrow0.
\]
The integrated Taylor remainder is bounded by a constant times $h^{3/2}$ uniformly in $x$. The Gaussian calculation in the preceding proof is therefore replaced here by uniform local diffusion estimates. Passing the discrete martingale identities to a subsequential limit identifies the unique martingale problem for the smooth elliptic generator $\mathcal B$ on the compact manifold. Reversing time returns the original stationary diffusion, exactly as in Theorem~\ref{WC-thm:OU-limit}. The martingale-problem and tightness criteria are standard \cite{EK1986}; all estimates required to apply them have been stated above.
\end{proof}
This gives Brownian motion on any compact connected manifold as a limit of stationary deterministic maps on that same manifold. It also covers smooth nonreversible elliptic diffusions with their own invariant densities. Four ambient polynomial orders suffice even when the generator preserves no finite polynomial space. The all-horizon exact moment identities of Section~\ref{WC-sec:gaussian} use the additional invariant-space property.


\Needspace{12\baselineskip}
\part{Gaussian geometry and canonical martingales}
\label{R11-part:Gaussian}
Gaussian domination controls convex containment and continuous martingales. We first recover actual sets and paths from that comparison, then select optimal couplings and determine their stability. Prescribing a general driving law raises the additional question of a common convex potential.

\section{From hitting to convex containment}
\label{R11-sec:hitting}
Banaszczyk's theorem and its subgaussian formulations connect vector balancing to hitting convex sets of large Gaussian measure~\cite{Banaszczyk,DGLN}. Hua--Song--Tudose represent every Gaussian-convex-dominated vector as a sum of three standard Gaussian vectors and resolve Talagrand's convexity conjecture~\cite{HST}. To extract a convex body inside an actual sumset, one must turn a family of hitting statements into a single containment statement at an exact measure level.

Suppose every compactly supported law dominated by a reference puts positive mass in an open set. A separating family of convex tests would produce a convex sublevel set missed by that support. A lower bound on the reference mass of the sublevel set excludes the separation. For a Gaussian, a convex function has median at most its mean, giving the half-measure threshold. The fixed-reference version identifies its exact analogue and the universal $1/e$ threshold within the log-concave class. Combined with the three-Gaussian representation, it gives three undilated summands; compact approximation keeps containment inside the actual Borel sum. The input representation and the extracted geometric statement are identified separately below.

\subsection{Compact Gaussian hitting and exact convex extraction}\label{GC-sec:hitting}
A hitting statement can select one law that controls every convex loss at the same Gaussian scale. The argument uses the convex median inequality, finite-dimensional separation and compactness of probability measures on the permitted support. It applies to compact sets without symmetry.

\begin{lemma}[Gaussian convex median inequality]\label{GC-lem:median}\label{a-lem:median}
If $f:\R^d\to\R$ is finite and convex with $\E|f(G)|<\infty$, $G\sim\gamma_d$, a lower median of $f(G)$ is at most $\E f(G)$. For every $\varepsilon>0$,
\begin{equation}\label{GC-eq:median}
\gamma_d\{f<\E f(G)+\varepsilon\}>\tfrac12.
\end{equation}
\end{lemma}
\begin{proof}
For a nonconstant $f$, set $F(t)=\gamma_d\{f\le t\}$. The inclusion
\[
(1-a)\{f\le s\}+a\{f\le t\}\subseteq\{f\le(1-a)s+at\}
\]
and Ehrhard's inequality~\cite{Ehrhard} show that $\Phi^{-1}\circ F$ is concave wherever it is finite. Consequently the increasing Gaussian quantile rearrangement $T$ of $f$ is convex. If $f(G)$ has an atom at its minimum, extend $T$ constantly below the corresponding quantile. Then $T$ is finite on $\R$, $T(g)$ has the law of $f(G)$ for $g\sim N(0,1)$, and $T(0)$ is a lower median. Jensen gives $T(0)\le\E T(g)=\E f(G)$. Continuity of $T$ provides some $s>0$ with $T(s)<\E f(G)+\varepsilon$; hence the probability in~\eqref{GC-eq:median} is at least $\Phi(s)>1/2$. Constant functions satisfy the assertion immediately. This is Kwapie\'n's median argument~\cite{Kwapien}.
\end{proof}

\begin{theorem}[Compact Gaussian hitting]\label{GC-thm:hitting}\label{a-thm:lossless}
Let $S\subset\R^d$ be compact. Suppose that $S$ meets every bounded open convex set $C$ with $\gamma_d(C)\ge1/2$. Then there exists $\mu\in\cP(S)$ such that $\mu\cx\gamma_d$.
If $S=-S$, it suffices to assume the hitting condition only for symmetric $C$; the resulting $\mu$ may then be symmetric.
\end{theorem}
\begin{proof}
For any finite convex $f$ with finite Gaussian expectation, its negative part is Gaussian integrable because $f$ admits an affine minorant. Lemma~\ref{GC-lem:median} and monotone convergence in the radius show that
\[
C=\{f<\E f(G)+\varepsilon\}\cap B(0,R)
\]
has Gaussian measure at least $1/2$ for all sufficiently large $R$. Thus $C$ meets $S$, and letting $\varepsilon\downarrow0$ yields
\begin{equation}\label{GC-eq:min-test}
\min_{s\in S}f(s)\le\E f(G).
\end{equation}
For a finite list $f_1,\ldots,f_m$ of such functions, form the compact convex set
\[
H=\operatorname{conv}\{(f_j(s)-\E f_j(G))_{j=1}^m:s\in S\}.
\]
If $H$ missed the closed negative orthant, strict separation would produce coefficients $c_j\ge0$, not all zero, for which
\[
\min_{s\in S}\sum_jc_jf_j(s)>\E\sum_jc_jf_j(G),
\]
contradicting~\eqref{GC-eq:min-test}. Therefore some probability measure on $S$ satisfies the given finite list of inequalities. The set of probability measures on compact $S$ is weakly compact, and every such test defines a closed subset, since $f_j$ is continuous on $S$. The finite intersection property therefore yields a measure in all these closed sets, satisfying every test. Tests with infinite Gaussian expectation impose no further restriction. Affine tests ensure that the selected law is centered.

For symmetric $S$, repeat the argument using even convex functions and symmetric sublevel sets, symmetrize each finite-test law, and take a symmetric weak limit. An arbitrary convex test then follows by replacing it with its even part.
\end{proof}

\begin{proposition}[Sharp Gaussian scale]\label{R10-prop:Gaussian-scale}
The scale one in Theorem~\ref{GC-thm:hitting} cannot be lowered uniformly in dimension, even for finite symmetric supports.
\end{proposition}
\begin{proof}
Let $r_d$ be a median of $\|G_d\|_2$ and $0<r<r'<r_d$.
A symmetric convex set of Gaussian measure at least $1/2$ has support
function at least $\Phi^{-1}(3/4)>2/3$ in every direction, since it is
contained in the slab of that width, and so contains $(2/3)B_2$; it also
contains a point $y$ with $\|y\|_2>r'$. By convexity it contains a ball
of radius $(2/3)(1-r/r')$ around a point of the sphere of radius $r$. A
sufficiently fine finite symmetric net $S$ on that sphere therefore
meets every such set, and every law on $S$ has $\E\|X\|_2=r$, so
domination by $cG_d$ forces $c\ge r/\E\|G_d\|_2$. Since
$\E\|G_d\|_2^2=d$ and $\Var(\|G_d\|_2^2)=2d$, $r_d/\E\|G_d\|_2\to1$.
Let $r\uparrow r_d$ and then $d\to\infty$.

\end{proof}

The next consequence is the exact geometric transfer needed below. Its quantifier concerns bounded convex-dominated laws; no assertion about all unbounded laws is needed.

\begin{proposition}[Extraction of a convex subset]\label{GC-prop:extraction}
Let $B\subset\R^d$ be bounded and open. If
\begin{equation}\label{GC-eq:positive-hit}
\Pp(X\in B)>0\quad\text{for every bounded random vector }X\cx\gamma_d,
\end{equation}
then $B$ contains a bounded open convex set $C$ with $\gamma_d(C)\ge1/2$.
For symmetric $B$, it suffices to test symmetric laws in~\eqref{GC-eq:positive-hit}, and $C$ may be symmetric.
\end{proposition}
\begin{proof}
Choose $r<R$ with $\overline B\subset B(0,r)$ and $\gamma_d(B(0,R))>1/2$. Suppose that $B$ contains no convex set of the required type. Put $S=\overline B(0,R)\setminus B$. We claim that $S$ meets every bounded open convex set $C$ of Gaussian measure at least $1/2$. Such $C$ has a point $z\in B(0,R)$ because the Gaussian mass outside that ball is less than $1/2$. If $C$ avoided $S$, then $z\in B\subset B(0,r)$. Any point of $C$ outside $\overline B(0,R)$ would, by the segment from $z$, force $C$ to meet the annulus $r<|x|\le R$, which lies in $S$. Therefore $C\subset\overline B(0,R)$, and avoidance of $S$ gives $C\subset B$, a contradiction.

Theorem~\ref{GC-thm:hitting} now gives a bounded $X\cx\gamma_d$ supported on $S\subset B^c$, contradicting~\eqref{GC-eq:positive-hit}. In the symmetric case all sets and laws in this argument can be restricted to their symmetric versions.
\end{proof}

\begin{proposition}[The half-measure level of extraction is sharp]\label{GC-prop:half-sharp}
The constant $1/2$ in Proposition~\ref{GC-prop:extraction} cannot be increased to any larger dimension-independent constant.
\end{proposition}
\begin{proof}
Let $m_d=\E|G_d|$ and let $B_d=B(0,m_d+d^{-1})$. If a convex-dominated law gave $B_d$ zero mass, then $\E|X|\ge m_d+d^{-1}>\E|G_d|$, a contradiction. Hence~\eqref{GC-eq:positive-hit} holds for $B_d$.
On the other hand, $m_d=\sqrt d+O(d^{-1/2})$ and the central limit theorem for $|G_d|^2$ gives
\[
\gamma_d(B_d)\longrightarrow\tfrac12.
\]
For completeness, the first estimate follows from the ratio
$m_d=\sqrt2\,\Gamma((d+1)/2)/\Gamma(d/2)$ and Stirling's formula; the standardized threshold $((m_d+d^{-1})^2-d)/\sqrt{2d}$ tends to zero. Every subset of $B_d$ has at most this measure.
\end{proof}

This sharpness concerns the general extraction criterion. It does not claim that the input condition $\sum_i\gamma_d(A_i)>2$ in Theorem~\ref{GC-thm:three} is an optimal mass threshold.

Applied to $S=\{As:s\in\{-1,1\}^n\}$, the theorem converts any hitting
statement for the signed sums into a law dominated by the standard
Gaussian. Its proof only needs open convex tests of Gaussian measure
strictly greater than $1/2$; each contains a closed convex body of measure
at least $1/2$, so Banaszczyk's closed-body theorem applies at the same scale.
Since $X\cx G_d$ gives $\E e^{\ip uX}\le e^{\|u\|_2^2/2}$, the
law is $1$-subgaussian; \cite[Theorem~3.4]{DGLN} obtains a
$\sqrt2$-subgaussian law from the same hypothesis. Banaszczyk's theorem
gives such a law for $\|a_j\|_2\le1/5$, hence $A\sigma\cx N(0,25I_m)$
for unit columns; the cosine reference gives
$A\sigma\cx N(0,\alpha R_A^2I_m)$ with $\alpha R_A^2<6.5$, together with
the hard bound.


\subsection{The exact threshold for a fixed reference}
\label{EXT-sec:fixed}
The Gaussian median inequality gives the threshold $1/2$ in the
preceding separation argument. For a general reference, the relevant
number is the smallest mass of a convex sublevel set immediately above
its mean. This identifies the precise geometric input to extraction. The numerical
$1/e$ for log-concave laws comes from the convex mean-sublevel inequality
proved below; a mean-one exponential already has mass $1/e$ above its mean.
The new conclusion is the exact characterization of when hitting by all
dominated laws forces a convex subset, with the largest possible threshold
for each fixed reference.

Let $\nu\in\mathcal P_1(\R^d)$. A finite convex function is
\emph{coercive} if it tends to $+\infty$ as $|x|\to\infty$. Set
\begin{equation}\label{EXT-eq:kappa}
 \kappa(\nu)=
 \inf_{\substack{f\text{ finite, convex, coercive}\ f\in L^1(\nu),\ \delta>0}}
 \nu\left\{f<\int f\,d\nu+\delta\right\}.
\end{equation}
The test class is nonempty: a first moment permits an integrable
convex function growing faster than linearly, obtained by choosing
successive slopes on sufficiently distant tails.

\begin{theorem}[Exact fixed-reference extraction]\label{EXT-thm:fixed}
Assume $\kappa(\nu)>0$. If a compact set $S$ meets every bounded open
convex set of $\nu$-measure at least $\kappa(\nu)$, then it supports
a probability measure $\mu\cx\nu$. Consequently, every bounded
open $B$ satisfying
\begin{equation}\label{EXT-eq:positive-hit}
 \mu(B)>0\quad\text{for every compactly supported }\mu\cx\nu
\end{equation}
contains a bounded open convex set $K$ with
$\nu(K)\ge\kappa(\nu)$. The constant $\kappa(\nu)$ is the largest
one for which this extraction conclusion holds for all such $B$.
\end{theorem}
\begin{proof}
We first verify the finite-test inequality used in
Theorem~\ref{GC-thm:hitting}. Let $F$ be a finite convex function
integrable under $\nu$. If $\min_S F>\int F\,d\nu$, choose an
integrable nonnegative convex function $\theta$ growing faster than
linearly and then choose $\varepsilon,\delta>0$ small enough that
\[
 \min_S(F+\varepsilon\theta)
 >\int(F+\varepsilon\theta)\,d\nu+\delta.
\]
A finite convex function has an affine minorant, so
$F+\varepsilon\theta$ is coercive. Its displayed strict sublevel
set is bounded, open and convex, has measure at least $\kappa(\nu)$
by \eqref{EXT-eq:kappa}, and misses $S$, a contradiction.
Thus $\min_S F\le\int F\,d\nu$.

Apply this inequality to every nonnegative linear combination of a
finite list of convex tests. Finite-dimensional separation and
compactness of probability measures on $S$ give one measure satisfying
all the test inequalities, exactly as in the compact Gaussian proof.
It suffices to use maxima of finitely many affine functions; these
are integrable under a first moment and approximate general convex
tests from below. Affine tests of both signs also retain the mean.
Hence the resulting measure satisfies $\mu\cx\nu$.

For extraction, choose a closed ball of radius $R$ containing
$\overline B$ in its interior and having outside mass less than
$\kappa(\nu)$. Put $S=\overline B_R\setminus B$. If $S$ met all
bounded open convex sets of mass at least $\kappa(\nu)$, the first
assertion would give a compactly supported $\mu\cx\nu$ with
$\mu(B)=0$. Thus some such convex set $K$ misses $S$. By connectedness
and the separation between $\overline B$ and the exterior of the ball,
$K$ either lies in $B$ or lies outside the ball. Its mass excludes the
latter possibility, so $K\subset B$.

Finally, for any test and $\delta>0$ in \eqref{EXT-eq:kappa}, the
bounded open convex set $B=\{f<\nu f+\delta\}$ satisfies
\eqref{EXT-eq:positive-hit}: otherwise $\mu f\ge\nu f+\delta$,
contradicting $\mu\cx\nu$. Its measure can approach $\kappa(\nu)$.
No larger uniform extraction constant is therefore possible.
\end{proof}

\begin{corollary}[Gaussian and log-concave thresholds]
\label{EXT-cor:logconcave}
For the standard Gaussian, $\kappa(\gamma_d)=1/2$.
For every full-dimensional log-concave probability measure $\nu$,
\begin{equation}\label{EXT-eq:logconcave}
 \kappa(\nu)\ge e^{-1}.
\end{equation}
The value $e^{-1}$ is sharp uniformly over dimensions and log-concave
references. In particular, a compact set meeting every bounded open
convex set of \emph{log-concave reference} mass at least $1/e$ supports
a law dominated by that reference. A bounded open set positively hit
by all compactly supported dominated laws contains a bounded open
convex subset of reference mass at least $1/e$.
\end{corollary}
\begin{proof}
Lemma~\ref{GC-lem:median} gives the Gaussian lower bound. The tests
$f_\varepsilon(x)=x_1+\varepsilon|x|^2$, with
$\varepsilon\downarrow0$ and $\delta\downarrow0$, give the reverse
bound in \eqref{EXT-eq:kappa}.

For a log-concave reference and a convex test $f$, the distribution
function $F(t)=\nu\{f\le t\}$ is log-concave by the Pr\'ekopa
inequality~\cite{EXT-Prekopa} applied to the convex sublevel epigraph. If $Z=f(X)$, the
random variable $F(Z)$ stochastically dominates a uniform variable on
$[0,1]$; this statement also allows atoms. Thus
$\E\log F(Z)\ge-1$. Jensen's inequality for the concave function
$\log F$, with truncation at the endpoints when needed, gives
$\log F(\E Z)\ge-1$. The strict sublevel at $\E Z+\delta$ has mass
at least $1/e$. This proves \eqref{EXT-eq:logconcave}.

For sharpness, let $\nu_d$ be uniform measure on the Euclidean unit
ball and use $f(x)=|x|$. Then
$\E_{\nu_d}|X|=d/(d+1)$, and
\[
 \nu_d\{|X|\le\E|X|\}
   =\left(\frac d{d+1}\right)^d\longrightarrow e^{-1}.
\]
The final hitting assertion follows directly by the finite-test
argument with $1/e$: a missed sublevel set would have mass at least
that value. Extraction then uses the same complement argument.
\end{proof}
The hypothesis of log-concavity is essential to the universal $1/e$
statement. For a general fixed reference the exact threshold is
\eqref{EXT-eq:kappa}. In the next application the reference is Gaussian,
so the stronger value $1/2$ governs the convex subset contained in the
Minkowski sum.



\subsection{Three summands, with exact containment}\label{GC-sec:convexity}
Hua, Song and Tudose resolved Talagrand's convexity conjecture by representing every Gaussian-convex-dominated law as a sum of three marginally standard Gaussians~\cite{HST}. Their geometric Corollary~1.3 places a symmetric half-measure body in $4(A+A+A)$ for $\gamma_d(A)>5/6$. Exact extraction gives the original Gaussian scale and the threshold $2/3$ in the following statement.

\begin{theorem}[Threefold Gaussian convexity]\label{GC-thm:three}\label{R10-intro:three}
Let $d\ge1$ and let $A_1,A_2,A_3$ be Borel subsets of $\R^d$ such that
\begin{equation}\label{GC-eq:mass}
\gamma_d(A_1)+\gamma_d(A_2)+\gamma_d(A_3)>2.
\end{equation}
There is a convex body $K$ satisfying
\begin{equation}\label{GC-eq:three}
K\subset A_1+A_2+A_3,\qquad \gamma_d(K)>\tfrac12.
\end{equation}
If every $A_i$ is symmetric about the origin, $K$ can be chosen symmetric. In particular, $\gamma_d(A)>2/3$ implies that $A+A+A$ contains such a body. The universal number of unweighted summands cannot be reduced to two.
\end{theorem}
A convex body is compact, convex and has nonempty interior. The number of unweighted summands is optimal by Talagrand's two-summand obstruction~\cite[Proposition~2.6]{Tal95}; this does not assert optimality of the input mass threshold.

For equal input sets, the undilated half-mass conclusion is also proved by Mazhar~\cite[Theorem~1.2]{R20-Mazhar}. The formulation here permits three different inputs, gives a strict finite-dimensional output inequality, and passes to the compact mass endpoint below.

The external representation used in the proof is the following theorem.

\begin{theorem}[Hua--Song--Tudose]\label{GC-input:three-gaussian}
If $X\cx\gamma_d$, there exists a coupling of three vectors $G_1,G_2,G_3$, each with law $\gamma_d$, whose sum has the law of $X$.
\end{theorem}
This is~\cite[Theorem~1.1]{HST}. The Gaussian vectors may be dependent.

\begin{proof}[Proof of Theorem~\ref{GC-thm:three}, initially with a non-strict output bound]
First assume the $A_i$ are bounded and open, and put $B=A_1+A_2+A_3$. For any bounded $X\cx\gamma_d$, Theorem~\ref{GC-input:three-gaussian} and the union bound yield
\begin{align*}
\Pp(X\in B)
&\ge\Pp(G_1\in A_1,\ G_2\in A_2,\ G_3\in A_3)\\
&\ge\gamma_d(A_1)+\gamma_d(A_2)+\gamma_d(A_3)-2>0.
\end{align*}
Proposition~\ref{GC-prop:extraction} gives a bounded open convex $C\subset B$ with Gaussian measure at least one half.

Next suppose each $A_i$ is compact. For $\varepsilon>0$, apply the open-set assertion to the open $\varepsilon$-neighborhoods $A_i^\varepsilon$. The resulting $C_\varepsilon$ satisfy
\[
\overline{C_\varepsilon}\subseteq
A_1+A_2+A_3+3\varepsilon\overline B(0,1),
\qquad \gamma_d(\overline{C_\varepsilon})\ge\tfrac12.
\]
These compact convex sets are uniformly bounded. Along a sequence $\varepsilon\downarrow0$, Blaschke compactness gives a Hausdorff limit $K$. It is compact and convex and lies in the exact compact sum $A_1+A_2+A_3$. For every $a>0$, eventually $\overline{C_\varepsilon}\subset K+a\overline B(0,1)$, so continuity of measure from above gives $\gamma_d(K)\ge1/2$. Positive Gaussian measure implies that $K$ has nonempty interior.

For Borel $A_i$, inner regularity gives compact subsets $D_i\subset A_i$ with $\sum_i\gamma_d(D_i)>2$. Apply the compact assertion to these sets. This avoids any need to assume that the original Minkowski sum is Borel. If each $A_i$ is symmetric, choose symmetric compact $D_i$ by replacing a compact subset by its union with its reflection; use the symmetric extraction assertion and symmetric neighborhoods at every stage.
\end{proof}

\begin{lemma}[Strict Gaussian mass from a variance perturbation]\label{GC-lem:strict}
The non-strict conclusion just proved implies the strict output bound in Theorem~\ref{GC-thm:three}.
\end{lemma}
\begin{proof}
Write $\gamma_{d,\sigma}$ for $N(0,\sigma^2 I_d)$. Its density varies continuously in $L^1$ as $\sigma\to1$. Because the sum in~\eqref{GC-eq:mass} is strictly greater than two, we can choose $\sigma>1$ sufficiently close to one that
\[
\sum_i\gamma_{d,\sigma}(A_i)>2.
\]
By a linear change of variables, the preceding non-strict assertion holds for this Gaussian law. It gives a compact convex $K\subset\sum_iA_i$ with $\gamma_{d,\sigma}(K)\ge1/2$.

Necessarily $0\in\intt K$. Otherwise separation places $K$ inside a halfspace through the origin; boundedness of $K$ makes its Gaussian measure strictly smaller than the halfspace's measure $1/2$. Since $0\in\intt K$ and $\sigma>1$, the contraction $\sigma^{-1}K$ lies in $\intt K$, and $K\setminus\sigma^{-1}K$ has positive Lebesgue measure. Thus $\gamma_d(K)>\gamma_d(\sigma^{-1}K)=\gamma_{d,\sigma}(K)\ge\tfrac12$. Symmetry is retained if it was requested.
\end{proof}

\begin{corollary}[Compact inputs at the mass endpoint]\label{GC-cor:compact-endpoint}
If the $A_i$ are compact and $\sum_i\gamma_d(A_i)\ge2$, then $A_1+A_2+A_3$ contains a compact convex set $K$ with $\gamma_d(K)\ge1/2$.
\end{corollary}
\begin{proof}
The three compact sets necessarily have positive measure. Every positive open neighborhood of a nonempty compact set has strictly greater Gaussian measure than the set itself. Thus the neighborhood inputs satisfy the strict mass condition. The uniformly bounded Hausdorff-limit argument in the proof of Theorem~\ref{GC-thm:three} applies without loss at the endpoint.
\end{proof}

\begin{corollary}[Anisotropy and finite certificates]\label{GC-cor:polytope}
Theorem~\ref{GC-thm:three} holds for every centered nondegenerate Gaussian measure, with that same measure used for the input and output. The body can be chosen to be a polytope; for standard Gaussian measure its vertices may be rational.
\end{corollary}
\begin{proof}
Apply an invertible linear map to the standard statement. A compact convex body of measure strictly greater than one half has inscribed polytopes whose measures approach its measure. In standard coordinates choose their vertices from the dense rational points in the interior. Symmetric rational polytopes are available in the symmetric case.
\end{proof}

The existence of a finite certificate is not an algorithm for locating it from an arbitrary Borel set. Talagrand's request for a constructive procedure remains a separate problem~\cite[Problem~2.3]{Tal26}.

\begin{corollary}[Symmetric output from a larger nonsymmetric input]\label{GC-cor:symmetrize}
If $\gamma_d(A)>5/6$, there is a symmetric convex body $K\subset A+A+A$ with $\gamma_d(K)>1/2$.
\end{corollary}
\begin{proof}
The symmetric Borel set $A\cap(-A)$ has measure at least $2\gamma_d(A)-1>2/3$. Apply Theorem~\ref{GC-thm:three} to it.
\end{proof}

\begin{corollary}[Solid sets]\label{GC-cor:solid}
Suppose the $A_i$ satisfy~\eqref{GC-eq:mass} and are solid, meaning that $x\in A_i$ and $|y_j|\le|x_j|$ for every $j$ imply $y\in A_i$. Then $A_1+A_2+A_3$ contains a solid convex body of Gaussian measure greater than one half.
\end{corollary}
\begin{proof}
Minkowski sums of solid sets are solid: coordinatewise reduce a representation of a sum, allocating each coordinate's required magnitude among the available summands. Start with the convex body $K$ from Theorem~\ref{GC-thm:three}. Choose an orthant $O$ for which $\gamma_d(K\cap O)$ is maximal. Its solidification
\[
C=\{y:\text{there exists }x\in K\cap O\text{ with }|y_j|\le|x_j|\ \forall j\}
\]
is compact and convex, because $x\mapsto(|x_j|)_j$ is linear on the fixed orthant. It is contained in the solid sum. Reflections of $K\cap O$ have pairwise null intersections, so
$\gamma_d(C)\ge2^d\gamma_d(K\cap O)\ge\gamma_d(K)>1/2$.
This solidification argument is the one used in~\cite[Appendix~A]{HST}.
\end{proof}

The existing two-summand counterexamples~\cite{Tal95,HST} prove that three is optimal here. The discrete $p$-smallness conjectures in those works additionally require combinatorial coverings; that further condition remains outside this Gaussian-containment theorem.



\subsection{Gaussian measures on separable Banach spaces}\label{GC-sec:gaussian-spaces}
The finite-dimensional geometric result extends to path spaces without any enlargement of the Minkowski sum. The output is a compact convex set. In infinite dimension such a set need not have nonempty interior.

\begin{lemma}[Finite-dimensional Gaussian approximation]\label{GC-lem:finite-rank}
Let $\gamma$ be a centered Radon Gaussian measure on a separable Banach space $E$. There are continuous finite-rank linear maps $P_n:E\to E$ such that
\[
P_nG\longrightarrow G\quad\text{almost surely and in }L^1(E),\qquad G\sim\gamma.
\]
For each $\delta>0$, there is a Borel set $D$ with $\gamma(D)>1-\delta$ on which $P_nx\to x$ uniformly.
\end{lemma}
\begin{proof}
Choose countably many continuous linear functionals $\ell_j$ generating the Borel sigma-field of $E$. A Radon Gaussian random element is Bochner integrable. Gaussian regression gives $P_nG=\E[G\mid\ell_1(G),\ldots,\ell_n(G)]$, where $P_n$ is a continuous finite-rank linear map; a singular covariance matrix is handled by its pseudoinverse. The closed-martingale convergence theorem gives almost sure and $L^1$ convergence to $G$. For this particular closed martingale, the Banach-valued conclusion follows by approximation of $G$ by simple random variables and Doob's scalar maximal inequality applied to the conditional expectation of the norm of the error. Egorov's theorem gives uniform convergence off a set of probability less than $\delta$.
\end{proof}

\begin{theorem}[Threefold convexity in a Gaussian Banach space]\label{GC-thm:banach}
Let $\gamma$ be a centered Radon Gaussian measure on a separable Banach space $E$. If Borel sets $A_1,A_2,A_3\subset E$ satisfy $\sum_i\gamma(A_i)>2$, there is a compact convex set
\[
K\subset A_1+A_2+A_3,\qquad\gamma(K)\ge\tfrac12.
\]
The same conclusion holds with the sum replaced by $A_0+(A_1+A_2)/2$ under the corresponding mass assumption. Symmetric inputs admit a symmetric output.
\end{theorem}
\begin{proof}
Let the mass excess be $a=\sum_i\gamma(A_i)-2>0$. Use Lemma~\ref{GC-lem:finite-rank} with $\delta<a/6$, and then inner regularity, to choose compact $D_i\subset A_i$ with $\sum_i\gamma(D_i)>2$ and with $P_n\to I$ uniformly on every $D_i$. We may intersect all these sets with the closed linear support of $\gamma$. For symmetric inputs choose a symmetric uniform-convergence set and symmetrize the compact subsets.

Put $\gamma_n=(P_n)_\#\gamma$. In its finite-dimensional linear support, this is a nondegenerate Gaussian measure (or a point mass, for which the conclusion is immediate). Since $\gamma_n(P_nD_i)\ge\gamma(D_i)$, Theorem~\ref{GC-thm:three}, after a linear change of variables, gives compact convex sets
\[
K_n\subset P_nD_1+P_nD_2+P_nD_3,\qquad\gamma_n(K_n)\ge\tfrac12.
\]
Uniform convergence on the $D_i$ implies Hausdorff convergence of their sums to the compact set $D_1+D_2+D_3$. Hence the union of the $K_n$ is totally bounded, and a subsequence converges in Hausdorff distance to a nonempty compact set $K$. Convexity passes to the limit, and $K\subset D_1+D_2+D_3$.

The measures $\gamma_n$ converge weakly to $\gamma$. For $\varepsilon>0$, let $K^{[\varepsilon]}=\{x:\dist(x,K)\le\varepsilon\}$, a closed set. Eventually $K_n\subset K^{[\varepsilon]}$, and Portmanteau gives
\[
\tfrac12\le\limsup_n\gamma_n(K^{[\varepsilon]})\le\gamma(K^{[\varepsilon]}).
\]
Let $\varepsilon\downarrow0$ to get $\gamma(K)\ge1/2$. The weighted assertion uses Theorem~\ref{GC-thm:weighted} and the same limit. Symmetry passes through all the compact approximations.
\end{proof}

\begin{corollary}[Wiener space]\label{GC-cor:wiener}
Let $\mathsf W$ be Wiener measure on $C_0([0,T];\R^d)$. If a Borel family of paths $A$ has $\mathsf W(A)>2/3$, then the set of pathwise sums of three members of $A$ contains a compact convex family $K$ with $\mathsf W(K)\ge1/2$.
The same is true for $A+(A+A)/2$. Three different path families are permitted when their probabilities sum to more than two.
\end{corollary}
\begin{proof}
Wiener measure is a centered Radon Gaussian measure on the indicated separable Banach space. Apply Theorem~\ref{GC-thm:banach}.
\end{proof}

The strict finite-dimensional output inequality was obtained by total-variation continuity under a small change of Gaussian variance. That argument does not extend to general infinite-dimensional Gaussian measures, which may be singular under variance changes. The conclusion above is therefore stated with $\gamma(K)\ge1/2$.

\begin{corollary}[Green's Gaussian product-space problem]\label{R20-cor:green}
Let $\gamma_\infty$ be the product standard Gaussian measure on $\R^{\N}$ with its product topology. If $A\subset\R^{\N}$ is Borel and $\gamma_\infty(A)>2/3$, its threefold coordinatewise sum contains a product-compact convex set $K$ with $\gamma_\infty(K)\ge1/2$. In particular this answers Problem~54 in Green's list~\cite{R20-Green}, which asks for ten summands and output mass at least $0.01$ from balanced compact input of mass at least $0.99$.
\end{corollary}
\begin{proof}
Choose $w_j>0$ with $\sum_jw_j<\infty$. The weighted Hilbert space
$H_w=\{x:\sum_jw_jx_j^2<\infty\}$ has full $\gamma_\infty$ measure and carries the centered Radon Gaussian with covariance eigenvalues $w_j$ in its orthonormal coordinates. The inclusion $H_w\to\R^{\N}$ is continuous, so $A\cap H_w$ is Borel in $H_w$. Theorem~\ref{GC-thm:banach} gives a norm-compact convex $K\subset3(A\cap H_w)$ of measure at least one half. Its continuous image is product-compact and lies in $A+A+A$. For balanced input, $0\in A$, so this threefold sum is contained in the tenfold sum requested in the problem.
\end{proof}


\section{Continuous martingales with bounded covariance}
\label{R11-sec:bounded}
The joint rounding law couples its output to a Gaussian reference by conditional means. To interpolate while controlling covariance after stopping times, first normalize terminal curvature by an affine map. Canonical contraction gives Lipschitz conditional endpoint maps; heat interpolation turns their Lipschitz bound into a covariance cap. Brownian completion gives the exact Gaussian residual.

The argument covers every law dominated by the uniformly log-concave terminal law, including discrete initial laws. Conditional log-concavity gives a separate finite strict-slack construction. The next section proves the canonical contraction used here and its sharp stability and future-path consequences.


\subsection{A continuous martingale with bounded instantaneous covariance}\label{GC-sec:martingales}
Terminal curvature controls how much covariance a martingale can accumulate in each unit of time. The canonical contraction theorem proved in Section~\ref{R11-sec:canonical} gives Lipschitz conditional maps for every initial law. Their heat interpolation proves the realization below. We use that single regularity proof for existence and for the later sharp variational comparison; the Brownian completion here then gives the weighted geometric consequences.

\begin{theorem}[A bounded-volatility realization]\label{GC-thm:volatility}\label{GC-prop:finite-mart}
Let $Q$ be a positive-definite $d\times d$ matrix and let
\[
\nu(dy)=Z^{-1}e^{-V(y)}\,dy,\qquad 0<Z<\infty,
\qquad y\longmapsto V(y)-\tfrac12\ip{Q^{-1}y}{y}\ \text{is convex}.
\]
The potential is lower semicontinuous and may take the value $+\infty$; this includes convex support constraints. For every integrable law $\mu\cx\nu$, there exists a continuous square-integrable martingale $(M_t)_{0\le t\le1}$ on a filtered probability space such that
\begin{equation}\label{GC-eq:volatility}
\law(M_0)=\mu,\qquad \law(M_1)=\nu,\qquad
0\preceq\frac{d\langle M\rangle_t}{dt}\preceq Q
\quad dt\otimes d\Pp\text{-almost everywhere}.
\end{equation}
Conversely, existence of a martingale with these two endpoint laws implies $\mu\cx\nu$. No irreducibility or absolute-continuity hypothesis on $\mu$ is required.
\end{theorem}

\begin{corollary}[Gaussian convex order as a continuous-time condition]\label{GC-cor:gauss-dynamic}\label{R10-eq:dynamic-Gaussian}
For an integrable law $\mu$ on $\R^d$,
\[
\mu\cx\gamma_d
\quad\Longleftrightarrow\quad
\exists M:\ \law(M_0)=\mu,\ \law(M_1)=\gamma_d,
\quad 0\preceq d\langle M\rangle_t/dt\preceq I_d.
\]
The bound $I_d$ cannot be replaced uniformly by $cI_d$ with $c<1$.
\end{corollary}
\begin{proof}[Proof of the Gaussian specialization]
Take $V(y)=|y|^2/2$ and $Q=I_d$ in Theorem~\ref{GC-thm:volatility}. For sharpness, a martingale from zero to a standard Gaussian has $\E\langle M\rangle_1=I_d$.
\end{proof}

\begin{lemma}[Conditional densities with retained curvature]\label{GC-input:conditional}
Let $\nu$ have a density proportional to $e^{-V}$, where $V(y)-|y|^2/2$ is lower semicontinuous and convex, possibly extended-valued, and let its mean be $m$. Let $\mu$ be finitely supported and suppose that for some $0<\varepsilon<1$,
\[
\mu\cx\law\bigl(m+(1-\varepsilon)(Y-m)\bigr),\qquad Y\sim\nu.
\]
There is a martingale coupling $(X,Y)$ of $\mu,\nu$ such that every conditional law $\law(Y\mid X=x_i)$ is $1$-uniformly log-concave.
\end{lemma}
\begin{proof}
Translate by $m$ and apply~\cite[Proposition~3.3]{HST}. Writing $p_i=\Pp(X=x_i)$, it gives conditional density ratios
\[
f_i(y)=\frac{\exp(a_i+\ip{b_i}{y})}
{\sum_jp_j\exp(a_j+\ip{b_j}{y})},\qquad
\int f_i\,d\nu=1,\quad\int y f_i(y)\,\nu(dy)=x_i,
\quad\sum_ip_if_i=1.
\]
Each $\log f_i$ is concave. The conditional potential is $V-\log f_i$ up to a constant, so its difference from $|y|^2/2$ is convex, including at the boundary of its support. The existence of the displayed ratios is exactly the cited external proposition.
\end{proof}

\begin{lemma}[Heat martingale of a Lipschitz transport]\label{GC-lem:heat}
Let $F:\R^d\to\R^d$ be $1$-Lipschitz and let $B$ be standard Brownian motion. Put
\[
H_t=(P_{1-t}F)(B_t)\quad(0\le t<1),\qquad H_1=F(B_1),
\]
where $P_sF(x)=\E F(x+\sqrt s\,G)$. Then $H$ is a continuous square-integrable martingale, $H_0=\E F(G)$, and
\[
\frac{d\langle H\rangle_t}{dt}
=\nabla P_{1-t}F(B_t)\nabla P_{1-t}F(B_t)^{\mathsf T}
\preceq I_d.
\]
\end{lemma}
\begin{proof}
The Markov property of Brownian motion gives $H_t=\E[F(B_1)\mid\cF_t^B]$. Heat smoothing preserves the Lipschitz constant, so the gradient has operator norm at most one. It\^o's formula gives the bracket identity on $[0,1)$. Continuity at one follows from
\[
|P_{1-t}F(B_t)-F(B_1)|
\le |B_t-B_1|+\sqrt{1-t}\,\E|G|.
\]
The Lipschitz growth of $F$ ensures square integrability. This heat representation also underlies Song's Gaussian Lipschitz-image lemma~\cite[Lemma~2.4]{Song}.
\end{proof}

\begin{proof}[Proof of Theorem~\ref{GC-thm:volatility}]
An affine minorant of the convex part of $V$ gives Gaussian tails. Thus both endpoint laws have second moments, by convex order. Whiten by $Q^{-1/2}$ and apply Theorem~\ref{CM-thm:contraction} with $L=1$. Its proof uses finite-source Bass representation, contraction of Gaussian-mixture transports and the fixed-terminal kernel limit; it does not use the present realization theorem. On $X\sim\mu$ and an independent Brownian motion $B$, the canonical conditional maps $T_x$ therefore give
\[
 M_t=Q^{1/2}(P_{1-t}T_{Q^{-1/2}X})(B_t).
\]
Lemma~\ref{GC-lem:heat} proves continuity and the bracket bound. The initial value is $X$, and the terminal law is $\nu$. Conversely, conditional Jensen for $M_0=\E[M_1\mid\cF_0]$ gives $\mu\cx\nu$. Under the finite strict-slack hypotheses, Lemma~\ref{GC-input:conditional} gives a second construction by applying Caffarelli contraction to its conditional densities and the same heat lemma.
\end{proof}
The realization can consequently be chosen canonical in the whitened coordinates. The bracket and conditional-reference conclusions below use only its covariance cap; they apply equally to any continuous martingale satisfying that cap.

\subsubsection{Brownian completion and the weighted theorem}\label{GC-sec:completion}
The matrix bracket bound gives a completion by two Brownian motions on the whole time interval. Their individual independence from the initial sigma-field survives, even though they can be mutually dependent.

\begin{lemma}[Two Brownian completions]\label{GC-lem:brown-complete}
Let $M$ be a continuous square-integrable martingale with
$0\preceq d\langle M\rangle_t/dt\preceq I_d$. On an enlargement of its filtered space, there are standard $d$-dimensional Brownian motions $U,V$, both Brownian relative to the enlarged filtration, such that
\begin{equation}\label{GC-eq:brown-completion}
M_t=M_0+\tfrac12(U_t+V_t)\qquad(0\le t\le1).
\end{equation}
Each of the entire paths $U,V$ is independent of $\cF_0$.
\end{lemma}
\begin{proof}
Choose a predictable version $a_t$ of the bracket density and add an independent standard Brownian motion $W$. Put $N_t=M_t-M_0$ and
\[
R_t=\int_0^t(I_d-a_s)^{1/2}\,dW_s,
\qquad U_t=N_t+R_t,\qquad V_t=N_t-R_t.
\]
The cross variation of $N$ with the added Brownian integral is zero. Therefore $\langle U\rangle_t=\langle V\rangle_t=tI_d$. L\'evy's characterization makes each a Brownian motion relative to the enlarged filtration, hence independent of its initial sigma-field. Identity~\eqref{GC-eq:brown-completion} is exact. The completion is the continuous-time form of the two-Gaussian argument in~\cite{Song}.
\end{proof}

\begin{theorem}[Conditional Gaussian residuals at stopping times]\label{GC-thm:residual}\label{R10-eq:stopped-reference}
The martingale in Theorem~\ref{GC-thm:volatility} may be realized so that, for every stopping time $0\le\tau\le1$,
\begin{equation}\label{GC-eq:residual}
\law(M_1-M_\tau\mid\cF_\tau)
\cx N\bigl(0,(1-\tau)Q\bigr)\quad\text{almost surely}.
\end{equation}
Here the Gaussian covariance on the right is a function of the conditioned value of $\tau$. For deterministic $s\le t$, the same conclusion holds for $M_t-M_s$ with covariance $(t-s)Q$. In particular,
\[
W_p\bigl(\law(M_s),\law(M_t)\bigr)
\le \sqrt{t-s}\,\bigl(\E|Q^{1/2}G|^p\bigr)^{1/p}\qquad(p\ge1).
\]
\end{theorem}
\begin{proof}
Whiten by $Q$ and use Lemma~\ref{GC-lem:brown-complete}. Given $\cF_\tau$, each Brownian increment $U_1-U_\tau$ and $V_1-V_\tau$ has law $N(0,(1-\tau)I_d)$ by the strong Markov property. Conditional Jensen at their midpoint proves~\eqref{GC-eq:residual}. This argument needs their two marginal conditional laws; it imposes no joint independence. The deterministic increment assertion follows in the same way. Apply it to the convex function $z\mapsto|z|^p$ and use the existing coupling $(M_s,M_t)$ for the Wasserstein bound.
\end{proof}

\begin{corollary}[A conditional strengthening of Gaussian convex order]\label{GC-cor:conditional-gaussian}
For every $X\cx\gamma_d$, there is a coupling $(X,G_0,G_1,G_2)$ such that
\begin{equation}\label{GC-eq:weighted-rep}
G_i\sim\gamma_d\ (i=0,1,2),\qquad
\law(G_j\mid X)=\gamma_d\ (j=1,2),\qquad
X=G_0+\tfrac12(G_1+G_2).
\end{equation}
Consequently,
\[
\E[G_0\mid X]=X,\qquad
\law(G_0-X\mid X)\cx\gamma_d.
\]
The bound on the conditional residual has optimal universal scale one.
\end{corollary}
\begin{proof}
Take the Gaussian-terminal martingale from Corollary~\ref{GC-cor:gauss-dynamic}, set $X=M_0$, $G_0=M_1$, and use $G_1=-U_1$, $G_2=-V_1$ in Lemma~\ref{GC-lem:brown-complete}. The paths $U,V$ are individually independent of $X$, so the conditional marginals are standard Gaussian. Conditional Jensen gives the residual comparison. When $X=0$, the residual is exactly standard Gaussian, proving the optimal scale claim.
\end{proof}

The weighted representation below is also a consequence of the finite strict-slack inputs in Hua--Song--Tudose~\cite[Proposition~3.3 and Lemma~2.3]{HST}. Their conditional uniform log-concavity gives a contracting Gaussian transport by Caffarelli's theorem; the centered image is a half-sum of two standard Gaussians. Conditioning on the initial state and then using finite conditional-mean approximation gives the general weighted representation by weak compactness. The canonical construction above selects the variational optimizer and provides its covariance cap and stopping-time residual control.

\begin{theorem}[Weighted three-set Gaussian convexity]\label{GC-thm:weighted}
Let $A_0,A_1,A_2$ be Borel subsets of $\R^d$ whose Gaussian measures sum to more than two. Then the set
\begin{equation}\label{GC-eq:weighted-set}
A_0+\tfrac12A_1+\tfrac12A_2
\end{equation}
contains a convex body $K$ with $\gamma_d(K)>1/2$. Symmetric inputs admit symmetric $K$.
\end{theorem}
\begin{proof}
For bounded open inputs, every $X\cx\gamma_d$ has representation~\eqref{GC-eq:weighted-rep}. The event $\{G_i\in A_i\text{ for }i=0,1,2\}$ has probability at least $\sum_i\gamma_d(A_i)-2>0$, and on it $X$ belongs to~\eqref{GC-eq:weighted-set}. Apply Proposition~\ref{GC-prop:extraction}. Compact neighborhoods, inner regularity and the variance perturbation in Lemma~\ref{GC-lem:strict} give the exact Borel and strict-mass conclusion. The sum of the coefficients is two; the set in~\eqref{GC-eq:weighted-set} is generally different from $A+A$ when all inputs coincide.
\end{proof}

\begin{corollary}[Four summands and higher output mass]\label{R15-cor:fourfold}
If $A\subset\R^d$ is symmetric Borel and $\gamma_d(A)>2/3$, then
$A+A+A+A$ contains a symmetric compact convex body $C$ with
\[
 \gamma_d(C)>\beta_4:=2\Phi(2\Phi^{-1}(3/4))-1
                         >0.822656449.
\]
For compact $A$ the input and output inequalities may both be non-strict.
More generally, $4j$ summands give mass greater than
$2\Phi(2j\Phi^{-1}(3/4))-1$. Hence
$O(\sqrt{\log(1/\delta)})$ summands suffice for output mass $1-\delta$,
and this order is necessary uniformly for fixed input mass below one.
\end{corollary}
\begin{proof}
Theorem~\ref{GC-thm:weighted} gives symmetric convex
$K\subset A+\tfrac12A+\tfrac12A$ of mass greater than $1/2$.
The Gaussian $S$-inequality~\cite{R15-LO} gives
$\gamma_d(2jK)\ge2\Phi(2j\Phi^{-1}((1+\gamma_d(K))/2))-1$.
Since $2jK\subset2jA+jA+jA\subset A+\cdots+A$ with $4j$ summands,
this proves the strict bounds. At the compact endpoint, apply the weighted
theorem to shrinking neighborhoods and take a bounded Hausdorff limit
before dilating, as in Corollary~\ref{GC-cor:compact-endpoint}.
A one-dimensional interval of fixed mass in $(2/3,1)$ proves the asserted
order by the Gaussian tail quantile.
\end{proof}
This answers the higher-mass question in~\cite[Problem~2.1]{Tal26} with four summands. The optimality theorem for three concerns the half-mass output in Theorem~\ref{GC-thm:three}.

\subsubsection{Three-point convex hulls and explicit summand counts}
\begin{corollary}[Three-point convex combinations with dilation two]
\label{R19-cor:conv-three}
For a Borel set $A\subset\R^d$ with $\gamma_d(A)>2/3$, the set $2\operatorname{conv}_3(A)$ contains a convex body of Gaussian measure greater than $1/2$. A symmetric input admits a symmetric output. If $A$ is compact and $\gamma_d(A)\ge2/3$, the same conclusion holds with a compact convex output of mass at least $1/2$.
\end{corollary}
\begin{proof}
The weighted theorem gives $K\subset A+\tfrac12A+\tfrac12A$. For $a,b,c\in A$,
\[
 \frac12\left(a+\frac12b+\frac12c\right)
 =\frac12a+\frac14b+\frac14c\in\operatorname{conv}_3(A).
\]
Theorem~\ref{GC-thm:weighted} proves the strict and symmetric assertions. For compact input at mass $2/3$, apply it to shrinking neighborhoods of $A$. The resulting convex bodies lie in a fixed bounded set; a Hausdorff-convergent subsequence gives a compact convex set inside $A+\tfrac12A+\tfrac12A$, of Gaussian mass at least $1/2$, as in Corollary~\ref{GC-cor:compact-endpoint}.
\end{proof}
Johnston's examples exclude a uniform undilated finite-point convex-hull conclusion, even for much larger input mass~\cite{R19-Johnston}. The displayed weights give a sufficient dilation of two for three points. The sharp undilated unweighted sum theorem also implies Johnston's Conjecture~1.1 with $k=3$, $\varepsilon=1/4$: balanced input of mass at least $3/4$ already gives a convex core of mass $>1/2$ in $A+A+A$, hence in the larger dilation allowed there.

\begin{corollary}[Explicit counts above Gaussian mass one half]
\label{R19-cor:explicit-counts}
Let $A\subset\R^d$ be Borel with $\gamma_d(A)\ge\alpha$, where $1/2<\alpha<1$. Put
\begin{equation}\label{R19-eq:count-choice}
 k_\alpha=\left\lfloor
       \frac{\Phi^{-1}(2/3)}{\Phi^{-1}(\alpha)}\right\rfloor+1.
\end{equation}
The sum of $3k_\alpha$ copies of $A$ contains a convex body of Gaussian mass $>1/2$. If $A$ is symmetric, $4k_\alpha$ copies contain a symmetric convex body of mass $>\beta_4$, with $\beta_4$ as in \eqref{R19-eq:intro-fourfold}. More generally, $4jk_\alpha$ copies contain such a body of mass
\[
 >2\Phi\!\left(2j\Phi^{-1}(3/4)\right)-1\qquad(j\ge1).
\]
For $\alpha=0.60,0.55,0.51$, the symmetric counts are respectively $8,16,72$. For $\alpha>2/3$ the count is four. At $\alpha=2/3$, compact input also permits four with the non-strict endpoint output in Corollary~\ref{R15-cor:fourfold}.
\end{corollary}
\begin{proof}
Choose a compact $D\subset A$ whose mass $\beta$ is sufficiently close to $\gamma_d(A)$ that
$k_\alpha\Phi^{-1}(\beta)>\Phi^{-1}(2/3)$. For symmetric $A$ choose $D$ symmetric. Borell's general Gaussian Brunn--Minkowski inequality~\cite[Theorem~2.3]{R20-Borell} permits all coefficients to equal one: their sum is at least one and at least twice the largest coefficient minus one. Hence, for $k_\alpha\ge2$,
\[
 \gamma_d(D^{+k_\alpha})
 \ge\Phi\!\left(k_\alpha\Phi^{-1}(\beta)\right)>2/3.
\]
For $k_\alpha=1$ the assertion is the identity. The sum is compact, so Theorem~\ref{GC-thm:three} applies and gives $3k_\alpha$ original summands. For symmetric input, Corollary~\ref{R15-cor:fourfold} and its higher-output iterates give the remaining conclusions. Substitution in \eqref{R19-eq:count-choice} gives the displayed counts. The mass-$2/3$ compact endpoint uses the separate fourfold endpoint theorem.
\end{proof}
The amplification depends only on the input threshold and the desired output mass. A direct Gaussian Steinhaus estimate, useful without the full Brunn--Minkowski inequality, is
\begin{equation}\label{R19-eq:translation}
 \gamma_d(B-x)\ge\Phi\!\left(\Phi^{-1}(\gamma_d(B))-|x|\right).
\end{equation}
Indeed, among sets of fixed Gaussian mass, the half-space opposite $x$ minimizes the translated mass by the monotone likelihood ratio in that direction. If $\gamma_d(B)>1/2$, this gives $B(0,2\Phi^{-1}(\gamma_d(B)))\subset B+B$ by intersecting $B$ with $x-B$. The full inequality in the proof permits every integer number of summands.

\begin{proposition}[Two Gaussian summands are insufficient]\label{GC-prop:two-gaussian}
There exists a centered two-point law dominated by $\gamma_1$ which cannot be represented as the sum of two centered Gaussian variables, even if their variances are allowed to differ.
\end{proposition}
\begin{proof}
For $G\sim N(0,1)$, set $a=\E|G|=\sqrt{2/\pi}$. The variable $X=a\,\operatorname{sign}(G)=\E[G\mid\operatorname{sign}(G)]$ satisfies $X\cx G$. If $Y+Z\in\{-a,a\}$ almost surely and $t=\pi/a$, then $e^{itY}=-e^{-itZ}$ almost surely. Taking expectations contradicts
\[
\E e^{itY}=e^{-t^2\operatorname{Var}(Y)/2}>0,
\qquad \E e^{-itZ}=e^{-t^2\operatorname{Var}(Z)/2}>0.
\]
No independence of $Y,Z$ was used. This is a scalar instance of the known two-Gaussian obstruction discussed in~\cite[Section~1]{HST}.
\end{proof}

\begin{corollary}[Bounded-time Skorokhod embedding]\label{GC-cor:skorokhod}
In dimension one, under the hypotheses of Theorem~\ref{GC-thm:volatility} with $Q=q>0$, there is a Brownian motion $B$, independent of $X\sim\mu$ at time zero, and a stopping time $\tau\le q$ in a filtration allowing the initial variable and auxiliary randomness, such that
\[
X+B_\tau\sim\nu.
\]
In particular, every $\mu\cx N(0,1)$ has such an embedding with $\tau\le1$.
\end{corollary}
\begin{proof}
Apply the Dambis--Dubins--Schwarz time change to $M-M_0$, extending the probability space beyond the terminal quadratic variation if necessary. The stopping time is $\tau=\langle M\rangle_1\le q$. The time-changed Brownian motion is Brownian relative to a filtration containing $\sigma(M_0)$ initially, so it is independent of $M_0$. Standard time-change theory is recorded in~\cite{KS}.
\end{proof}

The embedding starts with $X\sim\mu$ and ends at the prescribed uniformly log-concave law. Ankirchner--Strack and Ankirchner--Hobson--Strack study bounded-time embeddings from a point and give Gaussian-quantile sufficient conditions~\cite{R20-AS,R20-AHS}. Here Gaussian domination is an exact criterion for evolving a prescribed random initial law to a Gaussian in bounded time. The filtration may contain auxiliary randomness. In higher dimensions Theorem~\ref{GC-thm:volatility} gives a matrix covariance cap; it does not replace that cap by a single scalar stopping time.

\begin{corollary}[Convex support retained for the whole path]\label{GC-cor:path-support}
If the terminal reference $\nu$ in Theorem~\ref{GC-thm:volatility} is supported on a closed convex set $C$, the martingale can be chosen with
\[
\Pp(M_t\in C\text{ for every }0\le t\le1)=1.
\]
For the product cosine reference
\[
\nu(dy)=\prod_{j=1}^d\left[\frac1{R_j}
\cos^2\!\left(\frac{\pi y_j}{2R_j}\right)
\mathbf1_{\{|y_j|<R_j\}}\,dy_j\right],
\]
every $\mu\cx\nu$ therefore has a martingale realization confined to $\prod_j[-R_j,R_j]$, with
\[
\frac{d\langle M\rangle_t}{dt}
\preceq\operatorname{diag}\left(\frac{2R_j^2}{\pi^2}\right)_{j=1}^d.
\]
\end{corollary}
\begin{proof}
For every deterministic $t$, $M_t=\E[M_1\mid\cF_t]$ lies in $C$ by conditional Jensen applied to the convex distance function. Apply this at rational times and use path continuity and closedness of $C$. For the cosine density, the second derivative of its $j$th coordinate potential is
\[
\frac{\pi^2}{2R_j^2}\sec^2\!\left(\frac{\pi y_j}{2R_j}\right)
\ge\frac{\pi^2}{2R_j^2}.
\]
Extending the potential by $+\infty$ outside the box gives the uniformly log-concave reference required by Theorem~\ref{GC-thm:volatility}.
\end{proof}

\begin{corollary}[Application to the same signing law]\label{GC-cor:sign-law}
Suppose a signing law satisfies
\[
S=(A\sigma,\sigma)\cx N(0,Q),\qquad Q\succ0.
\]
The entire law of $\sigma$ has a coupling with a Gaussian $R\sim N(0,Q)$ and a continuous martingale starting at $S$ and ending at $R$ such that
\[
\E[R\mid\sigma]=S,\qquad
\law(R-S\mid\sigma)\cx N(0,Q),\qquad
0\preceq d\langle M\rangle_t/dt\preceq Q.
\]
Every signing probability, the hard discrepancy support, and the entropy of the signing law are unchanged. If $Q$ is block diagonal, the Gaussian reference blocks retain their marginal independence.
\end{corollary}
\begin{proof}
Apply Theorem~\ref{GC-thm:volatility} and Theorem~\ref{GC-thm:residual} to the existing law of $S$. Its second block determines $\sigma$, so conditioning on $S$ is equivalent to conditioning on $\sigma$. All target-law properties are therefore retained. Independence of Gaussian blocks follows from the block-diagonal covariance of the unchanged Gaussian marginal.
\end{proof}
The source--sign joint law may change. Its mutual information is not fixed by preserving the signing law and the Gaussian marginal.


\subsection{Five summands below half Gaussian mass}
\label{R22-sec:balanced}
Balancedness lets the weighted theorem reach input masses below one half. Here a set $A$ is \emph{balanced} if $tA\subset A$ for every real $t$ with $|t|\le1$, and $A^{+k}$ denotes its actual $k$-fold Minkowski sum. The radial part of this hypothesis has a specific use: a long chord in each direction fills a ball in $A-A$. Gaussian enlargement then increases the input mass, while the two coefficients $1/2$ in Theorem~\ref{GC-thm:weighted} keep the enlarged sets inside five copies of the original set.

\begin{theorem}[Fivefold Gaussian convexity below half mass]
\label{R22-thm:fivefold}
In every finite dimension, a balanced Borel set $A\subset\R^d$ with $\gamma_d(A)\ge5/12$ admits a compact symmetric convex set
\[
 K\subset A^{+5},\qquad \gamma_d(K)>\tfrac12.
\]
For a centered Radon Gaussian measure $\gamma$ on a separable Banach space $E$, the same input hypothesis gives compact symmetric convex $K\subset A^{+5}$ with $\gamma(K)\ge1/2$.
\end{theorem}

\begin{lemma}[Gaussian chords in a star-shaped set]
\label{R22-lem:chords}
Let $A$ be compact and star-shaped about zero, with Gaussian measure at least $\alpha\in(0,1)$. In $\R^d$,
\begin{equation}\label{R22-eq:chords}
 \rho_\alpha B_2^d\subset A-A,
 \qquad \rho_\alpha=2\Phi^{-1}((1+\alpha)/2).
\end{equation}
In a separable Gaussian Banach space the conclusion is $\rho_\alpha B_H\subset A-A$, where $B_H$ is the closed unit ball of the Cameron--Martin space. The radius is sharp for a centered interval in dimension one. For balanced $A$, the difference set is $A+A$.
\end{lemma}
\begin{proof}
Fix a unit direction $e$. Compactness gives a longest chord of $A$ parallel to $e$, of length $D_e$. Every parallel section lies in an interval of length at most $D_e$, whose one-dimensional Gaussian measure is at most $2\Phi(D_e/2)-1$. Gaussian Fubini gives $\alpha\le2\Phi(D_e/2)-1$, hence $D_e\ge\rho_\alpha$. The difference set is star-shaped, since $t(a-b)=ta-tb$ for $0\le t\le1$. It therefore contains every shorter vector in the direction $e$, proving \eqref{R22-eq:chords}.

For a unit $h\in H$, write a Gaussian random element as $\xi h+G'$, with $\xi$ standard normal and independent of the Gaussian residual $G'$. This is the orthogonal decomposition associated with the Gaussian linear functional of $h$. The same conditional section argument applies and proves the assertion for $B_H$. Its image is closed in $E$: a strongly $E$-convergent sequence in a bounded $H$-ball has a weakly $H$-convergent subsequence, and the continuous embedding identifies its limit. By the proved inclusion, that closed image lies in compact $A-A$, so the balls used below are compact in $E$. The centered interval $[-a,a]$ has measure $2\Phi(a)-1$ and difference interval $[-2a,2a]$, giving equality in the radius.
\end{proof}

\begin{proof}[Proof of Theorem~\ref{R22-thm:fivefold}]
Corollary~\ref{R23-cor:five-mass} proves the stronger uniform mass
$\Phi(1/16)>1/2$, including the Banach conclusion. Its proof below uses
only the chord lemma and weighted Gaussian extraction.
\end{proof}

The half-mass existence problem in Talagrand's prize formulation~\cite{R22-TalagrandPrize} was resolved by Hua--Song--Tudose~\cite{HST}. The weighted representation used here has the HST derivation recorded after Corollary~\ref{GC-cor:conditional-gaussian}. There is also a shorter comparison with Mazhar's undilated threefold theorem~\cite[Theorem~1.3]{R20-Mazhar}. At compact balanced input mass $1/2$, put $a=\Phi^{-1}(3/4)$ and $B=A+(2a/3)B_2^d$. Then $\gamma_d(B)\ge\Phi(2a/3)>2/3$ and
\[
 B+B+B=A+A+A+2aB_2^d\subset A^{+5}.
\]
Thus the chord lemma already converts that earlier threefold theorem into a fivefold half-mass conclusion. Theorem~\ref{R22-thm:fivefold} uses the weighted coefficients to lower the input threshold to $5/12$, with the strict finite-dimensional and compact Banach conclusions stated above. The count five is not asserted to be optimal. The bounds $a>2/3$ and $\Phi(4/9)>2/3$ follow by integrating the alternating Taylor bounds for $e^{-x^2/2}$ and using $25/8<\pi<22/7$.

\subsubsection{A Gaussian enlargement inside five summands}
The chord ball in the preceding proof has two uses. Part of its radius raises the input mass enough to apply the weighted theorem. The remaining radius enlarges the resulting convex core while staying inside the same five summands. Keeping these two radii separate gives a quantitative output bound in the original Banach space.

\begin{theorem}[Fivefold Gaussian enlargement]
\label{R23-thm:five-enlargement}
Let $\gamma$ be a centered Radon Gaussian measure on a separable Banach space $E$, and let $A\subset E$ be compact and balanced with $\gamma(A)\ge\alpha$, where $0<\alpha<1$. Set
\begin{equation}\label{R23-eq:five-radii}
 \rho_\alpha=2\Phi^{-1}((1+\alpha)/2),\qquad
 r_\alpha=\bigl[\Phi^{-1}(1-\alpha/2)-\Phi^{-1}(\alpha)\bigr]_+,
 \qquad d_\alpha=\rho_\alpha-r_\alpha.
\end{equation}
There is a compact symmetric convex set $K$ with
\begin{equation}\label{R24-eq:weighted-core}
 \gamma(K)\ge\tfrac12,\qquad
 K\subset A+\tfrac12A+\tfrac12A+r_\alpha B_H.
\end{equation}
If $d_\alpha\ge0$, the same set satisfies
\begin{equation}\label{R23-eq:five-enlargement}
 \gamma(K)\ge\tfrac12,\qquad
 K+d_\alpha B_H\subset A^{+5},\qquad
 \gamma(K+d_\alpha B_H)\ge\Phi(d_\alpha).
\end{equation}
The enlarged set is compact in $E$. For balanced Borel $A$ of mass at least $\alpha$, every $0\le d<d_\alpha$ admits the same conclusion with $d$ in place of $d_\alpha$.
\end{theorem}
\begin{proof}
For $\varepsilon>0$ put $B_\varepsilon=A+(r_\alpha+\varepsilon)B_H$. Gaussian isoperimetry gives
\[
 \gamma(A)+2\gamma(B_\varepsilon)
 \ge\alpha+2\Phi\bigl(\Phi^{-1}(\alpha)+r_\alpha+\varepsilon\bigr)>2.
\]
The inequality also holds when $r_\alpha=0$, since then $\alpha\ge2/3$. The weighted assertion of Theorem~\ref{GC-thm:banach} gives a compact symmetric convex set $K_\varepsilon$ of mass at least $1/2$ in
\[
 A+\tfrac12A+\tfrac12A+(r_\alpha+\varepsilon)B_H.
\]
For $0<\varepsilon\le1$ these sets lie in one compact sum. Compactness of $B_H$ follows from Lemma~\ref{R22-lem:chords}; the zero Gaussian is immediate. Hausdorff compactness and upper semicontinuity of measure give a compact symmetric convex limit $K$ with $\gamma(K)\ge1/2$ and
\[
 K\subset A+\tfrac12A+\tfrac12A+r_\alpha B_H.
\]
If $d_\alpha\ge0$, adding $d_\alpha B_H$ and applying the chord inclusion yields
\[
 K+d_\alpha B_H
 \subset A+\tfrac12A+\tfrac12A+\rho_\alpha B_H
 \subset A+\tfrac12A+\tfrac12A+A+A\subset A^{+5}.
\]
The half-copies of $A$ remain separate. Cameron--Martin isoperimetry gives the last inequality in \eqref{R23-eq:five-enlargement}; the sum is compact, symmetric and convex.

For Borel $A$, choose compact subsets with measures tending to $\gamma(A)$ and take their compact balanced hulls inside $A$. Continuity of $d_\alpha$ permits a compact hull of mass $\alpha'$ with $d_{\alpha'}>d$. Apply the compact assertion to that hull, retaining only radius $d$. This argument requires strict inequality in the radius and makes no compactness assertion about the original Borel set.
\end{proof}

\begin{corollary}[Uniform output mass and the half-input endpoint]
\label{R23-cor:five-mass}
Every balanced Borel set of mass at least $5/12$ in a separable Gaussian Banach space has a compact symmetric convex subset $C\subset A^{+5}$ satisfying
\begin{equation}\label{R23-eq:five-mass}
 \gamma(C)\ge\Phi(1/16)>0.524917.
\end{equation}
For compact balanced input of mass at least $1/2$, five summands contain such a set of mass at least $3/4$. More precisely, in this case \eqref{R23-eq:five-enlargement} holds with $d=\Phi^{-1}(3/4)$. These conclusions include Wiener measure and actual pathwise sums.
\end{corollary}
\begin{proof}
We verify the strict margin $d_{5/12}>1/16$ using rational bounds. Put $P_3(x)=x-x^3/6+x^5/40-x^7/336$ and
$P_4(x)=P_3(x)+x^9/3456$. Alternating exponential remainders give
\[
 \Phi(x)\le\tfrac12+\tfrac25P_4(x),\qquad
 \Phi(x)\ge\tfrac12+\tfrac{199}{500}P_3(x)
 \quad(0\le x\le5/6).
\]
Direct rational substitution proves
\[
 \Phi(35/64)<17/24,\qquad
 \Phi(3/14)>7/12,\qquad
 \Phi(49/60)>19/24.
\]
Consequently
\[
 d_{5/12}=2\Phi^{-1}(17/24)-\Phi^{-1}(19/24)+\Phi^{-1}(5/12)
 >\frac{35}{32}-\frac{49}{60}-\frac3{14}
 =\frac{211}{3360}>\frac1{16}.
\]
Theorem~\ref{R23-thm:five-enlargement} gives the Borel assertion. For the decimal lower bound, $\pi<355/113$ implies $1/\sqrt{2\pi}>398942/10^6$. Thus
\[
 \Phi(1/16)\ge\tfrac12+\frac{398942}{10^6}P_3(1/16)
 >\frac{524917}{10^6}.
\]
At $\alpha=1/2$ the radii are $\rho_\alpha=2\Phi^{-1}(3/4)$ and $r_\alpha=\Phi^{-1}(3/4)$, proving the compact endpoint.
\end{proof}

The profile also specifies its sufficient input threshold exactly. There is a unique $\alpha_*\in(0,1/2)$ with $d_{\alpha_*}=0$: $\rho_\alpha$ strictly increases, $r_\alpha$ decreases to zero, and $d_\alpha$ tends to $-\infty$ as $\alpha\downarrow0$, while $d_{1/2}>0$. Compact balanced input of mass at least $\alpha_*$ therefore gives fivefold half-mass output; Borel input of mass strictly greater than $\alpha_*$ gives output strictly greater than one half in every separable Gaussian Banach space. This is the threshold of the stated radius profile; no optimality among all fivefold constructions is asserted.

\begin{corollary}[Sixfold convexity and the high-mass profile]
\label{R22-cor:sixfold}
Let $A$ be compact and balanced in a separable Gaussian Banach space, with $\gamma(A)\ge1/2$. For every integer $j\ge1$ there is a compact symmetric convex set $C_j$ such that
\begin{equation}\label{R22-eq:six-profile}
 C_j\subset A^{+6j},\qquad
 \gamma(C_j)\ge2\Phi\!\left(2j\Phi^{-1}(3/4)\right)-1.
\end{equation}
In particular six summands give output measure at least $\beta_4>0.822656449$. For output mass $1-\delta$, the required number has order $\sqrt{\log(1/\delta)}$ as $\delta\downarrow0$, which is sharp in order already for intervals. These assertions include Wiener measure and actual pathwise sums, with compactness in the original Banach norm.
\end{corollary}
\begin{proof}
Put $a=\Phi^{-1}(3/4)$. At input mass $1/2$,
\eqref{R24-eq:weighted-core} gives $\gamma(K)\ge1/2$ and
$K\subset A+\tfrac12A+\tfrac12A+aB_H$.
The chord lemma gives $2aB_H\subset A+A$, so
$2K\subset2A+A+A+A+A\subset A^{+6}$.
The Gaussian $S$-inequality gives \eqref{R22-eq:six-profile} for
$C_j=2jK$. Its finite-dimensional form~\cite{R15-LO} extends to a
separable Banach space by decreasing cylinder approximations to the
closed symmetric convex set, using countably many continuous linear
functionals. Compactness holds by \eqref{R24-eq:weighted-core}.
The Gaussian tail asymptotic gives the asserted order; the interval
$A=[-a,a]$ has $A^{+k}=[-ka,ka]$, proving the matching order lower bound.
\end{proof}

\begin{corollary}[Every positive input mass]
\label{R22-cor:positive-mass}
For $0<\alpha<1$ define
\[
 R_\alpha=\max\{0,\Phi^{-1}(1-\alpha/2)-\Phi^{-1}(\alpha)\},\qquad
 N(\alpha)=3+2\left\lceil R_\alpha/\rho_\alpha\right\rceil.
\]
Every compact balanced Gaussian set of mass at least $\alpha$ in a separable Banach space has a compact symmetric convex half-mass subset of $A^{+N(\alpha)}$.
\end{corollary}
\begin{proof}
Use \eqref{R24-eq:weighted-core}, which holds for every $0<\alpha<1$.
The chord lemma gives
$R_\alpha B_H\subset (R_\alpha/\rho_\alpha)A+
(R_\alpha/\rho_\alpha)A$. For $c>0$, balancedness implies
$cA\subset A^{+\lceil c\rceil}$ by writing $ca$ as a sum of
$\lceil c\rceil$ equal points of $A$. The two half-copies and the full
copy contribute three summands; each ball contribution contributes at
most $\lceil R_\alpha/\rho_\alpha\rceil$. Zero coefficients contribute
only zero. The compact endpoint is already included in the core theorem.
\end{proof}


\section{Contraction and stability of stretched Brownian motion}
\label{R11-sec:canonical}
Stretched Brownian motion is the unique optimizer of the martingale Benamou--Brenier problem: among fixed-endpoint martingales, it minimizes quadratic diffusion distance from Brownian motion~\cite{CM-BBHK}. Its Bass representation describes conditional endpoints by translated convex gradients~\cite{CM-Bass}. The question here is quantitative. How does curvature of the terminal law bound those maps, and how does loss in the variational objective control an arbitrary competitor on the same noise?

Caffarelli's contraction theorem compares an upper bound on source log-density curvature with a lower bound on terminal curvature~\cite{Caffarelli,Kolesnikov}. A Gaussian-mixture source has upper curvature at most the identity, so the established Bass representation fits that comparison. A fixed-terminal approximation removes irreducibility while retaining conditional kernels and paths. Strong convexity of the conjugate then turns a Fenchel gap into the objective deficit: the terminal-law terms cancel globally and the martingale means cancel the conditional linear terms. This yields the sharp same-driver estimate, including the finite-source integrable-dual passage needed for arbitrary competitors. Conditional future paths inherit Gaussian inequalities through the resulting Lipschitz maps.

\subsection{The canonical optimizer and its conditional maps}
\label{CM-sec:setup}
Let $\mu\cx\nu$ belong to $\mathcal P_2(\R^d)$, and write
$\mathcal M(\mu,\nu)$ for their martingale couplings. For
$q\in\mathcal P_2(\R^d)$ define
\[
 \MCov(q,\gamma_d)
 =\sup\{\E\langle Y,G\rangle:Y\sim q,\ G\sim\gamma_d\}.
\]
For $\pi(dx,dy)=\mu(dx)\pi_x(dy)$ put
\begin{equation}\label{CM-eq:primal}
 J(\pi)=\int\MCov(\pi_x,\gamma_d)\,\mu(dx),
 \qquad \mathcal P(\mu,\nu)=\max_{\pi\in\mathcal M(\mu,\nu)}J(\pi).
\end{equation}
The martingale Benamou--Brenier theorem gives a unique maximizer
$\pi^*$ and identifies its heat interpolation as stretched Brownian
motion \cite{CM-BBHK}. If $T_x=\nabla v_x$ is the Brenier map from
$\gamma_d$ to $\pi_x^*$, jointly measurable versions give
\begin{equation}\label{CM-eq:heat}
 M_t^*=P_{1-t}T_X(B_t),\qquad
 P_s f(z)=\int f(z+\sqrt{s}\,g)\,\gamma_d(dg),
\end{equation}
where $X\sim\mu$ is independent of a standard Brownian motion $B$.
The equality $P_1T_x(0)=x$ gives $M_0^*=X$, and $M_1^*=T_X(B_1)$.

\begin{theorem}[Contraction of the canonical conditional maps]
\label{CM-thm:contraction}
Suppose
\begin{equation}\label{CM-eq:curvature}
 \nu(dy)=Z^{-1}e^{-V(y)}\,dy,\qquad
 V(y)-\frac{|y|^2}{2L^2}\ \text{is lower semicontinuous and convex},
 \quad 0<L<\infty.
\end{equation}
The potential $V$ may take the value $+\infty$. For every
$\mu\cx\nu$, the conditional maps of the unique optimizer in
\eqref{CM-eq:primal} have versions satisfying
\begin{equation}\label{CM-eq:contraction}
 \Lip(T_x)\le L,\qquad 0\preceq DT_x\preceq LI_d
 \quad\text{almost everywhere},\quad\mu\text{-almost every }x.
\end{equation}
No irreducibility hypothesis is imposed on $(\mu,\nu)$. In
\eqref{CM-eq:heat},
\begin{equation}\label{CM-eq:cap}
 dM_t^*=\sigma_t^*\,dB_t,\qquad
 \sigma_t^*=D(P_{1-t}T_X)(B_t),\qquad
 0\preceq\sigma_t^*\preceq LI_d.
\end{equation}
In particular $d\langle M^*\rangle_t/dt\preceq L^2I_d$.
\end{theorem}

The curvature hypothesis concerns the whole terminal law $\nu$.
The conditional laws $\pi_x^*$ can be singular. We first approximate
the initial law while keeping $\nu$ fixed. On each finite approximation,
the Bass representation gives a Gaussian-mixture source, to which
Caffarelli's contraction theorem applies. The same approximation will
make the dual calculation in the deficit proof integrable.

\subsection{Keeping the terminal law fixed}\label{CM-sec:stability}
Finite conditional means followed by contraction put a fixed positive part of $\nu$ into every conditional law. Strict concavity then recovers the complete kernel when that perturbation vanishes. The argument also works for a non-Gaussian reference; this will allow the same approximation to be used in Section~\ref{R13-sec:qbass}.

\begin{lemma}[Strict concavity of maximal covariance]\label{CM-lem:strict}
Let $\rho\in\cP_2(\R^d)$ give no mass to sets of Hausdorff dimension at most $d-1$. Then $C_\rho(q)=\MCov(q,\rho)$ is $W_2$-continuous and strictly concave. For probability mixtures with finite averaged second moment,
\[
 C_\rho\left(\int q\,\Lambda(dq)\right)\ge\int C_\rho(q)\,\Lambda(dq),
\]
and equality forces $\Lambda$ to be a point mass. This includes $C(q)=\MCov(q,\gamma_d)$.
\end{lemma}
\begin{proof}
The identity $2C_\rho(q)=m_2(\rho)+m_2(q)-W_2^2(q,\rho)$ gives continuity. Mix the optimal transports from $\rho$ to the component laws. Equality makes this mixture optimal for the averaged target. Uniqueness of the convex-gradient transport from $\rho$ forces the same map for almost every component, hence the same target. This argument applies to general mixtures by measurable selection of the unique optimal couplings; alternatively partition the component laws by a countable determining family and use the two-point assertion.
\end{proof}

\begin{theorem}[Fixed-terminal stability under conditional means and contraction]\label{CM-thm:stability}
Let $\mu\cx\nu$ belong to $\cP_2(\R^d)$ and put $m=\int x\,\mu(dx)$. On $X\sim\mu$, let $b_n(X)$ be finite-valued conditional means satisfying
\[
 \E[X\mid b_n(X)]=b_n(X),\qquad b_n(X)\longrightarrow X\text{ in }L^2.
\]
The choice $b_n(X)=X$ is also permitted. For $0<r_n<1$ with $r_n\to1$, set
$X_n=m+r_n(b_n(X)-m)$ and $\mu_n=\law(X_n)$.
Write $\pi^{(n)},\pi^*$ for the Gaussian-reference optimizers and $q_n(x)=\pi^{(n)}_{X_n(x)}$. Then
\begin{align}
 \mathcal P(\mu_n,\nu)&\to\mathcal P(\mu,\nu),\label{CM-eq:valueconv}\\
 \int W_2^2(q_n(x),\pi_x^*)\,\mu(dx)&\to0,\label{CM-eq:kernelconv}\\
 \int\!\int|T_{q_n(x)}(z)-T_{\pi_x^*}(z)|^2\,\gamma_d(dz)\mu(dx)&\to0.\label{CM-eq:mapconv}
\end{align}
Their conditional Bass processes on the same $X,B$ converge in $L^2$ of the uniform path norm. More generally, for every reference $\rho$ in Lemma~\ref{CM-lem:strict}, the value, kernel and map assertions hold with $P_\rho$, its optimizers and maps from $\rho$. Each approximating pair has a feasible kernel whose every conditional law dominates $(1-r_n)\nu$.
\end{theorem}
\begin{proof}
Coarsen an optimal kernel by conditioning on $b_n(X)$, then mix it with $\nu$ in proportions $r_n,1-r_n$. Concavity gives feasibility and
\begin{equation}\label{CM-eq:valuelower}
 P_\rho(\mu_n,\nu)\ge r_nP_\rho(\mu,\nu)+(1-r_n)C_\rho(\nu).
\end{equation}
Lift the optimal kernels to $\law(X,q_n(X))$. Their averaged terminal law is $\nu$, and their conditional mean is $X_n\to X$. Lemma~\ref{AS-lem:lift} gives compactness, moment control and feasible averaged limits. Continuity and concavity give the reverse value bound to~\eqref{CM-eq:valuelower}. Equality and Lemma~\ref{CM-lem:strict} make every limiting random kernel a point mass at the unique optimal kernel. Thus the convergence is along the same $X$ and is integrated in $W_2^2$. In particular, with $Z_n(x)=m_2(q_n(x))$, the fixed terminal law gives
\begin{equation}\label{CM-eq:UI}
 \int Z_n\mathbf1_{\{Z_n>K\}}\,d\mu
 \le2\int|y|^2\mathbf1_{\{|y|^2>K/2\}}\,d\nu\longrightarrow0.
\end{equation}
Optimal maps from a fixed nondegenerate source converge in $L^2$ when their target laws converge in $W_2$: limits of optimal graph couplings remain optimal, uniqueness identifies the graph, and convergence of moments gives strong convergence. Apply this fiberwise and use~\eqref{CM-eq:UI}. For $\rho=\gamma_d$, conditional expectation and Doob yield
\[
 \E\sup_{t\le1}|M_t^{(n)}-M_t^*|^2
 \le4\E|T_{q_n(X)}(B_1)-T_{\pi_X^*}(B_1)|^2\to0.
\]
\end{proof}

This fixed-terminal approximation works without the compact interior hypothesis of Theorem~\ref{AS-thm:main}; it specifies the allowed initial perturbations instead. General one-dimensional all-marginal stability is treated in~\cite{CM-BJMP}. Here the same reference input also retains the complete conditional maps in every dimension.

\subsection{Contraction and an integrable finite-source dual}\label{CM-sec:contraction}
For an irreducible pair, the Bass representation first transports a
Gaussian convolution to $\nu$. The negative logarithm of the convolution
density has Hessian at most $I_d$, regardless of whether the mixing law
is log-concave. The upper-source/lower-target form of the contraction
theorem therefore gives a common Lipschitz bound. Finite initial
support has a second role: it makes the latent law finite and the
conjugate potential integrable against $\nu$.

\begin{lemma}[Contraction from a Gaussian convolution]\label{CM-lem:mixture}
Let $\alpha$ be any probability law on $\R^d$, without a moment assumption, and let $\beta=\alpha*\gamma_d$. If $\nu$ satisfies~\eqref{CM-eq:curvature}, the unique gradient-of-convex transport from $\beta$ to $\nu$ has an $L$-Lipschitz representative.
\end{lemma}
\begin{proof}
Write the density of $\beta$ in the form
\[
 \beta(z)=(2\pi)^{-d/2}e^{-|z|^2/2}
     \int e^{\ip za-|a|^2/2}\,\alpha(da).
\]
The logarithm of the integral is convex. Thus the source potential $U=-\log\beta$ satisfies
\begin{equation}\label{CM-eq:sourcecurvature}
 D^2U\preceq I_d.
\end{equation}
For boundedly supported $\alpha$ and smooth finite target potential, the upper-source/lower-target form of Caffarelli's theorem gives $\Lip(\nabla v)\le L$; a precise formulation is~\cite[Theorem~2.2]{Kolesnikov}. Source log-concavity is unnecessary in this formulation.

For completeness, the limiting step handles both qualifications in the statement. Truncate $\alpha$ to large balls and renormalize, obtaining $\alpha_j$, so $\alpha_j*\gamma_d\Rightarrow\beta$. Approximate the convex part $W=V-|y|^2/(2L^2)$ by smooth finite convex Moreau regularizations. An affine minorant $W\ge\ip ay+b$ gives a common Gaussian-integrable bound after normalization. The target laws $\nu_j$ consequently converge to $\nu$ in every finite Wasserstein distance. Let $F_j$ be the $L$-Lipschitz transport maps for these smooth finite-moment approximations.

Their values at the origin are uniformly bounded. Choose $R,p>0$ so that every sufficiently large approximating source gives $B_R$ probability at least $p$. Integrating $|F_j(0)|\le|F_j(z)|+LR$ on that ball gives
\[
 |F_j(0)|\le LR+p^{-1}\int|y|\,\nu_j(dy).
\]
By Arzel\`a--Ascoli, a subsequence converges locally uniformly to an $L$-Lipschitz map $F$. Tightness of the source laws shows $F_\#\beta=\nu$. Each $F_j$ is a gradient of a convex function; local uniform convergence preserves this property, by convergence of normalized potentials or cyclic monotonicity. McCann's extension of uniqueness of monotone transport maps to sources without moment assumptions identifies $F$ with the given gradient transport; see the discussion in~\cite[Section~1.4]{CM-Bass}. This also proves the statement for extended-valued $V$.
\end{proof}

\begin{lemma}[The irreducible calculation]\label{CM-lem:Bass}
Under \eqref{CM-eq:curvature}, Theorem~\ref{CM-thm:contraction} holds
when $(\mu,\nu)$ is irreducible. There are a probability law $\alpha$,
a convex $v$ and a measurable $a(x)$ such that
\[
 (\nabla v)_\#(\alpha*\gamma_d)=\nu,
 \qquad T_x(g)=\nabla v(a(x)+g),\qquad \Lip(\nabla v)\le L.
\]
\end{lemma}
\begin{proof}
The Bass representation theorem~\cite[Theorem~1.3]{CM-Bass} gives
$F=\nabla v$ and a latent variable $A$ with $X=P_1F(A)$.
Lemma~\ref{CM-lem:mixture} gives $\Lip(F)\le L$.
The map $P_1F$ is injective. Indeed, if $P_1F(a)=P_1F(b)$ and
$h=b-a\ne0$, monotonicity gives
\[
 0=\int_0^1\!\int h^{\mathsf T}DF(a+th+g)h\,\gamma_d(dg)\,dt.
\]
Positivity of the integrand and of the Gaussian density imply
$DFh=0$ almost everywhere. Symmetry of $DF$ makes $h\cdot F$
constant, contradicting the full-dimensional terminal law. Thus
$A=a(X)$; conditional optimality and Brenier uniqueness give the
asserted maps. In score notation, the source curvature identity is
\begin{equation}\label{CM-eq:mixture-score}
 D^2(-\log(\alpha*\varphi_d))(z)
 =I_d-\Cov(A\mid A+G=z)\preceq I_d.
\end{equation}
The conditional moments in this formula are finite locally even when
$\alpha$ has no second moment, by the Gaussian likelihood.
\end{proof}

\begin{lemma}[Finite initial law and full conditional communication]\label{CM-lem:finiteBass}
Assume~\eqref{CM-eq:curvature}. Suppose $\mu$ is finitely supported and some martingale coupling of $\mu,\nu$ has each conditional law bounded below by a positive multiple of $\nu$. Then the optimal coupling has the representation
\begin{equation}\label{CM-eq:finiteBass}
 T_x(z)=F(a(x)+z),\qquad F=\nabla v,\qquad \Lip(F)\le L,
\end{equation}
with finitely many vectors $a(x)$. The conjugate $\psi=v^*$ is $1/L$-strongly convex and belongs to $L^1(\nu)$.
\end{lemma}
\begin{proof}
The common positive component of every conditional law makes the pair
irreducible. Lemma~\ref{CM-lem:Bass} gives $F=\nabla v$ and
$A=a(X)$. Because $X$ has finite support, $A$ has finite range.
The conjugate of a convex function with $L$-Lipschitz gradient is
$1/L$-strongly convex. Moreover, $A+G$ has finite second moment,
$F$ has linear growth and $v$ has at most quadratic growth. Fenchel
equality therefore makes
\[
 \psi(F(A+G))=\langle A+G,F(A+G)\rangle-v(A+G)
\]
absolutely integrable. Since $F(A+G)\sim\nu$, this proves
$\psi\in L^1(\nu)$.
\end{proof}

\begin{lemma}[Passing uniform regularity and its deficit to the limit]\label{R15-lem:transfer}
Fix a positive-density $\rho\in\cP_2$ and $\nu\in\cP_2$. Suppose every finite initial law with a positive $\nu$ component in each conditional law has an optimal representation
$T_x(w)=\nabla v(a(x)+w)$ with $\Lip(\nabla v)\le M$, uniformly in that law. Then every $\mu\cx\nu$ has $M$-Lipschitz conditional convex-gradient maps, and every feasible triple with $W\sim\rho$, $W\perp X$ obeys
\begin{equation}\label{R15-eq:transfer-deficit}
 P_\rho(\mu,\nu)-\E\langle Y,W\rangle
 \ge\frac1{2M}\E|Y-T_X(W)|^2.
\end{equation}
\end{lemma}
\begin{proof}
Use the finite approximations in Theorem~\ref{CM-thm:stability}. Their maps converge in $L^2(\mu\otimes\rho)$. On almost every fiber, their common Lipschitz bound and $L^1(\rho)$ bound control one point value. Local compactness gives the unique continuous limit with the same bound; positive density identifies it everywhere. Normalized convex potentials retain the convex-gradient structure, and values on a countable dense set give joint measurability.

For a finite pair, set $Z=a(X)+W$, $Y^*=\nabla v(Z)$. Finite shifts and quadratic growth make $v(Z)$ integrable. Fenchel equality gives $v^*\in L^1(\nu)$, so strong convexity of $v^*$ may be integrated against every competitor:
\[
 v^*(Y)-v^*(Y^*)-\langle Z,Y-Y^*\rangle
 \ge\frac{|Y-Y^*|^2}{2M}.
\]
The complete terminal marginals cancel the conjugate terms; exact conditional means cancel the term in $a(X)$. This proves~\eqref{R15-eq:transfer-deficit} for the finite pair. For an arbitrary feasible triple, independently replace $Y$ by $Y'\sim\nu$ with probability $1-r_n$, calling the result $Y_n$. Then $Y_n\sim\nu$, $W\perp X_n$, $\E[Y_n\mid X_n]=X_n$, and $Y_n\to Y$ in $L^2$. Apply the finite inequality and pass to the limit using the value and same-input map convergence of Theorem~\ref{CM-thm:stability}. All limiting operations retain the unchanged terminal law.
\end{proof}

\begin{proof}[Proof of Theorem~\ref{CM-thm:contraction}]
Lemma~\ref{CM-lem:finiteBass} verifies the preceding lemma for $\rho=\gamma_d$ and $M=L$. Heat smoothing preserves the Hessian bounds $0\preceq DT_x\preceq LI_d$. It\^o's formula gives~\eqref{CM-eq:cap}; the Lipschitz bound gives continuity at time one.
\end{proof}

\begin{corollary}[Arbitrarily small initial contraction]
\label{CM-lem:contraction-limit}
Under \eqref{CM-eq:curvature}, let $m=\int x\,\mu(dx)$,
$x_\delta=(1-\delta)x+\delta m$ and
$\mu_\delta=((1-\delta)\id+\delta m)_\#\mu$. Then
\begin{equation}\label{CM-eq:limit}
 \mathcal P(\mu_\delta,\nu)\longrightarrow\mathcal P(\mu,\nu),\qquad
 \int\!\int|T^\delta_{x_\delta}(g)-T_x(g)|^2
       \,\gamma_d(dg)\mu(dx)\longrightarrow0.
\end{equation}
All the maps have $L$-Lipschitz convex-gradient representatives.
\end{corollary}
\begin{proof}
Use $b_n(X)=X$ in Theorem~\ref{CM-thm:stability} and apply
Theorem~\ref{CM-thm:contraction}.
\end{proof}

\subsection{The variational deficit controls the coupling}
\label{CM-sec:deficit}
An $L$-Lipschitz convex gradient has a $1/L$-strongly convex conjugate.
At its contact point this gives
\begin{equation}\label{CM-eq:Fenchel}
 v(g)+v^*(y)-\langle g,y\rangle
 \ge\frac1{2L}|y-\nabla v(g)|^2.
\end{equation}
For the finite-source dual, every term in this inequality is
integrable. Its expectation becomes the objective deficit because
both the terminal marginal and the conditional linear term are fixed.
The finite-source approximation then retains this identity for every
competing coupling on the same Gaussian input.

\begin{theorem}[Sharp endpoint deficit]\label{CM-thm:deficit}
Under \eqref{CM-eq:curvature}, let $(X,Y)\in\mathcal M(\mu,\nu)$,
and let $G\sim\gamma_d$ be independent of $X$. Its dependence on
$Y$ is arbitrary. Then
\begin{equation}\label{CM-eq:deficit}
 \mathcal P(\mu,\nu)-\E\langle Y,G\rangle
 \ge\frac1{2L}\E|Y-T_X(G)|^2.
\end{equation}
The constant $1/(2L)$ is sharp uniformly over the class. In particular,
for every $\pi\in\mathcal M(\mu,\nu)$,
\begin{equation}\label{CM-eq:kernel-stability}
 \int W_2^2(\pi_x,\pi_x^*)\,\mu(dx)
 \le2L\{\mathcal P(\mu,\nu)-J(\pi)\}.
\end{equation}
\end{theorem}
\begin{proof}
Apply Lemma~\ref{R15-lem:transfer} with $\rho=\gamma_d$ and $M=L$. For~\eqref{CM-eq:kernel-stability}, couple each competing conditional law to $G$ by its Brenier map, attaining $J(\pi)$. For sharpness take $\mu=\delta_0$, $\nu=N(0,L^2I_d)$ and $Y=-LG$. Both sides of~\eqref{CM-eq:deficit} equal $2Ld$.
\end{proof}
The equality example proves sharpness for the full triple inequality.
It also determines the same-driver diffusion constant below; it makes
no separate sharpness assertion for the weaker kernel-$W_2$ bound.

\begin{corollary}[Same-driver diffusion stability]\label{CM-cor:dynamic}
On a filtration for which $B$ is Brownian and $X\in\mathcal F_0$ is
independent of $B$, let
$M_t=X+\int_0^t\sigma_s\,dB_s$ have initial law $\mu$ and terminal
law $\nu$. Let $M^*$ be \eqref{CM-eq:heat} on the same $(X,B)$, and
put
\[
 \Delta=\mathcal P(\mu,\nu)-\E\int_0^1\tr\sigma_t\,dt.
\]
Then
\begin{equation}\label{CM-eq:dynamic}
 \E\int_0^1\|\sigma_t-\sigma_t^*\|_{\rm HS}^2\,dt
 =\E|M_1-M_1^*|^2\le2L\Delta,
 \qquad
 \E\sup_{0\le t\le1}|M_t-M_t^*|^2\le8L\Delta.
\end{equation}
The factor $2L$ is sharp for a fixed Brownian driver.
\end{corollary}
\begin{proof}
It\^o isometry gives
$\E\langle M_1,B_1\rangle=\E\int_0^1\tr\sigma_t\,dt$.
Apply Theorem~\ref{CM-thm:deficit} with $Y=M_1$, $G=B_1$.
The equality in \eqref{CM-eq:dynamic} is It\^o isometry for
$M-M^*$, and the last estimate is Doob's $L^2$ inequality.
For sharpness use $X=0$, $M_t=-LB_t$, and $M_t^*=LB_t$.
The coefficient error is $4L^2d$ and the deficit is $2Ld$.
\end{proof}

For the action $\mathcal A(M,B)=\E\int_0^1\|\sigma_t-I_d\|_{\rm HS}^2dt$,
fixed marginals imply $\mathcal A-\mathcal A^*=2\Delta$. Thus
\begin{equation}\label{CM-eq:pythagorean}
 \mathcal A(M,B)\ge\mathcal A(M^*,B)
 +\frac1L\E\int_0^1\|\sigma_t-\sigma_t^*\|_{\rm HS}^2dt.
\end{equation}
The driver is part of this comparison: replacing it can change the
covariance objective even when the path law is unchanged.

\begin{corollary}[An a posteriori primal--dual certificate]
\label{CM-cor:certificate}
Let $U$ be a finite valid upper bound for $\mathcal P(\mu,\nu)$.
The deficit and diffusion estimates hold with $\mathcal P$ replaced
by $U$. In particular one may use a convex dual test $\psi$ whose
Fenchel dual value
\[
 U_\psi=\int\psi\,d\nu-\int(P_1\psi^*)^*(x)\,\mu(dx)
\]
is finite, with both integrals well-defined and integrable.
Evaluating or bounding these integrals and the feasible primal value
gives an error certificate without evaluating the optimizer.
\end{corollary}
\begin{proof}
For a feasible triple, Fenchel's inequality at $G+a(X)$ gives
\[
 \E\langle Y,G\rangle\le\int\psi\,d\nu+
 \int\{P_1\psi^*(a(x))-\langle a(x),x\rangle\}\,\mu(dx).
\]
Measurable approximate minimization over $a(x)$ gives
$\mathcal P(\mu,\nu)\le U_\psi$, the usual integrable dual upper
bound~\cite[Theorem~1.4]{CM-Bass}. Substitution in the established
deficit estimates proves the assertion. Any independently justified
finite dual upper bound can be used in the same way.
\end{proof}

\begin{proposition}[Bounds for the optimal action]\label{CM-prop:values}
Let $D=\Cov(\nu)-\Cov(\mu)\succeq0$. Under~\eqref{CM-eq:curvature},
\begin{equation}\label{CM-eq:valuebounds}
 \frac{\tr D}{L}\le\mathcal P(\mu,\nu)\le\tr\sqrt D.
\end{equation}
The lower constant $1/L$ is sharp. Consequently
\[
 d+\tr D-2\tr\sqrt D
 \le\mathcal A^*
 \le d+\tr D-\frac{2}{L}\tr D.
\]
For $\nu=\gamma_d$ and $L=1$, this gives $\mathcal A^*\le\tr\Cov(\mu)$.
\end{proposition}
\begin{proof}
The canonical diffusion satisfies $0\preceq\sigma_t^*\preceq LI_d$, so
$(\sigma_t^*)^2\preceq L\sigma_t^*$. Taking expected traces and integrating gives the lower bound. Also
\[
 D=\E\int_0^1(\sigma_t^*)^2dt.
\]
Concavity of $A\mapsto\tr\sqrt A$ gives
$\E\int\tr\sigma_t^*dt\le\tr\sqrt D$. For $\mu=\delta_0$ and $\nu=N(0,L^2I_d)$ both bounds are equal to $Ld$. The action statements follow from its fixed-marginal identity.
\end{proof}

\subsection{Conditional functional inequalities and future-path geometry}\label{CM-sec:functional}\label{CM-sec:path}
The optimal conditional laws are Lipschitz images of Gaussian measure. This representation gives conditional functional inequalities even though those laws are not asserted to be log-concave.

\begin{lemma}[Functional inequalities for Lipschitz Gaussian images]\label{CM-lem:functional}
Suppose $K=F_\#N(0,sI_d)$, where $F$ is $L$-Lipschitz and $s\ge0$. Then, for smooth test functions and with the usual extension by closure,
\begin{align}
 \Var_K(f)&\le L^2s\int|\nabla f|^2\,dK,\label{CM-eq:poincare}\\
 \Ent_K(f^2)&\le2L^2s\int|\nabla f|^2\,dK.\label{CM-eq:LSI}
\end{align}
For every probability $\eta$ with finite relative entropy against $K$,
\begin{equation}\label{CM-eq:T2}
 W_2^2(\eta,K)\le2L^2s\,D(\eta\Vert K).
\end{equation}
Translations of the Gaussian source give the same constants. At $s=0$ the assertions have their point-mass interpretation.
\end{lemma}
\begin{proof}
Apply the Gaussian Poincar\'e and logarithmic Sobolev inequalities to $f\circ F$. The almost-everywhere chain rule gives
$|\nabla(f\circ F)|\le L|\nabla f|\circ F$; approximation extends the argument to nonsmooth $F$. This proves~\eqref{CM-eq:poincare}--\eqref{CM-eq:LSI}; the Gaussian logarithmic Sobolev inequality is classical~\cite{Gross1975}.

For~\eqref{CM-eq:T2}, let $r=d\eta/dK$ and lift $\eta$ to the Gaussian source by the density $r\circ F$. The lift has exactly entropy $D(\eta\Vert K)$. Talagrand's Gaussian transport inequality~\cite{CM-TalagrandT2} bounds its squared $W_2$ distance from the Gaussian by $2sD(\eta\Vert K)$. Push an optimal coupling through $F$ in both coordinates; its cost grows by at most $L^2$. This also proves the result for singular image measures.
\end{proof}

\begin{theorem}[Conditional inequalities at stopping times]\label{CM-thm:conditional}
Realize $M^*$ by~\eqref{CM-eq:heat} in the filtration generated by $X$ and $B$. If $\tau\le t\le1$, where $\tau$ is a stopping time and $t$ is deterministic, the regular conditional law
\[
 K_{\tau,t}=\law(M_t^*\mid\cF_\tau)
\]
satisfies~\eqref{CM-eq:poincare}--\eqref{CM-eq:T2} almost surely with $s=t-\tau$. In particular every terminal kernel $\pi_x^*$ satisfies a logarithmic Sobolev inequality with constant $L^2$ and a $T_2$ inequality with coefficient $2L^2$.
\end{theorem}
\begin{proof}
Write $F_{x,t}=P_{1-t}T_x$, with $F_{x,1}=T_x$. Each is $L$-Lipschitz. Given $\cF_\tau$, the conditional law is the image of $N(0,(t-\tau)I_d)$ under
\[
 z\longmapsto F_{X,t}(B_\tau+z).
\]
The strong Markov property of Brownian motion gives the same conditional law at stopping times. Apply Lemma~\ref{CM-lem:functional} conditionally. The filtration and the chosen process are part of the conclusion; no arbitrary coarsening of a conditional law is used.
\end{proof}

\begin{theorem}[Gaussian isoperimetry and transportation for the whole future path]\label{CM-thm:path}
Let $\tau<1$ be a stopping time in the filtration of Theorem~\ref{CM-thm:conditional}. Conditional on $\cF_\tau$, let $\mathsf K_\tau$ be the law of the future path $(M_t^*)_{\tau\le t\le1}$ with its uniform norm. Then
\begin{equation}\label{CM-eq:path-T2}
 W_{2,\infty}^2(\eta,\mathsf K_\tau)
 \le 2L^2(1-\tau)D(\eta\Vert\mathsf K_\tau)
\end{equation}
for every $\eta$ with finite relative entropy. Furthermore, for every Borel path set $A$ and $r>0$,
\begin{equation}\label{CM-eq:isoperimetry}
 \mathsf K_\tau(A^r)
 \ge \Phi\left(\Phi^{-1}(\mathsf K_\tau(A))
                    +\frac{r}{L\sqrt{1-\tau}}\right),
\end{equation}
where $A^r$ is its open uniform-distance $r$-neighbourhood and $\Phi$ is the standard one-dimensional Gaussian distribution function. The same isoperimetric estimate applies to each kernel $K_{\tau,t}$ of Theorem~\ref{CM-thm:conditional}, with $1-\tau$ replaced by $t-\tau$ and Euclidean neighbourhoods.
\end{theorem}
\begin{proof}
Condition on $\cF_\tau$ and put $h=1-\tau$, $x=X$, $b=B_\tau$. The future law is the image of Wiener measure on $[0,h]$ under
\[
 \Psi(w)_s=(P_{1-\tau-s}T_x)(b+w_s).
\]
The maps are uniformly $L$-Lipschitz and continuous up to the terminal time. For a Cameron--Martin shift $u$,
\[
 \|\Psi(w+u)-\Psi(w)\|_\infty
 \le L\|u\|_\infty\le L\sqrt h\,\|u\|_H,
 \qquad \|u\|_H^2=\int_0^h|\dot u_s|^2ds.
\]
Use the entropy-preserving density lift of Lemma~\ref{CM-lem:functional}, now on Wiener space. Its Cameron--Martin $T_2$ inequality has coefficient two~\cite[Section~3]{CM-Wiener}; pushing the coupling through $\Psi$ gives~\eqref{CM-eq:path-T2}. Borell's Wiener-space isoperimetry~\cite{CM-Borell} applies to $\Psi^{-1}(A)$. A Cameron--Martin ball of radius $r/(L\sqrt h)$ maps into $A^r$, giving~\eqref{CM-eq:isoperimetry}. Inner approximation and radii increasing to $r$ handle arbitrary Borel sets. For $K_{\tau,t}$, the same argument uses its $L\sqrt{t-\tau}$-Lipschitz standard Gaussian representation from Theorem~\ref{CM-thm:conditional}.

For sharpness take $\mu=\delta_0$, $\nu=N(0,L^2I_d)$ and deterministic $\tau$. The future law is translated $L$-scaled Brownian motion of duration $h$. Shifting by $s\mapsto sa/h$ costs entropy $|a|^2/(2L^2h)$ and uniform squared transport cost $|a|^2$: translation attains the terminal-mean lower bound. Terminal-coordinate half-spaces attain the isoperimetric profile.
\end{proof}

\subsection{Canonical kernels for the same rounding law}\label{CM-sec:rounding}
\begin{corollary}[Optimal proximal kernels]\label{CM-cor:proximal}
When $L=1$, for $\mu$-almost every $x$ there is a proper lower semicontinuous convex function $h_x$ such that
\begin{equation}\label{CM-eq:proximal}
 T_x=\prox_{h_x},\qquad I-T_x=\prox_{h_x^*}.
\end{equation}
Here $\prox_h(z)$ is the unique minimizer of $h(y)+|z-y|^2/2$. Thus, with $Z\sim\gamma_d$ independent of $X$,
\[
 Y=\prox_{h_X}(Z),\qquad
 \law(Y)=\nu,\qquad \E[Y\mid X]=X,
\]
and this is the optimal conditional Gaussian-source realization. For general $L$, $T_x/L$ is a proximal map.
\end{corollary}
\begin{proof}
Choose $T_x=\nabla v_x$. Theorem~\ref{CM-thm:contraction} says that $v_x^*$ is $1$-strongly convex. Thus
$h_x=v_x^*-|\cdot|^2/2$ is proper lower semicontinuous and convex. Minimizing
\[
 h_x(y)+|z-y|^2/2=v_x^*(y)-\ip zy+|z|^2/2
\]
gives $y=\nabla v_x(z)$. Moreau's identity gives the complementary formula in~\eqref{CM-eq:proximal}; see~\cite{Moreau}. Apply the same argument to $v_x/L$ for general $L$.
\end{proof}

In particular the conditional maps satisfy the standard firm-nonexpansiveness estimate, when $L=1$,
\begin{equation}\label{CM-eq:firm}
 |T_x(z)-T_x(w)|^2
 +|(z-T_x(z))-(w-T_x(w))|^2\le|z-w|^2.
\end{equation}
The martingale variational problem therefore selects a specific proximal map, whose complementary proximal map describes the residual Gaussian displacement.

\begin{remark}[Matrices and the choice of metric]\label{CM-rem:metric}
If $V(y)-y^{\mathsf T}Q^{-1}y/2$ is convex with $Q\succ0$, whiten by $Q^{-1/2}$ and apply all the preceding theorems. The result is the optimizer for the \emph{whitened} martingale Benamou--Brenier objective and has bracket density at most $Q$ after transformation back. This is a weighted objective in the original coordinates. A general linear change of coordinates need not preserve the unweighted Euclidean objective, so the theorem does not identify that different optimizer with the same matrix cap.
\end{remark}

\begin{corollary}[The whole signing law and a canonical reference]\label{CM-cor:signing}
Suppose a signing law in this paper satisfies
\[
 S=(A\sigma,\sigma)\cx N(0,\Gamma),\qquad \Gamma\succ0.
\]
There is a coupling with $R\sim N(0,\Gamma)$ which retains every signing probability and satisfies
\[
 \E[R\mid\sigma]=S.
\]
In whitened coordinates it is the unique optimal martingale endpoint coupling. Conditional on each sign outcome of positive probability, $\Gamma^{-1/2}R$ is a proximal image of standard Gaussian measure. Consequently
\begin{align}
 \Ent(f(R)^2\mid\sigma)
   &\le2\E[\nabla f(R)^{\mathsf T}\Gamma\nabla f(R)\mid\sigma],\label{CM-eq:signLSI}\\
 W_{2,\Gamma^{-1}}^2(\eta,\law(R\mid\sigma))
   &\le2D(\eta\Vert\law(R\mid\sigma)).\label{CM-eq:signT2}
\end{align}
The whitened optimal path has the sharp endpoint and control-deficit estimates of Theorem~\ref{CM-thm:deficit} and Corollary~\ref{CM-cor:dynamic}.
\end{corollary}
\begin{proof}
Apply Theorem~\ref{CM-thm:contraction} with initial law $\law(\Gamma^{-1/2}S)$ and terminal law $\gamma_{m+n}$. The map $\sigma\mapsto S$ is injective, since its second block is $\sigma$. Thus the conditional-mean identity and the complete signing law are retained. Apply Corollary~\ref{CM-cor:proximal} and Theorem~\ref{CM-thm:conditional}, then undo the whitening. The cost $W_{2,\Gamma^{-1}}$ uses the squared norm $z^{\mathsf T}\Gamma^{-1}z$.
\end{proof}

Hard discrepancy, exact marginals and the Shannon entropy of the signing law are unchanged in Corollary~\ref{CM-cor:signing}. If $\Gamma$ is block diagonal, the reference blocks remain independent because their full joint law is the same Gaussian. The original reference--sign joint law, including its mutual information, can change.

For a cosine-product terminal reference, apply the curvature theorem in
the corresponding normalized metric. Its support is a closed box,
which every point of the martingale path retains by conditional
expectation. The conditional functional inequalities hold for the
variational optimizer in that metric as well.

\subsection{Dependence on the initial state}
\label{R22-sec:initial-kernels}
The canonical contraction bounds the dependence of each conditional map on its Gaussian driving variable. Dependence on the initial state has an additional obstruction, even when the terminal law is Gaussian and the endpoint pair has the positive product component used in the adapted-stability theorem.

\begin{proposition}[A Gaussian-terminal obstruction to initial-state contraction]
\label{R22-prop:initial-kernels}
There are finitely supported laws $\mu\cx\gamma_3$ admitting a martingale coupling with a positive product component, such that no martingale kernel $x\mapsto\pi_x$ from $\mu$ to $\gamma_3$ satisfies
\[
 W_1(\pi_x,\pi_{x'})\le |x-x'|
 \qquad(x,x'\in\operatorname{supp}\mu).
\]
\end{proposition}
\begin{proof}
Put $a=\sqrt{2/\pi}$ and let $\mu_1$ be uniform on $a\{-1,1\}^3$. The coupling $X=a\operatorname{sign}(Y)$, $Y\sim\gamma_3$, is a martingale. It is unique. Indeed, writing $\sigma=X/a$, every martingale coupling satisfies
\[
 \E[\sigma_iY_i]=a=\E|Y_i|.
\]
Equality in $\sigma_iY_i\le|Y_i|$ forces $\sigma_i=\operatorname{sign}(Y_i)$ for each coordinate.

Let $\pi_+,\pi_-$ be the conditional terminal laws at $x_+=a(1,1,1)$ and $x_-=-x_+$. If $W_1(\pi_+,\pi_-)=|x_+-x_-|$, equality in the mean-distance lower bound forces an optimal displacement to be a nonnegative multiple of $(1,1,1)$ almost surely. Thus every projection perpendicular to that vector would have the same law under $\pi_+$ and $\pi_-$. For the vector $(1,1,-2)$ these laws are those of $T=H_1+H_2-2H_3$ and $-T$, with independent standard half-normal variables $H_i$. They have mean zero and opposite nonzero third moments, since
\[
 \E T^3=-6\sqrt{2/\pi}\,(4/\pi-1)\ne0.
\]
Consequently $W_1(\pi_+,\pi_-)>|x_+-x_-|$.

For $0<\lambda<1$, put $\mu_\lambda=(\lambda\operatorname{id})_\#\mu_1$. The conditional laws $\lambda\pi_x+(1-\lambda)\gamma_3$ at $\lambda x$ have mean $\lambda x$ and average to $\gamma_3$. Their joint law dominates $(1-\lambda)\mu_\lambda\otimes\gamma_3$. Suppose kernels satisfying the asserted contraction existed along $\lambda_n\uparrow1$. Compactness and the fixed Gaussian terminal moments give a limiting martingale from $\mu_1$ to $\gamma_3$, which must be the unique orthant coupling. Since the eight separated initial atoms have fixed masses, each row converges in $W_1$: its tails are bounded by eight times the Gaussian tails. The row contraction inequalities would pass to the limit, contradicting the strict inequality above. Hence every $\lambda$ sufficiently close to one gives the required example.
\end{proof}
This proposition concerns Lipschitz dependence on the initial state. The driver contraction and the convergence of conditional laws in Corollary~\ref{R20-cor:intro-three-laws} keep their stated scopes. The example imposes no obstruction to existence of a Markov martingale, and its initial laws are discrete.

\subsection{Mixed-reference convex containment}
\label{CM-sec:mixed}
The covariance cap and fixed-reference extraction can also be used
with different reference measures in the same sum. The terminal law
sets the measure of the convex body; two Gaussian terms complete
the martingale noise.

\begin{theorem}[A log-concave reference and two Gaussian summands]
\label{CM-thm:mixed}
Let $\nu$ satisfy \eqref{CM-eq:curvature}. If Borel sets
$A_0,A_1,A_2\subset\R^d$ satisfy
\[
 \nu(A_0)+\gamma_d(A_1)+\gamma_d(A_2)>2,
\]
then $A_0+\frac L2A_1+\frac L2A_2$ contains a compact convex body
$K$ with $\nu(K)\ge1/e$. More precisely the mass can be taken at
least $\kappa(\nu)$ from Theorem~\ref{EXT-thm:fixed}.
\end{theorem}
\begin{proof}
Reduce to compact subsets while preserving the strict mass inequality.
For any bounded $X\cx\nu$, Theorem~\ref{CM-thm:contraction} gives
$M_0=X$, $M_1=Y\sim\nu$ and bracket density bounded by $L^2I_d$.
Apply Lemma~\ref{GC-lem:brown-complete} to $M/L$. Its two Brownian
completions give
\[
 X=Y+\frac L2(G_1+G_2),\qquad G_1,G_2\sim\gamma_d.
\]
The two Gaussians can be dependent. The union bound gives positive
probability that all three variables belong to their respective sets.
Thus every compactly supported law dominated by $\nu$ hits their
compact sum $S$.

Apply Theorem~\ref{EXT-thm:fixed} to $S+\varepsilon B_1$, obtaining
an open convex subset of $\nu$-mass at least $\kappa(\nu)$.
Its closure lies in $S+\varepsilon\overline B_1$. These compact convex
sets lie in a fixed bounded region. A Hausdorff limit as
$\varepsilon\downarrow0$ lies in $S$ and has mass at least
$\kappa(\nu)$, since every open neighborhood of the limit contains
all sufficiently late approximants. Absolute continuity gives
nonempty interior. The log-concave bound
$\kappa(\nu)\ge1/e$ follows from the same extraction theorem.
\end{proof}
The coefficient and mass of this mixed sum are stated without a
separate optimality claim. For a Gaussian reference the exact
half-measure extraction remains available. Whitening replaces the
coefficients $L I_d$ by $Q^{1/2}$ when the curvature matrix is $Q^{-1}$.

The canonical coupling now has a distance estimate and conditional
functional inequalities. We next establish the common-potential
representation for a prescribed non-Gaussian driver. The sharp
Gaussian deficit will then quantify how the optimal coupling changes
when that driver approximates the Gaussian.


\section{Arbitrary-reference dual attainment and exact \texorpdfstring{$q$}{q}-Bass existence}
\label{R13-sec:qbass}\label{QB-R11:section}
Tschiderer proves duality and primal uniqueness for nondegenerate references and poses globally finite $q$-Bass existence~\cite{QB-Tsch}; Acciaio--Marini prove one-dimensional first-moment existence~\cite{QB-AM}. Here the second-moment common-potential theorem holds in every dimension. Theorem~\ref{R24-thm:q-Bass-universal} gives the exact reference-support criterion for universal globally finite existence, including atomic references through subgradients.

Separate conditional Brenier maps need a common potential. A spanning martingale coupling controls their conditional Jensen gaps and normalizes one global minimizing sequence. Mean-exact recovery admits boundary and unbounded competitors, giving conditional dual attainment and one shifted subdifferential for every optimizer, including atomic references. The concluding corollaries keep the original signing law fixed: finite drivers approximate its canonical Gaussian coupling, while strict terminal inflation permits a common potential for every admissible driver.

\subsection{One potential for the entire optimal coupling}\label{QB-sec:results}
Let $\cP_2(\R^d)$ denote probability measures with finite second moment, and write
$m_2(p)=\int|y|^2p(dy)$. A coupling $\pi(dx,dy)=\mu(dx)\pi_x(dy)$ belongs to
$\mathcal M(\mu,\nu)$ if its marginals are $\mu,\nu$ and $\bary(\pi_x)=x$ almost surely.
The notation $\mu\cx\nu$ denotes convex order. For $p,q\in\cP_2(\R^d)$ put
\[
 \MCov(p,q)=\sup_{\lambda\in\Pi(p,q)}\int\ip{y}{w}\,\lambda(dy,dw),
 \qquad
 P_q(\mu,\nu)=\sup_{\pi\in\mathcal M(\mu,\nu)}
       \int\MCov(\pi_x,q)\,\mu(dx).
\]
A convex-order pair is \emph{irreducible} if for every Borel $A,B$ with
$\mu(A)>0$ and $\nu(B)>0$ some $\pi\in\mathcal M(\mu,\nu)$ satisfies
$\pi(A\times B)>0$. We work in the affine hull of the terminal support and
state the theorem when that hull is $\R^d$.

\begin{theorem}[General-reference existence]\label{QB-thm:main}
Let $\mu,\nu,q\in\cP_2(\R^d)$, let $\mu\cx\nu$ be irreducible, and assume
$\aff\supp\nu=\R^d$. Suppose $q$ gives no mass to any Borel set of Hausdorff
dimension at most $d-1$. There exist a proper lower-semicontinuous convex
function $v:\R^d\to(-\infty,+\infty]$ and a Borel map $a:\R^d\to\R^d$ such that
\begin{equation}\label{QB-eq:representation}
 \pi_x^q=\bigl(w\mapsto\nabla v(a(x)+w)\bigr)_\# q,
 \qquad
 \int\nabla v(a(x)+w)\,q(dw)=x
\end{equation}
for $\mu$-almost every $x$, where $\pi^q$ is the unique optimizer of $P_q$.
The gradient exists at $a(x)+w$ almost everywhere for this joint law and
has joint second moment $m_2(\nu)$.

Set $\alpha=a_\#\mu$ and
$g(z)=\int\nabla v(z+w)\,q(dw)$ on the $\alpha$-full set where this integral
exists. Then
\begin{equation}\label{QB-eq:Bass}
 g_\#\alpha=\mu,\qquad
 (\nabla v)_\#(\alpha*q)=\nu.
\end{equation}
The function $v$ is finite on all of $\R^d$ if either
$\overline{\conv}\supp q=\R^d$ or $\supp\nu$ is compact.
\end{theorem}

Absolute continuity of $q$ implies the stated Hausdorff-null condition.
Neither marginal needs a density; the terminal support may be unbounded.
The latent law $\alpha$ is a probability measure; its second moment
may be infinite. The square-integrable object is the terminal gradient,
whose law is exactly $\nu$.

\begin{theorem}[One subgradient representation for arbitrary references]\label{QB-thm:arbitrary}
Let $(\mu,\nu)$ be irreducible in $\cP_2(\R^d)$ with full terminal affine
span, and let $q\in\cP_2(\R^d)$ be arbitrary. There are a proper
lower-semicontinuous convex function $v$ and a Borel map $a$ such that
every optimal martingale coupling can be lifted to $(X,Y,W)$ with
\[
 (X,Y)\sim\pi,\qquad W\sim q,\qquad W\perp X,
 \qquad Y\in\partial v(a(X)+W)\quad\text{almost surely}.
\]
The potential and the shift can be chosen the same for all the optimal
couplings. The conditional relaxed dual is attained. The reference may
be atomic or concentrated on a lower-dimensional set.
\end{theorem}

The same potential thus describes every optimal conditional allocation.
Reference nondegeneracy makes each of its relevant subgradients a single
vector, producing the gradient representation of
Theorem~\ref{QB-thm:main}.

Tschiderer proved duality and primal uniqueness and conjectured
$q$-Bass existence for irreducible pairs under suitable reference
hypotheses~\cite{QB-Tsch}. Theorem~\ref{R24-thm:q-Bass-universal} below
identifies the exact universal hypothesis for his globally finite
potential formulation: full convex reference support. The theorem above
also gives extended-potential existence for every allowed support, and
global finiteness whenever the terminal support is compact.

Acciaio and Marini have proved the one-dimensional existence theorem for
$\mu,\nu\in\cP_1(\R)$ and every absolutely continuous reference, without a
moment assumption on the reference \cite[Theorem 1.1]{QB-AM}. Their
definition uses a monotone transport map on the relevant support. The
present proof gives arbitrary dimension and the quadratic variational
characterization; their one-dimensional integrability range is larger.
The semidiscrete precursor is \cite{QB-AMsemi}.

Every strict dilation of the terminal law makes an arbitrary
convex-order pair irreducible while retaining the initial law
(Corollary~\ref{QB-cor:inflation}). The construction will be used
below for the complete signing distribution.

The main limiting issue is the dual cost: its conditional differences can stay finite while its separate global integrals diverge. A spanning coupling bounds the Jensen gaps of one minimizing sequence and fixes its common affine normalization. The normalized limit controls compact interior competitors. To reach the full problem, each remaining competitor is contracted toward its prescribed mean and replaced by finite conditional expectations. Both operations keep that mean exact and recover its convex cost. Conditional optimality then determines the shift $a(x)$, and Fenchel equality yields the common potential. Section~\ref{QB-sec:dynamics} realizes the endpoints with a specified Brownian or L\'evy driver.

\subsection{Duality and the geometry of irreducibility}
Translation of $q$ changes $\MCov(p,q)$ by an affine function of
$\bary(p)$, so it changes $P_q$ by a constant and does not change its
optimizer. We henceforth center $q$; a translation of $a$ restores the
original reference at the end.

For a finite convex function $f$ of quadratic growth with a quadratic
lower bound, put
\begin{equation}\label{QB-eq:Phi}
 \Phi_f(x)=\inf_{\substack{p\in\cP_2(\R^d)\\\bary(p)=x}}
       \left\{\int f\,dp-\MCov(p,q)\right\},
 \qquad
 D_q(f)=\int f\,d\nu-\int\Phi_f\,d\mu.
\end{equation}
Tschiderer's duality theorem and convexification proposition give
\begin{equation}\label{QB-eq:known-duality}
 P_q(\mu,\nu)=\inf_f D_q(f),
\end{equation}
where the infimum may be taken over finite convex functions satisfying
$\ell+|y|^2/2\le f(y)\le A+K|y|^2$ for constants depending on $f$
\cite[Propositions 2.1 and 2.3]{QB-Tsch}. For nondegenerate $q$, the primal optimizer is unique.
We use this result, including its exact class of minimizing sequences,
as an input. Affine normalization preserves $D_q$ even when it changes
the constants in a particular quadratic lower bound.

\begin{proposition}[Duality before the nondegeneracy assumption]\label{QB-prop:all-duality}
For every centered $q\in\cP_2(\R^d)$, including atomic references, the
primal value is finite and attained, and \eqref{QB-eq:known-duality} holds
with the same convex quadratic class. Only the uniqueness assertion
requires the nondegeneracy assumption.
\end{proposition}
\begin{proof}
Use the weak cost
\[
 C(x,p)=\begin{cases}
 \tfrac12 W_2^2(p,q)-\tfrac12m_2(q),&\bary(p)=x,\\
 +\infty,&\bary(p)\ne x.
 \end{cases}
\]
It is bounded below, jointly lower semicontinuous on
$\R^d\times\cP_2(\R^d)$, and convex in $p$. The weak-transport existence
and duality theorem \cite[Theorems 1.2--1.3]{QB-BBP} therefore applies.
Subtracting the fixed terminal second moment gives exactly the duality
in \eqref{QB-eq:Phi}, initially for continuous quadratic functions. This
is the proof of \cite[Proposition 2.1]{QB-Tsch} with its strict-convexity
uniqueness step omitted.

Convexification also uses no nondegeneracy. If $f$ is in that quadratic
class and $h=\conv f$, choose measurably, for each $y$, a probability
$r_y$ with mean $y$ and
$\int f\,dr_y\le h(y)+\epsilon$. For a competitor $p$, its spread
$\bar p=\int r_y\,p(dy)$ has the same mean, belongs to $\cP_2$ by the
quadratic lower bound on $f$, and satisfies
$\int f\,d\bar p\le\int h\,dp+\epsilon$.
Also $\MCov(\bar p,q)\ge\MCov(p,q)$: couple $p$ optimally to $q$ and
then apply the mean-preserving kernel $r_y$ conditionally, preserving
its covariance with the reference. It follows that
$\Phi_f\le\Phi_h+\epsilon$, while $h\le f$ gives the reverse inequality.
Letting $\epsilon\downarrow0$ gives $\Phi_f=\Phi_h$ and
$D_q(h)\le D_q(f)$. This proves the convex restriction and the proposition.
\end{proof}

For clarity, the elementary properties used below are
\begin{equation}\label{QB-eq:basic}
 \Phi_f(x)\le f(x),\qquad
 D_q(f+\ell)=D_q(f),\qquad
 D_q((1-t)f+th)\le(1-t)D_q(f)+tD_q(h)
\end{equation}
for affine $\ell$. The first inequality uses $p=\delta_x$ and centered
$q$. The second uses the mean constraint. The third follows because
$\Phi_f$ is an infimum of functions linear in $f$.

\begin{lemma}[A common conditional support]\label{QB-lem:span}
For an irreducible pair there is $\bar\pi\in\mathcal M(\mu,\nu)$ such that
\[
 \supp\bar\pi_x=\supp\nu\quad\text{for }\mu\text{-almost every }x.
\]
Consequently, if $C=\overline{\conv}\supp\nu$ and $I=\operatorname{int}C$,
then $\mu(I)=1$ and $\overline{\conv}\supp\bar\pi_x=C$ almost surely.
\end{lemma}
\begin{proof}
Fix a Borel $B$ with $\nu(B)>0$. For $\pi\in\mathcal M(\mu,\nu)$ write
$A_\pi=\{x:\pi_x(B)>0\}$. The supremum of $\mu(A_\pi)$ is attained by
a countable convex mixture of a sequence approaching the supremum,
since the positive set of a mixture is the union of the positive sets.
If its complement $A$ had positive $\mu$-measure, irreducibility would
give another coupling with positive mass on $A\times B$. Mixing in
that coupling would increase the supremum, a contradiction. Thus there
is a coupling for which $\pi_x(B)>0$ almost surely.

Apply this construction to the countable family of rational open balls
having positive $\nu$-mass, and take another countable convex mixture.
For a common full-measure set of $x$, every such ball has positive
conditional mass. The conditional law is also concentrated on
$\supp\nu$ almost surely. This proves the support identity. A probability
with full convex support $C$ and finite first moment has barycenter in
$I$: otherwise a supporting hyperplane at that barycenter would force
its whole mass onto that hyperplane. The martingale identity gives
$\mu(I)=1$.
\end{proof}

The spanning-coupling formulation is consistent with the irreducible
geometry used in \cite{CM-Bass,DMT2019}. The direct proof above is enough
for the global irreducibility hypothesis in Theorem~\ref{QB-thm:main}.

\begin{lemma}[Continuity of maximal covariance]\label{QB-lem:MCov}
For $p,p',q,q'\in\cP_2(\R^d)$,
\begin{align*}
 |\MCov(p,q)-\MCov(p',q)|&\le\sqrt{m_2(q)}\,W_2(p,p'),\\
 |\MCov(p,q)-\MCov(p,q')|&\le\sqrt{m_2(p)}\,W_2(q,q').
\end{align*}
Moreover, $p\mapsto\MCov(p,q)$ is concave. It is strictly concave when
$q$ gives no mass to sets of Hausdorff dimension at most $d-1$.
\end{lemma}
\begin{proof}
Glue an optimal covariance coupling to an optimal $W_2$ coupling of the
changed marginal and use Cauchy--Schwarz; interchanging the two measures
gives the absolute-value bounds. Mixtures of covariance couplings prove
concavity. In an equality case for two distinct $p$'s, the mixture of
their optimal couplings from $q$ would again be optimal. Brenier's
theorem makes an optimal coupling from this $q$ a unique deterministic
gradient map. The two component maps must agree $q$-almost surely, which
forces the two $p$'s to coincide. This is also the strictness argument in
\cite[Proposition 2.1]{QB-Tsch}.
\end{proof}

\subsection{Compactness from the Jensen gap}
The next lemma isolates the part of convex dual compactness that uses
only the martingale geometry. Its argument follows the conditional
convex-function compactness mechanism in \cite[Section 7]{CM-Bass}.
We include it because the absence of any Gaussian hypothesis is essential.

\begin{lemma}[Interior compactness]\label{QB-lem:compactness}
Let $q\in\cP_2$ be centered, with no nondegeneracy assumption.
Let $(f_n)$ be a finite convex minimizing sequence for
\eqref{QB-eq:known-duality}. There are convex averages, a subsequence,
and affine functions which may be subtracted, after which
\[
 f_n\ge0,\qquad f_n(x_0)=0\quad\text{for one }x_0\in I,
\]
and $f_n$ converges uniformly on compact subsets of $I$ to a finite convex
function $\psi_I$. There is a proper lower-semicontinuous convex extension
$\psi$ with $\psi\ge0$, $\psi|_I=\psi_I$, and $\dom\psi\subset C$, such that
\begin{equation}\label{QB-eq:boundary-liminf}
 \psi(y)\le\liminf_n f_n(y)\qquad(y\in C).
\end{equation}
The modified sequence still satisfies $D_q(f_n)\to P_q$.
\end{lemma}
\begin{proof}
Let $\bar\pi$ be given by Lemma~\ref{QB-lem:span}. Its conditional Jensen
gaps are
\[
 J_n(x)=\int f_n(y)\,\bar\pi_x(dy)-f_n(x)\ge0.
\]
They satisfy
\[
 \int J_n\,d\mu
 =\int f_n\,d\nu-\int f_n\,d\mu
 \le D_q(f_n).
\]
Thus $(J_n)$ is bounded in $L^1(\mu)$. The Koml\'os subsequence theorem
\cite{QB-Komlos} gives a subsequence whose Ces\`aro means converge to
a finite value almost surely. Replace $f_n$ by the corresponding means.
Convexity of $D_q$ and weak duality show that the new sequence is still
minimizing. For almost every $x$, $\sup_n J_n(x)<\infty$.

Choose one such $x_0\in I$ with conditional convex support $C$, and choose
$s_n\in\partial f_n(x_0)$. Subtract
$f_n(x_0)+\ip{s_n}{y-x_0}$. Then $f_n\ge0$, $f_n(x_0)=0$, and
\[
 \sup_n\int f_n(y)\,\bar\pi_{x_0}(dy)=\sup_n J_n(x_0)<\infty.
\]
For a compact $K\subset I$, put
\[
 \delta_K=\inf_{|e|=1}\bar\pi_{x_0}
  \{y:\ip{e}{y}>\sup_{z\in K}\ip{e}{z}\}.
\]
This number is positive. Each direction has positive mass beyond the
supporting level of $K$, since $K$ is inside the full convex support.
The probability of this strict halfspace is lower semicontinuous in
$e$, so compactness of the sphere makes the infimum positive.

For $z\in K$ take a subgradient of $f_n$ at $z$. If it is nonzero, the
supporting affine function gives $f_n(y)\ge f_n(z)$ on the corresponding
halfspace in the preceding display. If it is zero, $f_n(z)=0$ because
$f_n(x_0)=0$. Hence
\[
 \sup_{z\in K}f_n(z)\le
 \delta_K^{-1}\int f_n\,d\bar\pi_{x_0}.
\]
Local bounds for convex functions give local equicontinuity on the
interior. A diagonal Arzel\`a--Ascoli argument proves local uniform
convergence on $I$.

Take the lower-semicontinuous convex extension of $\psi_I$ to $C$, with
value $+\infty$ outside $C$. For $y\in C$ and $0<t<1$, convexity gives $f_n(y)\ge t^{-1} f_n((1-t)x_0+ty)$. First let $n\to\infty$, then $t\uparrow1$. The radial limit is exactly the
lower-semicontinuous extension $\psi(y)$. This proves
\eqref{QB-eq:boundary-liminf}, including points at which the limit is
infinite. Affine invariance preserves the minimizing property.
\end{proof}

The use of the integral Jensen gap is important. It bounds the full
convex function after one affine normalization, although the dual
minimizers need not have an integrable global limit under $\nu$.
Conditional integrability of that limit will follow from the optimality
gap, before any conjugate is used.

\subsection{Conditional dual attainment}
The limiting dual must compare against laws on the boundary of the terminal convex hull and against unbounded laws, while retaining the prescribed mean exactly. Contracting a competitor toward its own mean moves its support into the interior. Conditional expectations on finite partitions then discretize it without changing that mean. Convexity bounds the dual integrals throughout, and lower semicontinuity identifies their limit. These two operations preserve the constraint needed to test conditional optimality.

\begin{lemma}[Recovery by compact interior competitors]\label{QB-lem:recovery}
Let $C$ be closed convex with nonempty interior $I$, let $x\in I$, and
let $\psi\ge0$ be proper lower-semicontinuous convex, finite on $I$, with
$\dom\psi\subset C$. If $p\in\cP_2(\R^d)$ has mean $x$ and
$\int\psi\,dp<\infty$, there are finitely supported probability measures
$p_j$ in $I$ with mean $x$ such that
\[
 W_2(p_j,p)\longrightarrow0,\qquad
 \int\psi\,dp_j\longrightarrow\int\psi\,dp.
\]
Consequently $\MCov(p_j,q)\to\MCov(p,q)$ for every $q\in\cP_2$.
\end{lemma}
\begin{proof}
For $Y\sim p$ and $0<\lambda<1$ set
$Y_\lambda=x+\lambda(Y-x)$. Its mean is $x$ and it belongs to $I$ almost
surely. Convexity gives $0\le\psi(Y_\lambda)\le(1-\lambda)\psi(x)+\lambda\psi(Y)$. As $\lambda\uparrow1$, radial lower semicontinuity and convexity give
$\psi(Y_\lambda)\to\psi(Y)$, while $Y_\lambda\to Y$ in $L^2$.
Dominated convergence gives convergence of the $\psi$ integrals.

For fixed $\lambda$, choose increasing finite measurable partitions
generating $\sigma(Y_\lambda)$ and let $Y_{\lambda,j}$ be the corresponding
conditional expectations. These have finite support, mean $x$, and
converge to $Y_\lambda$ in $L^2$ and almost surely. Each of their values
is in $I$. Indeed a barycenter of a probability carried by an open
convex set with finite first moment belongs to that set, by the
supporting-hyperplane argument. Conditional Jensen gives
$\E\psi(Y_{\lambda,j})\le\E\psi(Y_\lambda)$. Continuity of $\psi$ in $I$
and Fatou's lemma give the reverse limiting inequality. Diagonalization
in $\lambda,j$ proves the claim. The covariance assertion follows from
Lemma~\ref{QB-lem:MCov}.
\end{proof}

The same normalized minimizing sequence is used for every primal optimizer. Its conditional optimality gaps are nonnegative and their integral tends to zero. Their almost-sure vanishing, together with the recovery lemma, gives conditional attainment for any chosen optimizer. This is why the eventual potential and shift can be common to all optimizers even when the driver has atoms.

\begin{theorem}[Conditional attainment]\label{QB-thm:attainment}
Let $q\in\cP_2$ be centered and let $\pi^q$ be any primal optimizer.
The function $\psi$ obtained
in Lemma~\ref{QB-lem:compactness} satisfies, for $\mu$-almost every $x$,
\begin{equation}\label{QB-eq:local-attain}
 \int\psi\,d\pi_x^q<\infty,\qquad
 \Phi_\psi(x)=\int\psi\,d\pi_x^q-\MCov(\pi_x^q,q)\in\R.
\end{equation}
This assertion needs no nondegeneracy of $q$. Here $\Phi_\psi$ is the infimum in \eqref{QB-eq:Phi}, with competitors
of infinite $\psi$ integral assigned value $+\infty$.
Moreover $\Phi_\psi$ is finite and continuous on $I$, is convex, and has
a nonempty subdifferential at every point of $I$.
\end{theorem}
\begin{proof}
For the normalized minimizing sequence set
\[
 d_n(x)=\int f_n\,d\pi_x^q-\Phi_{f_n}(x)-\MCov(\pi_x^q,q)\ge0.
\]
Its integral is $D_q(f_n)-P_q$, which tends to zero. Pass to a further
subsequence so $d_n(x)\to0$ almost surely. Fix a good $x\in I$, and write
$p_* =\pi_x^q$. Since $\Phi_{f_n}(x)\le f_n(x)$,
\[
 \int f_n\,dp_*\le f_n(x)+\MCov(p_*,q)+d_n(x).
\]
The right side is bounded. The nonnegativity of $f_n$ and
\eqref{QB-eq:boundary-liminf} imply $\int\psi\,dp_*<\infty$.

If $p$ is any finitely supported law in $I$ with mean $x$, then
\[
 \int f_n\,dp_*-\MCov(p_*,q)
 \le\int f_n\,dp-\MCov(p,q)+d_n(x).
\]
Fatou's lemma on the left and uniform convergence on the finite support
on the right give
\begin{equation}\label{QB-eq:competitor}
 \int\psi\,dp_*-\MCov(p_*,q)
 \le\int\psi\,dp-\MCov(p,q).
\end{equation}
Lemma~\ref{QB-lem:recovery} extends this inequality to every finite-cost
competitor with mean $x$. Infinite-cost competitors do not lower the
infimum. This proves \eqref{QB-eq:local-attain}.

Concavity of $p\mapsto\MCov(p,q)$ shows that $\Phi_\psi$ is convex.
It is at most $\psi$ on $I$. It cannot be $-\infty$ at any point $y$ of
$C$: a point $x$ where \eqref{QB-eq:local-attain} holds lies on the open
segment from $y$ to some $z\in I$; convexity and
$\Phi_\psi(z)\le\psi(z)<\infty$ would then force
$\Phi_\psi(x)=-\infty$. Outside $C$, finite-cost competitors do not exist.
Thus $\Phi_\psi$ is proper and finite on $I$. A finite convex function
on an open convex set is continuous and has supporting subgradients
there. Its supporting inequalities hold on the entire domain by
convexity, including its boundary.
\end{proof}

\begin{remark}[The integrability convention]\label{QB-rem:relaxed}
The theorem asserts a finite conditional integral for almost every $x$.
The global integral $\int\psi\,d\nu$ can be infinite. The relaxed dual
identity is the well-defined conditional difference
\[
 \int\left[\int\psi(y)\,\pi_x^q(dy)-\Phi_\psi(x)\right]\mu(dx)=P_q.
\]
No subtraction of two infinite global integrals is used. This is the
same conditional interpretation required by the relaxed duality
formulation in \cite[Proposition 2.8 and Remark 1.6]{QB-Tsch}.
\end{remark}

\subsection{Recovering one translated-reference potential}
\begin{lemma}[Subgradient recovery]\label{QB-lem:subgradient}
Under Theorem~\ref{QB-thm:attainment}, put $v=\psi^*$. There is a Borel
selection $a(x)\in\partial\Phi_\psi(x)$ on $I$ such that, for every good
$x$ and every covariance-maximizing coupling $(W,Y)$ of $q,\pi_x^q$,
\begin{equation}\label{QB-eq:Fenchel}
 \psi(Y)+v(a(x)+W)=\ip{a(x)+W}{Y}\quad\text{almost surely}.
\end{equation}
If $q$ gives no mass to sets of Hausdorff dimension at most $d-1$,
then $Y=\nabla v(a(x)+W)$ almost surely.
\end{lemma}
\begin{proof}
The subdifferential of the finite convex function $\Phi_\psi$ on $I$ has
nonempty closed convex values and a Borel graph. Its minimum-norm element
is a Borel selection. Fix a good $x$ and write $a=a(x)$. For any
square-integrable coupling $(W,Y')$ with first marginal $q$ and finite
$\psi$ integral, let $x'=\E Y'$. Subgradient optimality and the definition
of $\Phi_\psi$ give
\begin{align*}
 \E\{\psi(Y')-\ip{W+a}{Y'}\}
 &\ge \Phi_\psi(x')-\ip{a}{x'}\\
 &\ge \Phi_\psi(x)-\ip{a}{x}
 =\E\{\psi(Y)-\ip{W+a}{Y}\}.
\end{align*}
All terms on the last line are finite by conditional attainment and
Cauchy--Schwarz.

Replace $Y$ by a fixed $b\in I$ on an arbitrary measurable event and
leave it unchanged elsewhere. This is an admissible square-integrable
competitor with finite $\psi$ cost. The preceding inequality therefore
implies
\[
 \psi(Y)-\ip{W+a}{Y}\le \psi(b)-\ip{W+a}{b}\quad\text{almost surely}.
\]
Take the common full-measure event for a countable dense set of $b$ in
$I$. Continuity on $I$ and radial approximation of the boundary imply
\[
 \sup_{b\in I}\{\ip{W+a}{b}-\psi(b)\}
 =\sup_{b\in\R^d}\{\ip{W+a}{b}-\psi(b)\}=v(W+a).
\]
The defining Fenchel inequality gives the reverse inequality at $Y$,
so equality holds. Equivalently, $Y\in\partial v(W+a)$.

Under this additional hypothesis on $q$, the effective domain of $v$
has nonempty interior. Otherwise its affine
hull would have dimension at most $d-1$, contradicting the fact that
$W+a\in\dom v$ almost surely. The boundary of a convex domain and the
nondifferentiability set of a finite convex function in its interior
are contained in Borel sets of Hausdorff dimension at most $d-1$.
The assumption on $q$ therefore makes $v$ differentiable at $W+a$
almost surely. Its subgradient there is the single vector $\nabla v(W+a)$.
\end{proof}

\begin{proof}[Proof of Theorem~\ref{QB-thm:arbitrary}]
Proposition~\ref{QB-prop:all-duality} gives primal attainment and a
convex minimizing sequence. Lemma~\ref{QB-lem:compactness} uses only the
Jensen gap, so it gives one limit $\psi$. For any fixed optimal $\pi$,
the nonnegative conditional gaps in the proof of
Theorem~\ref{QB-thm:attainment} have integrals tending to zero. A further
subsequence makes them tend to zero almost surely without changing the
local uniform limit $\psi$. The proof therefore gives local attainment
for this $\pi$ by the same $\psi$. This holds for each optimizer.

Set $v=\psi^*$ and choose the same Borel subgradient
$a(x)\in\partial\Phi_\psi(x)$. Choose measurably an optimal covariance
coupling of $q$ with each $\pi_x$. Such a selection exists by the usual
compactness and Borel graph of the optimal-coupling correspondence in
$\cP_2$. Its joint law has $W$ independent of $X$. The proof of
Lemma~\ref{QB-lem:subgradient}, up to and including Fenchel equality,
uses no regularity of $q$ and gives the subgradient inclusion. The final
step of that lemma is precisely where nondegeneracy makes the inclusion
a deterministic gradient identity. Centering and then restoring $q$
only translates $a$.
\end{proof}

\begin{proof}[Proof of Theorem~\ref{QB-thm:main}]
Apply Lemmas~\ref{QB-lem:compactness} and \ref{QB-lem:subgradient} and
Theorem~\ref{QB-thm:attainment}. The covariance coupling from $q$ to each
$\pi_x^q$ is the unique Brenier coupling. Its gradient agrees with
$w\mapsto\nabla v(a(x)+w)$ by \eqref{QB-eq:Fenchel}. The mean constraint
and the terminal marginal give \eqref{QB-eq:representation}. Measurability
of the selected $a$ and of a Borel version of $\nabla v$ makes the
representation joint. Its squared norm integrates to $m_2(\nu)$.

For $\alpha=a_\#\mu$, the identity
$g(a(x))=x$ holds almost surely, so $g_\#\alpha=\mu$ and
$(\nabla v)_\#(\alpha*q)=\nu$. Values of $g$ or the gradient on the
irrelevant null set can be chosen arbitrarily. The first identity also
shows that $a$ has a measurable inverse $g$ on full-measure subsets of
its source and image; in particular it does not merge distinct relevant
initial states.

If $\overline{\conv}\supp q=\R^d$, then $\dom v$ contains a translate of
a $q$-full set and its closed convex hull is $\R^d$. A convex set dense
in $\R^d$ equals $\R^d$ when it is an effective convex domain with full
interior; equivalently, the interior of a convex set equals the interior
of its closure. Thus $\dom v=\R^d$. If $\supp\nu$ is compact, then
$C$ is compact and $\psi\ge0$, so
\[
 v(z)\le\sup_{y\in C}\ip{z}{y}<\infty.
\]
It is also bounded below by one finite affine function because $\psi$
is finite at $x_0$. Thus again $v$ is finite everywhere. Translating the
reference back to its original mean completes the proof.
\end{proof}

\begin{remark}[The averaged conjugate]
For the selected $a(x)$, \eqref{QB-eq:Fenchel} implies
\[
 u(a(x)):=\int v(a(x)+w)\,q(dw)
 =\ip{a(x)}{x}-\Phi_\psi(x)\in\R.
\]
This integral is absolutely defined: its Fenchel representation contains
$\psi(Y)\in L^1$ and $\ip{W+a}{Y}\in L^1$. The theorem needs neither
finiteness of $u$ at every latent point nor a global interchange of
convolution and differentiation. The identity
$g(a(x))=\int\nabla v(a(x)+w)q(dw)=x$ is proved directly.
\end{remark}

\subsection{The potential domain and uniqueness for positive references}

\begin{theorem}[Universal finite potentials and the $q$-Bass existence criterion]
\label{R24-thm:q-Bass-universal}
Let $q\in\cP_2(\R^d)$ be arbitrary. The following are equivalent:
\begin{enumerate}[label=\textup{(\roman*)},leftmargin=*]
\item $\overline{\conv}\supp q=\R^d$.
\item For every irreducible pair $\mu\cx\nu$ in $\cP_2(\R^d)$ with full
terminal affine span, the common potential representing all optimizers
in Theorem~\ref{QB-thm:arbitrary} can be chosen finite everywhere.
\item There exist a globally finite convex $v$ and a joint law of
$(Z,W,Y)$ with $W\sim q$, $Y\sim\gamma_d$, $W\perp Z$,
$\E[Y\mid Z]=0$ and $Y\in\partial v(Z+W)$ almost surely.
\end{enumerate}
The latent law in \textup{(iii)} is unrestricted. If $q$ gives no mass
to sets of Hausdorff dimension at most $d-1$, these are also equivalent
to globally finite $q$-Bass existence for every such endpoint pair;
its coupling is then the unique maximal-covariance optimizer.
Thus the single Gaussian terminal test detects universal existence,
including for atomic references and randomized subgradient allocations.
\end{theorem}
\begin{proof}
Under \textup{(i)}, Theorem~\ref{QB-thm:arbitrary} gives a common potential
$v$. For one admissible initial point $x$, its domain contains
$a(x)+w$ for $q$-almost every $w$. This convex domain is dense in
$\R^d$ and therefore equals $\R^d$, proving \textup{(ii)}.
The pair $(\delta_0,\gamma_d)$ is irreducible, so \textup{(ii)} gives
\textup{(iii)} with a constant latent variable. Nondegeneracy turns
subgradients into gradients by Theorem~\ref{QB-thm:main}.

For the converse, disintegrate the law in \textup{(iii)} over $Z=z$.
For almost every $z$ the conditional law $\eta_z(dw,dy)$ has first
marginal $q$, finite second moment in $y$, and mean zero in $y$.
For two such $z,z'$, couple $\eta_z,\eta_{z'}$ over the same $w$.
The two supporting-plane gaps have sum
$\langle y-y',z-z'\rangle\ge0$ and expected sum zero. Hence each gap
vanishes, and $y\in\partial v(z'+w)$ for $\eta_z$-almost every $(w,y)$.
Fubini allows one fixed $z_0$ for almost every $z$; averaging over $z$
then gives a coupling of $q,\gamma_d$ supported on
\[
 y\in\partial f(w),\qquad f(w)=v(z_0+w).
\]
This removes all latent mixing, including randomized allocations at
atoms of $q$.

Put $h=f^*$. Fenchel equality makes $h$ finite at
$\gamma_d$-almost every point. Its convex domain has full Lebesgue
measure and hence equals $\R^d$. Thus $h$ is locally Lipschitz and
differentiable almost everywhere. The inverse subgradient relation gives
$(\nabla h)_\#\gamma_d=q$.
If \textup{(i)} fails, separation gives a unit vector $e$ and $c\in\R$
with $\langle e,w\rangle\ge c$ on $\supp q$, so
$\langle e,\nabla h(y)\rangle\ge c$ Lebesgue-almost everywhere.
Absolute continuity on almost every line parallel to $e$, followed by
continuity, yields $h(y-te)\le h(y)-ct$ for all $y$ and $t\ge0$.
For $\langle e,z\rangle<c$, biconjugacy gives
\[
 f(z)\ge\langle z,y-te\rangle-h(y-te)
 \ge\langle z,y\rangle-h(y)+t(c-\langle e,z\rangle)
 \longrightarrow+\infty,
\]
contradicting finiteness of $f$.
\end{proof}

For absolutely continuous second-moment references, the theorem gives
the exact universal support criterion for Tschiderer's globally finite
$q$-Bass formulation
\cite[Definition~1.3 and the conjecture after Definition~1.7]{QB-Tsch}.
The common subgradient criterion covers every second-moment reference.
Extended potentials retain existence for arbitrary support, while
compact terminal support ensures global finiteness for an individual pair.

\begin{proposition}[Compact driving support and a one-point initial law]
\label{R13-prop:finite-domain}
Let $q\in\cP_2(\R^d)$ have compact support and satisfy the reference
hypothesis of Theorem~\ref{QB-thm:main}. Let $\nu\in\cP_2(\R^d)$
have full affine span and mean $m$. A $q$-Bass representation of
$(\delta_m,\nu)$ by a globally finite convex potential exists if and
only if $\nu$ has compact support. Necessity holds even when its latent
law is allowed to be any probability measure.
\end{proposition}
\begin{proof}
Compact terminal support gives a globally finite potential by
Theorem~\ref{QB-thm:main}, since $(\delta_m,\nu)$ is irreducible.
Conversely, the supporting-plane argument in the proof of
Theorem~\ref{R24-thm:q-Bass-universal} reduces any latent representation to a coupling on
$\partial v(a_0+w)$; the common mean $m$ makes the gap sum zero.
Reference nondegeneracy makes this coupling the gradient law
$(\nabla v(a_0+\cdot))_\#q=\nu$. Subgradients of a globally finite
convex function are bounded on the compact set $a_0+\supp q$.
Hence this common law has compact support.
\end{proof}

\begin{example}[Why a finite potential cannot always be required]\label{QB-ex:domain}
Let $d=1$, $q=\operatorname{Unif}[-1,1]$, $\mu=\delta_0$, and
$\nu=N(0,1)$. The pair is irreducible and all three measures have every
finite moment. A $q$-Bass representation is obtained from
\[
 T(w)=\Phi^{-1}\left(\frac{w+1}{2}\right),\quad -1<w<1,
 \qquad \alpha=\delta_0,
\]
where $\Phi$ is the standard normal distribution function. Its convex
primitive, extended lower semicontinuously to $[-1,1]$ and by $+\infty$
outside, is an admissible $v$ in Theorem~\ref{QB-thm:main}.

Proposition~\ref{R13-prop:finite-domain} shows that a globally finite
potential is impossible even with a different latent law. The quantile
map above has an integrable derivative primitive, but its slope diverges
at the endpoints of the reference interval.

Thus absolute continuity and arbitrarily strong moment assumptions on
$q$ alone cannot force a globally finite potential. This is a domain
issue in the original definition, resolved by the map-on-support or
extended-potential convention used here and consistent with the
one-dimensional formulation in \cite[Definition 2.4]{QB-AM}.
\end{example}

\begin{proposition}[Irreducibility is the exact criterion for positive references]\label{QB-prop:positive}
Let $\mu\cx\nu$ in $\cP_2(\R^d)$, with full terminal affine span,
and let $q\in\cP_2(\R^d)$ be equivalent to Lebesgue measure.
A common globally finite convex-gradient representation exists if and
only if $(\mu,\nu)$ is irreducible. In this case the optimizer in
Theorem~\ref{QB-thm:main} satisfies
\[
 \pi_x^q\sim\nu\quad\text{for }\mu\text{-almost every }x,
 \qquad \pi^q\sim\mu\otimes\nu.
\]
Among representations of the unique optimizer, the latent map $a$
is unique up to one common translation and the potential up to that
translation and an additive constant.
\end{proposition}
\begin{proof}
Irreducibility gives existence and global finiteness by
Theorem~\ref{QB-thm:main}. Every translate of $q$ is equivalent to Lebesgue measure. The image of
such a translate by $\nabla v$ therefore has the same null sets for every
latent point at which it is used. Fixing one good point and integrating
the conditional laws shows that their common measure class is that of
$\nu$. This also gives the joint equivalence. Any common-gradient
representation with a positive reference has this property, which
implies positive mass on every rectangle of positive marginal mass.

For uniqueness, let $(v_1,a_1)$ and $(v_2,a_2)$ represent the same optimal
kernel. Conditional Brenier uniqueness gives
\[
 \nabla v_1(a_1(x)+w)=\nabla v_2(a_2(x)+w)
\]
for almost every $w$ in Lebesgue measure and almost every $x$. At one
fixed good $x_0$, this identifies $v_2$ with a translate of $v_1$, up to
an additive constant. At another $x$, the difference between the two
latent shifts would consequently be a period $h$ of $\nabla v_1$.

A nonzero period forces the whole gradient image into a hyperplane.
To see this without smoothness, take points $z,z'$ of differentiability
such that the gradient repeats along all integer translates of $z$ by
$h$. Monotonicity gives
\[
 \ip{\nabla v_1(z)-\nabla v_1(z')}{z+nh-z'}\ge0
 \qquad(n\in\mathbb Z).
\]
Hence $\ip{h}{\nabla v_1(z)-\nabla v_1(z')}=0$. This holds on a
Lebesgue-full set, so the terminal law would have constant projection
onto $h$, contrary to full affine span. Every such period is zero.
Thus the latent shift is common to almost every $x$, as asserted.
\end{proof}

\begin{remark}
The uniqueness statement concerns representations of the variational
optimizer. It does not assert a moment bound for the latent measure or
uniqueness among arbitrary formal fixed points lacking the integrability
and optimality conditions. Conditional-kernel uniqueness itself holds
under the weaker Hausdorff-null hypothesis of Theorem~\ref{QB-thm:main}.
\end{remark}

\subsection{Changing the reference preserves the full conditional kernel}
The next result concerns the optimizer, so its fixed-marginal stability
does not require irreducibility. It permits discrete approximations to
an absolutely continuous reference.

\begin{theorem}[Reference stability]\label{QB-thm:stability}
Fix $\mu,\nu\in\cP_2(\R^d)$ in convex order. Let $q_n\to q$ in $W_2$,
where $q$ gives no mass to sets of Hausdorff dimension at most $d-1$.
Let $\varepsilon_n\downarrow0$ and choose
$\pi^n\in\mathcal M(\mu,\nu)$ with
\[
 \int\MCov(\pi_x^n,q_n)\,\mu(dx)\ge P_{q_n}(\mu,\nu)-\varepsilon_n.
\]
Then, writing $\pi^q$ for the unique optimizer with reference $q$,
\begin{equation}\label{QB-eq:kernel-stability}
 \int W_2^2(\pi_x^n,\pi_x^q)\,\mu(dx)\longrightarrow0.
\end{equation}
The values satisfy the explicit estimate
\begin{equation}\label{QB-eq:value-stability}
 |P_{q_n}(\mu,\nu)-P_q(\mu,\nu)|
 \le\sqrt{m_2(\nu)}\,W_2(q_n,q).
\end{equation}
No regularity or nondegeneracy is imposed on the approximating $q_n$.
More precisely, write $m=\int x\,\mu(dx)$, $b_n=\int w\,q_n(dw)$,
$b=\int w\,q(dw)$, and let $q_n^0,q^0$ be the centered references.
Then
\begin{equation}\label{R13-eq:increment-value}
 \left|P_{q_n}(\mu,\nu)-P_q(\mu,\nu)-\langle m,b_n-b\rangle\right|
 \le \sqrt{m_2(\nu)-m_2(\mu)}\,W_2(q_n^0,q^0).
\end{equation}
The coefficient in this centered estimate is sharp.
\end{theorem}
\begin{proof}
Lemma~\ref{QB-lem:MCov}, integrated against $\mu$, gives
\eqref{QB-eq:value-stability} by Cauchy--Schwarz. For the stronger
estimate center both references and replace $\pi_x$ by the law of
$Y-x$. Its second moment is $m_2(\pi_x)-|x|^2$, whose $\mu$-integral
is $m_2(\nu)-m_2(\mu)$. The same Cauchy--Schwarz argument gives
\eqref{R13-eq:increment-value}; restoring the reference means adds
$\langle m,b_n-b\rangle$. For sharpness take $\mu=\delta_0$,
a centered nonzero $\nu$, $q=\nu$, and $q_n=(c\id)_\#\nu$ with $c>0$.
Both sides then equal $|c-1|m_2(\nu)$. In particular the original
$\pi^n$ are asymptotically optimal for the fixed reference $q$.
Apply Lemma~\ref{AS-lem:lift} to $\mu(dx)\delta_{\pi_x^n}(dp)$, whose averaged terminal law is fixed at $\nu$. The value estimate makes these kernels asymptotically optimal for $q$. Concavity and the fixed marginals identify the limiting optimum; strict concavity from Lemma~\ref{CM-lem:strict} forces its conditional-kernel graph. The same-initial-state conclusion of Lemma~\ref{AS-lem:lift} is precisely~\eqref{QB-eq:kernel-stability}.
\end{proof}

\begin{remark}
The convergence in \eqref{QB-eq:kernel-stability} uses the same initial
state $x$ on both sides. It is stronger than weak convergence of the
joint couplings. For general nondegenerate references this gives qualitative kernel
convergence. The next theorem obtains a rate at the Gaussian reference
under terminal curvature. Stability of the latent distributions has
additional normalization and tightness requirements; compare the
one-dimensional fixed-point theorem in \cite[Theorem 1.3]{QB-AM}.
\end{remark}

\subsubsection{Uniform conditional contraction and reference sensitivity}
The fixed-terminal approximation separates a conditional regularity statement from the existence of one global potential. Curvature of a finite mixture gives a bound independent of its number of shifts; Lemma~\ref{R15-lem:transfer} retains that bound for arbitrary initial laws. The common-potential theorem still has its own irreducibility hypothesis.

\begin{theorem}[Unrestricted curvature-controlled contraction and deficit]
\label{R14-thm:reference-curvature}
Let $\mu\cx\nu$ in $\cP_2(\R^d)$ and suppose $\nu$ satisfies~\eqref{CM-eq:curvature}. Let $q_0\in\cP_2$ be centered with positive density proportional to $e^{-U}$, where $U\in C^2(\R^d)$ and $D^2U\preceq KI_d$, $K>0$. Put $M=L\sqrt K$. Every conditional optimal map $T_x^{q_0}$ has an $M$-Lipschitz convex-gradient representative, and every feasible triple satisfies
\begin{equation}\label{R14-eq:baseline-deficit}
 P_{q_0}-\E\langle Y,W_0\rangle
 \ge\frac1{2M}\E|Y-T_X^{q_0}(W_0)|^2,
 \quad W_0\sim q_0,\quad W_0\perp X.
\end{equation}
The coefficient is sharp. No finite-support or irreducibility condition on the initial pair is required. When the pair is irreducible, its common potential from Theorem~\ref{QB-thm:main} has an $M$-Lipschitz gradient globally.

For arbitrary $q\in\cP_2$ and $J_q(\pi)\ge P_q-\eta$, there is a same-initial-state coupling, with $W_0\perp X$ and $\law(X,Y)=\pi$, such that, for $\delta=W_2(q^0,q_0)$,
\begin{equation}\label{R14-eq:reference-curvature-rate}
 \left(\int W_2^2(\pi_x,\pi_x^{q_0})\,\mu(dx)\right)^{1/2}
 \le\|Y-T_X^{q_0}(W_0)\|_2
 \le M\delta+\sqrt{M^2\delta^2+2M\eta}.
\end{equation}
In particular, exact optimizers have conditional-kernel error at most $2M\delta$. Both endpoint laws are exact, and $q$ may be atomic.
\end{theorem}
\begin{proof}
For a finite approximating pair from Theorem~\ref{CM-thm:stability}, Theorem~\ref{QB-thm:main} gives a common potential transporting
$r(z)=\sum_jp_jq_0(z-a_j)$ to $\nu$. With $\theta_j(z)=p_jq_0(z-a_j)/r(z)$,
\begin{equation}\label{R14-eq:mixture-curvature}
 D^2(-\log r)=\sum_j\theta_j D^2U(\,\cdot-a_j)
 -\Cov_\theta(\nabla U(\,\cdot-a_j))\preceq KI_d.
\end{equation}
Caffarelli's estimate~\cite{Caffarelli,Kolesnikov} gives the common bound $M$. Convex approximation handles extended terminal potentials as in Lemma~\ref{CM-lem:mixture}. Lemma~\ref{R15-lem:transfer} gives the conditional contraction and deficit for the original pair. In the irreducible case, the common potential agrees with a Lipschitz conditional map on any full-measure fiber. Positivity of $q_0$ then gives the global gradient representative. For $\mu=\delta_0$, $q_0=N(0,K^{-1}I_d)$ and $\nu=N(0,L^2I_d)$, $T(w)=Mw$ and expansion of the square gives equality in~\eqref{R14-eq:baseline-deficit} for every feasible coupling.

For the rate, center $q$, couple $W\sim q^0$ optimally to $W_0$, independently of $X$, and glue this with the conditional optimal covariance couplings of $W$ to $Y\sim\pi_x$. Put $Y^*=T_X^{q_0}(W_0)$ and $r_*=\|Y-Y^*\|_2$. Feasibility of $(X,Y^*,W)$ gives $\E\langle Y^*-Y,W\rangle\le\eta$. Thus
\[
 \frac{r_*^2}{2M}\le\E\langle Y^*-Y,W_0\rangle
 \le\eta+r_*\delta.
\]
Solving the quadratic proves~\eqref{R14-eq:reference-curvature-rate}; conditioning the same coupling on $X$ gives its kernel bound.
\end{proof}

For a centered $d$-dimensional Student reference with $r>2$ degrees of freedom,
\[
 D^2U(w)=\frac{d+r}{r+|w|^2}I_d-
 \frac{2(d+r)}{(r+|w|^2)^2}ww^{\mathsf T}
 \preceq\frac{d+r}{r}I_d.
\]
Hence $M=L\sqrt{(d+r)/r}$ controls finite-atomic approximation even without exponential moments. The bound is independent of the cardinalities of both laws; an atomic driver may need randomized conditional allocation. The sharpness assertion concerns the deficit $1/(2M)$, rather than the perturbation factor $2M$.

\begin{corollary}[Gaussian reference perturbation]\label{R13-thm:reference-rate}
For every $\mu\cx\nu$ satisfying~\eqref{CM-eq:curvature}, arbitrary $q\in\cP_2$ and $\eta$-optimal $\pi$, the same-initial-state coupling satisfies
\begin{equation}\label{R13-eq:reference-rate}
 \left(\int W_2^2(\pi_x,\pi_x^\gamma)\,\mu(dx)\right)^{1/2}
 \le\|Y-T_X(G)\|_2\le L\delta+\sqrt{L^2\delta^2+2L\eta},
 \quad\delta=W_2(q^0,\gamma_d).
\end{equation}
For an exact optimizer this is
\begin{equation}\label{R13-eq:reference-Lipschitz}
 \left(\int W_2^2(\pi_x^q,\pi_x^\gamma)\,\mu(dx)\right)^{1/2}
 \le2L W_2(q^0,\gamma_d).
\end{equation}
\end{corollary}
\begin{proof}
Take $q_0=\gamma_d$ and $K=1$ in Theorem~\ref{R14-thm:reference-curvature}. Its comparison is on the same Gaussian input, and its balance is
\begin{equation}\label{R13-eq:perturb-balance}
 r^2/(2L)\le\eta+r\delta.
\end{equation}
\end{proof}

\begin{corollary}[First variation in the reference]\label{R15-cor:variation}
Under Theorem~\ref{R14-thm:reference-curvature}, let $\bar T(w)=\int T_x^{q_0}(w)\,\mu(dx)$. For $v\in L^2(q_0;\R^d)$ and $q_t=(\id+tv)_\#q_0$,
\begin{equation}\label{R15-eq:variation}
 0\le P_{q_t}-P_{q_0}-t\int\langle\bar T(w),v(w)\rangle\,q_0(dw)
 \le\frac M2t^2\int|v|^2\,dq_0.
\end{equation}
The law $q_t$ may be atomic. More generally the remainder for any coupling of $W_0\sim q_0$ and $W_1\sim q_1$, independent of $X$, is at most $\tfrac M2\E|W_1-W_0|^2$ after subtracting $\E\langle T_X^{q_0}(W_0),W_1-W_0\rangle$. If both references have deficit constants $1/(2M_j)$, their canonical outputs on that coupling satisfy
\[
 \|T_X^{q_1}(W_1)-T_X^{q_0}(W_0)\|_2
 \le\frac{2M_0M_1}{M_0+M_1}\|W_1-W_0\|_2.
\]
\end{corollary}
\begin{proof}
Using the old optimal endpoint with $W_1$ proves the lower bound. Glue a new optimizer to the drivers. The upper remainder is bounded by
$r\|W_1-W_0\|_2-r^2/(2M)$, whose maximum is $M\|W_1-W_0\|_2^2/2$. For two regular references, add their deficits with the other's endpoint as competitor; the optimal values cancel. Cauchy--Schwarz gives the displayed harmonic coefficient. The common Lipschitz bound and finite moments make $\bar T$ and all integrals well defined.
\end{proof}

\subsection{Continuous-time variational realizations}\label{QB-sec:dynamics}
The continuous-time question raised after Definition~1.3 of
\cite{QB-Tsch} requires a choice of driving process. Its filtration
specifies what is observed before time one, and its terminal law
determines the covariance optimization. We prove one variational
identity for every square-integrable reference martingale, then give
continuous Brownian and c\`adl\`ag L\'evy realizations. Under
Theorem~\ref{QB-thm:main}, their endpoint maps use one common potential.

\begin{theorem}[Optimal realization for a specified reference martingale]
\label{R13-thm:reference-process}
Let $N$ be a centered square-integrable $\R^d$-valued martingale on
$[0,1]$, with $N_0=0$, in a filtration satisfying the usual conditions
and having trivial initial sigma-field. Let $q=\law(N_1)$, and take
$X\sim\mu$ independently of the reference space, where
$\mu\cx\nu$ in $\cP_2(\R^d)$. After adjoining one independent scalar
Brownian motion, every optimizer of $P_q(\mu,\nu)$ is the endpoint
coupling of a square-integrable martingale $M$ with $M_0=X$ and
$M_1\sim\nu$. It solves
\begin{align}
 \sup_{\widetilde M}\E\tr\langle\widetilde M,N\rangle_1
   &=P_q(\mu,\nu),\label{R13-eq:general-driver}\\
 \inf_{\widetilde M}\E\tr\langle\widetilde M-N\rangle_1
   &=m_2(\nu)-m_2(\mu)+m_2(q)-2P_q(\mu,\nu),
 \label{R13-eq:general-action}
\end{align}
over adapted martingales with the prescribed initial variable and
terminal law. No nondegeneracy or irreducibility is required for
existence. If $q$ satisfies the nondegeneracy hypothesis, the additional
Brownian motion is unnecessary and the optimizer in each permitted
filtration is unique. For an irreducible pair its terminal value is
$\nabla v(a(X)+N_1)$.
\end{theorem}
\begin{proof}
Choose a primal optimizer and, measurably at each $x$, an optimal
covariance coupling of its conditional terminal law to $q$.
Disintegrate the resulting triple to write
$Y=H(X,N_1,U)$, where $U$ is independent and uniform on $(0,1)$.
One independent scalar Brownian motion $\beta$ realizes this variable
as $U=\Phi(\beta_1)$. Conditional expectation of $Y$ in the enlarged
filtration gives a c\`adl\`ag square-integrable martingale $M$, with
$M_0=X$ and the prescribed endpoint coupling. Independent enlargement
preserves the martingale $N$.

For any competitor, conditional covariance gives
\[
 \E\langle\widetilde M_1,N_1\rangle
 \le\int\MCov(\law(\widetilde M_1\mid X=x),q)\,\mu(dx)
 \le P_q.
\]
The constructed triple attains both inequalities. The product and
square-martingale identities give
\eqref{R13-eq:general-driver}--\eqref{R13-eq:general-action}.
These identities include all continuous and jump covariations.

For nondegenerate $q$, equality determines the unique endpoint kernel
and then its unique optimal map from $q$. Thus $Y=T_X(N_1)$ is
measurable without $U$, and $M_t=\E[Y\mid\mathcal F_t]$ is unique.
Theorem~\ref{QB-thm:main} gives the common potential under
irreducibility.\end{proof}

The terminal law fixes the optimal coupling; the reference filtration
specifies its realization before time one. The following specializations
keep both features explicit. They address the prescribed-driver formulation
of~\cite[Section~1.3]{QB-Tsch}; the one-dimensional additive-process
formulation appears in~\cite[Section~7]{QB-AM}.

\subsubsection{Continuous Brownian and prescribed L\'evy references}
For a centered nondegenerate $q$, let $S_\#\gamma_d=q$ be the Brenier map,
let $P_t$ be the heat semigroup, and take $B$ independent of $X\sim\mu$.
Under the irreducible hypotheses of Theorem~\ref{QB-thm:main}, put
\begin{align}
 N_t&=P_{1-t}S(B_t),\label{QB-eq:Brownian-reference}\\
 M_t&=\E[\nabla v(a(X)+S(B_1))\mid X,B_t].\label{QB-eq:Brownian-Bass}
\end{align}
\begin{theorem}[Continuous variational realization]\label{QB-thm:Brownian}
The process in \eqref{QB-eq:Brownian-Bass} is a continuous square-integrable
martingale with $M_0=X$ and endpoint coupling $\pi^q$. In the augmented
filtration of $(X,B)$, and in every enlargement where $B$ remains Brownian
and $X$ is initially measurable, it is the unique optimizer of
\begin{align}
 \sup_{\widetilde M}\E\tr\langle\widetilde M,N\rangle_1
       &=P_q(\mu,\nu),\label{QB-eq:dynamic-cov}\\
 \inf_{\widetilde M}\E\tr\langle\widetilde M-N\rangle_1
       &=m_2(\nu)-m_2(\mu)+m_2(q)-2P_q(\mu,\nu),
       \label{QB-eq:dynamic-action}
\end{align}
where competitors have initial variable $X$ and terminal law $\nu$.
\end{theorem}
\begin{proof}
Brownian martingale representation gives continuity up to time one for
both square-integrable endpoint functions~\cite[Chapter~V]{QB-RY}.
Independent future increments identify the same conditional expectations
in every permitted enlargement. Apply Theorem~\ref{R13-thm:reference-process}.
\end{proof}
In the original filtration the action is
$\E\int_0^1\|\sigma_t-\eta_t\|_{\mathrm{HS}}^2dt$ for
$dM=\sigma\,dB$, $dN=\eta\,dB$. For $q=\gamma_d$, $S=\id$ and this is
the martingale Benamou--Brenier problem. Deterministic covariance caps require
the additional curvature hypotheses of Sections~\ref{R11-sec:bounded}
and~\ref{R11-sec:canonical}.

\begin{theorem}[L\'evy-reference realization]\label{QB-thm:Levy}
Let $L$ be a centered square-integrable L\'evy process, independent of $X$,
with $q=\law(L_1)$ satisfying Theorem~\ref{QB-thm:main}, and write
$q_t=\law(L_t)$. Then
\begin{equation}\label{QB-eq:Levy}
 M_t=\int\nabla v(a(X)+L_t+z)\,q_{1-t}(dz)
\end{equation}
has a c\`adl\`ag version with endpoint coupling $\pi^q$. It uniquely maximizes
$\E\tr\langle\widetilde M,L\rangle_1$ and minimizes
$\E\tr\langle\widetilde M-L\rangle_1$, with the values in
\eqref{QB-eq:dynamic-cov}--\eqref{QB-eq:dynamic-action}. This holds in the
filtration of $(X,L)$ and in every enlargement retaining independent future
increments and initial measurability of $X$.
\end{theorem}
\begin{proof}
Independent increments give \eqref{QB-eq:Levy} as the conditional expectation
of the optimal endpoint. Martingale regularization gives the c\`adl\`ag
version, and Theorem~\ref{R13-thm:reference-process} gives the claims,
including jump covariations.
\end{proof}
For $L=\Sigma B+\int z\,\widetilde{\mathcal N}(dt,dz)$, the action of a
competitor with Brownian coefficient $A_t$, compensated-jump coefficient
$H_t(z)$ and orthogonal martingale component $K$ is
\[
 \E\int_0^1\!\|A_t-\Sigma\|_{\mathrm{HS}}^2dt
 +\E\int_0^1\!\int|H_t(z)-z|^2\,\Pi(dz)dt
 +\E\tr\langle K\rangle_1.
\]
An atomic compound-Poisson terminal reference instead uses the one scalar
Brownian enlargement in Theorem~\ref{R13-thm:reference-process}; its
allocation may be randomized and nonunique. A nondegenerate Brownian
component restores the unique common-gradient representation.

\subsection{Preserving the entire signing law}
\begin{corollary}[Strict terminal inflation]\label{QB-cor:inflation}
Let $\mu\cx\nu$ be any pair in $\cP_2(\R^d)$ with common mean $m$ and
full-dimensional terminal affine hull. For $c>1$ define
\[
 \nu_c=\bigl(y\mapsto m+c(y-m)\bigr)_\#\nu.
\]
Then $(\mu,\nu_c)$ has the representation of Theorem~\ref{QB-thm:main}
for every reference $q$ satisfying its assumptions. The entire initial
law $\mu$ is unchanged. The globally finite-potential conclusion holds
under the same alternative conditions on $q$ or $\nu$.
\end{corollary}

\begin{proof}
Take $\pi\in\mathcal M(\mu,\nu)$ and form
$c^{-1}\pi+(1-c^{-1})\mu\otimes\nu$. Its conditional terminal mean is
$m+(x-m)/c$. Applying $y\mapsto m+c(y-m)$ gives a martingale coupling of
$\mu,\nu_c$ which dominates $(1-c^{-1})\mu\otimes\nu_c$. The pair is
therefore irreducible, and Theorem~\ref{QB-thm:main} applies.
\end{proof}

\begin{corollary}[A general driving reference for the entire signing law]\label{QB-cor:signing}
Suppose a prescribed sign distribution satisfies
\[
 S=(A\sigma,\sigma)\cx N(0,\Gamma),\qquad \Gamma\succ0.
\]
For every $c>1$ and every centered $q\in\cP_2(\R^{m+n})$ satisfying the
Hausdorff-null hypothesis, there exist $a$ and a proper convex potential
$v$ such that, with $W\sim q$ independent of $\sigma$,
\[
 R=\nabla v(a(S)+W),\qquad
 R\sim N(0,c^2\Gamma),\qquad \E[R\mid\sigma]=S.
\]
The law of $\sigma$ is exactly the original law. The coupling of $S,R$
maximizes the integrated maximal covariance with $q$ among all
martingale couplings of those endpoints. If the reference $q$ has full
convex support, $v$ is finite everywhere. It has the continuous
realizations of Section~\ref{QB-sec:dynamics} for the specified drivers.
\end{corollary}
\begin{proof}
Apply Corollary~\ref{QB-cor:inflation} to $\mu=\law(S)$ and
$\nu=N(0,\Gamma)$, then Theorem~\ref{QB-thm:main}. Since $S$ contains
$\sigma$, conditioning on either variable gives the same conditional
mean identity. Its second marginal is exactly the stated Gaussian.
\end{proof}

This corollary retains every atom of the signing law and hence every
hard discrepancy constraint, prescribed sign marginal, and entropy
statement that depends only on that law. Independent Gaussian terminal
blocks remain independent when $\Gamma$ is block diagonal. The new
reference--sign joint law generally differs from a prior construction,
so information quantities involving that former joint law are not
asserted to survive. The driver $q$ may itself have dependent
coordinates. Its whole joint law is kept fixed.

\begin{example}[An exact non-Gaussian convolution model]\label{QB-ex:quadratic}
Let $K\succ0$ be symmetric, let $q$ be centered with finite covariance
$Q$, and let $\nu=\mu*(K_\# q)$. Then
\[
 v(z)=\tfrac12\ip{z}{Kz},\qquad a(x)=K^{-1}x,
 \qquad Y=X+KW
\]
gives a common-gradient representation. Its optimality can be verified directly,
without an existence theorem. For
$\psi(y)=\tfrac12\ip{y}{K^{-1}y}$, completing the square gives
\[
 \Phi_\psi(x)=\tfrac12\ip{x}{K^{-1}x}-\tfrac12\tr(KQ),
 \qquad P_q(\mu,\nu)=\tr(KQ).
\]
For any admissible triple with the same endpoints, $W\sim q$ independent
of $X$, the exact deficit identity is
\[
 \tr(KQ)-\E\ip{Y}{W}
 =\frac12\E\ip{Y-X-KW}{K^{-1}(Y-X-KW)}.
\]
The identity uses only the martingale mean and the two endpoint second
moments. Every loss of covariance is exactly a squared deviation from
the same additive-noise coupling. This is the quadratic model for the
general conditional Fenchel argument above.
\end{example}

\begin{corollary}[Finite-reference calibration of the original signing law]
\label{R13-cor:signing-reference}
Let $S=(A\sigma,\sigma)\cx N(0,\Gamma)$ with $\Gamma\succ0$,
and keep the entire law of $\sigma$ fixed. In whitened coordinates
let $\pi^q$ be any optimal endpoint coupling with driving law
$q\in\cP_2$, and let $\pi^\gamma$ be the canonical Gaussian-reference
coupling. Both have terminal law $N(0,\Gamma)$, and
\[
 \left(\E W_{2,\Gamma^{-1}}^2
  \bigl(\law(R^q\mid\sigma),\law(R^\gamma\mid\sigma)\bigr)
 \right)^{1/2}
 \le 2W_2(q^0,\gamma_{m+n}).
\]
For an $\eta$-optimal endpoint coupling the right side is
$\delta+\sqrt{\delta^2+2\eta}$, where
$\delta=W_2(q^0,\gamma_{m+n})$ and the objective is whitened.
Finite atomic $q$ are permitted. Every signing probability, hard
constraint and entropy statement depending only on $\sigma$ remains
unchanged, as does the complete Gaussian terminal law.
\end{corollary}
\begin{proof}
Apply Theorem~\ref{R13-thm:reference-rate} to
$\mu=\law(\Gamma^{-1/2}S)$, $\nu=\gamma_{m+n}$ and $L=1$.
The inclusion of $\sigma$ in $S$ identifies the conditional kernels.
Undoing the whitening gives the stated metric.
\end{proof}

The terminal reference in this corollary is the original Gaussian.
The inflation in Corollary~\ref{QB-cor:signing} serves the separate
purpose of obtaining one common potential for every admissible driver.

For Gaussian-reference approximation, the whole selected coupling is
therefore controlled before a latent calibration scheme is specified.
The next section proves quantitative curvature and trajectory
confinement, which yield convergence rates for Bass calibration.


\section{Strong convexity from positive noise and Bass calibration}
\label{R12-sec:bass}
Calibration varies the latent law while transporting its noisy image to a fixed target. The transport itself changes with the latent perturbation, so its optimizing multipliers must be differentiated as well. After their contributions cancel, a conditional covariance remains: positive noise leaves uncertainty in every nonconstant centered perturbation. Averaging the transport's distributional Hessian converts that uncertainty into curvature. For atomic targets this Hessian lives on cell boundaries; full target affine span is exactly what makes the averaged matrix positive.

With finite second moments and a continuous everywhere-positive noise density, this proves strong convexity on fixed-mean slices of bounded $L^\infty$ balls. The modulus is uniform over $W_2$-compact full-span target families. For Gaussian noise, confinement bounds the full latent trajectory before this modulus is used. It gives the exponential rate in every dimension conjectured by Backhoff--Pammer--Schachermayer~\cite{BF-BPS}, retaining both endpoint marginals at finite time and preserving the probabilities of a calibrated signing law.

\subsection{Positive noise and full-dimensional targets}
\label{BF-sec:results}
For $\rho,\nu\in\mathcal P_2(\R^d)$, let
$\MCov(\rho,\nu)=\sup_{\pi\in\Pi(\rho,\nu)}\int x\cdot y\,\pi(dx,dy)$,
as in the preceding section.
For a square-integrable random vector $Z$ on an arbitrary probability space, and an independent noise vector $\xi$ with law $q$, put
\begin{equation}\label{BF-eq:F}
 F_{\nu,q}(Z)=\MCov(\law(Z+\xi),\nu),\qquad
 \Var_2(Z)=\E|Z-\E Z|^2.
\end{equation}
All Euclidean norms use the same fixed coordinate system. A law has full affine span when its support is contained in no proper affine hyperplane.

\begin{theorem}[Strong convexity under additive noise]\label{BF-thm:maincurvature}
Suppose $q\in\mathcal P_2(\R^d)$ has a continuous density which is strictly positive everywhere, and $\nu\in\mathcal P_2(\R^d)$ has full affine span. For every finite $R$ there is $\kappa(R,q,\nu)>0$ such that, whenever $|Z_0|,|Z_1|\le R$ almost surely and $0\le t\le1$,
\begin{align}\label{BF-eq:strong}
 F_{\nu,q}((1-t)Z_0+tZ_1)
 &\le (1-t)F_{\nu,q}(Z_0)+tF_{\nu,q}(Z_1)\notag\\
 &\quad-\frac{\kappa}{2}t(1-t)\Var_2(Z_1-Z_0).
\end{align}
The constant is independent of the probability space. It can be chosen uniformly when $\nu$ ranges over a $W_2$-compact family all of whose members have full affine span. Conversely, if $\nu$ is supported on a proper affine hyperplane, no positive constant works on all probability spaces and all nontrivial centered bounded perturbations.
\end{theorem}

The deterministic translations are the exact degeneracy in this statement:
\begin{equation}\label{BF-eq:translation}
 F_{\nu,q}(Z+a)=F_{\nu,q}(Z)+a\cdot\int y\,\nu(dy).
\end{equation}
Thus \eqref{BF-eq:strong} is strong convexity on each fixed-mean slice of an $L^\infty$ ball. The target may be discrete, singular or unbounded. The uniformity over $W_2$-compact full-span families allows finite approximations to retain a positive modulus as their number of atoms grows. The proof first establishes that uniformity geometrically, then passes from entropy-regularized finite targets to their $W_2$ limits.

\subsection{Curvature at transport-cell boundaries}
\label{BF-sec:curvature}
For an atomic target, the optimal convex potential is piecewise affine.
Its classical Hessian vanishes almost everywhere, although the
transport changes when mass crosses a cell boundary. In one dimension,
a symmetric two-point target and symmetric source have potential
$v(w)=|w|$. The distributional identity $D^2v=2\delta_0$ records
precisely the curvature lost by differentiating only inside the cells.

We use the positive matrix-valued Hessian measure of a convex
potential, averaged against a smooth positive function. Positive
additive noise bounds the posterior below by a multiple of the prior;
the resulting covariance estimate will multiply this Hessian measure.
Full affine span ensures that no direction has zero averaged curvature.
All transport maps in the proof are Brenier maps~\cite{Brenier}.
Write $g$ for the density of the noise $q$. For $R<\infty$, define
\[
 m_R(w)=\min_{|z|\le R}g(w-z).
\]
This is continuous and positive, and every density of the form $\alpha*g$, with $\alpha\in\mathcal P(B_R)$, is at least $m_R$.

\begin{lemma}[A smooth positive minorant]\label{BF-lem:minorant}
There is a smooth function $\eta_R>0$ such that $\eta_R\le m_R$ and
\begin{equation}\label{BF-eq:etaenergy}
 \int_{\R^d}\frac{\eta_R^2+|\nabla\eta_R|^2}{m_R}\,dw<\infty.
\end{equation}
For Gaussian noise one may take
\begin{equation}\label{BF-eq:gausseta}
 \eta_R(w)=(2\pi)^{-d/2}e^{-R^2}e^{-|w|^2}.
\end{equation}
\end{lemma}
\begin{proof}
Choose a smooth, nonnegative, locally finite partition of unity on annuli, with uniformly bounded overlap and derivatives. For its $n$th member $\chi_n$, choose a positive coefficient $a_n$ sufficiently small that $a_n\le\inf_{\supp\chi_n}m_R$ and
\[
 a_n^2\int_{\supp\chi_n}\frac{\chi_n^2+|\nabla\chi_n|^2}{m_R}\le 2^{-n}.
\]
Shrinking all coefficients by a fixed overlap constant, $\eta_R=\sum_n a_n\chi_n$ has the required properties. For \eqref{BF-eq:gausseta}, use $|w-z|^2\le2|w|^2+2R^2$ and $m_R\ge\eta_R$; the energy is bounded by a Gaussian integral.
\end{proof}

For a finite convex function $v$, its distributional Hessian $D^2v$ is a positive semidefinite matrix-valued Radon measure. Whenever the integral is finite, put
\begin{equation}\label{BF-eq:Kdef}
 K_R(v)=\int\eta_R(w)\,D^2v(dw).
\end{equation}

\begin{lemma}[Compactness and positive curvature]\label{BF-lem:positiveK}
Let $\mathcal N$ be a $W_2$-compact family of full-affine-span target laws. Normalize by $v(0)=0$ the Brenier potentials satisfying
\[
 (\nabla v)_\#(\alpha*q)=\nu,
 \qquad \alpha\in\mathcal P(B_R),\quad\nu\in\mathcal N.
\]
Their normalized potentials form a compact family for local uniform convergence. The matrices $K_R(v)$ are finite, depend continuously on this family, and satisfy
\begin{equation}\label{BF-eq:uniformK}
 \inf_{\alpha,\nu}\lambda_{\min}(K_R(v))>0.
\end{equation}
\end{lemma}
\begin{proof}
Set $M_2=\sup_{\nu\in\mathcal N}\int|y|^2\,d\nu<\infty$. The transport identity and the density lower bound give
\begin{equation}\label{BF-eq:weightedL2}
 \int|\nabla v(w)|^2m_R(w)\,dw\le M_2.
\end{equation}
This bounds gradients in $L^2$ on every compact ball. To obtain a pointwise local bound, take $p\in\partial v(x)$ with $x\in B_s$. If $p\ne0$, let $A_p=\{u\in B_1:u\cdot p/|p|\ge1/2\}$. Monotonicity of the subgradient implies
$\nabla v(x+u)\cdot u\ge p\cdot u\ge|p|/2$ for almost every $u\in A_p$. Hence
\[
 |p|^2\le \frac4{|A_p|}\int_{B_{s+1}}|\nabla v(w)|^2\,dw.
\]
The cap volume depends only on $d$. Thus the normalized potentials are locally uniformly bounded and locally equi-Lipschitz.

Given a sequence, extract $\alpha_n\Rightarrow\alpha$ in $\mathcal P(B_R)$, $\nu_n\to\nu$ in $W_2$, and $v_n\to v$ locally uniformly. The source densities converge locally uniformly, and the source laws converge in $W_2$. The transport plans have a subsequential limit with the prescribed marginals. Local uniform convergence of convex functions preserves their subgradient graphs at finite limit points. The limiting plan is therefore supported on $\partial v$. Since its first marginal has a positive density, it is the graph of $\nabla v$, so $(\nabla v)_\#(\alpha*q)=\nu$. This proves compactness.

The weighted Hessian admits the integration-by-parts formula
\begin{equation}\label{BF-eq:Kibp}
 (K_R(v))_{ij}=-\int \partial_i v\,\partial_j\eta_R\,dw.
\end{equation}
For completeness, apply distributional integration by parts with $\eta_R$ multiplied by a smooth cutoff on $B_s$, chosen increasing to one. The boundary terms vanish by \eqref{BF-eq:weightedL2}, \eqref{BF-eq:etaenergy} and Cauchy--Schwarz. Positivity and monotone convergence give finiteness of diagonal Hessian integrals; the off-diagonal entries then follow by polarization. The same estimates give the uniform tail bound
\begin{equation}\label{BF-eq:Ktail}
 \int_{|w|>s}|\nabla v|\,|\nabla\eta_R|
 \le M_2^{1/2}\left(\int_{|w|>s}\frac{|\nabla\eta_R|^2}{m_R}\right)^{1/2}\longrightarrow0.
\end{equation}
Local uniform convergence of finite convex functions implies almost-everywhere convergence of their gradients. Local gradient bounds, \eqref{BF-eq:Ktail} and \eqref{BF-eq:Kibp} prove continuity of $K_R$.

Finally, if $e^TK_R(v)e=0$ for a nonzero $e$, positivity of $\eta_R$ and of $D^2v$ imply $\partial_{ee}v=0$ as a measure. The matrix Cauchy--Schwarz inequality for positive semidefinite Hessian measures then gives $\partial_{ie}v=0$ for every $i$. Thus $\partial_ev$ is a constant distribution, and $e\cdot\nabla v$ is constant almost everywhere. The target law would lie in an affine hyperplane. Consequently $K_R(v)$ is positive definite, and compactness gives \eqref{BF-eq:uniformK}.
\end{proof}

\begin{remark}\label{BF-rem:uniformfamilies}
For example, fix $d,R,D,s>0$ and let targets satisfy $\supp\nu\subset B_D$ and $\Cov(\nu)\succeq sI_d$. This is a $W_2$-compact full-affine-span family. The constant in \eqref{BF-eq:uniformK}, and hence the strong-convexity constant, is uniform over the family. It has no dependence on an atomic target's number of atoms or minimum atom mass. For unbounded targets, the $W_2$-compactness hypothesis includes the required uniform integrability of second moments.
\end{remark}

The lower bound must survive finer and finer atomic approximations. Compactness gives a geometric lower bound: a vanishing directional Hessian would force the transported target into an affine hyperplane. Full affine span excludes that limit, so the smallest averaged Hessian eigenvalue stays positive on the compact target family. The next step couples this geometric bound to the uncertainty in the latent perturbation that remains after observing the noisy input.

\begin{lemma}[Posterior covariance minorization]\label{BF-lem:posterior}
Let $Z,H$ be random vectors on the base probability space, with $|Z|\le R$ and $H\in L^2$, and let $W=Z+\xi$ with independent noise. If $\rho$ is the density of $W$, then
\begin{equation}\label{BF-eq:posterior}
 \Cov(H\mid W=w)\succeq\frac{\eta_R(w)}{\rho(w)}\Cov(H)
 \quad\text{for almost every }w.
\end{equation}
\end{lemma}
\begin{proof}
The posterior distribution on the original probability space has density
$g(w-Z)/\rho(w)\ge\eta_R(w)/\rho(w)$ relative to the prior. Write it as a mixture of that multiple of the prior and a residual probability law. The covariance formula for mixtures gives \eqref{BF-eq:posterior}, because both the residual covariance and the between-means term are positive semidefinite.
\end{proof}

\subsubsection{An optimized second variation}
\label{BF-sec:entropy}
The optimized transport must continue to have the prescribed target masses as the latent input moves. Its dual multipliers therefore vary with that input. Differentiating after optimization cancels their second derivatives by the marginal equations. The remaining variance splits conditionally on the noisy observation: one nonnegative term contains the multipliers, while the other pairs the unresolved latent covariance with the transport Hessian. Positive noise and full affine span bound this second term uniformly from below, even as entropy regularization vanishes.

We first prove the curvature theorem for a finite target
$\nu=\sum_{j=1}^m p_j\delta_{y_j}$ with distinct atoms and positive weights. Let $\varepsilon>0$ and define
\begin{align}\label{BF-eq:softF}
 F_{\nu,q,\varepsilon}(Z)
 =\min_{b\in\R^m/\R\mathbf1}
 \left\{\E\left[\varepsilon\log\sum_{j=1}^m p_j
 e^{((Z+\xi)\cdot y_j-b_j)/\varepsilon}\right]
 +\sum_{j=1}^m p_jb_j\right\}.
\end{align}
The minimization is finite dimensional. Its objective is strictly convex on the quotient: its Hessian there is the expectation of a strictly positive multinomial covariance. Coercivity follows after normalizing $\min_j b_j=0$: if $D=\max_j|y_j|$ and $p_{\min}=\min_jp_j$, the objective is at least
\[
 p_{\min}\max_j b_j-D\E|Z+\xi|+\varepsilon\log p_{\min}.
\]
Hence the minimizing $b$ exists and is unique modulo constants. Its stationarity equations are
\begin{equation}\label{BF-eq:softmarginal}
 \E q_j(W)=p_j,\qquad
 q_j(w)=\frac{p_j e^{(w\cdot y_j-b_j)/\varepsilon}}
 {\sum_i p_i e^{(w\cdot y_i-b_i)/\varepsilon}},\quad W=Z+\xi.
\end{equation}
The finite entropy duality identity gives
\begin{equation}\label{BF-eq:softerror}
 F_{\nu,q}(Z)-\varepsilon H(p)
 \le F_{\nu,q,\varepsilon}(Z)\le F_{\nu,q}(Z).
\end{equation}
Indeed \eqref{BF-eq:softF} maximizes covariance minus $\varepsilon$ times relative entropy with respect to $\law(W)\otimes\nu$. An optimal unregularized coupling is a map from $W$ and has relative entropy $H(p)=-\sum_jp_j\log p_j$. This also proves \eqref{BF-eq:softerror} directly.

\begin{lemma}[A multidimensional second-variation bound]\label{BF-lem:softHessian}
Let $Z_t=Z_0+tH$ with $Z_0,H\in L^\infty$ and $|Z_t|\le R$ for $0\le t\le1$. Let $b(t)$ solve \eqref{BF-eq:softF}, and put
\[
 v_t(w)=\varepsilon\log\sum_jp_j e^{(w\cdot y_j-b_j(t))/\varepsilon}.
\]
Then
\begin{equation}\label{BF-eq:secondvariation}
 \frac{d^2}{dt^2}F_{\nu,q,\varepsilon}(Z_t)
 \ge \tr\bigl(K_R(v_t)\Cov(H)\bigr).
\end{equation}
\end{lemma}
\begin{proof}
The implicit function theorem on $\mathbf1^\perp$ makes $b(t)$ smooth. Differentiation under the expectation is justified by the bounded atoms and bounded $H$. The terms involving $b''$ cancel by \eqref{BF-eq:softmarginal}, leaving
\begin{equation}\label{BF-eq:fullsoftsecond}
 \frac{d^2}{dt^2}F_{\nu,q,\varepsilon}(Z_t)
 =\frac1\varepsilon\E\Var_{q(W_t)}\bigl(H\cdot Y-\dot b_Y(t)\bigr),
 \qquad W_t=Z_t+\xi.
\end{equation}
Here, after conditioning on $W_t=w$, $Y$ has probabilities $q_j(w)$ and is conditionally independent of $H$. Write $a=\E[H\mid W_t]$ and $H=a+r$. Expanding the conditional variance and averaging in $H$ gives the exact identity
\begin{align}\label{BF-eq:varianceidentity}
 &\E\!\left[\Var_{q(W_t)}(H\cdot Y-\dot b_Y)\mid W_t\right]\notag\\
 &\quad=\tr\bigl(\Cov(H\mid W_t)\Cov(Y\mid W_t)\bigr)
 +\Var_{q(W_t)}(a\cdot Y-\dot b_Y).
\end{align}
The mixed term vanishes because $\E[r\mid W_t]=0$. Also
$D^2v_t(w)=\varepsilon^{-1}\Cov(Y\mid W_t=w)$. Discard the last, nonnegative term of \eqref{BF-eq:varianceidentity} and use Lemma~\ref{BF-lem:posterior}. Multiplication by the density of $W_t$ cancels its denominator and yields \eqref{BF-eq:secondvariation}.
\end{proof}

\begin{lemma}[Zero-temperature averaged curvature]\label{BF-lem:zerotemp}
Fix a finite full-affine-span target $\nu$. If $c_0>0$ is a lower bound in \eqref{BF-eq:uniformK} for this target, then, for all sufficiently small $\varepsilon>0$, uniformly over $\alpha\in\mathcal P(B_R)$, the minimizing potential in \eqref{BF-eq:softF} satisfies
\[
 K_R(v_{\alpha,\varepsilon})\succeq(c_0/2)I_d.
\]
\end{lemma}
\begin{proof}
Suppose otherwise, and choose $\varepsilon_n\downarrow0$ and $\alpha_n\in\mathcal P(B_R)$ violating the assertion. Pass to $\alpha_n\Rightarrow\alpha$. After normalizing $\min b_j=0$, the coercivity estimate preceding \eqref{BF-eq:softmarginal} and the comparison with $b=0$ give a uniform bound on all $b_j$, because $\E|Z+\xi|\le R+\E|\xi|$. Pass to $b^n\to b$ and normalize the potentials at zero. They converge locally uniformly to the corresponding normalized maximum of the affine functions $w\cdot y_j-b_j$.

Except on finitely many hyperplanes, the maximizing atom is unique. The stationarity equations \eqref{BF-eq:softmarginal}, local uniform convergence of the source densities, and tightness of the source laws show that the limiting gradient pushes $\alpha*q$ to $\nu$. Every soft gradient is bounded by $D$. Formula \eqref{BF-eq:Kibp} and $\nabla\eta_R\in L^1$, a consequence of \eqref{BF-eq:etaenergy}, give convergence of the matrices $K_R$. Their limit is at least $c_0I_d$, a contradiction.
\end{proof}

\begin{proof}[Proof of Theorem~\ref{BF-thm:maincurvature}]
For a finite target, combine Lemmas~\ref{BF-lem:positiveK}, \ref{BF-lem:softHessian} and \ref{BF-lem:zerotemp} to obtain
\[
 \frac{d^2}{dt^2}F_{\nu,q,\varepsilon}(Z_t)
 \ge(c_0/2)\Var_2(H).
\]
Integrating twice gives \eqref{BF-eq:strong} with $\kappa=c_0/2$ for the regularized functional. Letting $\varepsilon\downarrow0$ and using \eqref{BF-eq:softerror} proves it for $F_{\nu,q}$.

For a general target, choose finite laws $\nu_n\to\nu$ in $W_2$, for example conditional-mean quantizations on refining finite partitions. The covariance matrices converge, so all sufficiently late $\nu_n$ have full affine span. The compact family consisting of these targets and their limit has a common positive $c_0$ in Lemma~\ref{BF-lem:positiveK}. Each finite-target proof therefore uses the same $\kappa=c_0/2$. Coupling an optimal target pair and applying Cauchy--Schwarz gives
\begin{equation}\label{BF-eq:targetstability}
 |F_{\nu_n,q}(Z)-F_{\nu,q}(Z)|
 \le(R+\|\xi\|_2)W_2(\nu_n,\nu),\qquad |Z|\le R.
\end{equation}
This proves \eqref{BF-eq:strong} in the limit.

For the stated uniformity over a compact target family $\mathcal N$, its covariance matrices have a common positive lower eigenvalue. Finite spatial quantization after truncation gives uniform $W_2$ approximation on $\mathcal N$. Take the approximations sufficiently fine that all their covariance matrices retain half that lower bound. The closure of the union of the resulting approximating families is $W_2$-compact: truncation errors are uniformly small by uniform integrability, and each fixed finite grid gives a compact family. Every member of this closure has full affine span. Lemma~\ref{BF-lem:positiveK} gives one $c_0$ for the complete family. The preceding proof then gives a common $\kappa$ on $\mathcal N$.

Finally, if $e\cdot y=a$ $\nu$-almost surely, replacing $Z$ by $Z+he$, for any bounded scalar random variable $h$, changes every admissible covariance by $a\E h$. Taking $\E h=0$ and $\Var(h)>0$ disproves positive curvature in that direction. This proves the converse.
\end{proof}

\begin{remark}[No pointwise positive Hessian hypothesis]\label{BF-rem:singular}
For a finite target, the limiting Brenier potential is piecewise affine and its classical Hessian vanishes almost everywhere. The positive matrix in Lemma~\ref{BF-lem:positiveK} integrates its \emph{distributional} Hessian, including the curvature on cell boundaries. Passing through \eqref{BF-eq:fullsoftsecond} before taking the limit retains that curvature. This is why discrete targets are covered by the same proof.
\end{remark}

\subsection{The lifted Bass flow and its exponential rate}\label{BF-sec:flow-statement}
For the Gaussian specialization, let $G\sim\gamma_d=N(0,I_d)$ be independent of a fixed $X\sim\mu$. On $\mathsf H=L^2(\sigma(X);\R^d)$, define the lifted Bass functional
\begin{equation}\label{BF-eq:Bass}
 V(Z)=F_{\nu,\gamma_d}(Z)-\E[Z\cdot X].
\end{equation}
If $T_Z=\nabla v_Z$ is the Brenier map from $\law(Z)*\gamma_d$ to $\nu$, its gradient is
\begin{equation}\label{BF-eq:gradient}
 DV(Z)=\E[T_Z(Z+G)\mid X]-X=(\nabla v_Z*\gamma_d)(Z)-X.
\end{equation}
Convexity, differentiability and this formula are established for the lifted Bass functional in \cite[Lemma~2.1]{BF-BPS}; a proof of the facts needed here is included in Lemma~\ref{BF-lem:gradient}.

A convex-order pair $(\mu,\nu)$ is \emph{irreducible} if every pair of measurable sets of positive $\mu$- and $\nu$-mass is joined with positive mass by some martingale coupling. This is the irreducibility used in the Bass representation theorem \cite{CM-Bass}.

\begin{theorem}[Exponential Bass flow in every dimension]\label{BF-thm:flow}
Let $\mu\preceq_{\mathrm{cx}}\nu$ be irreducible, let $\nu\in\mathcal P_2(\R^d)$ have full affine span, and suppose
\begin{equation}\label{BF-eq:interior}
 K:=\operatorname{conv}(\supp\mu)\Subset
 \operatorname{int}\overline{\operatorname{conv}}(\supp\nu).
\end{equation}
For every $Z_0\in L^\infty(\sigma(X);\R^d)$, the Bass gradient flow
\[
 \dot Z_t=-DV(Z_t)
\]
exists globally and stays uniformly bounded. It converges to a bounded minimizer $Z^*$ having $\E Z^*=\E Z_0$. This minimizer is unique with that mean. For some $\kappa>0$,
\begin{align}\label{BF-eq:rates}
 \|Z_t-Z^*\|_2&\le e^{-\kappa t}\|Z_0-Z^*\|_2,\\
 V(Z_t)-V(Z^*)&\le e^{-2\kappa t}\bigl(V(Z_0)-V(Z^*)\bigr).\notag
\end{align}
The terminal law is only required to have a finite second moment; compact support and a density are unnecessary.
\end{theorem}

Backhoff--Pammer--Schachermayer proved strong convergence in every
dimension and exponential convergence in one dimension under additional
target-density assumptions~\cite[Theorems~1.8--1.9]{BF-BPS}.
They conjectured the exponential conclusion in higher dimensions.
Corollary~\ref{BF-cor:published} below proves that conclusion from their
boundedness theorem and Theorem~\ref{BF-thm:maincurvature}.
Theorem~\ref{BF-thm:flow} also permits unbounded terminal support; its
separate confinement argument is in Section~\ref{BF-sec:unbounded}.

The rate constant is the strong-convexity modulus on a ball containing
the complete trajectory and its minimizer. Confinement is established
before this modulus is used. Uniform rates require a common latent
bound and a compact target family as in Theorem~\ref{BF-thm:maincurvature}.

\begin{corollary}[Arbitrarily small reference inflation]\label{BF-cor:inflation}
Let $\mu\preceq_{\mathrm{cx}}\nu$, with $\mu$ compactly supported and $\nu\in\mathcal P_2(\R^d)$ of full affine span. Write $m$ for the common mean. For every $c>1$, set
\[
 \nu_c=(y\mapsto m+c(y-m))_\#\nu.
\]
Then $(\mu,\nu_c)$ satisfies the hypotheses of Theorem~\ref{BF-thm:flow}. Consequently the lifted Bass flow for this pair converges exponentially from every bounded initial value.
\end{corollary}
\begin{proof}
Let $\pi$ be any martingale coupling of $(\mu,\nu)$ and set $\lambda=1/c$. The coupling
$\lambda\pi+(1-\lambda)\mu\otimes\nu$ has conditional terminal mean $m+\lambda(x-m)$. Dilating its terminal coordinate by $c$ gives a martingale coupling of $(\mu,\nu_c)$ containing a positive multiple of the product coupling. This proves irreducibility. Since $\nu$ has full affine span, $m$ lies in the interior of its closed convex hull. Strict dilation about $m$ carries every compact subset of the original closed convex hull into the interior of the dilated hull. This proves \eqref{BF-eq:interior}.
\end{proof}

\subsubsection{From strong convexity to convergence}
\label{BF-sec:basicflow}
We now specialize to Gaussian noise and the functional \eqref{BF-eq:Bass}. All random variables in its domain are measurable with respect to $\mathcal G=\sigma(X)$, while $G$ is independent of $\mathcal G$.

\begin{lemma}[Gradient and monotonicity]\label{BF-lem:gradient}
The functional $V$ is finite, convex and continuously Fr\'echet differentiable on $\mathsf H$, with gradient \eqref{BF-eq:gradient}. Its gradient satisfies
\[
 \|DV(Z)\|_2\le \left(\int|y|^2\,d\nu\right)^{1/2}+\|X\|_2.
\]
If $\E X=\int y\,d\nu$ then $\E DV(Z)=0$. On any bounded fixed-mean slice covered by Theorem~\ref{BF-thm:maincurvature},
\begin{align}\label{BF-eq:bregman}
 V(Z_1)&\ge V(Z_0)+\langle DV(Z_0),Z_1-Z_0\rangle_{\mathsf H}
       +\frac\kappa2\|Z_1-Z_0\|_2^2,\\
 \langle DV(Z_1)-DV(Z_0),Z_1-Z_0\rangle_{\mathsf H}
       &\ge\kappa\|Z_1-Z_0\|_2^2.\notag
\end{align}
\end{lemma}
\begin{proof}
The maximal covariance is the supremum of the affine functions
$Z\mapsto\E[(Z+G)\cdot Y]$, where $Y\sim\nu$ is measurable on the product space. Brenier's theorem gives its unique maximizer $Y_Z=T_Z(Z+G)$. If $Z_n\to Z$ in $L^2$, stability of optimal transport plans and uniqueness of the limiting graph imply $Y_{Z_n}\to Y_Z$ in probability. One can see the last assertion on the common probability space by first restricting $T_Z$ to a compact Lusin set. Since every $Y_{Z_n}$ has the same law $\nu$, uniform integrability of their squared norms improves this to $L^2$ convergence.

The supremum identity gives
\[
 0\le F_{\nu,\gamma_d}(Z+h)-F_{\nu,\gamma_d}(Z)-\E[h\cdot Y_Z]
 \le\|h\|_2\|Y_{Z+h}-Y_Z\|_2.
\]
This proves continuous Fr\'echet differentiability and the conditional-expectation formula. The gradient bound and its mean follow immediately. Finally, differentiation of \eqref{BF-eq:strong} at an endpoint proves the first inequality in \eqref{BF-eq:bregman}; adding it in both directions proves the second.
\end{proof}

\begin{proposition}[Exponential convergence on a bounded trajectory]\label{BF-prop:boundedflow}
Suppose $V$ has a minimizer $Z^*$, and a gradient trajectory and this minimizer lie in the same bounded fixed-mean slice. Then \eqref{BF-eq:rates} holds with the strong-convexity constant of that slice. In addition,
\begin{equation}\label{BF-eq:PL}
 \frac\kappa2\|Z-Z^*\|_2^2\le V(Z)-V(Z^*)
 \le\frac1{2\kappa}\|DV(Z)\|_2^2
\end{equation}
for every $Z$ in the slice.
\end{proposition}
\begin{proof}
The Bregman bound~\eqref{BF-eq:bregman} at $Z^*$ gives the lower estimate in~\eqref{BF-eq:PL}. Applying it at $Z$ and maximizing
$\langle DV(Z),h\rangle-\kappa\|h\|_2^2/2$ over $h$ gives the upper estimate. Along the flow, strong monotonicity and this upper estimate give
\[
 \frac d{dt}\|Z_t-Z^*\|_2^2\le-2\kappa\|Z_t-Z^*\|_2^2,
 \qquad
 \frac d{dt}(V(Z_t)-V(Z^*))
 =-\|DV(Z_t)\|_2^2\le-2\kappa(V(Z_t)-V(Z^*)).
\]
Integrating proves both rates.
\end{proof}

\begin{corollary}[The published higher-dimensional conjecture]\label{BF-cor:published}
Under Assumption~(A) of \cite{BF-BPS}, a bounded initial value gives the two exponential bounds \eqref{BF-eq:rates} in every dimension. The extra target-density assumptions of their Theorem~1.9 can be omitted.
\end{corollary}
\begin{proof}
Assumption~(A) gives an irreducible convex-order pair, a compact full-dimensional terminal convex hull, and an initial support compactly contained in its interior. Theorem~1.8 of \cite{BF-BPS} gives a bounded minimizer and a uniformly bounded trajectory. Theorem~\ref{BF-thm:maincurvature} and Proposition~\ref{BF-prop:boundedflow} therefore apply. In particular they cover every pair and every bounded initialization in the conjecture after \cite[Theorem~1.9]{BF-BPS}.
\end{proof}

The conjecture thus follows independently of the next argument,
which proves confinement when the terminal support is unbounded.

\subsection{Confinement for unbounded terminal laws}
\label{BF-sec:unbounded}
Strong convexity gives a rate only inside a bounded latent region.
The marginal constraint alone does not give such a region. We first
show that a Bass minimizer is bounded by interior coercivity of its
heat-smoothed potential. The flow then stays close to that minimizer
in $L^2$, which makes its latent laws tight. Compactness of the
corresponding heat potentials turns tightness into a uniform inward
drift outside one ball.

Let $\nu\in\mathcal P_2(\R^d)$ have full affine span, and set
$C=\overline{\operatorname{conv}}(\supp\nu)$. Throughout this section, $K$ is a fixed compact convex subset of $\operatorname{int}C$.

\begin{lemma}[Uniform compactness after Gaussian smoothing]\label{BF-lem:heatcompact}
Let $\mathcal A$ be a uniformly tight family of probability laws $\alpha$. For each $\alpha\in\mathcal A$, suppose a finite convex potential $v_\alpha$, normalized by $v_\alpha(0)=0$, satisfies
\[
 (\nabla v_\alpha)_\#(\alpha*\gamma_d)=\nu.
\]
Then $u_\alpha=v_\alpha*\gamma_d$ is finite and smooth on all of $\R^d$. The family of normalized potentials $v_\alpha$ is precompact locally uniformly; every limit retains a transport to $\nu$ from a Gaussian convolution of a subsequential weak limit of $\alpha$. The heat-smoothed functions $u_\alpha$ are likewise precompact locally uniformly, and convergence of the potentials gives local uniform convergence of $u_\alpha$ and their first derivatives.
\end{lemma}
\begin{proof}
Choose $K_0<\infty$ and $p>0$ with $\alpha(B_{K_0})\ge p$ for every $\alpha$. Their convolution densities are bounded below by
\begin{equation}\label{BF-eq:sharpminorant}
 \ell(w)=p(2\pi)^{-d/2}
 \exp\left(-\frac{(|w|+K_0)^2}{2}\right).
\end{equation}
Consequently
\begin{equation}\label{BF-eq:ellL2}
 \int|\nabla v_\alpha|^2\ell\le\int|y|^2\,d\nu.
\end{equation}
The cap argument in Lemma~\ref{BF-lem:positiveK} gives local gradient bounds and local compactness of the normalized potentials. Tightness of $\alpha$ and of the fixed target gives tightness of the associated transport plans. Passing to the limit in subgradient inequalities shows that each limiting plan is concentrated on the limiting subgradient. Its source is a Gaussian convolution, hence absolutely continuous; the limit is therefore a gradient map to $\nu$.

For $z$ in a fixed compact set, both
\[
 \frac{\varphi(w-z)^2}{\ell(w)}
 \quad\text{and}\quad
 \frac{|w|^{2j}\varphi(w-z)^2}{\ell(w)},\qquad j\ge0,
\]
are integrable in $w$, with uniformly vanishing tails; here $\varphi$ denotes the standard Gaussian density. Their exponents have the form $-|w|^2/2+O(|w|)$ uniformly for bounded $z$.

For almost every $w$, convexity and $v_\alpha(0)=0$ imply
\[
 |v_\alpha(w)|\le |w|\bigl(|\nabla v_\alpha(w)|+|p_\alpha|\bigr),
 \qquad p_\alpha\in\partial v_\alpha(0),
\]
and the $p_\alpha$ are uniformly bounded. Cauchy--Schwarz with \eqref{BF-eq:ellL2} gives finite Gaussian convolutions and uniform tails, including for all derivatives of the Gaussian density. It proves smoothness, identifies $\nabla u_\alpha=\nabla v_\alpha*\gamma_d$, and upgrades the local convergence of $v_\alpha$ and their gradients to the asserted convergence of the heat convolutions.
\end{proof}

\begin{lemma}[Interior coercivity]\label{BF-lem:coercivity}
For each potential $v$ in Lemma~\ref{BF-lem:heatcompact}, its conjugate $v^*$ is finite on $\operatorname{int}C$. The function $u=v*\gamma_d$ has positive definite Hessian everywhere. For each $x\in\operatorname{int}C$, the function $u(z)-x\cdot z$ is coercive and has a unique minimizer. The inverse image $(\nabla u)^{-1}(K)$ is bounded.
\end{lemma}
\begin{proof}
For $\nu$-almost every $y$ there is a finite $w$ with $y=\nabla v(w)$, and Fenchel equality gives $v^*(y)<\infty$. The convex domain of $v^*$ contains the convex hull of this full-measure set. The interior of that convex hull is $\operatorname{int}C$, so $v^*$ is finite and locally bounded there.

Choose $\delta>0$ such that $K+2\delta B_1\subset\operatorname{int}C$, and put $A_v=\sup_{x\in K+2\delta B_1}v^*(x)<\infty$. Jensen and Fenchel give
\begin{equation}\label{BF-eq:coercive}
 u(z)\ge v(z)\ge h_K(z)+2\delta|z|-A_v,
 \qquad h_K(z)=\sup_{x\in K}x\cdot z.
\end{equation}
This proves coercivity uniformly over $x\in K$ and bounds their minimizers.

The Hessian of $u$ is the convolution of the positive Hessian measure of $v$ with the strictly positive Gaussian density. If it vanished in a direction at one point, $v$ would be affine with constant slope in that direction by the last argument of Lemma~\ref{BF-lem:positiveK}. This would contradict full affine span of $\nu$. Thus $u$ is strictly convex with positive definite Hessian, which proves uniqueness.
\end{proof}

\begin{proposition}[A bounded Bass measure]\label{BF-prop:boundedmin}
Under the hypotheses of Theorem~\ref{BF-thm:flow}, a Bass measure is bounded. On $\mathsf H=L^2(\sigma(X);\R^d)$ there is a bounded minimizer of $V$ with any prescribed mean.
\end{proposition}
\begin{proof}
The Bass representation theorem \cite[Theorem~1.3]{CM-Bass} gives a probability law $\alpha$ and a finite convex $v$ such that
\[
 (\nabla v)_\#(\alpha*\gamma_d)=\nu,
 \qquad (\nabla v*\gamma_d)_\#\alpha=\mu.
\]
No moment assumption on $\alpha$ is needed for that representation. A single probability law is tight, so Lemmas~\ref{BF-lem:heatcompact}--\ref{BF-lem:coercivity} apply. They show that $\alpha$ is supported on the bounded inverse image $(\nabla u)^{-1}(K)$, where $u=v*\gamma_d$. In particular $\alpha\in\mathcal P_2$.

Set $Z^*=(\nabla u)^{-1}(X)$. Its law is $\alpha$ and $DV(Z^*)=0$ by \eqref{BF-eq:gradient}; convexity makes it a global minimizer. Translating $Z^*$ and translating the source variable of $v$ gives a minimizer with any prescribed mean. Uniqueness with that mean is also the general uniqueness assertion of \cite[Lemma~2.1(iv)]{BF-BPS}; it applies once the Bass minimizer is in $\mathcal P_2$.
\end{proof}

\begin{lemma}[Uniform inward drift]\label{BF-lem:inward}
For any uniformly tight family $\mathcal A$ as in Lemma~\ref{BF-lem:heatcompact}, there are constants $R_0<\infty$ and $\delta_0>0$ such that
\begin{equation}\label{BF-eq:inward}
 z\cdot(\nabla u_\alpha(z)-x)\ge\delta_0|z|
 \quad(|z|\ge R_0,\ x\in K,\ \alpha\in\mathcal A).
\end{equation}
\end{lemma}
\begin{proof}
Use $\delta$ from Lemma~\ref{BF-lem:coercivity}. By \eqref{BF-eq:coercive}, for each limit heat potential $u$ there is $r_u$ such that
\[
 \frac{u(r_ue)-u(0)}{r_u}\ge h_K(e)+\delta
 \quad\text{for every unit }e.
\]
This remains true with $\delta/2$ for every heat potential in a sufficiently small local-uniform neighborhood of $u$. The family of all limit heat potentials is compact by Lemma~\ref{BF-lem:heatcompact}; a finite cover therefore gives finitely many such radii. Set $R_0$ to their maximum. Secant slopes of a convex function on a ray are nondecreasing, and its radial derivative dominates every preceding secant slope. Hence
$e\cdot\nabla u_\alpha(re)\ge h_K(e)+\delta/2$ for $r\ge R_0$. This gives \eqref{BF-eq:inward} with $\delta_0=\delta/2$.
\end{proof}

\begin{proof}[Proof of Theorem~\ref{BF-thm:flow}]
A finite continuous convex functional on a Hilbert space has a unique subgradient flow from every starting point; here its gradient is continuous and uniformly bounded in $L^2$ by Lemma~\ref{BF-lem:gradient}, so the flow is global and locally absolutely continuous. This is the same Hilbert-space flow used in \cite{BF-BPS}.

Choose the bounded minimizer from Proposition~\ref{BF-prop:boundedmin} with $\E Z^*=\E Z_0$. Convexity gives
\[
 \frac d{dt}\frac12\|Z_t-Z^*\|_2^2
 =-\langle DV(Z_t),Z_t-Z^*\rangle\le0.
\]
Consequently the laws $\alpha_t=\law(Z_t)$ have uniformly bounded second moments and are uniformly tight. Apply Lemma~\ref{BF-lem:inward} to their Brenier potentials. The flow has versions with almost surely absolutely continuous sample paths on every finite time interval, by Fubini and its $L^2$ absolute continuity. Along these paths,
\[
 \frac d{dt}|Z_t|^2=-2Z_t\cdot(\nabla u_t(Z_t)-X)
 \le-2\delta_0|Z_t|\quad\text{whenever }|Z_t|\ge R_0.
\]
The positive-part chain rule proves
\[
 |Z_t|\le\max\{R_0,\|Z_0\|_\infty\}
 \quad\text{for all }t\ge0\text{ almost surely}.
\]
Mean preservation follows from $\E DV=0$. The complete trajectory and $Z^*$ are now in one bounded fixed-mean slice. Proposition~\ref{BF-prop:boundedflow} proves \eqref{BF-eq:rates}.
\end{proof}

\begin{remark}
The interior condition governs confinement of the latent Bass variable. The strong-convexity theorem itself needs neither a convex-order pair nor irreducibility. The positive rate can deteriorate as the initial support approaches the boundary or as targets approach a lower-dimensional law. Each of these dependencies is separate from the dimension in which the proof operates.
\end{remark}

\subsection{Iterations and stability on a bounded latent region}
\label{BF-sec:algorithms}
Fix a radius $R$ containing a normalized bounded minimizer $Z^*$, and put
\[
 \mathsf K_{R,m}=\{Z\in\mathsf H:|Z|\le R\text{ almost surely},\ \E Z=m\},
 \qquad m=\E Z^*.
\]
This is a closed convex subset of the Hilbert space. Let $\kappa$ be its strong-convexity constant.

\begin{proposition}[Proximal iteration]\label{BF-prop:proximal}
For $h>0$ and $Z_0\in\mathsf K_{R,m}$, define
\begin{equation}\label{BF-eq:proxiteration}
 Z_{n+1}=\mathop{\rm argmin}_{Z\in\mathsf K_{R,m}}
 \left\{V(Z)+\frac1{2h}\|Z-Z_n\|_2^2\right\}.
\end{equation}
The minimizer exists uniquely and
\begin{equation}\label{BF-eq:proxrate}
 \|Z_n-Z^*\|_2\le(1+h\kappa)^{-n}\|Z_0-Z^*\|_2.
\end{equation}
If the optimality condition at step $n$ has an additive error $e_n\in\mathsf H$, its distance satisfies the corresponding estimate
\[
 \|Z_{n+1}-Z^*\|_2
 \le\frac{\|Z_n-Z^*\|_2+h\|e_n\|_2}{1+h\kappa}.
\]
\end{proposition}
\begin{proof}
Weak lower semicontinuity, coercivity of the quadratic term, and strict convexity give existence and uniqueness. The optimality inclusion is
$(Z_n-Z_{n+1})/h\in DV(Z_{n+1})+N_{\mathsf K_{R,m}}(Z_{n+1})$, where $N$ is the normal cone. Compare with $0\in DV(Z^*)+N_{\mathsf K_{R,m}}(Z^*)$. Strong monotonicity of $DV$ and monotonicity of the normal cone give $(1+h\kappa)\|Z_{n+1}-Z^*\|_2^2 \le\langle Z_n-Z^*,Z_{n+1}-Z^*\rangle$. Cauchy--Schwarz proves \eqref{BF-eq:proxrate}; including $e_n$ proves the last estimate.
\end{proof}

\begin{proposition}[Lipschitz dependence on the initial marginal]\label{BF-prop:initialstability}
Let $(\mu_i,\nu)$, $i=0,1$, have Bass minimizers with latent variables normalized to the same mean. Assume these variables are bounded by $R$. If $\alpha_i$ are their Bass measures, then
\begin{equation}\label{BF-eq:bassstability}
 W_2(\alpha_0,\alpha_1)\le\kappa(R,\gamma_d,\nu)^{-1}W_2(\mu_0,\mu_1).
\end{equation}
More precisely, on any coupling of $X_0\sim\mu_0$ and $X_1\sim\mu_1$, their normalized optimal maps $Z_i=z_i(X_i)$ satisfy
$\|Z_1-Z_0\|_2\le\kappa^{-1}\|X_1-X_0\|_2$.
\end{proposition}
\begin{proof}
On the common probability space, let $A(Z)=DF_{\nu,\gamma_d}(Z)$. Optimality gives $A(Z_i)=X_i$; this remains true after enlarging the base sigma-field, because $A(Z_i)=(T_{Z_i}*\gamma_d)(Z_i)$. Strong monotonicity yields
\[
 \kappa\|Z_1-Z_0\|_2^2
 \le\langle A(Z_1)-A(Z_0),Z_1-Z_0\rangle
 =\langle X_1-X_0,Z_1-Z_0\rangle.
\]
Take an optimal coupling of the initial laws to obtain \eqref{BF-eq:bassstability}.
\end{proof}

The common latent bound is part of Proposition~\ref{BF-prop:initialstability}. It is not inferred solely from a common bound on the initial support: the latent variables may diverge as an irreducible pair approaches a reducible one.

There is also a measure-space interpretation of the curvature theorem. For the fixed noise law $q$, define
\[
 \mathcal V_q(\alpha)=\MCov(\alpha*q,\nu)-\MCov(\alpha,\mu).
\]

\begin{proposition}[Strong convexity on generalized geodesics]\label{BF-prop:geodesic}
Suppose $(X,Z_i)$ are optimal covariance couplings of $(\mu,\alpha_i)$, glued on a common space, and $|Z_i|\le R$. Put $\alpha_t=\law((1-t)Z_0+tZ_1)$. Then
\begin{align}\label{BF-eq:geodesic}
 \mathcal V_q(\alpha_t)
 &\le(1-t)\mathcal V_q(\alpha_0)+t\mathcal V_q(\alpha_1)\notag\\
 &\quad-\frac\kappa2t(1-t)\Var_2(Z_1-Z_0).
\end{align}
If the two latent means agree, the final variance can be replaced by $W_2^2(\alpha_0,\alpha_1)$.
\end{proposition}
\begin{proof}
The cyclic monotonicity inequalities for $(X,Z_0)$ and $(X,Z_1)$ add with weights $1-t,t$. Thus $(X,(1-t)Z_0+tZ_1)$ is also an optimal covariance coupling. Its covariance is affine in $t$, so Theorem~\ref{BF-thm:maincurvature} gives \eqref{BF-eq:geodesic}. Equal means give $\Var_2(Z_1-Z_0)=\E|Z_1-Z_0|^2\ge W_2^2(\alpha_0,\alpha_1)$.
\end{proof}

For Gaussian noise, ordinary convexity along these generalized geodesics is already part of the Bass-functional theory \cite{Bass}. The positive modulus in Proposition~\ref{BF-prop:geodesic} holds on the
specified bounded generalized geodesics and for the full noise class
of Theorem~\ref{BF-thm:maincurvature}.

\subsection{Exact marginals throughout calibration}
\label{BF-sec:couplings}
The latent convergence rate becomes a rate for the complete reference
coupling when the target has curvature. Both marginals are already
exact along the flow. The quantity tending to zero is the conditional
mean residual, equal to the lifted gradient. The same-Gaussian estimate
below controls these two errors with a single transport comparison.

Assume in this subsection that
\begin{equation}\label{BF-eq:uniformlc}
 \nu(dy)=c_Ve^{-V_\nu(y)}\,dy,\qquad
 V_\nu(y)-\frac{|y|^2}{2L^2}\ \text{is convex},\quad L>0.
\end{equation}
Extended-valued convex potentials, and thus convex support constraints, are allowed by approximation. Only this section uses the curvature assumption.

\begin{lemma}[Common-noise stability of optimal maps]\label{BF-lem:mapsmooth}
For arbitrary $Z_0,Z_1\in L^2$ on the same base space, set
\[
 W_i=Z_i+G,\qquad Y_i=T_{Z_i}(W_i),
\]
using the same independent Gaussian $G$. Then
\begin{equation}\label{BF-eq:mapcocoercive}
 \E[(Z_1-Z_0)\cdot(Y_1-Y_0)]\ge\frac1L\E|Y_1-Y_0|^2.
\end{equation}
In particular,
\begin{equation}\label{BF-eq:maplip}
 \|Y_1-Y_0\|_2\le L\|Z_1-Z_0\|_2,
 \qquad
 \|DV(Z_1)-DV(Z_0)\|_2\le L\|Z_1-Z_0\|_2.
\end{equation}
\end{lemma}
\begin{proof}
Lemma~\ref{CM-lem:mixture} gives $\Lip(T_{Z_i})\le L$. Apply the two-reference deficit calculation in Corollary~\ref{R15-cor:variation} with a one-point initial law of mean $\int y\,d\nu$ and driver laws $\law(W_i)$. These Gaussian mixtures have upper log-density curvature one. Centering them changes neither the maps' Lipschitz constants nor the comparison, because $\E(Y_1-Y_0)=0$. Before Cauchy--Schwarz, the sum of the two deficits gives
\[
 L^{-1}\E|Y_1-Y_0|^2
 \le\E[(W_1-W_0)\cdot(Y_1-Y_0)]
 =\E[(Z_1-Z_0)\cdot(Y_1-Y_0)].
\]
This proves~\eqref{BF-eq:mapcocoercive} and the first bound in~\eqref{BF-eq:maplip}. Conditional expectation proves the gradient bound and also
\begin{equation}\label{BF-eq:gradientcocoercive}
 \langle DV(Z_1)-DV(Z_0),Z_1-Z_0\rangle
 \ge L^{-1}\|DV(Z_1)-DV(Z_0)\|_2^2.
\end{equation}
\end{proof}

\begin{corollary}[Projected gradient iteration]\label{BF-cor:projected}
On $\mathsf K_{R,m}$ from Section~\ref{BF-sec:algorithms}, the iteration
\[
 Z_{n+1}=\operatorname{Proj}_{\mathsf K_{R,m}}\bigl(Z_n-L^{-1}DV(Z_n)\bigr)
\]
satisfies
\begin{equation}\label{BF-eq:projectedrate}
 \|Z_n-Z^*\|_2\le(1-\kappa/L)^{n/2}\|Z_0-Z^*\|_2.
\end{equation}
\end{corollary}
\begin{proof}
Combine \eqref{BF-eq:gradientcocoercive} with strong monotonicity to obtain
\[
 \|(Z_1-Z_0)-L^{-1}(DV(Z_1)-DV(Z_0))\|_2^2
 \le(1-\kappa/L)\|Z_1-Z_0\|_2^2.
\]
Projection is nonexpansive and $Z^*$ is fixed by the iteration. This proves \eqref{BF-eq:projectedrate}.
\end{proof}

\begin{theorem}[Exponential convergence of the reference coupling]\label{BF-thm:couplingrate}
Under Theorem~\ref{BF-thm:flow} and \eqref{BF-eq:uniformlc}, realize
\[
 Y_t=T_{Z_t}(Z_t+G),\qquad Y^*=T_{Z^*}(Z^*+G)
\]
using a single Gaussian independent of $X$. Then every $Y_t$ has law $\nu$, $\E[Y^*\mid X]=X$, and
\begin{align}\label{BF-eq:couplingrate}
 \|Y_t-Y^*\|_2&\le Le^{-\kappa t}\|Z_0-Z^*\|_2,\\
 \|\E[Y_t\mid X]-X\|_2&\le Le^{-\kappa t}\|Z_0-Z^*\|_2.\notag
\end{align}
Consequently the complete couplings $(X,Y_t)$ converge in $W_2$ at the first rate in \eqref{BF-eq:couplingrate}. The analogous assertions hold for the iterations in Propositions~\ref{BF-prop:proximal} and Corollary~\ref{BF-cor:projected}, with their respective geometric factors.
\end{theorem}
\begin{proof}
The exact target law follows from the defining transport constraint. Stationarity gives the conditional mean at the limit. Lemma~\ref{BF-lem:mapsmooth} and \eqref{BF-eq:rates} give the first estimate. For the second, use
$\E[Y_t\mid X]-X=DV(Z_t)$ and the Lipschitz gradient estimate relative to $Z^*$. The common-space coupling of the two pairs has zero error in its first coordinate, so it gives the stated $W_2$ bound.
\end{proof}

\begin{corollary}[Stability of the optimal martingale coupling]\label{BF-cor:martingalestability}
Under Proposition~\ref{BF-prop:initialstability} and \eqref{BF-eq:uniformlc}, let $\pi_i^*$ denote the canonical martingale couplings of $(\mu_i,\nu)$. Then, with the product Euclidean metric,
\[
 W_2(\pi_0^*,\pi_1^*)
 \le\sqrt{1+L^2/\kappa^2}\,W_2(\mu_0,\mu_1).
\]
\end{corollary}
\begin{proof}
Take an optimal coupling of $X_0,X_1$, realize their normalized latent maps on it, and use a common independent Gaussian. Proposition~\ref{BF-prop:initialstability} controls the latent difference, and Lemma~\ref{BF-lem:mapsmooth} controls the terminal difference. Add the squared initial and terminal differences. Each limiting coupling is the Bass coupling, hence is the canonical martingale Benamou--Brenier optimizer by \cite{CM-Bass}.
\end{proof}

\begin{corollary}[A whole signing law with exponentially calibrated reference]\label{BF-cor:signs}
Let a finite signing law satisfy
\[
 S=(A\sigma,\sigma)\preceq_{\mathrm{cx}}N(0,\Gamma),\qquad \Gamma\succ0.
\]
For every $c>1$, its law can be coupled through the above flow to references $R_t$ and $R^*$ such that
\[
 R_t\sim R^*\sim N(0,c^2\Gamma),\qquad
 \E[R^*\mid\sigma]=S,
\]
and, for constants $C,\kappa>0$ depending on this fixed problem,
\[
 \E\big|\Gamma^{-1/2}(R_t-R^*)\big|^2
 +\E\big|\Gamma^{-1/2}(\E[R_t\mid\sigma]-S)\big|^2
 \le Ce^{-2\kappa t}.
\]
The complete signing distribution and the complete Gaussian reference marginal are exact at every stage. Block independence in the Gaussian marginal is retained when $\Gamma$ is block diagonal.
\end{corollary}
\begin{proof}
Whiten $S$ by $\Gamma^{-1/2}$. Its finite law is compactly supported. Corollary~\ref{BF-cor:inflation} gives irreducibility for the target $N(0,c^2I)$, and Theorem~\ref{BF-thm:couplingrate} applies with $L=c$. Return to the original coordinates. The sign vector is part of $S$, so conditioning on $S$ is equivalent to conditioning on $\sigma$.
\end{proof}

The rates measure iterations of exact variational subproblems. Their
implementation cost includes representing the initial law and solving
the transport, proximal or projection subproblem; these rates alone
are not bit-complexity bounds.

The calibration changes the coupling while retaining both endpoint
laws. The next part changes the reference itself. Its score and
support determine which hard constraints can be attained and which
geometric region survives rounding.


\Needspace{12\baselineskip}
\part{Reference geometry and sharp comparisons}
\label{R6-part:geometry}
The retained geometry depends on the full initial reference. Optimizing its support and score leads to directional energy and spectral problems; compact group symmetry gives exact comparisons for specified discrete geometries.

\section{The geometry retained by the full score law}
\label{R11-sec:scores}
Integrated quantiles determine scalar support, and conditional averaging assembles a retained convex body. Even all score moments can agree while sharp radii differ. The realization and stability theorems identify the geometry determined by the complete score law; equal Fisher matrices and three-coordinate marginals permit order-$\sqrt m$ width differences.

\subsection{The geometry retained by the full score law}\label{S-sec:score}
Let $f$ be a probability density whose zero extension has an integrable
weak gradient. Its score is
$Z_f=\nabla f(X)/f(X)$, where $X\sim f$; its value on $\{f=0\}$ is
irrelevant. The score is integrable and centered. For a nonzero vector
$a$, the marginal density $g$ of $\langle a,X\rangle$ satisfies
\begin{equation}\label{S-eq:intro-score}
 \E[Z_f\mid\langle a,X\rangle=z]=a\,g'(z)/g(z).
\end{equation}
This is a vector conditional expectation. It constrains all upper tails
of the scalar score under the same coupling.

For a centered integrable scalar law $\nu$, let $Q_\nu^\downarrow$ be
its decreasing quantile and put
\begin{equation}\label{S-eq:radius}
 L_\nu(t)=\int_0^tQ_\nu^\downarrow(s)\,ds,\qquad
 r(\nu)=\frac12\int_0^1\frac{dt}{L_\nu(t)}.
\end{equation}
We allow the value $+\infty$. In particular $r(\delta_0)=+\infty$.
For a centered integrable vector law $\mu$, define
\begin{equation}\label{S-eq:r4-value}
 \mathfrak b_\mu(a)=\inf\{r(\nu):\law(aS)\cx\mu,\ S\sim\nu\},
 \qquad
 \mathcal K_\mu=\bigcap_{a\ne0}
       \{x:\langle a,x\rangle\le\mathfrak b_\mu(a)\}.
\end{equation}
For a closed convex set $L$, write
$h_L(a)=\sup_{x\in L}\langle a,x\rangle$ for its support function.
The constraint in \eqref{S-eq:r4-value} is equivalent to the existence of one martingale
coupling $\E[Z\mid S]=aS$, with $Z\sim\mu$. Reflection of $S$ shows
that $\mathfrak b_\mu$ is even; rescaling $S$ makes it positively
homogeneous in $a$. Thus $\mathcal K_\mu$ is a closed centrally
symmetric convex set. When $0<\E\|Z\|<\infty$, it contains the ball
of radius $1/\E\|Z\|$: indeed
$L_\nu\le\E|S|/2$, whence $r(\nu)\ge1/\E|S|$, and the martingale
constraint gives $\|a\|\E|S|\le\E\|Z\|$.

Write $W_p$ for the $p$-Wasserstein distance. Convergence in $W_p$
means weak convergence together with convergence of the $p$th absolute
moments, where $p\ge1$.

\begin{theorem}[Retained body under a first-moment score hypothesis]
\label{S-thm:main}
Let $K\subset\R^m$ be bounded, open, convex and centrally symmetric
about $p$. Let $f$ be a probability density supported in $K$, with mean
$p$ and zero extension in $W^{1,1}(\R^m)$. Suppose that each vector
$v_j$ has a centered auxiliary $q_j$ with $H_f(v_j;q_j)\ge1$.
Then, with $\mu_f=\law(Z_f)$,
\begin{equation}\label{S-eq:terminal}
 p+\operatorname{int}\mathcal K_{\mu_f}
 \subset\calS_{v_n}\cdots\calS_{v_1}K.
\end{equation}
The inclusion holds for every ordering of the vectors. Each point on
the left admits backward sign recovery into $K$.
Every finite minimum in \eqref{S-eq:r4-value} has a unique scalar
optimizer. When $\mu_f$ is symmetric, its finite dyadic optimizers converge in
$W_1$, and in every $W_p$ supported by a finite $p$th moment of $\mu_f$.
\end{theorem}

The proof separates the scalar support constraint from the score
contraction. The first uses only an integrable derivative. The second
preserves the same first-moment hypothesis through the complete lift.
Finite Fisher information is available when the score is square
integrable, in which case $J(f)=\E Z_fZ_f^{\mathsf T}$; it is not
needed for \eqref{S-eq:terminal}.

\subsubsection{The interval functional and its constrained extremizer}\label{S-sec:interval}

The elementary variational formula
\begin{equation}\label{S-eq:quantile-variational}
 L_\nu(t)=\inf_{a\in\R}\{ta+\E(Z-a)_+\}
 =\sup_{\substack{0\le U\le1\\\E U=t}}\E(ZU)
\end{equation}
holds with randomized selection at atoms. The last supremum allows $U$
to depend on $Z$ and auxiliary randomness. Selecting the upper $t$-tail
proves equality in both expressions.

\begin{lemma}\label{S-lem:radius-properties}
For centered nondegenerate integrable laws,
\begin{align}
 r(\law(cZ))&=|c|^{-1}r(\law Z)\quad(c\ne0),\label{S-eq:scale}\\
 \nu\cx\eta&\ \Longrightarrow\ r(\nu)\ge r(\eta),\label{S-eq:order}\\
 r\!\left(\int\nu_s\,d\lambda(s)\right)
 &\le\int r(\nu_s)\,d\lambda(s).\label{S-eq:mixture}
\end{align}
The radius is unchanged by reflection.
\end{lemma}
\begin{proof}
For $c>0$, the quantiles scale by $c$. Centering gives
$L_{\law(-Z)}(t)=L_{\law Z}(1-t)$, proving the remaining part of
\eqref{S-eq:scale} and reflection invariance. Convex order and
\eqref{S-eq:quantile-variational} give $L_\nu\le L_\eta$, hence
\eqref{S-eq:order}. If $\bar\nu=\int\nu_s\,d\lambda(s)$, the infimum
in \eqref{S-eq:quantile-variational} gives
\[
 L_{\bar\nu}(t)\ge\int L_{\nu_s}(t)\,d\lambda(s).
\]
Convexity of $x\mapsto1/x$ and Tonelli's theorem now prove
\eqref{S-eq:mixture}. The argument also applies to infinite values.
\end{proof}

\begin{theorem}[Sharp interval support]\label{S-thm:interval}
Let $g$ be a probability density supported in a finite interval of
length $2d$. Assume its zero extension is in $W^{1,1}(\R)$, and set
$S_g=g'/g$ on $\{g>0\}$. If $\law(S_g(X))\cx\nu$, $X\sim g$, then
\begin{equation}\label{S-eq:interval}
 d\ge r(\nu).
\end{equation}
Conversely, for every centered nondegenerate integrable law $\nu$ with
$r(\nu)<\infty$, there is a log-concave density on $(-r(\nu),r(\nu))$
whose score law is exactly $\nu$. It is symmetric when $\nu$ is symmetric.
\end{theorem}
\begin{proof}
First suppose $g$ is positive on the interior of its supporting interval.
Write $x(t)=F_g^{-1}(t)$ and $q(t)=g(x(t))$. The chain rule gives
\[
 x'(t)=\frac1{q(t)},\qquad q'(t)=S_g(x(t)).
\]
Zero trace gives $q(0)=q(1)=0$. Hence
\[
 q(t)=\int_0^t S_g(x(s))\,ds
 \le L_{\law(S_g(X))}(t)\le L_\nu(t).
\]
Integrating $x'=1/q$ proves \eqref{S-eq:interval}.
For a general absolutely continuous density, the same calculation holds
on its positive components. Passing from one component to the next adds
a nonnegative jump to the quantile function. The function $q$ has zero
trace at the intervening probability values and
$q'=S_g\circ x$ almost everywhere. Since
$\int_0^1|S_g(x(t))|dt=\int|g'|<\infty$, summing the component identities
is legitimate. The total span is at least $\int_0^1dt/q(t)$, which gives
the same bound.

For the converse, put $r=r(\nu)$ and define
\begin{equation}\label{S-eq:reconstruction}
 x(t)=-r+\int_0^t\frac{ds}{L_\nu(s)},\qquad
 g(x(t))=L_\nu(t),\qquad 0<t<1.
\end{equation}
The map $x$ is increasing onto $(-r,r)$. Change of variables gives
$g(x(t))x'(t)=1$, so $g$ has mass one. Also
\[
 \frac{d}{dx}\log g(x(t))=L'_\nu(t)=Q^\downarrow_\nu(t)
\]
almost everywhere. This derivative is nonincreasing; therefore $g$ is
log-concave and has the required score law. Its zero trace and
$\int|g'|=\int_0^1|Q^\downarrow_\nu|<\infty$ justify the weak derivative.
Reflection symmetry of $L_\nu$ gives symmetry of $g$ when $\nu$ is symmetric.
\end{proof}

For $g=e^{-V}$, the distribution of $V'$ under $g$ is its
one-dimensional moment measure.  Thus \eqref{S-eq:reconstruction} is an
explicit one-dimensional case of the moment-measure correspondence
\cite{CEK}.  The support inequality and the mixture argument above are
the forms needed for the section theorem.

\paragraph{Compactness in the moment topology}
\label{S-sec:r7-compactness}

The scalar constraint implies uniform integrability. This
removes the second-moment hypothesis from existence and gives stronger
convergence whenever higher source moments are available.

\begin{lemma}[Compactness of the admissible laws]\label{S-lem:r7-compact}
Let $\mu$ be a centered integrable law on $\R^m$ and $a\ne0$. The set
\[
 \mathcal A_\mu(a)=\{\nu:\law(aS)\cx\mu,\ S\sim\nu\}
\]
is compact in $W_1$. If $\int\|z\|^p\,d\mu(z)<\infty$ for some
$1\le p<\infty$, it is compact in $W_p$.
Consequently, every finite common-coupling minimum has a unique scalar
optimizer under the first-moment hypothesis alone. Every minimizing sequence converges to it in each such $W_p$.
For symmetric $\mu$, the same holds for the dyadic optimizing laws
defined in Section~\ref{S-sec:r6-dyadic}.
\end{lemma}
\begin{proof}
Put $Y=\langle a,Z\rangle/\|a\|^2$, $Z\sim\mu$. Projecting a feasible
martingale gives $\nu\cx\law(Y)$. Conditional Jensen therefore gives,
for every $R>0$,
\begin{equation}\label{S-eq:r7-ui}
 \int_{|s|>R}|s|^p\,d\nu(s)
 \le2^p\int (|s|-R/2)_+^p\,d\nu(s)
 \le2^p\E\bigl[|Y|^p\mathbf1_{\{|Y|>R/2\}}\bigr].
\end{equation}
The right side tends to zero, uniformly over all feasible laws.
This proves tightness and uniform integrability in the required
moment order. The feasible set is closed under $W_1$ convergence,
since convex Lipschitz tests pass to the limit and determine convex
order. It is therefore compact in $W_p$.

The integrated quantiles converge uniformly under $W_1$, so Fatou's
lemma makes $r$ lower semicontinuous. It attains every finite minimum.
The strict mixture inequality for $r$ gives uniqueness. Compactness
then identifies every subsequential limit of a minimizing sequence.
The dyadic lower sums have the same increasing-limit argument, now
in this compact space. Thus their optimizing laws have the same
$W_p$ limit.
\end{proof}

In particular, a finite second moment gives $W_2$ convergence of the
scalar extremizers. Their second moments, and hence the total
conditional covariance
\[
 \E\Cov(Z\mid S)=\Cov(Z)-aa^{\mathsf T}\E S^2,
\]
converge as well. 

\paragraph{Strict convexity and the optimizing law}

The mixture inequality has an explicit strict remainder. If
$L_i=L_{\nu_i}$ and $0<\theta<1$, the upper-tail variational formula
and the reciprocal identity give
\begin{align}
 r((1-\theta)\nu_0+\theta\nu_1)
 &\le(1-\theta)r(\nu_0)+\theta r(\nu_1)\notag\\[-2pt]
 &\quad-\frac{\theta(1-\theta)}2
 \int_0^1\frac{(L_1-L_0)^2}
 {L_0L_1((1-\theta)L_0+\theta L_1)}\,dt.
 \label{S-eq:r5-strict-mixture}
\end{align}
The underlying pointwise identity is
\[
 \frac{1-\theta}{A}+\frac\theta B
 -\frac1{(1-\theta)A+\theta B}
 =\frac{\theta(1-\theta)(A-B)^2}
 {AB((1-\theta)A+\theta B)}.
\]
Distinct laws have different continuous integrated quantiles on a set
of positive measure. Therefore $r$ is strictly convex on centered laws
with finite radius.

\begin{theorem}[Unique extremizer and stability]\label{S-thm:r5-stability}
Let $\mu$ be centered and integrable, let $a\ne0$, and
suppose $\mathfrak b_\mu(a)<\infty$. There is a unique admissible scalar
law $\nu_a$ attaining the infimum in \eqref{S-eq:r4-value}.
Write $L_a=L_{\nu_a}$. Every admissible law $\nu$ of finite radius
satisfies
\begin{equation}\label{S-eq:r5-stability}
 r(\nu)-\mathfrak b_\mu(a)
 \ge\frac12\int_0^1
       \frac{(L_\nu(t)-L_a(t))^2}{L_a(t)^2L_\nu(t)}\,dt.
\end{equation}
If the covariance is finite, $\Sigma=\Cov(\mu)\succ0$, and
$\sigma_a^2=(a^{\mathsf T}\Sigma^{-1}a)^{-1}$, then also
\begin{equation}\label{S-eq:r5-variance-stability}
 r(\nu)-\mathfrak b_\mu(a)
 \ge\frac1{2\sigma_a^3}\int_0^1
 \frac{(L_\nu(t)-L_a(t))^2}{[t(1-t)]^{3/2}}\,dt.
\end{equation}
Every minimizing sequence converges to $\nu_a$ in $W_p$ for every
$1\le p<\infty$ for which the source has a finite $p$th moment.
When $\mu$ is symmetric, $\nu_a$ is symmetric.
\end{theorem}
\begin{proof}
Lemma~\ref{S-lem:r7-compact} gives compactness in every available
moment topology. Uniform convergence of integrated quantiles under
$W_1$ and Fatou's lemma give attainment. Convex order is preserved
under mixtures, so \eqref{S-eq:r5-strict-mixture} gives uniqueness.

Apply \eqref{S-eq:r5-strict-mixture} with $\nu_0=\nu_a$ and
$\nu_1=\nu$. The mixture is feasible, whence its radius is at least
$\mathfrak b_\mu(a)$. Subtraction and division by $\theta$ give
\[
 r(\nu)-\mathfrak b_\mu(a)\ge
 \frac{1-\theta}2\int_0^1
 \frac{(L_\nu-L_a)^2}
 {L_aL_\nu((1-\theta)L_a+\theta L_\nu)}\,dt.
\]
Fatou's lemma as $\theta\downarrow0$ proves
\eqref{S-eq:r5-stability}.

The covariance inequality
$aa^{\mathsf T}\E S^2\preceq\Sigma$ implies
$\E S^2\le\sigma_a^2$. For an upper-tail selector $U$ of mass $t$,
Cauchy--Schwarz gives
$\E(SU)\le\sigma_a\sqrt{t(1-t)}$.
This bounds both $L_\nu$ and $L_a$ and proves
\eqref{S-eq:r5-variance-stability}.
Compactness and uniqueness imply convergence of every minimizing
sequence. Reflection preserves the feasible set and objective when
$\mu$ is symmetric, so uniqueness gives symmetry of $\nu_a$.
\end{proof}

\begin{corollary}[An extremizing interval density]\label{S-cor:r5-density}
Among probability densities $g$ with zero extension in $W^{1,1}$ and
$\law(aS_g)\cx\mu$, the least support span is
$2\mathfrak b_\mu(a)$. Up to translation, there is a unique extremizer.
It is log-concave and has score law $\nu_a$; for symmetric $\mu$ it is
symmetric about the midpoint of its support.
\end{corollary}
\begin{proof}
The interval theorem bounds every span below by
$2r(\law S_g)\ge2\mathfrak b_\mu(a)$, and reconstruction from
$\nu_a$ attains this bound. Equality in the interval proof forces
$g(F_g^{-1}(t))=L_a(t)$ almost everywhere and eliminates every gap
between positive components. Integration of $x'=1/L_a$ then fixes
$g$ up to translation.
\end{proof}

\subsubsection{Contraction with an integrable score}
\label{F-sec:score-order}
For a density $g(y,s)$ with integrable horizontal weak gradient, define
its horizontal score by
\[
 S_g=\frac{\nabla_y g(Y,S)}{g(Y,S)},\qquad (Y,S)\sim g.
\]
Its value where $g=0$ is immaterial. Its mean is zero. When its second
moment is finite, $S_g=2\nabla_y\sqrt g/\sqrt g$ almost everywhere
and its horizontal Fisher information matrix is
$J_y(g)=\E S_gS_g^{\mathsf T}$. Fiber rearrangement averages the
score at each level. This is a score formulation of the convex-integral
P\'olya--Szeg\H{o} mechanism; see~\cite{Capriani}.

\begin{theorem}[The score under coordinate rearrangement]\label{F-thm:score-cx}
Let $g_0$ be a density with integrable horizontal weak gradient, and
let $g_j$ be obtained by rearranging the $j$th auxiliary coordinate of
$g_{j-1}$. Every stage satisfies
\begin{equation}\label{F-eq:score-cx}
 S_{g_j}\cx S_{g_{j-1}}\cx S_{g_0}.
\end{equation}
These inequalities also hold after auxiliary folding. Finite Fisher
information, or a finite $r$th score moment for any $r>1$, is sufficient
but is not required.
\end{theorem}
\begin{proof}
First let $g(y,s',t)$ be smooth, and fix a regular positive fiber
level $\lambda$. Write
\[
 L(y,s',\lambda)=|\{t:g(y,s',t)>\lambda\}|,
 \qquad D=-\partial_\lambda L
       =\sum_{t_i:g(t_i)=\lambda}\frac1{|\partial_tg(t_i)|}.
\]
Coarea and implicit differentiation give
\[
 \nabla_yL=\sum_i\frac{\nabla_yg(t_i)}{|\partial_tg(t_i)|}.
\]
Since $g^*(y,s',\pm L/2)=\lambda$, differentiation at fixed
rearranged coordinate yields
\begin{equation}\label{F-eq:score-boundary-average}
 \nabla_y\log g^*(y,s',\pm L/2)
 =\sum_i\omega_i\nabla_y\log g(y,s',t_i),
 \qquad
 \omega_i=\frac{|\partial_tg(t_i)|^{-1}}{D}.
\end{equation}
The mass at $(y,s',\lambda)$, before and after rearrangement, is
$\lambda D\,dy\,ds'\,d\lambda$. Jensen applied to
\eqref{F-eq:score-boundary-average}, followed by integration, proves
\[
 \int g^*\Psi(\nabla_y\log g^*)\le
 \int g\Psi(\nabla_y\log g).
\]
This also identifies the exact quadratic decrease in the regular case:
\begin{equation}\label{F-eq:score-defect}
 J_y(g)-J_y(g^*)
 =\int\lambda D\,
   \operatorname{Cov}_{\omega}
       \bigl(\nabla_y\log g(y,s',t_i)\bigr)
       \,dy\,ds'\,d\lambda.
\end{equation}

For the general case, convolve $g$ with a Gaussian density in all
variables. The smoothed horizontal score is the conditional mean of
the original score given the noisy observation, so every convex score
integral contracts. Apply the smooth rearrangement argument to these
densities, obtaining $g_\varepsilon^*\to g^*$ in $L^1$.

An integrable random vector admits an increasing convex function
$\Theta$ with $\Theta(t)/t\to\infty$ and
$\E\Theta(\|S_g\|_2)<\infty$. One may construct it piecewise linearly
by choosing successive tail thresholds whose first-moment tails are
summable. Consequently
\[
 \sup_\varepsilon\int g_\varepsilon^*
 \Theta\!\left(\frac{\|\nabla_y g_\varepsilon^*\|_2}
                         {g_\varepsilon^*}\right)<\infty.
\]
This also makes the gradients uniformly integrable. For a measurable
set $E$, split at score magnitude $K$ to obtain
\[
 \int_E\|\nabla_y g_\varepsilon^*\|_2
 \le K\int_E g_\varepsilon^*
  +\left(\sup_{t>K}\frac{t}{\Theta(t)}\right)
     \E\Theta(\|S_g\|_2).
\]
Strong $L^1$ convergence of the densities makes the first term uniformly
small on small sets and outside sufficiently large compact sets.
First take $K$ large and then use this convergence. Weak compactness in
$L^1$ gives a weakly convergent subsequence of the gradients.
Distributional differentiation identifies its limit as the horizontal weak gradient of $g^*$.

For a polyhedral convex $\Psi$, its perspective
$(a,z)\mapsto a\Psi(z/a)$ is a supremum of affine functions. Integration
is therefore lower semicontinuous under this strong density and weak
gradient convergence. This proves the score inequality for polyhedral
tests. Subtract a supporting affine function and approximate every
finite convex test increasingly by polyhedral ones; monotone convergence
proves the full assertion. No moment beyond integrability of the score
is used.

The initial full shear has score $S_f(X)$ with the original score law.
Apply this argument successively to its existing coordinates. Folding
is another conditional averaging of the horizontal score, so the same
inequality survives folding.
\end{proof}

The integrable-gradient hypothesis permits scores with no $L^r$ moment
for any $r>1$. For example, normalize
$f(x)=(1-2x)/\log(1/x)$ on $(0,1/2)$ and extend it by zero.
Its weak derivative is integrable, while near zero its score is
asymptotic to $1/[x\log(1/x)]$. Hence
$\int f|S_f|^r=\infty$ for every $r>1$.

The score laws admit a reverse martingale coupling: consecutive scores
can be coupled with the later score equal to the conditional mean of
the earlier one. This follows from the convex-order characterization
\cite{Strassen}. Equation~\eqref{F-eq:score-defect} gives the corresponding
conditional covariance in the regular finite-second-moment case.
The coupling concerns score vectors; it does not specify a signing
algorithm.

\subsubsection{A barycentric section with a common score coupling}
Let $B\subset\R^m\times\R^k$ be bounded, open and convex, and suppose
$(y,s)\in B$ implies $(2p-y,s)\in B$. Let $F$ have mean $(p,\bar s)$
and an integrable horizontal weak gradient, with horizontal score
$Z_F=\nabla_yF/F$.

\begin{theorem}[First-moment section inequality]\label{S-thm:section}
If $\law(Z_F)\cx\mu$ for a centered integrable law $\mu$, then every
supporting inequality
$\langle a,y-p\rangle+\langle b,s\rangle<h$, $a\ne0$, satisfies
\begin{equation}\label{S-eq:r4-margin}
 h-\langle b,\bar s\rangle\ge\mathfrak b_\mu(a).
\end{equation}
Consequently
\begin{equation}\label{S-eq:r4-section}
 p+\operatorname{int}\mathcal K_\mu
 \subset\{y:(y,\bar s)\in B\}.
\end{equation}
\end{theorem}
\begin{proof}
Condition on $s$ and project onto $z=\langle a,y-p\rangle$.
Let $g_s$ denote the normalized marginal density and let $\nu_s$ be
its scalar score law. Slicing and marginalization preserve an integrable
weak derivative and zero extension. Distributional integration by parts
against a function of $(z,s)$ gives
\[
 \E[Z_F\mid z,s]=a\,g_s'(z)/g_s(z)=aQ.
\]
The random variable $Q$ has law $\int\nu_s\,d\law(s)$ and is integrable
by conditional Jensen. Thus $\law(aQ)\cx\mu$.
Horizontal symmetry bounds the support of $g_s$ by
$(-d(s),d(s))$, where $d(s)=h-\langle b,s\rangle$.
The interval theorem and convexity of $r$ under mixing give
\[
 h-\langle b,\bar s\rangle
   =\E d(s)\ge\E r(\nu_s)
   \ge r(\law Q)\ge\mathfrak b_\mu(a).
\]
This proves \eqref{S-eq:r4-margin}. An interior point of
$\mathcal K_\mu$ satisfies the corresponding inequality strictly.
A supporting inequality with $a=0$ is strict at the auxiliary mean
because the supporting set is open and convex. Separation proves the
section inclusion.
\end{proof}

\begin{proof}[Proof of Theorem~\ref{S-thm:main}]
Use the full shear from \eqref{R6-eq:full-shear}. Its horizontal score
has exactly law $\mu_f$. Rearrange all auxiliary coordinates and fold
them by absolute value. The resulting horizontal score is below
$\mu_f$ in convex order by Theorem~\ref{F-thm:score-cx}, and its
auxiliary mean $M$ satisfies $M_j\ge1$.
The original supporting set is invariant under simultaneous reflection
$(y,t)\mapsto(2p-y,-t)$. Coordinate symmetrization preserves this
invariance and adds separate reflections in every auxiliary coordinate.
Their composition gives horizontal symmetry about $p$ in the final set.
Theorem~\ref{S-thm:section} puts
$p+\operatorname{int}\mathcal K_{\mu_f}$ in its section at $M$.
Auxiliary reflection symmetry and convexity move this section to
$(1,\ldots,1)$ without changing the horizontal points.
The section identity \eqref{R6-eq:full-shear-section} and backward
recovery prove \eqref{S-eq:terminal}. Uniqueness and convergence of the
scalar optimizers are proved in
Theorems~\ref{S-thm:r5-stability} and \ref{S-thm:r6-dyadic}.
\end{proof}

An invertible linear map $L$ transports feasible couplings in both
directions, so
\[
 \mathfrak b_{L_\#\mu}(a)=\mathfrak b_\mu(L^{-1}a),\qquad
 \mathcal K_{L_\#\mu}=L^{-\mathsf T}\mathcal K_\mu.
\]
If $\mu\cx\eta$, the enlarged class of admissible scalar laws gives
$\mathcal K_\eta\subset\mathcal K_\mu$.
For later use, the scalar interval theorem and the Dirichlet inequality
also imply
\begin{equation}\label{S-eq:variance-radius}
 r(\nu)\ge\frac\pi{\sqrt{\E S^2}}
 \quad\text{for centered $\nu$ with finite second moment.}
\end{equation}
When $r(\nu)<\infty$, apply the interval Dirichlet inequality to the
square root of the reconstructed density on $(-r(\nu),r(\nu))$:
its derivative energy is $\E S^2/4$. The infinite case is immediate.

\subsubsection{Finite optimization of the common-coupling radius}
\label{S-sec:r6-dyadic}

The common scalar optimizer can also be approximated without guessing
its conditional hyperplanes.  Assume throughout this section that $\mu$ is centered, symmetric and
integrable, and $a\ne0$.  For $N=2^k$, set
\begin{equation}\label{S-eq:r6-dyadic-J}
 J_N(\nu)=\frac1N\sum_{i=1}^{N/2}\frac1{L_\nu(i/N)}.
\end{equation}
Let $B_N$ be its minimum over symmetric laws $\nu$ with
$\law(aS)\cx\mu$.  At each level it is enough to use $N$ equally
weighted atoms, with repetitions allowed.

\begin{theorem}[Dyadic reconstruction]\label{S-thm:r6-dyadic}
The minima $B_N$ are attained and satisfy
\begin{equation}\label{S-eq:r6-dyadic-limit}
 B_2\le B_4\le\cdots\uparrow\mathfrak b_\mu(a).
\end{equation}
When $\mathfrak b_\mu(a)<\infty$, choose any minimizing equal-mass law
at each level. These laws converge
in every $W_p$ for which the source has a finite $p$th moment to the
unique scalar extremizer $\nu_a$.
\end{theorem}
\begin{proof}
For a symmetric law $L$ is increasing on $(0,1/2)$, so the right
endpoint sums satisfy $J_N\le J_{2N}\le r$ and increase pointwise to
$r$.  Replacing $S$ by its conditional means on $N$ equal upper-quantile
bins preserves every $L(i/N)$, preserves symmetry, and decreases the
law in convex order.  Thus restricting to equal-mass atoms leaves the
minimum unchanged.  At each fixed $N$, their uniformly bounded first
moments bound all their atom locations, so compactness proves
attainment.

Let $\nu_N$ be minimizers. Lemma~\ref{S-lem:r7-compact} gives a
subsequential limit $\nu$ in every available $W_p$.  It remains
admissible.  For every fixed dyadic $D$, uniform convergence of
integrated quantiles and $J_N(\nu_N)\ge J_D(\nu_N)$ for $N\ge D$ give
$\lim_N B_N\ge J_D(\nu)$.  Taking $D\to\infty$ gives
$\lim_N B_N\ge r(\nu)\ge\mathfrak b_\mu(a)$.  The reverse inequality follows from $J_N\le r$ by taking infima.
When the limiting value is finite, uniqueness identifies every
subsequential limit and proves convergence in each available $W_p$. The same
compactness argument also gives $B_N\to+\infty$ when
$\mathfrak b_\mu(a)=+\infty$, since a finite limiting value would
produce an admissible law of finite radius.
\end{proof}

The finite values admit explicit certificates involving one convex
function of the complete vector.  Set $M=N/2$.  Choose numbers
$s_i\ge0$, $u_i\in\R$, and vectors $v_i\in\R^m$ satisfying
\begin{equation}\label{S-eq:r6-dual-constraint}
 2\langle a,v_i\rangle=\frac1N\sum_{j=i}^M s_j^2,
 \qquad 1\le i\le M.
\end{equation}
Then
\begin{equation}\label{S-eq:r6-dual}
 \boxed{\quad
 \mathfrak b_\mu(a)\ge B_N\ge
 \frac2N\sum_{i=1}^M(s_i+u_i)
 -\E_\mu\max_{1\le i\le M}
             \{u_i+|\langle v_i,Z\rangle|\}.
 \quad}
\end{equation}
All earlier quantiles are represented by the same maximum of affine
functions.  A certificate consists of these finitely many parameters
and an upper bound for this one expectation.

\begin{proposition}[Exact finite dual]\label{S-prop:r6-dual}
If the support of $\mu$ spans $\R^m$, the supremum of the right side of
\eqref{S-eq:r6-dual} over \eqref{S-eq:r6-dual-constraint} equals $B_N$.
Thus these finite max-affine certificates are complete as $N\to\infty$.
\end{proposition}
\begin{proof}
For an equal-mass target, write $m_j$ for the positive-bin scalar first
moments and $L_i=\sum_{j\le i}m_j$.  The reflected bins have moments
$-m_j$.  A common kernel into these $N$ labels gives, by
\eqref{S-eq:r6-dual-constraint},
\[
 \E\max_i\{u_i+|\langle v_i,Z\rangle|\}
 \ge\frac2N\sum_i u_i+\frac1N\sum_i s_i^2L_i.
\]
Combine this with $1/L_i\ge2s_i-s_i^2L_i$ to prove the inequality.

For equality, view the possible bin masses and vector first moments
as a compact convex set.  Its support function at affine labels is
exactly their expected maximum.  Symmetry pairs the labels into the
absolute-value form above.  Finite-dimensional separation therefore
gives the linear dual for maximizing $\sum_i\lambda_iL_i$.
There is no boundary obstruction to this separation: choose bounded
odd functions $h_1,\ldots,h_m$ so that
$\E[Zh^{\mathsf T}]$ is invertible, which is possible by truncating $Z$ and using that its support spans
$\R^m$.  Kernels of the form
$1/N\pm b_j^{\mathsf T}h(Z)$ then realize an open neighborhood of zero
in the allowable centered moment coordinates.  They retain all the
prescribed bin masses.  The affine constraint set therefore has the required relative interior
for the linear dual.

The bounded odd functions can be taken as
$h(Z)=Z\ind_{\{\|Z\|\le R\}}$ for sufficiently large $R$.
Their first-moment matrix is positive definite. Small coefficients in
the displayed kernels therefore realize arbitrary small bin moments,
including moments on the prescribed line. Small constant perturbations
also vary the paired masses around $1/N$. This proves the stated
relative-interior assertion without a regularity hypothesis on $\mu$.

Sorting the absolute scalar bin means into positive decreasing order
preserves feasibility and increases every $L_i$.  Therefore the same
moment set may be used without imposing an order on its labels.
At a minimizing positive vector $L^*$, the gradient of
$\sum_i1/(NL_i)$ gives a supporting linear objective
$\lambda_i=1/(N(L_i^*)^2)$.  Choose $s_i=1/L_i^*$ and apply the exact
linear dual just obtained.  The reciprocal inequalities are equalities
at $L^*$, proving that the supremum is $B_N$.
\end{proof}

\paragraph{Polyhedral extremizers and an exact gap identity}
\label{S-sec:r7-dual-cells}

The finite dual determines the optimizing kernel when the source is
absolutely continuous. It also separates the discrepancy of the
integrated quantiles from the discrepancy of the source allocation.

For $N=2^k$ and $M=N/2$, write the positive scalar bin means as
$q_1\ge\cdots\ge q_M\ge0$. A feasible kernel has labels $(i,\sigma)$,
$1\le i\le M$, $\sigma\in\{-1,1\}$, each of mass $1/N$, and
\[
 \E[Z\mid i,\sigma]=\sigma a q_i,
 \qquad L_i=\frac1N\sum_{j\le i}q_j.
\]
For dual parameters satisfying \eqref{S-eq:r6-dual-constraint}, set
\[
 \phi_{i,\sigma}(z)=u_i+\sigma\langle v_i,z\rangle,
 \qquad \Phi(z)=\max_{i,\sigma}\phi_{i,\sigma}(z),
\]
and write $D_N$ for the right side of \eqref{S-eq:r6-dual}.

\begin{theorem}[Unique finite partition and quantitative complementarity]
\label{S-thm:r7-cells}
Let $\mu$ be centered, symmetric and absolutely continuous, with finite
first moment. For each finite $N$, the minimum $B_N$ has a unique
ordered scalar target law and a unique labelled optimal coupling.
The dual supremum is attained. For any maximizing dual, the coupling
is the deterministic partition
\begin{equation}\label{S-eq:r7-cells}
 P_{i,\sigma}
 =\{z:\phi_{i,\sigma}(z)>\phi_{j,\tau}(z)
          \text{ for all }(j,\tau)\ne(i,\sigma)\},
\end{equation}
up to $\mu$-null sets. Each cell is the intersection of at most $N-1$
open halfspaces, has mass $1/N$, and has the prescribed vector
barycenter. The scalar means satisfy
$q_1>\cdots>q_M>0$.

For any feasible labelled kernel and any feasible dual, with its
random label denoted by $I$, the exact identity is
\begin{equation}\label{S-eq:r7-gap}
 J_N-D_N
 =\frac1N\sum_{i=1}^M\frac{(1-s_iL_i)^2}{L_i}
  +\E\bigl[\Phi(Z)-\phi_I(Z)\bigr].
\end{equation}
In particular, the probability assigned to labels whose affine score
is at least $\delta>0$ below the maximum is at most
$(J_N-D_N)/\delta$.
\end{theorem}
\begin{proof}
The possible labelled mass and moment vectors form a compact convex
set. Since the source support spans the ambient space, there are bounded odd functions
$h$ for which $\E[Zh^{\mathsf T}]$ is invertible. The kernels
$1/N\pm b_i^{\mathsf T}h(Z)$, together with small constant mass
perturbations, realize a relative neighborhood of zero in the
prescribed centered moment coordinates. Small strictly positive
ordered scalar means can be chosen in this neighborhood, with every
$L_i>0$. Finite-dimensional convex separation therefore gives an attained linear dual for every
supporting linear objective in these moments. This strengthens the
supremum statement in Proposition~\ref{S-prop:r6-dual} to attainment.

The vector $(L_i)$ ranges over a convex set when its labels are left
unordered. Sorting absolute scalar means preserves feasibility and
increases every partial sum. Thus a minimizer of the strictly convex
function $\sum_i1/(NL_i)$ has a unique vector $L^*$ and may be put in
positive decreasing order. The supporting objective at $L^*$ has
weights $1/(N(L_i^*)^2)$. Applying the attained linear dual with
$s_i=1/L_i^*$ gives an optimal affine maximum.

The projections of its slopes on $a$ are all distinct:
\[
 2\langle a,v_i-v_{i+1}\rangle=\frac{s_i^2}{N}>0,
 \qquad \langle a,v_M\rangle>0.
\]
Thus every equality between two different affine labels lies in a
proper hyperplane. Absolute continuity makes all ties null. Equality
in the dual forces every optimal kernel to use the unique maximizing
label at almost every source point. This proves uniqueness and
\eqref{S-eq:r7-cells}, including the mass and barycenter constraints.
If two scalar means coincided, exchanging their positive labels would
produce another optimal coupling. If $q_M=0$, exchanging its positive
and negative labels would do the same. The cells have positive mass,
so both possibilities contradict uniqueness.

Finally, the mass and moment constraints give
\[
 \E\phi_I(Z)=\frac2N\sum_i u_i
                +\frac1N\sum_i s_i^2L_i.
\]
Subtract this identity from the expected maximum, and use
$1/L-2s+s^2L=(1-sL)^2/L$. This proves \eqref{S-eq:r7-gap}.
Both terms are nonnegative, so Markov's inequality gives the stated
allocation bound.
\end{proof}

For an optimal dual, the first term in \eqref{S-eq:r7-gap} is
\[
 \frac1N\sum_i
 \frac{(L_i-L_i^*)^2}{(L_i^*)^2L_i}.
\]
The finite certificate therefore controls both the integrated
quantiles and the coupling's assignments. When $\mathfrak b_\mu(a)<\infty$, the optimizing scalar laws converge
as $N\to\infty$ to the unique common-coupling law
by Theorem~\ref{S-thm:r6-dyadic}. The finite partition theorem does not
require a choice of a limiting map.

\begin{example}[An exact half-square partition]\label{S-ex:r7-square}
For $Z$ uniform on $[-1,1]^2$, $a=(1,4/11)$ and $N=2$, the positive
cell is $x+y/2>0$. Its conditional mean is $(11/24,1/6)$, so the
scalar atoms are $\pm11/24$, $L_1=11/48$, and $B_2=24/11$.
A dual witness is
\[
 s_1=48/11,\qquad u_1=0,\qquad
 v_1=(576/143,288/143).
\]
Its expected maximum is $24/11$, which makes the dual value $24/11$.
Its positive cell has vertices $(-1/2,1),(1,1),(1,-1),(1/2,-1)$;
polygon integration gives the stated conditional mean. This is a finite-level
example; a bounded source has infinite full interval radius.
\end{example}

The scalar support theorem has now been transported through the full
shear. Its remaining input is the initial score distribution. The next
section constructs such distributions and compares the geometry they
retain; Section~\ref{INITOPT-sec} then optimizes the initial density
for the specified columns.


\subsection{Realized score laws and explicit retained regions}
\label{R6-sec:score-realization}
In every even dimension $m\ge4$, the following data can agree for two
even compactly supported log-concave densities:
\[
 J(f_m^{\rm B})=J(f_m^{\rm E})=I_m,\qquad
 (\mu_{f_m^{\rm B}})_J=(\mu_{f_m^{\rm E}})_J\quad(|J|\le3).
\]
Their retained bodies nevertheless satisfy
\[
 \frac{h_{\mathcal K_{\mu_{f_m^{\rm E}}}}(u_m)}
      {h_{\mathcal K_{\mu_{f_m^{\rm B}}}}(u_m)}
   \asymp\sqrt m,\qquad u_m=m^{-1/2}\mathbf1.
\]
Here $h_K(a)=\sup_{x\in K}\langle a,x\rangle$ is the support
function. Theorem~\ref{S-thm:r7-compact-pair} constructs these
densities with the same radial score profile. The difference lies in
the dependence of the score directions. Realizing the prescribed
joint score laws makes this separation explicit and also permits
computation of the retained body for the cosine reference.

\subsubsection{Compact realization of score laws}
\label{S-sec:r7-realization}

The common-coupling body is defined for a vector law. Its use as a score
invariant raises a separate realization question. The general existence
statement is already contained in the moment-measure theorem of
Cordero-Erausquin and Klartag \cite{CEK}: a centered probability law with
finite first moment whose support spans the space is the distribution
of a convex gradient under its own log-concave density, uniquely up to
translation. Applying that theorem to the reflected law gives the
following Sobolev consequence.

\begin{proposition}[Sobolev score realization]\label{S-prop:r7-realization}
Every centered law $\mu$ on $\R^m$ with finite second moment and
nonsingular covariance is the score law of a log-concave probability
density $f$ whose zero-extended amplitude belongs to $H^1(\R^m)$.
The essentially-continuous log-concave realization is unique up to
translation. If $\mu$ is centrally symmetric, $f$ can be centered to be
even. Its Fisher matrix is $\E_\mu ZZ^{\mathsf T}$.
\end{proposition}
\begin{proof}
Apply the moment-measure theorem to $(-\mathrm{id})_\#\mu$ and write
$f=e^{-V}$. Essential continuity gives zero trace on almost every
boundary slice and local absolute continuity along coordinate lines.
On the interior,
$\partial_i\sqrt f=-\tfrac12\sqrt f\,\partial_iV$.
The right side has squared integral
$\tfrac14\int z_i^2\,d\mu(z)<\infty$. Slicing and extension by zero
therefore give the asserted weak derivatives on the whole space.
Uniqueness and symmetry follow from the moment-measure theorem, after
fixing the barycenter of $f$. The gradient formula gives the Fisher
identity.
\end{proof}

For the explicit examples below, let $W$ denote the symmetric random
variable with density $2/[\pi(1+w^2)^2]$. It has variance one and
$\E|W|=2/\pi$; its absolute value has density
$4/[\pi(1+s^2)^2]$ on $(0,\infty)$.

Compactness can be characterized explicitly when the radial score
magnitude is prescribed. This also permits the same radial density
profile to be used on many different convex bodies.

Let $S>0$ have a positive continuous density $p$ on $(0,\infty)$ and
finite second moment. Write
\begin{equation}\label{S-eq:r7-tail}
 M(s)=\int_s^\infty t p(t)\,dt,
 \qquad \mathfrak r_S=\int_0^\infty\frac{p(s)}{M(s)}\,ds
                    =r(\law(\varepsilon S)),
\end{equation}
where $\varepsilon$ is a fair independent sign. Let $K$ be a centrally
symmetric convex body with interior containing zero, and let $g_K$ be
its Minkowski functional. For $X$ uniform in $K$, the direction
$X/g_K(X)$ has cone measure on $\partial K$; denote a random point with
this distribution by $\Theta_K$.

\begin{theorem}[A common radial profile]\label{S-thm:r7-profile}
For each dimension $m$ there is a convex increasing profile $V_m$,
unique up to an additive constant, such that
\begin{equation}\label{S-eq:r7-profile-density}
 f_K(x)=\frac{e^{-V_m(g_K(x))}}
 {m|K|\int_0^{R_m}r^{m-1}e^{-V_m(r)}\,dr},
 \qquad \{V_m<\infty\}=[0,R_m),
\end{equation}
has score law
\begin{equation}\label{S-eq:r7-score-polar}
 \law(Z_{f_K})=\law(-S\nabla g_K(\Theta_K)),
\end{equation}
where $S$ and $\Theta_K$ are independent. The radius $R_m$ is finite
if and only if $\mathfrak r_S<\infty$. When it is finite, $f_K$ is
positive on $R_m\operatorname{int}K$, has zero trace, and
$\sqrt{f_K}\in H^1(\R^m)$.

Put $x(s)=(V_m')^{-1}(s)$. There is a positive function $y$ satisfying
\begin{equation}\label{S-eq:r7-profile-ode}
 x'(s)=\frac{p(s)}{y(s)},\qquad
 y'(s)=p(s)\left(\frac{m-1}{x(s)}-s\right),
\end{equation}
and the exact identity
\begin{equation}\label{S-eq:r7-y-tail}
 y(s)=M(s)-(m-1)\int_s^\infty\frac{p(t)}{x(t)}\,dt.
\end{equation}
For $m=1$ the second term is zero. If
$p(s)\sim c s^{-q-1}$ for some $q>2$, then
\begin{equation}\label{S-eq:r7-boundary}
 R_m-x(s)\sim\frac{q-1}{s},
 \qquad e^{-V_m(r)}\sim c_m(R_m-r)^{q-1}
 \quad(r\uparrow R_m)
\end{equation}
for a constant $c_m>0$.
\end{theorem}
\begin{proof}
Apply the moment-measure theorem to the rotationally invariant law
$S U$, with $U$ uniform on the Euclidean unit sphere. Uniqueness makes
its centered potential radial, so it has the form $V_m(|x|)$.
Radial integration shows that its radius $R$ has density
$C r^{m-1}e^{-V_m(r)}$ and $V_m'(R)$ has density $p$.
Since $p$ is positive and continuous, the increasing derivative $V_m'$
has no jump and no interval of constancy on $(0,R_m)$. The distribution-function identity
$F_R(x(s))=F_S(s)$, whose two densities are positive and continuous,
makes the inverse $x$ continuously differentiable away from zero.
Writing
$y(s)=C x(s)^{m-1}e^{-V_m(x(s))}$ gives
\eqref{S-eq:r7-profile-ode}. The value at the upper endpoint is zero:
this follows from essential continuity when the endpoint is finite,
and from the exponential upper tail of an integrable log-concave
density when it is infinite. Integrating $y'$ gives
\eqref{S-eq:r7-y-tail}.

For all sufficiently large $s$, monotonicity of $x$ gives
$sx(s)\ge2(m-1)$. Since
$\int_s^\infty p(t)/x(t)\,dt\le M(s)/(sx(s))$, we obtain
\begin{equation}\label{S-eq:r7-comparison}
 \frac12M(s)\le y(s)\le M(s),
 \qquad \frac{p(s)}{M(s)}\le x'(s)\le\frac{2p(s)}{M(s)}.
\end{equation}
The integral of $p/M$ on a bounded initial interval is finite because
$M(0)=\E S>0$. Thus \eqref{S-eq:r7-comparison} proves the equivalence
between finite $R_m$ and finite $\mathfrak r_S$.

For a general $K$, integration on its dilates factors the probability
in \eqref{S-eq:r7-profile-density} into exactly the same radial law and
an independent cone direction. Homogeneity gives
$\nabla g_K(r\theta)=\nabla g_K(\theta)$ almost everywhere. This proves
\eqref{S-eq:r7-score-polar}. The gauge is Lipschitz and its gradient is
bounded, so the finite second moment of $S$ and zero trace give the
whole-space $H^1$ assertion.

Under the tail assumption, $M(s)\sim c s^{1-q}/(q-1)$, so
$\int^\infty p/M<\infty$. Now $x(s)\uparrow R_m>0$, and
\eqref{S-eq:r7-y-tail} gives $y(s)/M(s)\to1$. Therefore
$x'(s)\sim(q-1)s^{-2}$, proving the first part of
\eqref{S-eq:r7-boundary}. Substitution in
$y=Cx^{m-1}e^{-V_m(x)}$ proves the second part, including the finite
positive constant.
\end{proof}

For a polytope
\begin{equation}\label{S-eq:r7-polytope}
 K=\{x:\langle z_j,x\rangle\le1,\ \text{all }j\},
\end{equation}
let $F_j$ be its facets, $h_j=\|z_j\|^{-1}$ their distances from zero,
and
\begin{equation}\label{S-eq:r7-cone-weights}
 \omega_j=\frac{h_j|F_j|}{m|K|}.
\end{equation}
Then \eqref{S-eq:r7-score-polar} becomes
\begin{equation}\label{S-eq:r7-polytope-score}
 \mu_{f_K}=\sum_j\omega_j\law(-S z_j),
 \qquad J(f_K)=\E S^2\sum_j\omega_jz_jz_j^{\mathsf T}.
\end{equation}
Indeed, the cone from zero to $F_j$ has volume $h_j|F_j|/m$,
and the gauge gradient there is $z_j$. Thus the normalized cone-volume
weights prescribe the angular score law, while the same profile
prescribes its magnitude.

For the cosine magnitude $S=|W|$,
$p(s)=4/[\pi(1+s^2)^2]$ and $M(s)=2/[\pi(1+s^2)]$, so
$\mathfrak r_S=\pi$. Every dimension therefore has a compact profile
with quadratic vanishing at its boundary.

\subsubsection{Equal three-coordinate score marginals and different retained bodies}
\label{S-sec:r7-dependence}

The following construction realizes the dependence separation by
compact log-concave densities. Their scores have finite second moment;
the comparison concerns the joint score marginals, rather than the
spatial marginals of the densities.

Let $m\ge4$ be even, put $u_m=m^{-1/2}(1,\ldots,1)$, and let
$H=u_m^\perp$. Write
$B_1^m=\{x:\sum_i|x_i|\le1\}$, and let $P_H$ denote orthogonal
projection onto $H$. Consider
\begin{equation}\label{S-eq:r7-two-bodies}
 K_{\rm B}=B_1^m,
 \qquad
 K_{\rm E}=P_H B_1^m+[-m^{-1/2},m^{-1/2}]u_m.
\end{equation}
The second body is a cylinder with base $P_HB_1^m$. It has the
halfspace description
\begin{equation}\label{S-eq:r7-cylinder}
 K_{\rm E}=\left\{x:
 \left|\sum_i x_i\right|\le1,
 \ \langle\varepsilon,x\rangle\le1
 \text{ for every }\varepsilon\in\{-1,1\}^m
 \text{ with }\sum_i\varepsilon_i=0\right\}.
\end{equation}
To see this, the polar of the projected cross-polytope in $H$ is
$[-1,1]^m\cap H$, whose vertices for even $m$ are precisely the
balanced sign vectors. The two remaining inequalities give the
cylinder's end facets.

\begin{theorem}[Compact densities with matching score marginals]
\label{S-thm:r7-compact-pair}
Use $S=|W|$ in Theorem~\ref{S-thm:r7-profile}, and use its same profile
$V_m$ for $K_{\rm B}$ and $K_{\rm E}$. The resulting even log-concave
densities $f_m^{\rm B}$ and $f_m^{\rm E}$ have compact supports
$R_mK_{\rm B}$ and $R_mK_{\rm E}$, zero-extended amplitudes in $H^1$,
and Fisher matrix $I_m$. Their score laws satisfy
\begin{equation}\label{S-eq:r7-three-marginals}
 (\mu_{f_m^{\rm B}})_J=(\mu_{f_m^{\rm E}})_J
 \quad\text{for every }J\subseteq\{1,\ldots,m\},\quad |J|\le3.
\end{equation}
Every coordinate score has density $2/[\pi(1+t^2)^2]$, and
$h_{\mathcal K_{\mu_{f_m^{\rm B}}}}(e_i)
 =h_{\mathcal K_{\mu_{f_m^{\rm E}}}}(e_i)=\pi$.
In the unit diagonal direction,
\begin{align}
 h_{\mathcal K_{\mu_{f_m^{\rm E}}}}(u_m)
 &=\frac{\pi(m+1)}{2\sqrt m},\label{S-eq:r7-E-width}\\
 h_{\mathcal K_{\mu_{f_m^{\rm B}}}}(u_m)
 &\longrightarrow r_*:=r(\law(WG)),\label{S-eq:r7-B-width}
\end{align}
where $G$ is an independent standard Gaussian. Thus the support ratio
is $(\pi/(2r_*)+o(1))\sqrt m$.

The geometric supports satisfy $K_{\rm B}\subset K_{\rm E}$ and
\begin{equation}\label{S-eq:r7-volume-ratio}
 \frac{|K_{\rm E}|}{|K_{\rm B}|}
 =\frac{m}{2^m}\binom{m}{m/2}
 \sim\sqrt{\frac{2m}{\pi}}.
\end{equation}
\end{theorem}
\begin{proof}
Every facet normal of $K_{\rm B}$ is a sign vector, and all normalized
cone-volume weights equal $2^{-m}$. The two end facets of the cylinder
$K_{\rm E}$ have normals $\pm(1,\ldots,1)$ and combined cone weight
$1/m$. Coordinate permutations act transitively on its remaining
facets; their normals are the balanced sign vectors. Their total
weight is $(m-1)/m$. Equation~\eqref{S-eq:r7-polytope-score} therefore
identifies the two score laws as $S\varepsilon^{\rm B}$ and
$S\varepsilon^{\rm E}$, where $S$ is independent of the signs,
$\varepsilon^{\rm B}$ is uniform on all signs, and
$\varepsilon^{\rm E}$ is a constant sign vector with probability $1/m$
and is otherwise uniformly balanced.

Both sign laws are centrally symmetric. The balanced law has
$\E\varepsilon_i\varepsilon_j=-1/(m-1)$ for $i\ne j$, so the mixture
has zero pair correlations. Every product of an odd number of distinct
signs has zero expectation. The Fourier expansion on the Boolean cube
then proves that every at-most-three-coordinate sign marginal is
uniform, which proves \eqref{S-eq:r7-three-marginals}. Since $\E S^2=1$,
both Fisher matrices are $I_m$.

For a unit vector $a$ with
$\E[Z\mid\langle a,Z\rangle]=a\langle a,Z\rangle$, one has
\begin{equation}\label{S-eq:r7-regression-support}
 h_{\mathcal K_\mu}(a)=\mathfrak b_\mu(a)
   =r(\law(\langle a,Z\rangle)).
\end{equation}
Indeed the projection is admissible. For every other normal $b$, projection
of an admissible martingale onto $a$ and monotonicity of $r$ give
$\mathfrak b_\mu(b)\ge|\langle a,b\rangle|r(\law(\langle a,Z\rangle))$.
Thus $a\,r(\law(\langle a,Z\rangle))\in\mathcal K_\mu$, proving equality.

Pairwise independence of the signs gives the regression for each
coordinate, and permutation invariance gives it for $u_m$. The
E-projection is zero with probability $1-1/m$ and has law $\sqrt mW$
with probability $1/m$. For any symmetric nondegenerate $X$, $0<p\le1$
and $c>0$, direct integration of its quantile plateau gives
\begin{equation}\label{S-eq:r7-atom-radius}
 r((1-p)\delta_0+p\law(cX))
 =\frac{r(\law X)}c+\frac{1-p}{pc\,\E|X|}.
\end{equation}
Using $r(\law W)=\pi$ and $\E|W|=2/\pi$ proves
\eqref{S-eq:r7-E-width}. The B-projection has law $A_mW$, where
$A_m=|\sum_i\varepsilon_i^{\rm B}|/\sqrt m$. The limiting radius is
justified below.

The inclusion of the two supports follows by projecting a point of
$B_1^m$ onto $H$ and using $|\sum_i x_i|\le1$. To compute the volume,
each facet of $B_1^m$ has $(m-1)$-area $\sqrt m/(m-1)!$ and unit normal
$\varepsilon/\sqrt m$. Summing the projected areas of the outward
facets gives
\[
 |P_HB_1^m|_{m-1}
 =\frac{2^{m-1}}{(m-1)!\sqrt m}\,
   \E\left|\sum_i\varepsilon_i^{\rm B}\right|.
\]
Multiplying by the cylinder height $2/\sqrt m$ and dividing by
$|B_1^m|=2^m/m!$ leaves the expected absolute sign sum. For even $m$,
this expectation is $m2^{-m}\binom m{m/2}$. Stirling's formula proves
the asymptotic.
\end{proof}

\paragraph{The limiting radius and the fourth sign statistic}
For an independent nonnegative $A$, the continuous part of $AW$ has,
for $x>0$,
\[
 p_A(x)=\frac2\pi\E\frac{A^3}{(A^2+x^2)^2},\qquad
 M_A(x)=\frac1\pi\E\frac{A^3}{A^2+x^2}.
\]
If $p_0=\Pp(A=0)$, its full radius is
\begin{equation}\label{S-eq:r7-product-radius}
 r(\law(AW))=\frac{p_0}{\E A\,\E|W|}
             +\int_0^\infty\frac{p_A(x)}{M_A(x)}\,dx.
\end{equation}
For $A=A_m$, the central limit theorem gives $A_m\Rightarrow|G|$ and
uniform integrability of $A_m$. The central atom tends to zero.
For fixed $x>0$, both expectations in the integrand converge, and the
ratio is at most $2/x^2$ for large $x$. On $(0,\epsilon)$ its integral
is at most
$\Pp(0<A_mW<\epsilon)/M_{A_m}(1)$ for $\epsilon\le1$.
The denominators converge to a positive limit, while $WG$ has no atom
at zero. Truncation at both ends and convergence on the intervening
compact interval prove \eqref{S-eq:r7-B-width}.

Writing $y=x^2/2$ and $E_1(y)=\int_y^\infty e^{-t}\,dt/t$ gives
\begin{equation}\label{S-eq:r7-rstar}
 r_* =\int_0^\infty
 \frac{(1+y)e^yE_1(y)-1}{1-ye^yE_1(y)}\,dx
 =3.2491341681\ldots.
\end{equation}
The integral defines the constant; numerical quadrature is used only
for the displayed decimal approximation.

The first differing bounded coordinate statistic occurs at order four:
for distinct $i,j,k,l$,
\begin{equation}\label{S-eq:r7-four-signs}
 \E\prod_{a\in\{i,j,k,l\}}\operatorname{sgn}Z_a^{\rm E}
 =\frac1{m-3},\qquad
 \E\prod_{a\in\{i,j,k,l\}}\operatorname{sgn}Z_a^{\rm B}=0.
\end{equation}
For the balanced sign law the fourth product is
$3/[(m-1)(m-3)]$, which gives the displayed mixture value. These are
bounded sign statistics. Raw fourth score moments are infinite.

\paragraph{Four-wise independence in the common-magnitude class}
The order of marginal agreement has a quantitative consequence within
the common-magnitude, exchangeable class.

\begin{proposition}[Four-wise independence bounds the diagonal radius]
\label{S-prop:r7-fourwise}
Let $\varepsilon\in\{-1,1\}^m$ have an exchangeable, centrally
symmetric, four-wise independent law. Let $S=|W|$ be independent and
put $\mu=\law(S\varepsilon)$. Then
\begin{equation}\label{S-eq:r7-fourwise-bound}
 \pi\le h_{\mathcal K_\mu}(u_m)
 \le\pi\sqrt{3-2/m}.
\end{equation}
\end{proposition}
\begin{proof}
Set $C=m^{-1/2}\sum_i\varepsilon_i$. Exchangeability gives the
regression in \eqref{S-eq:r7-regression-support}, so the support is
$r(\law(|C|W))$. Four-wise independence gives
$\E C^2=1$ and $\E C^4=3-2/m$. H\"older's inequality yields
$\E|C|\ge(\E C^2)^{3/2}/(\E C^4)^{1/2}$.
Conditional Jensen gives
$(\E|C|)W\cx |C|W$, so monotonicity and scaling of $r$ imply the
upper bound. The projection has variance one; the sharp one-dimensional
Dirichlet inequality gives $r\ge\pi$.
\end{proof}

Thus four-wise independence rules out an unbounded diagonal separation
within this specified common-magnitude class. This statement places no
restriction on general pairs of vector laws agreeing on four-coordinate
marginals.

\subsubsection{The explicit score body of the cosine product}\label{S2-sec}
The abstract retained body becomes highly non-Euclidean for the cosine
product.  Its score law can be written explicitly:
\begin{equation}\label{S2-eq:scorelaw}
 Z_{f_R}\stackrel d=-\frac\pi R(T_1,\ldots,T_m),
 \qquad
 p_T(t)=\frac2{\pi(1+t^2)^2},
 \quad \E T=0,
 \quad \E T^2=1,
\end{equation}
with independent $T_i$.  For $w\ne0$ put
$\nu_w=\law(\sum_iw_iT_i)$.

\begin{theorem}[Directional cosine radii]\label{S2-thm:radii}
For every nonzero $w$,
\begin{equation}\label{S2-eq:radii}
 \frac\pi{\|w\|_2}\le r(\nu_w)\le\frac\pi{\|w\|_\infty},
\end{equation}
and equality in the lower bound occurs exactly when $w$ has one nonzero
coordinate.  For $a,b>0$,
\begin{equation}\label{S2-eq:two}
 r(\law(aT_1+bT_2))
 =\frac\pi{a+b}
 \left(2-\sqrt{\frac{a^2-ab+b^2}{a^2+ab+b^2}}\right).
\end{equation}
In particular, for $u=(e_1+e_2)/\sqrt2$ the physical retained radius is
\begin{equation}\label{S2-eq:twoequal}
 \frac R\pi r(\nu_u)
 =R\left(\sqrt2-\frac1{\sqrt6}\right)>R.
\end{equation}
For the equal $k$-coordinate unit vector
$u_k=k^{-1/2}(1,\ldots,1)$,
\begin{equation}\label{S2-eq:diffuse}
 r(\nu_{u_k})^2
 =\log k+2\log\log k+4+\log(\pi/2)
 +O\!\left(\frac{\log\log k}{\log k}\right).
\end{equation}
\end{theorem}
\begin{proof}
The variance bound \eqref{S-eq:variance-radius} and
$\Var(\sum_iw_iT_i)=\|w\|_2^2$ give the lower bound.  Conditioning on
$T_i$ gives
$\law(w_iT_i)\cx\nu_w$, and monotonicity of $r$ under convex order gives
the upper bound.  The characteristic function is
\begin{equation}\label{S2-eq:cf}
 \widehat\nu_w(t)=e^{-\|w\|_1|t|}\prod_i(1+|w_i||t|).
\end{equation}
Equality in the variance lower bound forces a scaled copy of $T$; the
exponential decay in \eqref{S2-eq:cf} then forces
$\|w\|_1=\|w\|_2$, hence one nonzero coordinate.

For two coordinates put $L=a+b$, $U=a^2-ab+b^2$, and
$W=a^2+ab+b^2$.  Fourier inversion of
$e^{-L|t|}(1+a|t|)(1+b|t|)$ gives, for the density $p$ and its tail first
moment $M(x)=\int_x^\infty yp(y)dy$,
\[
 \pi p(x)=\frac{2L[Ux^2+L^2(L^2+ab)]}{(x^2+L^2)^3},\qquad
 \pi M(x)=\frac{L[Ux^2+L^2W]}{(x^2+L^2)^2}.
\]
Thus
\[
 \frac{p(x)}{M(x)}=
 \frac4{x^2+L^2}-\frac2{x^2+L^2W/U}.
\]
The one-dimensional radius formula from
Section~\ref{S-sec:score} integrates these two rational terms and gives
\eqref{S2-eq:two}.  The diffuse asymptotic follows from the logarithmic
zero-bias representation of the radius applied to
\eqref{S2-eq:cf}; the complete saddle-point estimate is recorded in
Appendix~\ref{S2-app:diffuse}.
\end{proof}

The first inequality in \eqref{S2-eq:radii} already has a strong
geometric consequence.  By \eqref{S2-eq:scorelaw}, linear transport of
the retained body and \eqref{S-eq:variance-radius} give
\begin{equation}\label{S2-eq:ball}
 R B_2^m\subset\mathcal K_{\law(Z_{f_R})}.
\end{equation}
The inclusion is strict away from coordinate directions, and it grows
on diffuse directions according to \eqref{S2-eq:diffuse}.

\begin{corollary}[Inhomogeneous vector balancing]\label{S2-cor:displaced}
Assume the height conditions for the cosine product at radius $R$.
Then for every $x\in\R^m$ with $\|x\|_2<R$ there are signs satisfying
\begin{equation}\label{S2-eq:displaced}
 \left\|x+\sum_j\sigma_ja_j\right\|_\infty<R.
\end{equation}
At the profile-sharpened radius $R=R_A$ from
Theorem~\ref{MAIN-signing}, this holds simultaneously for every
$x\in R_AB_2^m$.
\end{corollary}
\begin{proof}
Equation \eqref{S2-eq:ball} places $x$ in the interior of the retained
score body.  Apply Theorem~\ref{S-thm:main} and recover signs backward
from the terminal section.
\end{proof}

The body is substantially larger than this Euclidean core.  Its volume
radius obeys
\begin{equation}\label{S2-eq:volume}
 \operatorname{vrad}(\mathcal K_{\law(T_1,\ldots,T_m)})
 \ge(1-o(1))\sqrt{\log m},
\end{equation}
so the physical retained displacement region has volume at least
\begin{equation}\label{S2-eq:physical-volume}
 |B_2^m|\left((1-o(1))\frac{R\sqrt{\log m}}\pi\right)^m.
\end{equation}
Indeed a uniform spherical direction has
$\|w\|_\infty\to0$ and
$\sum_i|w_i|^3=(\E|G|^3+o_{\Prb}(1))/\sqrt m$; the uniform form of the
zero-bias expansion gives $r(\nu_w)/\sqrt{\log m}\to1$ in probability.
The radial star body generated by these radii lies inside its convex
hull $\mathcal C$, and taking $m$th roots of volume proves
\eqref{S2-eq:volume}.

\subsubsection{The retained horizontal zonoid and matrix-adapted displacement}
\label{ZON-sec}
The score body depends only on the initial density.  A different invariant
keeps the accumulated column geometry.  For an integrable symmetric random
vector $Y$, write
\begin{equation}\label{ZON-eq:def}
 h_{Z_1(Y)}(a)=\E|a^{\mathsf T}Y|.
\end{equation}
Equivalently, $Z_1(Y)$ is the compact convex set of vectors
$\E[Y\theta(Y)]$ with measurable $|\theta|\le1$.

\begin{theorem}[Horizontal zonoid retained by every stage]
\label{ZON-thm:retained}
Let $K=-K$, let $f$ and every auxiliary density $q_j$ be even, and
suppose $H_f(v_j;q_j)\ge1$.  If $X\sim f$ and the $T_j\sim q_j$ are
independent, set
\[
 Y=X-\sum_jT_jv_j.
\]
Then
\begin{equation}\label{ZON-eq:inclusion}
 Z_1(Y)\subset
 \calS_{v_n}\cdots\calS_{v_1}K.
\end{equation}
The same inclusion holds after every prefix and for every ordering of the
vectors.  Every point of the zonoid admits backward sign recovery into $K$.
\end{theorem}
\begin{proof}
The full shear in \eqref{R6-eq:full-shear} has horizontal marginal $Y$,
and every auxiliary-coordinate rearrangement preserves that marginal.
The sheared source obeys $F(-y,t)=F(y,-t)$.  Once every auxiliary
coordinate has been symmetrized and folded, the resulting density is even
in $y$ at each fixed folded auxiliary vector $S$.

It is enough to consider an odd measurable $\theta$ with
$|\theta|\le1$, since replacing any $\theta$ by its odd part does not
change $\E[Y\theta(Y)]$.  Conditional horizontal symmetry gives
$\E[\theta(Y)\mid S]=0$.  Weight the folded density by
$1+\theta(Y)$.  The new density is nonnegative, has total mass one,
retains the complete $S$ marginal, and has horizontal mean
$\E[Y\theta(Y)]$.  Its support lies in the same open convex lifted set.
The barycenter argument from Section~\ref{sec:coordinate}, followed
by coordinatewise lowering of the folded auxiliary mean to
$(1,\ldots,1)$, places this horizontal mean in the terminal section.
Taking all $\theta$ proves \eqref{ZON-eq:inclusion}.  Reversing the
Steiner transforms gives the signs.
\end{proof}

For the optimized auxiliary this invariant has an explicit ellipsoidal
core.  Write
\begin{equation}\label{ZON-eq:sigma}
 \beta_R=R^2\left(\frac13-\frac2{\pi^2}\right),\qquad
 w=\E T_1^2,
 \qquad
 \Sigma_R=\beta_RI_m+wAA^{\mathsf T}.
\end{equation}
Write the polynomial in \eqref{R-eq:htrial} as
$h(x)=\sum_{r=1}^5c_{2r}x^{2r}$.  It is strictly concave on $[-1,1]$:
omitting the negative $c_4$ and $c_8$ contributions gives
\[
 h''(x)\le2c_2+30c_6+90c_{10}<-0.5527.
\]
Hence $q_\dagger$ is log-concave.  The cosine product is log-concave as
well, so $Y_R=X_R-AT$ is symmetric and log-concave.

\begin{corollary}[Column-covariance displacement ellipsoid]
\label{ZON-cor:ellipsoid}
At the profile radius $R_A$ of Theorem~\ref{MAIN-signing},
\begin{equation}\label{ZON-eq:ellipsoid}
 \frac1{\sqrt2}\Sigma_{R_A}^{1/2}B_2^m
 \subset Z_1(X_{R_A}-AT)
 \subset\calS_{a_n}\cdots\calS_{a_1}(-R_A,R_A)^m.
\end{equation}
Consequently every
$x\in2^{-1/2}\Sigma_{R_A}^{1/2}B_2^m$ admits a signing with
\begin{equation}\label{ZON-eq:displaced}
 \|x+A\sigma\|_\infty<R_A.
\end{equation}
The terminal transform also contains the score body from
Theorem~\ref{S-thm:main}, hence contains the convex hull of the two
retained regions.
\end{corollary}
\begin{proof}
For a symmetric real log-concave random variable $Z$, its survival
function $S(t)=\Prb(|Z|>t)$ is log-concave on $[0,\infty)$.  Therefore
$S(s+t)\le S(s)S(t)$ and
\[
 \E Z^2
 =2\int_0^\infty\!\int_0^\infty S(s+t)\,ds\,dt
 \le2(\E|Z|)^2.
\]
Apply this to $Z=a^{\mathsf T}Y_{R_A}$.  Independence in
\eqref{ZON-eq:sigma} gives $\Cov(Y_{R_A})=\Sigma_{R_A}$, so
\[
 h_{Z_1(Y_{R_A})}(a)=\E|a^{\mathsf T}Y_{R_A}|
 \ge\sqrt{a^{\mathsf T}\Sigma_{R_A}a/2}.
\]
The support-function inequality is exactly the first inclusion in
\eqref{ZON-eq:ellipsoid}.  The second is
Theorem~\ref{ZON-thm:retained}, using the strict heights from
Theorem~\ref{MAIN-signing}.  Backward recovery gives
\eqref{ZON-eq:displaced}.  The final assertion uses convexity of the
terminal transform and Theorem~\ref{S-thm:main}.
\end{proof}

This ellipsoid and the score body encode different information.  The
score body is determined by the initial reference law and can grow in
diffuse directions even when its Fisher matrix is fixed.  The zonoid
instead sees the accumulated covariance $AA^{\mathsf T}$ of the actual
columns.  Both survive the same sequence of transforms.


\subsection{What the score moments fail to retain}
\label{R9-sec:score-moments}
The body in Theorem~\ref{S-thm:main} is determined by the full score
law. Even its complete sequence of polynomial moments can lose the
sharp support radius. This follows by applying the scalar
reconstruction to a classical moment-indeterminate family.

\begin{theorem}[Equal score moments and different sharp supports]
\label{R9-thm:score-moment-indeterminacy}
There is a family of even compactly supported log-concave densities
$f_t$ on the line, smooth on the interiors of their supports, whose
zero extensions belong to $W^{1,1}$, such that all their score laws
have exactly the same finite moments of every order and Fisher
information one. Their support half-lengths attain the sharp scalar
score radius and are not constant in $t$.
\end{theorem}
\begin{proof}
For $|t|\le1/2$, let $R_t>0$ have density
\[
 g_t(r)=\frac{e^{-(\log r)^2/2}}{r\sqrt{2\pi}}
                   [1+t\sin(\pi\log r)].
\]
If $Z$ is standard Gaussian, then for every integer $k\ge0$,
\[
 \E e^{kZ}\sin(\pi Z)
       =e^{(k^2-\pi^2)/2}\sin(\pi k)=0.
\]
Thus the distinct laws have identical moments. Let $\nu_t$ be the
symmetric sign symmetrization of $R_t$ and set
$\widetilde\nu_t=(\nu_t+N(0,1))/2$. The score density is now smooth
and positive on the whole line. Its common variance is
$\sigma^2=(e^2+1)/2$.

The symmetric lognormal law has finite radius. Writing $\phi$ and
$\overline\Phi$ for the Gaussian density and upper tail gives
\[
 r(\nu_0)=\int_{-\infty}^{\infty}
   \frac{\phi(z)}{e^{1/2}\overline\Phi(z-1)}\,dz<\infty.
\]
The integrand has Gaussian decay at $-\infty$ and order $ze^{-z}$
at $+\infty$. Each $\widetilde\nu_t$ contains at least one quarter
of $\nu_0$. Concavity of the upper integrated quantile under mixing,
and nonnegativity of that quantile for centered laws, give finite
radii uniformly in $t$.

For distinct centered laws $\nu,\eta$ with finite radii,
\[
 L_{(1-s)\nu+s\eta}\ge(1-s)L_\nu+sL_\eta
 \quad\Longrightarrow\quad
 r((1-s)\nu+s\eta)<(1-s)r(\nu)+sr(\eta),\quad0<s<1.
\]
Strictness follows from strict convexity of $x\mapsto1/x$ and the
fact that distinct laws have different integrated quantiles on an
interval. The affine family $\widetilde\nu_t$ therefore has
nonconstant radius. Apply the scalar reconstruction of
Theorem~\ref{S-thm:interval} to the variance-one laws of
$S/\sigma$, $S\sim\widetilde\nu_t$:
\[
 x(u)=-r+\int_0^uL(v)^{-1}\,dv,\qquad f(x(u))=L(u).
\]
Its decreasing score quantile makes $f$ log-concave. The positive
smooth score density gives interior smoothness. The finite first
score moment gives the zero-extension $W^{1,1}$ property, and the
second moment gives Fisher information one. The constructed support
has half-length exactly $r$.
\end{proof}

Lognormal Stieltjes classes are classical; see \cite{Kleiber}.
Here their realization as scores converts moment indeterminacy into
a difference in sharp compact support. The entropy inverse in
Theorem~\ref{SP-thm:entropy-inverse} later gives an operational
transform which does determine the full symmetric law.

\begin{proposition}[Exponential score moments on bounded supports]
\label{R9-prop:score-exponential}
Let $f$ be a density of bounded support in $\R^D$ whose zero extension
belongs to $W^{1,1}$, and let $S=\nabla\log f(X)$, $X\sim f$.
For every nonzero $a\in\R^D$ and every $t\ne0$,
\[
 \E e^{t\langle a,S\rangle}=\infty.
\]
All polynomial score moments can nevertheless be finite.
\end{proposition}
\begin{proof}
For a scalar centered score $T$, suppose
$M=\E e^{\lambda T}<\infty$ for some $\lambda>0$. Chernoff's
quantile bound gives
\[
 L_{\law(T)}(u)\le\frac{u}{\lambda}
                  [1+\log(M/u)],\qquad0<u<1.
\]
The integral of its reciprocal diverges at zero, contradicting the
finite support inequality in Theorem~\ref{S-thm:interval}.
Reflection handles $\lambda<0$. For a unit direction $a$, the score
of the marginal $\langle a,X\rangle$ is
$\E[\langle a,S\rangle\mid\langle a,X\rangle]$.
Marginalization preserves bounded support and the zero-extension
$W^{1,1}$ property. Conditional Jensen would pass a finite
exponential moment to this scalar score, giving the contradiction.
Rescaling handles general $a$.

For the last claim, the symmetric score density
$q(s)=\tfrac14 e^{-\sqrt{|s|}}$ has moments
$\E S^{2k}=(4k+1)!$ and finite radius
\[
 r(q)=\int_0^\infty\frac{z\,dz}{z^3+3z^2+6z+6}.
\]
Its scalar reconstruction gives a compactly supported log-concave
physical density with all polynomial score moments finite.
\end{proof}


\section{Optimizing the reference and the universal overlap threshold}
\label{R11-sec:optimization}
For fixed directions, overlap has a nonlocal Cheeger dual and its small-translation limit is directional Fisher energy. The cube criterion determines exactly when the cosine product is optimal; a mixed derivative improves it when the criterion fails.

The universal problem selects a reference before an unknown diffuse direction, allowing arbitrary dependence. Ground-state scores give the upper comparison and product references attain its limit. This reference-before-direction threshold has different quantifiers from the all-signing-law tradeoff of Section~\ref{R11-sec:joint}.

\subsection{Optimizing the initial law}\label{INITOPT-sec}
\subsubsection{The full nonlocal problem}
Let $K\subset\R^m$ be bounded and open, and let $q_0=1/6$ on $(-3,3)$.
For a probability density $f$ supported in $K$ and a direction $v$, set
\begin{equation}\label{INITOPT-eq:T}
 \mathcal T_v(f)=\int_0^6(6-t)
       \|f-f(\cdot+tv)\|_{L^1(\R^m)}\,dt.
\end{equation}
A direct change of variables in \eqref{F-eq:H-weighted} gives
\begin{equation}\label{INITOPT-eq:uniform-height}
 H_f(v;q_0)=\frac32-\frac1{24}\mathcal T_v(f).
\end{equation}
Thus $\mathcal T_v(f)\le12$ is exactly the height-one condition for the
uniform auxiliary.

For a measurable $E\subset K$ of positive volume define
\begin{equation}\label{INITOPT-eq:perimeter}
 P_v(E)=\int_0^6(6-t)|E\mathbin\triangle(E-tv)|\,dt.
\end{equation}
For directions $V=(v_1,\ldots,v_n)$ and
$\lambda\in\Delta_n=
\{\lambda\in[0,1]^n:\sum_j\lambda_j=1\}$, put
\begin{equation}\label{INITOPT-eq:cheeger}
 h_K(\lambda)=\inf_{E\subset K,\ |E|>0}
 \frac{\sum_j\lambda_jP_{v_j}(E)}{|E|}.
\end{equation}

\begin{theorem}[Exact nonlocal Cheeger dual]\label{INITOPT-thm:cheeger}
Let $\mathcal D(K)$ denote the probability densities supported in $K$.
Then
\begin{equation}\label{INITOPT-eq:cheeger-dual}
 \inf_{f\in\mathcal D(K)}\max_j\mathcal T_{v_j}(f)
 =\max_{\lambda\in\Delta_n}h_K(\lambda).
\end{equation}
The maximum is attained.  If $K=-K$ is convex and the common value is
strictly below $12$, there is a centered density with
$H_f(v_j;q_0)>1$ for every $j$; hence Theorem~\ref{thm:affine}
produces a signing whose sum lies in $K$.
\end{theorem}
\begin{proof}
Each $\mathcal T_v$ is convex in $f$ and satisfies
$|\mathcal T_v(f)-\mathcal T_v(g)|\le36\|f-g\|_1$.
Sion's minimax theorem~\cite{Sion}, with the simplex as its compact
variable, therefore gives
\[
 \inf_f\max_{\lambda\in\Delta_n}\sum_j\lambda_j\mathcal T_{v_j}(f)
 =\max_{\lambda\in\Delta_n}\inf_f
        \sum_j\lambda_j\mathcal T_{v_j}(f).
\]
For the superlevel sets $E_a=\{f>a\}$, layer cake and Tonelli give
simultaneously
\[
 \mathcal T_{v_j}(f)=\int_0^\infty P_{v_j}(E_a)\,da,
 \qquad 1=\int_0^\infty |E_a|\,da.
\]
Hence the inner infimum is at least $h_K(\lambda)$, while the normalized
indicator of a measurable set $E$ proves the reverse inequality.
Upper semicontinuity in $\lambda$ gives attainment.  If $K=-K$,
reflection preserves every $\mathcal T_{v_j}$ and convexity permits
central symmetrization of an almost minimizer.  Equations
\eqref{INITOPT-eq:uniform-height} and
Theorem~\ref{thm:affine} finish the proof.
\end{proof}

The formula remains exact under a prescribed barycenter.  If $p\in K$
and
\[
 \tau_K(p)=\inf_{\substack{f\in\mathcal D(K)\\\int xf(x)dx=p}}
             \max_j\mathcal T_{v_j}(f),
\]
then the same supporting-hyperplane argument and coarea formula give
\begin{equation}\label{INITOPT-eq:barycenter-dual}
 \tau_K(p)=\max_{\lambda\in\Delta_n,\ \ell\in\R^m}
 \inf_{E\subset K,\ |E|>0}
 \frac{\sum_j\lambda_jP_{v_j}(E)
 +\int_E\ell\cdot(x-p)\,dx}{|E|}.
\end{equation}
If $r=\dist(p,\partial K)$, a maximizing multiplier may be chosen with
$\|\ell\|_2\le36/r$.  The proof uses that $0\le\tau_K\le36$ and tests a
supporting affine function at $p\pm r'u$, then lets $r'\uparrow r$.

\subsubsection{Directional Fisher energy}
For $\psi\in H_0^1(K)$ and $v\in\R^m$, write
\[
 \mathcal E_v(\psi)=\int_K|\partial_v\psi|^2,
 \qquad
 \kappa_K(V)=\inf_{\substack{\psi\in H_0^1(K)\\\|\psi\|_2=1}}
          \max_j\mathcal E_{v_j}(\psi).
\]
For a positive semidefinite matrix $M$, define
\[
 \lambda_1(K;M)=\inf_{\substack{\psi\in H_0^1(K)\\\|\psi\|_2=1}}
 \int_K\nabla\psi^{\mathsf T}M\nabla\psi.
\]

\begin{theorem}[Spectral minimax]\label{INITOPT-thm:spectral}
Suppose $K$ is bounded, open, and convex, and the directions $v_j$ span
$\R^m$.  Then
\begin{equation}\label{INITOPT-eq:spectral}
 \kappa_K(V)=\max_{\lambda\in\Delta_n}
 \lambda_1\!\left(K;\sum_j\lambda_jv_jv_j^{\mathsf T}\right).
\end{equation}
Both extrema are attained.  There is an optimizing amplitude
$\psi\ge0$ for which $\psi^2$ is log-concave; if $K=-K$, it can be
chosen even.  For primal and dual optimizers,
\begin{equation}\label{INITOPT-eq:complementarity}
 \lambda_j>0\Longrightarrow
 \mathcal E_{v_j}(\psi)=\kappa_K(V).
\end{equation}
\end{theorem}
\begin{proof}
Use $f=\psi^2$ as the variable.  Pointwise Cauchy--Schwarz gives
\[
 \mathcal E_v(\sqrt{\theta f+(1-\theta)g})
 \le\theta\mathcal E_v(\sqrt f)+(1-\theta)\mathcal E_v(\sqrt g),
\]
so the directional energies are convex in the density.  They are lower
semicontinuous in $L^1$, and the simplex is compact, which proves the
minimax identity.  Since the directions span, a bound on the maximum
energy controls the full $H_0^1$ norm; Rellich compactness yields a
primal optimizer.  Equality in the primal-dual chain gives
\eqref{INITOPT-eq:complementarity}.

For log-concavity, add $\varepsilon\int|\nabla\psi|^2$ to the objective.
The resulting optimizer is the positive first Dirichlet eigenfunction
of a uniformly elliptic constant-coefficient operator.  A linear change
of variables reduces it to the Laplacian on a convex domain, so the
Brascamp--Lieb theorem~\cite{BL} makes the ground state log-concave.  Letting
$\varepsilon\downarrow0$ preserves log-concavity under the resulting
$L^1$ limit.  Central symmetry gives evenness by uniqueness of the
regularized positive ground state.
\end{proof}

\subsubsection{When the cosine product is optimal}
Put $Q_m=(-1,1)^m$, $e_0=\pi^2/4$, and
\[
 \psi_0(x)=\prod_{i=1}^m\cos(\pi x_i/2).
\]
For unit directions, $\mathcal E_v(\psi_0)=e_0$.

\begin{theorem}[Exact criterion for cosine optimality]
\label{INITOPT-thm:cosine}
Suppose the unit directions $v_j$ span $\R^m$.  Then
\begin{equation}\label{INITOPT-eq:cosine-criterion}
 \kappa_{Q_m}(V)=e_0
 \quad\Longleftrightarrow\quad
 \conv\{v_jv_j^{\mathsf T}:1\le j\le n\}
 \text{ contains a diagonal matrix}.
\end{equation}
If the convex hull contains no diagonal matrix, an even log-concave
nonproduct density has strictly smaller maximum directional Fisher
energy than every product density on $Q_m$.
\end{theorem}
\begin{proof}
If $M=\sum_j\lambda_jv_jv_j^{\mathsf T}$ is diagonal, then
$\tr M=1$ and the one-dimensional Dirichlet inequalities give
\[
 \int\nabla\psi^{\mathsf T}M\nabla\psi
 =\sum_iM_{ii}\int|\partial_i\psi|^2
 \ge e_0\tr M=e_0.
\]
The cosine product attains equality, so
Theorem~\ref{INITOPT-thm:spectral} gives $\kappa=e_0$.

Conversely, if $\kappa=e_0$, then $\psi_0$ is primal optimal.  For a
dual optimizer $M$, the complementarity relation makes $\psi_0$ a weak
eigenfunction of $-\operatorname{div}(M\nabla)$.  In the cube interior,
\[
 \frac{-\operatorname{div}(M\nabla\psi_0)}{\psi_0}
 =e_0\tr M-2e_0\sum_{i<\ell}M_{i\ell}
 \tan\frac{\pi x_i}{2}\tan\frac{\pi x_\ell}{2}.
\]
The left side is constant.  Setting all but two coordinates to zero
forces every off-diagonal entry of $M$ to vanish.

For a product amplitude, all cross derivative integrals vanish and
each marginal Dirichlet energy is at least $e_0$.  Hence every product
density has maximum directional energy at least $e_0$.  When the convex
hull contains no diagonal matrix, the spectral value is strictly below
$e_0$ and Theorem~\ref{INITOPT-thm:spectral} gives the even
log-concave nonproduct optimizer.
\end{proof}

The failure of the diagonal condition can be quantified.  Write
\[
 b_j=(v_{ij}v_{\ell j})_{i<\ell},\qquad
 \delta=\dist(0,\conv\{b_1,\ldots,b_n\}).
\]
Then $\delta>0$ precisely when the cosine product fails to be optimal.

\begin{theorem}[Explicit nonproduct energy gap]\label{INITOPT-thm:gap}
For unit spanning directions,
\begin{equation}\label{INITOPT-eq:gap}
 \kappa_{Q_m}(V)\le
 \frac{5e_0-\sqrt{9e_0^2+(4096/81)\delta^2}}2.
\end{equation}
\end{theorem}
\begin{proof}
If $\delta=0$, the cosine trial suffices.  Otherwise let $b_*$ be the
closest point of the convex hull to zero and put $h=b_*/\|b_*\|_2$.
Then $h\cdot b_j\ge\delta$.  Define
\[
 \varphi_h(x)=\sum_{i<\ell}h_{i\ell}\sin(\pi x_i)\sin(\pi x_\ell)
   \prod_{r\ne i,\ell}\cos(\pi x_r/2).
\]
Its product modes are orthonormal, $\|\varphi_h\|_2=1$, and it is
orthogonal to $\psi_0$.  The exact one-dimensional integral
\[
 \int_{-1}^1\cos(\pi x/2)\frac d{dx}\sin(\pi x)\,dx=\frac43
\]
gives
\[
 \langle\partial_v\psi_0,\partial_v\varphi_h\rangle
 =-\frac{32}{9}\sum_{i<\ell}h_{i\ell}v_iv_\ell.
\]
A direct orthogonality calculation bounds the directional energy of
$\varphi_h$ by $4e_0$.  Therefore the normalized trial
$(\psi_0+t\varphi_h)/\sqrt{1+t^2}$ has every directional energy at most
\[
 \frac{e_0-(64/9)\delta t+4e_0t^2}{1+t^2}.
\]
The smaller eigenvalue of this two-dimensional quadratic form is
\eqref{INITOPT-eq:gap}.  Replacing the trial by its absolute value leaves
all directional energies unchanged.
\end{proof}

Theorems~\ref{INITOPT-thm:cheeger} and \ref{INITOPT-thm:spectral}
separate two optimization scales.  The nonlocal problem is exact for
the rounding criterion, while the spectral problem identifies when a
matrix admits a strictly better smooth initial law before any scalar
auxiliary is optimized.

\begin{example}[Two directions with a common mixed derivative]
Let $0<\theta<\pi/4$ and take
$v_1=(\cos\theta,\sin\theta)$ and
$v_2=(\sin\theta,\cos\theta)$ in the square.
Both off-diagonal products equal $\frac12\sin(2\theta)>0$, so the
convex hull in \eqref{INITOPT-eq:cosine-criterion} contains no diagonal
matrix. Here $\delta=\frac12\sin(2\theta)$, and the explicit trial is
\[
 \psi_t(x)=
 \frac{\cos(\pi x_1/2)\cos(\pi x_2/2)
       +t\sin(\pi x_1)\sin(\pi x_2)}{\sqrt{1+t^2}}.
\]
The common sign of the mixed derivatives makes the linear term in $t$
negative for both directions. Formula~\eqref{INITOPT-eq:gap} therefore
improves the energy of every product density. The log-concave optimizer
follows from Theorem~\ref{INITOPT-thm:spectral}; the explicit
perturbation proves the quantitative upper bound.
\end{example}

The finite-direction criterion and the diffuse-direction theorem in
Section~\ref{SPECTRAL-app} use different quantifiers. Here the density is
chosen for a fixed set of columns. There the reference is chosen before
the adversary's diffuse direction, and the competing references may
have arbitrary dependence between their blocks.


\subsection{Spectral minimax and information capacity on open domains}\label{SPECTRAL-app}
The cube construction uses a product cosine reference. Optimizing over all joint densities asks for the smallest possible largest directional Dirichlet energy. This problem makes sense on every nonempty bounded open $K\subset\R^d$, including domains with irregular boundary or several components. Put
\begin{equation}\label{SPECTRAL-eq:lambda}
 \Lambda(K)=\inf_{\substack{u\in H^1_0(K),\ u\ge0\\\int_Ku^2=1}}
 \lambda_{\max}\!\left(\int_K\nabla u\nabla u\trans\right).
\end{equation}
The space $H^1_0(K)$ is the $H^1$ closure of $C_c^\infty(K)$, so its zero boundary condition is defined without boundary regularity. The spectral dual and capacity identity below require only boundedness and openness. The later overlap limit states its convexity hypothesis separately and keeps the block dimension $d$ fixed while the number of factors grows.

\begin{proposition}[An exact spectral dual]\label{SPECTRAL-prop:dual}
The infimum in \eqref{SPECTRAL-eq:lambda} is attained and
\begin{equation}\label{SPECTRAL-eq:dual}
 \Lambda(K)=\sup_{\substack{M\succ0\\\tr M=1}}
        \lambda_1(-\operatorname{div}(M\nabla);K).
\end{equation}
If $K=-K$, an even minimizer can be chosen. Convexity and connectedness are unnecessary.
\end{proposition}
\begin{proof}
Write $E(u)=\int\nabla u\nabla u\trans$. For normalized nonnegative
$u,v$ and $0\le t\le1$, the function
$w=\sqrt{tu^2+(1-t)v^2}$ belongs to $H^1_0(K)$, is normalized, and
Cauchy--Schwarz applied in each direction gives
$E(w)\preceq tE(u)+(1-t)E(v)$. Thus
$\mathcal E=\{S:S\succeq E(u)\text{ for some admissible }u\}$ is
convex and upward closed. It is closed: bounded $S_l$ bound the $H^1$ norms of the corresponding $u_l$. Extend these functions by zero to an enclosing cube. Rellich compactness there gives a subsequence converging strongly in $L^2$ and weakly in $H^1$. The zero-extension image of $H^1_0(K)$ is a closed linear subspace, so the limit still belongs to $H^1_0(K)$. Normalization survives, and weak lower semicontinuity of every directional energy gives $E(u)\preceq\lim S_l$. This proves closure and attainment on arbitrary bounded open $K$.

The matrix $\Lambda(K)I$ lies on the boundary of $\mathcal E$. A
supporting linear functional has the form $S\mapsto\tr(MS)$ with
$M\succeq0$, because $\mathcal E$ is upward closed. Normalize
$\tr M=1$. Then $\inf_u\tr(ME(u))=\Lambda(K)$.
For $M_\delta=(1-\delta)M+\delta I/d$, the variational eigenvalue is
at least $(1-\delta)\Lambda(K)$ and at most $\Lambda(K)$.
Let $\delta\downarrow0$ to prove \eqref{SPECTRAL-eq:dual}; the reverse
inequality follows by testing any fixed $u$. Symmetrizing $u^2$ under
$z\mapsto-z$ preserves or decreases $E$ and gives the even case.
\end{proof}

\begin{corollary}[The Smirnov--Vershynin Fisher-information capacity]\label{R20-cor:capacity}
For every nonempty bounded open $K\subset\R^d$, let
\[
 I(K)=\inf_f\left\|\int_{\{f>0\}}
          \frac{\nabla f\nabla f\trans}{f}\right\|_{\mathrm{op}}^{1/2},
\]
where the infimum is over smooth probability densities compactly supported in $K$. Then
\begin{equation}\label{R20-eq:capacity}
 I(K)^2=4\Lambda(K),
\end{equation}
with the matrix-valued Dirichlet dual \eqref{SPECTRAL-eq:dual}. In particular,
\begin{equation}\label{R21-eq:capacity-values}
 I((-R,R)^d)=\frac\pi R,\qquad
 I(RB_2^d)=\frac{2j_{d/2-1,1}}{\sqrt d\,R},
\end{equation}
where $j_{\nu,1}$ is the first positive zero of the Bessel function $J_\nu$.
The relaxation to densities $f=u^2$ with $u\in H^1_0(K)$ attains the infimum. Smooth compactly supported densities approximate that value.
\end{corollary}
\begin{proof}
For every admissible smooth $f$ with finite Fisher information, $u=\sqrt f$ belongs to $H^1_0(K)$ and has $\int u^2=1$. Indeed, the compactly supported functions $\sqrt{f+\delta}-\sqrt\delta$ converge to $u$ in $H^1$ as $\delta\downarrow0$, by dominated convergence with gradient bound $|\nabla f|/(2\sqrt f)$ on $\{f>0\}$. Their gradients give
$\mathcal I_f=4\int\nabla u\nabla u\trans$, hence $I(K)^2\ge4\Lambda(K)$.
Conversely, approximate a nonnegative normalized minimizer in $H^1_0(K)$ by nonnegative $C_c^\infty(K)$ functions $u_j$, normalized in $L^2$. Such approximation follows by taking positive parts and smoothing inside the compact supports. The densities $u_j^2$ have Fisher matrices $4E(u_j)\to4E(u)$ in operator norm. This proves equality and attainment in the stated relaxation.

For the cube, its cosine ground state has energy matrix $(\pi^2/(4R^2))I_d$; the dual test $M=I_d/d$ gives the matching lower bound. The positive radial ground state on the ball has energy matrix $(\lambda_1/d)I_d$, with $\lambda_1=j_{d/2-1,1}^2/R^2$. The same test proves the ball value in \eqref{R21-eq:capacity-values}.
\end{proof}
This identifies exactly the operator-norm capacity introduced by Smirnov--Vershynin~\cite{R20-SV}. The square-root identity is also used there; \eqref{R20-eq:capacity} optimizes the full directional matrix through its spectral dual. The product formula below holds for arbitrary nonempty bounded open factors:
\[
 I(K\times L)=\max\{I(K),I(L)\}.
\]
For convex $K,L$, Theorem~\ref{GEOM-thm:bm} gives
$I((1-t)K+tL)^{-1}\ge(1-t)I(K)^{-1}+tI(L)^{-1}$.
The convex-domain volume comparison and its quantitative stability are in Theorem~\ref{GEOM-thm:product-ball}. These statements concern the capacity and its Sobolev relaxation. In the Smirnov--Vershynin thinning theorem the destination is $2K$ and the rejection coefficient is $I(K)/2$. Thus the cube coefficient $\pi/(2R)$ corresponds to destination $[-2R,2R]^d$; fixing destination $[-R,R]^d$ gives $\pi/R$.

\paragraph{Disconnected domains can lower the largest directional cost.}
For a finite disjoint union $K=\bigcup_{j=1}^s K_j$ of bounded open components, the Dirichlet space splits into the direct sum of the component spaces. Consequently
\begin{equation}\label{R21-eq:disconnected-dual}
 \Lambda(K)=\sup_{M\succ0,\,\tr M=1}
      \min_{1\le j\le s}\lambda_1(-\operatorname{div}(M\nabla);K_j).
\end{equation}
The minimizing density can distribute mass between components with different expensive directions. For an explicit example, take $0<a<b$ and two disjoint translates $K_1,K_2$ of
$(-a,a)\times(-b,b)$ and $(-b,b)\times(-a,a)$.
Put $\alpha=\pi^2/(4a^2)$ and $\beta=\pi^2/(4b^2)$. Their normalized ground states have energy matrices $\diag(\alpha,\beta)$ and $\diag(\beta,\alpha)$. Giving the two squared ground states equal probability produces energy $(\alpha+\beta)I_2/2$. Testing \eqref{R21-eq:disconnected-dual} with $M=I_2/2$ proves the matching lower bound. Hence
\[
 \Lambda(K_1\cup K_2)=\frac{\alpha+\beta}{2}<\alpha
       =\Lambda(K_1)=\Lambda(K_2),\qquad
 I(K_1\cup K_2)=\frac\pi{\sqrt2}\sqrt{a^{-2}+b^{-2}}<\frac\pi a.
\]
This is a gain from averaging anisotropic energy matrices across components. It preserves the distinction between optimizing the full Fisher matrix and minimizing its trace.

\begin{proposition}[The capacity reduction from disconnected components]
\label{R22-prop:components}
Let $K=\bigcup_{j=1}^sK_j\subset\R^d$ be a finite disjoint union of nonempty bounded open sets, and write $I_j=I(K_j)$. Then
\begin{equation}\label{R22-eq:components}
 \max\left\{\frac{\min_j I_j}{\sqrt d},
             \left(\sum_{j=1}^s I_j^{-2}\right)^{-1/2}\right\}
 \le I(K)\le\min_j I_j.
\end{equation}
In particular, the universal factor $1/\sqrt{\min(s,d)}$ in
$I(K)\ge\min_jI_j/\sqrt{\min(s,d)}$ is sharp.
\end{proposition}
\begin{proof}
For normalized $u\in H_0^1(K)$, restrict to the components, let $p_j$ be their squared $L^2$ masses and normalize each nonzero restriction to $u_j$. If $E_j=\int_{K_j}\nabla u_j\nabla u_j^{\mathsf T}$ and $D=\lambda_{\max}(\sum_jp_jE_j)$, then
\[
 D\ge p_j\Lambda(K_j),\qquad
 dD\ge\sum_jp_j\tr E_j\ge\min_j\Lambda(K_j).
\]
Summing $p_j\le D/\Lambda(K_j)$ gives
$D\ge(\sum_j\Lambda(K_j)^{-1})^{-1}$. Infimizing and using $I^2=4\Lambda$ proves the lower bounds. Concentrating on the component of least capacity gives the upper bound.

For sharpness put $r=\min(s,d)$. Take translated boxes with half-width $a$ in one of the first $r$ coordinate directions and half-width $b>a$ in all others, with at least one box of every orientation; repeat orientations when $s>r$. Put $\alpha=\pi^2/(4a^2)$ and $\beta=\pi^2/(4b^2)$. Each box has $\Lambda=\alpha$. Equal mixing of one squared ground state in each orientation gives largest directional energy $(\alpha+(r-1)\beta)/r$. The dual matrix with diagonal $1/r$ in these $r$ directions and zero elsewhere gives the matching lower bound, by positive-definite approximation when $r<d$. Indeed the first Dirichlet energy in those directions is the same on every component. Thus the ratio of the union capacity to each component capacity tends to $1/\sqrt r$ as $b/a\to\infty$.
\end{proof}
The gain comes from averaging differently oriented energy matrices. The number of components and the ambient dimension bound how far this averaging can lower the largest directional energy.

\subsubsection{A Gaussian comparison uniform over dependent inputs}
The statistical comparison used below is independent of the unknown
input law. Let $K_0\subset\R^d$ be compact and let $g$ be a $C^2$
vector field on a fixed neighborhood of $K_0$. Suppose that, for
$M\succeq0$, $\tr M=1$, and $\lambda>0$,
\begin{equation}\label{SPECTRAL-eq:score}
 \tr(MDg)=\lambda+g\trans Mg
\end{equation}
pointwise there. Let $U$ have a symmetric law on unit vectors with
$\E UU\trans=M$. Given known $p_1,\ldots,p_N\in\R^r$ satisfying
$\sum_i p_ip_i\trans=I_r$, write $\delta=\max_i\|p_i\|_2$.
An arbitrary random input $X\in K_0^N$ is independent of independent
copies $U_i$. The observed experiment is
\[
 (U_i,Y_i)_{i\le N},\qquad
 Y_i=X_i+U_i\langle p_i,\theta\rangle,\qquad\|\theta\|_2\le T.
\]
The joint coordinates of $X$ may be dependent or singular.

Put $q(x)=g(x)\trans Mg(x)$ and, using only the observations, define
\begin{align*}
 S&=\sum_i p_ip_i\trans q(Y_i),&
 L&=\sum_i p_i\langle U_i,g(Y_i)\rangle,\\
 \widehat\theta&=(\lambda I+S)^{-1}L,&
 V(S)&=(\lambda I+S)^{-1}S(\lambda I+S)^{-1}.
\end{align*}
For fixed $\eta>0$, set $\sigma^2=(4\lambda)^{-1}+\eta^2$ and let
$G_r$ be an independent standard Gaussian. The output of the kernel is
\begin{equation}\label{SPECTRAL-eq:channel}
 Z=\widehat\theta+[\sigma^2I-V(S)]^{1/2}G_r.
\end{equation}
The square root is defined because $s/(\lambda+s)^2\le1/(4\lambda)$.
Outside the observation neighborhood the kernel can be extended
arbitrarily; those values have zero probability for sufficiently small
$\delta$.

\begin{theorem}[A distribution-free Gaussian channel]\label{SPECTRAL-thm:channel}
For fixed $K_0,g,M,\lambda,r,T,\eta$ there are constants $C_1,C_2$ such
that, for all sufficiently small $\delta$,
\begin{equation}\label{SPECTRAL-eq:tv}
 \sup_{\substack{\law(X)\text{ on }K_0^N\\\|\theta\|_2\le T}}
 \|\law(Z)-N(\theta,\sigma^2I_r)\|_{\TV}
 \le C_1\delta+C_2\delta^2.
\end{equation}
The same kernel works for every input law and parameter. Here total
variation is $\sup_A|P(A)-Q(A)|$.
\end{theorem}
\begin{proof}
Condition on the entire input $X=x$. Write
$a_i=\langle p_i,\theta\rangle$,
$S_0=\sum_i p_ip_i\trans q(x_i)$,
$D_0=\lambda I+S_0$, and
$L_0=\sum_i p_i\langle U_i,g(x_i)\rangle$.
Then $\E L_0=0$ and $\Cov(L_0)=S_0$. Taylor expansion and
\eqref{SPECTRAL-eq:score} give
\[
 L=L_0+D_0\theta+R_L,\qquad
 \|R_L\|_{L^2}\le
 \delta(2B_1T+B_2T^2/2),
\]
where $B_1,B_2$ bound the first two derivatives of $g$. Indeed, the
centered first-order error is the sum of independent vectors
$p_i a_i(U_i\trans Dg(x_i)U_i-\tr(MDg(x_i)))$, whose squared
$L^2$ norm is at most $4B_1^2\delta^2T^2$. The remaining terms are
bounded by $(B_2/2)\sum_i\|p_i\|a_i^2\le B_2\delta T^2/2$.
Similarly, if $Q_1,Q_2$ bound the derivatives of $q$, symmetry of $U_i$
cancels the first-order mean and yields $\|S-S_0\|_{L^2(\mathrm F)}\le \delta^2(Q_1T+Q_2T^2/2)$. For the centered part use
$\sum_i\|p_i\|^4 a_i^2\le\delta^4 T^2$; for the remainder use
$\sum_i\|p_i\|^2 a_i^2\le\delta^2 T^2$.

All inverse matrices have norm at most $\lambda^{-1}$. The inverse
identity and Cauchy--Schwarz therefore imply
\[
 \E\|\widehat\theta-\theta-D_0^{-1}L_0\|=O(\delta),\qquad
 \E\|V(S)-V(S_0)\|_{\mathrm F}=O(\delta^2),
\]
uniformly in $x,\theta$. Gaussian convolution with covariance at least
$\eta^2I$ is Lipschitz in its mean and covariance in total variation,
with constants depending only on $r,\eta$ and fixed covariance bounds.
It remains to compare
$\theta+D_0^{-1}L_0+[\sigma^2I-V(S_0)]^{1/2}G_r$ with its Gaussian
counterpart. Conditional on $x$, the summands of $D_0^{-1}L_0$ are
independent, symmetric, centered and bounded by a constant times
$\|p_i\|$. Replace them in sequence by centered Gaussians with the
same covariance. The first three Taylor terms of a Gaussian-smoothed
test agree; its fourth directional derivatives have bounded $L^1$
norm, of order $\eta^{-4}$ times the fourth power of the direction
norm. The total replacement error is therefore bounded by
$C\sum_i\|p_i\|^4\le Cr\delta^2$. The limiting covariance is exactly
\[
 \Cov(D_0^{-1}L_0)+\sigma^2I-V(S_0)=\sigma^2I.
\]
The bound holds uniformly for every conditional input $x$, so averaging
proves \eqref{SPECTRAL-eq:tv} for every joint input law.
\end{proof}
For a loss in $[0,1]$, any decision rule for the Gaussian location
experiment may be composed with this kernel; its risk changes by at
most the right side of \eqref{SPECTRAL-eq:tv}, uniformly in the
parameter and the nuisance input law. This is a comparison of the
whole statistical experiment, rather than of one particular overlap.

\subsubsection{Optimal overlap against diffuse directions}
For this overlap theorem, assume again that $K$ is bounded, open and convex with nonempty interior. Fix the auxiliary window $(-L,L)$, where $L>0$. The block dimension
and the window remain fixed. Let $\varepsilon_N\downarrow0$, with $N\varepsilon_N^2\ge1$, and put
\[
 \mathcal D_N=\left\{v=(v_1,\ldots,v_N):
 \sum_i\|v_i\|_2^2=1,\quad \max_i\|v_i\|_2\le\varepsilon_N\right\}.
\]
For a joint probability density $f$ of $(X,T)$ on $K^N\times(-L,L)$,
define
\[
 H_f(v)=\frac14\int_{-L}^L\int_{-L}^L\int_{\mathbb R^{dN}}
 \min\{f(y+sv,s),f(y+tv,t)\}\,dy\,ds\,dt.
\]
For $I>0$, let $\varphi_I$ be the density of $N(0,I^{-1})$, and set
\[
 \mathcal G_I(q)=\frac14\int_{-L}^L\int_{-L}^L\int_{\mathbb R}
 \min\{q(s)\varphi_I(z+s),q(t)\varphi_I(z+t)\}\,dz\,ds\,dt,
 \qquad \mathcal S_L(I)=\sup_q\mathcal G_I(q),
\]
where $q$ ranges over probability densities on $(-L,L)$.

\begin{theorem}[The spectral limit for arbitrary joint references]
\label{SPECTRAL-thm:all-joint}\label{SPECTRAL-thm:minimax}
With these definitions,
\[
 \lim_{N\to\infty}\ \sup_f\ \inf_{v\in\mathcal D_N}H_f(v)
       =\mathcal S_L(4\Lambda(K)).
\]
The supremum permits arbitrary dependence among the blocks of $X$ and
between $X$ and $T$. Independent references of the form
$(u_*^2)^{\otimes N}\otimes q$, where $u_*$ minimizes the definition of
$\Lambda(K)$, asymptotically attain the value after optimization over $q$.
If $K=-K$, requiring both reference blocks to have mean zero leaves the
limit unchanged. Additional dependent auxiliary coordinates cannot
increase the individual height above the corresponding two-block height.
\end{theorem}

\begin{proof}
For the upper bound, choose an elliptic matrix $M$ with
$\operatorname{tr}M=1$ and a positive Dirichlet ground state on a fixed
slight enlargement of $K$. Write $\lambda>0$ for its eigenvalue and fix
$\eta>0$. The preceding Gaussian-channel theorem constructs a Markov
kernel, independent of the input law and the translation parameter, with
output variance $\sigma^2=(4\lambda)^{-1}+\eta^2$. For the scalar design $p_i=N^{-1/2}$, its total-variation error is a
number $e_N\to0$, uniformly over every law on $K^N$ and every translation
parameter in $[-L,L]$. Let $U_i$ be the independent symmetric unit
vectors used by this kernel, with
$\mathbb E U_iU_i^{\mathsf T}=M$. The direction
$v_U=(U_1,\ldots,U_N)/\sqrt N$ belongs to $\mathcal D_N$.

Let $q$ be the marginal density of $T$. For almost every $s$, conditional
on $T=s$, the input $X$ has an arbitrary probability law on $K^N$.
Applying the same kernel to $(U,X-sv_U)$ therefore gives a law $P_s'$
satisfying $\|P_s'-N(-s,\sigma^2)\|_{\rm TV}\le e_N$. This assertion uses the uniformity over the input law: the conditional
law is allowed to change with $s$.

For finite measures $\mu,\nu$, write
$O(\mu,\nu)=\int\min\{d\mu,d\nu\}$. Observing the independent direction
$U$ makes the original overlap the average of the overlaps at $v_U$.
A Markov kernel cannot decrease overlap. Moreover, perturbing probability
laws $P,Q$ by total-variation errors at most $e_N$ changes
$O(\alpha P,\beta Q)$ by at most $(\alpha+\beta)e_N$. Consequently
\[
 \mathbb E_U H_f(v_U)
 \le \mathcal G_{\sigma^{-2}}(q)
       +\frac{e_N}{4}\int_{-L}^L\int_{-L}^L(q(s)+q(t))\,ds\,dt
 =\mathcal G_{\sigma^{-2}}(q)+Le_N.
\]
The bound is uniform over all joint densities $f$. Hence
\[
 \limsup_N\sup_f\inf_{v\in\mathcal D_N}H_f(v)
 \le\mathcal S_L(\sigma^{-2}).
\]
Take $N\to\infty$ first. Then let $\eta\downarrow0$, let the enlarged
body decrease to $K$, and choose elliptic $M$ approaching the spectral
dual value
\[
 \Lambda(K)=\sup_{\substack{M\succ0\\\operatorname{tr}M=1}}
 \lambda_1(-\operatorname{div}(M\nabla);K).
\]
This gives $\sigma^{-2}\to4\Lambda(K)$. Total variation between the
corresponding centered Gaussian laws tends to zero. The same weighted
risk estimate, integrated over the fixed auxiliary window, makes the
convergence of $\mathcal G_I(q)$ uniform in $q$; thus the auxiliary
supremum also passes to the limit.

For the lower bound set $f_*=u_*^2$. Translations of this density are
differentiable in quadratic mean because $u_*\in H^1(\R^d)$ after
zero extension. Its information matrix is
$J_*=4E(u_*)\preceq4\Lambda(K)I$. For any sequence of directions in
$\mathcal D_N$, the scalar information is
$I_N=\sum_i v_i\trans J_*v_i\le4\Lambda(K)$. Along a subsequence let
$I_N\to I$. The product translation experiment converges to the
Gaussian location experiment of information $I$. To verify the needed
uniformity, the quadratic-mean remainder is
$o(\|h\|^2)$ uniformly over small vector translations $h$, and the sum
of these remainders is $o(1)$ because $\max_i\|v_i\|\to0$ and
$\sum_i\|v_i\|^2=1$. If $\ell$ is the score of $f_*$, the Lindeberg
remainder is bounded by
\[
 \E\bigl[\|\ell(X)\|^2
       \mathbf1_{\{\|\ell(X)\|>c/\varepsilon_N\}}\bigr]\longrightarrow0
\]
for each $c>0$, since $\ell\in L^2(f_*)$.
Expanding the square-root likelihoods, followed by this triangular-array
central limit theorem, gives the Gaussian likelihood ratios and
contiguity. Weighted binary testing risks therefore converge to the
Gaussian risks for each fixed pair of shifts. Their bound by
$\min(q(s),q(t))$ allows integration for a fixed bounded auxiliary.
Since Gaussian overlap decreases with information, every such
subsequence has height limit at least $\mathcal G_{4\Lambda(K)}(q)$.
The sequential argument proves uniformity over $\mathcal D_N$.
Bounded auxiliaries approximate arbitrary densities in $L^1$, and
$|\mathcal G_I(q)-\mathcal G_I(\widetilde q)|\le
L\|q-\widetilde q\|_1$. Taking their supremum proves the lower bound.

For symmetric $K$, the minimizing $u_*$ can be even. Reflecting and
averaging an auxiliary preserves support and makes it even; concavity
of the weighted overlap in $q$ shows that this cannot decrease its
value. Thus centered references attain the same lower limit, while the
upper bound already covered every law.
Finally, write a larger joint density as $f(x,w,t)$, where $w$ denotes
additional auxiliary coordinates. At each fixed $y,s,t$,
\[
 \int\min\{f(y+sv,w,s),f(y+tv,w,t)\}\,dw
 \le\min\left\{\int f(y+sv,w,s)\,dw,
                   \int f(y+tv,w,t)\,dw\right\}.
\]
Integrating gives the asserted comparison with the marginal reference
$(X,T)$. This completes the proof.
\end{proof}
The formula optimizes the initial overlap criterion with the reference
chosen before the direction. It leaves finite-dimensional discrepancy
optimization, where the density may depend on the whole matrix, as a
separate problem. The product conclusion allows fully dependent
competitors in its upper bound.


\subsubsection{The sharp reference threshold for Koml\'os rounding}
\label{R7-sec:threshold}
For $K=(-R,R)$ and auxiliary interval $(-3,3)$,
Theorem~\ref{SPECTRAL-thm:all-joint} evaluates the exact diffuse
minimax overlap as
\begin{equation}\label{R7-eq:diffuse-value}
 \mathcal S_3(\pi^2/R^2)=M(\pi/(2R)),\qquad
 M(k)=\sup_{q\text{ on }(-3,3)}\G_k(q).
\end{equation}
The notation $\G_k$ here is the Gaussian likelihood overlap of
Appendix~\ref{app:overlap}: its Fisher parameter is $4k^2$.
The supremum permits every auxiliary density. The competing initial
reference in the spectral theorem permits arbitrary dependence among
all physical blocks and the auxiliary, chosen before the direction.
Define
\[
 R_{\rm diff}=\inf\{R>0:M(\pi/(2R))\ge1\}.
\]
For a common independent auxiliary in all unit directions, define
\[
 A_m(R)=\sup_{f,q}\inf_{\|v\|_2=1}H_f(v;q),\qquad
 R_{\rm com}=\inf\{R>0:\limsup_{m\to\infty}A_m(R)\ge1\},
\]
where $f$ ranges over centered densities in $(-R,R)^m$ and $q$ over
centered densities in $(-3,3)$. Both may depend on $m$.

\begin{theorem}[Sharp calibration of the universal reference]
\label{R7-thm:method-constant}
The diffuse reference problem is exactly characterized by
\eqref{R7-eq:diffuse-value}. Its threshold and the common-auxiliary
threshold satisfy
\begin{equation}\label{R7-eq:threshold-bracket}
 6.8383231851<R_{\rm diff}\le R_{\rm com}<6.8383231852.
\end{equation}
For the fixed density $q_\dagger$, the cosine product attains the
limiting worst-direction value $\G_{\pi/(2R)}(q_\dagger)$ in the
range $6\le R\le7$. In particular, an arbitrary dependent physical
reference with that same auxiliary has no larger limiting value.
\end{theorem}
\begin{proof}
The spectral theorem with $K=(-R,R)$ gives
\eqref{R7-eq:diffuse-value}; its upper proof also fixes any chosen
auxiliary instead of taking the final supremum. Restricting to the
equal-magnitude diffuse directions gives
$\limsup_m A_m(R)\le M(\pi/(2R))$.
For the fixed $q_\dagger$, the forward likelihood comparison in
Appendix~\ref{app:overlap} gives the matching lower bound from the
cosine product in every dimension and every unit direction.

The Gaussian overlap is homogeneous and concave in its density
argument. Its supporting hyperplane at a positive density $q$ gives
\begin{equation}\label{R7-eq:dual-all-q}
 \G_k(p)\le\int p(s)D_q^k(s)\,ds\le\sup_sD_q^k(s),\qquad
 \G_k(q)=\int q(s)D_q^k(s)\,ds.
\end{equation}
Indeed, the derivative of the overlap of two translated normal
weights $A,B$ with respect to $A$ is
$\Phi(-a-\log(A/B)/(2a))$. Integrating its support-plane inequality
proves \eqref{R7-eq:dual-all-q}, including competing densities with
zeros. Also $M$ is nonincreasing and locally Lipschitz, with
$|M(k)-M(k')|\le9|k-k'|$, by the derivative bound for the normal
kernel.

Put $C_-=6.8383231851$, $C_+=6.8383231852$, and $k_\pm=\pi/(2C_\pm)$.
The rational whole-interval bounds of
Appendix~\ref{R7-app:dual-height} give
\[
 M(k_-)<1,\qquad \G_{k_+}(q_\dagger)>1+10^{-12}.
\]
Continuity preserves the first strict inequality at radii slightly
larger than $C_-$, and the second at radii slightly smaller than
$C_+$. The forward comparison remains valid there because the
logarithmic slope condition is strict. This proves
\eqref{R7-eq:threshold-bracket}.
\end{proof}
The thresholds refer to a universal reference selected before the
direction. The finite-direction optimization in
Theorem~\ref{INITOPT-thm:cosine} instead chooses the density for the
given columns; the cubic correction in
\eqref{MAIN-eq:radius} also uses their geometry.
The decimal $C=\Cstar$ is a certified finite construction radius.
The exact extremal object settled here is the diffuse reference
problem \eqref{R7-eq:diffuse-value}, with the enclosure
\eqref{R7-eq:threshold-bracket}; a globally optimal Koml\'os
constant would optimize over all signing constructions.

\subsubsection{Sharpness of the first-order variation certificate}
Let $|Df|$ be the total variation of the full distributional gradient of
a density extended by zero outside the cube. For a unit vector $v$,
write $|D_vf|$ for its directional variation. Define
\[
 \tau_m(R)=\inf_f\sup_{\|v\|_2=1}|D_vf|,
\]
where the infimum ranges over probability densities supported in
$[-R,R]^m$ whose zero extensions belong to $BV(\R^m)$.

\begin{theorem}[All-density first-order optimum]\label{M-first-order}
For each $R>0$,
\begin{equation}\label{M-first-order-limit}
 \lim_{m\to\infty}\tau_m(R)=\frac{\sqrt{2\pi}}R.
\end{equation}
More explicitly, if $0<k<\pi/(2R)$ and $m>\sec^2(kR)$, then
\begin{equation}\label{M-first-order-finite}
 \tau_m(R)\ge
 \frac{\Gamma(m/2)}{\sqrt\pi\Gamma((m+1)/2)}
       2k\sqrt{m-\sec^2(kR)}.
\end{equation}
Thus the sufficient condition
$H_f(v;q_0)\ge3/2-(3/2)|D_vf|$, with $q_0$ uniform on $(-3,3)$,
has sharp dimension-uniform radius $3\sqrt{2\pi}$ over all initial
densities, including dependent ones.
\end{theorem}
\begin{proof}
Put $G(x)=(k\tan(kx_i))_i$ and $M=k^2\sec^2(kR)$. Then
$\operatorname{div}G=mk^2+|G|^2$ and $0\preceq DG\preceq MI$.
For $\delta>0$, the vector field
$W=G/(|G|^2+\delta)^{1/2}$ has norm at most one and
\[
 \operatorname{div}W
 \ge\frac{|G|^2+mk^2-M}{\sqrt{|G|^2+\delta}}
 \ge2\sqrt{mk^2-M-\delta},
\]
whenever the last radicand is positive. Extend this field smoothly past
the cube and use a cutoff equal to one on it. Distributional integration
by parts against $f$ gives
$|Df|\ge\int f\operatorname{div}W\ge2\sqrt{mk^2-M-\delta}$.
Let $\delta$ decrease to zero.

The polar decomposition $Df=\zeta\,|Df|$ has $|\zeta|=1$ almost
everywhere for $|Df|$. Averaging over a uniform unit vector $v$ yields
\[
 \sup_{\|v\|_2=1}|D_vf|
 \ge\E_v|D_vf|
 =c_m|Df|,\qquad
 c_m=\E|v_1|=\frac{\Gamma(m/2)}{\sqrt\pi\Gamma((m+1)/2)}.
\]
This proves \eqref{M-first-order-finite}. Since
$c_m\sqrt m\to\sqrt{2/\pi}$, first send $m$ to infinity and then
$k$ to $\pi/(2R)$. The cosine score estimate
\eqref{F-eq:cosine-score-L1} gives the matching upper bound in every
dimension. The final assertion follows from the requirement
$\sup_v|D_vf|\le1/3$ in the first-order certificate.
\end{proof}


\subsubsection{Products, volume and Minkowski addition}
\label{GEOM-sec}

Products depend on the largest directional energy; volume comparison
uses its trace. These two calculations determine how support geometry
changes the optimal overlap.

\begin{theorem}[Products and the ball of prescribed volume]
\label{GEOM-thm:product-ball}
Let $K_i\subset\R^{d_i}$ be nonempty bounded open sets.  Then
\begin{equation}\label{GEOM-eq:product}
 \Lambda(K_1\times\cdots\times K_s)=\max_{1\le i\le s}\Lambda(K_i).
\end{equation}
For a bounded open convex $K\subset\R^d$ and a ball $B$ of the same
volume,
\begin{equation}\label{GEOM-eq:faber}
 \Lambda(K)\ge\frac{\lambda_1(K)}d
 \ge\frac{\lambda_1(B)}d=\Lambda(B).
\end{equation}
Equality between the two outer terms holds precisely when $K$ is a
translate of $B$.  If $d\ge2$, there is $c_d>0$ such that
\begin{equation}\label{GEOM-eq:quantitative}
 \Lambda(K)\ge\Lambda(B)\bigl(1+c_d\mathcal A(K)^2\bigr),\qquad
 \mathcal A(K)=\inf_{z\in\R^d}\frac{|K\mathbin\triangle(z+B)|}{|K|}.
\end{equation}
Consequently, among convex supports of a prescribed volume, a ball maximizes
the universal diffuse-block overlap limit in
Theorem~\ref{SPECTRAL-thm:all-joint}.
\end{theorem}
\begin{proof}
It suffices to prove the product formula for two factors.  For normalized
$u\in H^1_0(K_1\times K_2)$, let
\[
 u_1(x)=\left(\int_{K_2}u(x,y)^2\,dy\right)^{1/2}.
\]
Then $u_1\in H^1_0(K_1)$ and $\|u_1\|_2=1$.  To justify the boundary
statement, approximate $u$ in $H^1$ by compactly supported smooth
functions, use the contraction of the fibrewise $L^2$ norm in $L^2$,
and the derivative bound below to pass weakly in $H^1_0$.
For every fixed direction $a$ in the first factor, Cauchy--Schwarz gives
\[
 |\partial_a u_1(x)|^2
 \le\int_{K_2}|\partial_a u(x,y)|^2\,dy.
\]
Thus the first diagonal block of the energy matrix of $u$ dominates the
energy matrix of $u_1$.  Its largest eigenvalue is at least
$\Lambda(K_1)$; the second marginal gives $\Lambda(K_2)$.
For the reverse inequality take a product $u(x,y)=v(x)w(y)$ of
normalized minimizing functions.  The off-diagonal energy block vanishes,
since $\int v\nabla v=\int w\nabla w=0$ for zero-boundary Sobolev
functions.  The two diagonal blocks are exactly the separate energy
matrices.  This proves \eqref{GEOM-eq:product}.

For any normalized $u$ on $K$, the largest eigenvalue of its energy
matrix is at least $d^{-1}\int|\nabla u|^2$.  Infimizing proves the first
inequality in \eqref{GEOM-eq:faber}.  The second is Faber--Krahn.  On a
ball, the positive radial ground state has energy matrix
$\lambda_1(B)I_d/d$, giving equality.  The classical equality case in
Faber--Krahn characterizes the ball.  Convexity removes the distinction
between equality up to a null set and equality of these open domains.
The quantitative Faber--Krahn inequality of Brasco, De Philippis and
Velichkov \cite{BrascoDPV} gives \eqref{GEOM-eq:quantitative} after division
by $d$.  Finally the Gaussian testing functional $\mathcal S_L(I)$ is
nonincreasing in its information parameter $I$: adding Gaussian noise
is a Markov kernel and increases weighted overlap.  Apply the minimax
identity with $I=4\Lambda(K)$.
\end{proof}

The product formula allows full dependence in the density on the product.
Its lower bound marginalizes that density before making any product
choice.  In particular,
\[
 \Lambda\!\left(\prod_{i=1}^d(-R_i,R_i)\right)
       =\max_i\frac{\pi^2}{4R_i^2}.
\]
Together with \eqref{R21-eq:capacity-values}, this shows how the interval scales and ball radius enter the same capacity. Thin factors
therefore determine the universal information cost of a rectangular
support, whereas the ball spreads the minimal energy equally among all
directions.

\begin{theorem}[Brunn--Minkowski inequality for the reference energy]
\label{GEOM-thm:bm}
For bounded open convex $K,L\subset\R^d$ with nonempty interior and
$0\le t\le1$,
\begin{equation}\label{GEOM-eq:bm}
 \Lambda((1-t)K+tL)^{-1/2}
 \ge(1-t)\Lambda(K)^{-1/2}+t\Lambda(L)^{-1/2}.
\end{equation}
The quantity $\Lambda^{-1/2}$ is homogeneous of degree one and invariant
under translations and orthogonal transformations.
\end{theorem}
\begin{proof}
For $M\succ0$ with $\tr M=1$, write
$\lambda_M(K)=\lambda_1(-\operatorname{div}(M\nabla);K)$.
The linear change of variables $x=M^{1/2}y$ identifies this eigenvalue
with the ordinary Dirichlet eigenvalue of $M^{-1/2}K$.  The
Brascamp--Lieb Brunn--Minkowski inequality for the first eigenvalue
\cite{BL,Colesanti} therefore gives
\[
 \lambda_M((1-t)K+tL)^{-1/2}
 \ge(1-t)\lambda_M(K)^{-1/2}+t\lambda_M(L)^{-1/2}.
\]
The spectral duality theorem gives
$\Lambda(K)^{-1/2}=\inf_{M\succ0,\tr M=1}\lambda_M(K)^{-1/2}$.
Take the infimum on the left and bound each term on the right by its
own infimum.  This proves \eqref{GEOM-eq:bm}; the invariances follow by
change of variables in the energy.
\end{proof}

These conclusions apply to the reference-before-direction limit in
Theorem~\ref{SPECTRAL-thm:all-joint}.


\section{Asymptotic equality in Cheeger and subexponential spectra}
\label{R11-sec:spectra}
The Dirichlet Cheeger inequality gives $\lambda_{1,p}^{1/p}\ge h_D/p$. Its geometric side uses sets with their full boundary perimeter, including the part on the domain boundary; relative perimeter would describe a different problem. The divergence and max-flow descriptions of the Cheeger constant~\cite{Grieser} are the starting point for a matching lower estimate on large products.

For a weighted ground state $u$, the score $v=-\nabla\log u$ satisfies $\operatorname{div}_w v=\lambda_{1,2}+|v|^2$. Product factors add the eigenvalue term while their derivative bound stays fixed. Normalizing the resulting field therefore gives a sharp perimeter lower bound as the total energy grows. For the spectral upper bound, $\lceil\log_2 k_N\rceil$ factors encode disjoint support choices for $k_N$ trial functions; all other factors retain near-ground-state energy. This cost is negligible when $\log k_N=o(N)$, giving $h_D\sim2\sqrt{E_N}$ and $\lambda_{k_N,p}^{1/p}\sim(2/p)\sqrt{E_N}$ for fixed $p>1$. Rectangles permit changing side lengths when no coordinate dominates the energy. The subsequent heat-mode construction applies exact deterministic realization to the resulting finite observable spaces.

\subsection{Dirichlet spectra of weighted Cartesian products}
\label{FP-sec:weighted}\label{SP-sec:spectrum}
Briani--Buttazzo--Prinari prove that the optimal dimension-dependent ratio between the first $p$-Dirichlet eigenvalue and Cheeger constant tends to $1/p$ as dimension increases~\cite{R20-BBP}. Here the asymptotic holds along every sequence of products from a fixed weighted family and throughout all subexponential genus indices. The factor family fixes the uniform local score bounds; its proportions may vary arbitrarily.

Let $\Omega_j\subset\R^{d_j}$, $1\le j\le J$, be a fixed finite family
of bounded connected Lipschitz domains. Let $d\mu_j=w_j\,dx$, where
$w_j$ is smooth and strictly positive on a neighborhood of
$\overline\Omega_j$. No normalization of $w_j$ is required. For an
arbitrary sequence $a_{i,N}\in\{1,\ldots,J\}$, put
\[
 D_N=\prod_{i=1}^N\Omega_{a_{i,N}},\qquad
 \mu_N=\bigotimes_{i=1}^N\mu_{a_{i,N}},\qquad
 E_N=\sum_{i=1}^N\lambda_{1,2}(\Omega_{a_{i,N}},\mu_{a_{i,N}}).
\]
Here $\lambda_{1,2}$ is the first eigenvalue of the weighted Dirichlet
Laplacian. The full Dirichlet Cheeger constant is
\[
 h_{\mu,D}(\Omega)=\inf_{\mu(E)>0}
 \frac{\Per_\mu(E)}{\mu(E)},\qquad
 \Per_\mu(E)=\int_{\partial^*E}w\,d\mathcal H^{d-1}.
\]
The infimum is over finite-perimeter subsets of $\Omega$; perimeter
on $\partial\Omega$ is included. For $1<p<\infty$,
$\lambda_{k,p}(D,\mu)$ denotes the Krasnosel'skii-genus sequence
\[
 \inf_{\substack{A\subset W^{1,p}_0(D)\text{ compact, symmetric}\
             \|u\|_{L^p(\mu)}=1\ (u\in A),\ \gamma(A)\ge k}}
       \sup_{u\in A}\int_D|\nabla u|^p\,d\mu.
\]
The genus $\gamma(A)$ is the least $q$ admitting an odd continuous
map $A\to\R^q\setminus\{0\}$. In particular, a $k$-dimensional
subspace has a unit sphere of genus $k$.

\begin{theorem}[Mixed weighted products]
\label{R9-thm:mixed-weighted}
\label{FP-thm:weighted-cheeger}\label{FP-thm:weighted-p-spectrum}
\label{SP-thm:mixed-cheeger}\label{SP-thm:subexp-spectrum}
For the fixed family above and every sequence of factor types,
\begin{equation}\label{SP-eq:mixed-cheeger}
 \frac{h_{\mu_N,D}(D_N)}{2\sqrt{E_N}}\longrightarrow1.
\end{equation}
For each fixed $1<p<\infty$ and every sequence $k_N\ge1$ with
$\log k_N=o(N)$,
\begin{equation}\label{SP-eq:subexp-spectrum}
 \frac{\lambda_{k_N,p}(D_N,\mu_N)^{1/p}}{\sqrt{E_N}}
 \longrightarrow\frac2p.
\end{equation}
No limiting proportions of factor types are assumed. In particular,
for a single weighted domain,
$\lambda_{k_N,p}(\Omega^N,\mu^N)^{1/p}
\sim(2/p)\sqrt{N\lambda_{1,2}(\Omega,\mu)}$.
\end{theorem}

The index condition also gives a uniform statement: for every positive sequence $a_N\to0$, the asymptotic in \eqref{SP-eq:subexp-spectrum} is uniform over $1\le k\le\max\{1,\lfloor e^{a_NN}\rfloor\}$. Indeed, the first and last indices have the same asymptotic and trap every intermediate eigenvalue by monotonicity. The proof makes the scale visible: its lower bound is independent of $k$, while the upper construction changes only $O(\log k)$ factors. The remaining factors retain the full ground-state energy.

\subsubsection{Lower bounds from the ground-state field}
\begin{lemma}[Weighted divergence certificate]
\label{SP-lem:divergence}\label{FP-lem:weighted-field}
Let $d\mu=e^{-V}\,dx$ near the closure of a bounded domain $D$.
Suppose a smooth vector field $v$ satisfies
\[
 \diver_\mu v=a+|v|^2,\qquad Dv\preceq MI,\qquad a>M\ge0,
 \quad \diver_\mu v=e^V\diver(e^{-V}v).
\]
Then $h_{\mu,D}(D)\ge2\sqrt{a-M}$.
\end{lemma}
\begin{proof}
For $\tau>0$, set $s=|v|^2+\tau$ and $W=v/\sqrt s$. Then
$|W|\le1$ and
\[
 \diver_\mu W
 =\frac{a+|v|^2}{\sqrt s}
   -\frac{\langle Dv\,v,v\rangle}{s^{3/2}}
 \ge\sqrt s+\frac{a-M-\tau}{\sqrt s}
 \ge2\sqrt{a-M-\tau}.
\]
Weighted Gauss--Green bounds the integral over any finite-perimeter
$E\subset D$ by $\Per_\mu(E)$. Let $\tau\downarrow0$. This is the
continuous max-flow bound with the full ambient perimeter
\cite{Grieser,Leonardi}.
\end{proof}

\begin{lemma}[Outer approximation of the weighted eigenvalue]\label{FP-lem:outer-domain}
There are bounded smooth connected domains $\Omega_j\supset\overline\Omega$, contained in the neighborhood on which $w$ is defined, with
\[
 \lambda_{1,2}(\Omega_j,\mu)\longrightarrow\lambda_{1,2}(\Omega,\mu).
\]
\end{lemma}
\begin{proof}
Take smooth outer approximations contained in shrinking neighborhoods of $\overline\Omega$, all inside a fixed bounded smooth domain. Normalize their first eigenfunctions in $L^2(\mu)$ and extend them by zero to this fixed domain. Domain inclusion bounds their energies from above by $\lambda_{1,2}(\Omega,\mu)$. Since the weight is bounded above and below, weak $H^1$ compactness and strong $L^2$ compactness apply. Every subsequential limit has norm one and vanishes outside $\overline\Omega$. The Lipschitz-domain zero-extension characterization places its restriction in $H^1_0(\Omega)$. Lower semicontinuity of energy gives the reverse limiting eigenvalue inequality. Domain monotonicity gives the other direction.
\end{proof}

\begin{proof}[Proof of the Cheeger limit in Theorem~\ref{R9-thm:mixed-weighted}]
Fix $\eta>0$. By Lemma~\ref{FP-lem:outer-domain}, choose for each
factor a smooth outer domain with ground-state eigenvalue
$\lambda_j'\ge(1-\eta)\lambda_{1,2}(\Omega_j,\mu_j)$.
Write $w_j=e^{-V_j}$ and let $u_j>0$ be the outer ground state.
The field $v_j=-\nabla\log u_j$ is smooth on a neighborhood of
$\overline\Omega_j$ and satisfies
\[
 \diver_{\mu_j}v_j
 =-(\Delta-\nabla V_j\cdot\nabla)\log u_j
 =\lambda_j'+|v_j|^2.
\]
The finite family has a common bound $Dv_j\preceq M_\eta I$.
Concatenate these fields over the $N$ factors. The divergence identity
has $a\ge(1-\eta)E_N$, while the derivative bound remains
$M_\eta$. Lemma~\ref{SP-lem:divergence} gives
\[
 h_{\mu_N,D}(D_N)\ge2\sqrt{(1-\eta)E_N-M_\eta}.
\]
The energies $E_N$ grow linearly with $N$, so first let $N\to\infty$
and then $\eta\downarrow0$.

For the upper bound take normalized ground states $\psi_j$ on the
original domains and let $f_N=\prod_i\psi_{a_{i,N}}^2$.
Its integral against $\mu_N$ is one and its zero extension has
integrable weak gradient. Coarea and Cauchy--Schwarz give
\[
 h_{\mu_N,D}(D_N)\le\int|\nabla f_N|\,d\mu_N
 \le\left(\int f_N|\nabla\log f_N|^2\,d\mu_N\right)^{1/2}
 =2\sqrt{E_N}.
\]
Smooth $H^1_0$ approximation justifies the expression at boundary
zeros. This proves \eqref{SP-eq:mixed-cheeger}.
\end{proof}

\subsubsection{Subexponentially many disjoint trial functions}
\begin{lemma}[The $p$-Cheeger bound and smooth trial densities]
\label{SP-lem:p-cheeger}\label{SP-lem:trials}
For the weighted Dirichlet problem,
$\lambda_{k,p}^{1/p}\ge h_{\mu,D}/p$.
For every $\xi>0$, a factor ground state can be approximated by a
nonnegative normalized $\psi\in C_c^\infty(\Omega)$ satisfying
\[
 \int|\nabla\psi|^2d\mu\le\lambda_{1,2}(\Omega,\mu)+\xi,
 \qquad
 \int\psi^2|\nabla\log\psi|^q d\mu<\infty
 \quad(1\le q<\infty).
\]
Two such normalized densities with disjoint compact supports also
exist, without the energy approximation requirement.
\end{lemma}
\begin{proof}
Coarea applied to $|u|^p$, followed by H\"older, gives
\[
 h_{\mu,D}\int|u|^p d\mu
 \le p\left(\int|u|^p d\mu\right)^{(p-1)/p}
       \left(\int|\nabla u|^p d\mu\right)^{1/p}.
\]
Density in $W^{1,p}_0$ proves the lower bound. To obtain $\psi$,
approximate the nonnegative ground state in $H^1_0$, enclose the
approximant in a smooth compact interior domain, and add a small
positive bump there. Choose the bump proportional to $e^{-1/d(x)^2}$
near its boundary, with $d$ a smooth defining function. The sum is
positive in the interior and exponentially flat at the boundary.
Exponential decay dominates every polynomial in the logarithmic
derivatives, and all positive powers used below have the required
zero trace. Normalization and a sufficiently small addition retain
the energy bound. Taking two disjoint interior balls proves the last claim.
\end{proof}

\begin{proof}[Proof of the $p$-spectrum limit in Theorem~\ref{R9-thm:mixed-weighted}]
The lower bound follows from the preceding lemma and the Cheeger
limit. For the upper bound choose, for each factor type $j$, a
near-minimal trial $\psi_{j,0}$ and disjoint trials
$\psi_{j,+},\psi_{j,-}$ as in Lemma~\ref{SP-lem:trials}.
Let $\ell_N=\lceil\log_2k_N\rceil$, with $\ell_N=0$ if $k_N=1$.
For every word $w\in\{+,-\}^{\ell_N}$, put
\[
 F_w=\left(\prod_{i\le\ell_N}\psi_{a_{i,N},w_i}
             \prod_{i>\ell_N}\psi_{a_{i,N},0}\right)^{2/p}.
\]
These unit $L^p(\mu_N)$ functions have disjoint supports. If
$X_i$ has the corresponding density $\psi^2d\mu$ and
$Z_i=|\nabla\log\psi(X_i)|^2$, then
\begin{equation}\label{SP-eq:trial-energy}
 \int|\nabla F_w|^p d\mu_N
 =\left(\frac2p\right)^p
       \E\left(\sum_{i=1}^NZ_i\right)^{p/2}.
\end{equation}
The $\ell_N=o(N)$ modified factors contribute $o(N)$ in the required
moment, uniformly in $w$. The remaining variables belong to a fixed
finite family. Grouping them by type and applying the mean convergence
theorem shows that their centered sum divided by $N$ tends to zero in
$L^{p/2}$ for $p\ge2$, and in $L^1$ for $p<2$. Types occurring only
$o(N)$ times are controlled by truncation. Thus the largest quotient
in \eqref{SP-eq:trial-energy}, divided by $E_N^{p/2}$, has limsup
at most
\[
 \left(\frac2p\right)^p
 \left(1+\frac{\xi}{\min_j\lambda_{1,2}(\Omega_j,\mu_j)}\right)^{p/2}.
\]
On the span of the $F_w$, disjoint support makes every Rayleigh
quotient a weighted average of the individual quotients. The unit
sphere of any $k_N$-dimensional subspace has genus $k_N$ and is an
admissible test. Let $\xi\downarrow0$.
\end{proof}

\subsubsection{Rectangles with varying side lengths}
Identical or finitely many fixed factors are unnecessary for intervals. Let
\[
 Q_N=\prod_{i=1}^N(-R_{i,N},R_{i,N}),\quad
 e_{i,N}=\frac{\pi^2}{4R_{i,N}^2},\quad
 E_N=\sum_i e_{i,N},\quad \rho_N=\frac{\max_i e_{i,N}}{E_N}.
\]
\begin{theorem}[Anisotropic rectangles]\label{SP-thm:rectangles}
If $\rho_N\to0$, then $h_D(Q_N)\sim2\sqrt{E_N}$. More quantitatively,
\[
 2\sqrt{E_N}\,[1-O(\rho_N^{1/3})]\le h_D(Q_N)\le2\sqrt{E_N}.
\]
For fixed $1<p<\infty$, if $(1+\log k_N)\rho_N\to0$, then
\[
 \lambda_{k_N,p}(Q_N)^{1/p}/\sqrt{E_N}\longrightarrow2/p.
\]
\end{theorem}
\begin{proof}
Choose $0<\eta<1$ and put
$k_i=(1-\eta)\pi/(2R_{i,N})$. The vector field
$V_i=k_i\tan(k_i x_i)$ has
\[
 a=(1-\eta)^2E_N,\qquad
 M=(1-\eta)^2\max_i e_{i,N}\csc^2(\pi\eta/2).
\]
Lemma~\ref{SP-lem:divergence}, with $\eta$ a constant multiple of $\rho_N^{1/3}$, gives the lower bound; the product squared ground state gives the upper bound.

For the nonlinear upper bound, rescale the three fixed interval trial functions from the preceding proof to each interval. Their score variables are fixed-shape variables multiplied by $e_{i,N}$. Weighted averages of the unmodified score variables converge in the required moment whenever the largest weight $e_{i,N}/E_N$ tends to zero. To justify this with only the stated finite moment, truncate the fixed-shape variables; the bounded centered part has variance bounded by a constant times the largest weight, while Minkowski and the moment tail control the remainder. This proves $L^r$ convergence for every fixed required $r\ge1$; first-moment convergence suffices when $p<2$.

Use $\ell_N=\lceil\log_2 k_N\rceil$ coordinates for disjoint support choices. Their total relative energy is at most $\ell_N\rho_N=o(1)$. The genus test and \eqref{SP-eq:trial-energy} therefore give the matching upper bound. The matching lower bound follows from the $p$-Cheeger inequality.
\end{proof}

The parameter $\rho_N$ permits a changing family of interval scales: no single interval may contribute a nonvanishing fraction of the Dirichlet energy. The weighted product theorem uses positive smooth weights and a fixed finite family of Lipschitz domains. All limits keep $p$ fixed and count the exterior boundary in the perimeter; relative-perimeter and Neumann problems have different normalizations.



\subsection{Heat spectra on manifolds and the Dirichlet ground-state law}
\label{WC-sec:heat}
The interpolation theorem applies to eigenfunctions without requiring a polynomial algebra.

\begin{theorem}[Finite heat spectra realized by a measure-preserving map]
\label{WC-thm:heat}
Let $M$ be a compact connected smooth Riemannian manifold without boundary, of positive dimension, with normalized volume $\mu$. Let $-\Delta\phi_j=\lambda_j\phi_j$, where $\phi_0=1$ and the real $\phi_j$ form an orthonormal eigenbasis. For every $t>0$ and every finite set $J$ of nonconstant eigenfunctions, there is a $\mu$-preserving measurable map $T$ such that
\[
 \Lp_T\phi_j=e^{-t\lambda_j}\phi_j,
 \qquad
 \Lp_T^n\phi_j=e^{-nt\lambda_j}\phi_j\quad(j\in J,n\ge0).
\]
The pair law $(T(Y),Y)$ can be made arbitrarily $W_\infty$-close to the stationary heat-kernel pair at time $t$. For every $j\in J$ and every $\ell\ge0$,
\[
 \int\phi_j(y)\phi_\ell(T^ny)\mu(dy)=\ind_{\{j=\ell\}}e^{-nt\lambda_j}.
\]
\end{theorem}
\begin{proof}
The heat kernel is strictly positive on a connected compact manifold. Thus its conditional source laws are equivalent to $\mu$. The selected eigenfunctions together with the constant are linearly independent, so the conditional feature vector has full affine support. Apply Theorem~\ref{WC-thm:interpolation} with $\lambda=1$ and then Corollary~\ref{WC-cor:iterate}. Orthogonality proves the correlation identity. The heat-kernel and spectral facts used here are standard; see \cite{Davies1989}.
\end{proof}
For the unit sphere $S^{d-1}$, $d\ge2$, taking complete harmonic spaces through degree $k$ retains all the factors $e^{-t\ell(\ell+d-2)}$ for $\ell\le k$. For the circle, this reads
\[
 \Lp_T\cos(j\theta)=e^{-tj^2}\cos(j\theta),\qquad
 \Lp_T\sin(j\theta)=e^{-tj^2}\sin(j\theta),\quad1\le j\le k.
\]
These are exact prescribed finite spectra of a Perron--Frobenius operator, together with the invariant volume law.

\begin{corollary}[The ground-state reference law]
\label{WC-cor:ground}
Let $\Omega\subset\R^d$ be a bounded connected smooth domain. Let $-\Delta u_j=\lambda_j u_j$ be a real orthonormal Dirichlet eigenbasis, with $u_1>0$. Put $\mu(dx)=u_1(x)^2dx$. For every finite $k$ and $t>0$, there is a $\mu$-preserving measurable map $T$ satisfying
\begin{equation}\label{WC-eq:ground-spectrum}
 \Lp_T(u_j/u_1)=e^{-t(\lambda_j-\lambda_1)}u_j/u_1,
 \qquad 1\le j\le k.
\end{equation}
All iterates and the corresponding cross-correlations are exact.
\end{corollary}
\begin{proof}
Use the Dirichlet ground-state transform
\[
 Q_tf(x)=\frac{e^{t\lambda_1}}{u_1(x)}e^{t\Delta_D}(u_1f)(x).
\]
It is a Markov semigroup with invariant probability $u_1^2dx$. Its kernel is strictly positive, and $Q_t(u_j/u_1)=e^{-t(\lambda_j-\lambda_1)}u_j/u_1$. The ratios form an orthonormal family in $L^2(\mu)$, so in particular they are integrable and affinely independent after excluding $j=1$. Apply Theorem~\ref{WC-thm:interpolation}.
\end{proof}
The density $u_1^2$ is the reference selected by the directional-energy problem and used in Theorem~\ref{R9-thm:mixed-weighted}. Here it supports deterministic realizations of finite portions of the ground-state-transformed heat spectrum. The original Dirichlet eigenvalues remain unchanged; the transfer operator carries their differences from $\lambda_1$.


\section{Exact comparisons for compact group orbits}
\label{R11-sec:orbits}
State-dependent rotations of a reference generate a convex set with an expected orbit support function. An independent Haar rotation converts membership into a martingale coupling. The stabilizer and uniqueness assumptions give the sharp constant and endpoint rigidity. Eigenvalue and singular-value sums give matrix criteria; assignment and parity expose the mechanism. Classical random-assignment, polar and spin-glass inputs are credited separately.

\subsection{Sharp references on compact orbits}\label{sec:sharp-references}

For a prescribed orbit, symmetry makes every orbit support function
constant on the target. The same support functions describe all means obtainable
by rotating the reference before taking its expectation. Haar averaging turns
membership in that convex set into one martingale coupling.

\begin{theorem}[Comparison with an arbitrary compact orbit]\label{R15-thm:orbit}
Let a compact group $K$ act continuously and orthogonally on a finite-dimensional
Euclidean space $H$, and let $Z$ be integrable and $K$-invariant. Define
\[
 \mathcal A_Z=\{\E[RZ]:\ R\in K\text{ may have any joint law with }Z\}.
\]
This is a compact convex $K$-invariant set, with
\begin{equation}\label{R15-eq:orbit-support}
 h_{\mathcal A_Z}(u)=\E\max_{g\in K}\ip{u}{gZ}.
\end{equation}
For Haar $U$, $a\in H$ and $c>0$,
\begin{equation}\label{R15-eq:orbit-criterion}
 Ua\cx cZ
 \quad\Longleftrightarrow\quad a\in c\mathcal A_Z
 \quad\Longleftrightarrow\quad
 h_{Ka}(u)\le c\E h_{KZ}(u)\quad(u\in H).
\end{equation}
These conditions are also equivalent to domination of some law supported on
$Ka$. They require neither a restriction on stabilizers nor uniqueness of a
support maximizer. The optimal scale is the gauge of $a$ in $\mathcal A_Z$;
if finite and positive it is attained.
\end{theorem}
\begin{proof}
Joint laws of $(Z,R)$ with fixed $Z$-marginal are compact: $K$ is compact,
and $|RZ|=|Z|$ makes the expectation continuous by uniform truncation.
Mixing gives convexity and rotating $R$ gives invariance. Maximizing over
successively finer finite subsets of $K$ proves
\eqref{R15-eq:orbit-support} by dominated convergence. Separation now
identifies membership with the displayed support inequalities.

The convex function $z\mapsto\max_g\ip{u}{gz}$ is constant on $Ka$, so
its expectation proves necessity for every target law on that orbit.
Conversely choose $(Z,R)$ with $c\E[RZ]=a$ and an independent Haar $U$.
Put $Y=Ua$ and $Z'=URZ$. Haar averaging is unchanged by a rotation within
each orbit; hence $Z'\sim Z$. Moreover
\[
 \E[cZ'\mid U]=U c\E[RZ]=Ua=Y.
\]
Conditioning on $Y$ proves the comparison. Compactness gives attainment.
\end{proof}

The general criterion is a convex membership test. The next theorem gives a
single scalar and a unique endpoint coupling when the stabilizer and support
maximizer have the stated rigidity.

\begin{theorem}[Sharp reference for a compact homogeneous orbit]
\label{FC-thm:compact-orbit}\label{FORD-thm:orbit}
Let a compact group $K$ act continuously and orthogonally on a finite-dimensional real Euclidean space $H$. Let $\mathcal V=Kv_0$, where $\norm{v_0}^2=r^2>0$. Suppose that the fixed space of the stabilizer $K_{v_0}$ is exactly $\R v_0$. Let $Z$ be an integrable $K$-invariant random vector, and assume that
\[
 Q(Z)=\argmax_{v\in\mathcal V}\ip{Z}{v}
\]
is unique almost surely. If $w=\E\max_{v\in\mathcal V}\ip{Z}{v}>0$, then $Q(Z)$ has the invariant probability law on $\mathcal V$ and
\[
 \E[Z\mid Q(Z)]=\frac{w}{r^2}Q(Z),\qquad
 Q(Z)\cx\frac{r^2}{w}Z.
\]
The scalar $r^2/w$ is minimal among all probability laws supported on $\mathcal V$. At that scalar, every martingale coupling is induced by $Q$, and the target law is the invariant law.
\end{theorem}
\begin{proof}
Uniqueness makes $Q$ equivariant, so its law is the unique invariant probability on the compact homogeneous space. An equivariant version of the conditional-mean function is obtained by averaging a regular conditional mean over Haar measure on $K$. Its value at $v_0$ is fixed by $K_{v_0}$, so it is $\alpha v_0$. Equivariance makes the same $\alpha$ valid on the whole orbit, up to the harmless choice of versions. Taking an inner product with $Q(Z)$ gives $w=\alpha r^2$. Conditional Jensen gives the comparison.

For sharpness put $h(z)=\max_{v\in\mathcal V}\ip{z}{v}$. Equal norms imply $h(v)=r^2$ on $\mathcal V$. Thus $Y\cx cZ$ for any target $Y$ supported there forces $r^2\le cw$. At equality, in any martingale coupling,
\[
 \E\ip{cZ}{Y}=\E\norm{Y}^2=r^2=c\E h(Z),
 \qquad \ip{Z}{Y}\le h(Z).
\]
The nonnegative gap has mean zero. Therefore $Y=Q(Z)$ almost surely, proving the endpoint assertion.
\end{proof}
For a source with ties, Haar symmetrization of a measurable maximizing selection gives an equivariant randomized maximizer and the same barycentre. Endpoint uniqueness then need not hold.

\subsubsection{An optimal Gaussian reference for permutation matrices}
Write $J_n=\1\1^{\mathsf T}$ and $P_n=I-J_n/n$. Let $G$ have independent standard Gaussian entries and define
\[
 M_n=\max_{\pi\in S_n}\sum_{i=1}^nG_{i,\pi(i)},\qquad \mu_n=\E M_n.
\]
Let $\Pi$ denote a uniform permutation matrix.
\begin{theorem}\label{FORD-thm:assignment}
For $n\ge2$,
\[
 \boxed{\quad \Pi-J_n/n\cx s_nP_nGP_n,
 \qquad s_n=\frac{n-1}{\mu_n}.\quad}
\]
This is the smallest possible scalar for a law supported on the centred permutation matrices. At this scale the target law is necessarily uniform. A joint coupling is obtained by selecting the maximum-weight assignment of $G$:
\[
 \E[P_nGP_n\mid \Pi]=\frac{\mu_n}{n-1}(\Pi-J_n/n).
\]
Moreover $s_n\sim(2\log n)^{-1/2}$.
\end{theorem}
\begin{proof}
Work on $H=\{X:X\1=X^{\mathsf T}\1=0\}$ with its Frobenius inner product. Independent row and column permutations act transitively on the centred permutation matrices. At $I-J_n/n$, the stabilizer includes simultaneous conjugation by all permutation matrices. A matrix fixed by all these conjugations has a common diagonal entry and a common off-diagonal entry. Intersection with $H$ leaves precisely the line through $I-J_n/n$.

Every target has squared norm $n-1$. Its support optimizer for $P_nGP_n$ is the same assignment as for $G$: the centring terms add the same scalar to every assignment. Distinct assignment costs have a nondegenerate Gaussian difference, so ties have probability zero. The expected support function is $\mu_n$. Theorem~\ref{FORD-thm:orbit} proves the result.

The asymptotic $\mu_n\sim n\sqrt{2\log n}$ is the Gaussian assignment theorem of Mordant--Segers \cite{FORD-MS}. Their normalization divides assignment costs by $\sqrt n$. Substitution gives the claimed scale.
\end{proof}
A maximum-weight assignment gives a polynomial-arithmetic decoder
from ideal Gaussian inputs. The normalization is the exact expectation
$\mu_n$. The theorem concerns uniform permutation matrices.

\subsubsection{The Parisi formula gives an exact non-Gaussian reference}
Let $E=(E_{ij})$ have independent mean-one exponential entries and let $\Pi$ minimize its assignment cost. The output is again uniform.
\begin{corollary}\label{FORD-cor:parisi}
Set $H_n^{(2)}=\sum_{j=1}^n j^{-2}$ and
\[
 \alpha_n=\frac{n-H_n^{(2)}}{n-1}>0\qquad(n\ge2).
\]
Then
\[
 \E[J_n-E\mid\Pi]=\alpha_n(\Pi-J_n/n),\qquad
 \Pi-J_n/n\cx\alpha_n^{-1}(J_n-E).
\]
One may project the source on both sides by $P_n$ to obtain the stronger supported reference $-\alpha_n^{-1}P_nEP_n$. For this projected source the displayed scalar is optimal.
\end{corollary}
\begin{proof}
The Linusson--W\"astlund theorem \cite{FORD-LW} gives expected minimum assignment cost $H_n^{(2)}$. Conditional on the identity being the minimizer, symmetry makes the expected selected entry $a=H_n^{(2)}/n$ and every off-diagonal entry equal to some $b$. The total matrix sum is invariant under row and column permutations. Transitivity of the selected permutation therefore gives conditional expected total $n^2$, so $a+(n-1)b=n$.
Consequently
\[
 \E[E\mid\Pi]=J_n-\alpha_n(\Pi-J_n/n).
\]
Projection and conditional Jensen prove the claims. For sharpness apply the support-function bound of Theorem~\ref{FORD-thm:orbit} to $-P_nEP_n$, whose expected maximum assignment value is $n-H_n^{(2)}$.
\end{proof}
\subsubsection{Linear codes and parity constraints}
A binary linear code will be written multiplicatively as a subgroup $\mathcal C\subseteq\{-1,1\}^n$. Assume no coordinate is fixed on $\mathcal C$. For independent standard Gaussians let
\[
 Q(G)=\argmax_{c\in\mathcal C}\langle G,c\rangle.
\]
\begin{theorem}[Diagonal reference from any linear code]\label{FORD-thm:code}
The output $Q(G)$ is uniform on $\mathcal C$. Put
$\alpha_i=\E[G_iQ_i(G)]$. Every $\alpha_i$ is positive, and
\[
 \E[G\mid Q=c]=\diag(\alpha_i)c,
 \qquad U_{\mathcal C}\cx\diag(\alpha_i^{-1})G.
\]
If coordinate automorphisms of the code act transitively, all $\alpha_i$ equal
\[
 \alpha=\frac1n\E\max_{c\in\mathcal C}\langle G,c\rangle.
\]
In this case $1/\alpha$ is the optimal scalar Gaussian reference, even among all target laws supported on the code. Its endpoint target law is uniquely uniform.
\end{theorem}
\begin{proof}
Multiplication by a codeword preserves the independent Gaussian source and permutes all decoder cells transitively. If $\alpha=\E[G\mid Q=\1]$ is now viewed as a vector, equivariance gives $\E[G\mid Q=c]=\diag(c)\alpha$. Averaging $G_iQ_i$ identifies its $i$th entry with $\alpha_i$.

To see positivity, fix $G_{-i}$. Because both coordinate signs occur in the code, maximizing over the two classes shows that $Q_i(G)$ is the sign of $G_i-t_i(G_{-i})$ for a finite threshold. Gaussian integration gives $\alpha_i=2\E\phi(t_i(G_{-i}))>0$. Coordinate transitivity gives equality of the $\alpha_i$. The support-function proof of Theorem~\ref{FORD-thm:orbit} gives scalar sharpness and endpoint rigidity directly.
\end{proof}
The decoder is weighted syndrome decoding. Write $G_i=s_i|G_i|$. A correction set $F$ has loss $2\sum_{i\in F}|G_i|$ from the unrestricted optimum, and must bring $s$ into the code. For a graph cycle code, the syndrome is the set of vertices with odd incidence, so this is a minimum-weight $T$-join. The standard reduction to shortest paths and weighted perfect matching gives a polynomial-arithmetic decoder for that class. The general code theorem is an existence statement.

\subsubsection{The exact cost of a parity constraint}
Let $\mathcal C_b=\{y\in\{-1,1\}^b:\prod y_i=1\}$, $b\ge2$. Its Gaussian decoder takes coordinatewise signs and, if parity is wrong, flips the coordinate with smallest absolute Gaussian value. This takes $O(b)$ operations.
\begin{theorem}\label{FORD-thm:parityGaussian}
Let $a_0=\sqrt{2/\pi}$ and
\[
 m_b=\E\min_{1\le i\le b}|G_i|
     =\int_0^\infty \operatorname{erfc}(t/\sqrt2)^b\,dt.
\]
The exact smallest scalar Gaussian reference for a law supported on $\mathcal C_b$ is
\[
 s_b=\left(a_0-\frac{m_b}{b}\right)^{-1}.
\]
The uniform code law attains it, and
\[
 \delta_b:=\frac{s_b}{\sqrt{\pi/2}}-1
       \sim\frac{\pi}{2b^2}.
\]
For $n$ a multiple of $b$, the product of these parity laws is $(b-1)$-wise independent and satisfies
\[
 Y\cx(1+\delta_b)\sqrt{\pi/2}\,G_n,
 \quad D(\mathcal L(Y)\Vert U_n)=\frac{n\log2}{b},
 \quad \operatorname{TV}(\mathcal L(Y),U_n)=1-2^{-n/b}.
\]
In particular,
\[
 \frac{D(\mathcal L(Y)\Vert U_n)}{n\sqrt{\delta_b}}
       \longrightarrow (\log2)\sqrt{2/\pi}.
\]
\end{theorem}
\begin{proof}
Signs and magnitudes are independent. Wrong parity has probability $1/2$ and costs twice the minimum magnitude. Thus the expected support value is $ba_0-m_b$, and Theorem~\ref{FORD-thm:code} gives the exact scale.

The distribution function of $|G|$ satisfies $F(t)=a_0t+O(t^3)$ at zero. Hence $b\min_i|G_i|$ converges to an exponential variable of rate $a_0$. Convergence of means follows by splitting the survival integral at a fixed small $\epsilon$: on $[0,\epsilon]$, $1-F(t)\le e^{-ct}$; on $[\epsilon,\infty)$, factor out $(1-F(\epsilon))^{b-2}$ and integrate the square of the Gaussian tail. Thus $bm_b\to1/a_0$, proving the expansion for $\delta_b$.

Take independent block couplings. Each block is uniform on $2^{b-1}$ points and every proper coordinate marginal is uniform. The information and total-variation formulas follow by counting.
\end{proof}
The displayed limit measures the entropy loss at the square-root scale of the excess Gaussian comparison. Near-critical full Gaussian comparison permits a law very far from the uniform cube in total variation while hiding its dependence from all small coordinate marginals.


\subsection{Compact orbits and complete spectral comparison criteria}
\label{FC-sec:orbits}

For eigenvalue and singular-value orbits, the support function in
Theorem~\ref{R15-thm:orbit} is a trace rearrangement. Its expectation depends
only on the expected ordered spectrum, including multiplicities. The
majorization criterion below therefore characterizes the entire convex-order
comparison at once.

\subsubsection{Real and complex spectral orbits}
Throughout this section $\mathbb F$ is $\R$ or $\C$,
$K_n=O(n)$ or $U(n)$ respectively, and $A^*$ denotes the transpose or
conjugate transpose. Complex matrices are regarded as real Euclidean
vectors with inner product $\operatorname{Re}\tr(A^*B)$.
The majorization arguments therefore apply to every real convex function
of their entries.

\subsubsection{The majorization statements used below}
For decreasing vectors $a,m\in\R^n$ with the same sum, write $a\prec m$ if
\[
 \sum_{i=1}^ja_i\le\sum_{i=1}^jm_i\quad(1\le j<n).
\]
For nonnegative decreasing vectors $s,m\in\R^k$, write $s\prec_wm$ if these partial-sum inequalities hold for $1\le j\le k$, without requiring equal sums.

\begin{lemma}[Permutation convex hulls]\label{FC-lem:majorization}
The relation $a\prec m$ is equivalent to $a\in\conv\{P_\pi m:\pi\in S_n\}$. The relation $s\prec_wm$ is equivalent to $s$ belonging to the convex hull of the signed permutations of $m$.
\end{lemma}
\begin{proof}
The partial sums are convex, so the conditions are necessary. For sufficiency,
rearrange a separating test $u$ decreasingly and sum by parts:
\[
 \sum_i u_i^{\downarrow}(a_i-m_i)
 =\sum_{j<n}(u_j^{\downarrow}-u_{j+1}^{\downarrow})
                  \sum_{i\le j}(a_i-m_i)\le0.
\]
The total-sum term is zero. For signed permutations use the ordered absolute
values of $u$; the final total-sum term is then nonpositive. Separation
proves both assertions.
\end{proof}

These are the usual majorization descriptions underlying Schur--Horn and Fan orbitopes \cite{FC-Orbitopes}. The following probabilistic statement applies them to the \emph{expected ordered spectrum} of a random source.

\begin{theorem}[Exact comparison for a fixed eigenvalue orbit]
\label{FC-thm:eigen-orbit}
Let $Z$ be an integrable Hermitian $n\times n$ matrix over $\mathbb F$ whose law is invariant under conjugation by $K_n$. Write its ordered eigenvalues as $\gamma_1\ge\cdots\ge\gamma_n$ and put $m_i=\E\gamma_i$. Let $a_1\ge\cdots\ge a_n$, and let
\[
 X_a=U\diag(a)U^*,\qquad U\text{ Haar on }K_n.
\]
For $c\ge0$,
\begin{equation}\label{FC-eq:eigen-criterion}
 X_a\cx cZ\quad\Longleftrightarrow\quad a\prec cm.
\end{equation}
The relation includes equality of total sums, $\sum_i a_i=c\E\tr Z$. The same conditions are necessary whenever any probability law supported on the orbit of $\diag(a)$ is dominated by $cZ$.
\end{theorem}
\begin{proof}
For a Hermitian test matrix $H$, trace rearrangement gives
\[
 \E\max_U\operatorname{Re}\tr(HUZU^*)
 =\sum_i\lambda_i(H)m_i
 =h_{\operatorname{conv}\{U\operatorname{diag}(m)U^*\}}(H).
\]
Thus $\mathcal A_Z$ is the spectral orbitope on the right. Its membership
condition is $a\prec cm$: necessity follows from the convex Ky Fan sums
and the trace; sufficiency follows by conjugating the permutation
representation in Lemma~\ref{FC-lem:majorization}. Apply
Theorem~\ref{R15-thm:orbit}, treating $c=0$ directly.
\end{proof}

In the centered-trace case $\sum_i a_i=\E\tr Z=0$, if all partial sums $M_j=\sum_{i\le j}m_i$ are positive, the smallest scalar is explicitly
\begin{equation}\label{FC-eq:eigen-scale}
 c_* =\max_{1\le j<n}\frac{\sum_{i=1}^ja_i}{M_j}.
\end{equation}
If a denominator is zero, the corresponding partial-sum inequality in \eqref{FC-eq:eigen-criterion} is the unambiguous formulation. A one-dimensional stabilizer fixed space is unnecessary for Theorem~\ref{FC-thm:eigen-orbit}; the finite majorization constraints describe the entire feasible orbit.

\begin{corollary}[Haar projectors]\label{FC-cor:projector}
Let $P_r$ be a uniformly distributed rank-$r$ self-adjoint projector in $\mathbb F^n$, where $1\le r<n$. Assume in addition to Theorem~\ref{FC-thm:eigen-orbit} that $\E\tr Z=0$, and suppose $M_r=\E\sum_{i=1}^r\gamma_i>0$. Then the exact minimal scalar, even over arbitrary target laws on the centered projector orbit, is
\[
 P_r-\frac rnI_n\cx c_{n,r}Z,\qquad
 c_{n,r}=\frac{r(1-r/n)}{M_r}.
\]
If the $r$th and $(r+1)$st eigenvalues of $Z$ are distinct almost surely, the endpoint coupling selects the top-$r$ eigenspace, and the endpoint target law is Haar.
\end{corollary}
\begin{proof}
The ordered target eigenvalues are $1-r/n$ repeated $r$ times and $-r/n$ repeated $n-r$ times. Since $m_i$ decrease, $M_j/j\ge M_r/r$ for $j\le r$. For $j\ge r$, the remaining tail averages give $M_j/(n-j)\ge M_r/(n-r)$. Hence the maximal ratio in \eqref{FC-eq:eigen-scale} occurs at $j=r$.

For endpoint rigidity, use the support function
\[
 h(Z)=\max_{P:\,\rank P=r}\tr\bigl((P-rI_n/n)Z\bigr),
 \qquad \E h(Z)=M_r.
\]
Every centered target projector has squared Frobenius norm $r(1-r/n)$. Equality in the support-function bound forces every optimal martingale coupling to choose the top-$r$ projector. The spectral gap makes that projector unique, and conjugation invariance makes its law Haar. The majorization criterion establishes existence of this coupling.
\end{proof}

\paragraph{GOE normalization and the proportional-rank limit.}
Let $W_n$ have independent upper-triangular entries, with off-diagonal law $N(0,1)$ and diagonal law $N(0,2)$, and set $Z_n=W_n-(\tr W_n/n)I_n$. If $r/n\to\rho\in(0,1)$, let $q_\rho\in(-2,2)$ satisfy
\[
 \int_{q_\rho}^2\frac{\sqrt{4-x^2}}{2\pi}\,dx=\rho,
 \qquad J_\rho=\frac{(4-q_\rho^2)^{3/2}}{6\pi}.
\]
The semicircle law \cite{FC-Wigner} gives
\begin{equation}\label{FC-eq:GOE-constant}
 \sqrt n\,c_{n,r}\longrightarrow
 \frac{\rho(1-\rho)}{J_\rho}.
\end{equation}
Indeed the expected sum of the top $r$ eigenvalues is asymptotic to $n^{3/2}J_\rho$. The centering subtracts a trace term of mean zero. Convergence of expectations follows from convergence of the empirical spectral laws with first moments: the spectral first moments are uniformly integrable by the bounded expected normalized second moment, and the corresponding top-fraction integrals converge. At $\rho=1/2$, $q_\rho=0$ and the constant is $3\pi/16$. The same scalar is optimal for the uncentered source $W_n$, since its expected trace is zero and the expected spectral partial sums agree.

\subsubsection{Fixed singular values and the Stiefel manifold}
\begin{theorem}[Exact comparison for a fixed singular-value orbit]
\label{FC-thm:singular-orbit}
Let $Z\in\mathbb F^{n\times k}$, $n\ge k$, be integrable and invariant under $Z\mapsto UZV^*$ for $U\in K_n$ and $V\in K_k$. Let $m_i=\E\sigma_i(Z)$, in decreasing singular-value order. For decreasing $s_1,\ldots,s_k\ge0$, let $X_s$ have the invariant law on matrices with those singular values. Then
\[
 X_s\cx cZ\quad\Longleftrightarrow\quad s\prec_w cm.
\]
The inequalities remain necessary for any target law supported on that singular-value orbit.
\end{theorem}
\begin{proof}
Singular-value trace rearrangement gives
$h_{\mathcal A_Z}(H)=\sum_i\sigma_i(H)m_i$.
The corresponding convex hull has membership condition $s\prec_w cm$
by the signed-permutation statement of Lemma~\ref{FC-lem:majorization}.
Theorem~\ref{R15-thm:orbit} proves necessity for any orbit-supported law
and sufficiency for the invariant law, including both fields.
\end{proof}

\begin{corollary}[Sharp Gaussian reference for Haar frames]
\label{FC-cor:Stiefel}
For $1\le k\le n$, let $Q$ be Haar on $\St_{\mathbb F}(n,k)=\{Q\in\mathbb F^{n\times k}:Q^*Q=I_k\}$. Let $G$ have independent standard real Gaussian entries if $\mathbb F=\R$, or independent circular complex Gaussian entries with $\E|G_{ij}|^2=1$ if $\mathbb F=\C$. Then
\begin{equation}\label{FC-eq:Stiefel}
 Q\cx\frac{k}{\E\norm{G}_*}G,
\end{equation}
where $\norm{G}_*$ is the nuclear norm. This scalar is minimal among all target laws on $\St_{\mathbb F}(n,k)$. The endpoint law is uniquely Haar, and a martingale coupling uses the polar factor
\[
 Q(G)=G(G^*G)^{-1/2},\qquad
 \E[G\mid Q(G)]=\frac{\E\norm{G}_*}{k}Q(G).
\]
\end{corollary}
\begin{proof}
For $s=(1,\ldots,1)$, the mean $j^{-1}\sum_{i\le j}m_i$ decreases with $j$. Thus the largest necessary ratio is $k/\sum_{i=1}^km_i$. Theorem~\ref{FC-thm:singular-orbit} proves \eqref{FC-eq:Stiefel}. For the explicit coupling, use the polar decomposition $G=QR$. Gaussian invariance makes $Q$ Haar and independent of $R=(G^*G)^{1/2}$; right invariance gives $\E R=(\E\tr R/k)I_k$. Equivalently apply Theorem~\ref{FC-thm:compact-orbit}. The support maximizer is the polar factor almost surely, since $G$ has full column rank.
\end{proof}

Tropp's convex comparison theorem \cite{FC-Tropp} gives a QR-based mean-chi constant, with a bound $3/2$ for normalized square Gaussian matrices. Formula~\eqref{FC-eq:Stiefel} identifies the optimal scalar through the polar factor. For $k=n$, let $\widetilde G_n=G_n/\sqrt n$ and
\[
 \alpha_{\R}(n)=\frac1n\E\norm{\widetilde G_n}_*.
\]
Then
\begin{equation}\label{FC-eq:Grothendieck}
 Q_n\cx\alpha_{\R}(n)^{-1}\widetilde G_n,
 \qquad \alpha_{\R}(n)^{-1}\longrightarrow\frac{3\pi}{8}.
\end{equation}
For either field, with $\alpha_{\mathbb F}(n)=n^{-1}\E\|G_n/\sqrt n\|_*$, the constant $\alpha_{\mathbb F}(n)^2$ is precisely the approximation ratio, with matching semidefinite integrality gap, in the little Grothendieck problem over $O(n)$ or $U(n)$ studied by Bandeira, Kennedy and Singer \cite{FC-BKS}. Thus the optimal all-convex comparison factor is the inverse square root of that ratio. The limiting value follows from the quarter-circle singular-value law, whose mean is $8/(3\pi)$.

The Gaussian polar identity is classical
\cite[Proposition~9]{FC-Tropp}; see also
\cite[Lemmas~6 and~17]{FC-BKS}. The expected-spectrum criteria above
specify exactly which fixed orbits it can support.


\subsection{Cut codes and the spin-glass Parisi constant}
\label{FC-sec:spin}
Let $K_n$ be the complete graph, with $M=\binom n2$ edges. Its cut-sign code is
\[
 \mathcal C_n=\{(x_ix_j)_{i<j}:x\in\{-1,1\}^n\}\subset\{-1,1\}^M.
\]
For $n\ge2$, this is a binary linear code in multiplicative notation. It has $2^{n-1}$ elements, has no fixed coordinate, and is transitive on coordinates under vertex permutations. A uniform element is obtained by taking independent fair vertex signs and recording all their pairwise products.

\begin{theorem}[Optimal Gaussian comparison for cut signs]
\label{FC-thm:cut}
Let $Y=(\xi_i\xi_j)_{i<j}$ with independent fair $\xi_i$, and let $G=(G_{ij})_{i<j}$ be standard Gaussian in $\R^M$. Put
\[
 w_n=\E\max_{x\in\{-1,1\}^n}\sum_{i<j}G_{ij}x_ix_j,
 \qquad s_n=\frac{M}{w_n}.
\]
Then
\begin{equation}\label{FC-eq:cut-reference}
 Y\cx s_nG.
\end{equation}
The scalar $s_n$ is minimal among all probability laws supported on $\mathcal C_n$. At the minimum, the target law is uniquely uniform, and the coupling selects the Gaussian ground-state cut. In particular,
\[
 \E\left[G\,\middle|\,
       \argmax_{c\in\mathcal C_n}\ip{G}{c}\right]
 =\frac{w_n}{M}\argmax_{c\in\mathcal C_n}\ip{G}{c}.
\]
\end{theorem}
\begin{proof}
Apply Theorem~\ref{FORD-thm:code} to the coordinate-transitive code $\mathcal C_n$. For completeness, codeword multiplication preserves the Gaussian law and permutes the decoder cells transitively. The decoder is uniform, and its conditional mean has the form $\diag(\alpha_e)c$. Edge transitivity makes $\alpha_e=w_n/M$ for every edge. Conditional Jensen gives \eqref{FC-eq:cut-reference}.

The support function $h(z)=\max_{c\in\mathcal C_n}\ip{z}{c}$ equals $M$ on every codeword. Thus any dominated code-supported target forces $M\le sw_n$. At equality, the same support-function equality forces the target to be the maximizing codeword. Distinct codeword scores have nondegenerate Gaussian differences, so the maximizing codeword is almost surely unique. The two vertex sign vectors $x$ and $-x$ represent the same codeword and cause no ambiguity.
\end{proof}

\begin{corollary}[Parisi normalization]\label{FC-cor:SK}
Define
\[
 P_* =\lim_{n\to\infty}\frac1n\E\max_{x\in\{-1,1\}^n}
       \frac1{\sqrt n}\sum_{i<j}G_{ij}x_ix_j.
\]
Then $P_*$ is the zero-field SK ground-state constant for covariance function $\xi(q)=q^2/2$, and
\begin{equation}\label{FC-eq:SK-scale}
 s_n\sim\frac{\sqrt n}{2P_*}.
\end{equation}
\end{corollary}
\begin{proof}
For $H_n(x)=n^{-1/2}\sum_{i<j}G_{ij}x_ix_j$, writing $R(x,y)=n^{-1}\sum_ix_iy_i$ gives
\[
 \E H_n(x)H_n(y)=\frac n2 R(x,y)^2-\frac12.
\]
Adding one independent $N(0,1/2)$ variable, common to every spin configuration, makes the covariance exactly $nR(x,y)^2/2$. This is the normalization in Auffinger--Chen \cite{FC-AC}. The common centered variable does not alter the expected maximum. Their ground-state formula identifies the limit above. Since $w_n\sim P_*n^{3/2}$ and $M\sim n^2/2$, \eqref{FC-eq:SK-scale} follows.
\end{proof}

The established Parisi variational expression can be substituted directly. Let $\gamma:[0,1)\to[0,\infty)$ be nondecreasing, right-continuous, and integrable, and let $\Psi_\gamma$ solve
\[
 \partial_t\Psi_\gamma
 =-\tfrac12\bigl(\partial_{xx}\Psi_\gamma
                +\gamma(t)(\partial_x\Psi_\gamma)^2\bigr),
 \qquad \Psi_\gamma(1,x)=|x|.
\]
Then, with this normalization,
\[
 P_* =\inf_\gamma\left\{
 \Psi_\gamma(0,0)-\frac12\int_0^1t\gamma(t)\,dt\right\}.
\]
Thus the existing spin-glass ground-state variational problem determines an optimal all-convex Gaussian comparison for random cut matrices. This uses the Parisi theorem as an input. It is distinct from the exponential assignment identity also called a Parisi formula in Corollary~\ref{FORD-cor:parisi}.

\paragraph{Quadratic observables.}
For any deterministic linear map $T:\R^M\to\R^q$, the same law gives $TY\cx s_nTG$. In particular, for any collection of zero-diagonal symmetric matrices $B_1,\ldots,B_q$, one obtains a simultaneous comparison for the quadratic forms $\xi^{\mathsf T}B_j\xi$. The scalar in \eqref{FC-eq:SK-scale} controls the entire cut-code law. It need not be the best scalar after a particular projection, and it does not replace the finer coefficient-sensitive estimates for a single quadratic form.

\paragraph{Complete hypergraphs.}
For fixed $k\ge2$, the code
\[
 \mathcal C_{n,k}=\left\{\left(\prod_{i\in S}x_i\right)_{S\in\binom{[n]}k}
                    :x\in\{-1,1\}^n\right\}
\]
has the same coordinate transitivity and no fixed coordinate. The exact optimal Gaussian scalar is therefore
\begin{equation}\label{FC-eq:hypergraph-scale}
 \frac{\binom nk}
 {\E\max_x\sum_{S\in\binom{[n]}k}G_S\prod_{i\in S}x_i}.
\end{equation}
Writing $P_{*,k}$ for the limiting expected maximum per spin of
$n^{-(k-1)/2}\sum_SG_S\prod_{i\in S}x_i$, the scalar is asymptotic to
$n^{(k-1)/2}/(k!P_{*,k})$. The applicable pure-spin covariance is $q^k/k!$. To verify the normalization, start from the ordered-index model with coefficient $1/\sqrt{k!}$. Its distinct-index part has exactly this law. The repeated-index part has variance $O_k(1)$ at every configuration, so the expected supremum of its absolute value is $O_k(\sqrt n)$ by the Gaussian exponential bound and a union bound over $2^n$ configurations. Dividing the maximum by $n$ removes that difference. The general mixed-spin ground-state theorem in \cite{FC-AC} then identifies $P_{*,k}$.

Sampling a uniform cut codeword is elementary. Realizing the maximizing Gaussian coupling requires solving a ground-state optimization problem; the comparison theorem does not bound its computational cost.


\Needspace{12\baselineskip}
\part{Entropy, feasibility and dependence}
\label{R6-part:information}
Which information survives the construction? The exact entropy frontier quantifies the cost of conditional means and identifies the source. Fixing the common and coordinate noises determines feasible dependence; varying the common noise leads to a different limiting classification.

\section{Exact entropy rates and recovery of the source}
\label{R11-sec:entropy}
Shannon's rate--distortion theorem minimizes information subject to an expected fidelity constraint~\cite{R9Shannon}. Here the output is constrained at every codeword by exact vector conditional means, and the full independent product source must be retained. For a symmetric scalar source the Gibbs variational calculation gives a lower bound through $\E\log\cosh(tW)$. The coding, puncturing and exact calibration steps below attain its convex conjugate under all those conditional identities.

Attainment makes this transform operational: the optimal entropy is measurable from the family of calibrated problems. Exact values on scales accumulating at one positive point determine the analytic transform; conjugacy and Mellin--Fourier inversion then recover the whole normalized source, including its atoms, under only a first moment. The inverse is an exact-data uniqueness statement. Its endpoint asymptotics answer a different quantitative question: the small-value behavior of the source determines the power of the entropy deficit. The Gaussian square-root law is one instance. This minimum-entropy problem fixes the product reference; the next section fixes common and coordinate noises and asks for the maximum feasible entropy.

\subsection{An exact entropy frontier at every reference scale}
\label{R9-sec:frontier}
The Gaussian entropy asymptotic in
Theorem~\ref{VB-thm:entropy-endpoint} requires both a universal
inequality and laws attaining it. The coefficient reference in the
signing theorem already gives the inequality under a hard matrix
constraint. This section determines the best
lower bound imposed by a product reference alone. Exact attainment will
require a separate code construction, while the variational lower bound
applies to every admissible target, including a hard-constrained target.

An equivariant code maximizer aligns every conditional mean with its label, initially at unequal coordinate amplitudes. The entropy identity shows that only a vanishing fraction fall below the desired amplitude. After removing those coordinates, independent source symmetries set every remaining mean exactly, preserving the whole product law and the optimal rate.

For a symmetric integrable scalar source $W$ and $0<a<\E|W|$, the
rate is
\begin{equation}\label{R6-eq:binary-frontier}
 R_W(a)=\sup_{t\ge0}\{ta-\E\log\cosh(tW)\}.
\end{equation}
The minimum normalized entropy among laws satisfying
$\E[W_i\mid Y]=aY_i$ tends to $R_W(a)$, with the full independent
source retained. We give the simplex theorem because the same proof
handles every prime-power alphabet. All entropies below use natural
logarithms.

Let $q\ge2$ and put
\[
 E_q=\{z\in\R^q:\textstyle\sum_a z_a=0\},\qquad
 v_a=e_a-q^{-1}\mathbf1,\qquad s_q^2=\|v_a\|^2=(q-1)/q.
\]
We identify a label word $y$ with the corresponding vector of simplex
vertices $(v_{y_i})_i$. The uniform law on all label words of length $n$
is denoted $U_{q,n}$. Let $Z\in E_q$ have a law $\nu$ invariant under all
permutations of its coordinates, and assume
\begin{equation}\label{R9-eq:source-assumptions}
 \E\|Z\|<\infty,\qquad b:=\E\max_a Z_a>0.
\end{equation}
When $q$ is a prime power, full permutation invariance can throughout
be weakened to invariance under the affine field permutations
$a\mapsto ua+v$, $u\in\mathbb F_q^*$, $v\in\mathbb F_q$.
This group is two-transitive, and all field symmetries used in the proof
belong to it. The weaker option is included in every statement below.
Write $M(Z)$ for the number of coordinates attaining the maximum.
Source blocks $Z_1,\ldots,Z_n$ will be independent with law $\nu$.
Define the scalar functions
\begin{align}
 L(t)&=\E\log\left(q^{-1}\sum_a e^{tZ_a}\right),\quad t\ge0,
       & A(t)&=L'(t),\label{R9-eq:partition}\\
 R(a)&=\sup_{t\ge0}\{ta-L(t)\},\quad 0\le a\le b.
       \label{R9-eq:rate}
\end{align}
For $0<a<b$, there is a unique $t_a>0$ with $A(t_a)=a$.
Moreover, $R$ is continuous and strictly convex on $[0,b]$, is strictly
increasing on $(0,b)$, and satisfies
$R(0)=0$ and $R(b)=\log q-\E\log M(Z)>0$.
These assertions are proved in Lemma~\ref{R9-lem:partition}.

\begin{theorem}[Global product-reference entropy frontier]
\label{R9-thm:frontier}
Suppose $q$ is a prime power and \eqref{R9-eq:source-assumptions} holds.
For every $0<a<b$,
\begin{equation}\label{R9-eq:global-frontier}
 \inf_{n\ge1}\ \inf_{\mu\cx (s_q^2/a)\nu^{\otimes n}}
       \frac{H(\mu)}n=R(a).
\end{equation}
Here and throughout, the inner infimum is over laws on the simplex
vertices. Its normalized value also converges to $R(a)$ as $n\to\infty$
through all integers. Equivalently, for every $\varepsilon>0$, with
$\kappa=s_q^2/b$,
\begin{equation}\label{R9-eq:entropy-envelope}
 \sup_{n\ge1}\ \sup_{\mu\cx(1+\varepsilon)\kappa\nu^{\otimes n}}
 \frac{D(\mu\|U_{q,n})}{n}
 =\log q-R\left(\frac b{1+\varepsilon}\right).
\end{equation}
There are finite-field linear codes $\mathcal C_j\subseteq\mathbb F_q^{n_j}$,
$n_j\to\infty$, and exact couplings
\begin{equation}\label{R9-eq:exact-attainment}
 Y_j\sim U_{\mathcal C_j},\qquad
 Z^{(j)}\sim\nu^{\otimes n_j},\qquad
 \E[Z_i^{(j)}\mid Y_j]=\frac{a}{s_q^2}v_{Y_{j,i}},\qquad
 \frac{\log|\mathcal C_j|}{n_j}\longrightarrow R(a).
\end{equation}
The finite-dimensional lower entropy bound holds for every integer $q$.
It is strict at every finite dimension when $0<a<b$. If the maximum
coordinate of $Z$ is almost surely unique, the only feasible target at
$a=b$ is $U_{q,n}$.
\end{theorem}

The scalar function $R$ is the usual single-letter rate-distortion
function \cite{R9Shannon} for the permutation-symmetric distortion
$\max_c z_c-z_a$, parameterized by retained score $a$.
Indeed, permutation symmetrization makes the reproduction label uniform;
the Gibbs variational principle then gives \eqref{R9-eq:rate}.
The theorem proves that the same rate is achievable with every conditional
vector mean exact and with the entire product source unchanged.
The argument below proves the coding statement directly, using only
pairwise independence and a finite first moment.

\begin{corollary}[Minimum support and lower-order R\'enyi entropies]
\label{R9-cor:support}
Under the hypotheses of Theorem~\ref{R9-thm:frontier}, for every
$\theta\in[0,1]$,
\[
 \inf_{n\ge1}\ \inf_{\mu\cx(s_q^2/a)\nu^{\otimes n}}
       \frac{H_\theta(\mu)}n=R(a),
\]
where $H_0(\mu)=\log|\operatorname{supp}\mu|$ and $H_1=H$.
In particular, the optimal exponential support size is $e^{nR(a)+o(n)}$.
\end{corollary}
\begin{proof}
For $0\le\theta\le1$, $H_\theta\ge H$. Uniform code laws have
$H_\theta=\log|\mathcal C_j|$ for all orders, so
\eqref{R9-eq:exact-attainment} gives equality in the dimension-normalized
infimum. Concatenating any fixed feasible block permits arbitrary
multiples of its length.
\end{proof}

\begin{corollary}[Every symmetric scalar source]
\label{R9-cor:binary}
Let $W$ be symmetric with $0<m=\E|W|<\infty$. Put
\[
 L_W(t)=\E\log\cosh(tW),\quad A_W(t)=L_W'(t),\qquad
 R_W(a)=\sup_{t\ge0}\{ta-L_W(t)\}.
\]
For every $\varepsilon>0$,
\begin{align}
 &\sup_{n\ge1}\ \sup_{\mu\cx(1+\varepsilon)W^{(n)}/m}
 \frac{D(\mu\|U_n)}n \notag\\
 &\qquad=\log2-R_W\left(\frac m{1+\varepsilon}\right)
 =\inf_{t\ge0}\left\{
 \frac{m\varepsilon}{1+\varepsilon}t+
 \E\log(1+e^{-2t|W|})\right\}.
 \label{R9-eq:binary-envelope}
\end{align}
Uniform binary code laws attain this envelope asymptotically, with exact
$\E[W_i\mid Y]=mY_i/(1+\varepsilon)$.
Atoms, bounded support and heavy tails with finite first moment are all
allowed. The right limit of the entropy envelope at $\varepsilon=0$ is
$\Prb(W=0)\log2$; endpoint uniqueness holds when this atom is absent.
\end{corollary}
\begin{proof}
Use the $q=2$ source $(W/2,-W/2)$ and replace $t$ by $2t$.
The last equality follows from
$\log\cosh(tW)=t|W|-\log2+\log(1+e^{-2t|W|})$.
\end{proof}

\begin{corollary}[The sharp Gaussian entropy loss at the endpoint]
\label{R7-cor:Gaussian-entropy}
Let $\kappa=\sqrt{\pi/2}$ and define
\[
 \Delta_G(\varepsilon)=
 \sup_{n\ge1}\sup_{\mu\cx(1+\varepsilon)\kappa G_n}
       \frac{D(\mu\|U_n)}n,
 \qquad \mu\text{ a law on }\{-1,1\}^n.
\]
Then
\begin{align}
 \Delta_G(\varepsilon)
 &=\inf_{t\ge0}\left\{
 \sqrt{\frac2\pi}\frac{\varepsilon}{1+\varepsilon}t+
      \E\log(1+e^{-2t|G_1|})\right\},
      \label{R7-eq:Gaussian-entropy-exact}\\
 \Delta_G(\varepsilon)&\le
 \sqrt{\frac\pi3}\sqrt{\frac{\varepsilon}{1+\varepsilon}},
 \qquad
 \Delta_G(\varepsilon)\sim\sqrt{\frac\pi3}\sqrt\varepsilon.
      \label{R7-eq:Gaussian-entropy-asymptotic}
\end{align}
At $\varepsilon=0$, every feasible law is $U_n$.
For positive $\varepsilon$, uniform binary code laws attain the
supremum asymptotically with all conditional means exact.
\end{corollary}
\begin{proof}
Corollary~\ref{R9-cor:binary} gives the exact formula and attainment.
The half-normal density is at most $a_0=\sqrt{2/\pi}$, so
\[
 \E\log(1+e^{-2t|G_1|})\le\frac{a_0\pi^2}{24t}.
\]
Minimization proves the upper bound. The density is continuous and
positive at zero, so substitution $u=2tx$ and dominated convergence
give equality asymptotically as $t\to\infty$ in this estimate.
Every minimizing sequence in the exact formula tends to infinity
as $\varepsilon\downarrow0$. Comparing with the elementary minimum
of $at+b/t$ gives the stated asymptotic. At the endpoint send
$t\to\infty$ in the entropy inequality to obtain $D(\mu\|U_n)=0$.
\end{proof}
The scale $\kappa$ is therefore an independence threshold as well as
the least scalar Gaussian scale. The hard discrepancy law of
Theorem~\ref{R3-thm:hard-Gaussian} approaches it at the optimal
radius order. For PSD partitions,
Corollary~\ref{R3-cor:Weaver} replaces the ambient $n$ in the entropy
loss by the matrix dimension $d$. These conclusions retain the hard
constraint for every output; the unconstrained coding laws in
\eqref{R7-eq:Gaussian-entropy-exact} establish the sharp envelope.

\begin{example}[A finite source recovers the complete Hamming curve]
\label{R9-ex:discrete}
Let $Z=\ell v_U$, where $U$ is uniform on $\mathbb F_q$ and $\ell>0$.
Then $b=\ell s_q^2$. Set
$d=s_q^2\varepsilon/(1+\varepsilon)$ and
$h(d)=-d\log d-(1-d)\log(1-d)$.
The right side of \eqref{R9-eq:entropy-envelope} equals
\[
 h(d)+d\log(q-1).
\]
This is the classical $q$-ary Hamming rate-distortion expression, now
attained within exact multivariate convex order. For $q=2$ and a
Rademacher scalar source, $d=\varepsilon/[2(1+\varepsilon)]$.
\end{example}

\begin{remark}[Gaussian specialization]
For standard Gaussian $Z$ on $E_q$, $b=\E\max_{a\le q}G_a$ for
independent standard real Gaussians. Theorem~\ref{R9-thm:frontier}
therefore determines the entire simplex Gaussian entropy curve.
The near-endpoint behaviour is obtained by expanding this exact variational formula. The code families used
here need not be Hamming codes or coordinate-transitive.
\end{remark}

\subsubsection{An exact entropy identity and its rigidity consequences}
\label{R9-sec:identity}
The entropy argument permits arbitrary couplings and only aggregate
coordinate scores. This weaker hypothesis is what makes average coding
performance useful for full-reference construction.

\begin{lemma}[Scalar partition geometry]\label{R9-lem:partition}
Under \eqref{R9-eq:source-assumptions}, $L$ is finite and continuously
differentiable on $[0,\infty)$, with $L'(0)=0$. On $(0,\infty)$,
$L''>0$, $A(t)\uparrow b$, and
\[
 L(t)=tb-\log q+\E\log M(Z)+o(1)\qquad(t\to\infty).
\]
The assertions following \eqref{R9-eq:rate} hold.
\end{lemma}
\begin{proof}
The sample log partition and its first derivative are bounded by a
constant times $\|Z\|$. For a fixed $t>0$, write
$d_a=\max_c Z_c-Z_a$. The variance of the softmax score is at most
\[
 \sum_a d_a^2e^{-td_a}\le \frac{4q}{e^2t^2}.
\]
Differentiation on compact subintervals of $(0,\infty)$ is justified by
this bound. The variance is positive whenever the coordinates are not
all equal; this has positive probability because $b>0$.
Dominated convergence gives $A(t)\uparrow b$ and
$\E\log\sum_a e^{-td_a}\to\E\log M(Z)$. Legendre duality gives strict convexity
and the interior assertions for $R$. The endpoint values follow by
monotone limits; the supporting slope tends to infinity at $b$.
\end{proof}

For $t\ge0$ and label $a$, define the probability measure
\begin{equation}\label{R9-eq:backward}
 \Gamma_{t,a}(dz)
 =\frac{q e^{tz_a}}{\sum_c e^{tz_c}}\,\nu(dz).
\end{equation}
Permutation invariance makes its integral one. Let $P$ be any joint law
of $Z^{(n)}\sim\nu^{\otimes n}$ and a label word $Y$, and set
$g_i=\E\langle Z_i,v_{Y_i}\rangle$.

\begin{theorem}[Full posterior entropy identity]
\label{R9-thm:identity}
For any $t_1,\ldots,t_n\ge0$,
\begin{align}
 H(Y)-\sum_i\{t_i g_i-L(t_i)\}
 &=H(Y\mid Z^{(n)}) \notag\\
 &\quad+\E_YD\left(P_{Z^{(n)}\mid Y}\,
       \middle\|\,\bigotimes_i\Gamma_{t_i,Y_i}\right).
 \label{R9-eq:identity}
\end{align}
In particular, for $g_i\in[0,b]$,
\begin{equation}\label{R9-eq:entropy-trace}
 \sum_iR(g_i)\le H(Y),\qquad
 nR\left(n^{-1}\sum_i g_i\right)\le H(Y).
\end{equation}
For an exact martingale reference
$\E[Z_i\mid Y]=(a/s_q^2)v_{Y_i}$, this gives $H(Y)\ge nR(a)$.
\end{theorem}
\begin{proof}
Define a reference joint law with source $\nu^{\otimes n}$ and the
product softmax channel
\[
 K(y\mid z)=\prod_i\frac{e^{t_i z_{i,y_i}}}{\sum_c e^{t_i z_{i,c}}}.
\]
Its target is $U_{q,n}$ and its posterior source is the product in
\eqref{R9-eq:identity}. Expanding $D(P\|\nu^{\otimes n}K)$ first by
conditioning on the source and then by conditioning on the target gives
\[
 n\log q-H(Y\mid Z)-\sum_i(t_i g_i-L(t_i))
 =n\log q-H(Y)+\E_YD(P_{Z\mid Y}\|\Gamma_{t,Y}).
\]
This is the identity. The finite label set and integrable log density
ensure that all terms are finite. Optimize each $t_i$, then use Jensen.
\end{proof}

The lower bound in Theorem~\ref{R9-thm:frontier} follows already. For
$0<a<b$, equality would force both terms on the right of
\eqref{R9-eq:identity} to vanish. The first makes the target deterministic
given the source, whereas the second gives, for every target atom of
positive probability, a posterior equivalent to the full source law.
These two properties allow only one target atom. A positive exact
simplex regression with centered source has uniform single-coordinate
label marginals and therefore at least $q$ target atoms. This proves
finite-dimensional strictness. When the source maximum is almost surely
unique, equality at $a=b$ in the pointwise maximum test forces each label
to be that source winner; these
winners are independent and uniform. They also have the required
endpoint conditional means by permutation symmetry.

\begin{proposition}[A quantitative removal bound]
\label{R9-prop:puncture-certificate}
Fix $a\in(0,b)$, let $t_a=R'(a)$, and put
\[
 J_a(x)=R(x)-R(a)-t_a(x-a).
\]
For any coupling with $0\le g_i\le b$ and any $c\in(0,a)$,
\begin{equation}\label{R9-eq:bad-coordinates}
 \#\{i:g_i\le c\}
 \le\frac{H(Y)-t_a\sum_i g_i+nL(t_a)}{J_a(c)}.
\end{equation}
The denominator is positive. If $H(Y_n)/n\to R(a)$ and
$n^{-1}\sum_i g_i\to a$, then the fraction of coordinates outside
$(a-\delta,a+\delta)$ tends to zero for every $\delta>0$.
\end{proposition}
\begin{proof}
The identity gives
$\sum_iJ_a(g_i)\le H(Y)-t_a\sum_i g_i+nL(t_a)$.
Strict convexity implies $J_a>0$ away from $a$ and
$J_a(x)\ge J_a(c)$ for $x\le c$. Compactness of $[0,b]$ proves the
last assertion.
\end{proof}

\begin{corollary}[Posterior rigidity at the entropy frontier]
\label{R9-cor:rigidity}
For an exact reference at score $a$, write
$\Delta_n=H(Y)-nR(a)$. Then
\[
 H(Y\mid Z^{(n)})\le\Delta_n,
 \qquad
 \E_YD\left(P_{Z^{(n)}\mid Y}\,
 \middle\|\,\bigotimes_i\Gamma_{t_a,Y_i}\right)\le\Delta_n.
\]
If $I$ is an independent uniform $k$-subset of coordinates, then
\begin{equation}\label{R9-eq:local-rigidity}
 \E_{I,Y}D\left(P_{Z_I\mid Y}\,
 \middle\|\,\bigotimes_{i\in I}\Gamma_{t_a,Y_i}\right)
 \le\frac{k\Delta_n}{n}.
\end{equation}
Thus every fixed random source block approaches the product backward
Gibbs law whenever $\Delta_n=o(n)$. Its mean total variation error is
at most $\sqrt{k\Delta_n/(2n)}$.
\end{corollary}
\begin{proof}
The first two statements are \eqref{R9-eq:identity}. For a fixed target,
expand relative entropy of a coordinate subvector by the chain rule.
Conditioning on all earlier coordinates can only increase each expected
conditional relative entropy, by convexity. Averaging over uniform
$k$-subsets counts each full-chain term $k/n$ times. This proves
\eqref{R9-eq:local-rigidity} without any differential-entropy assumption.
Pinsker's inequality and Jensen give the final assertion.
\end{proof}
If the target coordinates in $I$ are themselves independent and uniform,
relative-entropy contraction also gives
\[
 \E_I D\bigl(P_{Z_I,Y_I}\,\big\|\,K_{t_a}^{\otimes k}\bigr)
 \le k\Delta_n/n,
\]
where $K_{t_a}(dz,a)=q^{-1}\Gamma_{t_a,a}(dz)$ is the single-block
joint Gibbs law. Hence the asymptotically optimal code constructions
below have product Gibbs limits on random fixed coordinate sets,
although their global target support is exponentially smaller than the
full cube or simplex product.

\subsubsection{Random linear codes, puncturing and exact calibration}
\label{R9-sec:codes}
We prove attainment in Theorem~\ref{R9-thm:frontier}. The code is chosen for its average score; exact vector calibration will follow from puncturing and law-preserving source randomization. Throughout, $q$ is a prime power and source coordinates are indexed by $\mathbb F_q$. No transitivity assumption is placed on the coordinate
positions of a code.

\begin{lemma}[A vector regression from field symmetries]
\label{R9-lem:code-regression}
Let $\mathcal C\le\mathbb F_q^n$. From $Z^{(n)}\sim\nu^{\otimes n}$,
select $Y$ uniformly among the maximizers of
\[
 \sum_i\langle Z_i,v_{y_i}\rangle,\qquad y\in\mathcal C.
\]
Then $Y\sim U_{\mathcal C}$ and, with
$g_i=\E\langle Z_i,v_{Y_i}\rangle$,
\begin{equation}\label{R9-eq:code-regression}
 0\le g_i\le b,\qquad
 \E[Z_i\mid Y]=\frac{g_i}{s_q^2}v_{Y_i}.
\end{equation}
The assertions allow atoms and ties among codewords.
\end{lemma}
\begin{proof}
Translation of labels by any codeword, accompanied by the corresponding
coordinatewise permutation of each source block, preserves the joint
experiment and the uniform tie rule. The action is transitive on code
words, so the output is uniform. Conditional on $Y=0$, simultaneous
multiplication of all labels by any nonzero element of $\mathbb F_q$
also preserves the experiment. The nonzero labels in each source block
are permuted transitively, while zero is fixed. Its conditional mean is
therefore a multiple of $v_0$. Translating by a codeword proves the
conditional identity, with multiplier determined by $g_i$.

For nonnegativity, give coordinate $i$ a variable nonnegative weight
$x$. Conditional on all other source blocks, the expected support
maximum is a convex function of $x$ with a minimum at zero: Jensen's
inequality applies to the centered, independent $i$th block. Every
selected support maximizer is a subgradient. Its expectation at
$x=1$ is $g_i$, hence is nonnegative. The pointwise maximum test gives
$g_i\le b$.
\end{proof}

Write
\[
 h_q(\delta)=-\delta\log\delta-(1-\delta)\log(1-\delta)
                +\delta\log(q-1),\quad 0<\delta<1-1/q.
\]
For a linear code, $d^\perp$ denotes the minimum Hamming weight of a
nonzero word in its dual, with $d^\perp=\infty$ for a zero dual.

\begin{theorem}[Optimal average score with a finite first moment]
\label{R9-thm:random-code}
Fix $r\in(0,R(b))$ and write $a_r=R^{-1}(r)$.
There exist linear codes $\mathcal C_n\le\mathbb F_q^n$ with
$\log|\mathcal C_n|/n\to r$ such that their support decoders satisfy
\begin{equation}\label{R9-eq:average-optimal}
 \frac1n\E\max_{y\in\mathcal C_n}\sum_i\langle Z_i,v_{y_i}\rangle
 \longrightarrow a_r.
\end{equation}
For any fixed $\delta$ with $h_q(\delta)<r$, the codes can simultaneously
be chosen with $d^\perp>\delta n$. They can also simultaneously have
minimum nonzero codeword weight greater than $\delta_0n$ for any fixed
$\delta_0$ with $h_q(\delta_0)<\log q-r$.
\end{theorem}
\begin{proof}
Take an $n\times k_n$ matrix whose entries are independent and uniform
in $\mathbb F_q$, where $k_n\log q/n\to r$; its image is the code.
For a fixed nonzero message the resulting word is uniform on
$\mathbb F_q^n$. Two linearly independent messages give independent
uniform words. Each projective line contains $q-1$ nonzero messages.

Fix $t>0$ with $R(A(t))<r$. Given a source realization $z$, tilt the
uniform word by the product channel
\[
 P_t(y\mid z)=\prod_i\frac{e^{tz_{i,y_i}}}{\sum_c e^{tz_{i,c}}}.
\]
The empirical conditional mean score converges in probability to $A(t)$,
and the empirical log partition converges to $L(t)$, by the law of large
numbers and integrability. The conditional variance of the total score
is at most $4qn/(e^2t^2)$, by the bound in
Lemma~\ref{R9-lem:partition}. Thus, for every sufficiently small fixed
$\eta>0$, on source environments of probability tending to one,
\[
 P_t\left(\left|n^{-1}\sum_i z_{i,Y_i}-A(t)\right|\le\eta
                 \middle|z\right)\ge\tfrac12,
 \quad
 n^{-1}\sum_i\log\left(q^{-1}\sum_c e^{tz_{i,c}}\right)
          \ge L(t)-\eta.
\]
Changing measure on this event shows that a uniform word has score at
least $n(A(t)-\eta)$ with probability
\begin{equation}\label{R9-eq:tilted-count}
 p_n(z)\ge\tfrac12\exp\{-n[R(A(t))+(t+1)\eta]\}.
\end{equation}
Choose $\eta$ so that the exponent is smaller than $r$.

Count the nonzero messages whose words exceed this score. Conditional
on $z$, their indicators are independent except possibly for messages
on the same projective line. Each indicator has mean $p_n(z)$, and
\[
 \operatorname{Var}(N_{\rm good}\mid z)
 \le(q-1)(q^{k_n}-1)p_n(z).
\]
Consequently,
\[
 \Prb(N_{\rm good}=0\mid z)
 \le\frac{q-1}{(q^{k_n}-1)p_n(z)}\longrightarrow0
\]
on the typical environments. Let $\eta\downarrow0$ and then let $A(t)$
increase to $a_r$. This proves the lower bound in probability over code
and source. The absolute score divided by $n$ is bounded by
$n^{-1}\sum_i\max_c|Z_{i,c}|$, a uniformly integrable family. The same
lower bound therefore holds in expectation.

For every fixed code, its support decoder has entropy at most
$k_n\log q$. Equations~\eqref{R9-eq:entropy-trace} and
\eqref{R9-eq:code-regression} give the matching deterministic upper bound
$R^{-1}(k_n\log q/n)$ for its average score. Thus the expected deficit
from this bound tends to zero.

The random matrix has full column rank with probability tending to one:
a union bound over nonzero kernel messages gives failure probability at
most $q^{k_n-n}$. Also, a fixed nonzero dual word has probability
$q^{-k_n}$ of annihilating the generator. Since
\[
 \sum_{j\le\delta n}\binom nj(q-1)^j
 \le\exp\{n h_q(\delta)+o(n)\},
\]
$h_q(\delta)<r$ makes the probability of a dual word of weight at most
$\delta n$ tend to zero. Conditioning on these two high-probability
properties leaves the expected nonnegative score deficit tending to
zero. A further union bound over nonzero messages bounds the probability
of a codeword of weight at most $\delta_0n$ by
$\exp\{n[r+h_q(\delta_0)-\log q]+o(n)\}$, which also tends to zero
under the stated inequality. Select a deterministic code with all these
properties and vanishing deficit. This proves the theorem.
\end{proof}

\begin{theorem}[From average score to exact product calibration]
\label{R9-thm:puncturing}
Suppose codes have rate $r$ and average score tending to $a_r$ as in
Theorem~\ref{R9-thm:random-code}. For every fixed $a<a_r$, remove the
coordinates with $g_i<a$, and let $\mathcal C'_n$ be the projected code.
The removed fraction tends to zero. On the retained coordinates there
is an explicit coupling with source $\nu^{\otimes n'}$, target
$U_{\mathcal C'_n}$, and exact conditional means
\[
 \E[Z_i'\mid Y']=\frac{a}{s_q^2}v_{Y_i'}.
\]
Moreover, $\log|\mathcal C'_n|/n'\to r$. Puncturing does not decrease
the absolute minimum distance of the dual code.
\end{theorem}
\begin{proof}
Proposition~\ref{R9-prop:puncture-certificate}, centered at $a_r$,
shows that the fraction with $g_i<a$ vanishes. Projection maps a uniform
linear code to a uniform linear code; the tower property preserves
\eqref{R9-eq:code-regression} on its retained coordinates. For coordinate
$i$, independently of everything else, retain the source block with
probability $a/g_i$, and otherwise permute its labels by an independent
uniform field translation. The average of these $q$ translations is
zero on $E_q$. The conditional mean is therefore exactly as stated.
Each block law is invariant under every translation, so the complete
independent source remains unchanged.

Removing $n-n'$ coordinates changes the code dimension by at most
$n-n'$, which proves the rate assertion. Finally, $(\mathcal C'_n)^\perp =\{u_J:(u_J,0)\in\mathcal C_n^\perp\}$. A nonzero shortened dual word has its original weight, proving the last
claim.
\end{proof}

\begin{proof}[Completion of Theorem~\ref{R9-thm:frontier}]
For the desired $a$, take $r\downarrow R(a)$ through rates larger than
$R(a)$. At each rate apply Theorems~\ref{R9-thm:random-code} and
\ref{R9-thm:puncturing}. A diagonal sequence has exact score $a$ and
entropy rate $R(a)$. This gives \eqref{R9-eq:exact-attainment} and the
upper bound in \eqref{R9-eq:global-frontier}. The lower bound and the
endpoint assertions were proved after Theorem~\ref{R9-thm:identity}.
To obtain convergence through all block lengths, fix a calibrated block
whose rate is within $\eta$ of $R(a)$, concatenate copies, and fill the
bounded remainder with the uniform product softmax channel. Letting
$\eta\downarrow0$ proves the all-length assertion.
\end{proof}

\begin{corollary}[Optimal references with linearly many exact product marginals]
\label{R9-cor:privacy-code}
For every $\delta\in(0,1-1/q)$ satisfying $h_q(\delta)<R(a)$,
the codes in \eqref{R9-eq:exact-attainment} can be chosen with
$d^\perp>\delta n_j$. Their uniform laws are therefore
$\lfloor\delta n_j\rfloor$-wise independent. Simultaneously, for every
$\delta_0$ with $h_q(\delta_0)<\log q-R(a)$, their minimum codeword
distance can exceed $\delta_0n_j$ for all sufficiently large $j$.
\end{corollary}
\begin{proof}
Choose $\delta'>\delta$ with $h_q(\delta')<R(a)$ and impose the dual
distance bound in the random-code theorem throughout the diagonal
construction. Shortening the dual preserves its weight bound. A
projection of a linear code onto a coordinate set is the full coordinate
space exactly when there is no nonzero dual word supported on that set.
Uniform code laws then project to uniform independent labels. For the
last assertion, impose a slightly stronger codeword-distance bound in
Theorem~\ref{R9-thm:random-code}. Puncturing removes $o(n)$ coordinates
and reduces codeword weights by at most $o(n)$.
\end{proof}

\begin{remark}[The role of the two symmetries]
Field translations provide uniform target probabilities. Multiplication
by a nonzero field element provides the full vector regression in each
simplex block. Coordinate transitivity would make all $g_i$ equal, but
is unnecessary: strict convexity controls their exceptional coordinates,
and puncturing followed by source randomization makes the conditional means exactly equal to their targets.
This separates source coding from full-reference calibration.
\end{remark}

\subsubsection{Families at the optimal entropy exponent}
The strict slack in a stronger code reference can also be used to vary
its target probabilities. We give the scalar form, which applies to the
bounded calibrated auxiliary.

\begin{theorem}[Positive families of asymptotically optimal reference laws]
\label{R9-thm:families}
Let $W$ be symmetric, $|W|\le B$ almost surely, with
$v=\E W^2>0$ and $m=\E|W|$. Fix $0<a<a_+<m$ and $\eta>0$.
There are binary codes $\mathcal C\subseteq\{-1,1\}^n$ of arbitrarily
large length, with $N=|\mathcal C|$, such that
\[
 nR_W(a_+)\le\log N\le n[R_W(a_+)+\eta],\qquad d^\perp>2.
\]
For $p=U_\mathcal C$, put
\[
 \alpha=1-a/a_+,\qquad L=1+anB/v.
\]
Every function $f:\mathcal C\to[-1,1]$ with
$\E_p f=0$ and $\E_p[fY]=0$ gives a law
\begin{equation}\label{R9-eq:family}
 p_t(y)=p(y)[1+t f(y)],\qquad |t|\le\frac{\alpha}{2L},
\end{equation}
with the exact entire-source coupling
\[
 T\sim\law(W)^{\otimes n},\qquad
 \E[T\mid Y]=aY,\qquad Y\sim p_t.
\]
All laws in this relative neighborhood satisfy
\[
 nR_W(a)\le H(p_t)\le\log N.
\]
The neighborhood has dimension $N-n-1$, and every law has full support
$\mathcal C$. Taking $a_+\downarrow a$ and $\eta\downarrow0$ gives
families whose entropy and support exponents converge to the optimum
$R_W(a)$.
\end{theorem}
\begin{proof}
Choose a uniform code coupling with conditional mean $a_+Y$ using
Theorem~\ref{R9-thm:frontier} and Corollary~\ref{R9-cor:privacy-code}.
Write its conditional label probabilities as $k_y^+(z)$. Mixing with
an independent source--target draw gives
\[
 k_y(z)=(1-\alpha)k_y^+(z)+\alpha p(y),
 \qquad \int z k_y(z)\,\nu^{\otimes n}(dz)=a p(y)y.
\]
Set $w_y(z)=1+(a/v)\langle z,y\rangle$. Its absolute value is at most
$L$, and its source mass and first moment are $1$ and $ay$.
The modified kernel $k'_y(z)=k_y(z)+t p(y)f(y)w_y(z)$ has sum one pointwise, since $\E_p f=\E_p[fY]=0$. It is nonnegative by
$k_y\ge\alpha p(y)$ and $|t|L\le\alpha/2$. Integrating gives the exact
target mass $p_t(y)$ and moment $a p_t(y)y$. The source remains fixed.
Entropy is at most $\log N$ and at least $nR_W(a)$.

The functions $1,Y_1,\ldots,Y_n$ are linearly independent on the code:
dual distance greater than two gives zero coordinate means and identity
covariance. Hence the probability perturbation space has dimension
$N-n-1$. The strict relative bound keeps all target probabilities positive.
\end{proof}

The correction uses the positive independent component from the positive source-preserving construction. Here it applies at the sharp
entropy exponent, with no dependence of the relative radius on the
smallest atom $1/N$. It has a sampling implementation that avoids
conditional likelihoods. With probability $1-\alpha$ use the stronger
joint sampler. On the independent branch, propose $z\sim\nu^{\otimes n}$
and $y\sim U_\mathcal C$ and accept with probability
\[
 \frac{1+(t/\alpha)f(y)w_y(z)}{1+|t|L/\alpha}.
\]
The numerator integrates to one and lies between $1/2$ and $3/2$.
This branch uses at most $3/2$ proposals in expectation. Uniform code
sampling uses its message bits. Any real acceptance decisions require
exact coins or certified interval evaluations; a finite rational implementation uses these operations as exact primitives.


\subsection{The entropy frontier determines the reference law}\label{SP-sec:inverse}
Let $W$ be symmetric, with $0<m=\E|W|<\infty$, and let $W^{(n)}$ have independent coordinates of law $W$. Write $U_n$ for the uniform sign law. Define
\begin{equation}\label{SP-eq:frontier-def}
 \Delta_W(\varepsilon)=\sup_{n\ge1}\ 
 \sup_{\substack{\mu\text{ on }\{-1,1\}^n\\
 \mu\cx (1+\varepsilon)W^{(n)}/m}}
 \frac{\KL(\mu\Vert U_n)}n,\qquad \varepsilon>0.
\end{equation}
The normalization makes this curve invariant under positive rescaling of $W$. Corollary~\ref{R9-cor:binary} gives
\begin{equation}\label{SP-eq:frontier-input}
\begin{aligned}
 \Delta_W(\varepsilon)&=\log2-R_X\bigl((1+\varepsilon)^{-1}\bigr),
 \qquad X=W/m,\\
 R_X(a)&=\sup_{t\in\R}\{ta-F_X(t)\},\qquad
 F_X(t)=\E\log\cosh(tX).
\end{aligned}
\end{equation}
The maximum normalized entropy deficit is thus a convex-dual measurement of the scalar source. The exact coding theorem is needed to identify the operational frontier with the right side; a lower entropy bound alone would not provide this identity.

\begin{theorem}[Identification from accumulating entropy scales]\label{SP-thm:entropy-inverse}
Let $W_1,W_2$ be symmetric, integrable and nonzero. If their entropy frontiers agree on any set of positive scales having an accumulation point in $(0,\infty)$, then
\[
 \frac{W_1}{\E|W_1|}\ \stackrel d=\ 
 \frac{W_2}{\E|W_2|}.
\]
Conversely, equality of these normalized laws gives equality of the entire frontiers. The frontier is real analytic on $(0,\infty)$, even when $W$ has atoms or lacks a second moment.
\end{theorem}

The proof passes from the measured frontier to its convex dual and then to the source law. At every positive transform parameter, exponential damping makes all derivatives finite even when the source has only a first moment. Analytic continuation therefore recovers the full transform from the accumulating exact measurements. Mellin inversion reads the remaining source distribution through the Fourier transform of a weighted law of $\log|X|$.

\begin{proof}
Normalize $\E|X|=1$, and put $R=|X|$. For $t>0$,
\begin{equation}\label{SP-eq:analytic-F}
 F_X(t)=t-\log2+\E\log(1+e^{-2tR}),\quad
 F_X'(t)=\E[R\tanh(tR)],\quad
 F_X''(t)=\E[R^2\operatorname{sech}^2(tR)]>0.
\end{equation}
The second derivative is finite for every positive $t$. To justify analyticity without moment assumptions, extend $t$ to $\Re t>0$. The principal logarithm of $1+e^{-2tR}$ is well defined there: for $R>0$ its real part before taking the logarithm is positive. On a compact subset of that half-plane the logarithm is uniformly bounded as $R$ varies. For small $R$ this follows from convergence to $\log2$, and for large $R$ from exponential decay. Dominated holomorphic integration proves analyticity. Also $F'_X(t)\uparrow1$ as $t\to\infty$ and $F'_X(t)\downarrow0$ as $t\downarrow0$, by dominated convergence using $R\in L^1$.

The analytic inverse-function theorem makes the maximizing $t=t(a)$ and $R_X(a)$ analytic on $0<a<1$. Thus \eqref{SP-eq:frontier-input} makes $\Delta_W$ analytic on positive scales. The identity theorem for real analytic functions shows that agreement on a set with an interior accumulation point gives agreement throughout $(0,\infty)$.

The curve determines $R_X(a)$ for $0<a<1$. Symmetry gives its values for $-1<a<0$, and $R_X(0)=0$. If $p_0=\Pp(R=0)$, then $R_X(1)=R_X(-1)=(1-p_0)\log2$. Indeed, $t-F_X(t)=\log2-\E\log(1+e^{-2tR})$ increases to this value. The same endpoint is the limit of $R_X(a)$ as $a\uparrow1$: use any fixed $t$ for the lower bound and monotonicity for the upper bound. Outside $[-1,1]$ the conjugate is infinite, since $F_X(t)/|t|\to1$. Biconjugacy therefore recovers $F_X$ on all of $\R$.

It remains to recover $R$ from $F_X$. Differentiating and using $\E R=1$, the function $F_X$ determines
\begin{equation}\label{SP-eq:inverse-g}
 g(t)=1-F_X'(t)=2\E\frac{R}{e^{2tR}+1}.
\end{equation}
For $0<\Re s<1$, absolute Fubini gives
\begin{equation}\label{SP-eq:mellin-inverse}
 \int_0^\infty t^{s-1}g(t)\dd t
 =2^{1-s}\Gamma(s)\eta(s)\E R^{1-s},
 \qquad \eta(s)=(1-2^{1-s})\zeta(s).
\end{equation}
The kernel identity is the Fermi--Dirac integral \cite{TOMO-dlmf}; alternatively it follows by integrating the exponential series on a positive real strip and then by analytic continuation. Absolute integrability on the line $\Re s=c\in(0,1)$ follows from $\E R^{1-c}<\infty$.

Fix such a $c$. The analytic multiplier in \eqref{SP-eq:mellin-inverse} is not identically zero, so its zeros on that vertical line are isolated. Equal $g$ functions yield equal $\E R^{1-c-i\omega}$ away from those zeros, and continuity fills the isolated points. These quantities are the Fourier transforms of the finite measures
\[
 B\longmapsto\E[R^{1-c}\ind_{\{R>0,\ \log R\in B\}}].
\]
Fourier uniqueness determines the law of $R$ on $(0,\infty)$ after removing the positive weight; total mass determines its atom at zero. Symmetry then determines $X$. The converse is immediate from the normalized definition.
\end{proof}

\begin{corollary}[A finite reference cannot match a nonatomic one locally]\label{SP-cor:atomic-inverse}
The entropy frontiers of a symmetric finite atomic reference and a symmetric nonatomic reference cannot agree on a set of scales with an accumulation point in $(0,\infty)$. This applies in particular to a finite reference and a Gaussian reference.
\end{corollary}
\begin{proof}
Positive rescaling preserves atomicity. Apply Theorem~\ref{SP-thm:entropy-inverse}.
\end{proof}

\begin{corollary}[Directional entropy measurements determine a vector law]\label{SP-cor:directional-inverse}
Let $X,Y$ be centrally symmetric integrable vectors in $\R^d$. Suppose for every $v\in\R^d$ their absolute first moments agree,
\[
 m_v=\E|\langle v,X\rangle|=\E|\langle v,Y\rangle|,
\]
and, whenever $m_v>0$, the entropy frontiers of those scalar projections agree on a set with an accumulation point in $(0,\infty)$, which may depend on $v$. Then $X\stackrel d=Y$.
\end{corollary}
\begin{proof}
Theorem~\ref{SP-thm:entropy-inverse} and the known $m_v$ identify every nonzero scalar projection. A projection with $m_v=0$ vanishes almost surely. Fourier transforms of the vector laws therefore agree at every $v$.
\end{proof}

\begin{remark}[Normalization and stability]
The absolute first moments in the directional statement are necessary data. All nondegenerate one-dimensional Gaussian projections have the same normalized frontier, regardless of their variance; normalized curves alone therefore do not identify the covariance of a Gaussian vector. The accumulating-scale result is an exact uniqueness statement. It provides no stable analytic-continuation estimate from noisy measurements.
\end{remark}

The proof joins two parts of the argument that use different observables: optimal high-dimensional coding entropy and the nonlinear cost of one fair binary calibration. Their convex duals agree. Exact values on a sequence of distinct positive scales converging to one positive scale consequently identify the entire normalized scalar reference, including its atoms and tails.

\par

\subsection{Small source values determine the critical entropy exponent}\label{SP-sec:critical}
Keep a symmetric reference $W$ with $m=\E|W|\in(0,\infty)$, and put $R=|W|$, $p_0=\Pp(R=0)$. It is convenient to write
\[
 \delta=\frac{m\varepsilon}{1+\varepsilon},\qquad
 K(t)=\E\log(1+e^{-2tR}).
\]
Corollary~\ref{R9-cor:binary} gives the exact positive-scale formula
\begin{equation}\label{SP-eq:critical-var}
 \Delta_W(\varepsilon)=\inf_{t\ge0}\{\delta t+K(t)\}.
\end{equation}
Only the distribution close to zero enters its leading near-endpoint asymptotic. The rest of the reference enters through $m$.

\begin{theorem}[Power-law small values]\label{SP-thm:power-entropy}
Suppose, for some $\alpha,c>0$,
\[
 \Pp(0<R\le x)\sim c x^\alpha\qquad(x\downarrow0).
\]
Set
\[
 A_\alpha=\frac{c\Gamma(\alpha+1)\eta(\alpha+1)}{2^\alpha},
 \qquad
 C_\alpha=(\alpha+1)\alpha^{-\alpha/(\alpha+1)}
                 A_\alpha^{1/(\alpha+1)}.
\]
Then, as $\varepsilon\downarrow0$,
\begin{equation}\label{SP-eq:power-entropy}
 \Delta_W(\varepsilon)
 =p_0\log2+C_\alpha
       \left(\frac{m\varepsilon}{1+\varepsilon}\right)^{\alpha/(\alpha+1)}
       (1+o(1)).
\end{equation}
For the term after $p_0\log2$, the displayed relative asymptotic is intended. The minimizing parameter satisfies
\[
 t_\delta\sim(\alpha A_\alpha/\delta)^{1/(\alpha+1)}.
\]
Every positive-scale frontier here is attained asymptotically by the exact reference couplings of Theorem~\ref{R9-thm:frontier}. No density assumption or moment beyond $m<\infty$ is needed.
\end{theorem}
\begin{proof}
Let $F_+(x)=\Pp(0<R\le x)$. Stieltjes integration by parts yields
\[
 K(t)-p_0\log2
 =2t\int_0^\infty\frac{F_+(x)}{e^{2tx}+1}\dd x
 =\int_0^\infty\frac{F_+(u/(2t))}{e^u+1}\dd u.
\]
On a fixed small interval, $F_+(x)\le Cx^\alpha$. Beyond that interval the denominator makes the integral exponentially small after multiplication by $t^\alpha$. Dominated convergence therefore gives
\[
 K(t)-p_0\log2\sim A_\alpha t^{-\alpha},
\]
since $\int_0^\infty u^\alpha/(e^u+1)\dd u=\Gamma(\alpha+1)\eta(\alpha+1)$.

The trial $t=\delta^{-1/(\alpha+1)}$ makes the infimum in \eqref{SP-eq:critical-var} tend to $p_0\log2$. Every minimizing sequence consequently tends to infinity, since $K(t)>p_0\log2$ at finite $t$. Bound its second term above and below by $(1\pm\zeta)A_\alpha t^{-\alpha}$ for large $t$, and minimize $\delta t+A t^{-\alpha}$ explicitly. This gives \eqref{SP-eq:power-entropy} on letting $\zeta\downarrow0$. Rescaling by $(\alpha A_\alpha/\delta)^{1/(\alpha+1)}$ gives a strictly minimized limiting function, which also proves the optimizer asymptotic.
\end{proof}

\begin{example}[Every exponent strictly between zero and one]
Let $R$ have distribution function $x^\alpha$ on $[0,1]$, and give it an independent fair sign. Then $m=\alpha/(\alpha+1)$ and $c=1$. Its frontier has critical exponent $\alpha/(\alpha+1)$. Varying $\alpha>0$ realizes every exponent in $(0,1)$ with bounded references and exact code-attainment laws. A reference with a continuous positive density $f$ at zero has $\alpha=1$, $c=2f(0)$, and
\[
 \Delta_W(\varepsilon)\sim
 \pi\sqrt{\frac{m f(0)}3}\,\sqrt\varepsilon.
\]
For $W\sim N(0,1)$ this becomes $\sqrt{\pi/3}\sqrt\varepsilon$.
\end{example}

\begin{theorem}[A gap above zero]\label{SP-thm:gapped-entropy}
Suppose $r_*:=\operatorname*{ess\,inf}(R\mid R>0)>0$. Then
\begin{equation}\label{SP-eq:gapped-entropy}
 \Delta_W(\varepsilon)-p_0\log2
 \sim \frac{\delta}{2r_*}\log(1/\delta),
 \qquad \delta=\frac{m\varepsilon}{1+\varepsilon}\downarrow0.
\end{equation}
\end{theorem}
\begin{proof}
The upper bound $\log(1+z)\le z$ gives
$K(t)-p_0\log2\le(1-p_0)e^{-2r_*t}$. For every $a>r_*$, the event $0<R\le a$ has positive mass, and $\log(1+z)\ge z/2$ on $0\le z\le1$ gives a positive multiple of $e^{-2at}$ as a lower bound. Hence
\[
 -t^{-1}\log(K(t)-p_0\log2)\longrightarrow2r_*.
\]
As in the previous proof, minimizing parameters tend to infinity. Sandwich the positive part of $K$ between exponentials with rates $2r_*\pm\zeta$ and minimize $\delta t+e^{-bt}$. Its leading value is $(\delta/b)\log(1/\delta)$. Let $\zeta\downarrow0$.
\end{proof}

\begin{remark}[The zero atom and the endpoint]
The term $p_0\log2$ is the right limit of the positive-scale frontier, given by Theorem~\ref{R9-thm:frontier}. The statements here concern $\varepsilon>0$ tending to zero. They do not add an attainment assertion for the exact critical scale in the presence of a zero atom. When $p_0=0$, equality in the absolute-value comparison forces the product sign law at the endpoint.
\end{remark}

\begin{proposition}[Top-gap version for a finite alphabet]\label{SP-prop:top-gap}
Let the permutation-symmetric simplex reference of Theorem~\ref{R9-thm:frontier} have a unique largest coordinate almost surely. Write its ordered coordinates as $Z_{(1)}>Z_{(2)}\ge\cdots$, put $D_j=Z_{(1)}-Z_{(j)}$, and $b=\E Z_{(1)}$. Suppose
\[
 \Pp(D_2\le x)\sim c x^\alpha,
 \qquad \E\sum_{j=3}^q e^{-tD_j}=o(t^{-\alpha}).
\]
Then its entropy deficit at score $b-\delta$ is
\[
 (\alpha+1)\alpha^{-\alpha/(\alpha+1)}
 [c\Gamma(\alpha+1)\eta(\alpha+1)]^{1/(\alpha+1)}
 \delta^{\alpha/(\alpha+1)}(1+o(1)).
\]
For prime-power $q$ this is an attained asymptotic within the exact source class of Theorem~\ref{R9-thm:frontier}; for other $q$ it remains the corresponding variational lower-entropy bound.
\end{proposition}
\begin{proof}
The remainder in the log partition function is
$K_q(t)=\E\log(1+\sum_{j=2}^q e^{-tD_j})$. Its difference from
$\E\log(1+e^{-tD_2})$ is nonnegative and at most
$\E\sum_{j=3}^q e^{-tD_j}$. Repeat the power-law calculation without the factor two in the exponential, then minimize $\delta t+K_q(t)$.
\end{proof}


\subsection{A second-order Gaussian simplex entropy expansion}\label{SP-sec:simplex}
Let $G_1,\ldots,G_q$ be independent standard normal variables, and put
\[
 b_q=\E\max_{i\le q}G_i,
 \qquad v_i=e_i-q^{-1}\mathbf1,
 \qquad s_q^2=(q-1)/q.
\]
The projected Gaussian $Z=G-q^{-1}(\sum_iG_i)\mathbf1$ is standard on $\mathbf1^\perp$. Its least simplex-reference scale is $\kappa_q=s_q^2/b_q$. For prime-power $q$, define the normalized entropy-deficit frontier over words of simplex vertices by
\[
 \Delta_q(\varepsilon)
 =\sup_{n\ge1}\sup_{\mu\cx(1+\varepsilon)\kappa_q Z^{(n)}}
       \frac{n\log q-H(\mu)}n.
\]
Theorem~\ref{R9-thm:frontier} gives
\begin{equation}\label{SP-eq:simplex-variational}
 \Delta_q(\varepsilon)=\inf_{t\ge0}\{\delta t+K_q(t)\},
 \quad \delta=b_q\frac{\varepsilon}{1+\varepsilon},
 \quad K_q(t)=\E\log\sum_i e^{tG_i}-tb_q.
\end{equation}
The projection does not affect this expectation. The first correction to a hard maximum is governed by near ties; a centered Gumbel perturbation evaluates two orders without separately treating all multiple-tie regions.

\begin{theorem}[Gaussian simplex frontier to second order]\label{SP-thm:simplex-expansion}
Fix a prime-power $q\ge2$ and set $u=\varepsilon/(1+\varepsilon)$. Let $\varphi,\Phi$ be the standard normal density and distribution function. For $q\ge3$ define
\[
 J_q=\int_{\R}\varphi(x)^3\Phi(x)^{q-3}\dd x,
 \qquad
 c_q=\frac{2\zeta(3)}{\pi^2}q(q-1)(q-2)J_q,
\]
and put $c_2=0$. Then
\begin{equation}\label{SP-eq:simplex-expansion}
 \boxed{\quad
 \Delta_q(\varepsilon)
 =\frac{\pi b_q}{\sqrt3}\sqrt u+c_q u+O_q(u^{3/2}).
 \quad}
\end{equation}
Uniform finite-field code laws attain the exact frontier asymptotically at every fixed positive scale. The expansion is for fixed $q$; the error term is not asserted uniform as $q$ grows.
\end{theorem}

\begin{lemma}[A smooth maximum and centered Gumbel noise]\label{SP-lem:gumbel}
Let $Y_i$ be independent standard Gumbel variables after subtracting Euler's constant. They satisfy
\[
 \E Y_i=0,\qquad \E Y_i^2=\pi^2/6,\qquad \E Y_i^3=2\zeta(3).
\]
For $F(h)=\E\max_i(G_i+h_i)$,
\begin{equation}\label{SP-eq:gumbel-identity}
 K_q(t)=t\{\E F(Y/t)-F(0)\}.
\end{equation}
Moreover,
\begin{equation}\label{SP-eq:K-expansion}
 K_q(t)=\frac{A_q}{t}+\frac{B_q}{t^2}+O_q(t^{-3}),
 \quad A_q=\frac{\pi^2b_q}{12},\quad
 B_q=\frac{\zeta(3)}6q(q-1)(q-2)J_q,
\end{equation}
where $B_2=0$.
\end{lemma}
\begin{proof}
For uncentered standard Gumbels $V_i$, multiplication of their distribution functions gives
\[
 \Pp\{\max_i(a_i+V_i)\le x\}
 =\exp[-e^{-x}\sum_i e^{a_i}].
\]
Their maximum has mean $\log\sum_i e^{a_i}+\gamma$. Centering and then averaging in $G$ proves \eqref{SP-eq:gumbel-identity}. The Gumbel moment-generating function is $\Gamma(1-z)$ for $z<1$, and the Taylor series of $\log\Gamma(1-z)$ gives the displayed centered moments. The maximum/log-partition relation is also used in \cite{HazanJaakkola}.

Gaussian convolution makes $F$ smooth. Every derivative of order at least one needed here is bounded: move all but one derivative onto the Gaussian density and use the a.e. bound on the gradient of $\max_i x_i$. Taylor's theorem through order three therefore has a global remainder $O_q(\|h\|^4)$. Independence and centering of $Y_i$ give
\[
 \E F(sY)=F(0)+\frac{s^2\pi^2}{12}\sum_i\partial_{ii}F(0)
                  +\frac{s^3\zeta(3)}3\sum_i\partial_{iii}F(0)
                  +O_q(s^4).
\]
To evaluate the derivatives, condition on coordinate $i$ being tied at the maximum. Equivalently, differentiate
\[
 \partial_iF(h)=\int_\R\varphi(x-h_i)
                         \prod_{j\ne i}\Phi(x-h_j)\dd x.
\]
At zero, integration by parts gives
\begin{align*}
 \sum_i\partial_{ii}F(0)
 &=q(q-1)\int\varphi(x)^2\Phi(x)^{q-2}\dd x=b_q,\\
 \sum_i\partial_{iii}F(0)
 &=q(q-1)\int x\varphi(x)^2\Phi(x)^{q-2}\dd x\\
 &=\frac{q(q-1)(q-2)}2\int\varphi(x)^3\Phi(x)^{q-3}\dd x.
\end{align*}
For the first equality with $b_q$, integrate the density formula
$b_q=q\int x\varphi(x)\Phi(x)^{q-1}\dd x$. The last integral in the second line vanishes at $q=2$ by oddness. The second identity can also be obtained by integrating the distributional third derivative of the maximum against the Gaussian density; this justifies the tie interpretation. Substitution with $s=1/t$ proves \eqref{SP-eq:K-expansion}.
\end{proof}

\begin{proof}[Proof of Theorem~\ref{SP-thm:simplex-expansion}]
Minimizing parameters in \eqref{SP-eq:simplex-variational} tend to infinity. The first term of \eqref{SP-eq:K-expansion} shows they have order $\delta^{-1/2}$. Uniformly on that scale,
\[
 \inf_{t>0}\{\delta t+A_q/t+B_q/t^2+O_q(t^{-3})\}
 =2\sqrt{A_q\delta}+\frac{B_q}{A_q}\delta+O_q(\delta^{3/2}).
\]
For completeness, first minimize $\delta t+A_q/t+B_q/t^2$: its minimizer is
$\sqrt{A_q/\delta}+B_q/A_q+O_q(\sqrt\delta)$. Restricting to a fixed multiplicative neighbourhood of $\sqrt{A_q/\delta}$ changes no minimizer for small $\delta$, and the error term is uniformly $O_q(\delta^{3/2})$ there. Since $\delta=b_qu$, the leading coefficient is $\pi b_q/\sqrt3$, and $b_qB_q/A_q=12B_q/\pi^2=c_q$.
\end{proof}

\begin{example}[Three labels]\label{SP-ex:q3}
For $q=3$,
\[
 b_3=\frac3{2\sqrt\pi},\qquad J_3=\frac1{2\pi\sqrt3},
\]
so
\[
 \Delta_3(\varepsilon)
 =\frac{\sqrt{3\pi}}2\sqrt{\frac\varepsilon{1+\varepsilon}}
 +\frac{2\sqrt3\,\zeta(3)}{\pi^3}
                      \frac\varepsilon{1+\varepsilon}
 +O(\varepsilon^{3/2}).
\]
The term of order $\varepsilon$ is positive; the binary expansion has no corresponding term when expressed in $u=\varepsilon/(1+\varepsilon)$.
\end{example}

\begin{remark}[Scope for other alphabets]
The Gaussian log-partition expansion holds for every integer $q\ge2$. For prime powers, the exact code construction in Theorem~\ref{R9-thm:frontier} identifies \eqref{SP-eq:simplex-expansion} with the operational frontier. For general integer alphabets, the same calculation determines the variational entropy bound; operational equality additionally requires exact attainment.
\end{remark}


\section{Fixed reference laws and the complete dependence boundary}
\label{R11-sec:dependence}
For fixed common and coordinate noises, the empirical directing mean must fit the common law, and the residual alphabet law must fit the coordinate law. These conditions characterize feasibility in every dimension and determine the maximal entropy rate.

Gaussian noises give an explicit boundary and minimum dependence. When the common noise may vary and depend on coordinate noises, cancellation leaves an exact permitted set of empirical means. Smooth coordinate-density perturbations can prescribe that set even though every fixed dimension still admits independent signs.

\subsection{An exact decomposition for exchangeable reference laws}\label{r3:sec:exchangeable}
Let $\mathcal D=\{d_1,\ldots,d_q\}\subset\mathbb R^d$ be finite. For $p\in\Delta_q$, write
\[
 m(p)=\sum_{i=1}^qp_i d_i,
 \qquad
 \nu_p^0=\sum_{i=1}^qp_i\delta_{d_i-m(p)},
 \qquad
 H(p)=-\sum_{i=1}^qp_i\log p_i.
\]
Let $U,V_1,V_2,\ldots$ be independent, centred, integrable random vectors, with the $V_i$ identically distributed according to $\nu$. The reference in dimension $n$ is
\[
 Y^{(n)}=(U+V_1,\ldots,U+V_n).
\]
A random probability vector $P\in\Delta_q$ is a directing law for an exchangeable output if, conditional on $P$, the outputs are independent with probabilities $P$. The existence and uniqueness of this description are the finite-alphabet case of de Finetti's theorem \cite{HewittSavage}.

\begin{theorem}\label{r3:thm:factorization}
Fix $\mu\in\operatorname{conv}(\mathcal D)$. There are $\mathcal D$-valued random vectors $X^{(n)}$, with mean $\mu$ in each coordinate, such that
\begin{equation}\label{r3:eq:alldimfinitealphabet}
 (X_1^{(n)}-\mu,\ldots,X_n^{(n)}-\mu)
       \preceq_{\rm cx}Y^{(n)}
 \quad\hbox{for every }n
\end{equation}
if and only if there is a random $P\in\Delta_q$ such that
\begin{equation}\label{r3:eq:factorconditions}
 \mathbb E m(P)=\mu,\qquad
 m(P)-\mu\preceq_{\rm cx}U,
 \qquad \nu_P^0\preceq_{\rm cx}\nu\quad\hbox{almost surely}.
\end{equation}
More precisely, \eqref{r3:eq:factorconditions} characterizes exactly the directing laws of exchangeable processes satisfying \eqref{r3:eq:alldimfinitealphabet} for every $n$.

Let $H_n^*$ be the largest entropy of an $n$-coordinate output satisfying \eqref{r3:eq:alldimfinitealphabet}. When the equivalent feasibility conditions hold,
\begin{equation}\label{r3:eq:entropyvariational}
 \lim_{n\to\infty}\frac{H_n^*}{n}
 =\sup\left\{\mathbb E H(P):P\text{ satisfies }\eqref{r3:eq:factorconditions}\right\}.
\end{equation}
The supremum is attained by a single directing law. The finite-dimensional maximizing laws need not be projectively consistent.
\end{theorem}

The two conditions in \eqref{r3:eq:factorconditions} separate the roles of the noises. The random barycenter must be dominated by the common noise, and the centred conditional output law must be dominated by the coordinate noise. In particular, increasing the common noise leaves the latter constraint unchanged.

\begin{proof}
Suppose first that \eqref{r3:eq:alldimfinitealphabet} holds in every dimension. Symmetrize a martingale coupling by simultaneous permutations of its coordinate pairs. The additive reference is exchangeable for arbitrary integrable common and iid coordinate laws. Tightness, a diagonal subsequence, and uniform integrability of each reference coordinate give an infinite exchangeable coupling $(X_i,Y_i)_{i\ge1}$ with
\begin{equation}\label{r3:eq:infinitemartingale}
 \mathbb E[Y_i\mid X_1,X_2,\ldots]=X_i-\mu.
\end{equation}
To justify the passage, first test a fixed coordinate against a bounded function of finitely many outputs. The identity holds in every sufficiently large finite coupling and passes to the limit. Conditional-expectation convergence then gives \eqref{r3:eq:infinitemartingale}.

Apply de Finetti's theorem to the pairs. Conditional on their directing measure, the pairs are iid. Its first marginal is a random probability vector $P$. By uniqueness of the directing measure of the reference sequence, its second marginal is a translate $u+\nu$, where $u$ has the law of $U$. The empirical means converge in $L^1$ to $m(P)$ and $u$, respectively. Averaging \eqref{r3:eq:infinitemartingale} therefore gives
\[
 \mathbb E[u\mid X_1,X_2,\ldots]=m(P)-\mu,
 \qquad \mathbb E[u\mid P]=m(P)-\mu.
\]
This proves the barycenter comparison.

Given $P$ and $u$, the law of $X_i$ is still $P$, so
$\mathbb E[u\mid X_i,P]=m(P)-\mu$. Also, conditional on the full directing measure, $Y_i-u$ has law $\nu$; it therefore has law $\nu$ conditional on $P$. Subtracting the preceding identity from \eqref{r3:eq:infinitemartingale}, after conditioning on $(X_i,P)$, yields $\mathbb E[Y_i-u\mid X_i,P]=X_i-m(P)$. This proves $\nu_P^0\preceq_{\rm cx}\nu$ almost surely. The same argument applies to any fixed consistent exchangeable output process and identifies its own directing law.

Conversely, couple $P$ and $U$ so that
$\mathbb E[U\mid P]=m(P)-\mu$. This is obtained by a martingale coupling with $m(P)$ and then disintegration over $m(P)$. Conditional on $P=p$, take a martingale coupling of a $\mathcal D$-valued $X$ and a $V$ of law $\nu$ satisfying $\mathbb E[V\mid X,P=p]=X-m(p)$. Choose independent copies of this conditional coupling for the coordinates, independently of $U$ given $P$. The marginal law of each $V_i$ is $\nu$ for every $p$, so the $V_i$ are iid and independent of $(P,U)$. Since $X_i$ has law $P$ independently of $U$ conditional on $P$, $\mathbb E[U+V_i\mid X_1,\ldots,X_n,P] =m(P)-\mu+X_i-m(P)=X_i-\mu$. Removing $P$ from the conditioning proves the required martingale coupling.

For completeness, the conditional couplings can be chosen measurably. A finite-alphabet coupling is represented by kernels $k_i(v)\in[0,1]$ with $\sum_i k_i=1$, prescribed masses $p_i$, and prescribed first moments $p_i(d_i-m(p))$. The unit ball is compact and metrizable in the weak-star topology of $L^\infty(\nu)$; the constraints are continuous because $\nu$ has a first moment. A measurable choice is obtained by minimizing the strictly convex lower-semicontinuous functional $\sum_i\int k_i^2\,d\nu$ on each feasible set. Its unique minimizer gives a measurable version after completion. This also covers singular coordinate laws.

It remains to prove \eqref{r3:eq:entropyvariational}. Symmetrization preserves feasibility and increases entropy. An exchangeable $k$-coordinate law is uniform conditional on its vector of alphabet counts. There are at most $(k+1)^q$ count vectors, so its entropy lies between $k\mathbb E H(P)$ and
\[
 q\log(k+1)+k\,\mathbb E H(\widehat P_k),
\]
where $\widehat P_k$ is the empirical probability vector. Conditional laws of large numbers and boundedness of $H$ give entropy rate $\mathbb E H(P)$.

For arbitrary maximizing finite-dimensional laws, pass to a subsequential exchangeable limit as above. For any fixed $k$, subadditivity gives
\[
 H(X_1^{(n)},\ldots,X_n^{(n)})
 \le\lfloor n/k\rfloor H(X_1^{(n)},\ldots,X_k^{(n)})+k\log q.
\]
First let $n$ tend to infinity along the subsequence and then let $k$ tend to infinity. The upper bound in \eqref{r3:eq:entropyvariational} follows. Every feasible directing law gives a consistent family attaining its own rate, which proves the reverse inequality. Finally, the set of allowed probability vectors is closed in the compact simplex, because convex order is closed under convergence of the bounded output laws. The constraint on the bounded random barycenter is also closed under weak convergence of its law. The feasible directing laws therefore form a compact set, and the continuous functional $\mathbb E H(P)$ attains its supremum.
\end{proof}

\begin{proposition}\label{r3:prop:finitecommon}
If the common noise $U$ has $r$ atoms, an entropy-maximizing directing law in Theorem~\ref{r3:thm:factorization} can be chosen with at most $r$ atoms. This support bound is sharp when the alphabet and its ambient dimension are allowed to vary.
\end{proposition}
\begin{proof}
Write $\mathbb P(U=u_j)=\pi_j>0$, $1\le j\le r$. A martingale coupling with a probability vector $p$ is described by a posterior $\lambda\in\Delta_r$ satisfying
\[
 \sum_j\lambda_j u_j=m(p)-\mu.
\]
Consider the compact set of pairs $(p,\lambda)$ satisfying this identity and $\nu_p^0\preceq_{\rm cx}\nu$. Feasible directing laws are obtained from probability measures on this set subject to the $r$ constraints $\mathbb E\lambda_j=\pi_j$. Maximizing $\mathbb E H(p)$ is a continuous linear optimization over a compact convex set of measures, so an extreme maximizing measure exists.

Such an extreme measure has at most $r$ atoms. Otherwise, select $r+1$ disjoint sets of positive measure. The $r+1$ vectors of integrated coordinates $\lambda_j$ lie in $\mathbb R^r$ and are linearly dependent. The resulting bounded signed perturbation, constant on each selected set, preserves all $r$ constraints; its total mass is also zero because $\sum_j\lambda_j=1$. Small perturbations of both signs contradict extremality. Projecting this measure onto $p$ gives the asserted support bound.

For sharpness, take $r$ affinely independent alphabet points, put $\mu$ in the relative interior of their simplex, let $V=0$, and let $U$ have the corresponding centred vertex law. The residual condition forces each conditional output law to be a point mass. The unique barycentric representation of $\mu$ uses every vertex, so every feasible directing law has at least $r$ atoms.
\end{proof}

\subsection{The complete Gaussian boundary for sign laws}\label{r3:sec:gaussianboundary}
Let $\phi$ and $\Phi$ be the standard normal density and distribution function, and put
\[
 I(p)=\phi(\Phi^{-1}(p)),\quad 0<p<1,
 \qquad I(0)=I(1)=0.
\]
The function $I$ is the Gaussian isoperimetric function. The scalar quantity needed below is
\begin{equation}\label{r3:eq:kappa}
 \kappa(m)=\frac{1-m^2}{2I((1+m)/2)},\qquad -1<m<1,
 \qquad \kappa(-1)=\kappa(1)=0.
\end{equation}

\begin{lemma}\label{r3:lem:bernoulli}
Let $S\in\{-1,1\}$ have mean $m$. Then
\[
 S-m\preceq_{\rm cx}N(0,b^2)
 \quad\Longleftrightarrow\quad b\ge\kappa(m).
\]
The function $\kappa$ is continuous and even, strictly decreasing on $[0,1]$, with $\kappa(0)=\sqrt{\pi/2}$.
\end{lemma}
\begin{proof}
Write $p=(1+m)/2$ and let $G$ be standard normal. Among events, including randomized events, of probability $p$, the largest value of $\mathbb E[G\mathbf1_E]$ is $I(p)$, attained by the upper Gaussian tail. If $\mathbb E[bG\mid S]=S-m$, then $\mathbb E[bG\mathbf1_{S=1}]=p(1-m)=2p(1-p)$. Thus $bI(p)\ge2p(1-p)$ is necessary. It is sufficient by mixing the upper-tail assignment of $S$ with an independent Bernoulli assignment of the same mean, using the mixing probability $\kappa(m)/b$. Both conditional moments then have the required values. The endpoint cases are deterministic.

Evenness and continuity are immediate, with the limits at the endpoints following from the Gaussian Mills ratio. We give the monotonicity argument because it determines the topology of the feasible set. Write $z=\Phi^{-1}((1+m)/2)\ge0$, $J(z)=\int_0^z\phi(t)\,dt$, and
\[
 f(z)=2\phi(z)J(z)-z\bigl(1/4-J(z)^2\bigr).
\]
Differentiation gives
\[
 f'(z)=2\phi(z)^2+J(z)^2-1/4,
 \qquad f''(z)=2\phi(z)\bigl(J(z)-2z\phi(z)\bigr).
\]
The derivative of $J(z)-2z\phi(z)$ is $(2z^2-1)\phi(z)$. Hence this expression is negative and then positive, with one positive zero. Since $f'(0)=1/\pi-1/4>0$ and $f'(z)\to0$, the function $f'$ is positive and then negative. As $f(0)=0$ and $f(z)\to0$, it follows that $f(z)>0$ for $z>0$. Finally,
\[
 \frac{d}{dz}\left(\frac{2\Phi(z)(1-\Phi(z))}{\phi(z)}\right)
       =-\frac{2f(z)}{\phi(z)}<0.
\]
The expression differentiated is $\kappa(m)$, proving strict monotonicity.
\end{proof}

For $b\ge0$, define $r_b\in[0,1]$ by
\begin{equation}\label{r3:eq:rb}
 r_b=0\quad\hbox{if }b\ge\sqrt{\pi/2},
 \qquad \kappa(r_b)=b\quad\hbox{otherwise}.
\end{equation}
In particular, $r_0=1$. Let $h(p)=-p\log p-(1-p)\log(1-p)$, with the usual continuous endpoint values.

\begin{theorem}\label{r3:thm:gaussianboundary}
Fix $a,b\ge0$ and $\mu\in[-1,1]$. In dimension $n$, let
\[
 G_i=aZ_0+bZ_i,\qquad 1\le i\le n,
\]
where all $Z_i$ are independent standard normals. There is a sign law with $\mathbb E S_i=\mu$ and
\[
 (S_1-\mu,\ldots,S_n-\mu)\preceq_{\rm cx}(G_1,\ldots,G_n)
\]
for every $n$ if and only if
\begin{equation}\label{r3:eq:boundary}
 a\ge
 \begin{cases}
 r_b\,\kappa(\mu/r_b),& |\mu|<r_b,\\
 0,&|\mu|\ge r_b.
 \end{cases}
\end{equation}
Under this condition, if $H_n^*$ denotes the maximum entropy of a feasible sign law in dimension $n$, then
\begin{equation}\label{r3:eq:entropyrate}
 \lim_{n\to\infty}\frac{H_n^*}{n}
 =h\!\left(\frac{1+\max\{|\mu|,r_b\}}2\right).
\end{equation}
A single consistent family attains this rate. It has an explicit martingale coupling which uses $n+1$ Gaussian draws and $O(n)$ scalar operations after calibrating its scalar parameters.
\end{theorem}

For unbiased signs, the all-dimensional boundary reduces to
\begin{equation}\label{r3:eq:unbiasedboundary}
 a\ge\sqrt{\pi/2}\,r_b,
 \qquad
 \lim_{n\to\infty}H_n^*/n=h((1+r_b)/2).
\end{equation}
Thus the coordinate standard deviation $b$ alone determines the largest entropy rate whenever the common component is large enough for feasibility.

\begin{proof}
Theorem~\ref{r3:thm:factorization} and Lemma~\ref{r3:lem:bernoulli} reduce feasibility to the existence of $M\in[-1,1]$ satisfying
\[
 \mathbb EM=\mu,\qquad M-\mu\preceq_{\rm cx}N(0,a^2),
 \qquad |M|\ge r_b\quad\hbox{almost surely}.
\]
If $|\mu|\ge r_b$, the constant choice $M=\mu$ works.

Suppose $|\mu|<r_b=r$. Among all laws supported outside $(-r,r)$ with mean $\mu$, the two-point law on $\{-r,r\}$ with that mean is minimal in convex order. Indeed, a convex function lies above the affine extension of its chord between $-r$ and $r$ on the complement of $(-r,r)$, and expectation of that chord depends only on the mean. Write this minimal random variable as $rB$, where $B$ is a sign with mean $\mu/r$. Lemma~\ref{r3:lem:bernoulli} says exactly that $rB-\mu$ is dominated by $N(0,a^2)$ when $a\ge r\kappa(\mu/r)$. This proves both necessity and sufficiency, including the boundary.

For the entropy, if $|\mu|\le r$, the constraint $|M|\ge r$ gives
$h((1+M)/2)\le h((1+r)/2)$ pointwise. The two-point minimizing law attains equality. If $|\mu|\ge r$, Jensen's inequality gives
$\mathbb E h((1+M)/2)\le h((1+\mu)/2)$, and the constant law attains it. The entropy-rate formula now follows from Theorem~\ref{r3:thm:factorization}.

Here is the coupling explicitly. In the interior-gap case, set $r=r_b$, $q=(1+\mu/r)/2$, and $\lambda_0=r\kappa(\mu/r)/a$. Draw $Z_0$. With probability $\lambda_0$, let $B=1$ precisely when $Z_0>\Phi^{-1}(1-q)$; with the remaining probability, take an independent sign $B$ with probability $q$ of being $1$. Put $M=rB$. Then $\mathbb E[aZ_0\mid M]=M-\mu$.

For every coordinate, conditional on $M=m$, let $p=(1+m)/2$ and $\lambda=\kappa(m)/b$. With probability $\lambda$, set $S_i=1$ precisely when $Z_i>\Phi^{-1}(1-p)$; otherwise use an independent sign of mean $m$. Take these choices independently across coordinates. The deterministic endpoint cases use constant signs. This gives
\[
 \mathbb E[bZ_i\mid S_i,M]=S_i-M,
 \qquad
 \mathbb E[aZ_0+bZ_i\mid S_1,\ldots,S_n,M]=S_i-\mu.
\]
When $|\mu|\ge r_b$, take $M=\mu$ and use only the coordinate construction; the common normal is independent. In the interior-gap construction, the entropy differs from $n h((1+r_b)/2)$ by the mutual information between $M$ and the first $n$ signs, which lies between $0$ and $\log2$. This also proves the claimed rate directly.
\end{proof}

\begin{corollary}\label{r3:cor:dependence}
Under the feasibility condition of Theorem~\ref{r3:thm:gaussianboundary}, let $M_*$ be constant $\mu$ when $|\mu|\ge r_b$, and otherwise let it be supported on $\{-r_b,r_b\}$ with mean $\mu$. For every feasible infinite exchangeable sign process with directing mean $M$,
\[
 M_*\preceq_{\rm cx}M,
 \qquad
 \operatorname{Cov}(S_i,S_j)=\operatorname{Var}(M)
 \ge (r_b^2-\mu^2)_+\quad(i\ne j).
\]
The construction of Theorem~\ref{r3:thm:gaussianboundary} attains equality. Moreover, its scalar empirical mean is smaller in convex order, at every finite sample size, than the empirical mean of every other such consistent exchangeable process.

Among arbitrary feasible families in increasing dimensions, the same lower bound holds for the limit inferior of the average off-diagonal covariance. In particular, a family with asymptotically vanishing average off-diagonal covariance exists if and only if $b\ge\kappa(\mu)$, independently of $a$.
\end{corollary}
\begin{proof}
The chord argument in the proof of Theorem~\ref{r3:thm:gaussianboundary} gives $M_*\preceq_{\rm cx}M$ in the interior-gap case. The constant choice is minimal for the other case by Jensen's inequality. Apply the convex function $x^2$ and use conditional independence to obtain the covariance bound.

For a convex function $f$ of a scalar empirical mean, conditional expectation given $M=m$ is a Bernstein polynomial in $(1+m)/2$, with values $f(2j/n-1)$ on its grid. Its second derivative is a positive combination of the nonnegative second differences of these values. It is therefore convex in $m$, so comparison of $M_*$ and $M$ proves the empirical-mean assertion.

For finite families, symmetrize and extract a subsequential infinite limit along a sequence realizing the limit inferior of the average covariance. The preceding proof applies to this limit and proves the lower bound. The lower bound is zero exactly when $|\mu|\ge r_b$, which is equivalent to $b\ge\kappa(\mu)$. In that case iid signs of mean $\mu$ are feasible by the scalar Gaussian coupling.
\end{proof}

\subsubsection{Rare marginals and entropy restrictions}
For a Bernoulli variable $B$ of mean $p$, the exact smallest Gaussian standard deviation is
\begin{equation}\label{r3:eq:rare}
 \inf\{s:B-p\preceq_{\rm cx}N(0,s^2)\}
       =\frac{p(1-p)}{I(p)},
 \qquad
 \left(\frac{p(1-p)}{I(p)}\right)^2
       \sim\frac1{2\log(1/p)}\quad(p\downarrow0).
\end{equation}
The Gaussian Mills ratio proves the asymptotic. In comparison, the exact variance is $p(1-p)$. Thus a covariance estimate proportional to the Bernoulli variance and a Gaussian convex-order estimate have quantitatively different dependence on rare marginals. Formula~\eqref{r3:eq:rare} determines the exact Gaussian cost.

There is also an entropy restriction which is uniform in the size of the common Gaussian component. Let $h_0\in(0,\log2]$, and let $r(h_0)\in[0,1)$ solve
$h((1+r(h_0))/2)=h_0$. For any family of unbiased sign laws with entropy rate at least $h_0$, an additive Gaussian reference of the above form must satisfy
\begin{equation}\label{r3:eq:entropyobstruction}
 b\ge\kappa(r(h_0)).
\end{equation}
This conclusion applies even when the finite-dimensional sign laws are not consistent or exchangeable, because symmetrization increases their entropy and preserves the reference comparison. Any positive lower bound on the entropy rate can be used as $h_0$ in this formula.

Near $b=\sqrt{\pi/2}$, Taylor expansion at the origin gives
\[
 \kappa(r)=\sqrt{\pi/2}\left(1-(1-\pi/4)r^2+O(r^4)\right),
 \qquad h((1+r)/2)=\log2-r^2/2+O(r^4).
\]
Consequently, for unbiased signs and feasible common noise,
\[
 \log2-\lim H_n^*/n
 =\frac{2}{\pi(4-\pi)}\left(\frac\pi2-b^2\right)
     +O\!\left(\left(\frac\pi2-b^2\right)^2\right)
\]
as $b^2$ increases to $\pi/2$ from below. This gives the first-order loss of entropy at the threshold for independent signs.


\subsection{Dependence forced by a reference}\label{PHASE-sec}
The fixed-noise classification specifies both the common and coordinate laws before the dimension varies. Here only the centered integrable coordinate law $\nu$ is fixed: the common noise may change with dimension and may depend on the coordinate noises. Subtracting the empirical average of the martingale identities cancels that common term exactly. The remaining scalar convex-order condition determines all limiting empirical means, and hence the maximum entropy under these broader choices. The coordinate noises remain iid throughout.

Write $h(p)=-p\log p-(1-p)\log(1-p)$ for binary entropy, with $0\log0=0$.

For the centered integrable law $\nu$, recall the integrated upper
quantile
\[
 L_\nu(p)=\sup\{\E[V k(V)]:0\le k\le1,\ \E k(V)=p\},\qquad V\sim\nu,
\]
and the closed set
\begin{equation}\label{PHASE-eq:allowed}
 K_\nu=\{m\in[-1,1]:L_\nu((1+m)/2)\ge(1-m^2)/2\}.
\end{equation}
A sign $S$ with mean $m$ satisfies $S-m\cx\nu$ exactly when
$m\in K_\nu$.  Necessity follows from the martingale moment on the event
$S=1$.  For sufficiency, mix an upper-quantile assignment of mass
$p=(1+m)/2$ with an independent assignment of the same mass.  Randomizing
at an atom gives the same argument for arbitrary $\nu$.

\begin{theorem}[Complete empirical-mean limit set]
\label{R8-thm:empirical-limits}
Fix a centred integrable coordinate law $\nu$ and $\mu\in[-1,1]$.
For each $n$, let $S^{(n)}$ be a sign vector with coordinate means
$\mu$ satisfying
\[
 S^{(n)}-\mu\mathbf1\cx U_n\mathbf1+(V_1,\ldots,V_n),
 \qquad V_i\text{ iid with law }\nu.
\]
Allow arbitrary centred integrable $U_n$, including dependence on the
coordinate noise. The possible weak subsequential limits of
$n^{-1}\sum_iS_i^{(n)}$ are exactly
\begin{equation}\label{R8-eq:empirical-limit-set}
 \left\{\eta\in\mathcal P([-1,1]):
       \eta(K_\nu)=1,\ \int m\,\eta(dm)=\mu\right\}.
\end{equation}
Every law in this set is attained by an infinite exchangeable
construction using one fixed common noise, uniform on $\{-2,2\}$
and independent of all coordinate noise. That construction has entropy
rate $\int h((1+m)/2)\,\eta(dm)$.

For every continuous $\varphi:[-1,1]\to\R$, the limiting supremum of
$\E\varphi(n^{-1}\sum_iS_i^{(n)})$ equals
\begin{equation}\label{R8-eq:observable-limit}
 \max_{\eta:\,\eta(K_\nu)=1,\ \int m\,d\eta=\mu}
       \int\varphi(m)\,\eta(dm).
\end{equation}
A maximizer can be chosen with at most two atoms. The limiting
infimum has the corresponding minimum formula.
\end{theorem}
\begin{proof}
Take martingale couplings and write
$M_n=n^{-1}\sum_iS_i$, $\bar V_n=n^{-1}\sum_iV_i$. Subtracting the
average of the conditional-expectation identities gives
\begin{equation}\label{PHASE-eq:common-cancel}
 \E[V_i-\bar V_n\mid S]=S_i-M_n.
\end{equation}
The common noise has disappeared. Randomly permute the coordinate
pairs, retaining $M_n$. Entropy cannot decrease, the $V$ marginal
remains iid, and \eqref{PHASE-eq:common-cancel} is unchanged.
For any sequence of dimensions tending to infinity, tightness and a
diagonal subsequence give an infinite exchangeable marked sequence
$(S_i,V_i,M)_{i\ge1}$. The fixed integrable law of $V_i$ gives uniform
integrability, and $\E|\bar V_n|\to0$. Testing
\eqref{PHASE-eq:common-cancel} against bounded continuous functions of
$M_n$ and finitely many signs and passing to the limit therefore gives $\E[V_i\mid M,S_1,S_2,\ldots]=S_i-M$. For a sample of $k$ coordinates from an exchangeable $n$-sign vector,
\[
 \E\left|k^{-1}\sum_{i=1}^kS_i-M_n\right|^2
   =\frac{n-k}{k(n-1)}\E(1-M_n^2)\le\frac1k.
\]
Thus $M$ is the empirical-mean limit of the infinite sign sequence.
Apply de Finetti to the pairs. Its first directing marginal is the
Bernoulli law of mean $M$, while its second directing marginal is
$\nu$ almost surely because the reference sequence is iid. In
particular $V_i$ has law $\nu$ conditional on $M$. The last conditional
expectation therefore proves $S_i-M\cx\nu$ conditional on $M$.
Hence $M\in K_\nu$ almost surely, and $\E M=\mu$.

For sufficiency, draw $M$ with law $\eta$. Conditional on $M$, draw
$U\in\{-2,2\}$ with
$\Prb(U=2\mid M)=1/2+(M-\mu)/4$. Then $U$ is fair and
$\E[U\mid M]=M-\mu$. For each $m\in K_\nu$ choose a scalar
martingale coupling with sign mean $m$ and
$\E[V\mid S,M=m]=S-m$, $\law(V\mid M=m)=\nu$.
Measurable selection gives such a family; alternatively, use the
upper-quantile construction in \eqref{PHASE-eq:allowed}.
Conditional on $(M,U)$, take independent copies of these pairs.
The coordinate noises have the fixed product law $\nu^{\otimes n}$
conditional on $(M,U)$, so they are independent of $U$. Conditional
expectation and then the tower property give $\E[U+V_i\mid S_1,\ldots,S_n]=S_i-\mu$. The conditional law of the signs is iid with mean $M$, so their
empirical means tend almost surely to $M$. Their entropy obeys
\[
 n\E h((1+M)/2)\le H(S_1,\ldots,S_n)
 \le n\E h((1+M)/2)+\log(n+1).
\]
For the upper inequality, the count of positive signs is sufficient
for $M$ and has at most $n+1$ values. This proves the rate.

Compactness of the empirical-mean laws gives the upper bound in
\eqref{R8-eq:observable-limit}, while the preceding construction gives
its lower bound. A continuous linear functional on the compact
moment-constrained probability set has an extreme maximizer. Such an
extreme law has at most two atoms: positive mass on three disjoint
neighborhoods permits a nonzero signed perturbation with zero total
mass and zero first moment, contradicting extremality. Apply the same
argument to $-\varphi$ for the infimum.
\end{proof}
\begin{corollary}[Exact entropy rate]\label{PHASE-thm:free-common}
For coordinate mean $\mu$, let $M_*=\mu$ if $\mu\in K_\nu$.
Otherwise let $M_*$ be supported on the nearest points $a<\mu<b$
of $K_\nu$, with mean $\mu$. Among all feasible finite laws, allowing
the common noise to vary with dimension, the maximum asymptotic
entropy per sign is $\E h((1+M_*)/2)$. The least limiting average
off-diagonal covariance is $\Var M_*$, and both optima are attained
with the fixed fair common source on $\{-2,2\}$.
\end{corollary}
\begin{proof}
Every law on $K_\nu$ with mean $\mu$ dominates $M_*$ in convex order:
outside $(a,b)$ a convex function lies above the extension of its chord
between $a,b$. Entropy counting gives
$H(S)\le\log(n+1)+n\E h((1+M_n)/2)$ for $M_n=n^{-1}\sum_iS_i$.
Use the limiting classification and concavity of binary entropy for the
upper bound. Its construction attains the bound. The average covariance
is $[n\Var M_n-(1-\mu^2)]/(n-1)$, giving the second conclusion.
\end{proof}

\subsubsection{Arbitrary closed sets of conditional means}
\begin{theorem}[Realization of permitted directing means]
\label{PHASE-thm:closed}
Every closed $K\subset[-1,1]$ containing $-1,1$ is $K_\nu$ for a centred
law $\nu$ with a strictly positive smooth density on $[-2,2]$.  For each
fixed $r$, the density can be chosen arbitrarily close to $1/4$ in
$C^r([-2,2])$ and in total variation.  Symmetric $K$ admit symmetric
choices of $\nu$.

With the independent common noise $U$ uniform on $\{-2,2\}$, all the
directing laws supported on $K$ with a specified mean $\mu$ are feasible.
If $a\le\mu\le b$ are the nearest points of $K$ bracketing $\mu$, the
maximum entropy rate is the linear interpolation between
$h((1+a)/2)$ and $h((1+b)/2)$, and the minimum off-diagonal covariance is
$(\mu-a)(b-\mu)$.
\end{theorem}
\begin{proof}
Put $K'=(1+K)/2$.  Choose a nonnegative $g\in C^\infty([0,1])$, flat on
$K'$, whose zero set is $K'$.  Such a function is the sum of positive
smooth bumps on the complementary intervals, with coefficients chosen
to make every derivative series uniformly convergent.  Average under
$p\mapsto1-p$ for symmetric $K'$.  For small $\delta>0$ put
\[
 L(p)=2p(1-p)-\delta g(p).
\]
Then $L(0)=L(1)=0$ and $L''<0$.  If $T$ is uniform on $(0,1)$, the
variable $V=L'(T)$ has mean zero, support $[-2,2]$, and density
\[
 f_\delta(L'(p))=\frac1{4+\delta g''(p)}.
\]
The inverse map is smooth, including at the endpoints, and tends to the
uniform-law inverse map in every fixed derivative norm.  Its integrated
upper quantile is $L$, so \eqref{PHASE-eq:allowed} holds exactly on $K$.
Smoothness here is up to the endpoints of the support; the zero
extension need not be smooth across them.

For any $M\in[-1,1]$ with mean $\mu$, prescribe
$\Pr(U=2\mid M)=1/2+(M-\mu)/4$.  This gives a fair marginal $U$ and
$\E[U\mid M]=M-\mu$.  The construction in Theorem~\ref{R8-thm:empirical-limits} therefore applies to every such directing law.  The entropy and covariance assertions follow
from the chord argument above, the entropy formula in Theorem~\ref{R8-thm:empirical-limits}, and
$\Cov(S_i,S_j)=\Var M$.
\end{proof}

In particular the permitted conditional means may form a Cantor set,
although the coordinate density is smooth and arbitrarily close to a
uniform density.  The irregularity is in the exact set of admissible
conditional means, rather than in the regularity of the noise density.

\subsubsection{Finite-dimensional independence and all-dimensional collapse}
Let $P_2(u)=(3u^2-1)/2$, and for $0<\varepsilon<1$ put
\begin{equation}\label{PHASE-eq:perturb}
 f_\varepsilon(v)=\tfrac14[1-\varepsilon P_2(v/2)]\ind_{(-2,2)}(v).
\end{equation}
This symmetric density is strictly log-concave on its support and tends
smoothly to the uniform density.  Fix independent common noise uniform
on $\{-2,2\}$, or $N(0,9)$.

\begin{theorem}[Noncommuting dimension and perturbation limits]
\label{PHASE-thm:collapse}
At $\varepsilon=0$, iid fair signs are feasible in every dimension.
For every $\varepsilon>0$, each feasible infinite exchangeable fair-sign
process is a single fair sign repeated in every coordinate.
Nevertheless, for every fixed $n$, iid fair $n$-signs remain feasible
for all sufficiently small $\varepsilon>0$.

Let $H_n^*(\varepsilon)$ be the maximum over all feasible fair-sign laws
in dimension $n$, without an extension requirement.  Then
\[
 \lim_{\varepsilon\downarrow0}\lim_{n\to\infty}
       H_n^*(\varepsilon)/n=0,\qquad
 \lim_{n\to\infty}\lim_{\varepsilon\downarrow0}
       H_n^*(\varepsilon)/n=\log2.
\]
For the quantitative statement, put
\[
 \delta_{n,\varepsilon}=\min\left\{1,
 \sqrt{\frac{32}{\varepsilon\sqrt n}},
 \frac{32}{\varepsilon(1-\varepsilon)\sqrt n}\right\},\qquad
 \psi(u)=h\!\left(\frac{1+\sqrt{1-u}}2\right).
\]
For any centred integrable common noise, even one changing with $n$ or
depending on the coordinate noise, every feasible fair-sign law satisfies
\begin{equation}\label{PHASE-eq:improved-entropy}
 \frac{H(S)}n\le\min\left\{\log2,\frac{\log(n+1)}n+
 \psi(\delta_{n,\varepsilon})\right\}.
\end{equation}
In particular the normalized entropy is
$O_\varepsilon(n^{-1/2}\log n)$.  The average off-diagonal covariance
converges to one, uniformly over the common reference.  An iid fair
block with $n\ge2$ can be feasible only if $n\le4096/\varepsilon^2$.
\end{theorem}
\begin{proof}
Writing $u=v/2$, let
\[
 A(u)=-u(1-u^2)/2,\qquad
 p_\varepsilon(u)=(1-u+\varepsilon A(u))/2,
 \qquad m=2p_\varepsilon(u)-1.
\]
Direct integration gives
\begin{equation}\label{PHASE-eq:defect}
 L_{\nu_\varepsilon}(p_\varepsilon(u))
        -2p_\varepsilon(u)(1-p_\varepsilon(u))
 =-\frac\varepsilon8(1-u^2)^2(1-\varepsilon u^2).
\end{equation}
It is negative for $-1<u<1$, so $K_{\nu_\varepsilon}=\{-1,1\}$.
Theorem~\ref{R8-thm:empirical-limits} gives the permitted empirical means; de Finetti gives the infinite-process assertion.
At zero perturbation take $S_i=\operatorname{sign}(V_i)$, which gives
the iid martingale coupling.

Here is why each fixed finite block has a margin.  In this zero-perturbation
coupling, the convex hull of the reference conditional on a sign pattern
$s$ is
\[
 C_s=[-2,2]\mathbf1+\prod_i J_{s_i},\qquad
 J_+=(0,2),\quad J_-=(-2,0).
\]
The interiors of $C_s,C_t$ overlap when $s,t$ differ in one coordinate.
Consider the finite-dimensional convex body of kernel masses and first
moments $(\int k_s\,d\rho,\int yk_s\,d\rho)_s$ with
$k_s\ge0$ and $\sum_s k_s=1$.  A supporting functional at the present
target consists of affine scores $\ell_s(y)$.  Optimality of the
coupling forces $\ell_s$ to dominate all the others on its conditional
support and hence on $C_s$.  On each open overlap the two affine scores
agree identically.  The connected sign-cube graph forces all scores to
be identical, which gives only the common affine relations on the moment
body.  Thus the target is in its relative interior.  The reference
support spans $\R^n$, so these relations are exactly the sums of masses
and first moments.  Keeping masses fixed, all the target barycentres
can therefore be enlarged by a factor $1+\eta_n>1$.  Equivalently,
$S\cx(1-\delta_n)Y_0$ for some $\delta_n>0$.

The perturbed reference law dominates
$c_n\law(Y_0)$ as a measure, where $c_n=(1-\varepsilon)^n$.
The residual measure is centred, so $c_nY_0\cx Y_\varepsilon$.
For sufficiently small $\varepsilon$, $c_n\ge1-\delta_n$, proving the
fixed-block assertion.  The common Gaussian case follows by dominating
the fair $\{-2,2\}$ variable with $N(0,9)$.

For the quantitative conclusions, put $M=n^{-1}\sum_iS_i$ and
$d_\nu(m)=(1-m^2)/2-L_\nu((1+m)/2)$.  In a martingale coupling, multiply
$\E[U+V_i\mid S]=S_i$ by $S_i-M$ and sum.  The common variable cancels
pathwise, and rearrangement of the $V_i$ gives
\[
 \E(1-M^2)\le2\E L_\nu((1+M)/2)
                 +4\E W_1(\widehat\nu_n,\nu).
\]
Integrated quantiles are uniformly $1$-Lipschitz in $W_1$ and
$|\overline V|\le W_1(\widehat\nu_n,\nu)$; these are the two error terms
in this bound.  For a law supported on $[-2,2]$,
\[
 \E W_1(\widehat\nu_n,\nu)
 \le n^{-1/2}\int_{-2}^2\sqrt{F(t)(1-F(t))}\,dt\le2n^{-1/2}.
\]
Thus $\E d_\nu(M)\le4/\sqrt n$.
In the change of variables in \eqref{PHASE-eq:defect},
\[
 1-m^2=(1-u^2)\left[1-\varepsilon u^2
             -\frac{\varepsilon^2u^2(1-u^2)}4\right].
\]
The bracket lies between zero and $1-\varepsilon u^2\le1$.
Thus $d_{\nu_\varepsilon}(m)\ge\varepsilon(1-m^2)^2/8$, whence
\begin{equation}\label{PHASE-eq:finite-deficit}
 \E(1-M^2)^2\le\frac{32}{\varepsilon\sqrt n}.
\end{equation}

A sharper estimate uses the small empirical minority fraction. Let
$D_n=\|\widehat F_n-F\|_\infty$ and
$a_\varepsilon=(1-\varepsilon)/4$, a lower bound for the density on
its support. The inverse distribution function is
$a_\varepsilon^{-1}$-Lipschitz, so the empirical and true quantiles
are uniformly within $D_n/a_\varepsilon$. For $p\le1/2$, integrate
this bound over the upper $p$-tail and subtract $p\overline V$.
For $p\ge1/2$, use the complementary lower tail. This gives
\[
 |L_{\widehat\nu_n}(p)-p\overline V-L_\nu(p)|
 \le\frac{2\min(p,1-p)}{a_\varepsilon}D_n
 \le\frac{4p(1-p)}{a_\varepsilon}D_n.
\]
The same rearrangement inequality, with $u=1-M^2$, now yields
\[
 \frac\varepsilon8\E u^2
 \le a_\varepsilon^{-1}\E(uD_n).
\]
The Dvoretzky--Kiefer--Wolfowitz inequality with Massart's constant
\cite{Massart} gives $\E D_n^2\le1/n$. Cauchy--Schwarz consequently
proves the additional bound
\begin{equation}\label{PHASE-eq:sharp-finite-deficit}
 \sqrt{\E(1-M^2)^2}
 \le\frac{32}{\varepsilon(1-\varepsilon)\sqrt n}.
\end{equation}
This estimate also gives an explicit $O_\varepsilon(n^{-1/2})$ bound
on the deficit of the average off-diagonal covariance from one.

For every sign law the count-entropy bound is
$H(S)\le\log(n+1)+n\E h((1+M)/2)$.
The function $\psi$ in \eqref{PHASE-eq:improved-entropy} is increasing
and concave: for $r=\sqrt{1-u}$,
\[
 \psi'(u)=\frac{\operatorname{arctanh}r}{2r},\qquad
 \psi''(u)=-\frac{r/(1-r^2)-\operatorname{arctanh}r}{4r^3}<0.
\]
Jensen and Cauchy--Schwarz therefore give
\[
 \E h((1+M)/2)=\E\psi(1-M^2)\le
 \psi\!\left(\min\{1,\sqrt{\E(1-M^2)^2}\}\right).
\]
Combining \eqref{PHASE-eq:finite-deficit} and
\eqref{PHASE-eq:sharp-finite-deficit} proves the displayed entropy
estimate. The average off-diagonal
covariance is $(n\E M^2-1)/(n-1)$, so \eqref{PHASE-eq:finite-deficit}
gives its convergence to one.  For iid fair signs,
$\E(1-M^2)^2=1-2/n+3/n^2-2/n^3\ge1/2$ for $n\ge2$;
substitution proves the dimension bound.  The two iterated entropy
limits now follow from the fixed-block assertion and the finite estimate.
\end{proof}

\begin{corollary}[A coupling to one repeated bit]\label{cor:one-bit}
Under the quantitative hypotheses of Theorem~\ref{PHASE-thm:collapse},
there is a fair sign $R$ on an extension of the same probability space
such that
\[
 \E d_H(S,R\mathbf1)=\frac n2\E(1-M_n^2)
 \le\frac n2\delta_{n,\varepsilon}.
\]
A uniformly sampled $k$-coordinate block has total-variation distance
at most $k\delta_{n,\varepsilon}/2$ from a repeated fair bit.
For an exchangeable law the bound applies to any specified block.
\end{corollary}
\begin{proof}
Choose $J$ independently and uniformly in $[n]$ and set $R=S_J$.
Given $S$, the mean number of disagreements is
$2n\,[\#\{S_i=1\}/n][\#\{S_i=-1\}/n]=n(1-M_n^2)/2$.
A union bound for the sampled block proves the second assertion.
\end{proof}


\section{Local dynamics under global reference guarantees}
\label{R11-sec:local}
High entropy, near-Gaussian comparison and many independent marginals can coexist with exponentially slow local movement. Code separation gives the obstruction for stationary processes and hidden-state lifts with the specified visible jump bound. The Gaussian and cosine--auxiliary constructions retain their different jump hypotheses; the latter also has a fast global sampler.

\subsection{Near-critical reference laws with slow linear-block dynamics}\label{SP-sec:codes}
Let $\mathcal C\le\mathbb F_2^n$ be a linear code of dimension $k$, size $M=2^k$, minimum nonzero distance $d$, and dual minimum distance $d^\perp$. Identify $x\in\mathbb F_2^n$ with the sign vector $((-1)^{x_i})_i$. The uniform laws on the code and on the full cube are $U_{\mathcal C}$ and $U_n$. Define
\[
 D_0=\KL(U_{\mathcal C}\Vert U_n)=(n-k)\log2,
 \qquad \mu_\eta=(1-\eta)U_{\mathcal C}+\eta U_n,
 \quad0<\eta<1.
\]
The classical dual-distance criterion identifies $U_{\mathcal C}$ as an orthogonal array of strength $d^\perp-1$. A direct proof is included below, so the dynamical conclusions depend only on linear algebra and Corollary~\ref{R9-cor:privacy-code}.

\begin{proposition}[What the full-support mixture preserves]\label{SP-prop:code-mixture}
Suppose $U_{\mathcal C}\cx sG_n$ with $s\ge\sqrt{\pi/2}$. Then
\[
 \mu_\eta\cx sG_n,
 \qquad
 (1-\eta)D_0-h_2(\eta)\le\KL(\mu_\eta\Vert U_n)\le(1-\eta)D_0.
\]
The law $\mu_\eta$ has full support and exactly uniform marginals on every set of fewer than $d^\perp$ coordinates. Its total-variation distance from $U_n$ is
\[
 \operatorname{TV}(\mu_\eta,U_n)=(1-\eta)(1-M/2^n).
\]
\end{proposition}
\begin{proof}
Conditioning $\sqrt{\pi/2}G_n$ on its coordinate signs gives $U_n$, so both mixture components are dominated by $sG_n$. Convexity of expectations proves the mixture comparison. Introduce the mixture indicator $J$. Relative-entropy decomposition gives
\[
 (1-\eta)D_0=\KL(\mu_\eta\Vert U_n)+I(J;X),\qquad0\le I(J;X)\le h_2(\eta).
\]
For a coordinate subset $S$, projection of $\mathcal C$ onto $\mathbb F_2^S$ is surjective exactly when no nonzero dual vector is supported on $S$. Uniform measure projects uniformly onto the image, proving the marginal assertion. Both components then have the same such marginals. The remaining assertions follow from the explicit two values of the mixture density.
\end{proof}

\begin{lemma}[Stationary exit with arbitrary memory]\label{R10-lem:stationary-flow}
Let $(X_t)_{t\ge0}$ be any discrete-time process with $\law(X_t)=\mu$ at
every time and Hamming jumps at most $r$. For $\mu(c)>0$, set
$N_r(c)=\{x:1\le d_H(x,c)\le r\}$ and
$\tau_c=\inf\{t\ge1:X_t\ne c\}$. Then
\begin{equation}\label{R15-eq:stationary-exit}
 \Prb(\tau_c\le T\mid X_0=c)\le
 T\frac{\mu(N_r(c))}{\mu(c)},\qquad
 \TV(\law(X_T\mid X_0=c),\mu)
 \ge1-T\frac{\mu(N_r(c))}{\mu(c)}-\mu(c).
\end{equation}
For a code $\mathcal C$ of distance greater than $r$ and
$\mu=\mu_\eta$, the ratio is at most
\begin{equation}\label{R10-eq:stationary-flow}
 q_r=\frac{\eta}{1-\eta}e^{-D_0}\sum_{j=1}^r\binom nj.
\end{equation}
Thus every stationary Markov kernel with such jumps has the former
mixing-time bound $t_{\rm mix}(1/4)\ge(4q_r)^{-1}$ when
$\mu(c)\le1/2$. The exit conclusion applies equally to stationary hidden-state
lifts, initialized conditionally on their visible state $c$.
\end{lemma}
\begin{proof}
On $\{X_0=c,\tau_c=j\}$, $X_j\in N_r(c)$. A union bound and one-time
stationarity give
$\Prb(X_0=c,\tau_c\le T)\le T\mu(N_r(c))$.
The atom $c$ tests total variation. No Markov assumption enters.
Since $N_r(c)$ contains no codeword,
$\mu(N_r(c))=\eta2^{-n}\sum_{j=1}^r\binom nj$ and
$\mu(c)\ge(1-\eta)/|\mathcal C|$, giving the ratio. A stationary lift has
these visible marginals. Its conditional initialization is a mixture of
hidden initial states, so at least one such state has no larger exit
probability up to any specified horizon.
\end{proof}

\subsubsection{Conditional influence eigenvalues}
For a sign law with nondegenerate individual marginals, define its off-diagonal influence matrix by
\[
 \Psi_{ij}=\Pp(X_j=1\mid X_i=1)-\Pp(X_j=1\mid X_i=-1)
 \quad(i\ne j),\qquad\Psi_{ii}=0.
\]
Spectral independence controls this matrix after every feasible coordinate conditioning, using the standard symmetrization when marginal variances differ \cite{AnariLiuOveis}. The conditional law used here is symmetric with fair marginals, so its influence matrix is already symmetric.

\begin{proposition}[A linear conditional eigenvalue with full support]\label{SP-prop:influence}
For the full-support mixture above there is a conditioning leaving $d$ coordinates free whose influence matrix has largest eigenvalue at least
\[
 (1-\eta)(d-1).
\]
The unperturbed code gives exactly $d-1$.
\end{proposition}
\begin{proof}
Let $S$ be the support of a minimum-weight nonzero codeword, and condition all coordinates outside $S$ to be zero in binary notation. The shortened code on $S$ consists of zero and that codeword. Any other nonzero word supported there would have smaller weight. Thus the conditional code law is the fair mixture of $\mathbf1_S$ and $-\mathbf1_S$ in sign notation.

Let $s=|S|=d$. The conditioning event has probability $2/M$ under the code and $2^{-(n-s)}$ under the cube. Its code probability is at least its cube probability: projection onto the complement has kernel dimension one and rank $k-1\le n-s$, so $s\le n-k+1$. Therefore the conditional mixture weight $w$ of the code component satisfies $w\ge1-\eta$.

The conditional full-support law is
\[
 w\,\tfrac12(\delta_{\mathbf1_S}+\delta_{-\mathbf1_S})+(1-w)U_S.
\]
It has covariance $w\mathbf1\mathbf1^{\mathsf T}+(1-w)I$. Its influence matrix is $w(\mathbf1\mathbf1^{\mathsf T}-I)$, whose largest eigenvalue is $w(s-1)$.
\end{proof}

\subsubsection{Heat-bath updates of a linear number of coordinates}
A block heat-bath update chooses a coordinate set $B$ and resamples its coordinates from $\mu_\eta$ conditional on the coordinates outside $B$. A fixed mixture of such updates is reversible with stationary law $\mu_\eta$. The pointwise escape estimate below also holds for any adaptive choice of blocks, provided each chosen block is resampled using this conditional law.

\begin{theorem}[Uniform escape bound for local block updates]\label{SP-thm:block-escape}
Suppose every allowed block has size at most $b<d$. At any codeword, the one-step probability of leaving that codeword is at most
\begin{equation}\label{SP-eq:escape}
 q_*:=\frac{\eta}{1-\eta}\exp(b\log2-D_0).
\end{equation}
For any fixed mixture of these block heat-bath kernels, if $q_*\le1/4$ and the stationary mass of one codeword is at most $1/4$, then
\[
 t_{\rm mix}(1/4)\ge\frac1{4q_*}
\]
up to the immaterial integer rounding of time. Its spectral gap is at most
\[
 \frac{q_*}{1-\mu_\eta(c)}
\]
for any codeword $c$ with $\mu_\eta(c)<1$.
\end{theorem}
\begin{proof}
If $|B|<d$, the only codeword agreeing with $c$ outside $B$ is $c$ itself. The conditional fibre has $2^{|B|}$ points. Its distinguished point has mass $(1-\eta)/M+\eta2^{-n}$, while each other point has mass $\eta2^{-n}$. Hence the escape probability after this block update is
\[
 \frac{\eta(2^{|B|}-1)2^{-n}}
      {(1-\eta)/M+\eta2^{|B|-n}}
 \le\frac{\eta}{1-\eta}(2^b-1)M2^{-n}\le q_*.
\]
This is uniform in the chosen block, so mixtures and adaptive selections satisfy the same bound as long as the chain remains at $c$. Starting there, a union bound gives probability at least $3/4$ of still being there at every time $t\le1/(4q_*)$. Its stationary mass is at most $1/4$, so total variation is at least $1/2$ at those times. For the spectral gap, use the indicator of $\{c\}$ in the reversible Rayleigh quotient: its Dirichlet form is $\mu_\eta(c)P(c,c^c)$ and its variance is $\mu_\eta(c)(1-\mu_\eta(c))$.
\end{proof}

\begin{theorem}[Entropy-extremizing laws with slow linear-block dynamics]\label{SP-thm:nearcritical-slow}
Fix $\varepsilon>0$, put $\kappa=\sqrt{\pi/2}$ and let
$\Delta=\Delta_G(\varepsilon)\in(0,\log2)$ be the exact Gaussian entropy frontier in Corollary~\ref{R7-cor:Gaussian-entropy}. Choose positive constants $\delta,\delta_0$ such that
\[
 h_2(\delta)<\log2-\Delta,
 \qquad h_2(\delta_0)<\Delta,
 \qquad \delta,\delta_0<1/2.
\]
There are a sequence of dimensions tending to infinity and full-support laws $\mu_n$ for which
\begin{gather*}
 \mu_n\cx(1+\varepsilon)\kappa G_n,
 \qquad n^{-1}\KL(\mu_n\Vert U_n)\longrightarrow\Delta,\\
 \text{every marginal on at most $\lfloor\delta n\rfloor$ coordinates is uniform},
 \qquad \operatorname{TV}(\mu_n,U_n)\longrightarrow1.
\end{gather*}
Their spectral-independence parameters are at least $(\delta_0+o(1))n$. For every fixed $0<\beta<\delta_0$, every fixed mixture of heat-bath blocks of size at most $\lfloor\beta n\rfloor$ has
\begin{equation}\label{SP-eq:exponential-mixing}
 t_{\rm mix}(1/4)\ge
 \exp\{n[\Delta-\beta\log2-o(1)]\}.
\end{equation}
The exponent is positive for these choices. The spectral gap is at most the reciprocal exponential rate.
For the same laws, every stationary process with visible jumps at most $\lfloor\beta n\rfloor$, including every stationary hidden-state lift, has conditional visible total variation bounded away from zero up to times of order
\begin{equation}\label{R10-eq:all-stationary}
 T_n=\exp\{n(\Delta-h_2(\beta)-o(1))\}.
\end{equation}
For stationary Markov kernels this gives the same worst-case mixing lower
bound. For a lift it holds from the stationary conditional hidden law given
a codeword, and from some individual hidden state for each specified horizon.

This exponent is positive for every $0<\beta<\delta_0$.
\end{theorem}
\begin{proof}
Use the uniform codes constructed in Corollary~\ref{R9-cor:privacy-code}, with
$d_n>\delta_0n$, $d_n^\perp>\delta n$ and $D_{0,n}/n\to\Delta$. Put
$\mu_n=(1-1/n)U_{\mathcal C_n}+n^{-1}U_n$. Proposition~\ref{SP-prop:code-mixture} gives the comparison, exact marginals, entropy limit and total variation. Proposition~\ref{SP-prop:influence} gives the conditional influence eigenvalue.

For $b=\lfloor\beta n\rfloor<d_n$, the escape bound is at most
$(n-1)^{-1}\exp\{-D_{0,n}+b\log2\}$. Individual codeword masses tend to zero because the code rate $\log2-\Delta$ is positive. Theorem~\ref{SP-thm:block-escape} therefore gives \eqref{SP-eq:exponential-mixing}. Finally,
$h_2(\delta_0)>\delta_0\log2$ for $0<\delta_0<1/2$, so
$\Delta>\beta\log2$.
For arbitrary stationary kernels, Lemma~\ref{R10-lem:stationary-flow} and the binomial estimate $\sum_{j\le\beta n}\binom nj\le\exp\{nh_2(\beta)+o(n)\}$ give~\eqref{R10-eq:all-stationary}. Since $h_2$ is strictly increasing on $(0,1/2)$, $h_2(\beta)<h_2(\delta_0)<\Delta$.
\end{proof}

As $\varepsilon\downarrow0$, the allowed independent-marginal fraction $\delta$ approaches $1/2$. Theorem~\ref{SP-thm:nearcritical-slow} keeps each positive $\varepsilon$ fixed while the dimension grows. It makes no claim of computational hardness for arbitrary samplers. The conclusion concerns conditional influence and every stationary bounded-jump dynamics of the stated linear size; the heat-bath subclass has the sharper exponent in~\eqref{SP-eq:exponential-mixing}. A sampling theorem based on a different mechanism is compatible with these examples.

\subsubsection{Slow local dynamics with the original independent references}
The preceding theorem concerns near-endpoint Gaussian references and the optimal entropy deficit. The following construction instead retains the actual independent cosine and optimized coefficient sources of the hard signing theorem, together with a positive conditional allocation kernel. These reference guarantees coexist with a different code geometry.

Let $q_\dagger$ be the density in~\eqref{R-eq:htrial}, rescaled to $(-3,3)$, let $T_1,\ldots,T_n$ be independent with this density, and let $X_R$ be the independent product cosine source on $(-R,R)^m$. Put
\[
 v_R=R^2\left(\frac13-\frac2{\pi^2}\right).
\]
\begin{theorem}[Exact references and a local stationary obstruction]\label{R10-thm:cosine-code}
Fix $R>0$. For $b\ge1$, let $n=7b$, let $\mathcal C_b$ be the $b$-fold product of the binary $[7,4,3]$ Hamming code, and put
\[
 \mu_b=\frac9{10}U_{\mathcal C_b}+\frac1{10}U_n.
\]
For every $A\in\R^{m\times n}$ with $\sum_{i,j}|A_{ij}|\le v_R/(42R)$ there is a coupling with
\[
 \law(\sigma)=\mu_b,\qquad X_R\perp T,\qquad
 \E[X_R\mid\sigma]=A\sigma,\quad \E[T\mid\sigma]=\sigma,
\]
and $\Pp(\sigma=s\mid X_R,T)\ge\mu_b(s)/42$ almost surely. In particular, $A$ can have every column nonzero. The law is symmetric, has full support and three-wise independent coordinates, and
\[
 \Cov(\sigma)=I_n,\qquad H(\sigma)\ge4.3b\log2.
\]
Every stationary process with one- or two-coordinate visible jumps has the corresponding conditional exit obstruction for $2^{3b}/\operatorname{poly}(b)$ steps. In particular stationary chains, including arbitrary hidden-state lifts, have this worst-case visible mixing obstruction. The random-scan single-coordinate heat-bath chain has exact exit probability $1/(9\,8^b+2)$ from each codeword. A coordinate conditioning gives an influence eigenvalue of order $n$, whereas independent global samples require expected $O(n)$ fair-bit operations.
\end{theorem}
\begin{proof}
Write $a=\E|T_1|$. The explicit rational polynomial $h$ in~\eqref{R-eq:htrial} satisfies $h(x)\le0$ for $0\le x\le1$: group each positive even-power term with a larger preceding negative coefficient. Thus $1+h\le e^h\le1+h+h^2/2$. Exact integration gives
\[
 a\ge\frac{3\int_0^1x(1+h(x))\,dx}
                {\int_0^1(1+h(x)+h(x)^2/2)\,dx}
 =\frac{156689961407407355794760643750000000000000000}
 {110986716051424670746467932873184051716949883}
 >\frac75.
\]
In each seven-coordinate block, correct the syndrome of $\operatorname{sgn}T$ by flipping its indexed coordinate when the syndrome is nonzero. Each codeword has exactly eight preimages, with errors $0,e_1,\ldots,e_7$. The output is uniform on $\mathcal C_b$, and its conditional auxiliary mean is $(3a/4)\sigma$. Keep that output with probability $4/(3a)$ and otherwise reset to an independent uniform codeword. The resulting regression is $\sigma$, and its kernel has floor $(1-4/(3a))U_{\mathcal C_b}\ge U_{\mathcal C_b}/21$. Keeping $\operatorname{sgn}T$ with probability $1/a$ and otherwise resetting to $U_n$ gives the cube law with the same regression and floor $(1-1/a)U_n\ge2U_n/7$. Their mixture therefore has a kernel $k_s^0(t)\ge\mu_b(s)/21$.

Independently sample $X_R$ and set
\[
 k_s(x,t)=k_s^0(t)+\mu_b(s)v_R^{-1}\langle x,As\rangle.
\]
Since $\E_{\mu_b}\sigma=0$, the kernels sum to one. The added term has absolute value at most $\mu_b(s)/42$ by the hypothesis on $A$. Centering of $X_R$ leaves label masses and auxiliary moments unchanged, while $\E X_RX_R^{\mathsf T}=v_RI_m$ gives the physical conditional mean exactly. Both prescribed sources, including their independence, retain their entire joint marginal.

The Hamming dual has minimum weight four, so every nonconstant Walsh coefficient of degree at most three vanishes; this passes to products and mixtures. The all-ones binary word belongs to the Hamming code, giving sign symmetry. Entropy concavity gives $H(\mu_b)\ge(9\cdot4+7)b\log2/10$.

For $r\le2$, every codeword has no other codeword within distance $r$. The stationary-flow calculation gives the sharper bound
\[
 P(c,c^c)\le
 \frac{\sum_{j=1}^r\binom{7b}{j}}{9\,2^{3b}+1}.
\]
Lemma~\ref{R10-lem:stationary-flow} then gives the mixing obstruction without reversibility. In a single-coordinate heat-bath fibre, the other point has mass $2^{-7b}/10$ and the codeword has mass $(9\,2^{-4b}+2^{-7b})/10$. Their ratio yields the asserted exit probability.

Pin the four information coordinates in each block to $+1$. Both mixture components give this event probability $2^{-4b}$, so the remaining $L=3b$ signs have law $(9/10)\delta_{\mathbf1}+(1/10)U_L$. Their covariance is $(1/10)I+(9/100)\mathbf1\mathbf1^{\mathsf T}$. All variances are $19/100$, so the largest eigenvalue of the off-diagonal influence matrix is $9(L-1)/19$. This is a conditioned-law statement, as required by spectral independence. A global sampler chooses the mixture component and then samples either $4b$ information bits followed by block encoding, or $7b$ independent bits, using $O(n)$ operations after the mixture choice. That Bernoulli choice takes a constant expected number of fair bits.
\end{proof}
The slack in Theorem~\ref{R10-thm:cosine-code} permits nonzero physical columns. Its global sampler is efficient, while the stationary local kernels have the stated exponential obstruction. The exact Gaussian comparisons follow from the same joint cosine--$q_\dagger$ coupling. They differ from the near-endpoint Gaussian scale of Theorem~\ref{SP-thm:nearcritical-slow}.

\subsubsection{Affine-code representations of hard-balanced laws}
The entropy-optimal codes in the product-reference problem and the
hard-balanced laws have different representation requirements.
\begin{proposition}[Number of affine-code components]
\label{R9-prop:affine-mixtures}
Suppose $\mu\cx(1+\varepsilon)\sqrt{\pi/2}\,G_n$ is a mixture of
$M$ uniform affine binary code laws, each supported on
$\{s:|\sum_i s_i|\le R\}$. Then
\[
 \log M\ge n\left[\tfrac12\log2-\Delta_G(\varepsilon)\right]
                          -\tfrac R2\log2.
\]
For sufficiently small fixed $\varepsilon$ and $R=O_\varepsilon(1)$,
this requires exponentially many components.
\end{proposition}
\begin{proof}
A dimension-$k$ affine code projects bijectively onto a full cube on
$k$ information coordinates. One codeword has all these signs positive,
so its total sum is at least $2k-n$. Hence $k\le(n+R)/2$.
If $J$ indexes the mixture,
\[
 H(\mu)\le H(J)+H(\sigma\mid J)
              \le\log M+\tfrac{n+R}{2}\log2.
\]
The Gaussian entropy envelope gives
$H(\mu)\ge n[\log2-\Delta_G(\varepsilon)]$. Combining the two
bounds proves the claim.
\end{proof}
The proposition concerns the size of this representation. An implicit
sampler can draw from an exponentially large mixture without listing
its components.


\Needspace{12\baselineskip}
\part{Calibration, privacy and comparison of experiments}
\label{R11-part:experiments}
We retain an entire observation--target law while adding an independent source. Its feasible conditional means lead to information costs, exact privacy and Blackwell comparison. Directed transport costs then reconstruct the posterior law itself.

\section{Extending an experiment and measuring its information cost}
\label{R11-sec:extensions}
Keeping the old observation--target law fixed gives an integrated-quantile constraint for each target-label subset. Conditional allocation attains every greedy vertex, proving completeness. For a fair added source, a nonlinear potential gives the minimum information cost; its Hessian identifies feasible directions. A bounded perturbation of the original kernel gives the vector extension and an exact sampler.

\subsection{Global extension of a retained coupling}\label{NLI-sec:global}
Let $Y\in\{1,\ldots,N\}$ and let $Z$ be standard Borel. Fix the entire law of $(Z,Y)$ and write
\[
 k_i(z)=\Pp(Y=i\mid Z=z),\qquad p_i=\E k_i(Z)>0.
\]
We prescribe a centered scalar source $X$, independent of $Z$, with $\E|X|<\infty$. The unknown coupling may correlate $X$ with $Y$. We seek its joint means
\[
 m_i=\E[X\ind_{\{Y=i\}}],
\]
so the conditional target means are $m_i/p_i$. Throughout this section the original $(Z,Y)$ law and the independent $(X,Z)$ law are retained.

\subsubsection{A submodular base polytope is the complete feasible region}
Write $\rho=\law(X)$. The upper integrated quantile is the same
function as in \eqref{S-eq:radius} and \eqref{PHASE-eq:allowed}:
\[
 L_\rho(s)=\int_{1-s}^1Q_X(v)\,dv,\qquad0\le s\le1.
\]
The function is concave and continuous, with $L_\rho(0)=L_\rho(1)=0$. For $A\subseteq[N]$ set
\begin{equation}\label{NLI-eq:submodular}
 k(A)=\sum_{i\in A}k_i,
 \qquad F(A)=\E L_\rho(k(A)).
\end{equation}
Concavity of $L_\rho$ implies that $F$ is submodular: $F(A)+F(B)\ge F(A\cup B)+F(A\cap B)$. Define its base polytope by
\[
 B(F)=\{m\in\R^N:m([N])=0,\ m(A)\le F(A)\text{ for every }A\subseteq[N]\}.
\]

\begin{theorem}\label{NLI-thm:global}
The complete set of feasible joint mean vectors is $B(F)$. Every feasible vector is attained by a mixture of at most $N$ conditional quantile couplings. The mixing variable can be chosen independently of $(X,Z)$.
\end{theorem}
\begin{proof}
For a fixed $z$, any selection event of conditional mass $s$ has $X$-integral at most $L_\rho(s)$ by the quantile rearrangement inequality. Since $X$ has its unchanged law conditional on $Z$, this applied to $Y\in A$ gives $m(A)\le F(A)$. Centering gives $m([N])=0$.

For sufficiency choose an order $\pi=(\pi_1,\ldots,\pi_N)$ and the nested sets $A_j=\{\pi_1,\ldots,\pi_j\}$. Given $Z=z$, draw a uniform $V\in(0,1)$ independently, put $X=Q_X(V)$, and allocate the highest $k_{\pi_1}(z)$ fraction of quantiles to label $\pi_1$, the next fraction to $\pi_2$, and so on. This is a measurable coupling and it preserves both required pair laws. Its mean vector is $m^\pi_{\pi_j}=F(A_j)-F(A_{j-1})$. Atoms of $X$ cause no difficulty because the uniform quantile variable randomizes within their quantile intervals.

The greedy theorem for submodular base polytopes says that these vectors have the support function of $B(F)$. Here it also follows directly: for $u_{\pi_1}\ge\cdots\ge u_{\pi_N}$, summation by parts gives
\[
 u\cdot m\le\sum_{j=1}^{N-1}(u_{\pi_j}-u_{\pi_{j+1}})F(A_j),
\]
and the displayed quantile coupling attains equality. Thus their convex hull is exactly $B(F)$. Carath\'eodory's theorem in the hyperplane $m([N])=0$ represents each vector using at most $N$ of them. Taking that mixture independently of $(X,Z)$ retains the two pair laws.
\end{proof}
The greedy base-polytope mechanism is classical \cite{NLI-Bach}. Here its rank function is the expected conditional quantile capacity \eqref{NLI-eq:submodular}, and its greedy allocations preserve both prescribed pair laws. No boundedness, density or atomlessness assumption is imposed on $X$.

For a fair sign $X\in\{-1,1\}$, the formula simplifies to
\begin{equation}\label{NLI-eq:fair-region}
 F(A)=\E\min\{k(A),1-k(A)\}.
\end{equation}
The formula determines every finite amplitude. The range of the Hessian below describes its feasible directions at the origin.

\subsubsection{The full binary potential contains the local information metric}
Write
\[
 j(t)=\frac{1+t}{2}\log(1+t)+\frac{1-t}{2}\log(1-t),\qquad -1\le t\le1,
\]
with $0\log0=0$. For fair $X$, set $t_i(z)=\E[X\mid Z=z,Y=i]$. The constraints become
\[
 |t_i|\le1,\qquad \sum_i k_it_i=0,\qquad m_i=\E k_it_i,
\]
and the added information is $I(X;Y\mid Z)=\E\sum_i k_ij(t_i)$.

\begin{theorem}\label{NLI-thm:potential}
Let $J_k(m)$ denote the minimum of this information, with value $+\infty$ outside the feasible region. Then
\begin{equation}\label{NLI-eq:potential}
 J_k=\Phi_k^*,\qquad
 \Phi_k(u)=\E\inf_{\alpha\in\R}\sum_i k_i(Z)\log\cosh(u_i-\alpha).
\end{equation}
Moreover
\begin{equation}\label{NLI-eq:hessian}
 \nabla^2\Phi_k(0)=G:=D_p-\E[kk^{\mathsf T}],\qquad D_p=\diag(p).
\end{equation}
On the range of $G$, the minimum cost has quadratic term $\tfrac12m^{\mathsf T}G^\dagger m$.
\end{theorem}
\begin{proof}
For a fixed posterior vector $k$, minimize the Lagrange multiplier for the constraint $\sum_i k_it_i=0$. Since $j^*(v)=\log\cosh v$, the conditional conjugate is the integrand of \eqref{NLI-eq:potential}. At finite $u$, the minimizing multiplier solves
\[
 \sum_i k_i\tanh(u_i-\alpha)=0
\]
and lies between the minimum and maximum of the active $u_i$. It is measurable in $k$. Pointwise optimization therefore commutes with integration. The feasible family of bounded kernels is weak-star compact and the finite moment constraints are continuous. Relative entropy is lower semicontinuous, so the minimum cost is a closed convex function. Taking its conjugate proves \eqref{NLI-eq:potential} and its reverse by biconjugacy.

At zero, the first-order multiplier is $\alpha=\sum_i k_iu_i$. The conditional Hessian is $D_k-kk^{\mathsf T}$. Integration gives \eqref{NLI-eq:hessian}. Its kernel consists of vectors constant on each connected component of the graph $\E k_ik_j>0$. Restricting to its orthogonal complement gives a positive-definite Hessian, so local Legendre inversion gives the cost expansion.
\end{proof}
The inverse conditional covariance is therefore the quadratic part of an explicit global convex potential. The exact feasible region of Theorem~\ref{NLI-thm:global} is the domain of its conjugate.

\subsubsection{Vector sources and the cost of exact augmentation}
The preceding Hessian also controls bounded vector sources. Keep the
whole experiment $(Z,Y)$, and prescribe $X\perp Z$ centered and bounded,
with covariance $\Sigma>0$. Let the desired means be $\tau b_i$, and
write $B$ for the matrix with rows $b_i^{\mathsf T}$.
The matrix $G$ is the one in \eqref{NLI-eq:hessian}; $G^\dagger$ denotes
its Moore--Penrose inverse. Every computational statement below is
conditional on access to the indicated baseline and source samplers.
\begin{theorem}[Exact augmentation of a fixed joint law]
\label{R10-thm:augmentation}
There is a neighborhood of $\tau=0$ on which couplings with
\begin{equation}\label{R10-eq:three-marginals}
 \law(X,Z)=\law(X)\otimes\law(Z),\quad
 \law(Z,Y)=\text{the prescribed law},\quad
 \E[X\mid Y=i]=\tau b_i
\end{equation}
exist if and only if every column of $D_pB$ belongs to $\ran G$.
In that case put
\[
 H=G^\dagger D_pB,\quad h_i=H_{i,\cdot}^{\mathsf T},\quad
 \bar h(z)=\sum_i k_i(z)h_i,\quad
 r_i(x,z)=\langle\Sigma^{-1}x,h_i-\bar h(z)\rangle.
\]
For any $\tau$ with $\eta=|\tau|\esssup_{x,z,i:k_i(z)>0}|r_i(x,z)|<1$,
the explicit kernel
\begin{equation}\label{R10-eq:kernel}
 K_\tau(i\mid x,z)=k_i(z)[1+\tau r_i(x,z)]
\end{equation}
satisfies \eqref{R10-eq:three-marginals}. It also has the exact posterior
\begin{equation}\label{R10-eq:physical-posterior}
 \frac{dP_{X\mid Y=i}}{dP_X}(x)
 =1+\tau\langle\Sigma^{-1}b_i,x\rangle.
\end{equation}
Given a sampler for $(Z,Y)$, a sampler for $X$, and evaluations of
$h_Y$ and $\bar h(Z)$, rejection with acceptance probability
$(1+\tau r_Y)/(1+\eta)$ samples the full coupling using exactly
$1+\eta$ expected proposals. Conditional likelihoods are not inputs.
\end{theorem}
\begin{proof}
For any scalar vector $u$,
\[
 u^{\mathsf T}Gu=\E\Var(u_Y\mid Z).
\]
Thus $G$ is a graph Laplacian, with edge weights $\E[k_i k_j]$.
Its kernel consists of the vectors for which $u_Y$ is measurable with
respect to $Z$. For any coupling satisfying \eqref{R10-eq:three-marginals}
and any such $u$, independence gives
$\tau\sum_i p_i u_i b_i=\E[Xu_Y]=0$. For nonzero $\tau$ this is precisely
the stated range condition.

Conversely $GH=D_pB$. Pointwise, $\sum_i k_i r_i=0$, so
\eqref{R10-eq:kernel} is a probability kernel. Since $\E X=0$,
$\E_XK_\tau(i\mid X,z)=k_i(z)$, which retains the whole $(Z,Y)$ law.
The source marginal $P_X\otimes P_Z$ is unchanged by construction.
Moreover,
\[
 \E[XK_\tau(i\mid X,Z)]
 =\tau\E_Z k_i(Z)(h_i-\bar h(Z))
 =\tau(GH)_i=\tau p_i b_i.
\]
Integrating the kernel over $Z$ gives \eqref{R10-eq:physical-posterior}.
The rejection distribution has density $1+\tau r_Y$ relative to
$P_X\otimes P_{ZY}$, and its total mass is one. This proves the cost
and correctness assertions. Boundedness of $X$ and finiteness of $H$
give a nonzero interval of admissible $\tau$.
\end{proof}

\begin{remark}[The compatibility condition has an exact meaning]
The connected components of the graph $\E[k_i k_j]>0$ are labels
recoverable from $Z$. The condition is $\sum_{i\in C}p_i b_i=0$ on each
component $C$. A deterministic baseline kernel has $G=0$ and permits
only zero physical means while its entire $(Z,Y)$ law is retained.
A positive baseline kernel has one component, so only $\sum_i p_i b_i=0$
is required. Singular physical covariances are handled by restricting to
$\ran\Sigma$ and using $\Sigma^\dagger$ there.
\end{remark}

\begin{theorem}[Sharp information cost with a fourth-order certificate]
\label{R10-thm:energy}
Assume in addition that $X$ is centrally symmetric. Define
\begin{equation}\label{R10-eq:energy}
 Q=\frac12\tr(\Sigma^{-1}B^{\mathsf T}D_pG^\dagger D_pB),\qquad
 M_4=\E_{P_X\otimes P_{ZY}}r_Y^4.
\end{equation}
Let $J_{\min}(\tau)$ be the minimum of $I(X;Y\mid Z)$ over
\eqref{R10-eq:three-marginals}. The explicit kernel above satisfies
\begin{align}
 \tau^2Q-\frac{e^\eta\tau^4M_4}{24}
 &\le J_{\min}(\tau)\le I_{K_\tau}(X;Y\mid Z)
 \label{R10-eq:cost-sandwich}\\
 &\le\tau^2Q+\frac{\tau^4M_4}{12(1-\eta^2)}.\notag
\end{align}
Consequently $J_{\min}(\tau)=\tau^2Q+O(\tau^4)$, and the explicit
rejection sampler is within
$\tau^4M_4[e^\eta/24+1/(12(1-\eta^2))]$ of the minimum.
\end{theorem}
\begin{proof}
Write $P_0=P_X\otimes P_{ZY}$. Direct calculation gives
$\E_0r_Y^2=\tr(\Sigma^{-1}H^{\mathsf T}GH)=2Q$.
For every feasible coupling $P$, the constraints imply
\[
 \E_P r_Y=\tau\sum_i p_i\langle\Sigma^{-1}b_i,h_i\rangle=2\tau Q.
\]
The term involving $\bar h(Z)$ vanishes by independence. Conditional
Gibbs duality therefore gives
\[
 I_P(X;Y\mid Z)\ge2\tau^2Q
 -\E_{X,Z}\log\sum_i k_i(Z)e^{\tau r_i(X,Z)}.
\]
Use $\log u\le u-1$ and the fourth-order exponential remainder. The
linear term vanishes pointwise; the cubic term vanishes after integration
by symmetry of $X$. This yields the lower bound in
\eqref{R10-eq:cost-sandwich}.
For the explicit kernel the information is
$\E_0(1+\tau r_Y)\log(1+\tau r_Y)$. Symmetry averages the integrand with
its value at $-\tau r_Y$, and
\[
 \frac{(1+z)\log(1+z)+(1-z)\log(1-z)}2
 =\frac{z^2}2+\sum_{j\ge2}\frac{z^{2j}}{(2j)(2j-1)}.
\]
The tail is at most $z^4/[12(1-\eta^2)]$. This proves the upper bound.
A minimizing kernel exists by weak-star compactness of finite-action
kernels with the stated integrable moment constraints and lower
semicontinuity of relative entropy. Strict convexity gives uniqueness
up to null sets. For nonsymmetric bounded centered $X$, the same argument
has a cubic remainder and gives $\tau^2Q+O(|\tau|^3)$.
\end{proof}

\begin{corollary}[An independent component gives a dimension-free gap]
\label{R10-cor:positive-component}
Suppose $k_i(z)\ge\alpha p_i$, $0<\alpha\le1$, and put
$\gamma=\alpha(2-\alpha)$. If $\sum_i p_i b_i=0$, then
\[
 Q\le\frac1{2\gamma}\sum_i p_i\|b_i\|_{\Sigma^{-1}}^2.
\]
With $B_X=\esssup_{x,i}|\langle\Sigma^{-1}b_i,x\rangle|$, the sufficient
condition $2|\tau|B_X/\gamma\le\eta<1$ guarantees the explicit lift.
The bound is independent of $\min_i p_i$.
\end{corollary}
\begin{proof}
Let $k'= (k-\alpha p)/(1-\alpha)$ when $\alpha<1$. The reversible
Markov transition $P=D_p^{-1}\E[kk^{\mathsf T}]$ has decomposition
\[
 P=(1-\gamma)P'+\gamma\Pi_p,\qquad
 P'=D_p^{-1}\E[k'k'^{\mathsf T}],\quad \Pi_p=\mathbf1p^{\mathsf T}.
\]
On the $p$-centered subspace, $(I-P)^{-1}$ has norm at most $1/\gamma$
in both $L^2(p)$ and $L^\infty(p)$. The Neumann series proves both claims;
adding a constant to all $h_i$ does not change the kernel. The case
$\alpha=1$ is immediate.
\end{proof}

The transition $P$ is the two-step Gibbs, or data-augmentation, transition
$Y\to Z\to Y'$. Its inverse Poisson operator $(I-P)^{-1}$ therefore gives the local quadratic information cost of calibration. It is also the conditional
covariance matrix underlying the finite exponential-family information
metric. In particular, spectral information sufficient for sampling comparisons
also quantifies the cost of adjoining independent physical means.


\subsection{Exact fair-bit privacy and the information cost of calibration}
\label{FC-sec:privacy}

\subsubsection{A finite formula for full-observation perfect privacy}
Fix a standard Borel random variable $Z$ and a label $Y\in[N]$ with $p_i=\PP(Y=i)>0$. Write
\[
 k_i(z)=\PP(Y=i\mid Z=z),\qquad
 k(A)=\sum_{i\in A}k_i(Z),\qquad
 F(A)=\E\min\{k(A),1-k(A)\}.
\]
A released bit $X\in\{-1,1\}$ is required to be fair and independent of $Z$. Its release mechanism may use the pair $(Z,Y)$, and the law of that pair must be preserved.

\begin{lemma}[Exact feasible correlation polytope]\label{FC-lem:bit-polytope}
The vectors $m_i=\E[X\one_{\{Y=i\}}]$ obtainable by such mechanisms are exactly
\begin{equation}\label{FC-eq:bit-base}
 \mathcal B(F)=\{m\in\R^N:m([N])=0,\ m(A)\le F(A)\text{ for all }A\subseteq[N]\}.
\end{equation}
This is a submodular base polytope. Its vertices are among the greedy vectors
\[
 m^\pi_{\pi_j}=F(A_j)-F(A_{j-1}),\qquad
 A_j=\{\pi_1,\ldots,\pi_j\},\quad\pi\in S_N.
\]
Every greedy vector is realizable by a conditional quantile coupling.
\end{lemma}
\begin{proof}
For a fair sign source, its upper integrated quantile is
$L(s)=\min\{s,1-s\}$. Theorem~\ref{NLI-thm:global} therefore gives
\eqref{FC-eq:bit-base}, and its conditional quantile allocations give the
greedy vertices displayed here.
\end{proof}

Because a fair-bit joint law is determined by $m$, the lemma determines the full utility optimization rather than only a moment relaxation.

\begin{theorem}[Optimal utility of a fair perfectly private bit]
\label{FC-thm:privacy}
With the model above, let
\[
 \mathcal I_{\rm bit}(Z,Y)
 =\sup\{I(X;Y):X\text{ fair binary},\ X\perp Z,\ \law(Z,Y)\text{ fixed}\}.
\]
Then
\begin{equation}\label{FC-eq:privacy-formula}
 \mathcal I_{\rm bit}(Z,Y)
 =\max_{\pi\in S_N}\sum_{j=1}^Np_{\pi_j}
  j\!\left(\frac{F(A_j)-F(A_{j-1})}{p_{\pi_j}}\right),
\end{equation}
where $j(t)=\tfrac{1+t}{2}\log(1+t)+\tfrac{1-t}{2}\log(1-t)$ on $[-1,1]$, with $0\log0=0$. A single greedy conditional quantile coupling attains the maximum.
\end{theorem}
\begin{proof}
Fairness gives
\[
 \PP(X=1,Y=i)=\frac{p_i+m_i}{2},\qquad
 \PP(X=-1,Y=i)=\frac{p_i-m_i}{2}.
\]
Consequently $I(X;Y)=\sum_i p_i j(m_i/p_i)$. This is a continuous convex function of $m$. A maximum over the polytope $\mathcal B(F)$ is attained at a vertex, because every point is a convex combination of vertices and convexity bounds its value by their maximum. Apply Lemma~\ref{FC-lem:bit-polytope}.
\end{proof}

\begin{corollary}[Binary utility and Bayes error]\label{FC-cor:Bayes}
If $Y$ is binary, put $p=\PP(Y=1)$ and
\[
 c=\E\min\{\PP(Y=1\mid Z),\PP(Y=0\mid Z)\}.
\]
Then
\[
 \mathcal I_{\rm bit}(Z,Y)
 =p\,j(c/p)+(1-p)\,j(c/(1-p)).
\]
If $Y$ is fair, $c=e_B$ is the Bayes error for predicting $Y$ from $Z$, and
\begin{equation}\label{FC-eq:Bayes-privacy}
 \mathcal I_{\rm bit}(Z,Y)=j(2e_B).
\end{equation}
\end{corollary}
\begin{proof}
The feasible polytope is $m=(a,-a)$, $|a|\le c$. Its utility is even and increases with $|a|$, so its maximum is at $|a|=c$. The Bayes error is precisely the expected smaller posterior mass.
\end{proof}

For example, if $Y$ is a fair bit observed through a binary symmetric channel of crossover probability $\tau\le1/2$, then $e_B=\tau$ and the optimal fair private output has utility $j(2\tau)$. More generally, all experiments with a fair binary useful label and the same Bayes error have the same value of this constrained full-observation privacy problem.

The full-observation and output-perturbation models are distinguished in Rassouli--G\"und\"uz \cite{RG2021}. In the latter model one additionally requires the Markov chain $Z-Y-X$; the formula above does not impose that restriction. It also fixes the output to be one fair bit. Allowing arbitrary output alphabets and marginals is a larger optimization problem. For finite alphabets, linear programs characterize perfect privacy in the earlier work; \eqref{FC-eq:privacy-formula} gives the explicit posterior-subset formula for this prescribed-output problem, valid for arbitrary standard Borel $Z$.

\subsubsection{Effective resistance and maximal correlation}
Define
\[
 D_p=\diag(p_1,\ldots,p_N),\qquad
 G=D_p-\E[kk^{\mathsf T}].
\]
It is a graph Laplacian with edge conductances $w_{ij}=\E[k_ik_j]$ for $i\ne j$. Its kernel consists of vectors constant on each connected component of that weighted graph. Write $G^\dagger$ for its Moore--Penrose inverse.

\begin{theorem}[Local conditional information cost]\label{FC-thm:resistance}
Let $v\in\operatorname{ran}G$. Among fair bits independent of $Z$ which preserve $(Z,Y)$ and satisfy $m=tv$, the minimum conditional mutual information obeys
\begin{equation}\label{FC-eq:local-information}
 \inf I(X;Y\mid Z)=\frac{t^2}{2}v^{\mathsf T}G^\dagger v+O_v(t^4)
 \qquad(t\to0).
\end{equation}
The same expression without the remainder is a lower bound at every feasible $t$. For $v=e_a-e_b$ in a connected component, the coefficient is half the effective resistance between $a$ and $b$.
\end{theorem}
\begin{proof}
Put $r_i(z)=\E[X\mid Z=z,Y=i]$. Fairness and independence require $\sum_i k_i(z)r_i(z)=0$, and $m_i=\E[k_ir_i]$. Since the conditional law of $X$ given $Z$ is fair,
\[
 I(X;Y\mid Z)=\E\sum_i k_i j(r_i)\ge\frac12\E\sum_i k_i r_i^2.
\]
Let $h=G^\dagger v$ and $r_i^0(z)=t(h_i-\sum_jk_j(z)h_j)$. Then $\sum k_i r_i^0=0$ and $\E[k_ir_i^0]=t(Gh)_i=tv_i$. Every other admissible $r$ satisfies
\[
 \E\sum_i k_i(r_i-r_i^0)r_i^0
 =t\sum_i h_i\E[k_i(r_i-r_i^0)]=0.
\]
Hence its quadratic energy is at least that of $r^0$, namely $t^2v^{\mathsf T}G^\dagger v$. For sufficiently small $|t|$, $|r_i^0|\le2|t|\norm{h}_\infty<1$, so this kernel is feasible. Expanding $j(r)=r^2/2+O(r^4)$ uniformly in that bounded range gives the upper bound in \eqref{FC-eq:local-information}. Effective resistance is $(e_a-e_b)^{\mathsf T}G^\dagger(e_a-e_b)$.
\end{proof}

Every finite weighted graph can occur in this way up to scale. If its conductances are $w_{ij}$ with $W=\sum_{i<j}w_{ij}$, choose $0<c<1/(4W)$ when $W>0$. Give a baseline atom for each edge probability $4cw_{ij}$ and posterior $(e_i+e_j)/2$. Allocate the remaining probability among singleton posteriors, positively at every vertex. The resulting matrix $G$ is exactly $c$ times the original Laplacian and all $p_i$ are positive. Thus the resistance coefficient can represent any finite electrical network.

For the operator interpretation, assume $N\ge2$ and set
\[
 C=D_p^{-1/2}\E[kk^{\mathsf T}]D_p^{-1/2}.
\]
The vector $\sqrt p$ is an eigenvector of eigenvalue one. On its orthogonal complement, the largest eigenvalue is the square $\rho_{\rm HGR}(Y,Z)^2$ of the Hirschfeld--Gebelein--R\'enyi maximal correlation. Indeed, for $\E b(Y)=0$,
\[
 \norm{\E[b(Y)\mid Z]}_2^2
 =b^{\mathsf T}\E[kk^{\mathsf T}]b,
 \qquad \norm{b(Y)}_2^2=b^{\mathsf T}D_pb.
\]
If $\rho_{\rm HGR}<1$, the worst coefficient in \eqref{FC-eq:local-information} among $v=D_pb$ with $\E b(Y)=0$ and $\E b(Y)^2=1$ is therefore
\begin{equation}\label{FC-eq:HGR-cost}
 \frac1{2(1-\rho_{\rm HGR}(Y,Z)^2)}.
\end{equation}
This follows by solving the Laplacian equation on the centered subspace and using $G=D_p^{1/2}(I-C)D_p^{1/2}$. The maximum is achieved at an eigenvector corresponding to the largest nonconstant eigenvalue of $C$. The reciprocal spectral gap thus measures the most expensive infinitesimal calibration direction. Section~\ref{NLI-sec:global} identifies the related two-step data-augmentation operator; \eqref{FC-eq:HGR-cost} and the resistance formula give explicit named interpretations of that operator.


\section{Nonlinear sparsification and exact latent realization}
\label{R11-sec:latent}
A posterior Laplacian describes the local information metric, while calibration costs also contain nonlinear finite-amplitude terms. The nonnegative cosine representation of $\log\cosh$ converts them into graph energies. A spectral comparison~\cite{FC-BSS} therefore controls the whole potential at every amplitude for every symmetric integrable source.

An exact latent experiment instead requires a completely positive posterior second moment. Its cp-rank is the least finite latent cardinality. The examples distinguish this exact realization complexity from the sparse approximation of all calibration costs.

\subsection{Nonlinear graph representations and sparse approximation}\label{NLI-sec:graphs}
An edge reference has latent states consisting of unordered pairs $e=\{i,j\}$ with masses $w_{ij}\ge0$ and individual vertices with masses $v_i\ge0$. Given an edge, the target is fair on its endpoints; given a vertex, it is deterministic. The target marginal is
\begin{equation}\label{NLI-eq:edge-marginal}
 p_i=v_i+\frac12\sum_{j\ne i}w_{ij},\qquad \sum_i v_i+\sum_{i<j}w_{ij}=1.
\end{equation}
The source is independent of the latent edge or vertex. Orient the edges arbitrarily and let $B$ be their incidence matrix, with column $e_i-e_j$ for an edge oriented from $j$ to $i$.

\subsubsection{The complete information cost is a convex flow problem}
Let $X$ have a symmetric integrable scalar law. Define
\[
 F_X(t)=\E\log\cosh(tX),\qquad \mathcal I_X=F_X^*.
\]
The symmetric binary duality in Section~\ref{TOMO-sec:binary} identifies $\mathcal I_X(b)$ with the minimum information needed for a fair sign $\sigma$ satisfying $\E[X\mid\sigma]=\sigma b$.

\begin{theorem}\label{NLI-thm:flow}
The edge reference has dual potential
\begin{equation}\label{NLI-eq:edge-potential}
 \Phi_w^X(u)=\sum_{i<j}w_{ij}F_X((u_i-u_j)/2).
\end{equation}
Its minimum conditional information at joint means $m$ is
\begin{equation}\label{NLI-eq:flow}
 J_w^X(m)=\min_{Bf=m}\sum_e w_e\mathcal I_X(2f_e/w_e),
\end{equation}
with zero-weight edges omitted and infeasible constraints assigned value $+\infty$.
\end{theorem}
\begin{proof}
On edge $e=(i,j)$ use the sign $\sigma=1$ at $i$ and $-1$ at $j$. If its conditional mean parameter is $b_e$, the edge contributes $w_eb_e/2$ to $m_i$ and its negative to $m_j$. Its information cost is $w_e\mathcal I_X(b_e)$. Setting $f_e=w_eb_e/2$ proves \eqref{NLI-eq:flow}. Deterministic vertex states contribute neither mean nor information because the centered source is independent of the latent state. Fenchel duality on the finite incidence constraint gives \eqref{NLI-eq:edge-potential}. Existence follows from compactness of the bounded feasible edge means and lower semicontinuity of the costs.
\end{proof}
Near zero this reduces to electrical energy. Formula~\eqref{NLI-eq:flow} retains the finite-amplitude entropy cost at every edge.

\subsubsection{A global comparison for an arbitrary baseline}
For the original posterior law $k(Z)$, define
\begin{equation}\label{NLI-eq:surrogate}
 w_{ij}=2\E[k_ik_j],\qquad v_i=\E[k_i^2],\qquad
 \Psi(u)=\sum_{i<j}w_{ij}\log\cosh((u_i-u_j)/2).
\end{equation}
The masses satisfy \eqref{NLI-eq:edge-marginal}, so this is a genuine edge reference with exactly the original target marginal.

\begin{theorem}\label{NLI-thm:graph-comparison}
For every finite-target baseline and every $u\in\R^N$,
\begin{equation}\label{NLI-eq:comparison}
 \Psi(u)\le\Phi_k(u)\le\Psi(2u),\qquad
 \frac12\Psi(2u)\le\Phi_k(u).
\end{equation}
Writing $J_{\rm edge}=\Psi^*$, one has, as extended-valued inequalities,
\begin{equation}\label{NLI-eq:primal-comparison}
 \boxed{2J_k(m)\le J_{\rm edge}(m)\le J_k(2m).}
\end{equation}
In particular, the edge feasible region contains half the original feasible region.
\end{theorem}
\begin{proof}
Work first with a fixed posterior $k$. Let
\[
 f_k(u)=\inf_\alpha\sum_i k_i\log\cosh(u_i-\alpha).
\]
Choosing $\alpha=u_j$ and averaging over $j\sim k$ proves
\[
 f_k(u)\le\sum_{i,j}k_ik_j\log\cosh(u_i-u_j).
\]
Convexity and evenness of $\log\cosh$ give, for every $\alpha$,
\[
 \log\cosh((u_i-u_j)/2)
 \le\tfrac12\log\cosh(u_i-\alpha)+\tfrac12\log\cosh(u_j-\alpha).
\]
Averaging proves the first lower bound.

For the stronger lower bound take the minimizing $\alpha$, put $a_i=u_i-\alpha$ and $t_i=\tanh a_i$. Then $\sum_i k_it_i=0$ and
\[
 \log\cosh(a_i-a_j)=\log\cosh a_i+\log\cosh a_j+\log(1-t_it_j).
\]
Jensen's inequality gives
\[
 \sum_{i,j}k_ik_j\log(1-t_it_j)
 \le\log\left(1-(\sum_i k_it_i)^2\right)=0.
\]
Thus the last pair sum is at most $2f_k(u)$. Integrating all three inequalities over $Z$ proves \eqref{NLI-eq:comparison}. Taking Fenchel conjugates, with the appropriate argument and scalar rescalings, proves \eqref{NLI-eq:primal-comparison}.
\end{proof}

\begin{theorem}[Lifting retains the entire old pair law]\label{NLI-thm:lift}
Every coupling for the unsparsified edge reference \eqref{NLI-eq:surrogate}, with any independent centered scalar source $X$, can be lifted to a coupling of $(X,Z,Y)$ that retains the original entire $(Z,Y)$ law and $X\perp Z$. It has the same joint means and no larger conditional information.
\end{theorem}
\begin{proof}
Draw $Z$ and $X$ independently. Conditional on $Z=z$, draw two independent labels $I,J$ from $k(z)$, independently of $X$. Let $E$ be their unordered pair, interpreted as a vertex if they coincide. Its law has exactly the masses in \eqref{NLI-eq:surrogate}. Apply the given edge coupling to choose $Y$ from $X,E$; it uses the fair endpoint marginal or deterministic vertex marginal of that edge model. Averaging over $X$ and the pair shows
\[
 \Pp(Y=i\mid Z=z)=k_i(z)^2+\sum_{j\ne i}k_i(z)k_j(z)=k_i(z).
\]
The $(X,E,Y)$ marginal is the given edge experiment, so its $(X,Y)$ law and all its means are preserved. Since $X$ is independent of $(E,Z)$ and the kernel depends on $E,X$ only, $I(X;Y\mid Z)\le I(X;E,Y\mid Z)=I(X;Y\mid E)$. \end{proof}
This lifting result is useful when $Z$ already includes every auxiliary variable of a previously constructed reference. All those variables remain in the retained pair law.

\begin{corollary}[The exact two-step transition]\label{NLI-cor:lazy}
Let $M=\E[kk^{\mathsf T}]$ and $P=D_p^{-1}M$. The edge reference in \eqref{NLI-eq:surrogate} has second-moment matrix $(D_p+M)/2$ and two-step transition $(I+P)/2$.
\end{corollary}
\begin{proof}
Its off-diagonal second moment is $w_{ij}/4=M_{ij}/2$. Its diagonal is $v_i+\sum_jw_{ij}/4=M_{ii}+(p_i-M_{ii})/2$. This is exactly the stated matrix.
\end{proof}

\begin{proposition}[The universal half-mean factor is sharp]\label{NLI-prop:half}
Take a constant posterior $k=p$ with $p_i=1/N$ and even $N$. Among all edge references with this same target marginal, the largest factor $c$ for which their feasible regions can contain $c$ times the original fair-binary region is
\[
 c_N=\frac{N}{2(N-1)}.
\]
Consequently the dimension-independent factor $1/2$ cannot be improved within edge references.
\end{proposition}
\begin{proof}
For a balanced set $A$, the original support bound is $F(A)=1/2$. In an edge reference, its corresponding maximum is half the total weight crossing $A$. A uniformly random balanced cut separates each pair with probability $N/[2(N-1)]$. Since total edge mass is at most one, some balanced cut has capacity at most $N/[4(N-1)]$. Containment therefore forces $c/2\le N/[4(N-1)]$.

For attainment use every edge with weight $1/\binom N2$ and no vertex states. This retains $p$. A set of size $a\le N/2$ has capacity $a(N-a)/(N(N-1))$, whereas its original bound is $a/N$. The ratio is $(N-a)/(N-1)$, minimized at $a=N/2$. The base-polytope description proves the desired containment for all sets.
\end{proof}

\subsubsection{Spectral approximation controls the whole nonlinear potential}
The identity
\begin{equation}\label{NLI-eq:levy}
 \log\cosh x=\int_0^\infty\frac{1-\cos(sx)}{s\sinh(\pi s/2)}\,ds
\end{equation}
follows by differentiating in $x$, using the Fourier sine integral for $\tanh x$, and matching the value at zero. Both sides are even. Equivalently it follows from the infinite product for $\cosh x$ by representing each $\log(1+x^2/a^2)$ as a cosine integral and summing the positive kernels. It is the L\'evy--Khintchine representation of $\log\cosh$, consistent with the classical infinite divisibility of hyperbolic secant laws \cite{NLI-JurekYor}.

At each frequency $s$, the edge contribution $1-\cos(s(u_i-u_j)/2)$ is half the sum of the squared differences of the cosine and sine coordinates. A Laplacian inequality therefore controls it for every phase vector. Integrating these inequalities against the positive kernel in \eqref{NLI-eq:levy} retains the same constants at arbitrary amplitude; averaging over the source retains them for every symmetric integrable law.

\begin{theorem}\label{NLI-thm:nonlinear-sparse}
Let $L_w$ be the weighted graph Laplacian of an edge reference. If another nonnegative edge system satisfies
\[
 cL_w\preceq L_{w'}\preceq L_w,\qquad 0<c\le1,
\]
then for every symmetric integrable source $X$ and every $u\in\R^N$,
\begin{equation}\label{NLI-eq:uniform-sparse}
 c\Phi_w^X(u)\le\Phi_{w'}^X(u)\le\Phi_w^X(u),
\end{equation}
and
\begin{equation}\label{NLI-eq:sparse-cost}
 J_w^X(m)\le J_{w'}^X(m)\le cJ_w^X(m/c).
\end{equation}
For each $0<\epsilon<1$, one can choose $O(N/\epsilon^2)$ positive edge states and at most $N$ vertex states, retaining the target marginal $p$ exactly, with
\[
 c=\frac{1-\epsilon}{1+\epsilon}.
\]
The same sparse graph works simultaneously for all such source laws and all amplitudes.
\end{theorem}
\begin{proof}
For any fixed frequency $s$ and source value $x$, write
\[
 1-\cos(sx(u_i-u_j)/2)
 =\frac12\left\|\begin{pmatrix}\cos(sxu_i/2)\\\sin(sxu_i/2)\end{pmatrix}
 -\begin{pmatrix}\cos(sxu_j/2)\\\sin(sxu_j/2)\end{pmatrix}\right\|^2.
\]
Apply the Laplacian inequalities to the two real coordinate vectors, then integrate the positive representation \eqref{NLI-eq:levy} and average over $X$. Tonelli's theorem applies; finiteness follows from $\log\cosh(tX)\le |t|\,|X|$. This proves \eqref{NLI-eq:uniform-sparse}. Conjugation gives \eqref{NLI-eq:sparse-cost}.

By the Batson--Spielman--Srivastava theorem \cite{FC-BSS}, there are $O(N/\epsilon^2)$ reweighted edges with $(1-\epsilon)L_w\preceq L_{\widetilde w}\preceq(1+\epsilon)L_w$. Set $w'=\widetilde w/(1+\epsilon)$. Diagonal entries of the upper Loewner bound imply that its degree at each vertex is at most the original degree. Hence
\[
 v_i'=p_i-\tfrac12\sum_jw_{ij}'\ge0
\]
completes the sparse graph to a probability reference with exactly the target marginal $p$. These new vertex masses do not affect its information potential.
\end{proof}
The Fourier proof also works for centrally symmetric vector sources, with vector-valued dual labels $u_i$, because each fixed source value produces scalar phases $\langle u_i,X\rangle$. It extends further to any edge loss with a nonnegative cosine-integral representation, with its quadratic term treated directly.

\begin{corollary}[Sparse approximation of an arbitrary fair baseline]\label{NLI-cor:general-sparse}
Every finite-target baseline admits a new edge-and-vertex reference with $O(N/\epsilon^2)$ latent states and the same target marginal such that
\[
 c\Phi_k(u/2)\le\Phi_{\rm sparse}(u)\le\Phi_k(u),
 \qquad c=(1-\epsilon)/(1+\epsilon).
\]
\end{corollary}
\begin{proof}
Use the surrogate of \eqref{NLI-eq:surrogate}. Theorem~\ref{NLI-thm:graph-comparison} gives $\Phi_k(u/2)\le\Psi(u)\le\Phi_k(u)$, and Theorem~\ref{NLI-thm:nonlinear-sparse} sparsifies $\Psi$.
\end{proof}

\begin{remark}[Which law is preserved]
Theorem~\ref{NLI-thm:lift} retains the actual entire $(Z,Y)$ law because its pair is sampled conditionally from $k(Z)$. Spectral sparsification changes the latent reference and retains $p$ and the stated potential inequalities. A sparse reweighting may increase individual edge weights, so it need not admit the same conditional pair construction. Applying a sparse solver inside a fixed original baseline requires an additional positive lifting representation. This distinction is part of the theorem statements.
\end{remark}

\subsection{Positive two-step realizations and exact latent complexity}\label{NLI-sec:cp}
A real symmetric matrix $M$ is completely positive if
\[
 M=\sum_{r=1}^R b_rb_r^{\mathsf T},\qquad b_r\in[0,\infty)^N.
\]
Its \emph{cp-rank} is the minimum number of vectors in such a factorization. This is an established positive-factorization problem; its graph restrictions are developed in \cite{NLI-Berman}. The following normalization identifies its exact probabilistic role for finite reference experiments.

\begin{theorem}\label{NLI-thm:cp}
Let $p_i>0$, $\sum_i p_i=1$, and let $P$ be a stochastic matrix reversible with respect to $p$. There is a baseline experiment $(Z,Y)$ with $Y\sim p$ whose two-step conditional transition $Y\to Z\to Y'$ is $P$ if and only if $D_pP$ is completely positive. The minimum finite latent cardinality is
\[
 \cprank(D_pP).
\]
\end{theorem}
\begin{proof}
Every baseline satisfies
\[
 D_pP=\E[k(Z)k(Z)^{\mathsf T}],
\]
which is completely positive. Finite-dimensional compactness of the normalized rank-one cone, or conic Carath\'eodory, gives a finite factorization even for continuous $Z$. A baseline on $R$ latent states gives a factorization with at most $R$ terms.

Conversely suppose $M=D_pP=\sum_r b_rb_r^{\mathsf T}$ with $b_r\ge0$. Discard zero vectors and put
\[
 s_r=\ind^{\mathsf T}b_r,\quad q_r=s_r^2,\quad k^{(r)}=b_r/s_r.
\]
Then
\[
 \sum_rq_r=\ind^{\mathsf T}M\ind=1,\qquad
 \sum_rq_r k^{(r)}=M\ind=p,
\]
and $\sum_rq_r k^{(r)}(k^{(r)})^{\mathsf T}=M$. Use $q_r$ as the latent probabilities and $k^{(r)}$ as the target posteriors. The two-step transition is the prescribed $P$. Minimizing the number of terms proves the cardinality assertion.
\end{proof}
Positive semidefiniteness of the reversible transition is a necessary spectral condition. Complete positivity incorporates the nonnegativity of the conditional probabilities in the factorization.

\begin{proposition}[A five-state obstruction and its exact correction]\label{NLI-prop:c5}
Let $A$ be the adjacency matrix of the five-cycle and set
\[
 P=\frac{5I+3A}{11},\qquad p_i=1/5.
\]
This is a positive-definite stochastic reversible transition. It has no two-step posterior realization. For
\[
 P_s=(1-s)I+sP,\qquad 0\le s\le1,
\]
a two-step posterior realization exists exactly when $s\le11/12$.
\end{proposition}
\begin{proof}
The smallest eigenvalue of $P$ is $(5-3(1+\sqrt5)/2)/11>0$. Its nonzero off-diagonal pattern is the five-cycle. In any completely positive factorization, each vector must be supported on a clique of that pattern, because a positive product on a nonedge cannot cancel. The graph is triangle-free, so supports have size at most two. For every such vector,
\[
 \sum_i b_i^2\ge2\sum_{i<j}b_ib_j.
\]
Consequently any completely positive matrix with this zero pattern satisfies trace at least twice the sum of its upper-triangular off-diagonal entries. For $P_s$, this requires $5(1-6s/11)\ge30s/11$, or $s\le11/12$. The same test excludes $s=1$.

For the converse, the explicit factorization is
\[
 P_s=\frac{3s}{11}\sum_{\{i,j\}\in E(C_5)}(e_i+e_j)(e_i+e_j)^{\mathsf T}
       +(1-12s/11)I.
\]
It is completely positive throughout the stated range. Multiplying by $1/5$ and applying Theorem~\ref{NLI-thm:cp} gives the desired baseline.
\end{proof}
Thus adding a holding probability of at least $1/12$ repairs this particular positive realization obstruction. The parameter $s$ is the weight of the original transition.

\subsubsection{Triangle-free second moments determine the posterior law}
\begin{theorem}\label{NLI-thm:rigidity}
Let a triangle-free graph carry a posterior law supported entirely on its fair edges,
\[
 k=\tfrac12(e_i+e_j),\qquad \Pp(k=\tfrac12(e_i+e_j))=w_{ij},\qquad \sum w_{ij}=1.
\]
Its matrix $M=\E[kk^{\mathsf T}]$ determines this entire probability law uniquely among all random vectors in the probability simplex.
\end{theorem}
\begin{proof}
If $\widetilde k$ has the same second moment, every nonedge satisfies $\E\widetilde k_i\widetilde k_j=0$. Nonnegativity implies $\widetilde k_i\widetilde k_j=0$ almost surely. Thus its support is almost surely a clique, and has size at most two. Every probability vector on at most two points has squared norm at least $1/2$, with equality exactly at a fair pair. But $\E\|\widetilde k\|^2=\tr M=1/2$. Therefore $\widetilde k$ is almost surely a fair pair. Finally $M_{ij}=w_{ij}/4$ identifies each edge mass.
\end{proof}

\begin{corollary}[Quadratic exact complexity and linear approximate complexity]\label{NLI-cor:complexity}
Put $N=2q$, $q\ge2$, and give the $q^2$ edges of $K_{q,q}$ equal mass $1/q^2$. Then
\[
 \rank M=2q-1,\qquad \cprank M=q^2.
\]
Every latent reference preserving its second moments needs at least $q^2$ states. The same lower bound holds for any latent reference with the same $p$ and the same full fair-binary calibration potential. For each fixed $0<\epsilon<1$, Theorem~\ref{NLI-thm:nonlinear-sparse} gives a reference with $O(q/\epsilon^2)$ states, the same $p$, and its uniform nonlinear potential comparison.
\end{corollary}
\begin{proof}
In the bipartition order,
\[
 M=\frac1{4q^2}\begin{pmatrix}qI&\ind\ind^{\mathsf T}\\
                       \ind\ind^{\mathsf T}&qI\end{pmatrix}.
\]
Its eigenvalues have one zero, so its rank is $2q-1$. Every factor in a completely positive representation has support of size at most two and can cover at most one of the $q^2$ positive off-diagonal edges. Hence at least $q^2$ factors are required, and the given fair-edge representation attains that number. Alternatively Theorem~\ref{NLI-thm:rigidity} fixes the whole posterior law and its $q^2$ distinct values.

The Hessian of the fair-binary potential is $D_p-M$ by Theorem~\ref{NLI-thm:potential}. Equality of potentials at the prescribed marginal therefore implies equality of $M$. The lower bound follows. The sparse assertion applies directly to the existing edge reference.
\end{proof}
In this class, quadratic response already determines all higher posterior information. The example also specifies an unavoidable exact representation cost even though its ordinary matrix rank grows linearly. The sparse bound measures approximation of the complete potential, uniformly over its arguments, rather than a fit to finitely many derivatives.


\section{Exact privacy capacity and smooth-data instability}
\label{R11-sec:privacy}
The observation model matters for perfect privacy. In the useful-data-only model a release is generated from $Y$; in the full-data model it may also use the private variable $Z$. Rassouli--G\"und\"uz study both and formulate the finite-alphabet Gaussian question in the latter model~\cite{RG2021}. The finite-alphabet and positive-budget theorems use the full-data model, with $U\perp Z$ and a prescribed output law. At zero conditional-information budget, the Markov constraint $Z-Y-U$ restricts the release to the useful-data-only model; the smooth-data and Fourier theorems state that restriction explicitly.

Conditional atomlessness gives labels that are almost conditionally balanced. Their error is controlled by the probability that two conditionally independent observations occupy the same fine cell. Conditional quantile allocation repairs the labels on a set of vanishing probability while making independence exact. A conditional atom gives the converse obstruction when the alphabet size ranges over all finite values. Compatible refinements with summable conditional information costs then produce one uniform release with infinite useful information. The Gaussian example keeps the logical conclusions separate: each finite-alphabet supremum is maximal and unattained, while a single infinite release attains infinite utility under any positive conditional-information budget.

Fixing that budget at zero raises a further question: how much of the Gaussian value is determined by the smooth geometry of the data law? Uniform masking and exact variance compensation keep the covariance fixed and the potential Hessian uniformly close to Gaussian, while producing infinite utility. Fourier completeness then distinguishes the zero values from the finite and infinite examples.

\subsection{Exact conditional marginal correction}
\label{FP-sec:repair}
The preceding fair-bit formula treats a prescribed finite useful
alphabet. We now allow arbitrary standard Borel private data $Z$ and
useful data $Y$. The released symbol $U\in[M]$ is produced by a kernel
observing the full pair $(Z,Y)$. Set
\begin{equation}\label{FP-eq:privacy-value}
 G^M_\epsilon(Z,Y)=\sup_{P_{U\mid Z,Y}:\ I(U;Z)\le\epsilon}I(U;Y).
\end{equation}
Define $h_b(t)=-t\log t-(1-t)\log(1-t)$ and, for $M\ge2$,
\[
 \omega_M(t)=h_b(s)+s\log(M-1),\qquad
 s=\min\{t,1-1/M\}.
\]
An approximately private symbol can be corrected to exact privacy by
changing it only where its conditional mass differs from its target
mass. The correction keeps the entire original experiment. Its error
is measured in total variation, so finite-alphabet information changes
by a controlled amount.

For finite probability vectors $q,p$, use $\TV(q,p)=\frac12\sum_i|q_i-p_i|$.

\begin{lemma}[Conditional maximal coupling]\label{FP-lem:repair}
Fix the complete law of $(Z,Y,U)$, where $U\in[M]$, and let $p$ be a probability vector. Write $q_i(z)=\Pp(U=i\mid Z=z)$. There is an extension by $U^*\in[M]$ such that
\[
 \Pp(U^*=i\mid Z)=p_i,
 \qquad
 \Pp(U^*\ne U)=\E\TV(q(Z),p).
\]
The disagreement probability is the least possible among extensions satisfying the prescribed conditional output law. The original $(Z,Y,U)$ law is unchanged.
\end{lemma}
\begin{proof}
For each $z$, set $d=\TV(q(z),p)$ and define the joint mass matrix
\begin{equation}\label{FP-eq:max-coupling}
 r_{ij}(z)=\ind_{\{i=j\}}\min(q_i(z),p_i)
 +\frac{(q_i(z)-p_i)_+(p_j-q_j(z))_+}{d},
\end{equation}
where the second term is zero when $d=0$. The positive and negative discrepancies both sum to $d$, so the row sums are $q_i(z)$ and the column sums are $p_j$. The product term vanishes on the diagonal. Given $(Z,Y,U)=(z,y,i)$, sample $U^*$ using the row $r_{ij}(z)/q_i(z)$; choose any row at $q_i(z)=0$. This is a measurable kernel. Its diagonal mass is $\sum_i\min(q_i,p_i)=1-d$. In any coupling, the $i$th diagonal mass is at most $\min(q_i,p_i)$, proving optimality after integration.
\end{proof}
The lemma is the classical maximal-coupling construction, applied to regular conditional probabilities. Its role here is to retain the original experiment exactly while imposing the new independent marginal.

\begin{lemma}[Finite-alphabet utility loss]\label{FP-lem:utility-loss}
For jointly defined $U,U^*\in[M]$, with disagreement probability $e$,
\[
 I(U;Y)\le I(U^*;Y)+h_b(e)+e\log(M-1).
\]
Consequently the same bound holds with the last two terms replaced by $\omega_M(t)$ whenever $e\le t$.
\end{lemma}
\begin{proof}
The chain rule gives $I(U;Y)\le I(U,U^*;Y) \le I(U^*;Y)+H(U\mid U^*)$. Writing $E_0=\ind_{\{U\ne U^*\}}$, the conditional entropy satisfies $H(U\mid U^*)\le H(E_0)+\Pp(E_0=1)\log(M-1)$. The right side increases until $e=1-1/M$, where it equals $\log M$, and never exceeds $\log M$. This proves the saturated form as well.
\end{proof}

\begin{theorem}[Uniform continuity at exact privacy]\label{FP-thm:privacy-continuity}
For every joint law of $(Z,Y)$ and every $M\ge2$,
\[
 0\le G^M_\epsilon(Z,Y)-G^M_0(Z,Y)
 \le\omega_M\!\left(\sqrt{\epsilon/2}\right).
\]
The same inequality holds when both suprema require one prescribed output law $p$ on $[M]$.
\end{theorem}
\begin{proof}
Take any admissible $U$ and set $p=\law(U)$. By conditional Pinsker and Jensen,
\begin{align*}
 e:=\E\TV(P_{U\mid Z},P_U)
 &\le\E\sqrt{\tfrac12D(P_{U\mid Z}\Vert P_U)}\\
 &\le\sqrt{\tfrac12I(U;Z)}
 \le\sqrt{\epsilon/2}.
\end{align*}
Lemma~\ref{FP-lem:repair} gives $U^*\perp Z$, with the same output law $p$ and disagreement $e$. Its kernel is available to a full-data encoder, by composing the old encoder with the correction depending on $(Z,U)$. Lemma~\ref{FP-lem:utility-loss} bounds the loss. Taking the supremum over $U$ proves the result. No optimizing encoder or compactness assumption is needed.
\end{proof}

\begin{corollary}[Bounded decision rewards]\label{FP-cor:reward}
Suppose $c(y,i)$ is measurable and
\[
 L=\sup_y\bigl(\max_i c(y,i)-\min_i c(y,i)\bigr)<\infty.
\]
The correction in Theorem~\ref{FP-thm:privacy-continuity} changes expected reward by at most $L\sqrt{\epsilon/2}$. In particular the full-data optimum for this reward is right-continuous at zero leakage, with this modulus.
\end{corollary}
\begin{proof}
The reward difference is zero on $\{U=U^*\}$ and has absolute value at most $L$ elsewhere.
\end{proof}

\begin{remark}[Scope of the correction]
The correction may use $Z$. It therefore preserves full-data observation, but does not in general preserve the output-perturbation restriction $Z-Y-U$. It can also change a pointwise input-dependent output constraint. The claims above permit every output symbol and constrain its distribution, privacy, and reward. They do not assert an automatic correction under arbitrary additional admissibility relations.
\end{remark}

\begin{corollary}[Short resolution of the Gaussian conjecture]\label{FP-cor:RG-short}
For a nondegenerate jointly Gaussian pair with $0<|\rho|<1$ and every $M\ge2$, $G_0^M=\log M$.
\end{corollary}
\begin{proof}
Rassouli--G\"und\"uz establish $G_\epsilon^M=\log M$ for every $\epsilon>0$ immediately before their Conjecture~1 \cite{RG2021}. Theorem~\ref{FP-thm:privacy-continuity} and $\epsilon\downarrow0$ give the lower bound $G_0^M\ge\log M$. The entropy bound $I(U;Y)\le H(U)\le\log M$ gives equality.
\end{proof}
Section~\ref{FP-sec:diffuse} gives a direct construction under substantially weaker distributional hypotheses and identifies Gaussian nonattainment.


\subsection{Exact privacy for conditionally atomless data}
\label{FP-sec:diffuse}
\begin{lemma}[Deterministic conditional correction]\label{FP-lem:deterministic-repair}
Suppose $B=b(Y)\in[M]$ and $P_{Y\mid Z=z}$ is atomless almost surely. The correction of Lemma~\ref{FP-lem:repair} from $B$ to a prescribed output law $p$ can be realized by $U=F(Z,Y)$, with the same conditional joint law of $(B,U)$ and the same minimum disagreement.
\end{lemma}
\begin{proof}
Choose a Borel embedding of the useful-data space into $\R$, and use its order. For each $i$, put
\[
 R_i(z,y)=P(Y'\le y,\ b(Y')=i\mid Z=z).
\]
This function is jointly measurable. If $q_i(z)>0$, then, conditional on $Z=z,B=i$, the variable $R_i(z,Y)$ is uniform on $[0,q_i(z)]$: the restricted distribution is atomless and its cumulative distribution function is continuous. Divide this interval into pieces of lengths $r_{ij}(z)$ from~\eqref{FP-eq:max-coupling}, in any fixed order with $j=i$ first. Output $j$ on its piece. Conditional masses are exactly $r_{ij}(z)$, and the diagonal piece has length $\min(q_i,p_i)$. All endpoints are measurable; zero-mass and endpoint events can be assigned arbitrarily.
\end{proof}

\begin{theorem}[Full finite-output entropy under conditional atomlessness]
\label{R9-thm:conditional-capacity}\label{SP-thm:diffuse}
Let $(Z,Y)$ be standard Borel and suppose the regular conditional law
$\nu_z=P_{Y\mid Z=z}$ is atomless for $P_Z$-almost every $z$.
For every finite probability vector $p=(p_1,\ldots,p_M)$,
\begin{equation}\label{R9-eq:conditional-capacity}
 \sup_{\substack{U=F(Z,Y),\ U\perp Z\,;\ \law(U)=p}} I(U;Y)=H(p).
\end{equation}
The encoders in this supremum are deterministic measurable functions.
For every $\epsilon>0$ there is also a deterministic function $B(Y)$
with law $p$, $I(B;Z)\le\epsilon$, and $I(B;Y)=H(p)$.

More precisely, let $\mathcal P_n$ be any sequence of finite refining
Borel partitions separating points of the useful-data space, and put
\[
 c_n=\E\sum_{C\in\mathcal P_n}\nu_Z(C)^2,
 \qquad e_n=\tfrac12\sqrt{(M-1)c_n}.
\]
Then $c_n\downarrow0$, and there are exactly private deterministic
encoders $U_n$, each with law $p$, such that
\begin{equation}\label{R9-eq:privacy-collision}
 H(p)-I(U_n;Y)\le\omega_M(e_n).
\end{equation}
\end{theorem}
\begin{proof}
Remove zero coordinates of $p$; the one-label case is immediate.
Choose a Borel embedding of the useful-data space into $[0,1]$.
Dyadic partitions of its image provide the required finite sequence.
If $Y',Y''$ are conditionally independent with law $\nu_Z$, then
$c_n$ is the probability that they belong to the same cell of
$\mathcal P_n$. These events decrease to $\{Y'=Y''\}$, whose
conditional probability is zero by atomlessness. Hence $c_n\downarrow0$.

Independently color each cell $C$ with a label $L_C$ of law $p$.
For a fixed coloring let $B_n(Y)=L_C$ on $C$, and write
$q_{n,i}(z)=P(B_n=i\mid Z=z)$. Averaging only over the independent
colors gives
\[
 \E_L q_{n,i}(z)=p_i,\qquad
 \E_L(q_{n,i}(z)-p_i)^2
       =p_i(1-p_i)\sum_{C\in\mathcal P_n}\nu_z(C)^2.
\]
Jensen and Cauchy--Schwarz therefore give
\begin{align*}
 \E_L\E_Z\TV(q_n(Z),p)
 &\le\tfrac12\sum_i\sqrt{p_i(1-p_i)}\sqrt{c_n}\\
 &\le\tfrac12\sqrt{(M-1)c_n}=e_n.
\end{align*}
There is a deterministic coloring with this bound. By
Lemma~\ref{FP-lem:deterministic-repair}, its label can be corrected to
$U_n=F_n(Z,Y)$ with $\law(U_n\mid Z)=p$ and
$P(U_n\ne B_n)\le e_n$. Since $B_n$ is determined by $Y$, $H(U_n\mid Y)\le H(U_n\mid B_n)\le\omega_M(e_n)$. This proves \eqref{R9-eq:privacy-collision} and the lower bound in
\eqref{R9-eq:conditional-capacity}. The upper bound is $I(U;Y)\le H(U)$.

For the positive-leakage assertion, the unconditional law of $Y$ is
also atomless. Write $r_n=\law(B_n)$. We have
$\TV(r_n,p)\le\E_Z\TV(q_n(Z),p)\le e_n$.
Splitting the label sets using the atomless law of $Y$ changes $B_n$
on exactly this much mass to a function $B'_n(Y)$ with law $p$.
The conditional total-variation triangle inequality gives $\E_Z\TV(\law(B'_n\mid Z),p)\le2e_n$. Put $p_* =\min_i p_i>0$. For any probability vector $q$,
\[
 D(q\|p)\le\sum_i\frac{(q_i-p_i)^2}{p_i}
       \le\frac{2}{p_*}\TV(q,p).
\]
Thus $I(B'_n;Z)\le4e_n/p_*\to0$, while
$I(B'_n;Y)=H(p)$ exactly. Take a sufficiently large $n$.
\end{proof}

\begin{corollary}[Attainment and an arbitrary entropy budget]
\label{R9-cor:privacy-attainment}\label{FP-cor:entropy-budget}
Under Theorem~\ref{R9-thm:conditional-capacity}, the supremum for a fixed
output law $p$ is attained exactly when there is a measurable $B(Y)$
of law $p$ independent of $Z$. For every $R\ge0$,
\[
 \sup\{I(U;Y):U\perp Z,\ U\text{ finite-valued},\ H(U)\le R\}=R.
\]
In particular these conclusions apply whenever $Y$ is atomless and
$I(Z;Y)<\infty$.
\end{corollary}
\begin{proof}
Attainment of $I(U;Y)=H(p)$ is equivalent to $H(U\mid Y)=0$,
which for a finite output makes $U$ a measurable function of $Y$.
The converse is immediate. Every $R\ge0$ is the entropy of a finite
probability vector, so the capacity theorem proves the second claim.
Finite mutual information implies
$P_{ZY}\ll P_Z\otimes P_Y$; almost every conditional law is then
absolutely continuous with respect to the atomless $P_Y$ and is atomless.
\end{proof}

\subsubsection{An explicit quantizer for dominated sources}
\begin{lemma}[Oscillating partitions]\label{FP-lem:oscillating}
Let $p=(p_1,\ldots,p_M)$ have positive entries and partition $[0,1)$ into consecutive intervals $J_i$ of lengths $p_i$. For each $f\in L^1([0,1])$,
\[
 \int_0^1 f(t)\ind_{\{\{nt\}\in J_i\}}\,dt
 \longrightarrow p_i\int_0^1 f(t)\,dt.
\]
Here $\{x\}$ denotes the fractional part.
\end{lemma}
\begin{proof}
For an interval indicator $f=\ind_{[a,b]}$, complete periods have exact proportion $p_i$ and the at most two incomplete periods contribute error at most $2/n$. The result follows for interval step functions. Approximate an arbitrary $L^1$ function by such step functions; the error is bounded by the $L^1$ distance because the oscillating indicators are bounded by one.
\end{proof}

\begin{corollary}[Periodic useful-data quantizers]
\label{FP-thm:diffuse-capacity}
Suppose $P_Y$ is atomless and
\begin{equation}\label{FP-eq:dominated-joint}
 P_{ZY}\ll P_Z\otimes P_Y.
\end{equation}
Let $T(Y)$ be a measure-preserving Borel coordinate with uniform law
on $[0,1]$, invertible modulo null sets. If $J_i$ are consecutive
intervals of lengths $p_i$, the deterministic quantizers
\[
 B_n(Y)=i\quad\Longleftrightarrow\quad\{nT(Y)\}\in J_i
\]
have law $p$ and $I(B_n;Z)\to0$. Their conditional corrections from
Lemma~\ref{FP-lem:deterministic-repair} have exact privacy, the same
law $p$, and utility tending to $H(p)$.
\end{corollary}
\begin{proof}
Condition \eqref{FP-eq:dominated-joint} implies the existence of conditional densities
$f_z$ of $T(Y)$ with respect to Lebesgue measure. The oscillating-partition
lemma gives $q_{n,i}(z)=P(B_n=i\mid Z=z)\to p_i$ almost surely.
Dominated convergence yields $\E\TV(q_n(Z),p)\to0$.
The deterministic correction and the entropy inequality prove exact
privacy with limiting utility $H(p)$. Finally
$D(q_n(z)\|p)\to0$ and is bounded by $\log(1/\min_i p_i)$,
so $I(B_n;Z)\to0$.
\end{proof}

\subsubsection{An explicit Gaussian construction}
Let the conditional Gaussian variance be $s^2=\operatorname{Var}(Y\mid Z)>0$. For a period $h>0$ define
\[
 B_h(Y)=1+\left\lfloor M\{Y/h\}\right\rfloor,
 \qquad a=\frac{2\pi^2s^2}{h^2}.
\]

\begin{proposition}[A quantitative exactly private Gaussian output]\label{FP-prop:gaussian-rate}
There is a deterministic full-data encoder $U_h=F_h(Z,Y)$ such that $U_h$ is uniform on $[M]$, $U_h\perp Z$, and
\begin{align*}
 \Pp(U_h\ne B_h)&\le e_h:=\sum_{\ell=1}^{\infty}e^{-a\ell^2}
 \le\frac{e^{-a}}{1-e^{-3a}},\\
 I(U_h;Y)&\ge\log M-\omega_M(e_h).
\end{align*}
For fixed $M$, the information gap is $O(ae^{-a})$ as $h\downarrow0$.
\end{proposition}
\begin{proof}
The wrapped conditional Gaussian density on the unit circle has Fourier expansion
\[
 1+2\operatorname{Re}\sum_{\ell\ge1}
 e^{-a\ell^2}e^{2\pi i\ell(t-m(Z)/h)},
 \qquad m(Z)=\E[Y\mid Z].
\]
Its total variation distance to uniform measure is at most $\sum_{\ell\ge1}e^{-a\ell^2}$, uniformly in $Z$. Passing to the $M$ bins contracts total variation. Conditional maximal coupling, realized by Lemma~\ref{FP-lem:deterministic-repair}, gives the first assertion. Since $\ell^2\ge1+3(\ell-1)$ for $\ell\ge1$, the series has the stated geometric bound. The entropy estimate follows as above. Finally $h_b(t)=O(t\log(1/t))$ as $t\downarrow0$.
\end{proof}

\begin{theorem}[Gaussian endpoint and nonattainment]\label{FP-thm:gaussian-nonattainment}\label{R9-intro:privacy}
Let $(Z,Y)$ be a scalar jointly Gaussian pair with positive marginal variances and correlation $0<|\rho|<1$. For every $M\ge2$,
\[
 G_0^M(Z,Y)=\log M,
\]
but no admissible encoder attains this value. Under the output-perturbation restriction $Z-Y-U$, the exact-privacy value is zero. Under that restriction every positive leakage level has value $\log M$, attained by a deterministic finite quantizer.
\end{theorem}
\begin{proof}
The capacity identity follows from Theorem~\ref{R9-thm:conditional-capacity} or Proposition~\ref{FP-prop:gaussian-rate}. Standardize the Gaussian coordinates. For a square-integrable function $f(Y)$, its expansion in the orthonormal Hermite basis has coefficients $a_k$, and
\[
 \E[f(Y)\mid Z]=\sum_{k\ge0}\rho^k a_k H_k(Z)
 \quad\text{in }L^2.
\]
The identity for polynomials follows by taking conditional expectations in the generating function $e^{tY-t^2/2}$; polynomial density extends it to $L^2$. Since $\rho\ne0$, constant conditional expectation forces $f$ to be constant almost surely.

If $I(U;Y)=\log M$, both $H(U)=\log M$ and $H(U\mid Y)=0$ hold. Hence $U=f(Y)$ is uniform. Exact privacy would give $\E[\ind_{\{f(Y)=i\}}\mid Z]=1/M$, contradicting injectivity for every $i$. For a randomized output-perturbation encoder, apply the same argument to $f_i(Y)=P(U=i\mid Y)$. Exact privacy makes each $f_i$ constant, so $U\perp Y$. The final assertion is the positive-leakage statement of Theorem~\ref{R9-thm:conditional-capacity} with uniform $p$.
\end{proof}
If $\rho=0$, an equal-probability quantizer attains exact utility $\log M$. If $|\rho|=1$, useful data are measurable from private data and the exact-privacy value is zero. These boundary cases explain the nondegeneracy and correlation conditions.


\subsection{Conditional atoms and one infinite private release}\label{R10-sec:privacy-characterization}
The capacity theorem has an exact converse when every finite alphabet is tested. It also gives one uniform output with infinite useful information, by retaining a summable conditional entropy budget along all its binary prefixes. The mechanism observes both $Z$ and $Y$ throughout.

Write $G_0^M(Z,Y)=\sup\{I(U;Y):U\perp Z,\ |\mathcal U|\le M\}$, allowing arbitrary randomized encoders. The prescribed-law version fixes $\law(U)=p$.

\subsubsection{Conditional atoms give a quantitative converse}
For a regular conditional distribution define its largest atom
\[
 a(z)=\sup_y\Pp(Y=y\mid Z=z),
 \qquad
 b_M=\E(a(Z)-1/M)_+.
\]
The function $a$ is measurable. For example it is the limit of the largest masses of the refining dyadic cells after a Borel embedding into $[0,1]$; the limit for any probability law on that interval equals its largest atom mass. To see the latter assertion, the maxima decrease and always dominate every atom. A subsequence of cells whose masses approach the limiting maximum has centers converging in the compact interval. Every neighborhood of the limiting center eventually contains these cells, so that center has an atom at least as large as the limit.

\begin{theorem}[A gap forced by conditional atoms]
\label{SP-thm:atoms}
For every exactly private uniform $M$-valued output,
\[
 I(U;Y)\le\log M-b_M.
\]
Allowing any output law on at most $M$ labels gives
\begin{equation}\label{SP-eq:atom-gap}
 G_0^M(Z,Y)\le\log M-\chi(b_M),
 \qquad
 \chi(b)=\frac{(\sqrt{1+8b}-1)^2}{8}.
\end{equation}
These bounds allow randomized encoders observing both variables.
\end{theorem}
\begin{proof}
Given an admissible $U$, let $B=b(Y)$ be a Bayes predictor of $U$ and let $e=\Pp(U\ne B)$. A finite probability vector $r$ has entropy at least $-\log\max_i r_i\ge1-\max_i r_i$, so
\begin{equation}\label{SP-eq:bayes-entropy}
 e\le H(U\mid Y).
\end{equation}
For almost every $z$, some conditional label mass of $B$ is at least $a(z)$. Therefore, for $u_M$ the uniform vector,
\[
 \TV(\law(B\mid Z=z),u_M)\ge(a(z)-1/M)_+.
\]
If $U$ is uniform and independent of $Z$, the coupling $(B,U)$ conditional on $Z$ yields $b_M\le e$. Combine this with \eqref{SP-eq:bayes-entropy} and $I(U;Y)=\log M-H(U\mid Y)$.

For general $p=\law(U)$, exact privacy and the triangle inequality instead give $b_M\le e+\TV(p,u_M)$. Put $\Delta=\log M-I(U;Y)=D(p\Vert u_M)+H(U\mid Y)$. Pinsker's inequality and \eqref{SP-eq:bayes-entropy} imply $b_M\le\Delta+\sqrt{\Delta/2}$. Solving this quadratic inequality in $\sqrt\Delta$ gives \eqref{SP-eq:atom-gap}. Padding the alphabet with zero-probability labels covers outputs with fewer than $M$ values.
\end{proof}

\begin{corollary}[Exact all-alphabet characterization]
\label{SP-cor:privacy-iff}\label{SP-thm:privacy-main}
For full-data encoders observing $(Z,Y)$, the following are equivalent:
\begin{enumerate}[label=(\roman*),nosep]
\item $\law(Y\mid Z=z)$ is atomless for almost every $z$;
\item $G_0^M(Z,Y)=\log M$ for every finite $M\ge2$;
\item for every finite probability vector $p$, the supremum of $I(U;Y)$ under $U\perp Z$ and $\law(U)=p$ is $H(p)$.
\end{enumerate}
Under these conditions deterministic encoders approach each supremum.
\end{corollary}
\begin{proof}
Theorem~\ref{SP-thm:diffuse} proves (i)$\Rightarrow$(iii), and uniform $p$ gives (iii)$\Rightarrow$(ii). If conditional atoms occur on a set of positive probability, $\Pp(a(Z)>0)>0$, so $b_M>0$ for some finite $M$. Equation~\eqref{SP-eq:atom-gap} then contradicts (ii). Thus (ii)$\Rightarrow$(i).
\end{proof}
For one fixed alphabet size, absence of large conditional atoms is only a necessary condition in this converse. The exact characterization quantifies over every finite alphabet.

\subsubsection{One exactly private uniform output with infinite utility}
\label{SP-sec:infinite}
Each finite private prefix restricts the conditional source to an atomless set of fixed positive conditional mass. The next bit can therefore be corrected within every prefix while retaining exact privacy. Choosing summable conditional entropy costs makes all prefixes parts of one infinite release. Its useful information grows with the prefix length, while its total conditional-information cost stays bounded.

\begin{theorem}[An infinite release with bounded conditional information cost]
\label{SP-thm:infinite}\label{R10-intro:privacy}
Assume $\law(Y\mid Z=z)$ is atomless almost surely. For every $b>0$, there exists a deterministic binary sequence $U=(U_j)_{j\ge1}$, measurable from $(Z,Y)$, such that its coordinates are independent fair bits, the whole sequence is independent of $Z$, and
\begin{equation}\label{SP-eq:prefix-budget}
 H(U_1,\ldots,U_n\mid Y)\le b
 \quad\text{for every }n.
\end{equation}
Consequently
\[
 I(U;Y)=\infty,\qquad I(Z;U\mid Y)\le b.
\]
The binary expansion maps this sequence to a deterministic $\Unif[0,1]$ output with the same independence and information conclusions.
\end{theorem}
\begin{proof}
Choose positive $b_j$ with $\sum_jb_j\le b$. Construct the bits inductively. Suppose the finite prefix $A=(U_1,\ldots,U_{j-1})$ is uniform and independent of $Z$, and is a deterministic function of $(Z,Y)$. For each of its finitely many values $a$, the conditional law of $Y$ given $(Z,A=a)$ is an atomless restriction of $\law(Y\mid Z)$; the conditional probability of $A=a$ given $Z$ is the positive constant $2^{-(j-1)}$.

Apply Theorem~\ref{SP-thm:diffuse}, with a fair-bit output, to the source pair $(Z,Y)$ under its conditional law given $A=a$. Choose a deterministic next bit in that branch with
\[
 \Pp(U_j=0\mid Z,A=a)=\Pp(U_j=1\mid Z,A=a)=1/2,
 \qquad H(U_j\mid Y,A=a)\le b_j.
\]
Combine the branch encoders. The enlarged prefix is uniform and independent of $Z$. The conditional entropy chain rule gives
\[
 H(U_1,\ldots,U_n\mid Y)
 =\sum_{j=1}^n H(U_j\mid Y,U_1,\ldots,U_{j-1})
 \le\sum_{j=1}^n b_j\le b.
\]
The infinite sequence is measurable from $(Z,Y)$. Independence of every finite prefix from $Z$ implies independence of the generated sigma-field. Since each prefix is uniform,
\[
 I(U_1,\ldots,U_n;Y)
 =n\log2-H(U_1,\ldots,U_n\mid Y)\ge n\log2-b.
\]
Data processing proves $I(U;Y)=\infty$.

The prefixes are deterministic given $(Z,Y)$, and hence $I(Z;U_1,\ldots,U_n\mid Y)=H(U_1,\ldots,U_n\mid Y)\le b$. Conditional mutual information is continuous from below along increasing finite partitions of the second variable that generate its sigma-field. Applying this fact to the prefixes gives $I(Z;U\mid Y)\le b$. Finally $\sum_{j\ge1}2^{-j}U_j$ is uniform on $[0,1]$ and determines its binary sequence off a null dyadic set, so it preserves the asserted information quantities.
\end{proof}
The bound concerns Shannon information; the theorem imposes no regularity or computational guarantee on the encoder. Its prefix bound also forces the conditional law of the whole release given $Y$ to be atomic. Selecting its largest atom gives a predictor of the entire output, despite the uniform unconditional marginal.

\begin{corollary}[Prediction of the entire private output]\label{R15-cor:private-prediction}
The output in Theorem~\ref{SP-thm:infinite} can be predicted from $Y$ as a
whole sequence: there are measurable maps $g_j(Y)$ whose graphs support its
joint law, with
\begin{equation}\label{R15-eq:private-prediction}
 \Prb\{U=g_1(Y)\}\ge e^{-b},\qquad
 \Prb\{U\notin\{g_1(Y),\ldots,g_L(Y)\}\}
       \le\frac{b}{\log(L+1)}.
\end{equation}
It still has the atomless uniform product law, $U\perp Z$,
$I(U;Y)=\infty$ and $I(Z;U\mid Y)\le b$. For a nondegenerate correlated
Gaussian scalar pair, the supremum exact whole-output prediction probability
is one, even under any fixed positive budget; it is not attained.
\end{corollary}
\begin{proof}
Let $\kappa_y=\law(U\mid Y=y)$ and $[u]_n$ be the length-$n$ cylinder.
Monotone convergence in the prefix bound gives
\[
 \E[-\log\kappa_Y(\{U\})]
 =\lim_n H(U_1,\ldots,U_n\mid Y)\le b.
\]
Thus $\kappa_y$ is purely atomic almost everywhere and its discrete entropy
$H(\kappa_y)$ has expectation at most $b$. Its atoms can be enumerated
measurably in decreasing mass: on each binary subtree the largest mass is
the limit of the maxima over its finite levels; choose the lexicographically
first maximizing branch, delete that atom, and repeat. Compactness of the
binary tree and continuity of a measure from above identify each limit
with an atom. These branches define $g_j$ and exhaust the mass.

The largest atom has mass at least $e^{-H(\kappa_y)}$, so Jensen gives
the first bound. After $L$ largest atoms each remaining mass is at most
$1/(L+1)$; its entropy contribution is at least its mass times
$\log(L+1)$. Integrating gives the second bound. Letting $b\downarrow0$
gives the Gaussian supremum. Attainment would make $U$ a function of $Y$
independent of $Z$; the Hermite conditional-expectation calculation in
Corollary~\ref{SP-cor:gaussian-budget} forces every bounded function of $U$
to be constant, contrary to its uniform law.
\end{proof}
The predictor itself is not asserted to be private or uniform. The equality
event in \eqref{R15-eq:private-prediction} concerns the complete release,
rather than any fixed finite prefix.

\begin{corollary}[A discontinuity in the Gaussian conditional-information budget]
\label{SP-cor:gaussian-budget}
Let $(Z,Y)$ be a nondegenerate jointly Gaussian scalar pair with correlation $0<|\rho|<1$. Define
\[
 \mathcal U(b)=\sup\{I(U;Y):\ U\perp Z,\ I(Z;U\mid Y)\le b\},
\]
where arbitrary standard Borel outputs and full-data mechanisms are allowed. Then
\[
 \mathcal U(0)=0,
 \qquad
 \mathcal U(b)=\infty\quad\text{for every }b>0.
\]
For every positive $b$, the infinite value is attained by a deterministic uniform $[0,1]$ output.
\end{corollary}
\begin{proof}
For $b>0$, the Gaussian conditional laws are atomless and Theorem~\ref{SP-thm:infinite} applies. At $b=0$, conditional independence gives the Markov chain $Z-Y-U$. For a bounded measurable test function $f$ on the output space, put $h(Y)=\E[f(U)\mid Y]$. Exact privacy implies $\E[h(Y)\mid Z]=\E f(U)$. After standardizing the Gaussian pair, expand $h$ in the Hermite basis of $L^2$ of the standard normal law. Conditional expectation multiplies the degree-$k$ Hermite coefficient by $\rho^k$. Since every $\rho^k$ is nonzero, the displayed identity forces every nonconstant coefficient of $h$ to vanish. Thus $\E[f(U)\mid Y]=\E f(U)$ for every bounded $f$, which gives $U\perp Y$ and $I(U;Y)=0$.
\end{proof}
The zero-budget condition gives the output-perturbation Markov chain, whose Gaussian zero-utility conclusion is established in \cite{RG2021}. The same paper proves unbounded full-observation utility. Here the intervening budget has its exact threshold: every positive $b$ admits one deterministic uniform infinite-utility output, with whole-output prediction probability at least $e^{-b}$.

\subsection{Exact privacy under smooth changes of the data law}
\label{R22-sec:privacy-instability}
The Gaussian profile fixes the data law and varies the conditional-information budget. We now fix that budget at zero and vary the data law, preserving its covariance. The zero-budget requirement makes the release a channel from $Y$ alone. An arbitrarily fine uniform component of $Y$ can nevertheless encode an exactly private continuous output. The construction below retains analytic densities, uniform log-concavity and convergence of the entire potential Hessian to the Gaussian precision matrix.

For a joint law $P$ of $(Z,Y)$ write
\[
 \mathcal U_P(b)=\sup\{I(U;Y):U\perp Z,\ I(Z;U\mid Y)\le b\},
\]
with the same full-data observation model and arbitrary standard Borel outputs as in Corollary~\ref{SP-cor:gaussian-budget}. Relative entropy is denoted by $D(P\Vert Q)$.

\begin{theorem}[Smooth covariance-preserving privacy instability]
\label{R22-thm:privacy-instability}
Fix $0<|\rho|<1$ and let $P_0=N(0,\Sigma_\rho)$, where
$\Sigma_\rho=\left(\begin{smallmatrix}1&\rho\\\rho&1\end{smallmatrix}\right)$.
There are centered laws $P_a$, $0<a\le1/2$, with covariance exactly $\Sigma_\rho$ and private marginal $Z\sim N(0,1)$, whose densities $p_a=e^{-V_a}$ are positive and real analytic and satisfy
\begin{equation}\label{R22-eq:privacy-curvature}
 (1-a^2)\Sigma_\rho^{-1}\preceq\nabla^2V_a
 \preceq(1+a^2)\Sigma_\rho^{-1}\qquad\hbox{everywhere}.
\end{equation}
As $a\downarrow0$,
\begin{equation}\label{R22-eq:privacy-convergence}
 W_2(P_a,P_0)\longrightarrow0,\qquad
 D(P_a\Vert P_0)+D(P_0\Vert P_a)\longrightarrow0.
\end{equation}
Nevertheless,
\begin{equation}\label{R22-eq:privacy-instability}
 \mathcal U_{P_0}(0)=0,\qquad
 \mathcal U_{P_a}(0)=\infty\quad(a>0).
\end{equation}
For every $a>0$ a deterministic function of $Y$ alone is uniform on $[0,1)$, independent of $Z$, and attains the infinite value. Its whole-output prediction probability is one. Each $P_a$ and $P_0$ has infinite utility at every positive budget.
\end{theorem}
\begin{proof}
Let $Z,G$ be independent standard normals, let $T_a$ be independent and uniform on $[-a,a]$, and put
\begin{equation}\label{R22-eq:privacy-construction}
 s_a^2=1-a^2/3,\qquad E_a=s_aG+T_a,\qquad
 Y_a=\rho Z+cE_a,\quad c=\sqrt{1-\rho^2}.
\end{equation}
The mean and covariance are as asserted. With $\{x\}=x-\lfloor x\rfloor$, define
\begin{equation}\label{R22-eq:modulo-release}
 U_a=\left\{\frac{Y_a}{2ac}\right\}.
\end{equation}
Conditional on $Z,G$, the variable $T_a/(2a)$ is uniform on an interval of length one. Translation modulo one preserves uniform measure, so $U_a$ is uniform and independent of $Z$. It is a function of $Y_a$, giving $I(Z;U_a\mid Y_a)=0$. The graph of this function has zero measure under $\law(Y_a)\otimes\law(U_a)$ because $U_a$ is atomless. Thus $I(U_a;Y_a)=\infty$. The value at $P_0$ is Corollary~\ref{SP-cor:gaussian-budget}.

Let $q_a$ be the density of $E_a$ and $v_a=-\log q_a$. Gaussian convolution of a compactly supported law is positive and real analytic. Differentiation under the integral gives
\begin{equation}\label{R22-eq:posterior-curvature}
 v_a''(x)=s_a^{-2}-s_a^{-4}\operatorname{Var}(T_a\mid E_a=x).
\end{equation}
Since $|T_a|\le a$,
\[
 \frac{1-4a^2/3}{(1-a^2/3)^2}\le v_a''(x)
 \le\frac1{1-a^2/3}.
\]
For $0\le a^2\le1/4$ the lower bound is at least $1-a^2$ and the upper bound at most $1+a^2$, by multiplication by their positive denominators. The joint potential is $z^2/2+v_a((y-\rho z)/c)$ up to a constant. The linear change $(z,y)\mapsto(z,(y-\rho z)/c)$ turns its Hessian into $\operatorname{diag}(1,v_a'')$ and the Gaussian precision into the identity. This proves \eqref{R22-eq:privacy-curvature}.

Using the same $Z,G$ for $P_a,P_0$ gives
\[
 W_2(P_a,P_0)^2\le c^2\bigl((s_a-1)^2+a^2/3\bigr).
\]
Convexity of relative entropy under mixing and $s_a^2+\E T_a^2=1$ give
\[
 D(q_a\Vert\gamma_1)\le-\tfrac12\log s_a^2.
\]
If $\phi_s$ is the $N(0,s^2)$ density, then
\[
 q_a(x)=\phi_{s_a}(x)\E\exp\!\left(\frac{xT_a}{s_a^2}
                       -\frac{T_a^2}{2s_a^2}\right)
 \ge\phi_{s_a}(x)e^{-a^2/(2s_a^2)},
\]
where Jensen uses $\E T_a=0$. Consequently
\[
 D(\gamma_1\Vert q_a)
 \le\tfrac12(s_a^{-2}-1+\log s_a^2)+\frac{a^2}{2s_a^2}.
\]
The two joint relative entropies equal these noise relative entropies under the fixed invertible change of coordinates. Both tend to zero, proving \eqref{R22-eq:privacy-convergence}; total variation convergence follows as well. The conditional laws of $Y_a$ given $Z$ are atomless, so Theorem~\ref{SP-thm:infinite} also applies at every positive budget.
\end{proof}

The masking step is the classical uniform-dither mechanism: its Fourier cancellation is the one used in independence criteria for quantization noise~\cite{R22-Schuchman}. In this theorem it produces an exactly private release while the full covariance and the stated analytic and curvature constraints are retained. The feature that persists under these constraints is visible in
\begin{equation}\label{R22-eq:sinc}
 \widehat q_a(t)=e^{-s_a^2t^2/2}\frac{\sin(at)}{at}.
\end{equation}
The nonzero Fourier zeros move to infinity as $a\downarrow0$. The release in \eqref{R22-eq:modulo-release} reads the corresponding periodic structure. The next criterion explains precisely what is lost when those zeros disappear.

\begin{theorem}[A Fourier criterion for zero-budget utility]
\label{R22-thm:privacy-Fourier}
Let $Z$ have full support on $\R$, let $E$ be independent of $Z$ with an everywhere-positive continuous density $q$, and put $Y=Z+E$. Then
\begin{equation}\label{R22-eq:Fourier-criterion}
 \mathcal U_{\law(Z,Y)}(0)=0
 \quad\Longleftrightarrow\quad
 \widehat q(t)\ne0\quad\hbox{for every real }t.
\end{equation}
A real Fourier zero gives strictly positive zero-budget utility. Every positive budget has infinite utility, attained by one deterministic uniform full-data release with the prediction guarantee of Corollary~\ref{R15-cor:private-prediction}.
\end{theorem}
\begin{proof}
At budget zero, $Z-Y-U$ is a Markov chain. For any bounded measurable output test $f$, put $h(y)=\E[f(U)\mid Y=y]-\E f(U)$. Privacy gives $h*\check q=0$ for $\law(Z)$-almost every argument, where $\check q(x)=q(-x)$. Translation continuity in $L^1$ makes this convolution continuous; full support makes it zero everywhere. If $\widehat q$ has no real zeros, Wiener's translation-density theorem~\cite{R22-Wiener} says that the translates of $\check q$ span a dense subspace of $L^1(\R)$. The bounded function $h$ annihilates this space, so $h=0$ Lebesgue-almost everywhere. The useful-data law has a positive density, hence $U\perp Y$.

If $\widehat q(t_0)=0$, then $t_0\ne0$ and the binary channel
\[
 \Prb\{U=1\mid Y=y\}=\frac{1+\cos(t_0y)}2
\]
is private: its conditional probability given $Z=z$ is $1/2$. The displayed probability is nonconstant on the positive density of $Y$, so $I(U;Y)>0$. This proves the converse. Conditional atomlessness and Theorem~\ref{SP-thm:infinite} prove the positive-budget assertion. Thus the harmonic-analytic input is classical bounded completeness; the criterion identifies its exact privacy consequence.
\end{proof}

\begin{corollary}[Every finite logarithmic zero-budget value]
\label{R22-cor:privacy-logm}
For every integer $m\ge2$ and every $0<|\rho|<1$ there are centered analytic uniformly log-concave laws with covariance $\Sigma_\rho$, converging to $P_0$ in all the senses of \eqref{R22-eq:privacy-convergence} and with their potential Hessians converging uniformly to $\Sigma_\rho^{-1}$, for which
\[
 \mathcal U_P(0)=\log m\quad\hbox{and}\quad
 \mathcal U_P(b)=\infty\quad(b>0).
\]
The zero-budget supremum is over every standard Borel output and is attained by a deterministic $m$-valued output.
\end{corollary}
\begin{proof}
Take $J$ uniform on $\{0,\ldots,m-1\}$ and independent of $Z,G$, and let
\[
 T=h(J-(m-1)/2),\qquad v=h^2(m^2-1)/12,
 \qquad E=\sqrt{1-v}\,G+T,
 \qquad Y=\rho Z+cE.
\]
For sufficiently small $h>0$, the posterior-variance identity \eqref{R22-eq:posterior-curvature}, with $s^2=1-v$ and $|T|\le h(m-1)/2$, gives uniform log-concavity and Hessian convergence. The same coupling and entropy bounds give both convergence statements, with $a=h(m-1)/2$ and $\E T^2=v$. The release $\lfloor Y/(ch)\rfloor\pmod m$ is uniform conditional on $Z,G$, since adding $J$ cycles through all residues. It is private, deterministic from $Y$, and has information $\log m$.

For the upper bound rescale to $\widetilde Y=Y/c$ and $\widetilde Z=\rho Z/c$, whose private law has full support. Let $K_y$ be any admissible zero-budget output kernel with marginal $p$. Privacy, followed by the nonvanishing-Fourier argument for the Gaussian convolution, gives
\begin{equation}\label{R22-eq:cyclic-privacy}
 \frac1m\sum_{j=0}^{m-1}K_{y+h(j-(m-1)/2)}=p
 \quad\hbox{for Lebesgue-almost every }y.
\end{equation}
Apply Gaussian deconvolution first to bounded output tests. A countable determining class on the standard Borel output space gives the measure identity simultaneously. Every nonnegative summand is bounded above by $mp$, so $K_y\le mp$ almost everywhere and $D(K_y\Vert p)\le\log m$. Integrating proves $I(U;Y)\le\log m$ for every output. The positive-budget claim again follows from conditional atomlessness.
\end{proof}

\begin{corollary}[Failure of both semicontinuity directions]
\label{R22-cor:privacy-tilt}
At fixed nonzero correlation and fixed covariance, zero-budget utility is neither upper nor lower semicontinuous in $W_2$ or total variation within positive analytic uniformly log-concave data laws.
\end{corollary}
\begin{proof}
Theorem~\ref{R22-thm:privacy-instability} gives failure of upper semicontinuity at the Gaussian. For the other direction fix small $a>0$ and give $T_{a,t}$ density proportional to $e^{tu}$ on $[-a,a]$. Center it and choose the independent Gaussian variance to be $1-\operatorname{Var}(T_{a,t})$. Before centering, its characteristic function is
\[
 \frac{t}{t+i\xi}\,
 \frac{\sinh(a(t+i\xi))}{\sinh(at)}\qquad(t\ne0).
\]
It has no real zero: the zeros of $\sinh$ have real part zero, whereas $at\ne0$. Centering only changes its phase. Theorem~\ref{R22-thm:privacy-Fourier} therefore gives zero utility for every $t\ne0$. As $t\to0$ these laws converge in $W_2$ and total variation to the uniform-mask example, whose utility is infinite. The Gaussian variance stays positive, and the posterior variance of a random variable in an interval of length $2a$ is at most $a^2$. Formula~\eqref{R22-eq:posterior-curvature} therefore gives a common positive curvature lower bound for $a\le1/2$. Means and covariances are fixed by the centering and variance choice. These examples give failure of lower semicontinuity.
\end{proof}

\subsubsection{Greatest and maximal private experiments}
Infinite mutual information leaves open which decision problems a release solves. Here the zero-budget class consists of useful-data-only channels $Z-Y-U$. A release is \emph{Blackwell-greatest} in this class if every such channel is a garbling of it: its output can be transformed by a Markov kernel independent of the useful variable $Y$ to reproduce the other channel. This property gives one release for every bounded decision problem about $Y$. A channel is \emph{maximal} if every private channel dominating it is Blackwell-equivalent to it. The uniform mask permits a complete channel factorization.

\begin{theorem}[All private channels factor through the modulo release]
\label{R23-thm:private-greatest}
Let $Z$ have full support on $\R$. Let $N$ have an everywhere-positive continuous density with no real Fourier zeros, and let $T$ be uniform on an interval of length $h>0$. Assume $Z,N,T$ are independent and put $Y=Z+N+T$. A standard Borel channel $K_y=\law(U\mid Y=y)$ is exactly private if and only if there is a Markov kernel $L$ on $[0,1)$ such that
\begin{equation}\label{R23-eq:private-factor}
 K_y=L_{\{y/h\}}\qquad\hbox{for Lebesgue-almost every }y.
\end{equation}
Thus $R=\{Y/h\}$ is uniform, independent of $Z$, and Blackwell-greatest among all useful-data-only perfectly private channels. In particular, for every law $P_a$ of Theorem~\ref{R22-thm:privacy-instability}, the release $U_a$ in \eqref{R22-eq:modulo-release} has this property, with every stated covariance, curvature and convergence constraint retained.
\end{theorem}
\begin{proof}
Fix a bounded output test $f$ and let $k(y)=\int f\,dK_y$ and $p_f=\E f(U)$. Write the uniform interval as $[b,b+h]$. Privacy and the full support of $Z$ imply
\[
 \E k(z+N+T)=p_f\qquad\hbox{for every }z.
\]
Indeed the left side is continuous by translation continuity of the noise density in $L^1$, and equality initially holds for $\law(Z)$-almost every $z$. Let $a(y)=h^{-1}\int_b^{b+h}k(y+t)\,dt$. Then $(a-p_f)*\check q_N=0$. Wiener's translation theorem~\cite{R22-Wiener} gives $a=p_f$ almost everywhere. The moving average is locally absolutely continuous, with
\[
 a'(y)=\frac{k(y+b+h)-k(y+b)}h
\]
almost everywhere. Hence $k(y+h)=k(y)$ almost everywhere. Apply this argument to a countable determining class on the standard Borel output space. Outside one Lebesgue null set, the kernels $K_y$ agree under every integer translation by $h$. Restricting a measurable version to $[0,h)$ gives $L$ and \eqref{R23-eq:private-factor}. The positive density of $Y$ makes these identities valid for its channel.

Conversely, conditional on $Z,N$, the fractional part $R$ is uniform. Any channel of the form \eqref{R23-eq:private-factor} is therefore independent of $Z$. The same formula is the required garbling of $R$. For $P_a$, rescale to $Y_a/c=(\rho/c)Z+s_aG+T_a$. The private variable $(\rho/c)Z$ has full support, the Gaussian convolution has no real Fourier zeros, and $h=2a$. This gives precisely $U_a$.
\end{proof}

The forward cancellation is the classical uniform-dither identity~\cite{R22-Schuchman}. Deconvolution and the moving-average derivative give the converse: privacy forces the whole output kernel to be periodic. Consequently the same deterministic output serves every zero-budget private decision problem. The discrete masks with finite capacity have a different decision order even under the same smooth covariance-preserving approximation.

\begin{theorem}[Incomparable maximal private releases near a Gaussian]
\label{R23-thm:private-maximal}
For every integer $m\ge2$, the analytic uniformly log-concave laws constructed in Corollary~\ref{R22-cor:privacy-logm} admit two deterministic, uniform $m$-valued private releases that attain the capacity $\log m$ and are incomparable maximal elements of the Blackwell order. They have no common private upper bound; in particular there is no Blackwell-greatest useful-data-only private channel. All covariance and convergence assertions of that corollary remain valid.
\end{theorem}
\begin{proof}
Use its notation and put
\[
 \begin{aligned}
 W&=Y/(ch),&\quad R_1&=\lfloor W\rfloor\pmod m,\\
 S&=\{W\},& B&=\mathbf1_{\{S\ge1/2\}},\qquad
 R_2=(R_1+B)\pmod m.
 \end{aligned}
\]
Conditional on $Z,G$, the mask $J$ cycles through the $m$ residues, while $S$ and $B$ are fixed. Both $R_1$ and $R_2$ are therefore uniform and private. Moreover $R_1$ is independent of $B$. The positive density of $Y$ gives $0<\Prb(B=1)<1$, and $B$ is recovered from $R_2-R_1$ modulo $m$. Hence
\begin{equation}\label{R23-eq:private-no-upper}
 H(R_1,R_2)=\log m+H(B)>\log m.
\end{equation}

Suppose a private channel $V$ dominates a deterministic function $R_i$ of $Y$. A garbling from $V$ to $R_i$ has conditional law given $Y$ equal to a point mass. Its conditional law given $V$ must consequently be a point mass almost surely, so $R_i$ is a measurable function of $V$. If $V$ dominates both releases, data processing and \eqref{R23-eq:private-no-upper} imply
\[
 I(V;Y)\ge H(R_1,R_2)>\log m,
\]
contradicting the capacity bound, which holds for every standard Borel output. Thus there is no common private upper bound, and the two releases are incomparable.

Finally, if a private $V$ dominates $R_i$, the same measurability and the chain rule give
\[
 \log m\ge I(V;Y)=\log m+I(V;Y\mid R_i).
\]
The conditional information vanishes, so $V$ is a garbling of $R_i$. Each $R_i$ is therefore maximal.
\end{proof}

\subsubsection{Greatest useful-data-only experiments and posterior simplices}
For finite useful data the same question has a geometric answer. The feasible posterior polytope and its vertex reduction for mutual information appear in Rassouli--G\"und\"uz~\cite[Lemma~1 and Theorem~1]{RG2021}. Combining that representation with Blackwell comparison~\cite{Blackwell,Strassen} identifies exactly when one private experiment serves every decision problem.

\begin{proposition}[Greatest private experiment and posterior simplex]
\label{R23-prop:privacy-simplex}
Let $Y\in\{1,\ldots,q\}$ have strictly positive probabilities $p_i$, and let $Z$ be standard Borel. Write $\ell_i=\law(Z\mid Y=i)$ and $\ell=\law(Z)$. Define the polytope
\begin{equation}\label{R23-eq:private-polytope}
 \mathcal C=\left\{r\in\Delta_q:\ \sum_{i=1}^q r_i\ell_i=\ell\right\}.
\end{equation}
There is a Blackwell-greatest useful-data-only private channel if and only if $\mathcal C$ is a simplex in its affine hull, including the one-point simplex. More generally, the maximal private experiments, up to Blackwell equivalence, are precisely the laws supported on the vertices of $\mathcal C$ with barycenter $p$. When $\mathcal C$ is a simplex, the greatest experiment has the unique such law.
\end{proposition}
\begin{proof}
Privacy and $Z-Y-U$ imply that each posterior $r=\law(Y\mid U)$ lies in $\mathcal C$, while averaging posteriors gives $p$. Conversely, every probability $\eta$ on $\mathcal C$ with mean $p$ defines a private channel by
\[
 \Prb\{U\in dr\mid Y=i\}=\frac{r_i}{p_i}\,\eta(dr).
\]
The constraint in \eqref{R23-eq:private-polytope} is a finite-dimensional linear constraint: the finite set of signed measures $\ell_i$ spans a finite-dimensional space. Thus $\mathcal C$ is a polytope. Positivity of $p_i$ puts $p$ in its relative interior.

Every point of $\mathcal C$ has a vertex decomposition, chosen measurably, for example by a finite triangulation. Refining each posterior into those vertices produces a dominating private experiment. If the original posterior law charges nonvertices, the mean squared norm strictly increases, so the refinement is not Blackwell-equivalent. A vertex-supported law is maximal: in any martingale refinement, a vertex can be the mean of a distribution on $\mathcal C$ only when that distribution is concentrated at the vertex. These observations prove the stated classification by the posterior characterization of Blackwell order.

If $\mathcal C$ is a simplex, its barycentric coordinates are unique, and every refinement just described has the same vertex law of mean $p$. That law dominates every private experiment. If $\mathcal C$ is not a simplex, its vertices are affinely dependent. The relative interior point $p$ has a representation assigning positive weight to every vertex. Perturbing those weights in a nonzero affine dependence gives two distinct vertex laws with mean $p$. Both are maximal, and no common private upper bound can dominate them, by the vertex argument. This excludes a greatest private experiment.
\end{proof}

The same mechanism also gives a density statement without curvature assumptions. If a centered noise $E$ has a density and finite variance $v>0$, then
$\sqrt{1-a^2/(3v)}E+T_a$ preserves its variance and converges to it in $W_2$ and total variation as $a\downarrow0$. Translation and dilation continuity in $L^1$ give the latter convergence. The modulo release has infinite zero-budget utility for every $a>0$. The analytic and uniform-curvature conclusions require the Gaussian component used in Theorem~\ref{R22-thm:privacy-instability}.


\begin{example}[A singular source with infinite source mutual information]
\label{SP-ex:singular}
Let $Z$ be uniform on the unit circle $\R/\mathbb Z$, let $C$ have the standard atomless Cantor distribution independently of $Z$, and put $Y=Z+C\pmod1$. Then $Y$ is uniform and $\law(Y\mid Z=z)$ is an atomless translate of the Cantor law. Theorem~\ref{SP-thm:diffuse} applies. The joint law is concentrated on the set $\{(z,y):y-z\in\operatorname{supp}C\pmod1\}$, which has zero measure under the product of the two uniform marginals. Hence the joint law is singular and $I(Z;Y)=\infty$. No finiteness assumption on $I(Z;Y)$ is needed for the privacy theorem.
\end{example}


\section{Greatest full-data private releases}\label{R27-sec:full-data}

The preceding section constructs exact privacy capacities and, in selected useful-data-only models, a greatest release. We now classify when one release is best for every decision about the useful variable $Y$ after the encoder is allowed to observe the full data. A full-data mechanism observes $(Z,Y)$ and may use auxiliary randomness; privacy remains $U\perp Z$. All source and output spaces in this section are standard Borel.

A release $V$ dominates $U$ in the Blackwell order if a Markov kernel $K$ satisfies
\begin{equation}\label{fdp:eq:garbling}
 \law(Y,U)(dy,du)=\int \law(Y,V)(dy,dv)K(v,du).
\end{equation}
A \emph{greatest private release} dominates every perfectly private full-data release. The useful-data-only class imposes the additional Markov constraint $Z-Y-U$. Its greatest modulo release and finite posterior-polytope criterion are Theorem~\ref{R23-thm:private-greatest} and Proposition~\ref{R23-prop:privacy-simplex}. The classification below concerns the larger full-data class.

Write $\rho=\law(Z)$ and let
\[
 \mu_z=\law(Y\mid Z=z).
\]
All statements about these conditional laws hold $\rho$-almost surely. Let
\[
 D=\{z:\mu_z\text{ is a point mass}\},\qquad N=D^c.
\]
The point-mass map is Borel on a standard Borel space, so these sets and the location $d(z)$ with $\mu_z=\delta_{d(z)}$ on $D$ can be chosen measurably.

\begin{theorem}[Greatest full-data private releases]\label{fdp:thm:classification}
Let $(Z,Y)$ take values in standard Borel spaces. A Blackwell-greatest perfectly private full-data release exists if and only if at least one of the following conditions holds:
\begin{enumerate}[label=(\roman*),nosep]
\item There is a probability law $\nu$ such that $\mu_z=\nu$ for almost every $z\in N$.
\item There are two distinct points $a,b$ such that $\mu_z(\{a,b\})=1$ for almost every $z\in N$.
\end{enumerate}
The conditional point masses on $D$ may have arbitrary locations and arbitrary total probability. If $\rho(N)=0$, every private release is independent of $Y$ and a constant release is greatest.
\end{theorem}

The theorem gives explicit mechanisms in both cases. The first uses one common realization of the common nondeterministic conditional law. The second uses one uniform variable to couple all binary conditional laws monotonically. These constructions appear in Subsection~\ref{fdp:sec:sufficiency}.

\begin{corollary}[Conditionally atomless data]\label{fdp:cor:atomless}
Suppose $\mu_z$ is atomless for almost every $z$. A greatest full-data private release exists if and only if $Y\perp Z$. More generally, the same conclusion holds if almost every $\mu_z$ is not concentrated on a set of at most two points.
\end{corollary}
\begin{proof}
There are no deterministic conditional laws. The two-point alternative in Theorem~\ref{fdp:thm:classification} is excluded. The remaining alternative says exactly that $\mu_z$ is almost surely constant, which is independence. At independence, the release $U=Y$ is private and dominates every experiment about $Y$.
\end{proof}

The obstruction is already visible in quadratic prediction. Fix a countable algebra $\mathcal A$ that generates the Borel sets of the useful space and separates its points. Let
\begin{equation}\label{fdp:eq:ternary-family}
 \mathcal F=\{\ind_B-\ind_C:B,C\in\mathcal A,\ B\cap C=\varnothing\}.
\end{equation}
Every member takes values in $\{-1,0,1\}$. For a bounded observable $f$, set
\[
 S_f(U)=\E\bigl[\E[f(Y)\mid U]^2\bigr].
\]
Maximizing $S_f$ is equivalent to minimizing squared prediction error for $f(Y)$.

\begin{theorem}[Quadratic prediction detects the greatest experiment]\label{fdp:thm:quadratic}
A private release is Blackwell-greatest if and only if it maximizes $S_f$ over all private full-data releases for every $f\in\mathcal F$. Consequently, a common optimizer for this countable family exists exactly in the cases of Theorem~\ref{fdp:thm:classification}. If $Y$ has $n$ possible values, it is enough to test the $3^n$ ternary observables.
\end{theorem}

The decisive necessity argument has two stages. Each scalar task can be optimized by aligning conditional quantiles. If one release solves all tasks, its pairwise conditional couplings must therefore be optimal for every ternary valuation of the useful states. Two support points of such a coupling can cross neither under one valuation nor under another. This forces the coupling onto a diagonal, a single row, a single column, or a triangular two-point support. That small support classification determines the full standard Borel theorem.

\subsection{A Gaussian obstruction witnessed by two polynomial tasks}\label{fdp:sec:gaussian}

Let $Z,E$ be independent standard Gaussian variables and set
\[
 Y=rZ+sE,\qquad 0<|r|<1,\qquad s=\sqrt{1-r^2}.
\]
Consider the deterministic full-data releases
\begin{equation}\label{fdp:eq:two-gaussian}
 U_1=E=\frac{Y-rZ}{s},\qquad
 U_2=\operatorname{sign}(Z)E.
\end{equation}
Both are standard Gaussian and independent of $Z$.

For a square-integrable target $f(Y)$, write
\[
 V_f(U)=\Var\bigl(\E[f(Y)\mid U]\bigr).
\]

\begin{theorem}[Gaussian private releases with no common upper bound]\label{fdp:thm:gaussian}
The releases in \eqref{fdp:eq:two-gaussian} satisfy
\begin{align*}
 V_Y(U_1)&=s^2,& V_Y(U_2)&=0,\\
 V_{Y^2}(U_1)&=2s^4,&
 V_{Y^2}(U_2)&=2s^4+\frac8\pi r^2s^2.
\end{align*}
Every private full-data release has $V_Y\le s^2$. Every release attaining this bound is Blackwell-equivalent, as an experiment about $Y$, to $U_1$. Consequently, $U_1$ and $U_2$ are incomparable and have no common private upper bound.
\end{theorem}
\begin{proof}
For a private $W$, independence gives $\E[Z\mid W]=0$, so
\[
 \E[Y\mid W]=s\E[E\mid W],\qquad
 V_Y(W)\le s^2.
\]
Equality forces $E$ to be measurable with respect to $W$. Because $Z\perp W$, the posterior law of $Y$ given $W$ is then $N(sE,r^2)$, exactly the posterior determined by $U_1=E$. Thus $W$ is Blackwell-equivalent to $U_1$.

For $U_1$, the conditional first and second moments are
\[
 \E[Y\mid U_1]=sU_1,\qquad
 \E[Y^2\mid U_1]=r^2+s^2U_1^2.
\]
For $U_2$, the pair $(Z,U_2)$ is independent and $E=\operatorname{sign}(Z)U_2$. Therefore
\[
 \E[Y\mid U_2]=0,\qquad
 \E[Y^2\mid U_2]=r^2+2rs\sqrt{2/\pi}\,U_2+s^2U_2^2.
\]
The odd and centered even Gaussian terms are orthogonal, which gives the stated variances.

A common private upper bound would dominate $U_1$ and therefore attain the maximal linear utility. It would be equivalent to $U_1$, but it could not dominate $U_2$, whose quadratic utility is strictly larger. The same strict inequalities can be witnessed by bounded truncations of the targets.
\end{proof}

\begin{corollary}[Infinite information and the Gaussian innovation]\label{fdp:cor:infinite-info}
Every perfectly private full-data release $V$ for the Gaussian source with $I(V;Y)=\infty$ is Blackwell-incomparable with the innovation $U_1$. The two have no common private upper bound. In particular, this applies to the infinite-information releases with arbitrarily small positive conditional-information budget in Theorem~\ref{SP-thm:infinite}.
\end{corollary}
\begin{proof}
The Gaussian innovation has $I(U_1;Y)=-\log|r|<\infty$, so it cannot dominate $V$. Any private release dominating $U_1$ is equivalent to $U_1$ by Theorem~\ref{fdp:thm:gaussian}, and therefore has the same finite mutual information. Such a release cannot be $V$ and cannot also dominate $V$.
\end{proof}

The obstruction is quantitative. Approaching the optimal linear prediction score constrains nonlinear prediction even before exact optimality is reached.

\begin{proposition}[A robust Gaussian utility gap]\label{fdp:prop:gaussian-gap}
For every private full-data release $W$, let
\[
 \delta=s^2-V_Y(W)\ge0.
\]
Then
\begin{equation}\label{fdp:eq:gaussian-gap}
 V_{Y^2}(W)\le\bigl(\sqrt2\,s^2+2|r|\sqrt\delta\bigr)^2.
\end{equation}
In particular, achieving at least the quadratic utility of $U_2$ requires
\begin{equation}\label{fdp:eq:positive-deficit}
 \delta\ge
 \left(\frac{\sqrt{2s^4+(8/\pi)r^2s^2}-\sqrt2\,s^2}{2|r|}\right)^2>0.
\end{equation}
\end{proposition}
\begin{proof}
Set $m=\E[E\mid W]$ and $v=\E[(E-m)^2\mid W]$. Then $\E v=\delta/s^2$. Since $Z\perp W$,
\[
 \E[ZE\mid W]=\E[Z(E-m)\mid W].
\]
Conditional Cauchy--Schwarz and $\E[Z^2\mid W]=1$ give
\[
 \|\E[ZE\mid W]\|_2\le\sqrt{\E v}=\sqrt\delta/s.
\]
Also, conditional expectation is an $L^2$ contraction, so
\[
 \|\E[E^2\mid W]-1\|_2\le\|E^2-1\|_2=\sqrt2.
\]
The identity
\[
 \E[Y^2\mid W]-1
 =2rs\E[ZE\mid W]+s^2\bigl(\E[E^2\mid W]-1\bigr)
\]
and the triangle inequality prove \eqref{fdp:eq:gaussian-gap}. Solving it for $\delta$ gives \eqref{fdp:eq:positive-deficit}.
\end{proof}

To determine when every task can have a common optimizer, we first identify the exact optimum for a single scalar observable. Equality in that calculation will constrain every pair of conditional useful laws.

\subsection{Every scalar task has an exact convex-order optimum}\label{fdp:sec:scalar}

A single task has a complete solution even when a common greatest experiment does not exist. Let $f(Y)$ be bounded and let $Q_z(t)$, $0<t<1$, be a nondecreasing quantile of its conditional law given $Z=z$. Define
\begin{equation}\label{fdp:eq:quantile-barycenter}
 m_f(t)=\int Q_z(t)\rho(dz).
\end{equation}
The map $(z,t)\mapsto Q_z(t)$ can be chosen measurable using rational thresholds. The function $m_f$ is nondecreasing and bounded.

\begin{theorem}[The greatest posterior mean for a scalar task]\label{fdp:thm:scalar}
For every perfectly private full-data release $U$,
\begin{equation}\label{fdp:eq:scalar-cx}
 \law\bigl(\E[f(Y)\mid U]\bigr)\preceq_{\rm cx}\law(m_f(T)),\qquad T\sim\mathrm{Unif}(0,1).
\end{equation}
There is a private release $T_f$ that attains equality of these laws. Therefore, for every convex function $\phi$ finite on a compact interval containing the range of $f$,
\begin{equation}\label{fdp:eq:scalar-opt}
 \sup_{U\perp Z}\E\phi\bigl(\E[f(Y)\mid U]\bigr)=\int_0^1\phi(m_f(t))\dd t.
\end{equation}
In particular,
\begin{equation}\label{fdp:eq:squared-opt}
 \sup_{U\perp Z} S_f(U)=\int_0^1m_f(t)^2\dd t.
\end{equation}
The same release attains all the convex objectives in \eqref{fdp:eq:scalar-opt} for this fixed $f$.
\end{theorem}
\begin{proof}
Write $X=\E[f(Y)\mid U]$. For an event $A\in\sigma(U)$, set $s=\Prb(A)$. Independence gives $\Prb(A\mid Z=z)=s$ almost surely. The upper-tail rearrangement bound, applied conditionally on $Z=z$, yields
\[
 \E[\ind_A f(Y)\mid Z=z]\le\int_{1-s}^1Q_z(t)\dd t.
\]
This bound follows by putting mass $s$ on the largest values of $f(Y)$, with randomization at a boundary atom. Consequently, for every real $c$,
\begin{align*}
 \E(X-c)_+
 &=\sup_{A\in\sigma(U)}\{\E[\ind_A f(Y)]-c\Prb(A)\}\\
 &\le\sup_{0\le s\le1}\left\{\int_{1-s}^1m_f(t)\dd t-cs\right\}
 =\int_0^1(m_f(t)-c)_+\dd t.
\end{align*}
The final identity uses monotonicity of $m_f$. The two random variables have the same mean, so the stop-loss characterization of convex order proves \eqref{fdp:eq:scalar-cx}.

To attain it, use the randomized conditional probability integral transform. With $F_z$ the conditional distribution function of $f(Y)$ and $V$ an independent uniform variable, set
\[
 T_f=F_Z(f(Y)-)+V\bigl(F_Z(f(Y))-F_Z(f(Y)-)\bigr).
\]
It is uniform conditionally on $Z$, hence private, and $f(Y)=Q_Z(T_f)$ almost surely. Since $Z\perp T_f$, conditioning on $T_f=t$ gives
\[
 \E[f(Y)\mid T_f=t]=\int Q_z(t)\rho(dz)=m_f(t).
\]
This proves attainment and all the stated formulas.
\end{proof}

The law in \eqref{fdp:eq:quantile-barycenter} is the one-dimensional quadratic Wasserstein barycenter of the conditional laws of $f(Y)$. Its importance here is the change of quantifiers: each observable has one release controlling every convex posterior-mean objective, whereas the existence of one release for all observables is governed by Theorem~\ref{fdp:thm:classification}.

\begin{corollary}[One rank release for all monotone tasks]\label{fdp:cor:monotone}
If $Y$ is real-valued, its randomized conditional rank
\[
 T=F_{Y\mid Z}(Y-)+V\bigl(F_{Y\mid Z}(Y)-F_{Y\mid Z}(Y-)\bigr)
\]
is one private release attaining \eqref{fdp:eq:scalar-opt} simultaneously for every bounded monotone observable $f$ and every convex posterior-mean objective in that formula.
\end{corollary}
\begin{proof}
For nondecreasing $f$, the function $f(Q_z(t))$ is a quantile of the conditional law of $f(Y)$. The same rank release therefore attains the proof of Theorem~\ref{fdp:thm:scalar} for all such $f$. For nonincreasing $f$, reverse the uniform parameter. The resulting posterior-mean distribution is unchanged by this reversal.
\end{proof}

This corollary explains where the obstruction to a greatest experiment appears. The rank release already solves every monotone scalar task. The full Blackwell order also compares nonmonotone decisions. The Gaussian example in Subsection~\ref{fdp:sec:gaussian} isolates the difference using $Y$ and $Y^2$.

\subsection{A coupling lemma that forces the classification}\label{fdp:sec:coupling}

The remaining question is whether the conditional-quantile optima for different observables can come from one release. Equality in each scalar problem constrains the same pairwise conditional coupling. For bounded real marginals, a coupling maximizes the product integral if and only if it is comonotone. One useful form is the two-copy criterion: if $(A,B)$ and $(A',B')$ are independent copies, then
\begin{equation}\label{fdp:eq:comonotone}
 (A-A')(B-B')\ge0\quad\text{almost surely}.
\end{equation}
The usual quantile coupling has this property. Conversely, a positive mass of crossing pairs permits a mass exchange that strictly increases the product integral. Choosing two separated rational rectangles makes this exchange quantitative and proves the criterion even with atoms.

\begin{lemma}[Simultaneous scalar transport]\label{fdp:lem:pair}
Let $\pi$ couple two laws $\mu,\nu$ on a standard Borel space $S$. Suppose that, for every $f\in\mathcal F$, $\pi$ maximizes
\[
 \int f(x)f(y)\pi(dx,dy)
\]
among couplings of $\mu$ and $\nu$. Then $\pi$ is supported on one of the following:
\begin{enumerate}[label=(\roman*),nosep]
\item the diagonal $\{(x,x):x\in S\}$;
\item one row $\{a\}\times S$;
\item one column $S\times\{b\}$;
\item $\{(a,a),(a,b),(b,b)\}$ for two distinct points $a,b$, with the orientation possibly reversed.
\end{enumerate}
In particular, either $\mu=\nu$, one marginal is a point mass, or both marginals are supported on one two-point set. Conversely, each of those marginal conditions admits a coupling that is optimal for every bounded scalar observable.
\end{lemma}
\begin{proof}
Apply \eqref{fdp:eq:comonotone} to each $f\in\mathcal F$. Countability gives one full-measure set on which all the inequalities hold. Because the algebra separates points, it can assign any prescribed values in $\{-1,0,1\}$ to any finite collection of distinct points. Consequently, for two independent pairs $(x,y),(x',y')$ drawn from $\pi$, the inequalities force
\begin{equation}\label{fdp:eq:compatibility}
 x=x'\quad\text{or}\quad y=y'\quad\text{or}\quad(x=y\text{ and }x'=y').
\end{equation}
Indeed, if neither coordinate agrees and both pairs are diagonal, every product is a square. In every other configuration, a ternary valuation makes the two differences have opposite signs. This can be checked by assigning opposite extreme values to the two ordered differences; a shared intermediate point receives zero.

If $\pi$ is diagonal, the conclusion follows. Otherwise Fubini's theorem allows an off-diagonal pair $(a,b)$, $a\ne b$, such that \eqref{fdp:eq:compatibility} holds against $\pi$-almost every other pair. Since $(a,b)$ is not diagonal, this gives
\[
 \pi\bigl((\{a\}\times S)\cup(S\times\{b\})\bigr)=1.
\]
If all mass is on one row or one column, we are done. Otherwise both parts outside their intersection have positive mass. Compare a pair $(a,y)$ with $y\ne b$ to a pair $(x,b)$ with $x\ne a$. Neither coordinate agrees. The last alternative in \eqref{fdp:eq:compatibility} therefore requires $y=a$ and $x=b$ almost surely. The support is the stated triangular set.

For the converse, use the diagonal coupling when $\mu=\nu$, the unique row or column coupling when a marginal is deterministic, and the monotone binary coupling when both marginals are supported on two points. Every scalar valuation preserves the two-copy condition.
\end{proof}

\subsection{Necessity on standard Borel spaces}\label{fdp:sec:necessity}

To apply the pairwise lemma, all conditional useful laws must be realized over one private parameter. The following refinement preserves the original source and release while making the useful variable a function of that parameter and the private data. Privacy then writes the loss of each scalar score as an integral of nonnegative pairwise transport gaps. A common optimizer forces those gaps to vanish almost everywhere.

\begin{lemma}[Private refinement]\label{fdp:lem:purify}
Every private full-data release $U$ has a private refinement $V$ such that $U$ is a function of $V$ and
\[
 Y=F(Z,V)
\]
for a measurable $F$.
\end{lemma}
\begin{proof}
Embed the useful standard Borel space into a Borel subset of $[0,1]$. Apply the randomized conditional probability integral transform to the law of $Y$ given $(Z,U)$. This constructs a uniform variable $T$ independent of $(Z,U)$, on an extension preserving the original $(Z,U,Y)$ law, such that $Y=F(Z,U,T)$. Set $V=(U,T)$. Since $U\perp Z$ and $T\perp(Z,U)$, the refinement remains private.
\end{proof}

Suppose $U$ maximizes all the scores in Theorem~\ref{fdp:thm:quadratic}. Its refinement $V$ does too, because conditional Jensen gives $S_f(V)\ge S_f(U)$ and the latter is already maximal. Write $Y=F(Z,V)$ and $\lambda=\law(V)$. For almost every $z$, the map $v\mapsto F(z,v)$ pushes $\lambda$ to $\mu_z$.

For a fixed $f\in\mathcal F$, define
\[
 H_f(z,z')=\int_0^1Q_z(t)Q_{z'}(t)\dd t
       -\int f(F(z,v))f(F(z',v))\lambda(dv).
\]
The rearrangement inequality gives $H_f\ge0$. Privacy and Fubini's theorem give
\begin{align*}
 S_f(V)
 &=\int\left(\int f(F(z,v))\rho(dz)\right)^2\lambda(dv),\\
 \int_0^1m_f(t)^2\dd t-S_f(V)
 &=\iint H_f(z,z')\rho(dz)\rho(dz').
\end{align*}
By scalar optimality, the left side is zero. Thus $H_f=0$ for almost every pair $(z,z')$. Intersecting over the countable family $\mathcal F$ shows that, for almost every pair, the coupling
\[
 \pi_{z,z'}=\law(F(z,V),F(z',V))
\]
is simultaneously optimal for every $f\in\mathcal F$. Lemma~\ref{fdp:lem:pair} applies.

It remains to convert the pairwise conclusion into a common conditional-law structure. Restrict to $N$. If the random conditional law $\mu_Z$ is almost surely constant there, condition (i) holds. Otherwise, choose a nondeterministic conditional law $\mu_*$ whose full-measure compatibility section contains another distinct nondeterministic law. Their union of supports has at most two points, so $\mu_*$ is supported on exactly two points $a,b$. Every other nondeterministic law in that section either equals $\mu_*$ or has union of support of size at most two with $\mu_*$. It must therefore be supported on the same $\{a,b\}$. This is condition (ii).

This proves necessity in Theorem~\ref{fdp:thm:classification}, and also proves that simultaneous quadratic optimization forces one of its two alternatives. Notice that no moments, densities, topological support assumptions, or finite alphabets were used.

\subsection{Explicit greatest releases and the converse}\label{fdp:sec:sufficiency}

Write
\[
 \ell=\rho(N),\qquad \alpha=\int_D\delta_{d(z)}\rho(dz).
\]
The measure $\alpha$ has mass $1-\ell$ and is fixed in every private posterior law: independence of the release from $Z$ prevents the deterministic part from changing its weight.

\subsubsection{A common nondeterministic conditional law}

Assume condition (i), with $\ell>0$. Take $V\sim\nu$ independent of $Z$ and set
\[
 Y=\begin{cases}V,&Z\in N,\\d(Z),&Z\in D.\end{cases}
\]
This has the required source law and gives
\begin{equation}\label{fdp:eq:common-posterior}
 \law(Y\mid V=v)=\alpha+\ell\delta_v.
\end{equation}

For any private release $U$, its posterior has the form
\[
 \law(Y\mid U=u)=\alpha+\ell\eta_u,
 \qquad\int\eta_u\law(U)(du)=\nu.
\]
Define a joint law of $(V,U)$ by $\law(U)(du)\eta_u(dv)$ and disintegrate it as $\nu(dv)K(v,du)$. Garbling the release in \eqref{fdp:eq:common-posterior} through $K$ produces exactly the required posterior of $Y$ given $U$. Thus $V$ dominates every private release.

\subsubsection{A common two-point support}

Assume condition (ii). Write $p(z)=\mu_z(\{b\})$ on $N$. Take $V\sim\mathrm{Unif}(0,1)$ independent of $Z$ and set
\[
 Y=\begin{cases}
 b,&Z\in N,\ V\le p(Z),\\
 a,&Z\in N,\ V>p(Z),\\
 d(Z),&Z\in D.
 \end{cases}
\]
The posterior is
\begin{equation}\label{fdp:eq:binary-posterior}
 \law(Y\mid V=t)=\alpha+(\ell-x(t))\delta_a+x(t)\delta_b,
 \qquad x(t)=\int_N\ind_{\{t\le p(z)\}}\rho(dz).
\end{equation}
The function $x$ is nonincreasing.

Every private release $U$ has a posterior
\[
 \law(Y\mid U=u)=\alpha+(\ell-v(u))\delta_a+v(u)\delta_b,
\]
where $v(u)=\int_N r_z(u)\rho(dz)$, $0\le r_z\le1$, and
$\int r_z(u)\law(U)(du)=p(z)$. For an event $A\in\sigma(U)$ of probability $s$,
\[
 \E[\ind_A r_z(U)]\le\min\{p(z),s\}.
\]
Consequently,
\begin{align*}
 \E(v(U)-c)_+
 &\le\sup_{0\le s\le1}\left\{\int_N\min(p(z),s)\rho(dz)-cs\right\}\\
 &=\sup_{0\le s\le1}\left\{\int_0^sx(t)\dd t-cs\right\}
 =\int_0^1(x(t)-c)_+\dd t.
\end{align*}
The means also agree. Hence $v(U)\preceq_{\rm cx}x(V)$.

Strassen's martingale coupling theorem~\cite{Strassen} gives a coupling of these two scalar laws with
\[
 \E[x(V)\mid v(U)]=v(U).
\]
Attach the original conditional law of $U$ given $v(U)$. The resulting coupling satisfies $\E[x(V)\mid U]=v(U)$. Disintegrating it produces a channel from $x(V)$, hence from $V$, to $U$. The affine posterior formula \eqref{fdp:eq:binary-posterior} then verifies \eqref{fdp:eq:garbling}. Thus $V$ is greatest.

\subsubsection{Every simultaneous quadratic optimizer is greatest}

We have proved existence of a greatest release whenever a simultaneous quadratic optimizer exists. To complete Theorem~\ref{fdp:thm:quadratic}, let $V$ be a greatest release and let $U$ maximize every score in $\mathcal F$. Realize $U$ as a garbling of $V$. For every $f\in\mathcal F$, conditional Jensen gives
\[
 S_f(V)-S_f(U)=\E\Var\bigl(\E[f(Y)\mid V]\mid U\bigr).
\]
Both scores are maximal, so this difference is zero. Countability implies that all the posterior integrals against $\mathcal F$ are determined by $U$. Since $\mathcal F$ contains a countable determining class of indicators, the whole posterior law of $Y$ given $V$ is determined by $U$. Disintegrating $V$ conditionally on that posterior gives a reverse garbling from $U$ to $V$. Hence $U$ is Blackwell-equivalent to the greatest release. The converse follows immediately from data processing.

\begin{corollary}[Posterior covariance matrices]\label{fdp:cor:covariance}
Suppose $Y$ takes $n$ values, with prior vector $p$, and write $R_U$ for the posterior probability vector given a private release $U$. There is an attained greatest matrix in the Loewner order among
\[
 \{\E[(R_U-p)(R_U-p)^{\mathsf T}]: U\perp Z\}
\]
if and only if the conditional laws have one of the two structures in Theorem~\ref{fdp:thm:classification}. Every release attaining that greatest matrix is a greatest private experiment.
\end{corollary}
\begin{proof}
For a valuation vector $t\in\R^n$, the quadratic form of this matrix is $\Var(\E[t_Y\mid U])$. An attained greatest matrix therefore gives a common optimizer for all ternary valuations. Apply Theorem~\ref{fdp:thm:quadratic}. The reverse implication follows from conditional Jensen under garbling.
\end{proof}

\subsection{Vector mean-square utility and Wasserstein barycenters}\label{fdp:sec:barycenter}

The barycenter connection extends to vector prediction. Wasserstein barycenters in quadratic transport~\cite{R27-AguehCarlier} also occur in fair regression under demographic parity; the scalar fair-regression connection is explicit in Chzhen, Denis, Hebiri, Oneto and Pontil~\cite{R27-Chzhen}. Theorem~\ref{fdp:thm:classification} determines when those task-specific optima have one common greatest experiment.

Assume now that $Y\in\R^d$ and $\E|Y|^2<\infty$. Let $\mathcal P_2(\R^d)$ denote the laws with finite second moment and define
\begin{equation}\label{fdp:eq:barycenter-functional}
 B=\inf_{\nu\in\mathcal P_2(\R^d)}\int W_2^2(\nu,\mu_z)\rho(dz).
\end{equation}

\begin{theorem}[All mean-square private optimizers]\label{fdp:thm:vector}
The infimum in \eqref{fdp:eq:barycenter-functional} is attained, and
\begin{equation}\label{fdp:eq:vector-score}
 \sup_{U\perp Z}\E\bigl|\E[Y\mid U]\bigr|^2=\E|Y|^2-B.
\end{equation}
For an optimal release $U$, put $X=\E[Y\mid U]$. Then $\law(X)$ is a minimizer of \eqref{fdp:eq:barycenter-functional}, and $\law(X,Y\mid Z=z)$ is an optimal quadratic transport coupling for almost every $z$.

Conversely, every minimizer $\nu$ and a measurable family of optimal couplings between $\nu$ and $\mu_z$ produce a private release $X$ satisfying
\[
 \law(X)=\nu,\qquad \E[Y\mid X]=X,
\]
and attaining \eqref{fdp:eq:vector-score}.
\end{theorem}
\begin{proof}
For a private $U$, its posterior mean $X$ is also private. The conditional law of $(X,Y)$ given $Z=z$ couples the same law $\nu=\law(X)$ to $\mu_z$. Since $\E[Y\mid X]=X$,
\begin{equation}\label{fdp:eq:orthogonality}
 \int W_2^2(\nu,\mu_z)\rho(dz)
 \le\E|Y-X|^2=\E|Y|^2-\E|X|^2.
\end{equation}
This proves the upper bound in \eqref{fdp:eq:vector-score}.

Existence in \eqref{fdp:eq:barycenter-functional} follows from the direct method. For every $\nu$,
\[
 \int|x|^2\nu(dx)\le2\int W_2^2(\nu,\mu_z)\rho(dz)+2\E|Y|^2.
\]
A minimizing sequence therefore has uniformly bounded second moments and is tight in $\R^d$. A weak limit has finite second moment. Lower semicontinuity of quadratic transport cost for each $z$, followed by Fatou's lemma, proves that this limit minimizes the functional.

Fix a minimizer $\nu$. Optimal coupling sets for fixed marginals are nonempty and compact in the weak topology. Their graph is Borel when the marginals range over $\mathcal P_2$, and standard measurable selection gives a family $\pi_z$ of optimal couplings, up to completion of $\rho$. Assemble the joint law
\[
 \rho(dz)\pi_z(dx,dy).
\]
Its $(Z,Y)$ marginal is the prescribed source and $X\perp Z$, because the first marginal of every $\pi_z$ is $\nu$. Let $b(X)=\E[Y\mid X]$. The orthogonal decomposition gives
\[
 B=\E|Y-X|^2=\E|Y-b(X)|^2+\E|b(X)-X|^2.
\]
The pair $(b(X),Y)$ conditionally on $Z=z$ is a feasible transport coupling from $\law(b(X))$ to $\mu_z$. Therefore
\[
 B\le\int W_2^2(\law(b(X)),\mu_z)\rho(dz)
 \le B-\E|b(X)-X|^2.
\]
It follows that $b(X)=X$ almost surely. This proves attainment of \eqref{fdp:eq:vector-score}.

For an arbitrary optimal release, both inequalities in \eqref{fdp:eq:orthogonality} must be equalities. Its posterior mean law is therefore a barycenter, and the nonnegative conditional transport-cost gaps integrate to zero. The conditional couplings are optimal almost everywhere.
\end{proof}

\begin{remark}
For scalar laws, the barycenter is unique and has quantile $m_f$. In higher dimensions, uniqueness requires additional assumptions. The theorem represents all posterior-mean laws without imposing such assumptions. It does not claim that all mean-square optimal releases are Blackwell-equivalent.
\end{remark}

\subsection{A two-by-three prediction gap}\label{fdp:sec:certificate}

Let $Z$ be equiprobable on two states and let $Y$ take three values. Its conditional probability vectors are
\[
 p=(1/2,1/3,1/6),\qquad q=(1/6,1/3,1/2).
\]
The unconditional useful law is uniform. Take the target valuations
\[
 f=(-1,0,1),\qquad g=(0,1,0).
\]

\begin{proposition}[A two-by-three simultaneous utility gap]\label{fdp:prop:finite-gap}
For this source,
\[
 \sup_{U\perp Z}S_f(U)=\frac12,\qquad
 \sup_{U\perp Z}S_g(U)=\frac13,
\]
whereas every private full-data release satisfies
\begin{equation}\label{fdp:eq:finite-gap}
 S_f(U)+S_g(U)\le\frac23.
\end{equation}
The sum of the deficits from the two individual optima is therefore at least $1/6$. The bound is attained.
\end{proposition}
\begin{proof}
Private refinement reduces the upper bound to releases specifying a pair $(i,j)$ whose distribution couples $p$ and $q$: when the private state is the first or second one, the useful value is respectively $i$ or $j$. The posterior is $(\delta_i+\delta_j)/2$. Define
\[
 c_{ij}=\left(\frac{f_i+f_j}{2}\right)^2+
         \left(\frac{g_i+g_j}{2}\right)^2.
\]
The rational row and column potentials
\[
 a=(0,1/2,1),\qquad b=(1,1/2,0)
\]
satisfy $c_{ij}\le a_i+b_j$ for all nine pairs. Their dual value is
\[
 \sum_i p_i a_i+\sum_jq_jb_j=\frac23,
\]
which proves \eqref{fdp:eq:finite-gap} for every refined release, and hence for every release by conditional Jensen.

The monotone coupling for $f$ has masses $1/6,1/3,1/3,1/6$ on
\[
 (1,1),(1,2),(2,3),(3,3).
\]
It gives $S_f=1/2$ and $S_g=1/6$. The coupling with masses $1/6,1/3,1/6,1/3$ on
\[
 (1,1),(2,2),(3,3),(1,3)
\]
gives $S_f=1/3$ and $S_g=1/3$. Scalar quantile optimality proves the first individual bound. The second is bounded by $\E g(Y)^2=1/3$ and is attained by the second coupling. Both couplings attain the joint upper bound.
\end{proof}

\subsection{Observation models and universal fair representations}\label{fdp:sec:consequences}

\subsubsection{Access to the private data can destroy existence of a greatest release}

The smooth additive-noise examples of Theorem~\ref{R22-thm:privacy-instability} have positive conditional densities on the whole real line and retain nonzero dependence between $Z$ and $Y$. Their useful-data-only private channels factor through a single modulo release. The classification gives a different answer as soon as the encoder may use the full pair.

\begin{corollary}[The two observation models have different greatest elements]\label{fdp:cor:two-models}
For every dependent source in the smooth full-support additive-noise family of Theorem~\ref{R23-thm:private-greatest}, a greatest useful-data-only perfectly private release exists, while no greatest full-data perfectly private release exists.
\end{corollary}
\begin{proof}
The useful-data-only assertion is Theorem~\ref{R23-thm:private-greatest}. Every conditional useful law has a positive density and is therefore atomless. The source is dependent, so Corollary~\ref{fdp:cor:atomless} excludes a greatest full-data release.
\end{proof}

This conclusion includes the covariance-matched analytic uniformly log-concave approximations of Theorem~\ref{R22-thm:privacy-instability}, for which the greatest useful-data-only release already has infinite mutual information. Infinite utility for that one numerical objective therefore coexists with the absence of a greatest experiment in the enlarged full-data class.

\subsubsection{Universal representations satisfying demographic parity}

In statistical learning, the condition $U\perp Z$ is demographic parity for a representation $U$ and a protected variable $Z$. Taking $Y$ to be the useful feature vector, a downstream task is a decision problem about $Y$. Theorem~\ref{fdp:thm:classification} classifies when one representation that may use the full data can be optimal for every such task. Conditional atomlessness and dependence exclude this possibility. Each specified regression task still has its own optimal quantile or Wasserstein barycenter construction.

The task-specific fair-regression construction is credited to~\cite{R27-Chzhen}. The universal classification follows from simultaneous optimality of the ternary tests, and applies to the full standard Borel source.


\section{Calibration tests and Blackwell comparison}
\label{R11-sec:Blackwell}
Blackwell compares experiments by all decision problems, equivalently by an observation kernel that garbles one experiment into the other~\cite{Blackwell}. For a fixed prior, the posterior law carries this comparison. The preceding section compares private releases by the garbling identity~\eqref{fdp:eq:garbling}. The tests here ask whether an independent source can be added with specified conditional means while the old observation--target law remains fixed. Scalar binary costs identify a symmetric source; comparison of experiments requires the full family of vector-source tests and independent target refinements.

One unobserved fair bit permits events of posterior probability $k_2+\tfrac12k_1$. It changes the available calibration constraints without changing the observed information. The three-state example displays the resulting failure of unrefined comparison. Cyclic transport encodes arbitrary finite decision functions, while transport costs have only finitely many walls at every bounded refinement level; hiding a discrepancy between those walls proves strictness. Quantitative completion then invokes barycentric weak-transport duality~\cite{GRST} to produce one simulation kernel for every later bounded decision problem. The kernel also preserves the specified predictive observation marginal exactly.

\subsection{A full binary calibration cost determines the law}\label{TOMO-sec:binary}
The preceding information potential records more than its local Hessian.
Even with one fair target bit, the complete cost as the requested mean
varies determines a symmetric reference. The proof recovers every
one-dimensional projection from a single integral transform.

\subsubsection{The potential under a finite first moment}
Let $\mu$ be a Borel probability measure symmetric under $x\mapsto-x$ on $\R^d$, with $\int\|x\|\mu(\dd x)<\infty$. Given a fair sign $\sigma\in\{-1,1\}$, define
\[
 \mathcal I_\mu(b)=\inf I(X;\sigma),
 \qquad X\sim\mu,\quad
 \E[X\mid\sigma]=\sigma b,
\]
and set the value to $+\infty$ if the constraint is infeasible. All information is measured in nats. The dual potential is
\begin{equation}\label{TOMO-eq:calpotential}
 F_\mu(\theta)=\E\log\cosh\ip{\theta}{X},\qquad
 \mathcal I_\mu=F_\mu^*.
\end{equation}
For completeness, write $m(x)=\E[\sigma\mid X=x]$. The information integrand is
\[
 h(m)=\tfrac12[(1+m)\log(1+m)+(1-m)\log(1-m)].
\]
Antisymmetrizing $m$ preserves $\E[Xm(X)]$ and reduces its convex integral, while ensuring a fair sign. Pointwise duality gives the optimizer $m(x)=\tanh\ip{\theta}{x}$ and the potential in \eqref{TOMO-eq:calpotential}. The feasible channels form a weak-star compact subset of $L^\infty(\mu)$; the moment constraints are continuous because $X\in L^1$. Lower semicontinuity of the convex entropy integral proves that the extended cost is closed. Thus the conjugacy also controls boundary values.

\begin{theorem}[Identification from fair-binary calibration]\label{TOMO-thm:binary}
If $\mu$ and $\nu$ are centrally symmetric probability measures on $\R^d$ with finite first moments and
\[
 \mathcal I_\mu(b)=\mathcal I_\nu(b)\quad\text{for every }b\in\R^d,
\]
then $\mu=\nu$.
\end{theorem}
\begin{proof}
Biconjugacy gives $F_\mu=F_\nu$. On a direction $v$, this determines
$f(t)=\E\log\cosh(tR)$ for $R=|\ip{v}{X}|$, and
$m=\lim_{t\to\infty}f'(t)=\E R$. If $m=0$ the projection vanishes.
Otherwise the Mellin inversion in the proof of
Theorem~\ref{SP-thm:entropy-inverse} applies after normalization. In its
unnormalized form it reads
\begin{equation}\label{TOMO-eq:mellin}
 \int_0^\infty t^{s-1}(m-f'(t))\,dt
 =2^{1-s}\Gamma(s)\eta(s)\E R^{1-s},\qquad0<\Re s<1.
\end{equation}
The first moment ensures absolute integrability; the isolated zeros of the
analytic multiplier are filled by continuity. Fourier uniqueness on the
logarithmic coordinate recovers $R>0$, and normalization recovers the zero
atom. Symmetry gives each signed projection, so Cram\'er--Wold proves the claim.
\end{proof}

\subsubsection{Consequences for finite representations and for processes}
\begin{corollary}\label{TOMO-cor:binarycompression}
A finite atomic centrally symmetric reference cannot have the same complete fair-binary calibration cost as a nonatomic centrally symmetric integrable reference.
\end{corollary}
The assertion follows from Theorem~\ref{TOMO-thm:binary}.

For a process whose finite-dimensional distributions are invariant under simultaneous sign reversal, the complete binary costs of every finite observation vector determine all its finite-dimensional distributions. Thus this nonlinear information observable can determine the process law. The observation consists of the full nonlinear cost for every finite vector.


\subsection{Deterministic garblings with an exact output law}
\label{SP-sec:experiments}

Let $\Theta\in\{1,\ldots,q\}$ have prior $p$, and let $Y$ be an observation with atomless marginal $\nu$. Its posterior vector is
\[
 F(Y)=\bigl(\Pp(\Theta=i\mid Y)\bigr)_{i=1}^q.
\]
A randomized garbling of $Y$ gives an output $X$ with law $\mu$ and posterior $g(x)=\law(\Theta\mid X=x)$. The relation between these posterior vectors is $\E[F(Y)\mid X]=g(X)$, which is the usual martingale description of Blackwell comparison \cite{Blackwell}.

\begin{theorem}[Deterministic garbling after strict attenuation]
\label{SP-thm:blackwell}
For every randomized garbling with output law $\mu$ on a Polish space and every $0\le\lambda<1$, there is a Borel statistic $T_\lambda(Y)$ with law $\mu$ such that
\[
 \Pp(\Theta=i\mid T_\lambda(Y)=x)
   =p_i+\lambda(g_i(x)-p_i),\qquad i=1,\ldots,q.
\]
Equivalently, the entire state--output law is
\begin{equation}\label{SP-eq:attenuated-experiment}
 \law(\Theta,T_\lambda(Y))
 =\lambda\law(\Theta,X)+(1-\lambda)p\otimes\mu.
\end{equation}
For $0<\lambda<1$, $\law(T_\lambda(Y),Y)$ can be chosen arbitrarily close in $W_\infty$ to
\[
 \lambda\law(X,Y)+(1-\lambda)\mu\otimes\nu.
\]
\end{theorem}
\begin{proof}
Apply Theorem~\ref{WC-thm:attenuate} to the posterior vector $F$ and the garbling coupling. Its profile is $p+\lambda(g-p)$ and its output marginal is $\mu$. Conditional expectation from the state through the observation gives the displayed posterior identity; integration gives the full state--output law. The $W_\infty$ comparison follows from the same theorem, including a constant posterior in relative dimension zero.
\end{proof}
The source experiment and the full output marginal are exact. The state--output law differs from the requested one by at most $1-\lambda$ in total variation. For any payoff bounded in absolute value by one, the corresponding expected payoff changes by at most $2(1-\lambda)$. The construction gives a measurable statistic; computing it has no asserted running-time bound.

\begin{example}[Strict attenuation can be necessary for an exact output law]
\label{SP-ex:blackwell-endpoint}
Let $Y\sim\Unif[0,1]$ and $\Theta\mid Y\sim\operatorname{Bernoulli}(Y)$. A randomized garbling appends an independent $V\sim\Unif[0,1]$, giving target output $(Y,V)$ with the uniform square law and posterior equal to its first coordinate. A deterministic statistic with this exact output law and posterior would attain the zero cost in Example~\ref{R10-ex:weak-nonattainment}. Its second-moment argument forces the first output coordinate to equal $Y$, so the second cannot have the required independent uniform law. Thus the endpoint $\lambda=1$ fails. 
\end{example}
The randomized target experiment in this example is Blackwell-equivalent to the original observation; the obstruction concerns retaining its particular prescribed output law as a deterministic statistic.


\subsection{Reference refinement and Blackwell comparison}
\label{BW-sec}

Section~\ref{NLI-sec:global} determines the feasible extension set for
one fixed experiment, and Section~\ref{TOMO-sec:binary} identifies a
source from the full information cost of binary calibration. Here we
compare two retained experiments using feasibility of added independent
sources. The vector-source support function is an averaged conditional
transport value. Target refinement makes these tests complete, even
though it changes none of the observed information.

Let $(Z,Y)$ be an experiment with $Y\in[q]$, and let
$K_i(Z)=\Prb(Y=i\mid Z)$ be its posterior vector. Write $\nu$ for the
law of $K$ on $\Delta_q$ and $p=\E K$ for its prior. Blackwell's order
on experiments with this prior is the convex order of their posterior
laws \cite{Blackwell}. Every statement below concerns finite label sets;
the variable $Z$ and the independent source may be arbitrary standard
Borel random variables with the stated integrability.

The private experiments in Theorems~\ref{R23-thm:private-greatest}--\ref{R23-thm:private-maximal} give explicit endpoints for this order: one class has a greatest release, while another has incomparable maximal releases at the same information capacity. Proposition~\ref{R23-prop:privacy-simplex} identifies the finite-state geometric obstruction. Here the question concerns arbitrary pairs of experiments: how much independent-source calibration data suffice to determine their order? The refinement below increases the range of tests while preserving the observed experiment exactly.

\subsubsection{Vector-source extension sets}\label{BW-sec:calibration}
\begin{definition}
For a centered integrable law $\rho$ on $\R^d$, define
\[
 \mathcal C_\nu(\rho)=\left\{(\E[X\ind_{Y=i}])_{i=1}^q:
 \law(X)=\rho,\ X\perp Z,\ \law(Z,Y)\text{ is prescribed}\right\}.
\]
The elements are vectors $m=(m_1,\ldots,m_q)$ with $m_i\in\R^d$ and $\sum_i m_i=0$.
\end{definition}
Every admissible coupling is represented by kernels $h_i(x,z)\ge0$ with
\begin{equation}\label{BW-eq:kernel}
 \sum_i h_i(x,z)=1,\qquad \int h_i(x,z)\rho(dx)=K_i(z)
 \quad\text{for almost every }z.
\end{equation}
Its labelwise moments and conditional information are
\begin{align}
 m_i&=\int xh_i(x,z)\rho(dx)\Pp_Z(dz),\label{BW-eq:moments}\\
 \mathcal I(h)&=\int\sum_i h_i(x,z)\log\frac{h_i(x,z)}{K_i(z)}
                 \,\rho(dx)\Pp_Z(dz).\label{BW-eq:information}
\end{align}
The usual conventions apply at zero probabilities. Feasibility forces $h_i=0$ when $K_i=0$.

\begin{lemma}\label{BW-lem:posterior}
The calibration set depends on $(Z,Y)$ only through $\nu$. It is compact and convex. For $u=(u_1,\ldots,u_q)\in(\R^d)^q$, its support function is
\begin{equation}\label{BW-eq:support}
 h_{\mathcal C_\nu(\rho)}(u)=\int V_{\rho,u}(k)\nu(dk),
\end{equation}
where
\[
 V_{\rho,u}(k)=\sup\left\{\sum_i\int \langle u_i,x\rangle\gamma_i(dx):
 \sum_i\gamma_i=\rho,\ \gamma_i(\R^d)=k_i\right\}.
\]
The function $V_{\rho,u}$ is concave and continuous on $\Delta_q$.
\end{lemma}
\begin{proof}
Conditioning a kernel in \eqref{BW-eq:kernel} on $(X,K(Z))$ gives a kernel depending on $z$ only through $K(z)$, with the same moments. Conversely, any such kernel may be applied to the original $Z$, so that its entire pair law with $Y$ is restored. This proves the first assertion.

The feasible kernels form a weak-star compact convex subset of the unit ball of $L^\infty(\rho\otimes\Pp_Z)^q$. The conditional constraints are closed because they may be tested against every bounded function of $z$. The moment maps are continuous in that topology because $x\in L^1(\rho)$. Their image is therefore compact and convex.

For finite $\rho$, conditional optimization is a finite transportation problem. Its value is concave in the prescribed target masses, and measurable optimal choices give \eqref{BW-eq:support}. Approximate a general $X$ in $L^1$ by finite-valued $X_n$. Replacing $X$ by $X_n$ changes every coupling objective by at most $\max_i\|u_i\|\E\|X-X_n\|$. This gives the same formula by uniform convergence of the finite-source value functions. Their limit is concave and continuous on the compact simplex.
\end{proof}

Write $\mu\cx\nu$ when $\int f\,d\mu\le\int f\,d\nu$ for every continuous convex $f$ on the simplex. Since the simplex is compact, no tail convention is needed. The common prior is the equality of barycenters.

\begin{corollary}\label{BW-cor:blackwell-forward}
If $\mu\cx\nu$, then $\mathcal C_\nu(\rho)\subseteq\mathcal C_\mu(\rho)$ for every centered integrable finite-dimensional source.
\end{corollary}
\begin{proof}
Apply the convex-order inequality to the negative of the concave function in \eqref{BW-eq:support}, and use separation of compact convex sets.
\end{proof}
The direction is worth noting: a more informative old experiment imposes more restrictions on what an independent new source can do.

\subsubsection{An unobserved fair bit changes the comparison}
\label{R8-sec:first-refinement}
Consider three target states and a fixed observation $Z$. Appending a
fair bit $B$, independent of $(Z,Y)$, replaces a posterior
$(k_1,k_2,k_3)$ by
\[
 (k_1/2,k_1/2,k_2/2,k_2/2,k_3/2,k_3/2).
\]
The observation is unchanged. A new source $X$, required to be
independent of $Z$, can nevertheless be correlated with the new target
bit. In particular, the event
\begin{equation}\label{R8-eq:refined-event}
 E=\{Y=2\}\cup\{Y=1,B=1\}
 \quad\hbox{has posterior probability }k_2+k_1/2.
\end{equation}
Its coefficients define a new direction in the posterior simplex,
which no unrefined subset of the three labels realizes.

Here is an explicit pair for which that new direction matters. On the
unit disk in $\R^2$, put
\begin{align*}
 \psi(x)&=(1-\|x\|_2^2)^4_+,\qquad
 Q=\frac1{18}\begin{pmatrix}71&-38\\-38&14\end{pmatrix},
 \qquad \eta=\frac9{9016},\\
 f_\mu(x)&=\frac3\pi(1-\|x\|_2^2)^2_+,\qquad
 f_\nu(x)=f_\mu(x)+\frac{3\eta}{\pi}Q:D^2\psi(x).
\end{align*}
Push the densities forward by
\[
 k=(1/3+x_1/12,\ 1/3+x_2/12,\ 1/3-(x_1+x_2)/12).
\]
Both posterior laws have the uniform prior and are supported in the
simplex interior. They satisfy
\begin{equation}\label{R8-eq:first-bit-comparison}
 \mathcal C_\nu(\rho)\subseteq\mathcal C_\mu(\rho)
 \quad\hbox{for every centred integrable finite-dimensional source }\rho.
\end{equation}
After the fair-bit refinement this inclusion fails for the scalar
source $X=\pm1/2$ with equal probabilities. The subset
\eqref{R8-eq:refined-event} has a strictly reversed support bound, with
gap
\begin{equation}\label{R8-eq:first-bit-gap}
 \frac{2\sqrt5}{71001\pi}>0.
\end{equation}

The proof in Section~\ref{BW-sec:hierarchy} explains all the constants.
For this first level, unrefined transportation walls have normal
directions $(1,0),(0,1),(1,1)$. Their quadratic forms under $Q$ are
$71/18$, $7/9$ and $1/2$, respectively. The new direction is $(1,2)$,
whose quadratic form is $-25/18$. Integration by parts transfers these
signs to every old transport value and to the new subset test. The
relative density perturbation has absolute value at most $1/2$, so
positivity is explicit. The general construction below produces
$C^\infty$ versions arbitrarily close to one another; the displayed
polynomial densities are $C^1$ across the disk boundary.

\begin{proposition}[The information retained by an independent refinement]
\label{R8-prop:information-invariance}
For any finite independent target refinement $B$ with positive masses
$\lambda_b$, the information density is unchanged:
\[
 \log\frac{\Prb(Y,B\mid Z)}{\Prb(Y,B)}
 =\log\frac{K_Y(Z)}{p_Y}\quad\hbox{almost surely}.
\]
In particular $I((Y,B);Z)=I(Y;Z)$. More generally, every
$f$-divergence between the joint law and the product of its marginals
is unchanged, as is every R\'enyi divergence between those two laws
whenever defined. The prior-dependent decision problem for any loss
which depends only on $Y$ is also unchanged.
\end{proposition}
\begin{proof}
The joint posterior and prior coordinates are $K_i\lambda_b$ and
$p_i\lambda_b$, so their ratio is $K_i/p_i$. This identity holds under
both the joint law and the product law. Summing out $b$ therefore
preserves the distribution of the likelihood ratio under either
measure, which proves the divergence assertions. The observation and
the joint law of $(Z,Y)$ have been retained, proving the decision-risk
assertion. Zero-prior labels can simply be removed.
\end{proof}
Thus the strict increase in calibration testing power occurs under an
operation preserving the full information-density law. Its effect is
to enlarge the available target events for the independent source.
The completion theorem asks how much refinement makes those events
sufficient for every decision problem.

\subsubsection{Independent refinements recover Blackwell order}\label{BW-sec:completion}
A uniform $M$-refinement replaces $Y$ by $(Y,B)$, where $B$ is independent of $(Z,Y)$ and uniform on $[M]$. The new posterior law is
\[
 \nu^{[M]}=(R_M)_\#\nu,\qquad (R_Mk)_{i,b}=k_i/M.
\]
When comparing calibration sets, the new physical source may be coupled to the refined label; the auxiliary variable is independent of the old experiment, not required to remain independent of the new source.

\begin{theorem}[Refinement completion]\label{BW-thm:completion}
For posterior laws $\mu,\nu$ with the same prior, the following are equivalent:
\begin{enumerate}
\item $\mu\cx\nu$;
\item $\mathcal C_{\nu^{[M]}}(\rho)\subseteq\mathcal C_{\mu^{[M]}}(\rho)$ for every $M\ge1$ and every centered integrable finite-dimensional $\rho$;
\item the same inclusions hold for $M=2^b$, $b=0,1,2,\ldots$.
\end{enumerate}
\end{theorem}
A finite decision problem has a Bayes value which is a maximum of
finitely many affine functions of the posterior. The proof encodes that
maximum as the largest cumulative imbalance in a directed transport
cycle. Refinement implements the required rational grouping of labels.
The affine part of the cycle cost has the same expectation under both
posterior laws, so only the decision value remains in the comparison.
We first establish the grouping and transport identities.

\begin{lemma}[Grouping refined labels]\label{BW-lem:grouping}
Let $T$ be a $q\times m$ stochastic matrix whose entries are integer multiples of $1/M$. Group the $M$ copies of each label $i$ into $m$ classes with respective counts $MT_{ij}$. The image of the refined calibration set under summing moments within these classes is exactly $\mathcal C_{T_\#\nu}(\rho)$, where $T_\#\nu$ denotes the law of the row vector $KT$.
\end{lemma}
\begin{proof}
Every refined coupling produces a grouped coupling with conditional masses $(KT)_j$. Conversely, split a grouped label $j$, conditionally on $Z$, among its constituent copies $(i,b)$ in proportions $(K_i/M)/(KT)_j$. Perform this split independently of $X$ given $(Z,j)$. The required refined masses and grouped moments are then exact. Classes of zero conditional mass can be defined arbitrarily.
\end{proof}

\begin{lemma}[The cycle transportation value]\label{BW-lem:cycle}
On $m$ cyclically ordered sites put $a_i=1/m$ and transport cost $c_{ij}=(j-i)\bmod m$. For target masses $b$, set
\[
 S_j=\sum_{i=0}^{j-1}(b_i-a_i),\quad 1\le j<m,\qquad S_0=0.
\]
Then the minimum transportation cost is
\begin{equation}\label{BW-eq:cycle}
 W(a,b)=m\max_{0\le j<m}S_j-\sum_{j=1}^{m-1}S_j.
\end{equation}
Moreover, for any two target distributions,
\begin{equation}\label{BW-eq:transport-lipschitz}
 |W(a,b)-W(a,b')|\le\frac{m-1}{2}\|b-b'\|_1.
\end{equation}
\end{lemma}
\begin{proof}
Route each transported unit clockwise along the cycle. If $F_i$ is the flow on arc $i\to i+1$, conservation gives $F_i=t-S_{i+1}$, with $S_m=0$. Its least possible nonnegative value of $t$ is $\max_jS_j$, and the resulting total arc cost is \eqref{BW-eq:cycle}. At least one arc then has zero flow. Cut there and decompose the remaining acyclic flow into paths from surplus to deficit sites, allowing mass to stay at its site. Each path has length less than $m$, so this decomposition realizes the original cyclic transport cost and proves equality.

To change the demand from $b$ to $b'$, reassign at most $\|b-b'\|_1/2$ units from surplus destinations to deficit destinations while keeping their source sites fixed. Each unit changes its cost by at most $m-1$. Apply this argument in both directions to prove \eqref{BW-eq:transport-lipschitz}.
\end{proof}

\subsubsection{A finite quantitative witness}
\begin{theorem}[Linear-size refinement witness]\label{BW-thm:quantitative}
Suppose the refined calibration inclusion holds at a fixed integer $M$
for every centered integrable finite-dimensional source. Let
\[
 f(k)=\max_{0\le j<m}a_j\cdot k,\qquad
 a_j\in[-1,1]^q,\quad m\ge2.
\]
Then
\begin{equation}\label{BW-eq:quantitative}
 \int f\,d\mu-\int f\,d\nu\le\frac{2m}{M}.
\end{equation}
In fact it suffices to assume the inclusion for the uniform centered
regular-simplex source with atoms $e_\ell-m^{-1}\mathbf1$, supported on
$\mathbf1^\perp\subset\R^m$.
If the displayed expectation difference is $\Delta>0$, every
$M>2m/\Delta$ admits a failing comparison, witnessed by this source with
$m$ atoms in dimension $m-1$. If $M\ge4m/\Delta$, its support-function
violation in the cyclic cost field is at least $\Delta/4$.
\end{theorem}
\begin{proof}
Put $\delta=1/(2m)$ and encode the affine pieces by cumulative rows:
\begin{equation}\label{BW-eq:encoded-T}
 C_{i0}=0,\qquad
 C_{ij}=\frac jm+\delta(a_j-a_0)_i\quad(1\le j<m),
 \qquad C_{im}=1.
\end{equation}
Every consecutive increment is nonnegative, since its deviation from
$1/m$ is at most $2\delta=1/m$, including the first and last increments.
Thus these are the cumulative rows of a stochastic matrix $T$.
Round each $C_{ij}$ to the nearest multiple of $1/M$, using a common
monotone rounding map. The resulting $C^M_{ij}$ stay nondecreasing,
retain endpoints zero and one, and differ from $C_{ij}$ by at most
$1/(2M)$. Their increments define a stochastic matrix $T_M$ with
entries that are integer multiples of $1/M$.

For the target masses $b=kT_M$, the cycle prefixes satisfy
\[
 S_j^M(k)=\delta(a_j-a_0)\cdot k+e_j(k),\qquad
 |e_j(k)|\le\frac1{2M},\qquad e_0=0.
\]
Every $e_j$ is linear in $k$. Therefore
\[
 \left|\max_j S_j^M(k)-\delta(f(k)-a_0\cdot k)\right|
 \le\frac1{2M}.
\]
The sum of the prefixes in Lemma~\ref{BW-lem:cycle} is affine in $k$;
its expectations cancel because the posterior laws have the same
barycenter. The same is true of $a_0\cdot k$. Consequently
\begin{equation}\label{BW-eq:witness-gap}
 \E_\mu W(a,KT_M)-\E_\nu W(a,KT_M)
 \ge\frac\Delta2-\frac mM.
\end{equation}
This uses the exact prefix formula rather than an entrywise transport
Lipschitz bound.

Let the source atoms be $x_\ell=e_\ell-a$ with uniform masses
$a_\ell=1/m$, and take $u_j=-c_{\cdot j}$, projected onto
$\mathbf1^\perp$ if desired. Every cost column has average $(m-1)/2$,
so the conditional calibration support function is $V(kT_M)=(m-1)/2-W(a,kT_M)$. Grouping the $M$ refined copies as in Lemma~\ref{BW-lem:grouping}
realizes this as a projection of the refined calibration set.
Inclusion forces the left side of \eqref{BW-eq:witness-gap} to be
nonpositive, proving all the assertions.
\end{proof}

\begin{proof}[Proof of Theorem~\ref{BW-thm:completion}]
Convex order is preserved by a linear refinement map, so Corollary~\ref{BW-cor:blackwell-forward} proves (i)$\Rightarrow$(ii)$\Rightarrow$(iii). Assume (iii). Every maximum of finitely many affine functions, after addition of an affine function and a scalar normalization if needed, is covered by Theorem~\ref{BW-thm:quantitative}. Let $M=2^b\to\infty$. The resulting inequality holds for every convex polyhedral function. Continuous convex functions on the compact simplex are uniformly approximable by such functions; equivalently, use supporting planes after a convex Lipschitz approximation from an interior homothety. Thus $\mu\cx\nu$.
\end{proof}

\subsubsection{Quantitative recovery from exact or approximate tests}
Define the convex Lipschitz deficit by
\begin{equation}\label{R8-eq:metric-definition}
 d_{\mathrm{cx}}(\mu,\nu)=\sup\left\{\int f\,d\mu-\int f\,d\nu:
 f\text{ convex and }1\text{-Lipschitz in }\ell_1\right\}.
\end{equation}
The exact comparison theorem has a stable form. Its useful error
parameter is the normalized failure of the cyclic tests themselves.

\begin{lemma}[Approximate cyclic comparisons]
\label{R7-cor:approximate-blackwell}
Let $\mu,\nu$ have the same prior on $\Delta_q$. For an integer $M\ge1$,
define
\[
 e_M=\sup_{m\ge2}\ \sup_T
 \frac{[\E_\mu W(a,KT)-\E_\nu W(a,KT)]_+}{m-1},
 \qquad a=m^{-1}\mathbf1,
\]
where $T$ ranges over stochastic $q\times m$ matrices with entries in
$M^{-1}\mathbb Z$, and $W$ uses the directed $m$-cycle cost. Every
$f(k)=\max_{0\le j<m}a_j\cdot k$, with $a_j\in[-1,1]^q$, satisfies
\begin{equation}\label{R8-eq:approximate-polyhedral}
 \E_\mu f-\E_\nu f\le2(m-1)e_M+2m/M.
\end{equation}
Exact refined inclusion makes $e_M=0$.
\end{lemma}
\begin{proof}
The encoded matrix in Theorem~\ref{BW-thm:quantitative} belongs to
the defining family for $e_M$. Substitute the upper bound
$(m-1)e_M$ into \eqref{BW-eq:witness-gap}. Exact inclusion makes every
cyclic deficit nonpositive by the same support-function argument.
\end{proof}

\paragraph*{Convex approximation with a controlled number of planes}
\label{R8-sec:convex-approximation}
The cycle witness charges for the number of affine pieces in a test.
Convexity permits a more economical approximation than a uniform grid
of evaluation points. The relevant parameter is a point of the graph
of the subdifferential, written in the coordinates $x+g$.
Nearby parameters give a quadratic error in the supporting planes.
This is the elementary monotonicity underlying the proximal
parametrization of a convex subdifferential \cite{Moreau}; it gives the
classical convex-approximation exponent \cite{Dudley} with the coefficient
control needed by our cycle construction.

\begin{lemma}[Polyhedral approximation from the subdifferential]
\label{R8-lem:polyhedral-approximation}
Let $q\ge2$ and let $f$ be convex and $1$-Lipschitz in $\ell_1$ on
$\Delta_q$, normalized by $\min f=0$. Put
\[
 B_q=2+2\sqrt q,\qquad
 N_q(\epsilon)=\left(1+\frac{1+\sqrt q}{\sqrt\epsilon}\right)^{q-1}.
\]
For $0<\epsilon\le1$ there are $2\le m\le N_q(\epsilon)$ affine
forms $a_j\cdot k$, with $\|a_j\|_\infty\le B_q$, whose maximum
$f_\epsilon$ satisfies
\[
 0\le f-f_\epsilon\le\epsilon\quad\hbox{on }\Delta_q.
\]
Repeated forms are allowed when one supporting plane suffices.
\end{lemma}
\begin{proof}
Let $p_0=q^{-1}\mathbf1$ and $H=\mathbf1^\perp$. The function
\[
 F(x)=\inf_{k\in\Delta_q}\{f(k)+\|p_0+x-k\|_1\},\qquad x\in H,
\]
is finite, convex and $\sqrt q$-Lipschitz in the Euclidean norm, and
$F(k-p_0)=f(k)$. Its subgradients therefore have Euclidean norm at most
$\sqrt q$. Consider the compact subset of $H$
\[
 \mathcal Z=\{k-p_0+g:k\in\Delta_q,\ g\in\partial F(k-p_0)\}.
\]
It is contained in the ball of radius $1+\sqrt q$. Choose a maximal
$\delta$-separated subset $z_j=k_j-p_0+g_j$ of $\mathcal Z$, where
$\delta=2\sqrt\epsilon$. Disjoint balls of radius $\delta/2$ give
$m\le(1+2(1+\sqrt q)/\delta)^{q-1}=N_q(\epsilon)$.
The selected points form a $\delta$-net.

For $k\in\Delta_q$, choose $g\in\partial F(k-p_0)$ and a selected
$z_j$ within $\delta$ of $z=k-p_0+g$. The supporting plane
$\ell_j(k)=f(k_j)+\langle g_j,k-k_j\rangle$ obeys
\begin{align*}
 0\le f(k)-\ell_j(k)
 &\le\langle g-g_j,k-k_j\rangle\\
 &=\langle z-z_j,k-k_j\rangle-\|k-k_j\|_2^2
 \le\tfrac14\|z-z_j\|_2^2\le\epsilon.
\end{align*}
Thus $f_\epsilon=\max_j\ell_j$ has the required error. Since $0\le
f\le2$, the identity $\sum_i k_i=1$ writes $\ell_j(k)=a_j\cdot k$
with
\[
 (a_j)_i=f(k_j)+(g_j)_i-\langle g_j,k_j\rangle,
 \qquad |(a_j)_i|\le2+2\sqrt q.
\]
The displayed upper bound on $N_q(\epsilon)$ exceeds two, so a single
plane may be repeated if necessary.
\end{proof}

\begin{theorem}[Quantitative Blackwell recovery]
\label{R8-thm:blackwell-rate}
For posterior laws with the same prior and normalized cyclic deficit
$e_M$ as in Lemma~\ref{R7-cor:approximate-blackwell},
\begin{equation}\label{R8-eq:blackwell-rate}
 d_{\mathrm{cx}}(\mu,\nu)\le
 \inf_{0<\epsilon\le1}
 \left\{\epsilon+2B_qN_q(\epsilon)(e_M+M^{-1})\right\}.
\end{equation}
Consequently
\[
 d_{\mathrm{cx}}(\mu,\nu)
 =O_q\bigl((e_M+M^{-1})^{2/(q+1)}\bigr).
\]
Under exact refined calibration inclusion, $e_M=0$, so the rate is
$O_q(M^{-2/(q+1)})$. With $b$ fair auxiliary bits it is
$O_q(2^{-2b/(q+1)})$. These are upper bounds; optimality of the
recovery exponent is not asserted.
\end{theorem}
\begin{proof}
Apply the cycle witness inequality, including its error term from
Lemma~\ref{R7-cor:approximate-blackwell}, to $f_\epsilon/B_q$:
\[
 \E_\mu f_\epsilon-\E_\nu f_\epsilon
 \le2B_q\{(m-1)e_M+m/M\}.
\]
The one-sided approximation in Lemma~\ref{R8-lem:polyhedral-approximation}
adds at most $\epsilon$ to this deficit. Take the supremum over $f$.
For $e_M+M^{-1}\le1$, choose $\epsilon$ of order
$(e_M+M^{-1})^{2/(q+1)}$ and use
$N_q(\epsilon)=O_q(\epsilon^{-(q-1)/2})$. The uniform bound
$d_{\mathrm{cx}}\le2$ covers the other case. Exact inclusion makes
every cyclic deficit nonpositive, so $e_M=0$.
\end{proof}

\begin{corollary}[Size of a decision witness]
\label{R8-cor:witness-size}
For fixed $q\ge3$, a convex Lipschitz deficit $\eta>0$ is exposed by
one uniform centred simplex source with
$O_q(\eta^{-(q-1)/2})$ atoms and a uniform refinement of size
$O_q(\eta^{-(q+1)/2})$. Equivalently,
\[
 b\le\frac{q+1}{2}\log_2(1/\eta)+O_q(1)
\]
fair target bits suffice. The statement concerns existence and size of
a witness, and makes no oracle-complexity assertion for finding it.
\end{corollary}
\begin{proof}
Normalized convex $1$-Lipschitz tests form a compact family, so a
maximizer for the deficit exists. Approximate it with error
$\epsilon=\min\{\eta/2,1\}$. Its polyhedral deficit is at least
$\eta/2$, and it uses $m=O_q(\eta^{-(q-1)/2})$ planes of bounded
coefficients. Theorem~\ref{BW-thm:quantitative}, after scaling by
$B_q$, produces a violation for every $M>4B_qm/\eta$. Choosing an
integer $M$ or the next power of two proves the claims.
\end{proof}

\begin{corollary}[Canonical source tests]\label{BW-cor:canonical-sources}
Let $\rho_m$ be uniform on the centered regular simplex
$\{e_\ell-m^{-1}\mathbf1:0\le\ell<m\}$. For posterior laws with the
same prior, $\mu\cx\nu$ is equivalent to
\[
 \mathcal C_{\nu^{[2^b]}}(\rho_m)
 \subseteq\mathcal C_{\mu^{[2^b]}}(\rho_m)
 \quad\text{for all }m\ge2\text{ and }b\ge0.
\]
\end{corollary}
\begin{proof}
The forward implication follows from
Corollary~\ref{BW-cor:blackwell-forward}. Conversely, the strengthened
source-specific assertion in Theorem~\ref{BW-thm:quantitative}, with
$b\to\infty$, gives every convex polyhedral test. Affine normalization
and uniform convex approximation then give $\mu\cx\nu$.
\end{proof}

\begin{proposition}[The binary exception]\label{BW-prop:binary-unrefined}
For $q=2$, unrefined calibration inclusion for every centered scalar
two-point source is equivalent to Blackwell order.
\end{proposition}
\begin{proof}
Let $X$ take $1-t$ with probability $t$ and $-t$ with probability
$1-t$, where $0<t<1$. Its upper integrated quantile is
\[
 Q_t(s)=\min\{(1-t)s,t(1-s)\}=(1-t)s-(s-t)_+.
\]
The largest first-label calibration moment is $\E Q_t(K_1)$.
The inclusion gives $\E_\nu Q_t(K_1)\le\E_\mu Q_t(K_1)$.
Equal means cancel the linear term, leaving every stop-loss inequality
$\E_\mu(K_1-t)_+\le\E_\nu(K_1-t)_+$. These characterize convex order
on $[0,1]$. Since $K_2=1-K_1$, this is full posterior convex order.
The converse follows from the general forward implication.
\end{proof}

\subsubsection{What scalar sources see}
For a centered real source define its upper integrated quantile
\[
 Q_\rho(s)=\sup\{\E[Xa(X)]:0\le a\le1,\ \E a=s\},\qquad 0\le s\le1.
\]
This is a concave function with $Q_\rho(0)=Q_\rho(1)=0$. Conditional rearrangement shows that the scalar calibration set is the base polytope with subset bounds
\begin{equation}\label{BW-eq:scalar-subsets}
 \sum_{i\in A}m_i\le\E Q_\rho(K(A)),\qquad A\subseteq[q],
 \qquad \sum_i m_i=0.
\end{equation}
For clarity, these bounds are sufficient as well as necessary. Given a sorted field $u$, successive upper level sets form a nested chain. A single quantile allocation realizes every bound along that chain. Its objective is the usual telescoping sum of the field increments times the right sides of \eqref{BW-eq:scalar-subsets}. These are exactly the support functions of the base polytope, proving equality.

\begin{proposition}\label{BW-prop:linear-order}
All-refinement calibration comparison restricted to scalar sources is equivalent to linear convex order of the posterior laws: every one-dimensional projection of $\mu$ is below the corresponding projection of $\nu$ in convex order. Fair scalar signs suffice for this equivalence.
\end{proposition}
\begin{proof}
For a refined subset, its posterior mass is $a\cdot K$ with rational $a_i\in[0,1]$. Linear convex order and concavity of $Q_\rho$ give every comparison in \eqref{BW-eq:scalar-subsets}. Conversely, for a fair sign $Q_\rho(s)=\min(s,1-s)$. Given a field $v$ and threshold $t$, choose a small positive $\epsilon$ so $a_i=1/2+\epsilon(v_i-t)$ lies in $[0,1]$, and approximate these coefficients rationally. The resulting subset comparison is exactly $\E_\nu|v\cdot K-t|\ge\E_\mu|v\cdot K-t|$. Together with equal means, these absolute-value tests characterize every scalar convex-order comparison.
\end{proof}
Thus the dimension of the physical source is part of the comparison problem. Increasing only the number of scalar refinements does not silently replace linear convex order by full Blackwell order.

\subsubsection{No fixed finite auxiliary family suffices}\label{BW-sec:hierarchy}
The completion proof needs new directions as the target is refined.
For a fixed refinement, all physical sources and all cost matrices still
produce transport walls from one finite set of normal directions. The
counterexample chooses an indefinite quadratic form positive on that
entire set. A smooth second-derivative perturbation then passes every
old comparison while reversing the sign along one new direction.

A general finite independent refinement has a law $\lambda=(\lambda_1,\ldots,\lambda_s)$ and posterior coordinates $k_i\lambda_b$. Zero auxiliary masses may be discarded.

\begin{theorem}[Smooth finite-family obstruction]\label{BW-thm:finite-family}
Let $q\ge3$, let $p$ be an interior point of $\Delta_q$, and fix any finite collection $\lambda^{(1)},\ldots,\lambda^{(L)}$ of finite auxiliary laws. For every neighborhood $U$ of $p$ in the simplex interior there are distinct posterior laws $\mu,\nu$ with barycenter $p$ and $C^\infty$ densities supported in $U$ such that
\begin{equation}\label{BW-eq:old-family}
 \mathcal C_{\nu^{\lambda^{(\ell)}}}(\rho)
 \subseteq\mathcal C_{\mu^{\lambda^{(\ell)}}}(\rho)
\end{equation}
for every $\ell$ and every centered integrable finite-dimensional source $\rho$, but another finite uniform refinement fails this inclusion, already for a scalar two-point source. For every fixed integer $s\ge0$ the two densities can be chosen arbitrarily close in $C^s$ norm.
\end{theorem}
\begin{proof}
Use affine coordinates $x=(k_1,\ldots,k_{q-1})$ on the simplex and put $d=q-1\ge2$. Consider any finite physical source and any cost matrix on the refined output. As a function of its prescribed target masses, the minimum transportation value is convex and piecewise affine. Its chamber walls are among the balance hyperplanes
\[
 \sum_{(i,b)\in A}k_i\lambda_b=\sum_{j\in B}a_j,
\]
where $a_j$ are the physical masses. Indeed, the feasibility inequalities for a basic transportation tree are the mass balances obtained by removing one tree edge. Degenerate optimal faces are covered by the same finite chamber complex.

The right side may vary with the physical source, but the normal of the left side belongs to a fixed finite set. In the affine coordinates it has the form
\[
 n=(a_1-a_q,\ldots,a_{q-1}-a_q),\qquad
 a_i=\sum_{b\in A_i}\lambda_b.
\]
Taking all subsets and all the prescribed auxiliary laws gives a finite set $\mathcal N$ of nonzero old normals, independent of the physical source and cost.

Choose a nonzero rational vector $v\in\R^d$ parallel to none of these normals. With $\widehat v=v/\|v\|$, choose $\delta>0$ small enough that
\begin{equation}\label{BW-eq:Q-signs}
 Q=I-(1+\delta)\widehat v\widehat v^{\mathsf T}
 \quad\text{satisfies}\quad n^{\mathsf T}Qn>0\ (n\in\mathcal N),
 \qquad v^{\mathsf T}Qv<0.
\end{equation}
This is possible because every old normal has a nonzero component perpendicular to $v$.

Choose a ball about $p$ inside $U$, a nonnegative smooth bump $\psi$ positive near $p$ and compactly supported in that ball, and a smooth probability density $f$ of mean $p$ which is strictly positive on $\supp\psi$ and is supported in a slightly larger concentric ball. For example, a radial smooth density about $p$ has the required mean. Define
\begin{equation}\label{BW-eq:smooth-perturb}
 d\mu=f(x)\,dx,\qquad d\nu=(f(x)+t\,Q:D^2\psi(x))\,dx.
\end{equation}
For all sufficiently small $t>0$ the second density is nonnegative. Integration by parts shows that the perturbation has zero mass and zero first moments. Thus both laws have barycenter $p$, and the densities approach one another in every fixed $C^s$ norm as $t\downarrow0$.

Let $W$ be any of the old finite-source transportation values. The distributional Hessian of a convex piecewise-affine function is supported on its walls. On each open facet it is a nonnegative scalar measure times $nn^{\mathsf T}$, where $n$ is a wall normal. This follows either by differentiating the jump of its gradient, which is normal to a common affine facet, or by its convexity on lines. Lower-dimensional facet intersections add no separate term to this derivative measure. Hence \eqref{BW-eq:Q-signs} and integration by parts give
\[
 \int W\,d\nu-\int W\,d\mu
   =t\int \psi\,Q:D^2W\ge0.
\]
Calibration support functions are negative transportation values after a choice of cost signs. This proves \eqref{BW-eq:old-family} for finite physical sources; the $L^1$ approximation in Lemma~\ref{BW-lem:posterior} extends it to every integrable source.

It remains to exhibit a failure. Since $v$ is rational, choose rational coefficients $a_i\in(0,1)$ with $(a_1-a_q,\ldots,a_{q-1}-a_q)=s_0v$ for a positive rational $s_0$. They are realized as a subset of copies in some uniform $M$-refinement. Let $r_0=a\cdot p$ and take the centered Bernoulli source $X=B-r_0$, with $\Pp(B=1)=r_0$. Its upper integrated quantile is
\[
 Q_X(s)=\min(r_0,s)-r_0s.
\]
The distributional Hessian of $Q_X(a\cdot k)$ is a negative multiple of $vv^{\mathsf T}$ on the hyperplane $a\cdot k=r_0$, which passes through the positive part of $\psi$. Therefore
\[
 \int Q_X(a\cdot k)\,d\nu-\int Q_X(a\cdot k)\,d\mu>0.
\]
This violates the necessary subset support inequality for the new refinement. The perturbation is consequently nonzero and the laws are distinct.
\end{proof}

\begin{corollary}[Every next level is strict, with smooth densities]\label{BW-cor:bounded-bits}
For every $M\ge1$, $q\ge3$, interior prior $p$, and neighborhood $U$ of $p$, the smooth local examples in Theorem~\ref{BW-thm:finite-family} can be chosen to pass every uniform refinement $1\le r\le M$ and fail at $M+1$. They can still be arbitrarily close in any fixed $C^s$ norm. The hierarchy by fair auxiliary bits is strict after each additional bit as well.
\end{corollary}
\begin{proof}
For a uniform $r$-refinement the old normals, after multiplication by $r$, are integer vectors $(c_1-c_q,\ldots,c_{q-1}-c_q)$ with $0\le c_i\le r$. Their coordinate range, including zero, is at most $r$. Choose the new direction $v=(M,M+1,0,\ldots,0)$. Since its two nonzero coordinates are coprime, a nonzero integer vector parallel to $v$ is an integer multiple of it, and its coordinate range is at least $M+1$. Thus none of the old normals is parallel to $v$. Use precisely this direction in the proof of Theorem~\ref{BW-thm:finite-family}.

At refinement $M+1$, take an event containing $M$ copies of label 1, every copy of label 2 and no others. Its affine normal is proportional to $v$, so the strict failing Bernoulli test in that proof occurs at this exact next size. For $M=2^b$, the same normal is realized with $2M$ copies by taking $M,M+1,0,\ldots,0$ copies of the labels. Thus the failure can be placed after exactly one additional fair bit.
\end{proof}

For $q=3$ and the uniform prior, the same construction gives an explicit polynomial-density example. Put
\[
 v=(M,M+1),\quad w=(M+1,-M),\quad
 \epsilon_M=\frac1{2M^2(2M+1)^2},\quad Q=ww^{\mathsf T}-\epsilon_Mvv^{\mathsf T}.
\]
Every nonzero old normal has $(w\cdot n)^2\ge1$ and $|v\cdot n|\le M(2M+1)$, so $n^{\mathsf T}Qn\ge1/2$, while $v^{\mathsf T}Qv<0$. With $\psi=(1-\|x\|^2)^4_+$, $L=\sum_{ij}|Q_{ij}|$ and $\eta=(112L)^{-1}$, its planar densities are
\[
 f_\mu=\frac3\pi(1-\|x\|^2)^2_+,\qquad
 f_\nu=f_\mu+\frac{3\eta}{\pi}Q:D^2\psi.
\]
The bracketed relative perturbation is bounded by $56\eta L=1/2$, so positivity is explicit. After the map $k=(1/3+sx_1,1/3+sx_2,1/3-s(x_1+x_2))$ with $s=1/12$, the next-level support gap is
\[
 \frac{256\eta s\epsilon_M\|v\|^3}{105\pi(M+1)}>0.
\]
This follows by the same integration by parts and $\int_{-1}^1(1-t^2)^4dt=256/315$. The polynomial densities are $C^1$ across their support boundary. The smooth local construction above strengthens regularity and permits arbitrary interior priors; the polynomial formula gives an explicit elementary gap.

The obstruction has a geometric interpretation with a precise content. A finite target refinement determines only finitely many transportation-wall directions. An indefinite perturbation of the Hessian can be positive on all of those directions and negative on a new one. The growing family in Theorem~\ref{BW-thm:completion} is what makes those directions sufficiently rich to test arbitrary convex functions.

\subsubsection{One simulator for all bounded decision problems}
For a posterior law $\mu$ of prior $p$, let its canonical experiment
have joint law $P_\mu(\{i\}\times dx)=x_i\mu(dx)$.
A kernel $L$ from the observation of $\nu$ to that of $\mu$ preserves
the target's prior-predictive marginal exactly when $\nu L=\mu$.

\begin{theorem}[Approximate simulation with the exact observation marginal]
\label{thm:blackwell-simulator}
With total variation defined by $\TV(P,Q)=\sup_A|P(A)-Q(A)|$,
\begin{equation}\label{eq:exact-deficiency}
 \inf_{L:\nu L=\mu}\TV(P_\mu,P_\nu L)
       =\tfrac12 d_{\rm cx}(\mu,\nu).
\end{equation}
The infimum is attained. Consequently the refinement comparisons in
Theorem~\ref{R8-thm:blackwell-rate} construct one kernel with joint error
$O_q((e_M+M^{-1})^{2/(q+1)})$. Every decision rule with loss in $[0,1]$
transfers through that kernel with at most this change in prior-averaged
risk. The kernel is chosen before the loss and the decision rule.
\end{theorem}
\begin{proof}
Barycentric weak-transport duality \cite{GRST} on the compact simplex gives
\[
 d_{\rm cx}(\mu,\nu)=\min_{\pi\in\Pi(\mu,\nu)}
       \E_\pi\|X-\E_\pi[Y\mid X]\|_1.
\]
Disintegrate $\pi(dx,dy)=\nu(dy)L(y,dx)$ and put
$b(x)=\E[Y\mid X=x]$. Then $\nu L=\mu$ and
$P_\nu L(\{i\}\times dx)=b_i(x)\mu(dx)$, so
\[
 \TV(P_\mu,P_\nu L)=\tfrac12\int\sum_i|x_i-b_i(x)|\,\mu(dx).
\]
Conversely every such kernel determines this coupling. The displayed
identities prove equality and attainment. Total variation bounds the
change in expectation of every function valued in $[0,1]$, proving
the risk assertion. For a positive prior, worst-state error is at most
the joint error divided by $\min_i p_i$.
\end{proof}


\section{Exact posterior reconstruction and the sharp directed inverse}
\label{R11-sec:inverse}
Varying source weights on one fixed directed cycle produces a cost profile whose derivative recovers the posterior distribution function, including boundary mass. The three-state example shows why orientation is essential: symmetric costs agree while the directed costs differ. The chamber identity gives uniqueness and the sharp inverse modulus on the closed simplex, separating full reconstruction from extension-order comparison.

\subsection{Reconstructing the reference experiment}\label{sec:posterior-reconstruction}
The cyclic transport values used to compare experiments also determine
an experiment exactly. This section fixes one directed cost matrix and
varies only the independent source weights. Its one-sided directional
derivative is a posterior distribution function. Keeping that derivative
at atoms will be essential for recovery on the closed simplex.

Let $q\ge2$, $d=q-1$, and
\[
 \Delta_q=\{k\in[0,1]^q:\textstyle\sum_i k_i=1\},\qquad
 C=\{x\in\R^d:0<x_1<\cdots<x_d<1\}.
\]
The cumulative map and its inverse are
\begin{equation}\label{CYC-eq:cumulative}
 x_i=\sum_{l=1}^i k_l,\qquad
 k=(x_1,x_2-x_1,\ldots,x_d-x_{d-1},1-x_d).
\end{equation}
The closed simplex maps to $\overline C$. All experiments have finitely many states and standard Borel observation spaces. When two experiments are compared, their prior $p$ is the same and $p_i>0$.

We use total variation distance in the convention $\TV(P,Q)=\frac12\|P-Q\|_1$ and Wasserstein distance $W_1$ on the simplex with cost $\|k-k'\|_1$. Constants depending only on $q$ may change from line to line.

\subsubsection{One directed cycle}
Fix the cyclic order $1,\ldots,q,1$, with ground cost
\[
 c(l,i)=(i-l)\bmod q\in\{0,\ldots,q-1\}.
\]
For $a,k\in\Delta_q$, put
\begin{equation}\label{CYC-eq:transport}
 W(a,k)=\min\left\{\sum_{l,i}c(l,i)\pi_{li}:\pi\ge0,
              \ \pi\1=a,\ \pi^T\1=k\right\}.
\end{equation}
If $\nu$ is a probability law on $\Delta_q$, our observable is
\begin{equation}\label{CYC-eq:J}
 J_\nu(a)=\int W(a,k)\nu(\dd k),\qquad a\in\Delta_q^\circ.
\end{equation}
No unknown posterior can be selected by the observer; only the source $a$ is selected.

\begin{lemma}[The cycle formula]\label{CYC-lem:cycle}
Let $x$ and $y$ be the cumulative vectors of $k$ and $a$. For
$f(v)=\max(0,v_1,\ldots,v_d)$,
\begin{equation}\label{CYC-eq:cycle}
 W(a,k)=qf(x-y)-\sum_{i=1}^d(x_i-y_i).
\end{equation}
This identity holds on the closed simplex.
\end{lemma}
\begin{proof}
Apply Lemma~\ref{BW-lem:cycle} with the source and target prefixes in
the displayed cumulative coordinates. Its nonnegative cyclic flow is
$g_i=g_q+y_i-x_i$, so the least $g_q$ is $\max(0,x_i-y_i)$.
Summing the flows gives \eqref{CYC-eq:cycle}. The path decomposition
in that lemma permits zero source and target masses, hence the identity
holds on every face of the simplex.
\end{proof}

\subsubsection{A three-state example}
\label{R6-sec:three-state-cycle}
The scalar extension inequalities in Theorem~\ref{NLI-thm:global}
depend only on the distributions of posterior subset sums. A directed
cycle can distinguish laws having all those distributions in common.
Let $p=(1/3,1/3,1/3)$ and $0<h<1/3$. Give equal mass to the three
points in each row:
\[
\begin{array}{c|ccc}
 \nu_+&p+h(-1,1,0)&p+h(0,-1,1)&p+h(1,0,-1)\\
 \nu_-&p+h(-1,0,1)&p+h(0,1,-1)&p+h(1,-1,0).
\end{array}
\]
Both priors are $p$. Each individual posterior coordinate has the same
three-point law; each two-coordinate sum is the complement of the
remaining coordinate. Thus every scalar extension set agrees for the
two experiments. Their symmetrized cyclic costs also agree.

For the directed cost at source $p$, the three values for $\nu_+$ are
all $h$, and the three values for $\nu_-$ are all $2h$. Hence
\[
 J_{\nu_+}(p)=h,\qquad J_{\nu_-}(p)=2h.
\]
Indeed, in cumulative coordinates the positive-row differences are
$(-h,0),(0,-h),(h,h)$; applying \eqref{CYC-eq:cycle} gives $h$ in
each case. Their negatives give $2h$. The equality of the symmetrized
costs at every source, and the higher-alphabet construction that explains
it, follow from the finite-difference argument below.

\subsubsection{The derivative of the directed cost}
The cycle formula expresses the observed profile as an expected maximum plus an affine term. Moving source mass from the last state to the first raises every cumulative source threshold equally. The right derivative of the positive maximum records whether any posterior coordinate exceeds its threshold. Its complement is the joint lower-orthant distribution function, including equality at atoms. This recovers dependence between coordinates as well as their separate laws.

Write $\lambda$ for the law of the cumulative vector $X$. Define
\[
 H_\lambda(y)=\int f(x-y)\lambda(\dd x),\qquad
 L=\left(\sum_{i=1}^d\partial_i\right)\prod_{i=1}^d\partial_i.
\]
The differential operator has order $q$, not $d$.

\begin{lemma}[Directional and distributional inverses]\label{CYC-lem:inverse}
For every $y\in\R^d$,
\begin{equation}\label{CYC-eq:cdf}
 \partial_{\1,+}H_\lambda(y)=\Pp(X\le y)-1.
\end{equation}
Consequently $LH_\lambda=\lambda$ as distributions on $\R^d$. If $a$ is the source corresponding to $y\in C$, then
\begin{equation}\label{CYC-eq:direct-cdf}
 \Pp(X\le y)=\frac1q\left(1+
 \lim_{t\downarrow0}\frac{J_\nu(a+t(e_1-e_q))-J_\nu(a)}t\right).
\end{equation}
In particular $L(J_\nu/q)=\lambda$ inside $C$, regarding $J_\nu$ as a function of $y$.
\end{lemma}
\begin{proof}
For $M=\max_i(x_i-y_i)$, the right derivative at zero of $\max(0,M-t)$ is $-\1_{\{M>0\}}$. The difference quotient is bounded by one, so dominated convergence gives \eqref{CYC-eq:cdf}, with the inclusive lower orthant at its atoms. Applying $\partial_1\cdots\partial_d$ to the distribution function gives $\lambda$. Formula~\eqref{CYC-eq:cycle} gives
$J_\nu(a)=qH_\lambda(y)-\sum_i(\E X_i-y_i)$; the directional derivative of its affine term is $d$. Thus $q(F-1)+d=qF-1$, proving \eqref{CYC-eq:direct-cdf}. The affine term is annihilated by $L$.
\end{proof}

\begin{lemma}[Extension across the chamber boundary]\label{R15-lem:chamber}
For $y\in\R^d$ set $t_i=\min_{j\ge i}y_j$, $b=(-t_1)_+$ and
$z_i=\min(1,t_i+b)$. Then $z\in\overline C$ and
\begin{equation}\label{R15-eq:chamber}
 f(x-y)=b+f(x-z)\quad(x\in\overline C).
\end{equation}
Consequently the values of $H_\lambda$ on the closed chamber determine its
values on all of $\R^d$.
\end{lemma}
\begin{proof}
Suffix minima preserve the maximum: if $t_i=y_j$ for $j\ge i$, then
$x_i-t_i\le x_j-y_j$. If $t_1<0$, the maximum is at least
$x_1-t_1\ge b$; subtracting $b$ therefore moves it inside the positive part.
For $t_1\ge0$ this step changes nothing. Finally $t_i+b\ge0$ is increasing,
and clipping coordinates above one leaves their nonpositive contributions
unchanged. This proves \eqref{R15-eq:chamber}.
\end{proof}

\subsubsection{Closed-simplex recovery without posterior positivity}
\begin{theorem}[Complete recovery from one fixed cycle]\label{CYC-thm:boundary}
For arbitrary probability laws $\nu,\nu'$ on the closed simplex $\Delta_q$,
\[
 J_\nu(a)=J_{\nu'}(a)\quad\hbox{for every }a\in\Delta_q^\circ
 \quad\Longleftrightarrow\quad \nu=\nu'.
\]
The same fixed cost matrix is used throughout. Neither coordinate positivity nor absolute continuity nor finite support is required.
\end{theorem}
\begin{proof}
The cost profile is continuous up to the boundary. At the source $e_1$,
whose cumulative vector is $\1$, $H_\lambda(\1)=0$. Hence for two laws
\[
 q(H_\lambda-H_{\lambda'})(y)
   =J_\nu(a(y))-J_{\nu'}(a(y))-J_\nu(e_1)+J_{\nu'}(e_1)
 \quad(y\in\overline C).
\]
Equality of the profiles, Lemma~\ref{R15-lem:chamber} and the global
identity $LH_\lambda=\lambda$ prove equality of the laws.
\end{proof}
The chamber identity determines the values outside the observable source
region without losing facial atoms. The one-sided derivative remains useful
for the finite exact algorithm below.

\subsubsection{A finite exact reconstruction algorithm}
Let
\[
 \Delta_q[M]=\{k\in\Delta_q:Mk_i\in\mathbb Z\ \hbox{for all }i\}.
\]
An exact cost oracle returns $J_\nu(a)$ at the requested rational interior source. For rational posterior weights, these outputs and all operations below are rational.

\begin{theorem}[Boundary atoms with at most $2^q$ queries]\label{CYC-thm:grid}
Suppose $\nu$ is supported on $\Delta_q[M]$. For every specified candidate $k\in\Delta_q[M]$, its mass $\nu\{k\}$ can be recovered with at most $2^q$ exact queries, all at strictly interior rational sources of the one fixed directed cycle. This includes vertices and all proper faces. Recovering all atoms uses at most
\[
 2^q\binom{M+q-1}{q-1}
\]
queries. The bound is an oracle/output bound, not a polynomial-in-$q$ claim.
\end{theorem}
\begin{proof}
Choose a label $j$ with $k_j>0$, rotate it to the final position, and write $x$ for the cumulative candidate. Then $x_d\le1-1/M$. The atom is the mass of
\[
 \prod_{i=1}^d\left(x_i-\frac1{2M},x_i+\frac1{2M}\right],
\]
which is an alternating sum of $2^d$ lower-orthant probabilities. For each corner $z_i=x_i\pm1/(2M)$, return zero when some $z_i<0$. Otherwise set
\[
 t_i=\min_{r\ge i}z_r,\qquad
 \epsilon=\frac1{8qM},\qquad y_i=t_i+i\epsilon.
\]
Both $y$ and $y+\epsilon\1$ lie in $C$. Every $t_i$ is a half-grid point. Moving it by at most $(d+1)\epsilon=1/(8M)$ crosses no possible cumulative posterior coordinate. Hence the lower-orthant probability is constant on the directional segment and
\[
 F_X(z)=F_X(y)=\frac{J_\nu(a(y+\epsilon\1))-J_\nu(a(y))+\epsilon}{q\epsilon}.
\]
Two calls give each corner probability. Cyclically unrotate each source before sending it to the fixed oracle. Inclusion--exclusion gives the candidate mass exactly.
\end{proof}

\subsubsection{A metric on posterior laws and its sharp inverse modulus}
Define
\begin{equation}\label{CYC-eq:metric}
 d_{\rm cyc}(\nu,\nu')=\sup_{a\in\Delta_q}|J_\nu(a)-J_{\nu'}(a)|.
\end{equation}
Continuity makes the supremum over the open simplex equivalent.

\begin{proposition}[Forward stability and topology]\label{CYC-prop:topology}
The function $d_{\rm cyc}$ is a metric on $\cP(\Delta_q)$, and
\begin{equation}\label{CYC-eq:forward}
 d_{\rm cyc}(\nu,\nu')\le\frac{q-1}{2}W_1(\nu,\nu').
\end{equation}
It induces the weak topology, equivalently the $W_1$ topology on this compact simplex.
\end{proposition}
\begin{proof}
Given an optimal coupling of source $a$ and target $k$, minimally recouple its target to $k'$. At most $\TV(k,k')$ mass changes label, and each changed cost changes by at most $q-1$. Applying the argument in both directions gives
$|W(a,k)-W(a,k')|\le(q-1)\TV(k,k')$.
Integrate against a coupling of $\nu$ and $\nu'$ to get \eqref{CYC-eq:forward}. Theorem~\ref{CYC-thm:boundary} proves definiteness. The map from the compact space of probability laws, in the weak topology, into this metric space is continuous and injective; its inverse on its image is continuous. Alternatively take a weakly convergent subsequence and use \eqref{CYC-eq:forward} and uniqueness of the limit.
\end{proof}

The inverse differentiates the cost profile $q$ times. A Lipschitz Wasserstein test already has one bounded weak derivative; smoothing it at scale $h$ pays for the remaining $q-1$. A profile error $\eta$ therefore contributes $\eta h^{1-q}$, alongside smoothing error $h$. Balancing these terms gives $h=\eta^{1/q}$. Oscillations with zero mass and first moments attain this balance, preserving a common prior and smooth interior densities. The sharpness construction thus tests the same inverse mechanism as the upper bound.

\begin{theorem}[Sharp inverse on the closed simplex]\label{CYC-thm:sharp}
For all probability laws on $\Delta_q$,
\begin{equation}\label{CYC-eq:holder}
 W_1(\nu,\nu')\le C_q\,d_{\rm cyc}(\nu,\nu')^{1/q}.
\end{equation}
No common prior, interior support or density is required. The exponent is
optimal even for smooth bounded densities with a common prior and common
compact support in the open simplex.
\end{theorem}
\begin{proof}
Put $\eta=d_{\rm cyc}(\nu,\nu')$. The preceding proof and
Lemma~\ref{R15-lem:chamber} give
$\|H_\lambda-H_{\lambda'}\|_{L^\infty(\R^d)}\le2\eta/q$;
for equal priors the factor is $1/q$. Extend a normalized Lipschitz
Kantorovich test from $\overline C$ and multiply by one smooth cutoff equal
to one on a neighborhood of the whole chamber. Its mollification $g_h$
satisfies
\[
 \|g_h-g\|_{L^\infty(\overline C)}\le C_qh,
 \qquad \|L^*g_h\|_1\le C_qh^{1-q}.
\]
Differentiate once on the Lipschitz test and $q-1$ times on the mollifier.
The global distributional inverse now gives an error at most
$C_q(h+\eta h^{1-q})$. Choose $h=\eta^{1/q}$ for $\eta\le1$, and use
bounded diameter otherwise. The cumulative coordinates are bi-Lipschitz.

For sharpness choose nonzero $\chi\ge0$ smooth with compact support in $C$,
and a smooth probability density $r$ positive near that support. Set
$\phi_n(y)=n^{-q}\chi(y)\cos(n\1\cdot y)$ and $\rho_n=L\phi_n$.
These have uniformly bounded supremum norm and zero mass and first moments.
For a fixed sufficiently small $\tau>0$, $r\pm\tau\rho_n$ are probability
densities with a common prior. Integration by parts in the fundamental
solution gives
\begin{equation}\label{CYC-eq:fundamental}
 \int f(x-y)L\phi_n(x)\,dx=\phi_n(y).
\end{equation}
The $q$ sign changes from $x$-differentiation cancel those from integration
by parts. Thus the profile difference is $2\tau q\phi_n$, of size
$O(n^{-q})$. The leading term of $\rho_n$ is
$d\chi(y)\cos(n\1\cdot y+q\pi/2)+O(n^{-1})$.
A constant multiple of $n^{-1}\cos(n\1\cdot y+q\pi/2)$ is Lipschitz and
pairs with $2\tau\rho_n$ to give at least $c/n$, since the squared cosine
has limiting average $1/2$ on $\chi$. No exponent greater than $1/q$ can hold.
\end{proof}

\begin{corollary}[Weighted cycles, faces and endpoint laws]\label{R15-cor:faces}
Positive cycle lengths $\ell_i$ of total length $A$ replace the cost by
$A f(x-y)-\sum_{i<q}\ell_i(x_i-y_i)$. Theorem~\ref{CYC-thm:sharp} holds
with a constant depending on $q,A$; for equal priors the profile difference
divided by $A$ is independent of the edge lengths. If one law lies on a
known face with $r$ active states, the other may be arbitrary and the sharp
exponent is $1/r$. In particular, writing $\eta=d_{\rm cyc}(\nu,\nu')$,
\[
 W_1(\delta_a,\nu')\le2d_{\rm cyc}(\delta_a,\nu'),\qquad
 W_1\!\left(\sum_i p_i\delta_{e_i},\nu'\right)\le6\eta,
\]
where the last constant is $4$ for equal priors. These endpoint exponents
are sharp; optimality of their constants is not asserted.
\end{corollary}
\begin{proof}
The nonnegative edge-flow proof of Lemma~\ref{BW-lem:cycle} gives the
weighted formula. Vertex queries give, with cyclic indices,
\[
 J_\nu(e_{i+1})-J_\nu(e_i)=qp_i-1,
 \qquad |p_i-p_i'|\le2\eta/q.
\]
If $\nu$ lies on an $r$-state face, the expected mass of $\nu'$ outside it
is at most $2\eta$. Moving that mass to one face vertex costs at most
$4\eta$ in $W_1$. On the face the original cycle is the induced weighted
$r$-cycle. Apply the theorem there and the forward Lipschitz estimate.
The case $r=1$ follows from the vertex identity.

For a point law, $W(a,k)\ge\|a-k\|_1/2$ gives the displayed bound at
source $a$. For the fully revealing law use $u=\1/q$:
\[
 (q-1)/2-W(u,k)=\min_j\sum_i k_i c(i,j)\ge1-\max_i k_i.
\]
Hence $\E_{\nu'}(1-\max_i k_i)\le\eta$. Choosing a vertex $e_I$ with
probabilities $k_i$ costs
$2\E(1-\sum_i k_i^2)\le4\eta$ and produces prior $p'$; changing it to
$p$ costs $\|p-p'\|_1\le2\eta$. Oscillatory examples on the face give
sharpness for $r\ge2$, and mixtures of distinct point or vertex laws give
the endpoint exponent one.
\end{proof}

\subsubsection{Orientation, parity, and hidden dependence}
The symmetrized observable below is the sum of two separately optimized transport costs. It is not optimal transport for the entrywise symmetrization of the cost matrix:
\[
 S_\nu(a)=\int\bigl[W(a,k)+W(k,a)\bigr]\nu(\dd k).
\]
By the cycle formula,
\begin{equation}\label{CYC-eq:sym}
 S_\nu(a(y))=q\int\bigl[f(x-y)+f(y-x)\bigr]\lambda(\dd x).
\end{equation}

\begin{proposition}[Parity of the inverse]\label{CYC-prop:parity}
As distributions on $\R^d$,
\[
 L\int\bigl[f(x-y)+f(y-x)\bigr]\lambda(\dd x)
       =\bigl(1+(-1)^q\bigr)\lambda.
\]
For even $q$, $S_\nu$ determines every posterior law carried by the open simplex. For odd $q$, the right side vanishes; the following theorem exhibits actual nonidentifiability, not merely failure of this differential formula.
\end{proposition}
\begin{proof}
Reflection contributes $(-1)^q$ to a derivative of total order $q$. For even $q$, apply $L/(2q)$ to \eqref{CYC-eq:sym} on $C$. If the posterior has no boundary mass, that recovers its entire law.
\end{proof}

\begin{theorem}[Explicit odd-alphabet indistinguishability]\label{CYC-thm:parity}
Let $q\ge3$ be odd and $0<h<1/q$. There are two different finitely supported posterior laws $\nu_+,\nu_-$, each with $2^{q-1}-1$ atoms and uniform prior, satisfying all of the following.

Every posterior coordinate is at least $1/q-h$. Every subset sum $\sum_{i\in A}K_i$ has the same distribution under the two laws. All polynomial moments of total degree less than $q$ agree. The symmetrized costs $S_{\nu_+}(a)$ and $S_{\nu_-}(a)$ agree for every $a\in\Delta_q$. Nevertheless, at the uniform source $p=(1/q,\ldots,1/q)$,
\begin{equation}\label{CYC-eq:gap}
 J_{\nu_+}(p)-J_{\nu_-}(p)=-\frac{qh}{2^{q-1}-1}.
\end{equation}
In particular, equality of every scalar calibration body and equality of all these low-order moments do not imply equality of vector calibration bodies.
\end{theorem}
\begin{proof}
Let $v_0=\1\in\R^d$, $v_i=e_i$ for $1\le i\le d$, and set
\[
 (x_0)_i=i/q-h,\qquad
 \rho=\sum_{b\in\{0,1\}^q}(-1)^{q-|b|}
     \delta_{x_0+h\sum_{r=0}^d b_rv_r}.
\]
Only one point is represented twice: $x_0+h\1$ arises from
$(b_0,b_1,\ldots,b_d)=(1,0,\ldots,0)$ and $(0,1,\ldots,1)$. Since $q$ is odd, the two coefficients cancel. All other points are distinct and have coefficients $\pm1$. The positive and negative masses both equal
$m=2^{q-1}-1$. Divide the Jordan parts by $m$ and map back to the simplex to define $\nu_+$ and $\nu_-$.

Each posterior coordinate is $1/q$ plus an element of $\{-h,0,h\}$. The order-$q$ difference $\rho$ annihilates all polynomials of lower degree; symmetry of the full binary cube then gives each law mean $(1/q,\ldots,1/q)$.

For a subset $A$, put $s_i=\1_{\{i\in A\}}$. Its sum is $s_q+\sum_{i<q}(s_i-s_{i+1})x_i$, whose increments along the stencil are $s_1-s_q$ and $s_i-s_{i+1}$. On an odd binary cycle some increment is zero. Pairing the corresponding two stencil terms makes the entire pushforward of $\rho$ vanish, proving equality of the subset-sum distributions.

By Proposition~\ref{CYC-prop:parity}, the mixed directional derivative of $f(x-y)+f(y-x)$ along $v_0,\ldots,v_d$ vanishes. Convolution with the $q$ interval measures on $[0,h]$ gives its zero mixed finite difference, pointwise by continuity. This proves equality of the symmetrized costs. At $y_i=i/q$, however, $x_i-y_i=h(b_0+b_i-1)$. The maximum is zero for $b_0=0$ and is $h$ for $b_0=1$ except when every other bit is zero. Its signed sum is $-h$. The affine term integrates to zero, proving~\eqref{CYC-eq:gap}.

Finally, for a scalar source with quantile $Q_X$, order any scalar coefficients as $u_{\pi(1)}\le\cdots\le u_{\pi(q)}$. Conditional rearrangement gives the maximal reward
\[
 \sum_{i=1}^q u_{\pi(i)}
 \int_{t_{i-1}}^{t_i}Q_X(s)\dd s,\qquad
 t_i=\sum_{l\le i}k_{\pi(l)}.
\]
After telescoping, its expectation depends only on distributions of subset sums. Thus all scalar bodies agree. The directed transport witness is the support of a centered simplex-valued source of dimension at most $q-1$, so some vector bodies differ.
\end{proof}

\begin{corollary}[Smooth indistinguishability]\label{CYC-cor:smooth}
For every odd $q\ge3$ there are distinct smooth posterior densities with compact interior support having the same subset-sum distributions, moments through degree $q-1$, and symmetrized cyclic costs. They can also be made strictly positive throughout the open simplex by adding the same smooth positive baseline density to both laws.
\end{corollary}
\begin{proof}
Convolve both cumulative-coordinate laws with the same sufficiently small smooth compactly supported probability kernel. Their supports stay in $C$. Projection equalities, moments through degree $q-1$ and equality of the symmetrized potentials survive convolution. The directed gap converges to~\eqref{CYC-eq:gap}, so the laws remain distinct. Mixing both with the same positive smooth density preserves these conclusions and makes them positive throughout $C$.
\end{proof}



\appendix
\clearpage

\section{Scalar overlap and the column correction}\label{app:overlap}
This appendix proves the scalar height input \eqref{eq:scalar-input}.
Write $\phi,\Phi$ for the standard normal density and distribution
function. The weighted overlap of two normal laws is
\begin{align}\label{R-eq:gaussianrisk}
 \mathcal R(t,b)&=e^b\Phi(-\sqrt t-b/\sqrt t)
                  +e^{-b}\Phi(-\sqrt t+b/\sqrt t),\\
 \Pi(A,B;a)&=\sqrt{AB}\,\mathcal R(a^2,\tfrac12\log(A/B)).\notag
\end{align}
Set $\mathcal R(0,b)=e^{-|b|}$ and $\Pi=0$ when $AB=0$. For a density
$q$ on $(-3,3)$ put
\begin{equation}\label{R-eq:G}
 \G_k(q)=\frac14\iint\Pi(q(s),q(t);k|s-t|)\,ds\,dt.
\end{equation}
Differentiation gives
\begin{align}\label{R-eq:derivatives}
 \partial_t\mathcal R(t,b)&=-\frac{\phi(\sqrt t)}{\sqrt t}
                   e^{-b^2/(2t)},\\
 \partial_a\Pi(A,B;a)&=-2\sqrt{AB}\,\phi(a)
                   e^{-\log^2(A/B)/(8a^2)}.\notag
\end{align}

\subsection{The likelihood transform}
For a pair $P,Q$ with densities $p,q$, define
\[
 R_{P,Q}(b)=\int\min\{e^bp,e^{-b}q\}.
\]
On the common support put $\ell=\tfrac12\log(p/q)$ and let $\nu$ be the image of $\sqrt{pq}\dd x$ under $\ell$. Its Fourier transform is
\[
 J(\omega)=\int e^{i\omega\ell}\nu(\dd\ell).
\]
The finite measure $\nu$ has mass at most one; parts of $P$ or $Q$ outside the common support contribute zero to every risk. The pointwise identity
\[
 \min\{e^bp,e^{-b}q\}=\sqrt{pq}\,e^{-|b+\ell|}
\]
gives
\begin{equation}\label{R-eq:fourier-risk}
 R_{P,Q}(b)=\frac1\pi\int_{\R}\frac{e^{i\omega b}J(\omega)}{1+\omega^2}\dd\omega.
\end{equation}
All integrals here converge absolutely. Independence multiplies the transforms because half log likelihoods add and Hellinger measures convolve. For the translated pair $f(x),f(x+w)$ write its transform as $J_{f,w}$ and its weighted overlap as
\[
 O_f(w;A,B)=\int\min\{Af(x),Bf(x+w)\}\dd x.
\]

For the translated cosine pair with angular separation $s\in(0,\pi)$,
solving the half-log-likelihood equation and taking the Jacobian gives
\begin{equation}\label{R-eq:likelihood}
 \nu_s(\ell)=\frac{\sin^3s}{2\pi(\cosh\ell-\cos s)^2}.
\end{equation}
The elementary product
\begin{equation}\label{R-eq:trigproduct}
 \frac{\cosh x-\cos s}{1-\cos s}
 =\prod_{j\in\mathbb Z}\left(1+\frac{x^2}{(s+2\pi j)^2}\right)
\end{equation}
follows by factoring $T_N(z)-\cos s$ at
$z=\cos((s+2\pi j)/N)$ and evaluating its ratio at
$z=\cosh(x/N)$ and $z=1$. Pairing indices near the two ends gives
locally uniform convergence, since the logarithms have a summable
$O((1+|j|)^{-2})$ majorant from the elementary lower bound for sine.

Let $G_j$ be independent gamma variables of shape two and rate one,
and put $Y=\sum_{j\in\mathbb Z}G_j/(s+2\pi j)^2$.
The gamma integral gives $\nu_s(\ell)=\nu_s(0)\E e^{-Y\ell^2}$.
Reweight the law of $Y$ by
$\nu_s(0)\sqrt{\pi/Y}e^{1/(4Y)}$ and set $T=1/(2Y)$.
The resulting measure $\mu_s$ satisfies
\begin{equation}\label{R-eq:clock}
 \nu_s(d\ell)=\int e^{-t/2}N(0,t)(d\ell)\,\mu_s(dt).
\end{equation}
Its mass is $\int e^\ell\nu_s(d\ell)\le1$. Put the missing mass at
$t=\infty$, where the overlap is zero. Thus the weighted cosine
risk equals $\E\mathcal R(T,b)$; independent coordinates add their
times, with infinity absorbing.
Direct integration gives
\begin{equation}\label{R-eq:Js}
 J_s(\omega)=
 \frac{\omega\sin s\cosh((\pi-s)\omega)+\cos s\sinh((\pi-s)\omega)}
 {\sinh(\pi\omega)}.
\end{equation}
The value at $\omega=0$ is interpreted continuously. One way to verify the integral without contour methods is to differentiate \eqref{R-eq:likelihood} under the integral and integrate twice by parts in $\ell$. Both sides then solve
$J''-2\cot s\,J'=(1+\omega^2)J$ with $J(0)=1,J(\pi)=0$. The boundary limits follow directly from \eqref{R-eq:likelihood}. The elementary maximum principle gives uniqueness on $(0,\pi)$, proving the formula.

For $0<s\le\pi/2$,
\[
 -\partial_s\log J_s(\omega)
 =\frac{1+\omega^2}{\cot s+\omega\coth((\pi-s)\omega)}
 \le(1+\omega^2)s.
\]
Indeed, $x\coth x\ge1$, and
$\cot s+1/(\pi-s)\ge1/s$. For the latter inequality multiply by
$s(\pi-s)\sin s$; the resulting numerator vanishes at $0,\pi/2$, and its derivative is $[2-s(\pi-s)]\sin s$, with one sign change. Integrating proves
\begin{equation}\label{R-eq:laplace}
 J_s(\omega)\ge e^{-(1+\omega^2)s^2/2}.
\end{equation}
Through \eqref{R-eq:clock}, a product with $a^2=\sum_i s_i^2\le\pi^2/4$ consequently satisfies
\begin{equation}\label{R-eq:time-laplace}
 \E e^{-\lambda T}\ge e^{-\lambda a^2},\qquad \lambda\ge\tfrac12.
\end{equation}

\subsection{Restricted prior comparison}
\begin{lemma}\label{R-lem:tangent}
If $0<a\le\pi/2$ and $|b|\le11a/25$, there is $\lambda\ge1/2$ such that
\begin{equation}\label{R-eq:tangent}
 \mathcal R(t,b)\ge\mathcal R(a^2,b)e^{-\lambda(t-a^2)}
 \qquad(t\ge0).
\end{equation}
\end{lemma}
\begin{proof}
Let $g_b(t)=(2\pi t)^{-1/2}e^{-t/2-b^2/(2t)}$. Equation~\eqref{R-eq:derivatives} gives
$\mathcal R(t,b)=\int_t^\infty g_b(u)\dd u$. Write
$\eta_b=g_b/\mathcal R$ and
$\xi_b=-(\log g_b)'=1/2+1/(2t)-b^2/(2t^2)$.
Then $\eta_b'=\eta_b(\eta_b-\xi_b)$. For $t\ge2b^2$, $\xi_b$ is positive and decreasing. Hence
\[
 \frac{\mathcal R(t,b)}{g_b(t)}
 =\int_0^\infty\exp\left(-\int_t^{t+u}\xi_b(w)\dd w\right)\dd u
 \ge\frac1{\xi_b(t)},
\]
so the hazard decreases there. Below $2b^2$, any zero of
$\eta_b-\xi_b$ is a strict downward crossing, since $\xi_b'>0$.
The hazard therefore has one maximum and then decreases to $1/2$.

Take $\lambda=\eta_b(a^2)$. With $\kappa=|b|/a\le11/25$, the Mills integral gives
\[
 \mathcal R(a^2,\kappa a)
 =\phi(a)e^{-\kappa^2/2}\{m(a+\kappa)+m(a-\kappa)\},
 \quad m(x)=\int_0^\infty e^{-xu-u^2/2}\dd u.
\]
The inequality $\cosh(\kappa u)\le e^{u^2/2}$ implies that the braces are at most $2/a$, so $\lambda\ge1/2$. Since $a^2>2b^2$, the only candidates for the global minimum of
$\log\mathcal R(t,b)+\lambda t$ are $t=0$ and $t=a^2$. The latter has the smaller value precisely when
\begin{equation}\label{R-eq:endpoint}
 \log\mathcal R(a^2,\kappa a)+
 \frac{a\phi(a)e^{-\kappa^2/2}}{\mathcal R(a^2,\kappa a)}+\kappa a\le0.
\end{equation}
The rational checker in Appendix~\ref{R-app:arithmetic} encloses this expression divided by $a$ on
$0<a\le11/7$, $0\le\kappa\le11/25$. It covers the rectangle by $7939$ accepted rational subrectangles and gives a largest upper endpoint less than $-0.0000021803913799612835$. The arithmetic and treatment of $a=0$ are described in Appendix~\ref{R-app:arithmetic}. Thus \eqref{R-eq:endpoint} holds on a rectangle containing the required domain. When $b=0$, the same conclusion follows directly from log convexity of the Gaussian tail as a positive exponential mixture. The limit $t=\infty$ is harmless.
\end{proof}

\begin{proposition}\label{R-prop:forward}
Suppose $R\ge6$, $k=\pi/(2R)$, and $q$ is a positive auxiliary density with
$\Lip(\log q)\le22k/25$. Then for every $m$ and every $\norm v_2\le1$,
\begin{equation}\label{R-eq:forward}
 H_{f_{R,m}}(v;q)\ge\G_k(q).
\end{equation}
\end{proposition}
\begin{proof}
For $\norm v_2=1$, the pair of positions $s,t$ gives cosine shifts
$s_i=k|s-t||v_i|$, with $a^2=\sum s_i^2=k^2(s-t)^2\le\pi^2/4$.
The prior imbalance satisfies
$|b|=\tfrac12|\log(q(s)/q(t))|\le11a/25$.
Integrating \eqref{R-eq:tangent} against the law of the random time and using
\eqref{R-eq:time-laplace} gives $\E\mathcal R(T,b)\ge\mathcal R(a^2,b)$.
Multiply by $\sqrt{q(s)q(t)}$ and integrate.

For shorter vectors, weighted translation overlap of a log-concave density increases when the translation shrinks. This fact has a one-dimensional proof. For a log-concave line density $g$ and $\delta>0$, the ratio $g(x)/g(x+\delta)$ is nondecreasing in $x$. At the crossing $c$ of $Ag(x)$ and $Bg(x+\delta)$,
\[
 \int\min\{Ag(x),Bg(x+\delta)\}\dd x
 =A\int_{-\infty}^{c}g(x)\dd x+B\int_c^\infty g(x+\delta)\dd x.
\]
Its derivative in $\delta$ is $-Bg(c+\delta)\le0$, because the terms from the moving crossing cancel. Boundary crossings follow by limits. Approximation covers nonsmooth log-concave functions. Applying this along lines parallel to $v$ and integrating the transverse coordinate proves the claim for the cosine product. This reduces $\norm v_2<1$ to the unit case.
\end{proof}

\subsection{The fixed auxiliary density}
All coefficients in the polynomial
\begin{align}\label{R-eq:htrial}
 h(x)={}&-0.2765022959752712062x^2
 -0.0005935651463657652x^4\notag\\
 &+0.0000085009936360202x^6
 -0.0000001379669116675x^8\notag\\
 &+0.0000000033324446571x^{10}.
\end{align}
Set
\begin{equation}\label{R-eq:qtrial}
 q_\dagger(s)=\frac{e^{h(s/3)}}{\int_{-3}^3e^{h(t/3)}dt},
 \qquad |s|<3.
\end{equation}

The coefficients are exact rationals. The density is even, with
$\Lip(\log q_\dagger)\le\sum_j j|h_j|/3<22\pi/(50R)$ for
$6\le R\le7$. Its normalizer is denoted by $Z$.
For a positive $q$, the first variation is
\begin{equation}\label{R-eq:gradient}
 D_q^k(s)=\frac12\int_{-3}^3
 \Phi\left(-k|s-t|-\frac{\log q(s)-\log q(t)}{2k|s-t|}\right)dt.
\end{equation}
Homogeneity of $\Pi$ gives
\begin{equation}\label{R-eq:dual} \G_k(q)=\int q(s)D_q^k(s)\,ds.
\end{equation}
Using this identity, the rational integration in
Appendix~\ref{NEW-app:primal} proves
\begin{equation}\label{eq:primal-height}
 \G_{\pi/(2C)}(q_\dagger)
 >1.00000000000000639518571466,
 \qquad C=\Cstar.
\end{equation}
The zero-column height also exceeds one: $h(1)>-9/32$, $Z\le6$
and $e^{-9/32}>3/4$ give $q_\dagger>1/8$, hence
$\frac14\iint\min(q_\dagger(s),q_\dagger(t))\,ds\,dt>9/8$.

\subsection{The cubic remainder}
For $\lambda\ge1/2$, write
\begin{equation}\label{NEW-eq:Delta}
 \omega_\lambda=\sqrt{2\lambda-1},\qquad
 \Delta_\lambda(z)=\log J_z(\omega_\lambda)+\lambda z^2,
 \qquad 0<z\le\pi/2.
\end{equation}
The value at $\lambda=1/2$ is defined by continuity in
\eqref{R-eq:Js}. The function $\Delta_\lambda$ measures the part of the
cosine likelihood transform discarded by its quadratic bound.

\begin{lemma}[A cubic comparison for the likelihood transform]
\label{NEW-lem:cubic}
For every $\lambda\ge1/2$, $\Delta_\lambda(z)$ is nonnegative and
increasing in $z$, and $\Delta_\lambda(z)/z^3$ is decreasing on
$(0,\pi/2]$. For each fixed $z$, $\Delta_\lambda(z)$ is increasing in
$\lambda$. Consequently, if $\|u\|_2=1$ and $0<a\le\pi/2$, then
\begin{equation}\label{NEW-eq:product-profile}
 \prod_iJ_{a|u_i|}(\omega_\lambda)
 \ge e^{-\lambda a^2}
       \exp\!\left\{\Delta_\lambda(a)\sum_i|u_i|^3\right\}.
\end{equation}
Equality holds in \eqref{NEW-eq:product-profile} for a direction supported
on one coordinate.
\end{lemma}
\begin{proof}
Put
\[
 D_\omega(z)=\cot z+\omega\coth((\pi-z)\omega),\qquad
 E_\omega(z)=D_\omega(z)-\frac1z.
\]
At $\omega=0$, the second term in $D_\omega$ is $1/(\pi-z)$.
The partial-fraction identity for $\csc^2$ and $\sinh x\ge x$ give
\begin{align*}
 E_\omega'(z)
 &= -\csc^2z+\omega^2\operatorname{csch}^2((\pi-z)\omega)
       +z^{-2}\\
 &\le-\sum_{k\in\mathbb Z}(z-k\pi)^{-2}
       +(\pi-z)^{-2}+z^{-2}<0.
\end{align*}
Also
$E_\omega(\pi/2)=\omega\coth(\pi\omega/2)-2/\pi\ge0$.
Thus $E_\omega$ is nonnegative on the required interval. The logarithmic
derivative of $J$ yields
\[
 \Delta_\lambda'(z)
 =2\lambda z^2\frac{E_\omega(z)}{1+zE_\omega(z)}.
\]
The function $U(z)=E_\omega(z)/(1+zE_\omega(z))$ satisfies
\[
 U'(z)=\frac{E_\omega'(z)-E_\omega(z)^2}
              {(1+zE_\omega(z))^2}\le0.
\]
Since $\Delta_\lambda(0)=0$, integration gives
\[
 \frac{\Delta_\lambda(z)}{z^3}
 =2\lambda\int_0^1t^2U(zt)\,dt.
\]
This proves monotonicity in $z$ and nonnegativity. The function
$\omega\coth(c\omega)$ increases in $\omega\ge0$ for each $c>0$;
indeed its derivative has the sign of $\sinh(2c\omega)-2c\omega$.
Therefore $E_\omega$, $U$, and $2\lambda U$ increase with $\lambda$.
The integral formula proves the remaining monotonicity.

For each nonzero coordinate,
$\Delta_\lambda(a|u_i|)\ge |u_i|^3\Delta_\lambda(a)$.
Sum these inequalities and use $\sum_i u_i^2=1$. A one-coordinate
direction gives equality directly.
\end{proof}

Let $p=\sum_i|u_i|^3$ and $\rho=\|u\|_\infty$ for a unit vector $u$.
Since $\rho^3\le p\le\rho$, the same monotonicity gives
\begin{equation}\label{PROFILE-eq:max}
 \sum_i\Delta_\lambda(a|u_i|)
 \ge\frac p{\rho^3}\Delta_\lambda(a\rho)
 \ge\Delta_\lambda(ap^{1/3})\ge p\Delta_\lambda(a).
\end{equation}
Indeed, apply the decrease of $\Delta_\lambda(z)/z^3$ to each
$a|u_i|\le a\rho$, sum, and apply it once more at $a\rho\le ap^{1/3}$.
The first inequality is an equality when the nonzero absolute
coordinates are equal.

For $6\le R\le7$, put
\begin{equation}\label{NEW-eq:profile-parameters}
 a=\frac{\pi|s-t|}{2R},\quad
 b=\frac12\log\frac{q_\dagger(s)}{q_\dagger(t)},\quad
 \lambda(a,b)=\frac{\phi(a)e^{-b^2/(2a^2)}}{a\mathcal R(a^2,b)}.
\end{equation}
This is the exponent in Lemma~\ref{R-lem:tangent}. Define
\begin{equation}\label{PROFILE-eq:one}
 \widetilde{\mathcal B}_R(p)=\frac14\iint
 \Pi(q_\dagger(s),q_\dagger(t);a)
 \exp\{\Delta_{\lambda(a,b)}(ap^{1/3})\}\,ds\,dt.
\end{equation}
At $p=0$ this is $\G_{\pi/(2R)}(q_\dagger)$.
Integrating the exponential minorant against the time law, and then
using \eqref{PROFILE-eq:max}, proves
\begin{equation}\label{PROFILE-eq:height}
 H_{f_{R,m}}(u;q_\dagger)
 \ge\widetilde{\mathcal B}_R(S_3(u)).
\end{equation}
The finite interval estimate in Appendix~\ref{NEW-app:profile} is
\begin{equation}\label{NEW-eq:profile-certificate}
 \widetilde{\mathcal B}_{C-\gamma p}(p)
 \ge\G_{\pi/(2C)}(q_\dagger)>1
 \quad(0\le p\le1),\qquad\gamma=\frac{167}{200}.
\end{equation}
Equations \eqref{PROFILE-eq:height}--\eqref{NEW-eq:profile-certificate}
prove \eqref{eq:scalar-input}. The exact overlap remains the hypothesis
of Theorem~\ref{MAIN-signing}; retaining the first bound in
\eqref{PROFILE-eq:max} gives a finer directly evaluable height whenever
both the cubic profile and the largest coordinate are known.


\section{Scalar Gaussian and information comparisons}\label{app:reference-inequalities}
\begin{lemma}[The scalar Gaussian envelopes]\label{GAUSS-lem:scalar}
For every $R>0$,
\[
 (X_R)_i\cx N(0,\alpha R^2),\qquad T_j\cx N(0,\tau).
\]
The variance on either right-hand side cannot be decreased.
\end{lemma}
\begin{proof}
For a symmetric integrable variable $Z$, comparison with $N(0,s^2)$
requires $\E|Z|\le s\sqrt{2/\pi}$. Thus the candidate standard deviation
is $s_Z=\sqrt{\pi/2}\,\E|Z|$. We check that this candidate suffices for
both reference densities.

We use a two-crossing criterion. Suppose symmetric densities $f,g$ have
the same first absolute moment and $g-f$ has successive signs $+,-,+$
on the positive half-line. Put
$D(t)=\int_t^\infty(x-t)(g(x)-f(x))\,dx$ for $t\ge0$.
Then $D(0)=D'(0)=0$ and $D''=g-f$. On the first interval $D'$ increases
from zero. On the middle interval it decreases, and on the final
interval it increases to zero from below. Therefore $D'$ has at most
one change from positive to negative, and $D(t)\ge0$ because
$D(0)=\lim_{t\to\infty}D(t)=0$. Symmetry gives all stop-loss tests,
and hence convex order.

For the unit cosine density, $\E|Z|=1/2-2/\pi^2$.
Its density divided by the $N(0,\alpha)$ density is less than one at
zero, tends to zero at the endpoint one, and has a unique interior
maximum. To verify the last assertion, the derivative of its logarithm
is $x/\alpha-\pi\tan(\pi x/2)$; after division by $x$, this expression
strictly decreases. It starts positive. Normalization forces the ratio
to exceed one between its endpoints. Thus there are two crossings of
the required type.

For $q_\dagger$ and $N(0,\tau)$, the logarithmic density ratio has
positive derivative on $(0,3)$, since
\[
 \frac1\tau-\frac19\sum_{k\in\{2,4,6,8,10\}}k|c_k|>0,
 \qquad h(x)=\sum_kc_kx^k.
\]
The ratio is less than one at zero because
$Z>\pi\E|T_1|$. It exceeds one near three by normalization and then
drops to zero outside the support of $q_\dagger$. These are again the
required two crossings. The rational inequalities just used, together
with enclosures for the constants, are certified in
Appendix~\ref{GAUSS-app}. Matching first absolute moments proves
optimality in both cases.
\end{proof}

\subsection{The constants}\label{GAUSS-app}
The same density polynomial gives
\[
 \E|T_1|=\frac{6\int_0^1xe^{h(x)}\,dx}
 {\int_{-1}^1e^{h(x)}\,dx}.
\]
Its numerator and denominator are integrated with exact rational
coefficients and the density remainder used above. The program
\texttt{certify\_gaussian\_envelopes.py} gives
\begin{align*}
 \E|T_1|&\in(1.43211362569999658,1.43211362569999659),\\
 \alpha&\in(0.13889226438133246,0.13889226438133247),\\
 \tau&\in(3.22162384194907020,3.22162384194907021).
\end{align*}
It also certifies
$1/\tau-\sum_k k|c_k|/9>0.24868$,
$Z>\pi\E|T_1|$, $1-\pi(1/2-2/\pi^2)>0$, and
$1/\alpha-\pi^2/2>0$. These are precisely the scalar inequalities
used in Lemma~\ref{GAUSS-lem:scalar}.

\subsection{The variance proxy and the rate function}\label{BAL-app:scalar}
If an even density on $(-L,L)$ is proportional to $e^{-W(x^2)}$ with $W$ convex, its squared variable is convex-order dominated by $s^2G^2$, where $s^2$ is its variance. To prove this, the logarithm of the squared-density ratio against that Gaussian is a constant plus $y/(2s^2)-W(y)$, a concave function. The density difference has signs $-,+,-$; equal mass and first moment rule out a single crossing. A convex test lies below its secant on the middle interval and above it outside, proving the convex order. Applying the squared comparison to the convex function $y\mapsto\cosh(t\sqrt y)$ gives the optimal subgaussian proxy $s^2$.

For the auxiliary density, direct differentiation of the rational
polynomial proves convexity of $-h(\sqrt y)$ on $[0,1]$. Thus, with
$\psi(t)=\log\E e^{tT_1}$,
\begin{equation}\label{T-eq:subg}
 \psi(t)\le wt^2/2,\qquad w=\E T_1^2<2.784448.
\end{equation}
The cosine comparison is already \eqref{eq:cosine-laplace}.
\subsection{The optimized auxiliary and finite certificates}\label{BAL-app:constants}
Here $\mathfrak i(s)=d((1+s)/2\Vert1/2)$, $\psi=\log\E e^{(\cdot)T_1}$, $I=\psi^*$, $w=\E T_1^2$, and $\alpha_{\rm rate}=I(1)/\log2=\kappa^{-1}$. Use $h_0=h$ and $q=q_\dagger$ from \eqref{R-eq:htrial}--\eqref{R-eq:qtrial}.
For the scalar comparisons, set $t_*=(\psi')^{-1}(1)$. Exact integration encloses
\[
 0.3861088059<t_*<0.3861088060,\quad
 0.26<\alpha_{\rm rate}<0.28,\quad w<14/5.
\]
Let $y_0=23/25$. The series
$\ri(y)=\sum_{k\ge1}y^{2k}/(2k(2k-1))$
makes $\ri(y)/y^2$ increasing, and elementary logarithm bounds give
\[
 \ri(y_0)/y_0^2<63/100,\qquad\operatorname{artanh}(y_0)>79/50.
\]
For $0<y\le y_0$, \eqref{T-eq:subg} gives
$I(y)/y^2\ge1/(2w)>5/28>(28/100)(63/100)$.
For $y_0\le y<1$, the derivative of $I(y)-\alpha_{\rm rate}\ri(y)$ is at most
$2/5-(26/100)(79/50)<0$.
The difference vanishes at one. This proves $I(s)\ge\kappa^{-1}\mathfrak i(s)$ on $[-1,1]$ with its exact equality set.

The rational calculation integrates the degree-$26$ Taylor polynomial of
$e^{h_0(x)+3tx}$ on $[-1,1]$. For $|t|\le2/5$ its exponent has absolute value below $3/2$, so the uniform remainder is bounded by
\[
 \frac{5(3/2)^{27}}{27!}.
\]
Integration against $1,x,x^2,x^4$ encloses the normalizer, the needed tilted means and the fourth moment. Logarithms are enclosed by the series
\[
 \log z=2\sum_{k=0}^{N-1}\frac{u^{2k+1}}{2k+1}+R_N,
 \quad u=\frac{z-1}{z+1},\quad
 0\le R_N\le\frac{2u^{2N+1}}{(2N+1)(1-u^2)}
\]
after scaling to $1\le z\le2$. The implementation uses $N=55$ and outward rational rounding to mesh $10^{-21}$.

\subsection{The rational information certificate}\label{INFO-app}
For odd $N$, the polynomial
\[
 P_N(x)=\sum_{k=1}^N
 \frac{2^{2k}(2^{2k}-1)B_{2k}}{2k(2k)!}\,x^{2k}
\]
is an upper bound for $\log\cosh x$ on the entire real line. Here
$B_{2k}$ are the Bernoulli numbers. To verify the inequality, apply the
odd truncation of the alternating series for $\log(1+y)$ to each
factor of
\[
 \cosh x=\prod_{j\ge0}\left(1+
                   \frac{4x^2}{\pi^2(2j+1)^2}\right).
\]
For $y\ge0$, the odd truncation is a global upper bound, as follows by
integrating the finite geometric identity for $(1+y)^{-1}$.
The coefficient sums converge absolutely and give the displayed
Bernoulli coefficients.

Use $N=41$ and $a=39/80$. The coefficients of $h$ in
\eqref{R-eq:htrial} are rational and $\sup_{[-1,1]}|h|<3/10$.
The degree-$28$ Taylor polynomial for $e^h$ has uniform error at most
\[
 2(3/10)^{29}/29!.
\]
Integrating this polynomial against $x^{2k}$ gives rational enclosures
for all moments needed in $\E P_{41}(aT_1)$. Division by the enclosed
positive normalizer is outward rounded. The resulting lower bound is
\[
 a-\E P_{41}(aT_1)
 \in[0.206576817856581526048,\,
      0.206576817856581526049].
\]
The upper endpoint here encloses the polynomial trial, which is itself
a lower bound for $d_q$; it is not an upper bound for $d_q$.
The rational logarithm series
$\log z=2\sum_{k\ge0}u^{2k+1}/(2k+1)$, where $u=(z-1)/(z+1)$,
with its geometric tail bound certifies that the lower endpoint
exceeds $\log(61473/50000)$. The file
\texttt{code/certify\_information.py} reproduces these rational
inequalities and the auxiliary variance enclosure. The global height
certificate and this information certificate use exactly the same
polynomial $h$.


\section{Certificates and numerical constants}\label{NEW-app:certificates}
The density and radius coefficients are rational. The accompanying computational source
and endpoint records implement the following directed enclosures over
whole parameter rectangles.

\subsection{The restricted-prior comparison}\label{R-app:arithmetic}
The scalar inequality in Lemma~\ref{R-lem:tangent} is certified by \texttt{comparison\_certificate.py}, together with \texttt{rational\_bounds.py}. It starts from a $64\times16$ rational partition of
$[0,11/7]\times[0,11/25]$. A rectangle is accepted when an outward interval for the left side of \eqref{R-eq:endpoint}, divided by $a$, has strictly negative upper endpoint. Otherwise it is subdivided into four rectangles. The computation terminates after $2305$ subdivisions, with $7939$ accepted rectangles.

On rectangles with $a\le3/25$, write
\[
 R=\mathcal R(a^2,\kappa a)=1+a d.
\]
The code bounds $d$ using integral mean-value formulas for the two exponential differences and the two differences of $\Phi$. It evaluates
\[
 d\,\frac{\log(1+ad)}{ad}+\frac{\phi(a)e^{-\kappa^2/2}}{1+ad}+\kappa
\]
by its continuous extension at $a=0$. Thus it never divides an interval containing zero by $a$. On the remaining rectangles it evaluates \eqref{R-eq:endpoint} directly and divides by a strictly positive interval.

Every interval endpoint is an integer multiple of $2^{-128}$. The routines for $\exp$, $\Phi$, $\log$, and $\log(1+x)/x$ use rational Taylor coefficients and explicit remainders on their declared compact ranges. Endpoint monotonicity is used where available. The checker also verifies a strict global margin of $2\cdot10^{-6}$ for the divided endpoint expression. Its output is \texttt{results/comparison\_bounds.json}.

\subsection{Direct integration of the optimized height}\label{NEW-app:primal}
Write $s=3x$, $t=3y$, and $q_\dagger(3x)=e^{h(x)}/Z$. The homogeneity identity \eqref{R-eq:dual} gives
\[
 G(C):=\mathcal G_{\pi/(2C)}(q_\dagger)
 =\frac{\int_{-1}^1e^{h(x)}D_{q_\dagger}^{\pi/(2C)}(3x)\,dx}
        {\int_{-1}^1e^{h(x)}\,dx}.
\]
The identity reduces the weighted overlap to polynomial integration
against the fixed density.

Let $k$ be the rational number obtained by rounding the upper endpoint
of $\pi/(2C)$ upward to denominator $10^{50}$. For $x\ne y$, the normal
argument occurring in the first variation, with its sign separated at
$x=y$, is the polynomial
\begin{equation}\label{NEW-eq:normal-poly}
 z(x,y)=-3k(x-y)-\frac{h(x)-h(y)}{6k(x-y)}.
\end{equation}
Its coefficient $\ell^1$ norm is less than two. Replace $\Phi(z)$ by
\[
 \frac12+\frac1{\sqrt{2\pi}}
 \sum_{j=0}^{30}\frac{(-1)^jz^{2j+1}}{2^jj!(2j+1)}.
\]
Taylor's theorem applied to $e^{-t^2/2}$ bounds the omitted integral by
\[
 \frac1{\sqrt{2\pi}}
 \frac{2^{63}}{2^{31}31!\,63}.
\]
The integration across the sign change is exact on every monomial:
\[
 \left(\int_{-1}^x-\int_x^1\right)x^ay^b\,dy
 =\frac{2x^{a+b+1}-(1+(-1)^{b+1})x^a}{b+1}.
\]
The resulting first-variation polynomial has degree $550$.

For the density weight use $\rho(x)=\sum_{j=0}^{30}h(x)^j/j!$.
The rational inequalities $-1/3<h\le0$ give
$|e^h-\rho|\le(1/3)^{31}/31!$. Multiplying the two polynomials and
integrating every even monomial exactly gives the numerator. The
normalizer is enclosed in the same way. Fixed-point polynomial
operations use scale $2^{224}$ with an explicit coefficient-norm
error; the final interval arithmetic is directed at scale $2^{128}$.
If $p,q$ are integer coefficient arrays with error bounds
$e_p/S,e_q/S$ at scale $S=2^{224}$, multiplication has error at most
\[
 S^{-1}\left\lceil\frac{\|p\|_1e_q+\|q\|_1e_p+e_pe_q}{S}\right\rceil
 +\frac{\#\{\text{output coefficients}\}}{S}.
\]
Each output coefficient is rounded once. Rational scaling and exact
monomial integration propagate the analogous coefficient-norm bounds.
A one-sided allowance of $9|k-\pi/(2C)|$ accounts for rounding $k$.
The constant nine follows from
$|\partial_k\mathcal G_k(q)|\le\frac12\phi(0)\iint
\sqrt{q(s)q(t)}|s-t|\,dsdt<9$.

For $C=\Cstar$, the program \texttt{certify\_universal.py} returns
\begin{align*}
 1.00000000000000639518571466107790&\le G(C),\\
 G(C)&\le1.00000000000000639518571469690118,\\
 5.4894231280195092062761207232756114&\le Z,\\
 Z&\le5.4894231280195092062761207232756115.
\end{align*}
The stored rational endpoints determine the enclosures; the displayed
decimals are outward roundings of those endpoints.

\subsection{A uniform profile certificate}\label{NEW-app:profile}
Put $R=C-\gamma p$, $\gamma=167/200$. On a parameter interval
$P=[p_-,p_+]\subset[0,1]$, let $c_P\ge p_+^{1/3}$ be the outward
upper endpoint of its cube-root enclosure, and define
\[
 D_P(a,b)=\frac{\Delta_{\lambda(a,b)}(ac_P)}{c_P^3}.
\]
The monotonicity in Lemma~\ref{NEW-lem:cubic} gives
$\Delta_{\lambda(a,b)}(ap^{1/3})\ge pD_P(a,b)$ for every $p\in P$.
Since $D_P\ge0$ and $e^x\ge\sum_{j=0}^4x^j/j!$ for $x\ge0$, it
suffices to prove
\begin{equation}\label{NEW-eq:certified-integral}
 \frac14\iint\Pi(q(s),q(t);a)
 \left(D_P+\frac p2D_P^2+\frac{p^2}{6}D_P^3+
                    \frac{p^3}{24}D_P^4\right)dsdt
 -\gamma\int_0^1G'(R+t(C-R))\,dt>0.
\end{equation}
For $p>0$, multiplication by $p$ then gives
$\widetilde{\mathcal B}_R(p)>G(C)$; at $p=0$ the two heights agree.

Use $d=|s-t|$ and $x=(s+t)/2=(3-d/2)y$, with $0\le y\le1$.
Symmetry removes the factor $1/4$ and leaves weight $w=3-d/2$ on
$0\le d\le6$, $0\le y\le1$. With $v=d/2$ write
\begin{align*}
 H_0&=\frac{h((x+v)/3)+h((x-v)/3)}2
 =\sum_k\frac{c_k}{3^k}\sum_{\ell\ \mathrm{even}}
          \binom{k}{\ell}x^{k-\ell}v^\ell,\\
 P_0&=\frac{h((x+v)/3)-h((x-v)/3)}d
 =\sum_k\frac{c_k}{3^k}\sum_{\ell\ \mathrm{odd}}
          \binom{k}{\ell}x^{k-\ell}v^{\ell-1}.
\end{align*}
Thus $\sqrt{q(s)q(t)}=e^{H_0}/Z$, $b=dP_0/2$ up to an immaterial
sign, and $a=\pi d/(2R)$. These polynomial formulas remove division
by $d$ from the prior slope. Apart from $we^{H_0}/Z$, the integrand of
$G'(R)$ is
\begin{equation}\label{NEW-eq:gprime}
 g_R(d,y)=\frac{2a}{R}\phi(a)e^{-b^2/(2a^2)}.
\end{equation}

We bound its radius integral by Simpson quadrature with an explicit
fourth-derivative remainder. As a function of $r$, the kernel is
$c r^{-2}e^{-A/r^2-Br^2}$. Set $X=A/r^2=a^2/2$ and
$Y=Br^2=b^2/(2a^2)$. The polynomials
$P_k=r^kg_r^{(k)}/g_r$ satisfy
\[
 P_0=1,\qquad
 P_{k+1}=(-2+2X-2Y-k)P_k-2X\partial_XP_k+2Y\partial_YP_k.
\]
On $0\le X\le5/4$, $0\le Y\le1/10$, the bidegree-$(4,4)$ Bernstein
coefficients of $P_4$ lie between
\[
 -\frac{771487}{5000}\quad\text{and}\quad\frac{82886}{625}.
\]
The exact rational program \texttt{check\_second\_pass.py} verifies
these coefficients. In particular $|P_4|<155$. The stated rectangle
contains every kernel in the certificate because $r\ge r_0=C-\gamma>6$
and $|b|/a<11/25$. Hence the radius average is at most
\begin{equation}\label{PROFILE-eq:simpson}
 \frac{g_R+4g_{(R+C)/2}+g_C}{6}
 +\frac{155(C-R)^4}{2880r_0^4}
                       \max_{R\le r\le C}g_r.
\end{equation}
All quantities on the right are enclosed by directed interval arithmetic.

The program \texttt{certify\_profile\_graded.cpp} integrates the
remaining two variables by midpoint quadrature with interval second
derivatives. If $F_P(p,d,y)$ is the integrand in
\eqref{NEW-eq:certified-integral} with the Simpson term replacing the
radius average, a rectangle of side lengths $h_d,h_y$ has average
quadrature error at most
\[
 \frac{h_d^2}{24}\sup|\partial_d^2 F_P|
       +\frac{h_y^2}{24}\sup|\partial_y^2 F_P|.
\]
The program propagates values and first and second derivatives using
the product and chain rules, with outward intervals at every operation.
For the parameter variable, it evaluates the spatial integral at the
midpoint of $P$, encloses the integral of $\partial_pF_P$ over all of
$P$, and subtracts half the length of $P$ times the largest absolute
endpoint of this integrated derivative. Thus it certifies every
parameter in the interval, while retaining cancellation in the
parameter derivative. The Simpson remainder in
\eqref{PROFILE-eq:simpson} is subtracted separately.

For the spatial mesh take $d_0=1/32$. Between $d_0$ and $1/2$ successive
nodes are multiplied by $1+8/n_d$, with the last node set to $1/2$;
use $n_d$ equal intervals on $[1/2,6]$ and $n_y$ equal intervals on
$[0,1]$. The initial run uses $128$ equal parameter intervals and
$(n_d,n_y)=(128,24)$. Each of the final three parameter intervals is
split into four equal pieces and evaluated with $(n_d,n_y)=(192,32)$.
The resulting $137$ certified intervals partition $[0,1]$ exactly.

On $0\le d<d_0$, discard the positive transform terms. Since
$q\le1/Z$, the omitted negative contribution is bounded by
\[
 \frac{\gamma\phi(0)\pi}{Zr_0^2}\frac{3d_0^2}{2}.
\]
All $137$ final intervals have a positive lower endpoint; the smallest
exceeds $1.0189937442658\cdot10^{-5}$ on $[31/32,125/128]$.
Directed MPFR arithmetic at $96$ bits is used for normal tails,
trigonometric functions and cube roots, with every argument in its
stated domain. The stated interval endpoints and exact coverage prove
\eqref{NEW-eq:certified-integral} on the full parameter range.

\subsection{The Gaussian-height certificate}\label{R7-app:dual-height}

Set $C_-=6.8383231851$, $C_+=6.8383231852$ and
$k_\pm=\pi/(2C_\pm)$. Use the rational enclosures
\[
 \underline k_-=
 \frac{459409795143201068077}{2000000000000000000000}<k_-,\qquad
 \overline k_+=
 \frac{2297048975682414518287}{10000000000000000000000}>k_+.
\]
The fixed trial of \eqref{R-eq:htrial} has, for all $s\in[-3,3]$,
\begin{align}
 0.9999999999940804374730475
 &\le D_{q_\dagger}^{\underline k_-}(s)
 \le0.9999999999985087548717213,\label{R-eq:certlower}\\
 1.0000000000016799598810746
 &\le D_{q_\dagger}^{\overline k_+}(s)
 \le1.0000000000061276852238026.\label{R-eq:certupper}
\end{align}
These enclose the first variation on the entire interval. The
following rational polynomial calculation proves both inequalities.

\begin{proof}[The whole-interval dual bounds]
The polynomial \eqref{R-eq:htrial} has rational coefficients, so the slope bound follows by summing $j|h_j|/3$. The normalizing integral in \eqref{R-eq:qtrial} cancels from \eqref{R-eq:gradient}. Put $s=3x,t=3y$ and define the polynomial
\[
 z(x,y)=-3k(x-y)-\frac{h(x)-h(y)}{6k(x-y)}.
\]
The divided difference is polynomial. On $y<x$ it is the Gaussian argument in \eqref{R-eq:gradient}, while on $y>x$ that argument is $-z$. Therefore
\begin{equation}\label{R-eq:Dpoly}
 D_{q_\dagger}^k(3x)=\frac32\left\{
 \int_{-1}^x\Phi(z(x,y))\dd y+
 \int_x^1\Phi(-z(x,y))\dd y\right\}.
\end{equation}
At both rational choices of $k$, the sum of the absolute coefficients of $z$ is less than two. Replace $\Phi(z)$ by
\[
 \frac12+\phi(0)\sum_{j=0}^{20}\frac{(-1)^jz^{2j+1}}{2^j j!(2j+1)}.
\]
Taylor's remainder for $e^{-u}$ with $u\ge0$, followed by integration, bounds the resulting error in \eqref{R-eq:Dpoly} by
\begin{equation}\label{R-eq:taylortail}
 \frac32\frac{2^{43}}{2^{21}21!\,43}.
\end{equation}
Here $\phi(0)<1/2$ has been used. Signed integration of a monomial is exact:
\[
 \left(\int_{-1}^x-\int_x^1\right)x^a y^b\dd y
 =\frac{2x^{a+b+1}-(1+(-1)^{b+1})x^a}{b+1}.
\]
It produces a univariate polynomial of degree $370$. Convert it to the Chebyshev basis $\sum c_jT_j(x)$. Since $|T_j(x)|\le1$ on $[-1,1]$, its entire range is enclosed by
$c_0\pm\sum_{j\ge1}|c_j|$, enlarged by \eqref{R-eq:taylortail} and the arithmetic errors.

The code performs the bivariate polynomial arithmetic with $224$-bit fixed-point coefficients and a separately propagated coefficient-$\ell^1$ error. If $P,Q$ have integer coefficient arrays $p,q$ and errors $e_p/S,e_q/S$, where $S=2^{224}$, the error after multiplication is at most
\[
 S^{-1}\left\lceil\frac{\|p\|_1e_q+\|q\|_1e_p+e_pe_q}{S}\right\rceil
 +\frac{\#\{\text{output coefficients}\}}S.
\]
Each output coefficient is rounded once. Rational scaling and signed integration carry analogous explicit rounding bounds. Chebyshev conversion and the constants $\pi,\phi(0)$ use outward $128$-bit dyadic intervals. The former constant is enclosed through Machin's arctangent formula and rational remainders; the latter uses integer square roots. These operations give \eqref{R-eq:certlower}--\eqref{R-eq:certupper} without evaluating a grid or running an optimizer.

Inequality~\eqref{R7-eq:dual-all-q} gives $M(\underline k_-)<1$ and $\G_{\overline k_+}(q_\dagger)>1+10^{-12}$. Monotonicity in $k$, rather than any monotonicity of the individual derivative $D_q^k(s)$, transfers these two bounds to $k_-$ and $k_+$. This proves the two inequalities used in Theorem~\ref{R7-thm:method-constant}.
\end{proof}

The code needs no approximation to the normalizer of $q_\dagger$. Its global upper certificate applies to densities with zeros, unbounded logarithmic derivatives, or arbitrarily narrow peaks. The fixed rational trial was selected numerically; its selection plays no role in the proof of the displayed inequalities.


\section{Uniform zero-bias analysis of diffuse cosine scores}
\label{S2-app:diffuse}
This appendix proves the uniform expansion used in
\eqref{S2-eq:diffuse} and the volume argument.  Let $T$ have density
$2/[\pi(1+t^2)^2]$.  For a unit vector $w$, put
\[
 a_i=|w_i|,\qquad S_w=\sum_i a_iT_i,
 \qquad \beta=\sum_i a_i^3,
 \qquad H=-2\log\beta.
\]
The regime $\max_i a_i\to0$ is equivalent to $\beta\to0$.

\subsection{The logarithmic zero-bias formula}
Let $X$ have a positive symmetric density $p$, mean zero, and variance
$\sigma^2\in(0,\infty)$.  Put
\[
 M(x)=\int_x^\infty yp(y)\,dy,
 \qquad p^*(x)=M(x)/\sigma^2.
\]
Then $p^*$ is the zero-bias density, characterized by
$\E[X\phi(X)]=\sigma^2\E[\phi'(X^*)]$.

\begin{proposition}[Logarithmic radius formula]
\label{S2D-prop:log}
With infinite values allowed,
\begin{equation}\label{S2D-eq:log}
 r(\law X)=\int_0^\infty
 \frac{\log p^*(0)-\log p^*(x)}{x^2}\,dx.
\end{equation}
\end{proposition}
\begin{proof}
For $t=\Prb(X>x)$, symmetry and the quantile definition
\eqref{S-eq:radius} give $L_{\law X}(t)=M(x)$ and $dt=-p(x)dx$.
Hence
\[
 r(\law X)=\int_0^\infty\frac{p(x)}{M(x)}\,dx.
\]
Since $M'(x)=-xp(x)$,
\[
 \log\frac{M(0)}{M(x)}=\int_0^x\frac{tp(t)}{M(t)}\,dt.
\]
Tonelli's theorem and $\int_t^\infty x^{-2}dx=t^{-1}$ give
\eqref{S2D-eq:log}.
\end{proof}

Direct integration shows that the zero-bias law of $T$ is standard
Cauchy.  If $C$ is standard Cauchy, independent of the $T_i$, and
$I$ is independent with $\Prb(I=i)=a_i^2$, then the independent-sum
zero-bias identity gives
\begin{equation}\label{S2D-eq:cauchy}
 S_w^*\stackrel d=\sum_{j\ne I}a_jT_j+a_IC.
\end{equation}
Writing $M_w$ for the density of $S_w^*$ and
$L_1=\sum_i a_i$, its Fourier transform is therefore
\begin{equation}\label{S2D-eq:cf}
 \widehat M_w(t)=e^{-L_1|t|}
 \sum_i a_i^2\prod_{j\ne i}(1+a_j|t|).
\end{equation}
For equal coefficients this reads
\[
 S_{u_k}^*\stackrel d=k^{-1/2}(C+T_1+\cdots+T_{k-1}).
\]

We shall also need a uniform tail floor.
\begin{lemma}\label{S2D-lem:global}
For $x\ge2$,
\begin{equation}\label{S2D-eq:global}
 M_w(x)\ge\frac{3\beta}{17\pi x^2}.
\end{equation}
Moreover $M_w(0)\le1/2$, and therefore for an absolute constant $C$,
\begin{equation}\label{S2D-eq:global-log}
 0\le\log\frac{M_w(0)}{M_w(x)}
 \le\log(1/\beta)+2\log x+C.
\end{equation}
\end{lemma}
\begin{proof}
Condition on $I=i$ in \eqref{S2D-eq:cauchy} and put
$Y_i=\sum_{j\ne i}a_jT_j$.  Since $\E Y_i^2\le1$,
$\Prb(|Y_i|\le x)\ge3/4$ for $x\ge2$.  On that event the conditional
Cauchy density at $x$ is at least $4a_i/(17\pi x^2)$.  Multiplication by
$\Prb(I=i)=a_i^2$ and summation proves \eqref{S2D-eq:global}.
Also $M_w(0)=\E|S_w|/2\le1/2$.  The last claim follows by taking
logarithms and using monotonicity of $M_w$.
\end{proof}

\subsection{A profile-uniform Fourier expansion}
Let $\varphi(x)=(2\pi)^{-1/2}e^{-x^2/2}$ and define
\begin{equation}\label{S2D-eq:h}
 \mathfrak h(x)=\frac1\pi\int_0^\infty e^{-t^2/2}
 \left(\frac{t^3}{3}-t\right)\cos(tx)\,dt.
\end{equation}
Put $\gamma=\max_i a_i$ and $\delta_4=\sum_i a_i^4$.

\begin{lemma}[Uniform central and tail expansion]\label{S2D-lem:fourier}
For all sufficiently small $\gamma$,
\begin{equation}\label{S2D-eq:fourier}
 \|M_w-\varphi-\beta\mathfrak h\|_\infty
 +\|M_w''-\varphi''-\beta\mathfrak h''\|_\infty
 \le C\delta_4\le C\beta^{4/3}.
\end{equation}
Furthermore
\begin{equation}\label{S2D-eq:h-tail}
 \mathfrak h(x)=\frac1{\pi x^2}
 \left(1+\frac5{x^2}+O(x^{-4})\right),
 \qquad \mathfrak h(0)=-\frac1{3\pi}.
\end{equation}
\end{lemma}
\begin{proof}
For $t\ge0$ write
\[
 P(t)=\prod_i(1+a_it)e^{-a_it},\qquad
 B(t)=\sum_i\frac{a_i^2}{1+a_it},
\]
so that \eqref{S2D-eq:cf} is $\widehat M_w=PB$.  Taylor expansion with
integral remainder gives
\begin{align*}
 \log P(t)&=-t^2/2+\beta t^3/3+R(t),
 &|R(t)|&\le\delta_4t^4/4,\quad R(t)\le0,\\
 B(t)&=1-\beta t+R_B(t),
 &0\le R_B(t)&\le\delta_4t^2.
\end{align*}
Also
\[
 P(t)\le\exp\left[-\frac{t^2}{2(1+\gamma t)}\right],
 \qquad 0<B(t)\le1.
\]
For $0\le t\le(2\gamma)^{-1}$ these inequalities yield
\[
 \left|P(t)B(t)-e^{-t^2/2}
 \{1+\beta(t^3/3-t)\}\right|
 \le C\delta_4(t^2+t^4+t^6)e^{-t^2/4}.
\]
Here $\beta\le\gamma$ and $\beta^2\le\delta_4$; the latter follows
from Cauchy--Schwarz with weights $a_i^2$.  On
$[(2\gamma)^{-1},\gamma^{-1}]$ the exponential bound is at most
$e^{-t^2/4}$, while for $t\ge\gamma^{-1}$ it is at most
$e^{-t/(4\gamma)}$.  The tails, with or without a factor $t^2$, are
$O(\delta_4)$.  Fourier inversion proves the first inequality in
\eqref{S2D-eq:fourier}.  Finally
$\delta_4\le\gamma\beta\le\beta^{4/3}$.

Repeated integration by parts in \eqref{S2D-eq:h} gives
\eqref{S2D-eq:h-tail}.  The endpoint expansion
\[
 e^{-t^2/2}(t^3/3-t)=-t+\frac56t^3-\frac7{24}t^5+O(t^7)
\]
produces the first reciprocal powers, while direct integration at zero
gives $\mathfrak h(0)=-1/(3\pi)$.
\end{proof}

\subsection{The transition scale}
Let $b_\beta>\sqrt2$ be the larger solution of
\begin{equation}\label{S2D-eq:b}
 b_\beta^2-4\log b_\beta=H+\log(\pi/2).
\end{equation}
For all sufficiently small $\beta$ it is unique.  The next theorem is
uniform over the entire coefficient profile.

\begin{theorem}[Sharp diffuse expansion]\label{S2D-thm:sharp}
As $\max_i|w_i|\to0$,
\begin{align}
 r(\nu_w)&=b_\beta+\frac2{b_\beta}+O(b_\beta^{-3}),
 \label{S2D-eq:r-b}\\
 r(\nu_w)^2
 &=H+2\log H+4+\log(\pi/2)
 +O\left(\frac{\log H}{H}\right).
 \label{S2D-eq:entropy}
\end{align}
The constants in the error terms are independent of dimension and of
$w$.
\end{theorem}
\begin{proof}
Put $L=\log(1/\beta)$, $b=b_\beta$, $X=\beta^{-1/12}$, and
$F_w(x)=\log(M_w(0)/M_w(x))$.  Proposition~\ref{S2D-prop:log} gives
$r(\nu_w)=\int_0^\infty F_w(x)x^{-2}dx$, while
$b^2=2L+O(\log L)$.

On $[0,1]$, the second-derivative part of
Lemma~\ref{S2D-lem:fourier}, evenness, and uniform positivity on compact
intervals give
\[
 F_w(x)-x^2/2=O(\beta x^2).
\]
On $[1,b/2]$ the same lemma compares $M_w$ with $\varphi$ with relative
error at most $C\beta e^{b^2/8}$, which tends to zero faster than any
inverse power needed below.  These two intervals contribute
$o(b^{-3})$ to the difference from the Gaussian model.

On $[b/2,X]$ use \eqref{S2D-eq:h-tail} and set
\[
 M_0(x)=\varphi(x)+\frac{\beta}{\pi x^2}.
\]
Then
\[
 \frac{|M_w(x)-M_0(x)|}{M_0(x)}
 \le Cb^{-2}+C\frac{\delta_4X^2}{\beta}
 \le Cb^{-2}+C\beta^{1/6},
\]
Moreover,
\[
 \log M_w(0)=\log\varphi(0)+O(\beta).
\]
After integration against $x^{-2}$, replacing $F_w$ by
$\log\{\varphi(0)/M_0(x)\}$ costs $O(b^{-3})$.  For $x\ge X$,
Lemma~\ref{S2D-lem:global} bounds both tails by
$O((L+\log X+1)/X)=o(b^{-3})$.

Put $c=\tfrac12\log(\pi/2)$ and
\[
 A(x)=x^2/2,\qquad B(x)=L+2\log x+c.
\]
The two functions meet at $x=b$.  On $[b/2,\infty)$,
\[
 \log\frac{\varphi(0)}{M_0(x)}
 =\min\{A(x),B(x)\}
 -\log(1+e^{-|A(x)-B(x)|}).
\]
For small $\beta$, $A-B$ is increasing there with derivative at least
$x/2$.  Changing variables to $A-B$ shows that the final logarithmic
term contributes $O(b^{-3})$, because
$\int_\R\log(1+e^{-|u|})du<\infty$.  Hence
\begin{align*}
 r(\nu_w)
 &=\int_0^b\frac12\,dx
 +\int_b^\infty\frac{L+2\log x+c}{x^2}\,dx+O(b^{-3})\\
 &=\frac b2+\frac{L+2\log b+c+2}{b}+O(b^{-3})
 =b+\frac2b+O(b^{-3}),
\end{align*}
using $b^2/2=L+2\log b+c$.  This proves \eqref{S2D-eq:r-b}.
Writing $y=b^2$, equation \eqref{S2D-eq:b} becomes
$y-2\log y=H+\log(\pi/2)$, whence
\[
 y=H+2\log H+\log(\pi/2)
 +O((\log H)/H).
\]
Squaring \eqref{S2D-eq:r-b} adds $4+O(1/H)$ and proves
\eqref{S2D-eq:entropy}.
\end{proof}

For $u_k=k^{-1/2}(1,\ldots,1)$, one has
$\beta=k^{-1/2}$ and $H=\log k$, so
Theorem~\ref{S2D-thm:sharp} is exactly \eqref{S2-eq:diffuse}.  The
uniformity in the coefficient profile is what justifies the spherical
volume argument leading to \eqref{S2-eq:volume}.




\begingroup
\makeatletter\let\addcontentsline\@gobblethree\makeatother
\pdfbookmark[1]{Acknowledgment}{R27-acknowledgment}

\section*{Acknowledgment}
The author used ChatGPT for mathematical exploration and proof development, literature searches, computational experiments, and code for numerical certificates. Proofs developed with this assistance were carefully checked, adapted, and rewritten by the author, who also rewrote the exposition. The author takes full responsibility for the mathematical content, computations, and any errors.


\endgroup
\end{document}